\documentclass{article}

\usepackage{tipa}

\usepackage{times}
\usepackage{graphicx} 
\usepackage{subfigure}

\usepackage{natbib}

\usepackage{algorithm}
\usepackage{algorithmic}
\usepackage[algo2e,lined,boxed,vlined]{algorithm2e}

\usepackage{amsfonts}
\usepackage{amssymb}
\usepackage{booktabs}

\usepackage{hyperref}
\usepackage{xurl}

\usepackage[accepted]{icml2016}

\usepackage{comment}
\usepackage{amsthm}
\usepackage{mathtools}
\usepackage{tabularx}
\usepackage{multirow}

\newtheorem{lemma}{Lemma}
\newtheorem{proposition}{Proposition}
\newtheorem{corollary}{Corollary}
\newtheorem{definition}{Definition}

\newtheorem{remark}{Remark}

\def\b{\ensuremath\boldsymbol}

\usepackage{cancel}

\usepackage{lastpage}
\usepackage{fancyhdr}
\usepackage{url}
\icmltitlerunning{Foundations of Riemannian Geometry for Riemannian Optimization: A Monograph with Detailed Derivations}

\begin{document}


\twocolumn[
\icmltitle{Foundations of Riemannian Geometry for Riemannian Optimization: \\A Monograph with Detailed Derivations}


\icmlauthor{Benyamin Ghojogh}{bghojogh@uwaterloo.ca}
\icmladdress{Waterloo, Ontario, Canada}


\vskip 0.3in
]

\begin{abstract}
Riemannian geometry provides the fundamental framework for optimization on nonlinear spaces such as matrix manifolds, which arise in machine learning, signal processing, and robotics. While the underlying theory is classical, existing literature often presents results at a high level of abstraction, omitting the detailed coordinate-level derivations required for implementation and algorithm development.

This work provides a self-contained and rigorous treatment of the foundations of Riemannian geometry, with a focus on explicit derivations tailored to Riemannian optimization. We systematically develop the key geometric structures---including tangent and cotangent spaces, tensor calculus, metric tensors, Levi-Civita connections, curvature, and geodesics---emphasizing step-by-step derivations in coordinates and matrix form.

Building on these foundations, we derive the Riemannian gradient, Hessian, exponential map, and retraction in a form suitable for numerical computation. We further specialize these constructions to important matrix manifolds, including the Stiefel, Grassmann, and SPD (Symmetric Positive Definite) manifolds, providing explicit formulas widely used in optimization and geometric machine learning.

This monograph develops a unified and implementation-oriented treatment of Riemannian geometry for optimization on manifolds. Its main contribution is the systematic organization and detailed derivation of classical geometric constructions in forms directly usable for algorithm design and numerical implementation. By connecting coordinate-level differential geometry with matrix-manifold formulas, the monograph bridges the gap between abstract theory and practical computation, and provides a reference for researchers and practitioners working in Riemannian optimization and related fields.
\end{abstract}

\hfill\break
\hfill\break

\hfill\break
\hfill\break







\onecolumn
{\small
\tableofcontents
}
\twocolumn

\section{Introduction}

\subsection{Motivation and Purpose of This Monograph}

Optimization on nonlinear spaces has become increasingly important in modern scientific and engineering applications, including machine learning, signal processing, computer vision, and robotics. In many of these problems, the variables naturally lie on smooth manifolds rather than Euclidean spaces. Examples include orthogonality constraints (Stiefel manifold), subspace representations (Grassmann manifold), and positive definite matrices (SPD manifold). Riemannian geometry provides the appropriate mathematical framework for handling such problems.

Over the past decades, a substantial body of literature has developed around Riemannian manifolds and optimization on manifolds. Foundational treatments of differential geometry can be found in classical texts such as \cite{lee2013smooth, do1992riemannian}, while modern treatments of Riemannian optimization (or optimization on smooth manifolds) are presented in works such as \cite{absil2008optimization, boumal2023introduction}. These references provide deep theoretical insights and elegant formulations of geometric concepts.

However, a practical gap remains between abstract theory and implementation. Many results are presented in a coordinate-free or highly compact form, which, while mathematically elegant, can obscure the detailed steps required for numerical implementation. In particular, explicit coordinate-level derivations of key objects—such as the Levi-Civita connection, Riemannian gradient, and Riemannian Hessian—are often omitted or scattered across different sources. This creates a barrier for researchers and practitioners who seek to develop algorithms or implement Riemannian optimization methods from first principles.

The goal of this work is to bridge this gap by providing a self-contained and detailed exposition of the foundations of Riemannian geometry, with an emphasis on explicit derivations suitable for computation. Rather than introducing new theoretical results, we focus on systematically deriving classical concepts in a way that makes every step transparent and implementable.

We begin by developing the fundamental structures of differential geometry, including smooth manifolds, tangent and cotangent spaces, tensor calculus, and coordinate systems. We then introduce the metric tensor and derive the Levi-Civita connection and curvature in explicit coordinate form. Building on these foundations, we derive key tools for Riemannian optimization, including the Riemannian gradient, Riemannian Hessian, exponential and logarithm maps, retractions, and vector transport.

Finally, we specialize these constructions to important matrix manifolds, including the Stiefel and Grassmann manifolds. For these manifolds, we provide explicit formulas for tangent spaces, projections, gradients, Hessians, and geometric operations, which are widely used in optimization and geometric machine learning.

This monograph is intended as an expository and implementation-oriented reference. Its contribution is the systematic consolidation of classical Riemannian-geometric constructions, with detailed coordinate-level and matrix-level derivations, for researchers implementing optimization algorithms on matrix manifolds.
By consolidating these derivations into a single reference, this work aims to serve as a practical bridge between differential geometry and modern computational applications.

\subsection{Notational Conventions in This Monograph}

Before entering the technical development, we summarize the
main notational conventions used throughout this monograph.
Because the presentation moves from general differential
geometry to Riemannian matrix-valued optimization, some
symbols are reused in increasingly specialized settings.
For example, in the general sections, $\mathcal{M}$ denotes a
smooth manifold and $\b{p}\in\mathcal{M}$ denotes a point on
the manifold, whereas in the matrix-manifold sections, the
point is often represented by a matrix such as $\b{X}$.
Likewise, tangent vectors are denoted in abstract form by
$\b{\xi}\in T_{\b{p}}\mathcal{M}$ and in matrix form by
$\b{\Delta}\in T_{\b{X}}\mathcal{M}$.

Table~\ref{table_notations} collects the most frequently used
symbols and their meanings. Its purpose is to make the
monograph easier to read by giving the reader a single place
to check notation before and during the derivations. Whenever
a symbol acquires a more specialized meaning in a later
section, that specialization is stated again locally in the
text.

\begin{table*}[!t]
\centering
\scriptsize
\caption{Notational conventions in this monograph}
\label{table_notations}
\renewcommand{\arraystretch}{1.12}
\setlength{\tabcolsep}{3pt}

\newcommand{\symwa}{0.12\linewidth}
\newcommand{\deswa}{0.35\linewidth}
\newcommand{\symwb}{0.12\linewidth}
\newcommand{\deswb}{0.35\linewidth}

\begin{tabular}{p{\symwa} p{\deswa} | p{\symwb} p{\deswb}}
\toprule
\textbf{Symbol$\quad\quad\quad\quad$ (general)} & \textbf{Description} & \textbf{Symbol$\quad\quad\quad\quad$ (matrix manifold)} & \textbf{Description} \\
\midrule

$\tau$ & Topology 
&  & \\

\hline
$\mathcal{M}$ & A smooth manifold or Riemannian manifold. 
& $\mathrm{St}(n,d)$ & Stiefel manifold of $n\times d$ matrices with orthonormal columns. \\

\hline
&
& $\mathrm{Gr}(n,d)$ & Grassmann manifold of $d$-dimensional subspaces of $\mathbb{R}^n$. \\

\hline
&
& $\mathbb{S}_{++}^n$ & Symmetric positive definite (SPD) manifold. \\

\hline
&
& $\mathbb{S}^n$ & Vector space of $n\times n$ symmetric matrices. \\

\hline
$\b{p}, \b{q}$ & Point on a general manifold $\mathcal{M}$. 
& $\b{X},\b{Y},\b{Z}$ & Matrix-valued point on matrix manifolds or matrices in the ambient space (in Section~\ref{section_important_Riemannian_matrix_manifolds}). \\

\hline
 & 
& $[\b{X}]$ & Equivalence class of $\b{X}$ (i.e., a point) in the Grassmann manifold. \\

\hline
$\b{V}, \b{W}, \b{X}, \b{Y}, \b{T}$ & Vectors, tensors, or vector fields in the general differential-geometric sections, depending on context (especially up to Section~\ref{section_important_Riemannian_matrix_manifolds}).  
&  &  \\

\hline
$\b{\xi}, \b{\eta}$ & Tangent vector at a general manifold $\mathcal{M}$.
& $\b{\Delta},\b{\Delta}_1,\b{\Delta}_2,\b{\Xi}$ & Tangent vectors (matrices) at a matrix-manifold point. \\

\hline
$T_{\b{p}}\mathcal{M}$ & Tangent space of $\mathcal{M}$ at $\b{p}$.
& $T_{\b{X}}\mathcal{M}$ & Tangent space at $\b{X}$ on a matrix manifold. \\

\hline
$T_{\b{p}}^*\mathcal{M}$ & Cotangent space of $\mathcal{M}$ at $\b{p}$. 
&  &  \\

\hline
$x^i$ & Local coordinates on the manifold. 
&  &  \\

\hline
$\partial_i,\b{e}_i$ & Coordinate basis vectors of the tangent space. 
&  &  \\

\hline
$dx^i$ & Dual coordinate basis covectors. 
&  &  \\

\hline
$V^i, V_i$ & Contravariant and covariant components of a vector/tensor $\b{V}$. 
&  &  \\

\hline
&   
& $\Pi_{\b{X}}^{\mathrm{St}}(\cdot)$ & Projection onto the tangent space of Stiefel manifold at $\b{X}$. \\

\hline
 &  
& $\Pi_{[\b{X}]}^{\mathrm{Gr}}(\cdot)$ & Projection onto the tangent space of Grassmann manifold at $[\b{X}]$. \\

\hline
$\b{X}(f)(\b{p}) = Df(\b{p})[\b{X}(\b{p})] $ & Directional derivative of the function $f$ along the vector field $\b{X}$ at point $\b{p}$. 
& $Df(\b{X})[\b{\Delta}] = D\bar f(\b{X})[\b{\Delta}]$ & ambient directional derivative of function $f$, or smooth local extension $\bar{f}$, at $\b{X}$ along (in the direction of) vector $\b{\Delta}$. \\

\hline
$D\b{Y}(\b{p})[\b{X}(\b{p})]$ & Coordinate directional derivative of the vector field $\b{Y}$ along (in the direction of) the vector field $\b{X}$ at
$p$.
&  $D\b{\Delta}_2(\b{X})[\b{\Delta}_1]$ & Ambient (Euclidean) directional derivative of vector $\b{\Delta}_2$ along (in the direction of) the vector $\b{\Delta}_1$ at point $\b{X}$ \\

\hline
$g_{\b{p}}$ & Riemannian metric at point $\b{p}$.
& $g^E_{\b{X}}, g^E_{[\b{X}]}$ & Euclidean metric on a matrix manifold. \\

\hline
$g_{ij}$ & Components of the metric tensor in local coordinates. 
& $g^C_{\b{X}}$ & Canonical metric on the Stiefel manifold. \\

\hline
$g^{ij}$ & Components of the inverse metric tensor. 
& $g^\alpha_{\b{X}}$ & $\alpha$-metric on the Stiefel manifold. \\

\hline
&
& $g^{AI}_{\b{X}}$ & Affine-invariant metric on the SPD manifold. \\

\hline
&
& $g^{LE}_{\b{X}}$ & Log-Euclidean metric on the SPD manifold. \\

\hline
&
& $g^{BW}_{\b{X}}$ & Bures--Wasserstein metric on the SPD manifold. \\

\hline
&
& $\operatorname{qf}(\cdot)$ & The $Q$ factor of QR decomposition. \\

\hline
$\Gamma^{k}_{ij}$ & Christoffel symbols of the second kind. 
&  &  \\

\hline
$\Gamma_{ijk}$ & Christoffel symbols of the first kind. 
&  &  \\

\hline
$\partial_i$ & Partial derivative with respect to coordinate $x^i$.
& & \\

\hline
$\nabla$ & Connection or covariant derivative operator.  
& $\nabla \bar f(\b{X})$  & Euclidean derivative of the smooth local extension of function $f$, at point $\b{X}$, on the matrix submanifold. \\

\hline
$\nabla_i$ & Covariant derivative with respect to $\partial_i$ or $\b{e}_i$.
& & \\

\hline
$\nabla_{\b{X}} \b{Y}$ & Covariant derivative of vector field $\b{Y}$ in the direction of vector field $\b{X}$.
& $(\nabla_{\b{\Delta}_1} \b{\Delta}_2)(\b{X})$ & Covariant derivative of vector field $\b{\Delta}_2$ in the direction of vector field $\b{\Delta}_1$ at point $\b{X}$. \\

\hline
$(\cdot)_{,k}$ & Partial derivative in Ricci calculus notation. 
& & \\

\hline
$(\cdot)_{;k}$ & Covariant derivative in Ricci calculus notation. 
& & \\

\hline
$T(\b{X},\b{Y})$ & Torsion tensor. 
& &  \\

\hline
$R^\ell{}_{ijk}$ & Components of the Riemann curvature tensor. 
& & \\

\hline
$\b{\gamma}(t)$ & Geodesic (curve), with parameter $t$, on a general manifold $\mathcal{M}$.
& $\b{X}(t)$ & Geodesic (curve), with parameter $t$, on a matrix manifold. \\

\hline
$\operatorname{grad} f$ & Riemannian gradient of a smooth function $f$. 
& & \\

\hline
$\operatorname{Hess} f$ & Riemannian Hessian of a smooth function $f$. 
& & \\

\hline
$\operatorname{Exp}_{\b{p}}(\b{\xi})$ & Exponential map at $\b{p}$, applied on tangent vector $\b{\xi}$. 
& $\operatorname{Exp}_{\b{X}}(\b{\Delta})$ & Exponential map at matrix point $\b{X}$, applied on tangent vector (matrix) $\b{\Delta}$. \\

\hline
$\operatorname{Log}_{\b{p}}(\b{q})$ & Logarithm map at $\b{p}$, applied on point $\b{q}$. 
& $\operatorname{Log}_{\b{X}}(\b{Y})$ & Logarithm map at matrix point $\b{X}$, applied on matrix point $\b{Y}$. \\

\hline
 &  
& $\operatorname{exp}(\cdot),\operatorname{log}(\cdot)$ & Matrix exponential and matrix logarithm. \\

\hline
$\operatorname{Ret}_{\b{p}}(\b{\xi})$ & Retraction map at $\b{p}$, applied on tangent vector $\b{\xi}$. 
& $\operatorname{Ret}_{\b{X}}(\b{\Delta})$ & Retraction map at matrix point $\b{X}$, applied on tangent vector (matrix) $\b{\Delta}$. \\

\hline
$\mathcal{T}_{\b{\eta}}(\b{\xi})$ & Vector transport of tangent vector $\b{\xi}$ along tangent vector $\b{\eta}$ in a general manifold $\mathcal{M}$. 
& $\mathcal{T}_{\b{\Delta}}(\b{\Xi})$ & Vector transport of tangent vector (matrix) $\b{\Xi}$ along tangent vector (matrix) $\b{\Delta}$ in a matrix manifold. \\

\hline
$\mathcal{O}(.)$ & Big-O complexity notation
& $\mathrm{O}(d)$ & Orthogonal group (set of orthogonal $d \times d$ matrices) \\

\hline
$\Vert\cdot\Vert_2$ & $\ell_2$ norm of either a vector or a matrix. &
$\Vert\cdot\Vert_F$ & Frobenius norm of a matrix. \\

\hline
 &   
& $\operatorname{tr}(\cdot)$ & Trace of a square matrix. \\

 &   
&  &  \\



\bottomrule
\end{tabular}
\end{table*}

\subsection{Organization of This Monograph}

The remainder of this monograph is organized as follows.
Section \ref{section_definitions} introduces the definition of a Riemannian manifold,
starting from topology, topological manifolds, and smooth
manifolds. Section \ref{section_intrinsic_flatness_curvature} discusses intrinsic flatness and
curvature of manifolds and contrasts intrinsic and extrinsic
viewpoints. Section \ref{section_tensor_algebra} reviews the main ingredients of tensor
algebra (tensor calculus), including tangent and cotangent spaces, vector
fields, tensors, tensor products, and coordinate
transformations. Section \ref{section_metric_tensor} introduces the metric tensor and
its role in measuring lengths, angles, and raising and
lowering indices. Section \ref{section_connection_covariant_derivative} develops connections, covariant
derivatives, Christoffel symbols, torsion, and the
Levi-Civita connection. Section \ref{section_ricci_calculus_notation} explains Ricci calculus
notation, including comma and semicolon derivatives.
Section \ref{section_riemannian_curvature_quantities} presents the main Riemannian curvature quantities, including
the Riemann curvature tensor, sectional curvature, Ricci
curvature, scalar curvature, Gaussian curvature, and some
other advanced curvature measurements. Section \ref{section_ricci_flow} gives an
introduction to Ricci flow and its geometric interpretation.
Section \ref{section_curves_parallel_transport_geodesics} studies curves, absolute differentiation, parallel
transport, and geodesics. Section \ref{section_manifold_valued_optimization} develops the
Riemannian optimization framework, including the
Riemannian gradient, Riemannian Hessian, exponential and
logarithm maps, retractions, vector transport, and first-order
and second-order optimization methods on manifolds.
Section \ref{section_important_Riemannian_matrix_manifolds} specializes the general theory to important matrix
manifolds, namely the Stiefel manifold, the Grassmann
manifold, and the symmetric positive definite manifold,
deriving their main geometric and optimization-related
objects explicitly. Section \ref{section_toolboxes_for_riemannian_optimization} briefly introduces some
important software toolboxes and textbooks related to
Riemannian geometry and Riemannian optimization.
Finally, Section \ref{section_conclusion} concludes the monograph and outlines
possible directions for future development.

\section{Definition of Riemannian Manifold: From Topology to Smooth Manifold}\label{section_definitions}

In this section, we define the Riemannian manifold. We build gradually from topology and topological space. Then, we define the topological manifold. After introducing some characteristics of a topological manifold---including chart, smooth atlas, and maximal atlas---we define the smooth manifold. Finally, we define the Riemannian manifold. 
For more information on topological manifolds, smooth manifolds, and Riemannian manifolds, the reader can refer to \cite{lee2010topological}, \cite{lee2013smooth}, and \cite{do1992riemannian,lee2006riemannian2,lee2018riemannian}, respectively. 

\subsection{Topology and Topological Space}

\begin{definition}[Topology and topological space \cite{lee2010topological,kelley2017general}]
Let $\mathcal{X}$ be a set. A \textbf{topology} on $\mathcal{X}$ is a collection $\tau$ of subsets of $\mathcal{X}$, called open sets, satisfying:
\begin{itemize}
\item $\varnothing, \mathcal{X} \in \tau$ 
\item If $U_1, \dots, U_k \in \tau$, then $\bigcap_{j=1}^k U_j \in \tau$. In other words, finite intersections of open sets are open. 
\item If $U_\alpha \in \tau, \forall \alpha \in A$ (where $A$ is the index set of topology), then $\bigcup_{\alpha \in A} U_\alpha \in \tau$. In other words, arbitrary unions of open sets are open. 
\end{itemize}
The pair $(\mathcal{X}, \tau)$ is called a \textbf{topological space} associated with the topology $\tau$. 
\end{definition}

\begin{figure}[!h]
\centering
\includegraphics[width=3.2in]{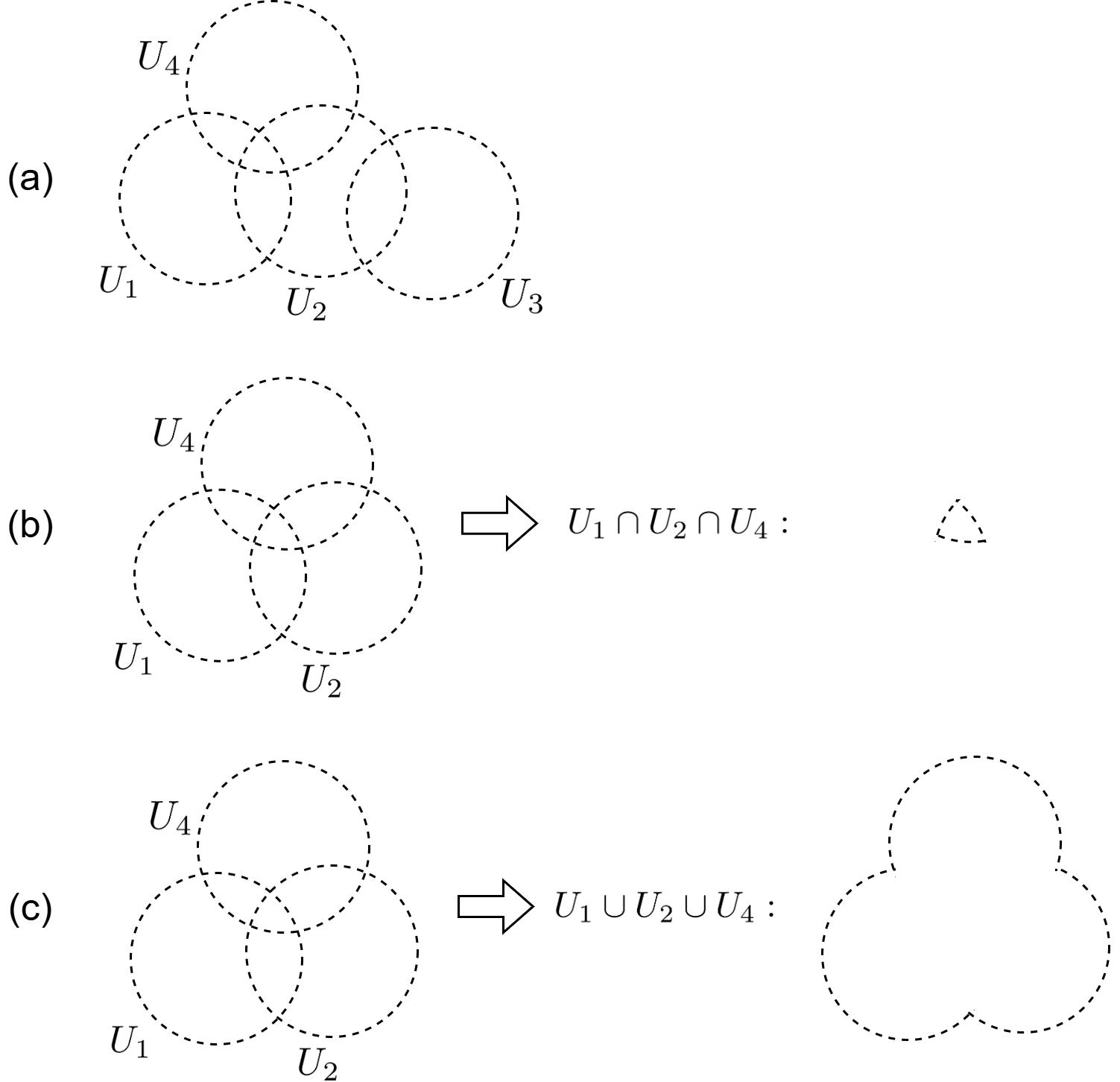}
\caption{Topology: (a) multiple open sets, (b) a finite intersection of some of the open sets is also an open set, and (c) union of some of the open sets is also an open set.}
\label{figure_topology}
\end{figure}

Intuitively, a topology is a collection (or set) of open sets. As shown in Fig. \ref{figure_topology}-a, these open sets may have overlap with each other. Any finite intersection of these open sets is another open set (see Fig. \ref{figure_topology}-b). The union of these open sets is also another open set (see Fig. \ref{figure_topology}-c). A topological space is a set equipped with a topology.

\begin{definition}[Hausdorff space \cite{lee2010topological,kelley2017general}]
A topological space $(\mathcal{X}, \tau)$ is \textbf{Hausdorff} if and only if for $x_1, x_2 \in \mathcal{X}$, $x_1 \neq x_2$, we have:
\begin{align}\label{equation_Hausdorff}
\exists\, \text{open sets } U,V \text{ such that } x_1 \in U,\, x_2 \in V,\, U \cap V = \varnothing.
\end{align}
\end{definition}

\begin{figure}[!h]
\centering
\includegraphics[width=2.5in]{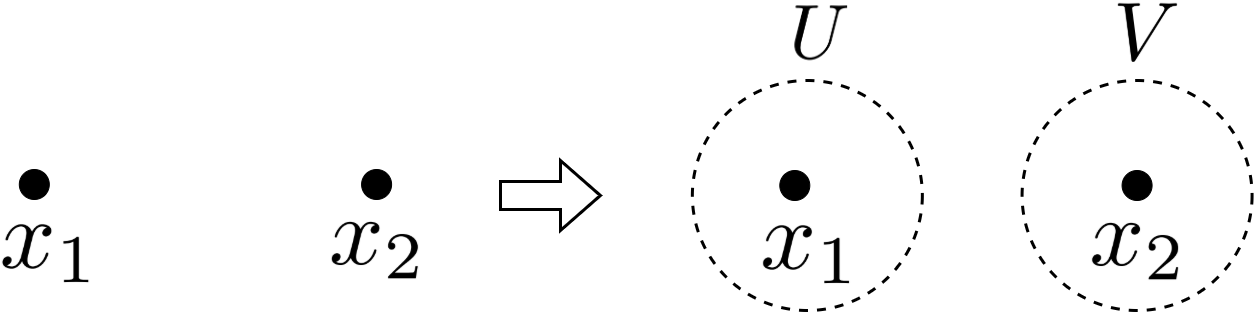}
\caption{A Hausdorff topological space where the points are distinguishable.}
\label{figure_Hausdorff}
\end{figure}

Intuitively, the points of a Hausdorff topological space are separable and distinguishable. Equation (\ref{equation_Hausdorff}) means that the two points $x$ and $y$ have neighborhoods or open sets $U$ and $V$ which do not overlap (see Fig. \ref{figure_Hausdorff}). Thus, they can be separated and distinguished.

\subsection{Homeomorphism and Diffeomorphism}

\begin{definition}[Isomorphism]
An \textbf{isomorphism} is a bijective mapping between two mathematical structures that preserves the relevant structure of the objects. In other words, if two objects are isomorphic, they are considered equivalent from the viewpoint of the structure being studied, even though they may look different geometrically or algebraically.
\end{definition}

\begin{definition}[Homeomorphism]
A \textbf{homeomorphism} is a structure-preserving map between topological spaces.
It is a transformation between two topological spaces
without tearing or gluing the topology (without tearing or gluing the space). It is studied in algebraic topology. The two topologies before and after a homeomorphism transformation are called homeomorphic to each other. The homeomorphic symbol is usually denoted by $\cong$.
\end{definition}

Intuitively, if you can transform a topology without tearing it, or without making any holes in it, and without gluing parts, this transformation is called homeomorphism. Note that the number of holes in a topology is important in algebraic topology. 

\begin{figure}[!h]
\centering
\includegraphics[width=3.2in]{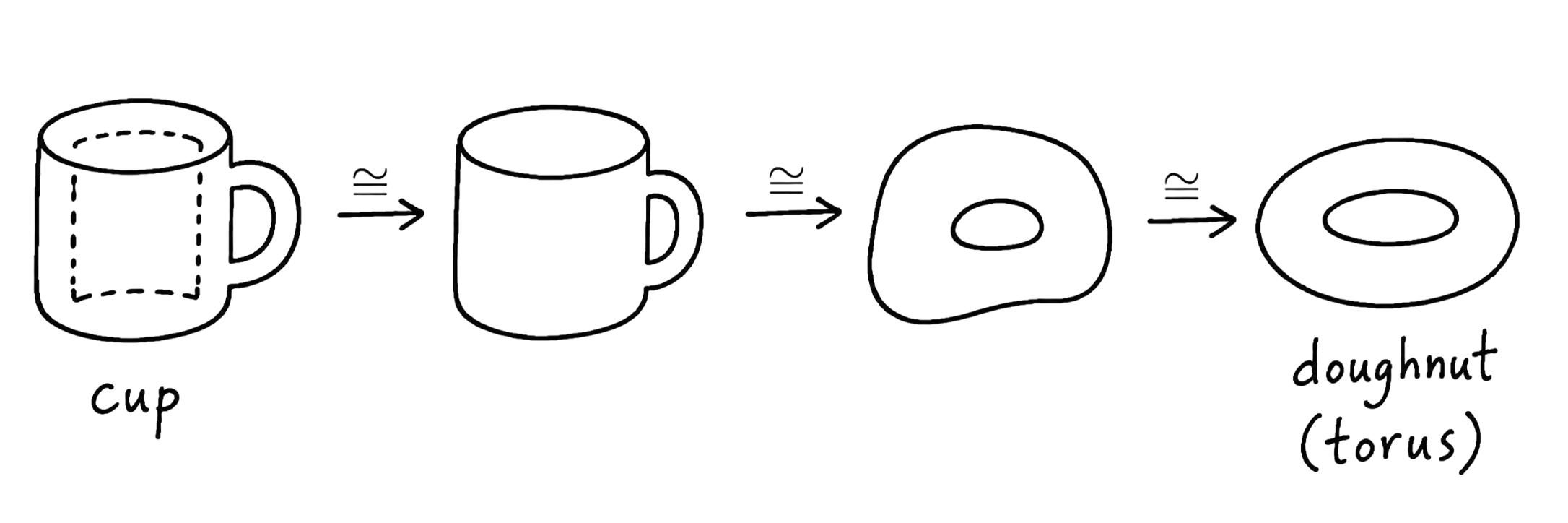}
\caption{Showing that a cup and a doughnut (genus-1 torus) are homeomorphic by molding the cup gradually to a doughnut. First, we fill up the cup by increasing the height of its bottom. Then, we shrink the cup except its handle. Then, we mold the cup to become a clean doughnut. The hole in the doughnut corresponds to the hole in the handle of the cup.}
\label{figure_homeomorphism}
\end{figure}

For example, a well-known example is that a cup and a doughnut are homeomorphic to each other. A doughnut is called a genus-1 torus in topology, where genus-1 refers to one hole in it. Therefore, a cup and a genus-1 torus are homeomorphic. This is illustrated in Fig. \ref{figure_homeomorphism} where a cup is molded gradually to become a genus-1 torus. Note that the hole in the handle of the cup is equivalent to the central hole in the torus. 

\begin{remark}[Difference between isomorphism and homeomorphism]
The following explains the difference of isomorphism and homeomorphism. 
In simple terms, homeomorphism is a special type of isomorphism in topology.
Isomorphism is a very general concept. It means two mathematical structures are equivalent because there is a bijective mapping between them that preserves the structure being studied. The exact meaning of “structure” depends on the field, such as algebra, graph theory, topology, etc.
Homeomorphism is a topological isomorphism. It is a bijective, continuous mapping with a continuous inverse that preserves topological properties such as connectedness and number of holes.

In other words, if two structures are isomorphic, they are considered the same with respect to the mathematical structure being studied, but the meaning depends on the context.
If two spaces are homeomorphic, you can deform one into the other without tearing, gluing, or creating holes.
Isomorphism is structure-preserving equivalence in general, and homeomorphism is topology-preserving equivalence, as a special case in topology.
\end{remark}

\begin{definition}[Diffeomorphism]
A \textbf{diffeomorphism} is a smooth bijection whose inverse is also smooth.
It can also be said that a diffeomorphism is a homeomorphism transformation which is smooth and differentiable.
\end{definition}

\subsection{Dimension of Topological Space}

\begin{definition}[Embedded submanifold {\citep[Chapter 5]{lee2013smooth}}]
Let $\mathcal{M}$ be a smooth manifold. An embedded submanifold of $\mathcal{M}$ is a subset $\mathcal{S} \subseteq \mathcal{M}$ which is itself a manifold endowed with a smooth structure where the inclusion map $\mathcal{S} \hookrightarrow \mathcal{M}$ is a smooth embedding.
The quantity $\mathrm{dim}(\mathcal{M}) - \mathrm{dim}(\mathcal{S})$ is called the codimension of $\mathcal{S}$ in $\mathcal{M}$, where $\mathrm{dim}(\mathcal{M})$ and $\mathrm{dim}(\mathcal{S})$ denote the dimensionalities of $\mathcal{M}$ and $\mathcal{S}$, respectively. 
\end{definition}

\begin{lemma}[Whitney embedding theorem \cite{whitney1936differentiable,whitney1944self}]
Every $d$-dimensional differentiable manifold can be embedded in $\mathbb{R}^{2d+1}$ \cite{whitney1936differentiable}. In some cases, it can be embedded in $\mathbb{R}^{2d}$ \cite{whitney1944self}.
\end{lemma}

\begin{definition}[$n$-sphere $S^n$]\label{definition_n_sphere}
An \textbf{$n$-sphere}, denoted by $S^n$, with radius $r > 0$, is defined as:
\begin{align}\label{equation_n_sphere}
\boxed{
S^n := \{\b{x} \in \mathbb{R}^{n+1} \mid \Vert\b{x}\Vert = r\}.
}
\end{align}
If $r=1$, then $S^n$ is called a unit sphere. 
\end{definition}

\begin{definition}[$n$-ball $B^n$]
An \textbf{$n$-ball}, denoted by $B^n$, with radius $r > 0$, is defined as:
\begin{align}
\boxed{
B^n := \{\b{x} \in \mathbb{R}^{n} \,|\, \Vert \b{x} \Vert \leq r\}.
}
\end{align}
If $r=1$, then $B^n$ is called a unit ball. 
\end{definition}

Intuitively, the sphere is the boundary only but the ball contains both boundary and the interior. 
That is why, for norms, we have unit balls and not unit spheres.

\begin{figure}[!h]
\centering
\includegraphics[width=3.2in]{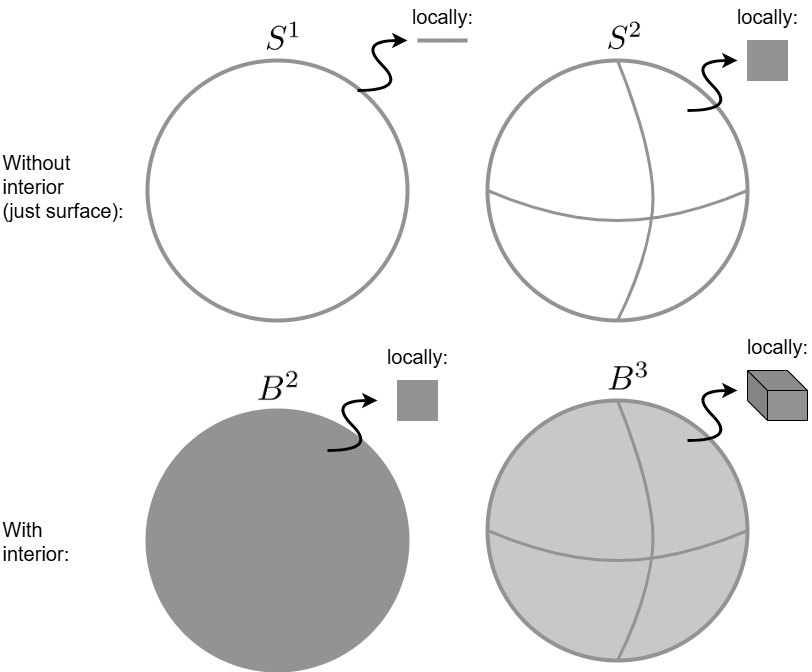}
\caption{Some examples of $S^1$, $S^2$, $B^2$, and $B^3$. The $S^1$, $S^2$, $B^2$, and $B^3$ are locally one-dimensional (embedded in 2D), two-dimensional (embedded in 3D), two-dimensional (embedded in 2D), and three-dimensional (embedded in 3D), respectively.}
\label{figure_sphere_ball}
\end{figure}

Note that an $n$-dimensional ball $B^n$ is locally $n$-dimensional, and it can be embedded in an $n$-dimensional space. 
But an $n$-dimensional sphere $S^n$ is locally $n$-dimensional, and it can be embedded in an $(n+1)$-dimensional space, according to the Whitney embedding theorem. 
To imagine this, consider a football filled with mud. It is a three-dimensional ball $B^3$ (because if we pick a part of the filled football, it is 3D) in a three-dimensional space we live in.
However, an empty football is a two-dimensional sphere $S^2$ (because if we pick a part of the skin of the football, it is 2D) in a three-dimensional space we live in.
Some examples of $S^1$, $S^2$, $B^2$, and $B^3$ are illustrated in Fig. \ref{figure_sphere_ball}.

\subsection{Topological Manifold}\label{section_topological_manifold}

\begin{definition}[Topological manifold \cite{lee2010topological}]
A topological space $(\mathcal{X}, \tau)$ is a \textbf{topological manifold} of dimension $n$, for $n \in \mathbb{Z}_{\geq 0}$, also called a topological $n$-manifold, if all the following conditions hold:
\begin{itemize}
\item $(\mathcal{X}, \tau)$ is Hausdorff. 
\item $(\mathcal{X}, \tau)$ has a countable basis. 
\item $(\mathcal{X}, \tau)$ is locally homeomorphic to $n$-dimensional Euclidean space, $\mathbb{R}^n$. 
\end{itemize}
\end{definition}

\begin{definition}[Chart \cite{lee2010topological}]
Consider a topological manifold $\mathcal{M} := (\mathcal{X}, \tau)$. It is locally homeomorphic to $\mathbb{R}^n$, meaning that for all $x \in X$, there exists an open set $U$ containing $x$ and a homeomorphism $\varphi: U \rightarrow \varphi(U)$ where $\varphi(U)$ is an open subset of $\mathbb{R}^n$. Such a mapping is denoted by $\varphi: U \overset{\cong}{\longrightarrow} \varphi(U)$ and the tuple $(U, \varphi)$ is called a \textbf{coordinate chart}, or a \textbf{chart} in short, for $\mathcal{M}$. 
\end{definition}

\begin{figure}[!h]
\centering
\includegraphics[width=2.5in]{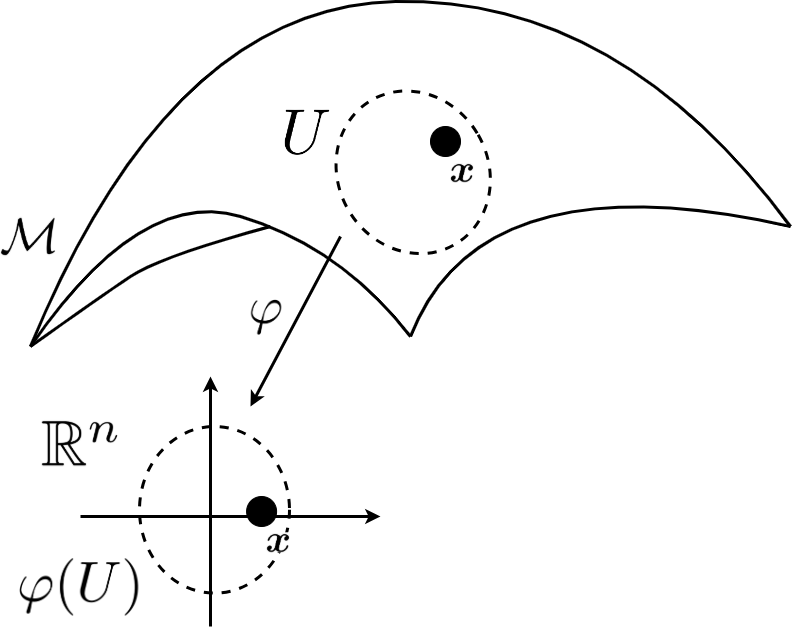}
\caption{A chart $(U, \varphi)$ on the manifold $\mathcal{M}$. The mapping $\varphi(U)$ approximates the open set $U$ locally to a local flat Euclidean space $\mathbb{R}^n$. In other words, $\varphi(U)$ and the flat Euclidean space $\mathbb{R}^n$ are homeomorphic.}
\label{figure_chart}
\end{figure}

An example chart is illustrated in Fig. \ref{figure_chart}.
As this figure shows, in a chart $(U, \varphi)$, the mapping $\varphi(U)$ approximates the open set $U$ locally to a local flat Euclidean space $\mathbb{R}^n$. In other words, $\varphi(U)$ and the flat Euclidean space $\mathbb{R}^n$ are homeomorphic. 

\begin{figure}[!h]
\centering
\includegraphics[width=3in]{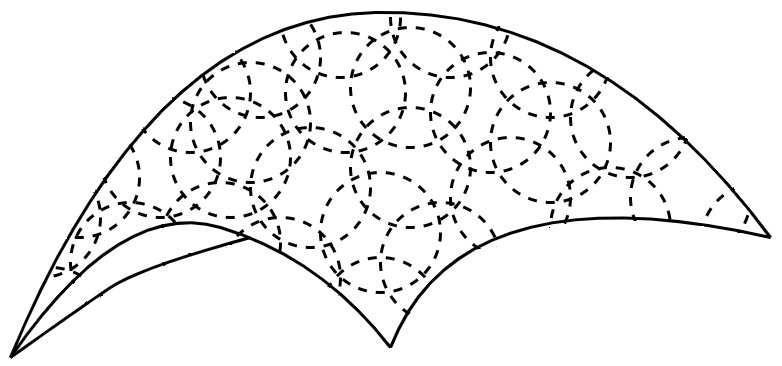}
\caption{A smooth atlas for a topological manifold.}
\label{figure_atlas}
\end{figure}

\begin{figure}[!h]
\centering
\includegraphics[width=3in]{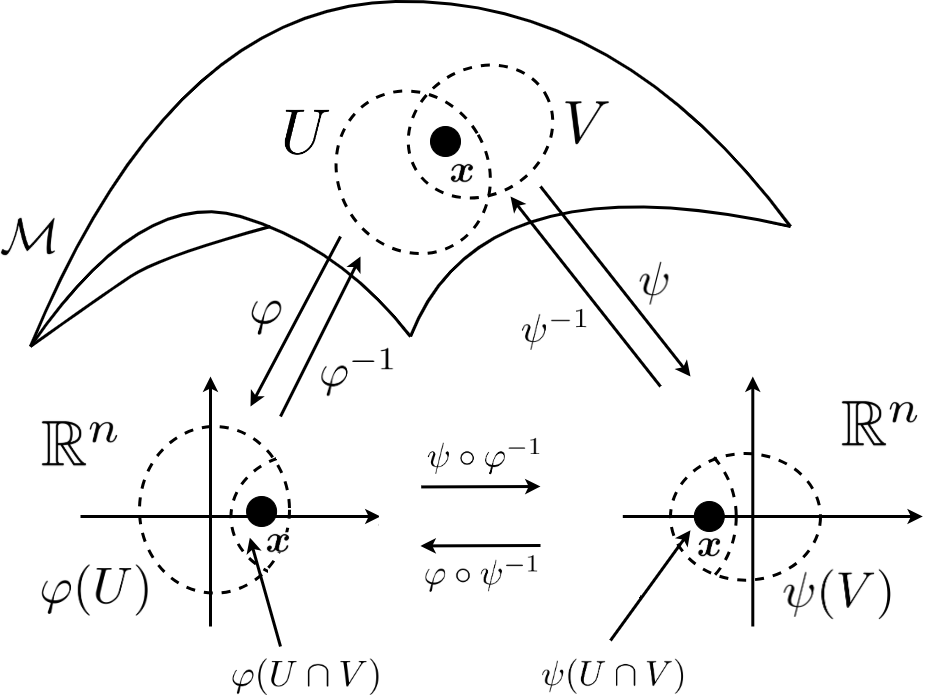}
\caption{Two charts $(U, \varphi)$ and $(V, \psi)$ are smoothly compatible if the mapping $\psi \circ \varphi^{-1}$ is a diffeomorphism. Note that the mapping $(\psi \circ \varphi^{-1})(U \cap V)$ maps $\varphi(U \cap V)$ back from $\mathbb{R}^n$ to the manifold $\mathcal{M}$ and then maps from the manifold $\mathcal{M}$ to $\mathbb{R}^n$ in $\psi(U \cap V)$.}
\label{figure_chart2}
\end{figure}

\begin{definition}[Smooth atlas \cite{lee2013smooth}]
A \textbf{smooth atlas} $\mathcal{A}$ for a topological $n$-manifold $\mathcal{M}$ is a collection of charts $(U_\alpha, \varphi_\alpha)$ for $\mathcal{M}$ such that:
\begin{itemize}
\item They cover $\mathcal{M}$, i.e., $\bigcup_{\alpha \in A} U_\alpha = \mathcal{M}$. For example, see Fig. \ref{figure_atlas} illustrating some open sets of the charts covering the manifold.
\item Any two charts in this collection are smoothly compatible. Note that two charts $(U, \varphi)$ and $(V, \psi)$ are smoothly compatible if the mapping $\psi \circ \varphi^{-1}$ is a diffeomorphism. The mapping $\psi \circ \varphi^{-1}$ is illustrated in Fig. \ref{figure_chart2}. 
\end{itemize}
\end{definition}

\begin{definition}[Maximal atlas \cite{lee2013smooth}]
A smooth atlas $\mathcal{A}$ for a topological $n$-manifold $\mathcal{M}$ is \textbf{maximal} if it is not contained in any other smooth atlas for $\mathcal{M}$. 
\end{definition}

\subsection{Smooth Manifold and Riemannian Manifold}\label{section_smooth_manifold_riemannian_manifold}

\begin{definition}[Smooth manifold \cite{lee2013smooth}]
A \textbf{smooth manifold} $\mathcal{M}$ of dimension $n$, also called a smooth $n$-manifold, is a topological $n$-manifold together with a choice of maximal smooth atlas $\mathcal{A}$ on $\mathcal{M}$. 
\end{definition}

\begin{definition}[Point on a manifold]
Let \(\mathcal{M}\) be a manifold. A \textbf{point} on
\(\mathcal{M}\) is simply an element of the set
\(\mathcal{M}\). In other words, if \(\b{p}\in\mathcal{M}\),
then \(\b{p}\) is a point of the manifold.

In differential geometry, the word ``point'' refers to a
location on the manifold itself, whereas tangent vectors,
cotangent vectors, and tensors are objects attached to that
point.
\end{definition}



\begin{definition}[Tangent space on a smooth or Riemannian manifold]\label{definition_tangent_space_no_math}
Let $\mathcal{M}$ be a smooth manifold and let $\b{p} \in \mathcal{M}$. The \textbf{tangent space} $T_{\b{p}}\mathcal{M}$ is the vector space consisting of all tangent vectors at point $\b{p}$, which represent possible directions in which one can pass through $\b{p}$ along smooth curves on the manifold. 

Intuitively, the tangent space is a linear approximation of the manifold around the point $\b{p}$. For a Riemannian manifold, the metric induces an inner product structure on $T_{\b{p}}\mathcal{M}$, allowing measurement of lengths and angles between tangent vectors.
\end{definition}

\begin{figure}[!h]
\centering
\includegraphics[width=3in]{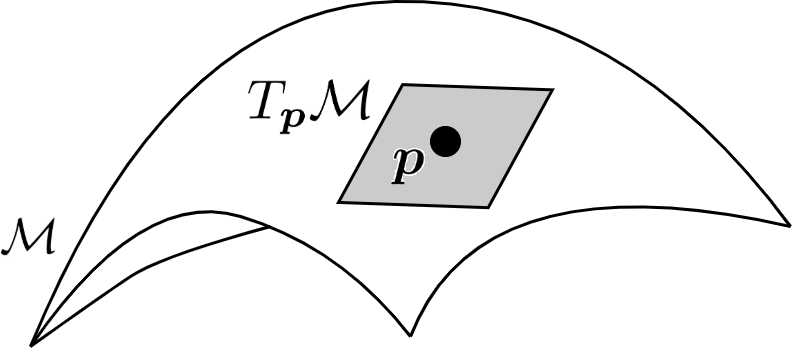}
\caption{A locally two-dimensional tangent space $T_{\b{p}}\mathcal{M}$ at point $\b{p}$ in a locally two-dimensional manifold $\mathcal{M}$.}
\label{figure_tangent_space}
\end{figure}

Intuitively, in an $n$-manifold, every point $\b{p} \in \mathcal{M}$ can have an $n$-dimensional tangent space $T_{\b{p}}\mathcal{M}$, where the tangent space is a flat Euclidean space. For example, consider Fig. \ref{figure_tangent_space} which illustrates a locally two-dimensional smooth manifold embedded in 3D. The tangent space at a point on the manifold is a two-dimensional flat Euclidean space. 

\begin{definition}[Riemannian manifold {\cite{do1992riemannian,lee2006riemannian2,lee2018riemannian}}]
Let $\mathcal{M}$ be a smooth manifold. A \textbf{Riemannian metric}
on $\mathcal{M}$ is a family of inner products:
\begin{align}
g_{\b{p}} : T_{\b{p}}\mathcal{M} \times T_{\b{p}}\mathcal{M} \to \mathbb{R},
\qquad \b{p}\in\mathcal{M},
\end{align}
such that:
\begin{enumerate}
    \item for every $\b{p}\in\mathcal{M}$, the map
    $g_{\b{p}}(\cdot,\cdot)$ is an inner product on the tangent space
    $T_{\b{p}}\mathcal{M}$, and
    \item the metric varies smoothly with $\b{p}$, meaning that for any
    smooth vector fields $\b{X},\b{Y}$ on $\mathcal{M}$, the function:
    \begin{align}
    \b{p} \mapsto g_{\b{p}}(\b{X}(\b{p}),\b{Y}(\b{p})),
    \end{align}
    is smooth.
\end{enumerate}
A \textbf{Riemannian manifold} is a smooth manifold $\mathcal{M}$
equipped with a Riemannian metric $g$.
\end{definition}

\begin{remark}[Information provided by Riemannian metric]
A Riemannian metric allows us to measure lengths of tangent vectors,
angles between tangent vectors, lengths of curves, and hence distances
on the manifold. It also provides the geometric structure needed to
define geodesics, curvature, gradients, Hessians, and other
Riemannian objects.
In other words, knowing the Riemannian metric at every point of a Riemannian manifold provides all needed information about the geometry of the manifold. 
\end{remark}

\begin{remark}[Coordinate-free and coordinate-based differential geometry]
There are two common styles of presentation in differential
geometry: the \textbf{coordinate-free} style and the
\textbf{coordinate-based} style.

In the coordinate-free style, geometric objects are defined
intrinsically, without fixing a local coordinate system.
For example, a tangent vector is treated as a geometric object
in $T_{\b{p}}\mathcal{M}$, a metric is treated as a bilinear map
$g_{\b{p}} : T_{\b{p}}\mathcal{M}\times T_{\b{p}}\mathcal{M}\to\mathbb{R}$,
and a connection is defined by its structural properties.
This style is used mostly in modern differential geometry,
global analysis, geometric mechanics, and much of contemporary
Riemannian optimization because it makes the geometric meaning
transparent and keeps the statements independent of the choice
of coordinates \cite{lee2013smooth, lee2018riemannian, tu2017differential, absil2008optimization, boumal2023introduction}.

In the coordinate-based style, one chooses a local coordinate
system $(x^1,\dots,x^n)$ and expresses geometric objects in
terms of their components, such as $g_{ij}$,
$\Gamma^k_{ij}$, and $R^\ell{}_{ijk}$.
This style is used mostly in tensor calculus, continuum
mechanics, general relativity (in physics), engineering, and explicit calculations,
because it turns abstract geometric definitions into formulas
that can be computed directly \cite{do1992riemannian, itskov2007tensor, bishop2012tensor, hartle2021gravity}.

Both styles are important. The coordinate-free viewpoint is
better for understanding the intrinsic meaning of definitions,
theorems, and geometric relations. The coordinate-based
viewpoint is better for deriving explicit formulas, performing
symbolic manipulations, and implementing numerical algorithms.
Therefore, the coordinate-free style emphasizes
\emph{what} a geometric object is, while the coordinate-based
style emphasizes \emph{how} that object is written and computed
in a chosen chart.

These two styles are not different theories; rather, they are
two equivalent ways of describing the same geometry.
A coordinate-based formula is obtained by expressing an
intrinsic coordinate-free object in a local basis, and a
coordinate-free statement can often be recovered by recognizing
which geometric object the coordinate expression represents.
For this reason, a complete understanding of differential
geometry benefits from both viewpoints.

In this monograph, we use both styles deliberately.
Whenever possible, we first introduce the geometric objects in
a coordinate-free manner and then derive their coordinate
expressions in detail. This is especially important for
Riemannian optimization, where the coordinate-free
viewpoint clarifies the underlying geometry, while the
coordinate-based viewpoint is often the one needed for
practical derivations and implementation.
\end{remark}

\begin{remark}[Meaning of ``endowed with'' and ``equipped with'' in differential geometry]
In differential geometry, the phrases ``\textbf{endowed with}"
and ``\textbf{equipped with}" usually mean that an already
existing mathematical object is considered together with
some additional structure\footnote{These phrases may be obvious to geometers but we explain them here for the people new to this field.}.

For example, a topological manifold \(\mathcal{M}\)
\emph{endowed with} a smooth atlas becomes a smooth
manifold, and a smooth manifold \(\mathcal{M}\)
\emph{equipped with} a Riemannian metric \(g\) becomes a
Riemannian manifold. Likewise, a vector bundle can be
equipped with a connection, and a manifold can be endowed
with a symplectic form, a complex structure, or other
geometric structures.

In most mathematical writing, the phrases
``\emph{endowed with}" and ``\emph{equipped with}" have essentially
the same meaning and are often used interchangeably.
Both indicate that the object itself is not changed as an
underlying set or manifold, but rather that extra structure
is assigned to it.

The slight difference is mostly stylistic. The phrase
``\emph{endowed with}" often sounds a bit more formal and is
frequently used when emphasizing that a structure is given
to an object as part of its definition. The phrase
``\emph{equipped with}" is also standard and often sounds a bit
more direct or concrete. In this monograph, both phrases
refer to the same idea unless stated otherwise.
\end{remark}

\section{Intrinsic Flatness and Curvature of Manifold}\label{section_intrinsic_flatness_curvature}

\subsection{Intrinsic versus Extrinsic Curvature}\label{section_intrinsic_extrinsic_curvature}



As also discussed in Section \ref{section_topological_manifold}, consider a piece of paper. It is flat. If an ant is on the paper, it can only traverse the paths on the two-dimensional paper. Therefore, the ant feels the paper as a flat two-dimensional space. We can get the paper and fold it smoothly. The ant still feels the paper as a flat two-dimensional space. The ant does not understand folding of the paper because it is too small and the paper is still two-dimensional locally. 

The paper is locally flat; we say that it is flat intrinsically. However, the paper is curved (folded) in the three-dimensional space. We say that the two-dimensional paper is \textit{embedded} in the three-dimensional space and the folded paper is curved extrinsically. 
In summary, there are two types of curvature, i.e., \textit{intrinsic curvature} and \textit{extrinsic curvature}. 

\begin{figure}[!h]
\centering
\includegraphics[width=3in]{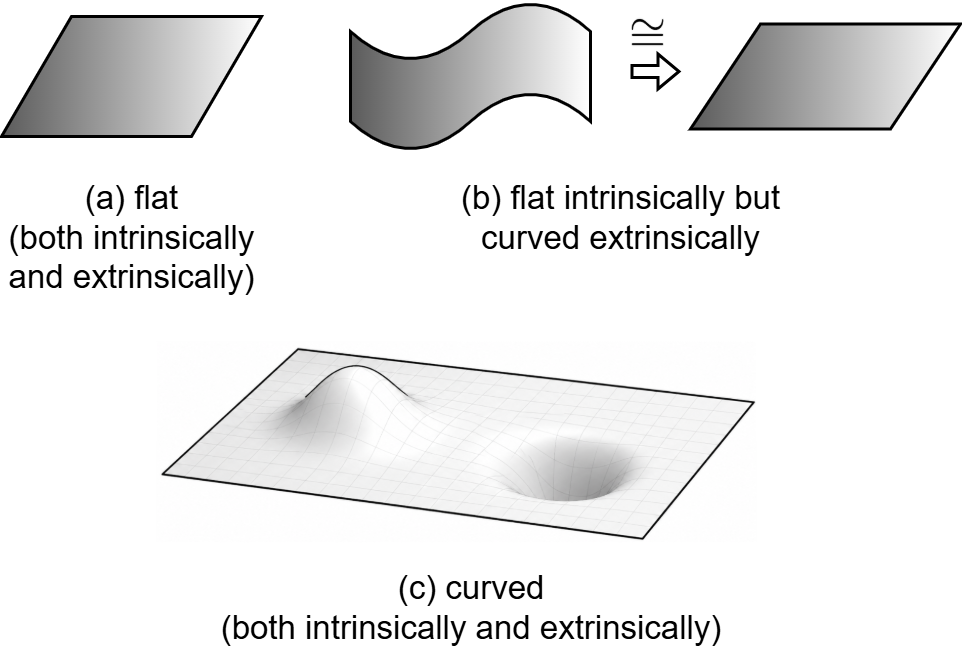}
\caption{Intrinsic and extrinsic curvature: (a) a manifold being flat both intrinsically and extrinsically, (b) a manifold being flat intrinsically (i.e., is homeomorphic to a flat manifold) but curved extrinsically, and (c) a manifold being curved both intrinsically and extrinsically.}
\label{figure_intrinsic_curvature}
\end{figure}

As illustrated in Figs. \ref{figure_intrinsic_curvature}-a and \ref{figure_intrinsic_curvature}-b, an unfolded or smoothly folded paper is intrinsically flat. This is because the folded paper can be unfolded and put on the two-dimensional surface. The folded paper is flat intrinsically but curved extrinsically. The unfolded paper is flat both intrinsically and extrinsically. 

However, consider one or multiple bumps in a locally two-dimensional manifold (see Fig. \ref{figure_intrinsic_curvature}-c). For example, assume landscape of ground to be a two-dimensional manifold where the mountains are positive bumps and valleys are negative bumps. This landscape can never be unfolded in a way to put on a two-dimensional surface. Therefore, this landscape is curved both intrinsically and extrinsically. 

The intrinsic curvature is the curvature felt by a small bug or ant on the manifold. For example, when a small bug is put on Fig. \ref{figure_intrinsic_curvature}-a or Fig. \ref{figure_intrinsic_curvature}-b, the bug feels it as a flat surface\footnote{That is why we humans also feel the earth (especially the flat plain) as a flat surface; this is because we are small in the scale of earth.}.
Therefore, Figs. \ref{figure_intrinsic_curvature}-a and \ref{figure_intrinsic_curvature}-b are flat intrinsically. 
However, when a small bug is put on Fig. \ref{figure_intrinsic_curvature}-c, it feels the bumps upwards and downwards. Thus, Fig. \ref{figure_intrinsic_curvature}-c has intrinsic curvature. 

In differential geometry, the term curvature typically denotes intrinsic curvature, unless explicitly specified otherwise\footnote{This is because \textit{Carl Friedrich Gauss}, who discussed curvature of a two-dimensional manifold for the first time in 1828 \cite{gauss1828disquisitiones}, considered the intrinsic curvature regardless of how it is embedded or folded in the three-dimensional space.}.
Similarly, in physics, curvature is the intrinsic curvature because we are already in the space-time manifold, i.e., the universe, and not outside of it. 
The notion of intrinsic curvature was developed by \textit{Carl Friedrich Gauss} in 1828 \cite{gauss1828disquisitiones}.


\subsection{Cartesian, Affine, and Curvilinear Coordinates}

\begin{figure}[!h]
\centering
\includegraphics[width=3.2in]{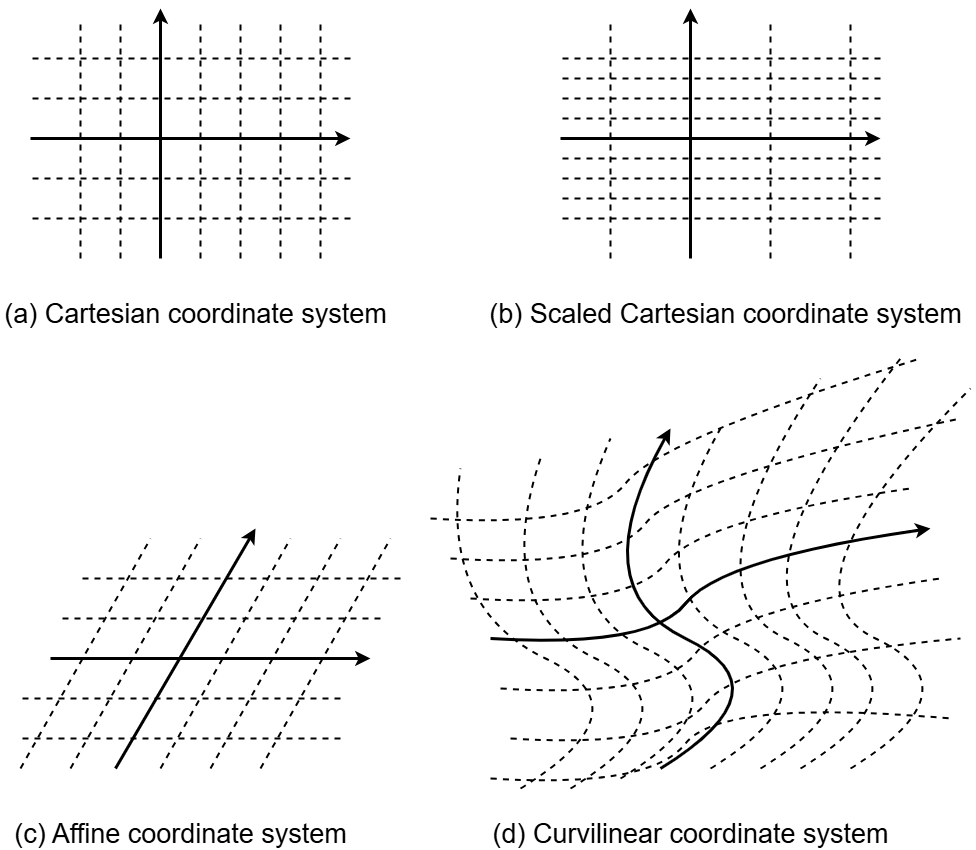}
\caption{Coordinate systems: (a) Cartesian coordinate system, (b) scaled Cartesian coordinate system, (c) affine coordinate system, and (d) curvilinear coordinate system.}
\label{figure_coordinate_system}
\end{figure}

The regular coordinate system that we learn in high school is the \textit{Cartesian coordinate system} where the axes are perpendicular to one another. A two-dimensional Cartesian coordinate system is depicted in Fig. \ref{figure_coordinate_system}-a.

The axes can be scaled (e.g., see Fig. \ref{figure_coordinate_system}-b).
Moreover, the axes can have angles more or less than $90$ degrees. If the axes remain straight, although they may be scaled and may have non-perpendicular angle with each other, the coordinate system is called \textit{affine coordinate system}. A two-dimensional affine coordinate system is illustrated in Fig. \ref{figure_coordinate_system}-c.

\begin{figure}[!h]
\centering
\includegraphics[width=3.2in]{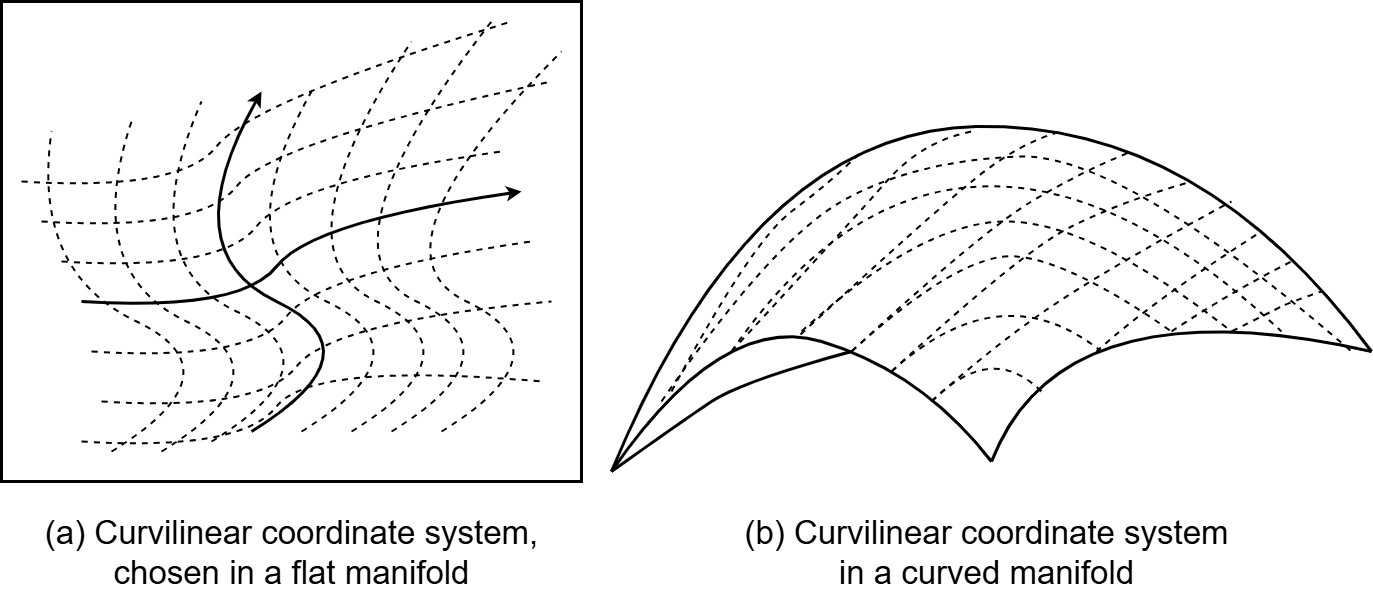}
\caption{Two occasions where curvilinear coordinate system may appear: (a) Curvilinear coordinate system chosen in an intrinsically flat manifold, (b) Curvilinear coordinate system in an intrinsically curved manifold.}
\label{figure_curvilinear_coordinates}
\end{figure}

If at least one of the axes is curved coordinate lines instead of straight affine axes, the coordinate system is called \textit{curvilinear coordinate system}.
A two-dimensional curvilinear coordinate system is illustrated in Fig. \ref{figure_coordinate_system}-d.
Curvilinear coordinate system occurs in two scenarios:
\begin{enumerate}
\item When the space is intrinsically flat but coordinates are chosen to be curvy (e.g., see Fig. \ref{figure_curvilinear_coordinates}-a), or
\item When the space is intrinsically curved, i.e., it has intrinsic curvature (e.g., see Fig. \ref{figure_curvilinear_coordinates}-b). 
\end{enumerate}

\section{Tensor Algebra (Tensor Calculus)}\label{section_tensor_algebra}

\textit{Tensor algebra}, also called the \textit{tensor calculus}, is one of the important backbones of differential geometry\footnote{There are many books on tensor calculus that the reader can also refer to. Some of the introductory books on tensor calculus are \cite{sochi2016tensor,grinfeld2013introduction,itskov2007tensor} and some of the advanced books are \cite{kobayashi1996foundations,bishop2012tensor}.}. Here, we briefly introduce the building blocks of tensor calculus, including tangent space, tangent vector, cotangent space, tangent bundle, cotangent bundle, vector, covector, vector field, tensor, tensor product, contravariant components, covariant components, and transformation of coordinates. 

\subsection{Preliminaries}

\textit{Einstein Convention}---also known as the \textit{Einstein notation}, \textit{Einstein summation convention}, and \textit{Einstein summation notation}---was proposed by Einstein in 1916 in his paper of general relativity \cite{einstein1916foundation}. It is defined in the following. 

\begin{definition}[Einstein summation convention \cite{einstein1916foundation}]
In differential geometry, there are usually many summations where the summation indices appear in both upper index and lower index of the summed expressions. 
According to \textbf{Einstein summation convention} \cite{einstein1916foundation}, when an index variable appears twice in a single term (once as a superscript and once as a subscript), it implies summation over all possible values of that index (e.g., $V^i \partial_i = \sum_{i=1}^n V^i \partial_i$).
In this way, the expressions are simplified by not writing the summations. Some examples of Einstein convention are:
\begin{align}
&A^j B_j^i := \sum_{j=1}^n A^j B_j^i, \\
&A^j_i B_j^i := \sum_{i=1}^n \sum_{j=1}^n A^j_i B_j^i, 
\end{align}
where summations are over the dimensions of indices. 
\end{definition}

Henceforth in this paper, we use Einstein convention, unless mentioned otherwise. 

\begin{definition}[Dummy variable]\label{definition_dummy_variable}
The summation variables, which the summation sums over, are \textbf{dummy variables}. They are called dummy variables because their name is not important and you can replace their name with anything. For example, these expressions are equivalent:
\begin{align*}
A^j B_j = A^i B_i = A^{\text{Benyamin}} B_{\text{Benyamin}},
\end{align*}
where Einstein convention has been used. In these expressions, $i$, $j$, and `Benyamin' are all dummy variables because one can change the name of summation variables:
\begin{align*}
\sum_{j=1}^n A^j B_j = \sum_{i=1}^n A^i B_i = \sum_{\text{Benyamin}=1}^n A^{\text{Benyamin}} B_{\text{Benyamin}}.
\end{align*}
\end{definition}

\begin{remark}[Naming of dummy variables]
According to Definition \ref{definition_dummy_variable}, the names of dummy variables are not important at all. However, in some applications of differential geometry, people may use some specific names for dummy variables by convention. For example, in general relativity in physics, people usually use Greek letter $\mu$ for indexing of space-time manifold (i.e., $\mu \in \{0, 1, 2, 3\}$ for time and space components) while they usually use Latin letter $i$ for indexing of three-dimensional space manifold (i.e., $i \in \{1, 2, 3\}$ for space components) \cite{susskind2025general}. 
\end{remark}

\begin{definition}[Kronecker delta]\label{definition_kronecker_delta}
The \textbf{Kronecker delta} acts as the identity tensor\footnote{Tensor will be defined later, but we define Kronecker delta here because it is required in many topics of differential geometry.}. 
The Kronecker delta with lower indices, denoted by $\delta_{ij}$, is defined as:
\begin{align}\label{equation_Kronecker_delta}
\boxed{
\delta_{ij} := 
\left\{
    \begin{array}{ll}
        1 & \mbox{if } i = j, \\
        0 & \mbox{if } i \neq j.
    \end{array}
\right.
}
\end{align}
The Kronecker delta with lower and upper indices, denoted by $\delta_i^j$, is defined as:
\begin{align}\label{equation_Kronecker_delta_lower_upper}
\boxed{
\delta_i^j := 
\left\{
    \begin{array}{ll}
        1 & \mbox{if } i = j, \\
        0 & \mbox{if } i \neq j.
    \end{array}
\right.
}
\end{align}
\end{definition}

\subsection{Smooth Functions on Smooth Manifold}

\begin{definition}[Function on a manifold]
Let $\mathcal{M}$ be a smooth manifold. A function:
\begin{align}
\boxed{
f : \mathcal{M} \to \mathbb{R},
}
\end{align}
is a \textbf{function on the manifold}, which maps every point of manifold, $\b{p} \in \mathcal{M}$, to a real-valued number. 
\end{definition}

\begin{definition}[Smooth real-valued function on a manifold]
Let $\mathcal{M}$ be a smooth manifold. A function:
\begin{align}
f : \mathcal{M} \to \mathbb{R},
\end{align}
is called a \textbf{smooth function}, if for every chart $(U, \varphi)$ of $\mathcal{M}$ with $\b{p} \in U \subset \mathcal{M}$,
the composition:
\begin{align}\label{equation_f_o_varphi_inverse}
f \circ \varphi^{-1} : \varphi(U) \subset \mathbb{R}^n \to \mathbb{R},
\end{align}
is a smooth function in the usual sense on $\mathbb{R}^n$ (i.e., all partial derivatives of all orders exist and are continuous)\footnote{Note that Eq. (\ref{equation_f_o_varphi_inverse}) means that a point from the chart $\varphi(U)$ is taken and passed through the function $\phi^{-1}(.)$ to obtain the point in the manifold. Then the point in the manifold is passed through the function $f$ to output a real-valued function.}. 

Intuitively, a smooth function is a function acting on the manifold locally where the function is smooth, i.e., it has continuous derivatives of all orders (denoted by $C^\infty$).
\end{definition}

\begin{definition}[The set of smooth real-valued functions on manifold]
Let $\mathcal{M}$ be a smooth manifold. Then:
\begin{align}
\boxed{
C^\infty(\mathcal{M}) := \{ f : \mathcal{M} \to \mathbb{R} \mid f \text{ is infinitely differentiable} \}.
}
\end{align}
The set $C^\infty(\mathcal{M})$ is called the \textbf{algebra (or set) of smooth real-valued functions} on $\mathcal{M}$.
In other words, $C^\infty(\mathcal{M})$ is the set of all smooth (infinitely differentiable) functions on the manifold.
As a result, a smooth function is denoted by $f \in C^\infty(\mathcal{M})$.
\end{definition}

\subsection{Tangent Space and Tangent Vector}\label{section_tangent_space_tangent_vector}

We know that the earth is like a ball but when we stand on it and look around, it seems that it is flat. It is because we are very small compared to the earth so we see it flat locally. Likewise, at every point on a curvy manifold, we can locally consider it flat. In other words, as illustrated in Fig. \ref{figure_tangent_space} (see Section \ref{section_smooth_manifold_riemannian_manifold}), we can consider a Euclidean (flat) tangent space at every point $\b{p}$ on the manifold $\mathcal{M}$. If the manifold is locally $n$-dimensional, the tangent space is also $n$-dimensional. For example, in the Fig. \ref{figure_tangent_space}, the manifold is locally two-dimensional (embedded in three dimensions), so its tangent space is a two-dimensional tangent plane. 
The tangent space of the manifold $\mathcal{M}$ at point $\b{p}$ is denoted by $T_{\b{p}}\mathcal{M}$.

\begin{figure}[!h]
\centering
\includegraphics[width=3in]{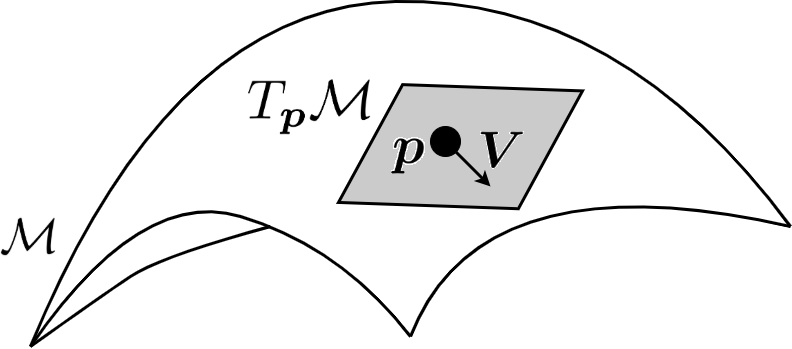}
\caption{A two-dimensional tangent vector $\b{V}$ existing in the locally two-dimensional tangent space $T_{\b{p}}\mathcal{M}$ at point $\b{p}$ in a locally two-dimensional manifold $\mathcal{M}$. The tangent vector lies in the tangent space, i.e., $\b{V} \in T_{\b{p}}\mathcal{M}$.}
\label{figure_tangent_vector}
\end{figure}

\begin{definition}[Tangent vector]
Let $\mathcal{M}$ be a smooth manifold and consider the point $\b{p}$ in the manifold, i.e., $\b{p} \in \mathcal{M}$.  
A \textbf{tangent vector} at $\b{p}$ is a map:
\begin{align}
\boxed{
\b{V} : C^\infty(\mathcal{M}) \to \mathbb{R},
}
\end{align}
such that for all $f,g \in C^\infty(\mathcal{M})$ and $a,b \in \mathbb{R}$ the following hold:
\begin{itemize}
\item Linearity:
\begin{align}
\b{V}(af + bg) = a\,\b{V}(f) + b\,\b{V}(g).
\end{align}
\item Leibniz rule (product rule):
\begin{align}
\b{V}(fg) = f(\b{p})\b{V}(g) + g(\b{p})\b{V}(f).
\end{align}
\end{itemize}
Such a map $\b{V}$ is called a derivation or tangent vector at $\b{p}$.

In other words, a tangent vector at point $\b{p} \in \mathcal{M}$ is a vector which is tangent to the manifold $\mathcal{M}$ at point $\b{p}$. 
In an $n$-dimensional space or manifold, a tangent vector is an $n$-dimensional vector.
An example two-dimensional tangent vector in a locally two-dimensional manifold is illustrated in Fig. \ref{figure_tangent_vector}.
\end{definition}

Recall Definition \ref{definition_tangent_space_no_math} in Section \ref{section_smooth_manifold_riemannian_manifold} for a descriptive definition of tangent space. In the following, we provide a more detailed definition of tangent space. 

\begin{definition}[Tangent space]
Let $\mathcal{M}$ be a smooth manifold and consider the point $\b{p}$ in the manifold, i.e., $\b{p} \in \mathcal{M}$.   
The \textbf{tangent space} at $\b{p}$, denoted by $T_{\b{p}}\mathcal{M}$, is the set of all tangent vectors at $\b{p}$:
\begin{align}
\boxed{
T_{\b{p}}\mathcal{M} := \{ \b{V} : C^\infty(\mathcal{M}) \to \mathbb{R} \mid \b{V} \text{ is a derivation at } p \}.
}
\end{align}
Suppose the function $f: \mathcal{M} \rightarrow \mathbb{R}$ be a smooth function on manifold $\mathcal{M}$, i.e., $f \in C^\infty(\mathcal{M})$.
The set $T_{\b{p}}\mathcal{M}$ forms a vector space with operations:
\begin{align}
&(\b{V}+\b{W})(f) = \b{V}(f) + \b{W}(f),\,\, \forall \b{V}, \b{W} \in T_{\b{p}}\mathcal{M}, \\ 
&(a\b{V})(f) = a\,\b{V}(f), \quad \forall \b{V} \in T_{\b{p}}\mathcal{M}, \forall a \in \mathbb{R}.
\end{align}
In an $n$-dimensional space or manifold, a tangent space is an $n$-dimensional space.
\end{definition}

\begin{figure}[!h]
\centering
\includegraphics[width=2.7in]{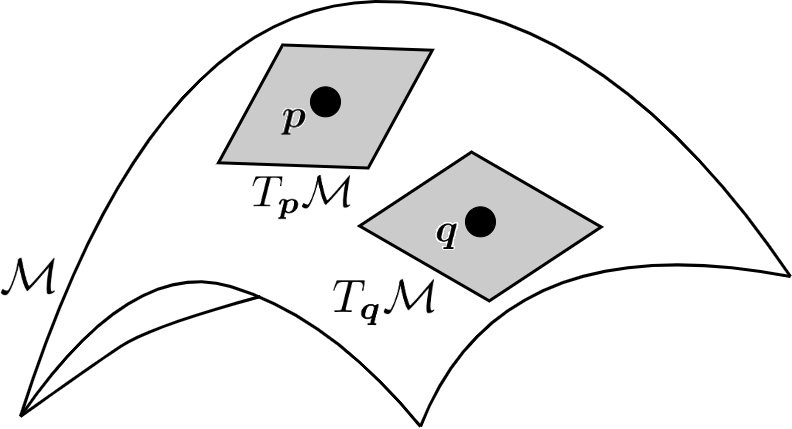}
\caption{Two tangent spaces $T_{\b{p}}\mathcal{M}$ and $T_{\b{q}}\mathcal{M}$ for two different points $\b{p}$ and $\b{q}$.}
\label{figure_tangent_spaces_different}
\end{figure}

\begin{remark}[Tangent spaces at different points]
Different points have different tangent spaces.
If the manifold is not flat, then:
\begin{align}
\b{p}, \b{q} \in \mathcal{M}, \b{p} \neq \b{q} \implies T_{\b{p}}\mathcal{M} \neq T_{\b{q}}\mathcal{M}.
\end{align}
An example of two different tangent spaces is illustrated in Fig. \ref{figure_tangent_spaces_different}.
\end{remark}

\subsection{Cotangent Space and Covector}

\begin{definition}[Covector (cotangent vector)]\label{definition_covector}
A \textbf{covector}, also called a \textbf{cotangent vector} or a \textbf{dual vector}, is a linear map that takes a vector and returns a real number. A covector acts on a vector and is defined as:
\begin{align}
\boxed{
W: T_{\b{p}}\mathcal{M} \rightarrow \mathbb{R},
}
\end{align}
such that:
\begin{align}
W(a\b{V} + b\b{Z}) = aW(\b{V}) + bW(\b{Z}),
\end{align}
where $\b{V}, \b{Z} \in T_{\b{p}}\mathcal{M}$ and $a,b \in \mathbb{R}$.
\end{definition}

\begin{definition}[Dual space]
Suppose we have a vector space $\mathcal{V}$.
The \textbf{dual space} of $\mathcal{V}$, denoted by $\mathcal{V}^*$, is the set of all linear maps from $\mathcal{V}$ to real numbers:
\begin{align}
\mathcal{V}^* := \{f: \mathcal{V} \rightarrow \mathbb{R} \mid f \text{ is linear}\}.
\end{align}
Every element of $\mathcal{V}^*$ takes a vector and returns a scalar number. 
\end{definition}

\begin{definition}[Cotangent space]
Let $\mathcal{M}$ be a smooth manifold and consider the point $\b{p}$ in the manifold, i.e., $\b{p} \in \mathcal{M}$.   
The \textbf{cotangent space} at $\b{p}$, denoted by $T_{\b{p}}^*\mathcal{M}$, is the set of all covectors at $\b{p}$. The cotangent space is the dual space of the tangent space:
\begin{align}
\boxed{
T_{\b{p}}^*\mathcal{M} := (T_{\b{p}}\mathcal{M})^*.
}
\end{align}
Every covector $W$ belongs to the cotangent space:
\begin{align}
\boxed{
W \in T_{\b{p}}^*\mathcal{M}.
}
\end{align}
\end{definition}

\subsection{Tangent and Cotangent Bundles}

\begin{definition}[Fiber]
Let \(\mathcal{M}\) be a manifold and let \(E\) be a space
together with a projection map:
\begin{align}
\pi : E \to \mathcal{M}.
\end{align}
For a point \(\b{p} \in \mathcal{M}\), the \textbf{fiber over}
\(\b{p}\) is defined as:
\begin{align}
\boxed{
F_{\b{p}} := \pi^{-1}(\b{p}) = \{e \in E \mid \pi(e)=\b{p}\}.
}
\end{align}
In other words, the fiber over \(\b{p}\) is the set of all
elements in the total space \(E\) that are attached to the
base point \(\b{p}\).
\end{definition}

\begin{definition}[Bundle (fiber bundle) {\cite{lee2013smooth, lee2018riemannian}}]
A \textbf{fiber bundle}, also called a \textbf{bundle} in short, is a quadruple:
\begin{align}
(E,\mathcal{M},\pi,F),
\end{align}
where:
\begin{itemize}
    \item \(E\) is the \textbf{total space},
    \item \(\mathcal{M}\) is the \textbf{base manifold},
    \item \(\pi : E \to \mathcal{M}\) is the \textbf{projection map},
    \item \(F\) is the \textbf{typical fiber}.
\end{itemize}
For every \(\b{p}\in\mathcal{M}\), the set:
\[
F_{\b{p}}:=\pi^{-1}(\b{p}),
\]
is called the fiber over \(\b{p}\), and each fiber
\(F_{\b{p}}\) is isomorphic to the typical fiber \(F\).

Intuitively, a fiber bundle is a space where, at every
point \(\b{p}\) of the base manifold \(\mathcal{M}\), another
space \(F_{\b{p}}\) is attached.
\end{definition}



\begin{definition}[Section of bundle]
A \textbf{section of bundle} is a rule that chooses one element from each fiber. Formally, a section is a map:
\begin{align}
s: \mathcal{M} \rightarrow E,
\end{align}
such that:
\begin{align}
\pi\big(s(\b{p})\big) = \b{p},
\end{align}
meaning that for every point $\b{p} \in \mathcal{M}$, the $s(\b{p})$ lies in the fiber over $\b{p}$.
The set of smooth sections of a bundle $E$ is denoted by:
\begin{align}\label{equation_Gamma_smooth_sections_of_bundle}
\Gamma(E) := \{ s : \mathcal{M} \to E \mid s \text{ is smooth} \}.
\end{align}
\end{definition}


\begin{definition}[Tangent bundle]
For a manifold $\mathcal{M}$, at each point $\b{p} \in \mathcal{M}$, we attach the tangent space $T_{\b{p}}\mathcal{M}$. The \textbf{tangent bundle}, denoted by $T\mathcal{M}$ is the union of all tangent spaces of all the points of manifold:
\begin{align}
\boxed{
T\mathcal{M} := \bigcup_{\b{p} \in \mathcal{M}} T_{\b{p}}\mathcal{M}. 
}
\end{align}
So, the tangent bundle is a fiber bundle where the base space is the manifold $\mathcal{M}$, the fiber at point $\b{p}$ is the tangent space at that point, $T_{\b{p}}\mathcal{M}$, and the total space is $T\mathcal{M}$. 
Each element of the tangent bundle is:
\begin{align}
(\b{p},\b{V}), \text{ where }\,\, \b{p} \in \mathcal{M}, \quad \b{V} \in T_{\b{p}}\mathcal{M}. 
\end{align}
The projection map of the tangent bundle is a map from the tangent bundle to the manifold:
\begin{align}
&\pi: T\mathcal{M} \to \mathcal{M}, \quad \pi(\b{p}, \b{V}) = \b{p}, \\
&\pi^{-1}(\b{p}) = T_{\b{p}} \mathcal{M}.
\end{align}
\end{definition}

\begin{definition}[Cotangent bundle]
For a manifold $\mathcal{M}$, at each point $\b{p} \in \mathcal{M}$, we attach the cotangent space $T_{\b{p}}^*\mathcal{M}$. The \textbf{cotangent bundle}, denoted by $T^*\mathcal{M}$ is the disjoint union\footnote{The regular union $\bigcup$ of sets just merges elements together, but the disjoint union $\bigsqcup$ keeps track of which set each element belongs to. Even if two sets share the same element, they are treated as different copies in disjoint union. For example, two different points could have vectors that \textit{look the same} as numbers, but they are different vectors in different spaces.} of all cotangent spaces of all the points of manifold:
\begin{align}
\boxed{
T^*\mathcal{M} := \bigsqcup_{\b{p} \in \mathcal{M}} T_{\b{p}}^* \mathcal{M} = \{ (\b{p}, W) \mid \b{p} \in \mathcal{M}, \, \omega \in T_{\b{p}}^* \mathcal{M} \}.
}
\end{align}
So, the cotangent bundle is a fiber bundle where the base space is the manifold $\mathcal{M}$, the fiber at point $\b{p}$ is the cotangent space at that point, $T_{\b{p}}^*\mathcal{M}$, and the total space is $T^*\mathcal{M}$. 
Each element of the cotangent bundle is:
\begin{align}
(\b{p},\omega), \text{ where }\,\, \b{p} \in \mathcal{M}, \quad \omega \in T_{\b{p}}^*\mathcal{M}. 
\end{align}
The projection map of the cotangent bundle is a map from the cotangent bundle to the manifold:
\begin{align}
&\pi: T^*\mathcal{M} \to \mathcal{M}, \quad \pi(\b{p}, \omega) = \b{p}, \\
&\pi^{-1}(\b{p}) = T_{\b{p}}^* \mathcal{M}.
\end{align}
\end{definition}

\subsection{Vector and Vector Field}\label{section_vector_and_vector_field}

\begin{definition}[Vector]\label{definition_vector}
For a smooth manifold $\mathcal{M}$, every point $\b{p} \in \mathcal{M}$ has its own tangent space $T_{\b{p}}\mathcal{M}$. A \textbf{vector} at point $\b{p}$ is a tangent vector at point $\b{p}$, so it belongs to the tangent space at point $\b{p}$:
\begin{align}
\boxed{
\b{V} \in T_{\b{p}}\mathcal{M}.
}
\end{align}
So, the vector lives in the tangent space at that specific point.
In an $n$-dimensional space or manifold, a vector is an $n$-dimensional vector.
\end{definition}

Intuitively, a vector is an $n$-dimensional arrow, with some length and direction, starting from point $\b{p}$ on the manifold $\mathcal{M}$ where the vector is tangent to the manifold at $\b{p}$. An example vector at point $\b{p} \in \mathcal{M}$ is illustrated in Fig. \ref{figure_tangent_vector} (see Section \ref{section_tangent_space_tangent_vector}).

\begin{definition}[Identity Map]
Let $\mathcal{M}$ be a smooth manifold. The \textbf{identity map} on $\mathcal{M}$ is the function:
\begin{align}
\mathrm{id}_\mathcal{M} : \mathcal{M} \to \mathcal{M},
\end{align}
defined by:
\begin{align}
\boxed{
\mathrm{id}_\mathcal{M}(\b{p}) = \b{p}, \quad \forall \b{p} \in \mathcal{M}.
}
\end{align}
\end{definition}

\begin{definition}[Vector field]\label{definition_vector_field}
Let $\mathcal{M}$ be a smooth manifold. A \textbf{smooth vector field} on $\mathcal{M}$ is a smooth map:
\begin{align}
\boxed{
\b{X} : \mathcal{M} \to T\mathcal{M},
}
\end{align}
such that:
\begin{align}
\pi \circ \b{X} = \mathrm{id}_\mathcal{M},
\end{align}
where $\pi : T\mathcal{M} \to \mathcal{M}$ is the projection from the tangent bundle to manifold $\mathcal{M}$.
The vector field assigns a vector to every point $\b{p} \in \mathcal{M}$:
\begin{align}
\boxed{
\b{p} \mapsto \b{X}(\b{p}) \in T_{\b{p}}\mathcal{M}, \quad \forall \b{p} \in \mathcal{M}.
}
\end{align}

The \textit{space of smooth vector fields on $\mathcal{M}$} is denoted by $\mathfrak{X}(\mathcal{M})$ and is defined as:
\begin{equation}
\boxed{
\begin{aligned}
\mathfrak{X}(\mathcal{M}) := \{\, &\b{X} : \mathcal{M} \to T\mathcal{M} \mid \b{X} \text{ is smooth and } \\
&\b{X}(\b{p}) \in T_{\b{p}}\mathcal{M} \text{ for all } \b{p} \in \mathcal{M} \,\}.
\end{aligned}
}
\end{equation}
\end{definition}

\begin{remark}[Vector field versus tangent vector]
Note that $\b{X} \in \mathfrak{X}(\mathcal{M})$ is a vector field, but when we apply the vector field on a point, it becomes a tangent vector at that point. In other words, $\b{X}(\b{p}) \in T_{\b{p}}\mathcal{M}$ is a tangent vector at point $\b{p}$ (existing in the tangent space at point $\b{p}$). Therefore, in summary:
\begin{align}
\text{Vector field: } &\boxed{\b{X} \in \mathfrak{X}(\mathcal{M}),} \\
\text{Tangent vector: } &\boxed{\b{X}(\b{p}) \in T_{\b{p}}\mathcal{M}.}
\end{align}
\end{remark}

A vector field is a mathematical concept where a vector (an arrow with magnitude and direction) is assigned to every point in a space or manifold. 
In an $n$-dimensional space or manifold, a vector field assigns an $n$-dimensional vector to each $n$-dimensional point. 

\begin{figure}[!h]
\centering
\includegraphics[width=2.7in]{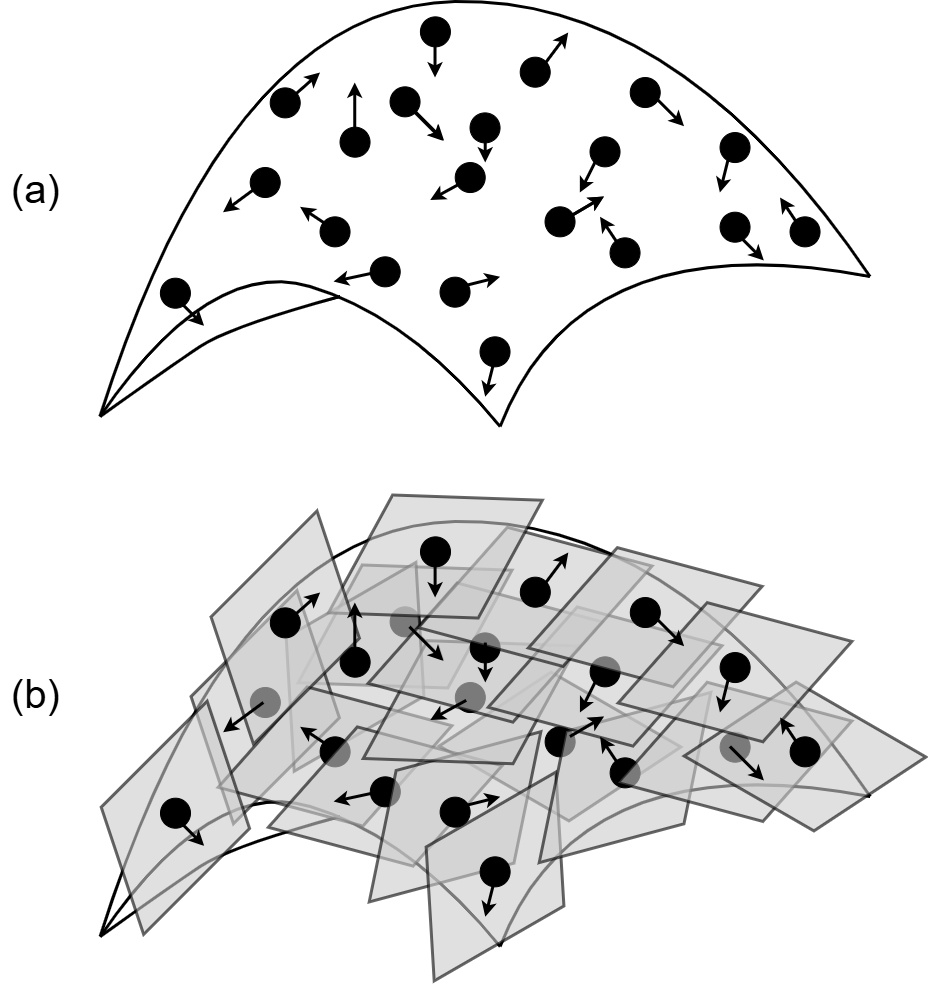}
\caption{Vector field: (a) vectors at different points of manifold where the vectors may change from point to point, and (b) the tangent spaces at different points of the manifold where each tangent vector lies in its tangent space. Note that for the sake of visualization, we are showing a few points on the manifold. The vector field has a vector measurement at ``every" point on the manifold and every point has its own tangent space.}
\label{figure_vector_field}
\end{figure}

Intuitively, if we have a vector measurement on different points of the manifold, where the vectors may change from point to point, we have a vector field on the manifold. For example, Fig. \ref{figure_vector_field} depicts a two-dimensional vector field on a locally two-dimensional manifold. 

\begin{remark}[Relation of vector field and tangent bundle]
A vector field $\b{X}: \mathcal{M} \rightarrow T\mathcal{M}$, defined in Definition \ref{definition_vector_field}, assigns a tangent vector to each point: $\b{X}: \b{p} \mapsto X(\b{p}) = \b{V} \in T_{\b{p}}\mathcal{M}$.
So, a vector field $\b{X}(\mathcal{M})$ belongs to the set of smooth sections of the tangent bundle:
\begin{align}
\boxed{
\b{X}(\mathcal{M}) \in \Gamma(T\mathcal{M}),
}
\end{align}
where $\Gamma(.)$ is defined in Eq. (\ref{equation_Gamma_smooth_sections_of_bundle}). 
\end{remark}

\begin{definition}[Pointwise definition of vector field]\label{definition_poinstwise_vector_field}
Let $\mathcal{M}$ be a smooth manifold, $f \in C^\infty(\mathcal{M})$, and $\b{X} \in \mathfrak{X}(\mathcal{M})$. 
The product $f\b{X}$ is the \textbf{vector field defined pointwise} by:
\begin{align}\label{equation_poinstwise_vector_field}
(f\b{X})_{\b{p}} := f(\b{p})\,\b{X}({\b{p}}), \quad \forall \b{p} \in \mathcal{M},
\end{align}
where $\b{X}({\b{p}}) \in T_{\b{p}}\mathcal{M}$ and the multiplication on the right-hand side 
is the usual scalar multiplication in the tangent space $T_{\b{p}}\mathcal{M}$.
\end{definition}

\subsection{Tensor}



\begin{definition}[Covariant and contravariant indices]\label{definition_covariant_and_contravariant_indices}
In differential geometry, when an index is upper (superscript) index, it is called a \textbf{contravariant index}.
When an index is lower (subscript) index, it is called a \textbf{covariant index}.
\end{definition}

\begin{definition}[Linear map]
A map or function $f$ is \textbf{linear} if:
\begin{align}
f(a\b{V} + b\b{Z}) = af(\b{V}) + bf(\b{Z}),
\end{align}
for vectors $\b{V}, \b{Z}$ and scalars $a, b$.
\end{definition}

\begin{definition}[Multiplinear map]
A map or function is \textbf{multilinear} if it is linear in each argument separately.
Let $\mathcal{V}_1, \mathcal{V}_2, \dots, \mathcal{V}_k, \mathcal{V}_{k+1}$ be vector spaces.
A map:
\begin{align}
\b{T} : \mathcal{V}_1 \times \mathcal{V}_2 \times \cdots \times \mathcal{V}_k \to \mathcal{V}_{k+1},
\end{align}
is called a multilinear map if it is linear in each argument separately.  
That is, for every $i \in \{1,\dots,k\}$, for all $\b{V}_j \in \mathcal{V}_j$ $(j \neq i)$, for all $\b{U}_i, \b{Z}_i \in \mathcal{V}_i$, and for all scalars $a,b \in \mathbb{R}$, we have:
\begin{equation}
\begin{aligned}
T(\b{V}_1,\dots,&\, a \b{U}_i + b \b{Z}_i,\dots, \b{V}_k) \\
&= a\,T(\b{V}_1,\dots, \b{U}_i,\dots, \b{V}_k) \\
&~~~~+ b\,T(\b{V}_1,\dots, \b{Z}_i,\dots, \b{V}_k).
\end{aligned}
\end{equation}
\end{definition}


\begin{definition}[Tensor]\label{definition_tensor}
A \textbf{tensor} is a multilinear object that can take several vectors and/or covectors as input and return a number \cite{kobayashi1996foundations,bishop2012tensor}.
More formally, a tensor is a multilinear map built from:
\begin{itemize}
\item vectors in the tangent space $T_{\b{p}}\mathcal{M}$, and/or
\item covectors in the cotangent space $T_{\b{p}}^*\mathcal{M}$.
\end{itemize}
So, tensors generalize both vectors and covectors.

A tensor of type $(r,s)$, also called $(r,s)$ tensor, has $r$ contravariant (upper) indices and $s$ covariant (lower) indices (see Definition \ref{definition_covariant_and_contravariant_indices}). So, it acts on $r$ covectors and $s$ vectors and produces a number:
\begin{align}
\boxed{
\b{T}: (T_{\b{p}}^*\mathcal{M})^r \times (T_{\b{p}}\mathcal{M})^s \rightarrow \mathbb{R},
}
\end{align}
which means\footnote{It feels logical to match ``vector indices" with ``vector inputs", but it actually works the opposite way because tensors are functions whose output is a scalar number. To get a scalar out of the tensor, we have to cancel out the basis vectors. So, $r$ contravariant (upper) indices are needed to cancel out $r$ covectors in the cotangent space to output a scalar number. Likewise, to get a scalar out of the machine, we have to cancel out the basis vectors. So, $s$ covariant (lower) indices are needed to cancel out $s$ vectors in the tangent space to output a scalar number.}:
\begin{align*}
\b{T}: \underbrace{T_{\b{p}}^*\mathcal{M} \times \cdots \times T_{\b{p}}^*\mathcal{M}}_{r \text{ times}} \times \underbrace{T_{\b{p}}\mathcal{M} \times \cdots \times T_{\b{p}}\mathcal{M}}_{s \text{ times}} \rightarrow \mathbb{R}.
\end{align*}
A $(r,s)$ tensor is denoted by:
\begin{align}
\boxed{
T_{j_1, \dots, j_s}^{i_1, \dots, i_r},
}
\end{align}
having $r$ contravariant (upper) indices and $s$ covariant (lower) indices\footnote{In simple words, a tensor of type $(r,s)$ has $r$ contravariant (upper) indices and $s$ covariant (lower) indices.}. 
\end{definition}

\begin{remark}[Geometrical interpretation of $(r,s)$ tensor]
A $(r,s)$ tensor has $r$ contravariant slots and $s$ covariant slots. One can imagine that a $(r,s)$ tensor, in an $n$-dimensional manifold, is a geometrical object with $(r+s)$ grids, where each grid is an $n$-dimensional array. 
So, a $(r,s)$ tensor can be imagined to have $n^{(r+s)}$ cells where each cell has a number in it. 
For example, in a four-dimensional space-time manifold of general relativity, a metric tensor $g_{ij}$, which can be introduced later, is a tensor of type $(0,2)$. This tensor can be imagined as a grid with $4^{(0+2)} = 16$ cells. Such a grid can be represented as a $(4 \times 4)$ matrix. 
\end{remark}

\begin{remark}[Interpretation of $(r,s)$ tensor in pure mathematics]
Pure mathematics interprets a $(r,s)$ tensor as follows. If we feed $r$ covectors $\omega_1, \dots, \omega_r \in T_{\b{p}}^*\mathcal{M}$ and $s$ vectors $\b{V}_1, \dots, \b{V}_s \in T_{\b{p}}\mathcal{M}$ to a $(r,s)$ tensor, you get a scalar:
\begin{align}
\b{T}(\omega_1, \dots, \omega_r, \b{V}_1, \dots, \b{V}_s) \in \mathbb{R}.
\end{align}
In other words, $(r,s)$ tensor $T_{j_1, \dots, j_s}^{i_1, \dots, i_r}$ takes $r$ covectors and $s$ vectors and outputs a scalar number.
Some examples are:
\begin{itemize}
\item $(0,0)$ tensor $T$ has no input and is just a scalar number. 
\item $(r,0)$ tensor $T{i_1, \dots, i_r}$ takes $r$ covectors and outputs a scalar number.
\item $(0,s)$ tensor $T_{j_1, \dots, j_s}$ takes $s$ vectors and outputs a scalar number.
\end{itemize}
\end{remark}

\begin{remark}[Interpretation of $(r,s)$ tensor in physics]
Physics interprets a $(r,s)$ tensor as follows. 
In physics, people often ``freeze" some slots and see what the tensor does in the remaining slots. 
In other words, in physics, we often interpret contravariant (upper) indices as outputs (vectors) and covariant (lower) indices as inputs (vectors that we act on).

In a $(r,s)$ tensor $T_{j_1, \dots, j_s}^{i_1, \dots, i_r}$, if we fix some inputs, the remaining inputs can produce vectors/covectors in a natural way.
For example, a $(r,s)$ tensor $T_{j_1, \dots, j_s}^{i_1, \dots, i_r}$ can act on $s$ vectors to produce a tensor of type $(r,0)$, i.e., something with $r$ contravariant slots; in fact, it takes $s$ vectors and produces an object with $r$ vector-like directions. 
Some examples are:
\begin{itemize}
\item $(1,1)$ tensor $T_{j}^{i}$ can act on a vector to produce another vector. An example in physics is a linear transformation, e.g., stress tensor or the Jacobian of a map.
\item $(0,1)$ tensor $T_{j}$ or $(1,0)$ tensor $T^{i}$ can act on a vector to produce a scalar. An example in physics is that the input is a vector $\b{V}$ like motion vector, and the output is a scalar rate of change along $\b{V}$. 
\item $(1,1)$ tensor $T_{j}^i$ can act on a vector to produce another vector. An example in physics is that the input is a vector $\b{V}$ like direction of surface normal, and the output is a force vector. Another example is the Jacobian $J_j^i = \partial x^i / \partial y^j$ where the input is a vector in $y$ coordinates and the output is a vector in $x$ coordinates. 
\item $(0,2)$ tensor $T_{j_1, j_2}$ is a bilinear map which can act on two vectors to produce a scalar. An example is the inner product (or metric tensor) that the inputs are two vectors, and the output is the scalar output of inner product. 
\end{itemize}
\end{remark}

\begin{remark}[Examples for tensor]
Some examples of tensor are:
\begin{itemize}
\item $(1,0)$ tensor in three-dimensional manifold:
\begin{align*}
\b{T} = T^i =
\begin{bmatrix}
4\\
5\\
3
\end{bmatrix},
\end{align*}
where $T^1 = 4, T^2 = 5, T^3 = 3$.
\item $(0,1)$ tensor in three-dimensional manifold:
\begin{align*}
\b{T} = T_i = [4, 5, 3],
\end{align*}
where $T_1 = 4, T_2 = 5, T_3 = 3$.
\item $(1,1)$ tensor in three-dimensional manifold:
\begin{align*}
\b{T} = T^i_j = 
\begin{bmatrix}
4 & 2.5 & 6\\
5 & 1 & 4\\
3 & 2 & 7
\end{bmatrix},
\end{align*}
where $T^1_1 = 4, T^1_2 = 5, T^1_3 = 3, T^2_1 = 2.5, T^2_2 = 1, T^2_3 = 2, T^3_1 = 6, T^3_2 = 4, T^3_3 = 7$.
\end{itemize}
\end{remark}

\begin{definition}[Rank of tensor]
\textbf{Rank} of a tensor of type $(r,s)$ is defined as:
\begin{align}
\boxed{
\text{rank}(T_{j_1, \dots, j_s}^{i_1, \dots, i_r}) = r + s.
}
\end{align}
\end{definition}

\begin{definition}[$k$-form]\label{definition_k_form}
Let \(\mathcal{M}\) be a smooth manifold and let
\(\b{p}\in\mathcal{M}\). A \textbf{\(k\)-form} at \(\b{p}\) is a
totally antisymmetric covariant tensor of rank \(k\), that
is, a multilinear map (a $(0,k)$ tensor):
\begin{align}
\omega: \underbrace{T_{\b{p}}\mathcal{M} \times \cdots \times T_{\b{p}}\mathcal{M}}_{k \text{ times}} \rightarrow \mathbb{R},
\end{align}
or equivalently:
\begin{align}
\boxed{
\omega: (T_{\b{p}}\mathcal{M})^k \rightarrow \mathbb{R},
}
\end{align}
where antisymmetric means that swapping any two input vectors flips the sign:
\begin{equation}
\begin{aligned}
\omega(\b{V}_1, \dots, &\b{V}_i, \dots, \b{V}_j, \dots, \b{V}_k) = \\
&- \omega(\b{V}_1, \dots, \b{V}_j, \dots, \b{V}_i, \dots, \b{V}_k).
\end{aligned}
\end{equation}
In other words, A $k$-form on a smooth manifold is a totally antisymmetric covariant tensor of rank $k$.
\end{definition}

\begin{remark}[Covector as $1$-form]
According to Definitions \ref{definition_covector} and \ref{definition_k_form}, a covector is a $1$-form. 
\end{remark}

\begin{remark}[Tensors are coordinate-independent geometric objects]\label{remark_tensor_coordinate_independent} 
A key property of a tensor is that it represents a geometric object that is independent of the choice of coordinate system (or the reference frame in the language of physics). When the coordinate system changes, the components of the tensor transform according to specific rules so that the underlying tensor itself remains unchanged. In other words, as long as the manifold is not altered, the tensor represents the same geometric object regardless of the coordinate system, although its components may vary depending on the chosen coordinates. 

In summary, the tensor itself is a geometrical object that does not depend on the coordinate system. However, the components of the tensor depend on the chosen coordinate system and transform according to the tensor transformation laws---which will be presented later in Corollary \ref{corollary_tensor_transformation_laws}.

To better understand, think of tensor as a ``real" object in the manifold independent of any coordinates. For example, the metric tensor (which we will introduce later) defines distances and angles on the manifold.
Or the curvature tensor $R$ (which we will introduce later) defines how the manifold bends. However, components of a tensor are numbers that describe the tensor relative to a specific basis (coordinate system).
If you change the basis (coordinate system), the numbers change, but the underlying object does not.
For an analogy, assume that the tensor is like the actual arrow in space, the basis are the coordinate axes, and the components are projections of the arrow onto the axes. 
\end{remark}

\begin{remark}[Use of tensors in general relativity]
According to Remark \ref{remark_tensor_coordinate_independent}, tensors are useful for representing the space-time manifold in general relativity because, for example, gravity depends only on the curvature of space-time manifold and not the choice of reference frame or coordinate system \cite{susskind2025general}. That is why Einstein used tensors in his gravitational field equations of general relativity \cite{einstein1915feldgleichungen}.
\end{remark}

\subsection{Tensor Product}

\begin{definition}[Tensor product \cite{ryan2002introduction}]
Suppose we have two tensors, i.e., tensor $T$ of type $(r,s)$ and tensor $S$ of type $(k,l)$:
\begin{align*}
&T: (T_{\b{p}}^*\mathcal{M})^r \times (T_{\b{p}}\mathcal{M})^s \rightarrow \mathbb{R}, \\
&S: (T_{\b{p}}^*\mathcal{M})^k \times (T_{\b{p}}\mathcal{M})^l \rightarrow \mathbb{R}.
\end{align*}
Then, their \textbf{tensor product}, denoted by $\otimes$, is defined as:
\begin{align}
\boxed{
T \otimes S: (T_{\b{p}}^*\mathcal{M})^{r+k} \times (T_{\b{p}}\mathcal{M})^{s+l} \rightarrow \mathbb{R},
}
\end{align}
which is a tensor of type $(r+k, s+l)$, with $r+k$ covariant indices and $s+l$ contravariant indices. 

For covectors $\omega_1, \dots, \omega_{r+k}$ and vectors $\b{V}_1, \dots, \b{V}_{s+l}$, there is:
\begin{align*}
(T \otimes S)&(\omega_1, \dots, \omega_{r+k}, \b{V}_1, \dots, \b{V}_{s+l}) \\
&= T(\omega_1, \dots, \omega_{r}, \b{V}_1, \dots, \b{V}_{s})\, \cdot \\
&~~~~~T(\omega_{r+1}, \dots, \omega_{r+k}, \b{V}_{s+1}, \dots, \b{V}_{s+l})
\end{align*}
\end{definition}

\begin{remark}[Geometrical interpretation]
We can think of tensor product as ``stacking grids":
the tensor $T$ of type $(r,s)$ is a $(r+s)$ dimensional grid and the tensor $S$ of type $(k,l)$ is a $(k+l)$ dimensional grid. 
The tensor product $T \otimes S$ is a new $(r+s+k+l)$ dimensional grid.
The new tensor has all input directions from both tensors, with outputs combined multiplicatively.
\end{remark}

\begin{remark}[Components view of tensor product]
If tensor $T$ of type $(r,s)$ has components $T_{j_1, \dots, j_s}^{i_1, \dots, i_r}$ and tensor $S$ of type $(k,l)$ has components $T_{q_1, \dots, q_l}^{p_1, \dots, p_k}$, then the components of the tensor product $T \otimes S$ is:
\begin{align}\label{equation_tensor_product_coordinate_based}
\boxed{
(T \otimes S)_{j_1, \dots, j_s, q_1, \dots, q_l}^{i_1, \dots, i_r, p_1, \dots, p_k} = T_{j_1, \dots, j_s}^{i_1, \dots, i_r}\, S_{q_1, \dots, q_l}^{p_1, \dots, p_k}.
}
\end{align}
\end{remark}

\begin{remark}[Example for tensor product]
An example of the tensor product, in a two-dimensional manifold, is as follows:
\begin{align*}
&\b{A} = A^i_j =
\begin{bmatrix}
a_{11} & a_{12}\\
a_{21} & a_{22}
\end{bmatrix}
, \quad \b{B} = B^k =
\begin{bmatrix}
b_1 \\
b_2 
\end{bmatrix}
,
\end{align*}
where $i,j,k \in \{1,2\}$ in the two-dimensional space. 
The tensor product of these tensors is:
\begin{align*}
\b{C} &= C^{i,k}_j = \b{A} \otimes \b{B} = A^i_j B^k \\
&=
\begin{bmatrix}
a_{11} 
\begin{bmatrix}
b_1 \\
b_2 
\end{bmatrix}
& a_{12}
\begin{bmatrix}
b_1 \\
b_2 
\end{bmatrix}
\\ \\
a_{21} 
\begin{bmatrix}
b_1 \\
b_2  
\end{bmatrix}
& a_{22}
\begin{bmatrix}
b_1 \\
b_2 
\end{bmatrix}
\end{bmatrix}
=
\begin{bmatrix}
a_{11} b_1 & a_{12} b_1\\
a_{11} b_2 & a_{12} b_2\\
a_{21} b_1 & a_{12} b_1\\
a_{21} b_2 & a_{12} b_2
\end{bmatrix},
\end{align*}
where $i,j,k \in \{1,2\}$ in the tensor product $C^{i,k}_j = A^i_j B^k$ so the result tensor $\b{C}$ has $2\times 2 \times = 8$ elements, no matter in what arrangement we put these elements in a matrix. 
\end{remark}

\begin{definition}[Contraction of indices]
In tensor product, when an index appears once as upper index and once as lower index---so it is summed over---that index goes away in the output of summation. This procedure is called \textbf{contraction of index}; in other words, that index is contracted. For example, in the following expression, the we have contraction of indices $j$ and $\ell$:
\begin{align*}
A_{ij\ell} B^{jk} C_m^{\ell} = \sum_{j} \sum_{\ell} A_{ij\ell} B^{jk} C_i^{\ell} = D^{k}_{im},
\end{align*}
where Einstein convention is used. 
\end{definition}


\begin{lemma}[Index substitution by Kronecker delta]
Let $A^i$ and $A_j$ be the components of a vector and a covector, respectively. The Kronecker delta $\delta_i^j$ acts as an \textbf{index substitution} operator such that:
\begin{equation}\label{equation_index_substitution_delta}
\boxed{
\delta_i^j A^i = A^j \quad \text{and} \quad \delta_i^j A_j = A_i.
}
\end{equation}
\end{lemma}

\begin{proof}
By the definition of the Einstein summation convention and the Kronecker delta:
\begin{equation*}
    \delta_i^j A^i = \sum_{i=1}^n \delta_i^j A^i = \delta_1^j A^1 + \dots + \delta_j^j A^j + \dots + \delta_n^j A^n.
\end{equation*}
Since $\delta_i^j = 1$ if $i=j$ and $0$ otherwise, all terms in the sum vanish except for the term where $i=j$:
\begin{equation*}
    \delta_i^j A^i = 1 \cdot A^j = A^j.
\end{equation*}
The proof for $A_j$ follows an identical substitution logic.
\end{proof}

\subsection{Coordinate Basis for Vectors and Covectors}

\begin{definition}[Coordinate system on manifold]
Let $\mathcal{M}$ be an $n$-dimensional smooth manifold. A \textbf{coordinate system} or \textbf{curvilinear coordinate system} on $\mathcal{M}$ is a collection of smooth functions:
\begin{align}
x^1, x^2, \dots, x^n : \mathcal{M} \to \mathbb{R},
\end{align}
such that the map:
\begin{align}
\boxed{
x = (x^1,\dots,x^n) : \mathcal{M} \to \mathbb{R}^n,
}
\end{align}
is a diffeomorphism onto an open subset of $\mathbb{R}^n$. Each $x^i$ is called a coordinate function. 

Intuitively, a coordinate system is generally a curvilinear coordinate system which covers the entire manifold. 
\end{definition}

\begin{definition}[Regular partial derivatives]
It is possible to calculate the rate of change of a vector or a function along the direction of a coordinate function, using \textbf{regular partial derivatives}. For example:
\begin{itemize}
\item $\frac{\partial f}{\partial x^i}$ is the derivative of smooth function $f$ with respect to the coordinate $x^i$. It gives the rate of change of function $f$ along the coordinate $x^i$.
\item $\frac{\partial V}{\partial x^i}$ is the derivative of vector $V$ with respect to the coordinate $x^i$. It gives the rate of change of vector $V$ along the coordinate $x^i$.
\end{itemize}
We define the notation:
\begin{align}\label{equation_partial_i}
\boxed{
\partial_i := \frac{\partial }{\partial x^i}.
}
\end{align}
\end{definition}

\begin{definition}[Coordinate bases for vectors and covectors]
Consider a coordinate system $x = (x^1, \dots, x^n)$ on an $n$-dimensional smooth manifold $\mathcal{M}$.

\begin{itemize}

\item \textbf{Coordinate basis} for \underline{vectors}: the associated coordinate vector fields:
\begin{align}\label{equation_bases_for_vectors}
\boxed{
\frac{\partial}{\partial x^1}, \dots, \frac{\partial}{\partial x^n},
}
\end{align}
form a basis of the tangent space $T_{\b{p}}\mathcal{M}$ at each point $\b{p} \in \mathcal{M}$. 
In other words, the bases of vectors (in the tangent space) are listed in Eq. (\ref{equation_bases_for_vectors}). 
The partial derivative $\partial / \partial x^i$ means the rate of change along the coordinate $x^i$.
According to Eq. (\ref{equation_partial_i}), the coordinate basis for vectors can be stated as:
\begin{align}\label{equation_bases_for_vectors_2}
\boxed{
\partial_1, \dots, \partial_n.
}
\end{align}

\item \textbf{Coordinate basis} for \underline{covectors}: The differentials of the coordinate functions:
\begin{align}\label{equation_bases_for_covectors}
\boxed{
dx^1, \dots, dx^n,
}
\end{align}
form the dual basis of the cotangent space $T_{\b{p}}^*\mathcal{M}$, satisfying:
\begin{align}\label{equation_dx_partial_partialx_delta}
\boxed{
dx^i\!\left(\frac{\partial}{\partial x^j}\right)=\delta^i_j.
}
\end{align}
In other words, the bases of covectors (in the cotangent space) are listed in Eq. (\ref{equation_bases_for_covectors}). 
The differential $dx^i$ means the infinitesimal change along the coordinate $x^i$.
Note that the basis for covectors is also called the \underline{dual basis}.

\end{itemize}
\end{definition}
\begin{proof}
Equation (\ref{equation_dx_partial_partialx_delta}) is by the chain rule:
\begin{align*}
& dx^i\!\left(\frac{\partial}{\partial x^i}\right) = 1, \\
& dx^i\!\left(\frac{\partial}{\partial x^j}\right) = 0, \quad \forall j \neq i.
\end{align*}
\end{proof}

Let $\{\b{e}_1, \b{e}_2, \dots, \b{e}_n\}$ denote the basis vectors of the $n$-dimensional coordinate system on the manifold. 
The coordinate system can be any coordinate system, such as a curvilinear coordinate system.  

\begin{definition}[Coordinate basis for the coordinate system on manifold]\label{definition_coordinate_basis}
In a coordinate system on a manifold, the coordinate basis vectors $\{\b{e}_i\}_{i=1}^n$ are coordinate basis for vectors, as defined in Eq. (\ref{equation_bases_for_vectors}). 
Therefore, the \textbf{coordinate basis vectors} $\{\b{e}_i\}_{i=1}^n$ are defined as the partial derivative with respect to the coordinates:
\begin{align}\label{equation_coordinate_basis_vectors}
\boxed{\b{e}_i := \frac{\partial}{\partial x^i} \overset{(\ref{equation_partial_i})}{=} \partial_i,} \quad \forall i \in \{1, \dots, n\}.
\end{align}
In other words, every basis vector $\b{e}_i$ determines how much of change some quantity can have with respect to the coordinate $x^i$. 
\end{definition}

\begin{lemma}[Coordinate basis transformation]
Consider a coordinate system $x$ (with basis vectors $\{\b{e}_i\}_{i=1}^n$) on a manifold. If we use another coordinate system $y$ (with basis vectors $\{\widetilde{\b{e}}_i\}_{i=1}^n$) on the same manifold, the relation of $i$-th basis vector in coordinate system $y$ and the $j$-th basis vector in coordinate system $x$ is:
\begin{align}\label{equation_coordinate_basis_transformation}
\boxed{
\widetilde{\b{e}}_i = \frac{\partial x^j}{\partial y^i} \b{e}_j.
}
\end{align}
\end{lemma}
\begin{proof}
According to Eq. (\ref{equation_coordinate_basis_vectors}), in the coordinate systems $x$ and $y$, we have:
\begin{align}
&\b{e}_i = \frac{\partial}{\partial x^i}, \label{equation_ej_partial} \\
&\widetilde{\b{e}}_i = \frac{\partial}{\partial y^i}, \label{equation_ej_tilde_partial}
\end{align}
respectively. 
According to the chain rule in derivatives, we have:
\begin{align*}
\frac{\partial}{\partial y^i} \overset{(a)}{=} \frac{\partial x^j}{\partial y^i} \frac{\partial}{\partial x^j} \overset{(b)}{\implies} \widetilde{\b{e}}_i = \frac{\partial x^j}{\partial y^i} \b{e}_j,
\end{align*}
where $(a)$ is because of the chain rule in derivatives and $(b)$ is because of Eqs. (\ref{equation_ej_partial}) and (\ref{equation_ej_tilde_partial}).
\end{proof}

\subsection{Contravariant and Covariant Components}

\subsubsection{Definition of Contravariant and Covariant Components}

We defined contravariant and covariant indices as upper (superscript) indices and lower (subscript) indices, respectively. Now, we define the contravariant and covariant components. Contravariant and covariant components are the components of a tensor expressed with respect to a basis. We will clarify this in the following. 


\begin{definition}[Contravariant components]
Consider an $n$-dimensional vector $\b{V} := [V^1, \dots, V^n]^\top \in T_{\b{p}}\mathcal{M}$ in an $n$-dimensional manifold $\mathcal{M}$. The $\{V^1, \dots, V^n\}$ are called the \textbf{contravariant components}. The contravariant components are denoted by upper (superscript) indices, as also mentioned in Definition \ref{definition_covariant_and_contravariant_indices}. The vector $\b{V} \in T_{\b{p}}\mathcal{M}$ can be stated as a linear combination of the basis vectors:
\begin{align}\label{equation_contravariant_components}
\boxed{
\b{V} = V^i \b{e}_i = V^1 \b{e}_1 + V^2 \b{e}_2 + \dots + V^n \b{e}_n,
}
\end{align}
with contravariant components as the coefficients. 
According to Eq. (\ref{equation_coordinate_basis_vectors}), the Eq. (\ref{equation_contravariant_components}) can be stated as:
\begin{align}\label{equation_contravariant_components_2}
\boxed{
\b{V} = V^i \frac{\partial}{\partial x^i} = V^i \partial_i = V^1 \partial_1 + V^2 \partial_2 + \dots + V^n \partial_n.
}
\end{align}
\end{definition}

\begin{remark}
To better understand---roughly speaking---the regular components of a vector, that we learned in high school, are the contravariant components. 
\end{remark}

\begin{remark}[Notation for coordinates and components]
We denote coordinate functions by lowercase letters, such as $x^i$ or $y^i$. Components of vectors and tensors are denoted by uppercase (or sometimes lowercase) symbols, such as $V^i$ or $W^i$.
In particular, coordinates $x^i$ should not be confused with components of vector fields.
\end{remark}

\begin{remark}[Coordinate basis for the tangent space]
The basis vectors in Eq. (\ref{equation_bases_for_vectors}) or (\ref{equation_bases_for_vectors_2}) or (\ref{equation_coordinate_basis_vectors}) exist in the tangent space at every point $\b{p} \in \mathcal{M}$. Thus, $\{\partial_1, \dots, \partial_n\}$ are also the basis vectors for the tangent space. As a result, a vector $\b{V} \in T_{\b{p}}\mathcal{M}$ can be stated as in Eq. (\ref{equation_contravariant_components_2}).
\end{remark}


\begin{definition}[Dual basis for the coordinate system on a manifold]
Let $\mathcal{M}$ be a smooth manifold and let $\{\b{e}_i\}_{i=1}^n$ be a basis for the tangent space $T_{\b{p}} \mathcal{M}$ at a point $\b{p} \in \mathcal{M}$.  
The \textbf{dual basis} $\{\b{e}^i\}_{i=1}^n$ is the set of covectors in the cotangent space $T_{\b{p}}^* \mathcal{M}$ satisfying:
\begin{align}
\boxed{
\b{e}^i(\b{e}_j) = \delta^i_j,
}
\end{align}
where $\delta^i_j$ is the Kronecker delta, defined in Eq. (\ref{equation_Kronecker_delta_lower_upper}).  

For any vector $\b{V} \in T_{\b{p}} \mathcal{M}$, we can write $\b{V} = V^i \, \b{e}_i$ according to Eq. (\ref{equation_contravariant_components}), and the contravariant components $V^i$ are recovered using the dual basis as:
\begin{align}
\boxed{
V^i = \b{e}^i(\b{V}).
}
\end{align}
\end{definition}

\begin{remark}[Coordinate expression of a $(r,s)$ tensor]
Let $\b{T}$ be a tensor of type $(r, s)$. It is expressed in coordinate bases as:
\begin{align}
\boxed{
\b{T} = T^{i_1 \dots i_r}_{j_1 \dots j_s} \b{e}_{i_1} \otimes \dots \otimes \b{e}_{i_r} \otimes \b{e}^{j_1} \otimes \dots \otimes \b{e}^{j_s},
}
\end{align}
where $\otimes$ denotes tensor product and $\{\b{e}_i\}_{i=1}^n$ and $\{\b{e}^i\}_{i=1}^n$ are basis vectors and dual basis vectors, respectively. 
The $\b{T} = T^{i_1 \dots i_r}_{j_1 \dots j_s}$ are the coordinates of this tensor, where $i_1 \dots i_r, j_1 \dots j_s \in \{1, \dots, n\}$.
\end{remark}

\begin{definition}[Covariant components]
Consider a vector $\b{V} \in T_{\b{p}}\mathcal{M}$ in an $n$-dimensional manifold $\mathcal{M}$. 
The inner product of the vector $\b{V}$ and the $i$-th basis vector of the coordinate system $\b{e}_i$ is defined as the $i$-th \textbf{covariant component}, denoted by $V_i$:
\begin{align}\label{equation_covariant_components}
\boxed{
V_i = \langle \b{V}, \b{e}_i \rangle,
}
\end{align}
where $\langle \cdot, \cdot \rangle$ denotes the inner product. 
The $\{V_1, \dots, V_n\}$ are called the covariant components. The covariant components are denoted by lower (subscript) indices, as also mentioned in Definition \ref{definition_covariant_and_contravariant_indices}. 
\end{definition}
\begin{proof}
The following is the proof for Eq. (\ref{equation_covariant_components}) in a coordinate system with orthonormal basis vectors. 
For general curvilinear coordinate system, the Eq. (\ref{equation_covariant_components}) is the definition of covariant components and the definition does not require a proof. 

Let $\{\b{e}_1, \b{e}_2, \dots, \b{e}_n\}$ denote the basis vectors of the $n$-dimensional (scaled) Cartesian coordinate system. In such coordinate system, the basis vectors are orthonormal, meaning that:
\begin{align}\label{equation_basis_orthonormal}
\langle \b{e}_i, \b{e}_j \rangle = 
\left\{
    \begin{array}{ll}
        1 & \mbox{if } i = j, \\
        0 & \mbox{if } i \neq j.
    \end{array}
\right.
\end{align}
In other words, in such coordinate system, the basis vectors are perpendicular (orthogonal) while each basis vector has a unit length. 

We start from the right hand side and obtain the left hand side of that equation: 
\begin{align*}
&\langle \b{V}, \b{e}_i \rangle \overset{(\ref{equation_contravariant_components})}{=} \langle (V^1 \b{e}_1 + V^2 \b{e}_2 + \dots + V^n \b{e}_n), \b{e}_i \rangle \\
&= (V^1 \underbrace{\langle \b{e}_1, \b{e}_i \rangle}_{0}) + \dots \\
&~~+ (V^{i-1} \underbrace{\langle \b{e}_{i-1}, \b{e}_{i} \rangle}_{0}) + (V^i \underbrace{\langle \b{e}_i, \b{e}_i}_{1}\rangle) \\
&~~+ (V^{i+1} \underbrace{\langle\b{e}_{i+1}, \b{e}_{i+1}\rangle}_{0}) + \dots + V^n \underbrace{\langle\b{e}_n, \b{e}_i\rangle}_{0} \overset{(\ref{equation_basis_orthonormal})}{=} V^i.
\end{align*}
\end{proof}

\begin{definition}[Coordinate basis for the cotangent space]
The basis vectors in Eq. (\ref{equation_bases_for_covectors}) exist in the cotangent space at every point $\b{p} \in \mathcal{M}$. Thus, $\{dx^1, \dots, dx^n\}$ are also the \textbf{basis vectors for the cotangent space} (also called the \textbf{dual basis}). As a result, a covector $\omega \in T_{\b{p}}^*\mathcal{M}$ can be stated as:
\begin{align}\label{equation_covector_Wi_dXi}
\boxed{
\omega = \omega_i\, dx^i.
}
\end{align}
\end{definition}

\begin{lemma}[Extracting the $i$-th contravariant component from vector]
Consider a vector $\b{V}$ with contravariant components $\{V^1, \dots, V^n\}$. 
We can extract the $i$-th contravariant component from vector as:
\begin{align}\label{equation_dX_V_Vi}
\boxed{
dx^i(\b{V}) = V^i.
}
\end{align}
\end{lemma}
\begin{proof}
According to Eq. (\ref{equation_contravariant_components_2}), we have:
\begin{align*}
\b{V} = V^1 \partial_1 + V^2 \partial_2 + \dots + V^n \partial_n.
\end{align*}
We have:
\begin{align*}
&dx^i(\b{V}) \overset{(\ref{equation_contravariant_components_2})}{=} dx^i(V^1 \partial_1 + V^2 \partial_2 + \dots + V^n \partial_n) \\
&\overset{(a)}{=} V^1 dx^i(\partial_1) + V^2 dx^i(\partial_2) + \dots + V^n dx^i(\partial_n) \\
&\overset{(\ref{equation_dx_partial_partialx_delta})}{=} V^1 \delta_1^i + V^2 \delta_2^i + \dots + V^n \delta_n^i \overset{(\ref{equation_Kronecker_delta_lower_upper})}{=} V^i,
\end{align*}
where $(a)$ is because of linearity of $dx^i$.
\end{proof}

\begin{lemma}[Extracting the $i$-th covariant component from a covector]
Consider a covector (1-form) $\omega$ with covariant components $\{\omega_1, \dots, \omega_n\}$.
We can extract the $i$-th covariant component from $\omega$ as:
\begin{align}
\boxed{
\omega(\partial_i) = \omega_i.
}
\end{align}
\end{lemma}
\begin{proof}
By definition of the covector components in a coordinate basis, we have:
\begin{align*}
\omega \overset{(\ref{equation_covector_Wi_dXi})}{=} \omega_1 dx^1 + \omega_2 dx^2 + \dots + \omega_n dx^n.
\end{align*}
Then, we evaluate $\omega$ on the basis vector $\partial_i$:
\begin{align*}
\omega(\partial_i) 
&= (\omega_1 dx^1 + \dots + \omega_n dx^n)(\partial_i) \\
&\overset{(a)}{=} \omega_1 dx^1(\partial_i) + \dots + \omega_n dx^n(\partial_i) \\
&\overset{(\ref{equation_dx_partial_partialx_delta})}{=} \omega_1 \delta^1_i + \dots + \omega_n \delta^n_i \overset{(\ref{equation_Kronecker_delta_lower_upper})}{=} \omega_i,
\end{align*}
where $(a)$ is because of linearity of $dx^i$.
\end{proof}

\begin{figure}[!h]
\centering
\includegraphics[width=3.2in]{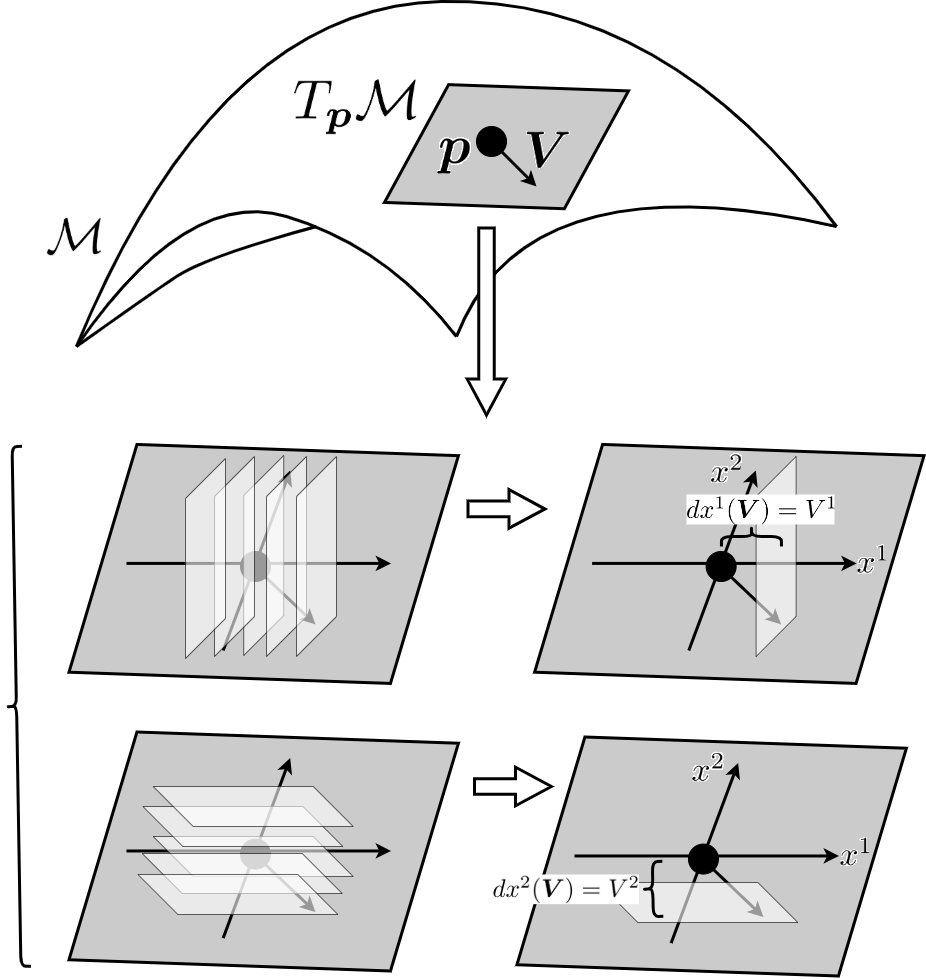}
\caption{Cotangent space as parallel planes along $x^1$ and $x^2$ where these planes are orthogonal to the tangent space. Here, the manifold is locally two-dimensional so its covectors are $w_1$ and $w_2$ where, according to Eq. (\ref{equation_dX_V_Vi}), we have $dx^1(\b{V}) = V^1$ and $dx^2(\b{V}) = V^2$.}
\label{figure_covector}
\end{figure}

\begin{remark}[Geometrical interpretation of vector and covector]
Intuitively, a vector is a direction, with some length, such as an arrow. According to Definition \ref{definition_vector}, a vector is in the tangent space. An example vector, in the tangent space of point $\b{p} \in \mathcal{M}$, is shown in Fig. \ref{figure_tangent_vector} (see Section \ref{section_tangent_space_tangent_vector}).

A covector is something that measures vectors. For example, a covector $\omega_i$ can give the $i$-th component of vector, i.e., $V^i$. This is because, on the one hand, $dx^i, \forall i \in \{1, \dots, n\}$ are basis vectors of the covectors, according to Eqs. (\ref{equation_bases_for_covectors}) and (\ref{equation_covector_Wi_dXi}). On the other hand, according to Eq. (\ref{equation_dX_V_Vi}), $dx^i$---which is the basis vector for covector---can extract the $i$-th component of vector, i.e., $V^i$. 

Thus, intuitively, consider parallel planes along $x^i$ where these planes are orthogonal to the tangent space; this explains why these planes are named ``cotangent space" containing the covectors or cotangent vectors. As illustrated in Fig. \ref{figure_covector}, a covector $\omega_i$ is like scalar values as intersection of these parallel planes along $x^i$. 
\end{remark}

\subsubsection{Geometric Interpretation of Contravariant and Covariant Components}

\begin{figure}[!h]
\centering
\includegraphics[width=3.2in]{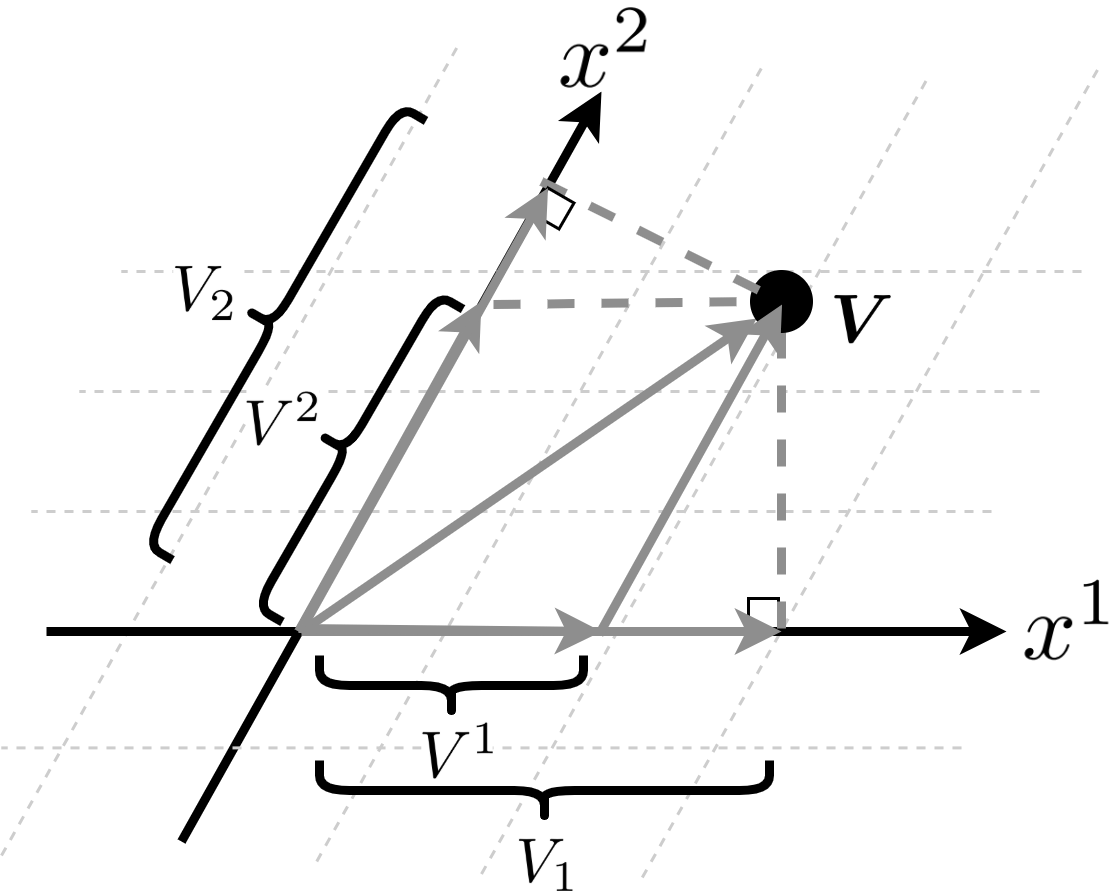}
\caption{Geometric interpretation of contravariant components $\{V^1, V^2\}$ and covariant components $\{V_1, V_2\}$ of a vector $\b{V}$ in a two-dimensional coordinate system $(x^1, x^2)$.}
\label{figure_contravariant_covariant_geometric}
\end{figure}

Consider a vector $\b{V}$ in an $n$-dimensional coordinate system. Although the coordinate system can be any coordinate system---including Cartesian, affine, and curvilinear---consider a two-dimensional affine coordinate system for simplicity, as illustrated in Fig. \ref{figure_contravariant_covariant_geometric}.

\begin{remark}[Geometric interpretation of contravariant components]
Consider the point $\b{V}$ in the coordinate system. Assume we draw lines, from the point, \underline{parallel} to the axes, as shown in Fig. \ref{figure_contravariant_covariant_geometric}. The length of line fragment from origin to the intersection of the parallel line and the $i$-th axis is the $i$-th contravariant component $V^i$. 
The $\{V^1, \dots, V^n\}$ are the contravariant components. 
\end{remark}

\begin{remark}[Geometric interpretation of covariant components]
Consider the point $\b{V}$ in the coordinate system. Assume we draw lines, from the point, \underline{perpendicular} to the axes, as illustrated in Fig. \ref{figure_contravariant_covariant_geometric}. The length of line fragment from origin to the intersection of the perpendicular line and the $i$-th axis is the $i$-th contravariant component $V_i$. 
The $\{V_1, \dots, V_n\}$ are the covariant components. 
\end{remark}

\begin{lemma}
In Cartesian coordinate system, the contravariant and covariant components are equivalent. 
\end{lemma}
\begin{proof}
In Cartesian coordinate system---scaled or non-scaled---the coordinate axes are orthogonal. Therefore, as shown in Fig. \ref{figure_contravariant_covariant_geometric} (assume that the axes $x^1$ and $x^2$ are perpendicular in this figure), the lines drawn as parallel or perpendicular are the same. So, the contravariant and covariant components are the same.
\end{proof}

\subsection{Transformation of Coordinate Systems}

\subsubsection{Transformation of Coordinates for Contravariant and Covariant Components}

\begin{proposition}[Transformation of coordinates for contravariant components]
Consider a coordinate system $x$ on a manifold, where the vector $\b{V} = [V^1, \dots, V^n]^\top \in T_{\b{p}}\mathcal{M}$ is represented in this coordinate system. If we use another coordinate system $y$ on the same manifold, the same vector in the new coordinate system is obtained as:
\begin{align}\label{equation_transformation_coordinates_contravariant}
\boxed{
\widetilde{V}^i = \frac{\partial y^i}{\partial x^j} V^j,
}
\end{align}
where $[\widetilde{V}^1, \dots, \widetilde{V}^n]^\top \in T_{\b{p}}\mathcal{M}$ is the vector in the new coordinate system and $j$ is the summation index in Einstein convention\footnote{As $\partial x^j$ is in the denominator, its index $j$ can be considered as a lower (subscript) index.}. 
\end{proposition}
\begin{proof}
Let $x=(x^1,\dots,x^n)$ and $y=(y^1,\dots,y^n)$ be two coordinate systems on a manifold.
In coordinate systems $x$ and $y$, the coordinate basis vectors are as in Eqs. (\ref{equation_ej_partial}) and (\ref{equation_ej_tilde_partial}), respectively.
In the $x$-coordinate basis, the vector $\b{V}$ can be represented as:
\begin{align}\label{equation_V_Vj_partialX}
\b{V} \overset{(\ref{equation_contravariant_components})}{=} V^j \b{e}_j \overset{(\ref{equation_ej_partial})}{=} V^j \frac{\partial}{\partial x^j}.
\end{align}
In the $y$-coordinate basis, the same vector is written as:
\begin{align}\label{equation_V_U_partialY}
\b{V} \overset{(\ref{equation_contravariant_components})}{=} \widetilde{V}^i\, \widetilde{\b{e}}_i \overset{(\ref{equation_ej_tilde_partial})}{=} \widetilde{V}^i \frac{\partial}{\partial y^i}.
\end{align}
Using the chain rule in derivatives, we have:
\begin{align}\label{equation_chain_rule_derivative}
\frac{\partial}{\partial x^j}
=
\frac{\partial y^i}{\partial x^j}
\frac{\partial}{\partial y^i}.
\end{align}
We have:
\begin{align}
\b{V}
&\overset{(\ref{equation_V_Vj_partialX})}{=}
V^j \frac{\partial}{\partial x^j}
\overset{(\ref{equation_chain_rule_derivative})}{=}
V^j
\left(
\frac{\partial y^i}{\partial x^j}
\frac{\partial}{\partial y^i}
\right) \nonumber
\\
&\overset{(a)}{=}
\left(
\frac{\partial y^i}{\partial x^j} V^j
\right)
\frac{\partial}{\partial y^i}. \label{equation_V_partialY_partialX_V_partialY}
\end{align}
where $(a)$ is because of rearranging. 

Comparing Eqs. (\ref{equation_V_U_partialY}) and (\ref{equation_V_partialY_partialX_V_partialY}) gives Eq. (\ref{equation_transformation_coordinates_contravariant}).
\end{proof}

\begin{proposition}[Transformation of coordinates for covariant components]
Consider a coordinate system $x$ on a manifold, where the covariant components of a vector $\b{V}$ are $\{V_1, \dots, V_n\}$ in this coordinate system. If we use another coordinate system $Y$ on the same manifold, the same covariant components in the new coordinate system are obtained as:
\begin{align}\label{equation_transformation_coordinates_covariant}
\boxed{
\widetilde{V}_i = \frac{\partial x^j}{\partial y^i} V_j,
}
\end{align}
where $\{\widetilde{V}_1, \dots, \widetilde{V}_n\}$ are the covariant components in the new coordinate system and $j$ is the summation index in Einstein convention\footnote{As $\partial x^j$ is in the numerator, its index $j$ can be considered as an upper (superscript) index.}. 
\end{proposition}
\begin{proof}
Let $\{\b{e}_i\}_{i=1}^n$ be the basis vectors associated with the coordinate system $X$, and $\{\widetilde{\b{e}}_i\}_{i=1}^n$ be the basis vectors associated with the coordinate system $y$.
The bases are related through the coordinate transformation, stated in Eq. (\ref{equation_coordinate_basis_transformation}). 

Let $\b{W}$ be an arbitrary vector on the manifold. According to Eq. (\ref{equation_transformation_coordinates_contravariant}), its contravariant components in the
two coordinate systems satisfy:
\begin{align}\label{equation_W_tilde_partialX_partialY_W}
W^j
\overset{(\ref{equation_transformation_coordinates_contravariant})}{=}
\frac{\partial x^j}{\partial y^i} \widetilde{W}^i,
\end{align}
where $W^j$ and $\widetilde{W}^i$ are contravariant components of $\b{W}$ in the $x$ and $y$ coordinates, respectively. 

The inner product of components of $\b{V}$ and $\b{W}$ in $x$-coordinates and the inner product of components of $\widetilde{\b{V}}$ and $\widetilde{\b{W}}$ in $y$-coordinates are:
\begin{align}
&\langle\b{V},\b{W}\rangle = V_j W^j, \label{equation_V_W_inner} \\
&\langle\b{V},\b{W}\rangle = \widetilde{V}_i \widetilde{W}^i, \label{equation_V_W_inner2}
\end{align}
respectively. 

We have:
\begin{align}
\langle\b{V},\b{W}\rangle &\overset{(\ref{equation_V_W_inner})}{=} V_j W^j \overset{(\ref{equation_W_tilde_partialX_partialY_W})}{=}
V_j
\left(
\frac{\partial x^j}{\partial y^i} \widetilde{W}^i
\right) \nonumber
\\
&\overset{(a)}{=} 
\left(\frac{\partial x^j}{\partial y^i} V_j\right) \widetilde{W}^i, \label{equation_V_dot_W_partialX_partialY_V_Wtilde}
\end{align}
where $(a)$ is because of rearranging the parentheses and rearranging the elements in multiplication. 

Comparing Eqs. (\ref{equation_V_W_inner2}) and (\ref{equation_V_dot_W_partialX_partialY_V_Wtilde}) gives Eq. (\ref{equation_transformation_coordinates_covariant}).
\end{proof}

\begin{remark}[The naming reason of contravariant and covariant components]
The following justifies the reason for the names of contravariant and covariant components:
\begin{itemize}
\item Equation (\ref{equation_transformation_coordinates_contravariant}) shows that for converting contravariant components in coordinates $x$ to contravariant components in coordinates $y$, we should use $\partial y / \partial x$ which is the inverse Jacobian. This is why they are called contravariant components, where ``contra" refers to ``inverse".
\item Equation (\ref{equation_transformation_coordinates_covariant}) shows that for converting covariant components in coordinates $x$ to covariant components in coordinates $y$, we should use $\partial x / \partial y$ which is the Jacobian. This is why they are called covariant components, where ``co" refers to ``same".
\end{itemize}
\end{remark}
\begin{proof}
Another way to prove that the transformation of coordinates in contravariant and covariant components use inverse Jacobian and Jacobian, respectively, is provided below:

\textbf{Proof for why vector components (contravariants) transform oppositely:}
According to Eq. (\ref{equation_contravariant_components_2}), a vector is:
\begin{align*}
\b{V} = V^i \frac{\partial}{\partial x^i}.
\end{align*}
But it must represent the same geometric vector, even in the new coordinates $y$. So, we also write:
\begin{align*}
\b{V} = \widetilde{V}^i \frac{\partial}{\partial y^i}.
\end{align*}
Substitute the basis transformation using the chain rule:
\begin{align*}
\frac{\partial }{\partial y^i} = \frac{\partial x^j}{\partial y^i} \frac{\partial }{\partial x^j}.
\end{align*}
Then:
\begin{align*}
\b{V} = \widetilde{V}^i \frac{\partial x^j}{\partial y^i} \frac{\partial }{\partial x^j}.
\end{align*}
Comparing with:
\begin{align*}
\b{V} = V^j \frac{\partial}{\partial x^j},
\end{align*}
gives:
\begin{align*}
V^j = \frac{\partial x^j}{\partial y^i} \widetilde{V}^i.
\end{align*}
So, vector components use the inverse Jacobian for transformation.
This is why they are called contravariant.

\textbf{Proof for why covector components (covariant) transform the same way:}
According to Eq. (\ref{equation_covector_Wi_dXi}), a covector is:
\begin{align*}
W = W_i\, dx^i.
\end{align*}
But it must represent the same geometric vector, even in the new coordinates $y$. So, we also write:
\begin{align*}
W = \widetilde{W}_i\, dy^i.
\end{align*}
According to chain rule in differentiation, we have:
\begin{align*}
dy^i = \frac{\partial y^i}{\partial x^j} dx^j.
\end{align*}
By substituting $dy^i$, we obtain:
\begin{align*}
W = \widetilde{W}_i\, \frac{\partial y^i}{\partial x^j} dx^j.
\end{align*}
Comparing with:
\begin{align*}
W = W_j\, dx^j,
\end{align*}
gives:
\begin{align*}
W^j = \widetilde{W}_i\, \frac{\partial y^i}{\partial x^j}.
\end{align*}
So, covector components use the Jacobian for transformation.
This is why they are called covariant.
\end{proof}

\subsubsection{Tensor Transformation Laws}

\begin{corollary}[Tensor transformation laws]\label{corollary_tensor_transformation_laws}
The Eqs. (\ref{equation_transformation_coordinates_contravariant}) and (\ref{equation_transformation_coordinates_covariant}) can be generalized to tensors and any number of lower and upper indices. 
Let $\b{T}$ denote a tensor represented in $x$ coordinate system and $\widetilde{\b{T}}$ be the same tensor in $y$ coordinate system. Suppose $T^i$ and $T_j$ denote the contravariant and covariant components of tensor $\b{T}$, respectively, and suppose $\widetilde{T}^i$ and $\widetilde{T}_j$ denote the contravariant and covariant components of tensor $\widetilde{\b{T}}$, respectively

For the contravariant components, we have:
\begin{align}
&\boxed{
\widetilde{T}^i = \frac{\partial y^i}{\partial x^j} T^j,
} \\
&\boxed{
\widetilde{T}^{ik} = \frac{\partial y^i}{\partial x^j} \frac{\partial y^k}{\partial x^\ell} T^{j\ell},
} \\
&\boxed{
\widetilde{T}^{ikp} = \frac{\partial y^i}{\partial x^j} \frac{\partial y^k}{\partial x^\ell} \frac{\partial y^p}{\partial x^q} T^{j\ell q},
}
\end{align}
and so on. 

For the covariant components, we have:
\begin{align}
&\boxed{
\widetilde{T}_i = \frac{\partial x^j}{\partial y^i} T_j,
} \\
&\boxed{
\widetilde{T}_{ik} = \frac{\partial x^j}{\partial y^i} \frac{\partial x^\ell}{\partial y^k} T_{j\ell},
} \label{equation_transformation_coordinates_covariant_Tik} \\
&\boxed{
\widetilde{T}_{ikp} = \frac{\partial x^j}{\partial y^i} \frac{\partial x^\ell}{\partial y^k} \frac{\partial x^q}{\partial y^p} T_{j\ell q},
} 
\end{align}
and so on.

For a combination of the contravariant and covariant components, we have:
\begin{align}
&\boxed{
\widetilde{T}_{i}^k = \frac{\partial y^k}{\partial x^\ell} \frac{\partial x^j}{\partial y^i} T_{j}^\ell,
} \\
&\boxed{
\widetilde{T}_{i}^{k p} = \frac{\partial y^k}{\partial x^\ell} \frac{\partial x^j}{\partial y^i} \frac{\partial y^p}{\partial x^q} T_{j}^{\ell q},
} \\
&\boxed{
\widetilde{T}_{i p}^{k} = \frac{\partial y^k}{\partial x^\ell} \frac{\partial x^j}{\partial y^i} \frac{\partial x^q}{\partial y^p} T_{j q}^{\ell},
} 
\end{align}
and so on.
\end{corollary}

Remembering the equations of tensor transformation laws is easy considering the contraction of indices in the Einstein convention. The lower indices of tensors correspond to indices in the denominator of derivatives, and the upper indices of tensors correspond to indices in the numerator of derivatives. 

Also note that, in the tensor transformation laws, the number of lower and upper indices of the tensor should be the same in the two coordinate systems because when a tensor is represented in different coordinate systems, its type (see Definition \ref{definition_tensor}) does not change.

\subsection{Directional Derivative}\label{section_directional_derivative}

Recall Section \ref{section_vector_and_vector_field} about vector field, which is required for understanding the following concepts. 

\begin{definition}[Directional derivative of a smooth function at a point]
Let $\mathcal{M}$ be a smooth manifold, let
$f \in C^\infty(\mathcal{M})$, let $\b{p} \in \mathcal{M}$, and let
$\b{\xi} \in T_{\b{p}}\mathcal{M}$.
According to Eq. (\ref{equation_contravariant_components_2}), in a local coordinate chart $(x^1,\dots,x^n)$ around $\b{p}$, we have:
\[
\b{\xi} = \xi^i \partial_i|_{\b{p}},
\]
where $|_{\b{p}}$ means ``at point $\b{p}$".

The \textbf{directional derivative} of $f$ at $\b{p}$ in the direction
$\b{\xi}$ is defined as:
\begin{equation}\label{equation_directional_derivative_tangent_vector}
\boxed{
\b{\xi}(f)
:=
\xi^i \frac{\partial f}{\partial x^i}(\b{p})
\overset{(\ref{equation_partial_i})}{=}
\xi^i \partial_i f(\b{p}).
}
\end{equation}
The quantity $\b{\xi}(f) \in \mathbb{R}$ measures the rate of change of function $f$ along the direction of tangent vector $\b{\xi}$ at point $\b{p}$.
\end{definition}

\begin{definition}[Differential of a smooth function at a point]
Let $\mathcal{M}$ be a smooth manifold, let
$f \in C^\infty(\mathcal{M})$, and let $\b{p} \in \mathcal{M}$.
The \textbf{differential} of $f$ at $\b{p}$, denoted by
$Df(\b{p})$, is the linear map:
\begin{equation}\label{equation_differential_map}
Df(\b{p}) : T_{\b{p}}\mathcal{M} \to \mathbb{R},
\end{equation}
defined by:
\begin{equation}\label{equation_differential_acts_on_vector}
\boxed{
\begin{aligned}
Df(\b{p})[\b{\xi}]
&:=
\b{\xi}(f) \\
&\overset{(\ref{equation_directional_derivative_tangent_vector})}{=} \xi^i \partial_i f(\b{p}),
\end{aligned}
}
\end{equation}
for $\b{\xi} \in T_{\b{p}}\mathcal{M}$.
The quantity $Df(\b{p})[\b{\xi}] \in \mathbb{R}$ measures the rate of change of function $f$ along the direction of tangent vector $\b{\xi}$ at point $\b{p}$.
\end{definition}

\begin{definition}[Directional derivative of a smooth function along a vector field]\label{definition_directional_derivative_function_along_vector_field}
Let $\mathcal{M}$ be a smooth manifold, let
$f \in C^\infty(\mathcal{M})$, and let
$\b{X} \in \mathfrak{X}(\mathcal{M})$ be a smooth vector field.
The \textbf{directional derivative of the function $f$ along the vector field $\b{X}$} at point $\b{p} \in \mathcal{M}$ is defined by:
\begin{equation}\label{equation_directional_derivative_function_along_vector_field}
\boxed{
\b{X}(f)(\b{p})
:=
Df(\b{p})[\b{X}(\b{p})] \in \mathbb{R}.
}
\end{equation}
Therefore, $\b{X}(f)$ is a scalar function on $\mathcal{M}$.
The quantity $\b{X}(f)(\b{p}) \in \mathbb{R}$ measures the rate of
change of function $f$ along the direction of tangent vector
$\b{X}(\b{p})$ at point $\b{p}$.

In local coordinates, according to
Eq. (\ref{equation_contravariant_components_2}), the vector field
$\b{X}$ can be written as:
\begin{align}\label{equation_X_Xi_partiali_in_directional_derivative}
\b{X}
=
X^i \frac{\partial}{\partial x^i}
\overset{(\ref{equation_partial_i})}{=}
X^i \partial_i,
\end{align}
where $X^i$ is the $i$-th contravariant component of $\b{X}$ in the
coordinate basis. Hence:
\begin{equation}\label{equation_X_f_X_partial_f_partial_x}
\boxed{
\b{X}(f)(\b{p})
=
X^i \frac{\partial f}{\partial x^i}\big|_{\b{p}}
=
X^i \partial_i f\big|_{\b{p}}.
}
\end{equation}

Thus, a vector field acts as a first-order differential operator on
smooth scalar functions: at each point $\b{p} \in \mathcal{M}$, the
vector $\b{X}(\b{p}) \in T_{\b{p}}\mathcal{M}$ specifies a direction,
and $\b{X}(f)(\b{p})$ measures the rate of change of $f$ at $\b{p}$
along that direction. Equation
(\ref{equation_X_f_X_partial_f_partial_x}) says that, in local
coordinates, $\b{X}(f)$ is obtained by taking the partial derivatives of $f$ along the coordinate basis vectors
$\partial_i=\frac{\partial}{\partial x^i}$ at point $\b{p}$, weighting them by the components $X^i$, and summing over $i$.
\end{definition}

\begin{remark}[Repeated action of vector fields on a smooth function]\label{remark_XYf}
Let $\b{X}, \b{Y} \in \mathfrak{X}(\mathcal{M})$ be smooth vector
fields, and let $f \in C^\infty(\mathcal{M})$. Since $\b{Y}(f)$ is a
scalar function on $\mathcal{M}$, the vector field $\b{X}$ can act on
it again. The map:
\begin{align}
\b{X}(\b{Y}(f)): \mathcal{M} \to \mathbb{R},
\end{align}
means the $\b{Y}(f): \mathcal{M} \to \mathbb{R}$ is scalar function on the manifold, and the vector field $\b{X}(\cdot)$ acts on this function. 
The expression $\b{X}(\b{Y}(f))$ means the directional derivative of the scalar function $\b{Y}(f)$
along the vector field $\b{X}$.

At a point $\b{p} \in \mathcal{M}$, this is:
\begin{equation}\label{equation_X_of_Y_of_f}
\boxed{
\b{X}(\b{Y}(f))(\b{p})
:=
D(\b{Y}(f))(\b{p})[\b{X}(\b{p})] \in \mathbb{R}.
}
\end{equation}
Hence, $\b{X}(\b{Y}(f))$ is again a scalar function on
$\mathcal{M}$.

In local coordinates, if:
\[
\b{X} = X^i \partial_i,
\qquad
\b{Y} = Y^j \partial_j,
\]
then:
\begin{equation}\label{equation_Yf_coordinate}
\boxed{
\b{Y}(f)=Y^j\partial_j f,
}
\end{equation}
and therefore:
\begin{equation}\label{equation_XYf_coordinate}
\boxed{
\b{X}(\b{Y}(f))
=
X^i \partial_i \!\left( Y^j \partial_j f \right).
}
\end{equation}
By the product rule, this can be expanded as:
\begin{equation}\label{equation_XYf_coordinate_expanded}
\boxed{
\b{X}(\b{Y}(f))
=
X^i \frac{\partial Y^j}{\partial x^i}\frac{\partial f}{\partial x^j}
+
X^i Y^j \frac{\partial^2 f}{\partial x^i \partial x^j}.
}
\end{equation}
\end{remark}

\begin{definition}[Coordinate directional derivative of a vector field along another vector field]\label{definition_directional_derivative_vector_field_along_vector_field}
Let $\b{X},\b{Y} \in \mathfrak{X}(\mathcal{M})$ be smooth vector
fields. In a local coordinate chart, according to Eq. (\ref{equation_contravariant_components_2}), we can write:
\[
\b{Y}
=
Y^k \frac{\partial}{\partial x^k}
\overset{(\ref{equation_partial_i})}{=}
Y^k \partial_k.
\]
The \textbf{coordinate directional derivative of the vector field $\b{Y}$ along (in the direction of) the vector field $\b{X}$} at
$p \in \mathcal{M}$ is a map:
\begin{align}
D\b{Y}(\b{p})[\b{X}(\b{p})] : T_{\b{p}}\mathcal{M} \to T_{\b{p}}\mathcal{M},
\end{align}
defined by:
\begin{align}\label{equation_DYpXp_XYk_p_partialkp}
\boxed{
D\b{Y}(\b{p})[\b{X}(\b{p})]
:=
\b{X}(Y^k)(\b{p})\,\partial_k\big|_{\b{p}} \in T_{\b{p}}\mathcal{M}.
}
\end{align}
Equivalently, in local coordinates, if $\b{X} = X^i\partial_i$ according to Eq. (\ref{equation_contravariant_components_2}), then:
\begin{align}
\boxed{
D\b{Y}(\b{p})[\b{X}(\b{p})]
=
\left(\!
X^i \frac{\partial Y^k}{\partial x^i}
\right)\!(\b{p})\,\partial_k\big|_{\b{p}}.
}
\end{align}
This derivative differentiates the coordinate component
functions of $\b{Y}$ along $\b{X}$.
In other words, the quantity
$D\b{Y}(\b{p})[\b{X}(\b{p})]$ represents the rate of change of the
coordinate components of $\b{Y}$ along $\b{X}$ at point $\b{p}$.
\end{definition}

Note that the notation we use for directional derivative of a vector field along another vector field follows this notational rule:
\begin{align}
\boxed{
D(\text{field being differentiated})(\text{base point})[\text{direction}].
}
\end{align}




\begin{remark}[Type of directional derivatives of a function or a vector field]
The differential of a smooth scalar function $f$ at point
$\b{p} \in \mathcal{M}$ is the linear map:
\[
Df(\b{p}) : T_{\b{p}}\mathcal{M} \to \mathbb{R}.
\]
Hence, for a vector field $\b{X} \in \mathfrak{X}(\mathcal{M})$, we have:
\[
\b{X}(f)(\b{p}) = Df(\b{p})[\b{X}(\b{p})] \in \mathbb{R}.
\]
Therefore, $\b{X}(f)$ is a scalar function on $\mathcal{M}$.

In contrast, the coordinate directional derivative of a vector field
$\b{Y} \in \mathfrak{X}(\mathcal{M})$ at $\b{p}$ is, in a chosen local
coordinate chart, the map:
\[
D\b{Y}(\b{p}) : T_{\b{p}}\mathcal{M} \to T_{\b{p}}\mathcal{M}.
\]
Therefore:
\[
D\b{Y}(\b{p})[\b{X}(\b{p})] \in T_{\b{p}}\mathcal{M}.
\]

In other words, the directional derivative of a scalar function along a
vector field is a scalar, while the coordinate directional derivative of
a vector field along another vector field is a tangent vector. The
former acts on a scalar-valued function $f:\mathcal{M}\to\mathbb{R}$,
whereas the latter differentiates the coordinate component functions of
the vector field $\b{Y}$ along the direction $\b{X}$.
\end{remark}

\begin{remark}[Difference between dimensionalities of directional derivatives]
Note that $\b{X}(\b{Y}(f))$, discussed in
Remark \ref{remark_XYf}, is a scalar function on $\mathcal{M}$,
whereas $D\b{Y}(\b{p})[\b{X}(\b{p})]$, discussed in
Definition \ref{definition_directional_derivative_vector_field_along_vector_field},
is the coordinate directional derivative of the vector field $\b{Y}$
along $\b{X}$ at $\b{p}$, which is a tangent vector at $\b{p}$ in the
chosen coordinate chart. In summary:
\begin{align*}
& \b{X}(\b{Y}(f)) : \mathcal{M} \to \mathbb{R}, \\
& \b{X}(\b{Y}(f))(\b{p}) \in \mathbb{R}, \\
& D\b{Y}(\b{p})[\b{X}(\b{p})] \in T_{\b{p}}\mathcal{M}.
\end{align*}
\end{remark}

\begin{remark}[Directional derivative in local coordinates versus intrinsic derivative]
The quantity $D\b{Y}(\b{p})[\b{X}(\b{p})]$ is not an intrinsic derivative of vector fields on $\mathcal{M}$; it is the ordinary derivative of the coordinate representation of $\b{Y}$ in a chosen chart.
That is, it corresponds to the standard Euclidean derivative of the component functions of $\b{Y}$ 
when the manifold is locally identified with $\mathbb{R}^n$.
Therefore, the directional derivative depends on the choice of coordinates and is not an intrinsic geometric object 
on the manifold.

To define a coordinate-invariant (intrinsic) derivative of vector fields on a manifold, one needs a connection. 
This leads to the notion of the covariant derivative, which will be introduced in Section \ref{section_connection_covariant_derivative}.
\end{remark}

Here, we introduced general directional derivative. 
The ambient directional derivative for matrix manifolds will be discussed in Section \ref{section_ambient_directional_derivative_matrix_manifolds}.

\subsection{Lie Bracket}

\begin{definition}[Lie bracket]\label{definition_Lie_bracket}
Let $\b{X} = X^i \frac{\partial}{\partial x^i}$ and 
$\b{Y} = Y^i \frac{\partial}{\partial x^i}$ be vector fields on a manifold $\mathcal{M}$.
The \textbf{Lie bracket} of $\b{X}$ and $\b{Y}$, denoted by $[\b{X},\b{Y}]$, is defined as:
\begin{align}\label{equation_Lie_bracket}
\boxed{
[\b{X},\b{Y}](f) := \b{X}(\b{Y}(f)) - \b{Y}(\b{X}(f)),
}
\end{align}
for all smooth functions $f: \mathcal{M} \rightarrow \mathbb{R}$ on the manifold.

In other words, the Lie bracket measures how much the flows of $\b{X}$ and $\b{Y}$ fail to commute.
\end{definition}

\begin{proposition}[Coordinate expression of Lie bracket]\label{proposition_Lie_bracket_coordinates}
Let $\b{X} = X^i \frac{\partial}{\partial x^i}$ and 
$\b{Y} = Y^i \frac{\partial}{\partial x^i}$ be vector fields on a manifold.
The Lie bracket is given by:
\begin{align}\label{equation_Lie_bracket_coordinate_expression}
\boxed{
[\b{X},\b{Y}] = \left( X^i \partial_i Y^k - Y^i \partial_i X^k \right) \partial_k.
}
\end{align}
\end{proposition}

\begin{proof}
By definition of the Lie bracket, for any smooth function $f \in C^\infty(\mathcal{M})$, we have:
\begin{align*}
[\b{X},\b{Y}](f) = \b{X}(\b{Y}(f)) - \b{Y}(\b{X}(f)).
\end{align*}
First compute $\b{Y}(f)$:
\begin{align}\label{equation_Y_f_Y_partialf_partialx}
\b{Y}(f) \overset{(\ref{equation_X_f_X_partial_f_partial_x})}{=} Y^j \frac{\partial f}{\partial x^j}.
\end{align}

Applying $\b{X}$ to $\b{Y}(f)$ gives:
\begin{align}
\b{X}(\b{Y}(f)) 
&\overset{(a)}{=} X^i \frac{\partial}{\partial x^i} \left( Y^j \frac{\partial f}{\partial x^j} \right) \nonumber \\
&\overset{(b)}{=} X^i \left( \frac{\partial Y^j}{\partial x^i} \frac{\partial f}{\partial x^j}
+ Y^j \frac{\partial^2 f}{\partial x^i \partial x^j} \right), \nonumber
\end{align}
where $(a)$ is because of Eqs. (\ref{equation_X_f_X_partial_f_partial_x}) and (\ref{equation_Y_f_Y_partialf_partialx}), and $(b)$ is because of product rule of derivatives. 

Similarly:
\begin{align*}
\b{Y}(\b{X}(f)) 
&= Y^j \frac{\partial}{\partial x^j} \left( X^i \frac{\partial f}{\partial x^i} \right) \\
&= Y^j \left( \frac{\partial X^i}{\partial x^j} \frac{\partial f}{\partial x^i}
+ X^i \frac{\partial^2 f}{\partial x^j \partial x^i} \right).
\end{align*}
Subtracting, we obtain:
\begin{align*}
[\b{X},\b{Y}](f)
&= X^i \frac{\partial Y^j}{\partial x^i} \frac{\partial f}{\partial x^j}
- Y^j \frac{\partial X^i}{\partial x^j} \frac{\partial f}{\partial x^i}.
\end{align*}

In the first term, we rename the dummy variables $j$ to $k$, and in the second term, we rename the dummy variables $j$ to $i$ and $i$ to $k$:
\begin{align*}
[\b{X},\b{Y}](f)
&= X^i \frac{\partial Y^k}{\partial x^i} \frac{\partial f}{\partial x^k} 
- Y^i \frac{\partial X^k}{\partial x^i}
\frac{\partial f}{\partial x^k} \\
&\overset{(a)}{=} \left( X^i \frac{\partial Y^k}{\partial x^i}
- Y^i \frac{\partial X^k}{\partial x^i} \right)
\frac{\partial f}{\partial x^k}.
\end{align*}
where $(a)$ is because of factoring out $\partial f / \partial x^k$. 

Since this holds for all smooth functions $f$, it follows that
\begin{align*}
[\b{X},\b{Y}]
=
\left( X^i \frac{\partial Y^k}{\partial x^i}
- Y^i \frac{\partial X^k}{\partial x^i} \right)
\frac{\partial}{\partial x^k}.
\end{align*}
According to Eq. (\ref{equation_partial_i}), this equation is euqal to Eq. (\ref{equation_Lie_bracket_coordinate_expression}). 
\end{proof}

\begin{lemma}[Lie bracket of coordinate basis vectors]
Let $\{\b{e}_i\}_{i=1}^n$ be the coordinate basis vectors associated with local coordinates $(x^1, \dots, x^n)$ on a smooth manifold. Then, for any $i,j$, the Lie bracket of coordinate basis vectors is zero:
\begin{align}\label{equation_Lie_bracket_coordinate_basis}
\boxed{
[\b{e}_i, \b{e}_j] = 0.
}
\end{align}
\end{lemma}
\begin{proof}
By definition, the Lie bracket of two vector fields $\b{X}$ and $\b{Y}$ is Eq. (\ref{equation_Lie_bracket}) for any smooth function $f$. If we consider $\b{e}_i$ and $\b{e}_j$ as vector fields $\b{X}$ and $\b{Y}$ in Eq. (\ref{equation_Lie_bracket}), we have:
\begin{align*}
[\b{e}_i, \b{e}_j](f) := \b{e}_i(\b{e}_j(f)) - \b{e}_j(\b{e}_i(f)). 
\end{align*}

According to Eq. (\ref{equation_coordinate_basis_vectors}), the basis vectors are $\b{e}_i = \frac{\partial}{\partial x^i}$ and $\b{e}_j = \frac{\partial}{\partial x^j}$. Then, for any smooth function $f$, we have:
\begin{align*}
\b{e}_i(\b{e}_j(f)) &= \frac{\partial}{\partial x^i} \left( \frac{\partial f}{\partial x^j} \right) = \frac{\partial^2 f}{\partial x^i \partial x^j}, \\
\b{e}_j(\b{e}_i(f)) &= \frac{\partial}{\partial x^j} \left( \frac{\partial f}{\partial x^i} \right) = \frac{\partial^2 f}{\partial x^j \partial x^i}.
\end{align*}

Since mixed partial derivatives commute for smooth functions:
\begin{align*}
\frac{\partial^2 f}{\partial x^i \partial x^j} = \frac{\partial^2 f}{\partial x^j \partial x^i},
\end{align*}
therefore,
\begin{align*}
[\b{e}_i, \b{e}_j](f) &= \b{e}_i(\b{e}_j(f)) - \b{e}_j(\b{e}_i(f)) = 0,
\end{align*}
for all smooth functions $f$.
\end{proof}

\begin{lemma}[Coordinate expression of Lie bracket]
The Lie bracket of vector fields $\b{X}, \b{Y} \in \mathfrak{X}(\mathcal{M})$ at point $\b{p} \in \mathcal{M}$ can be written as:
\begin{align}\label{equation_lie_bracket_directional_derivative}
\boxed{
[\b{X}, \b{Y}](\b{p})
=
D\b{Y}(\b{p})[\b{X}(\b{p})]
-
D\b{X}(\b{p})[\b{Y}(\b{p})],
}
\end{align}
where $D\b{Y}(\b{p})[\b{X}(\b{p})]$ and $D\b{X}(\b{p})[\b{Y}(\b{p})]$ are directional derivatives of vector field along vector field, defined in Definition \ref{definition_directional_derivative_vector_field_along_vector_field}.
\end{lemma}
\begin{proof}
Let $(x^1,\dots,x^n)$ be a local coordinate system around $\b{p} \in \mathcal{M}$. 
Write the vector fields in coordinates as:
\[
\b{X} = X^i \partial_i,
\qquad
\b{Y} = Y^i \partial_i.
\]
According to Proposition \ref{proposition_Lie_bracket_coordinates}, the Lie bracket has the coordinate expression:
\begin{equation}\label{equation_lie_bracket_coordinate_components_proof}
[\b{X},\b{Y}]
=
\left(
X^i \frac{\partial Y^k}{\partial x^i}
-
Y^i \frac{\partial X^k}{\partial x^i}
\right)\partial_k.
\end{equation}

Now we compute the right-hand side of Eq. \eqref{equation_lie_bracket_directional_derivative}.

Firstly, the vector field $\b{Y}$ can be viewed locally as the map:
\[
\b{p} \mapsto \b{Y}(\b{p}) = Y^k(\b{p})\,\partial_k|_{\b{p}}.
\]
According to Eq. (\ref{equation_DYpXp_XYk_p_partialkp}), we have:
\[
D\b{Y}(\b{p})[\b{X}(\b{p})]
=
\b{X}(Y^k)(\b{p})\,\partial_k|_{\b{p}}.
\]
By Definition \ref{definition_directional_derivative_function_along_vector_field}, the directional derivative of the scalar function $Y^k$ along $\b{X}$ is:
\[
\b{X}(Y^k)
=
X^i \partial_i Y^k.
\]
Therefore:
\begin{equation}\label{equation_directional_derivative_Y_along_X_components}
D\b{Y}(\b{p})[\b{X}(\b{p})]
=
\left(X^i \partial_i Y^k\right)(\b{p})\,\partial_k|_{\b{p}}.
\end{equation}

Similarly, we have:
\begin{equation}\label{equation_directional_derivative_X_along_Y_components}
D\b{X}(\b{p})[\b{Y}(\b{p})]
=
\left(Y^i \partial_i X^k\right)(\b{p})\,\partial_k|_{\b{p}}.
\end{equation}

Subtracting Eqs. \eqref{equation_directional_derivative_Y_along_X_components} and
\eqref{equation_directional_derivative_X_along_Y_components}, we obtain:
\begin{align*}
D\b{Y}(\b{p})[\b{X}(\b{p})]
-
&D\b{X}(\b{p})[\b{Y}(\b{p})]
= \\
&\left(
X^i \partial_i Y^k
-
Y^i \partial_i X^k
\right)(\b{p})\,\partial_k|_{\b{p}}.
\end{align*}
By Eq. \eqref{equation_lie_bracket_coordinate_components_proof}, this is exactly the coordinate expression of $[\b{X},\b{Y}](\b{p})$. Hence:
\[
[\b{X}, \b{Y}](\b{p})
=
D\b{Y}(\b{p})[\b{X}(\b{p})]
-
D\b{X}(\b{p})[\b{Y}(\b{p})].
\]
\end{proof}

\section{Metric Tensor}\label{section_metric_tensor}




\subsection{Definition of Metric Tensor}

\begin{definition}[Metric tensor --- coordinate-free definition]
Consider a Riemannian manifold $\mathcal{M}$. For each point $\b{p} \in \mathcal{M}$, the \textbf{Riemannian metric tensor} $g$ on the manifold is a bilinear, symmetric, and positive-definite map:
\begin{align}
\boxed{
g_{\b{p}} : T_{\b{p}}\mathcal{M} \times T_{\b{p}}\mathcal{M} \to \mathbb{R}.
}
\end{align}
For any two tangent vectors $\b{V}, \b{W} \in T_{\b{p}}\mathcal{M}$, the \textbf{inner product} is denoted by $g(\b{V}, \b{W})$, or $\langle \b{V}, \b{W} \rangle$, or $\langle \b{V}, \b{W} \rangle_g$, or $\langle \b{V}, \b{W} \rangle_{\b{p}}$.
\end{definition}

\begin{remark}[Notation of Riemannian manifold with a metric]
A Riemannian manifold is fully understandable by knowing its metric tensor at different points of manifold. That is why a Riemannian manifold $\mathcal{M}$ with a metric tensor $g$ is denoted by $(\mathcal{M}, g)$.
\end{remark}

\begin{definition}[Metric tensor --- coordinate-based definition]
Given a local coordinate system $\{x^1, \dots, x^n\}$ with the associated basis vectors $\{\partial_1, \dots, \partial_n\}$, the \textbf{components of the metric tensor $g$} are defined by:
\begin{align}\label{equation_g_components}
\boxed{
g_{ij} := g(\partial_i, \partial_j),
}
\end{align}
where $g(\cdot, \cdot)$ denotes:
\begin{align}\label{equation_g_inner_product}
\boxed{
g(\b{V}, \b{W}) = \langle \b{V}, \b{W} \rangle,
}
\end{align}
for two vectors $\b{V}$ and $\b{W}$, where $\langle \cdot, \cdot \rangle$ denotes inner product. 

In this basis, the inner product of two vectors $\b{V} = V^i \partial_i$ and $\b{W} = W^j \partial_j$ is computed as:
\begin{align}\label{equation_g_V_W}
\boxed{
g(\b{V}, \b{W}) := g_{ij} V^i W^j,
}
\end{align}
using the Einstein summation convention.
\end{definition}

According to Eqs. (\ref{equation_coordinate_basis_vectors}) and (\ref{equation_g_inner_product}), the Eq. (\ref{equation_g_components}) can be stated as:
\begin{align}\label{equation_g_e_e}
\boxed{
g_{ij} := \langle \partial_i, \partial_j \rangle \overset{(\ref{equation_coordinate_basis_vectors})}{=} \langle \b{e}_i, \b{e}_j \rangle,
}
\end{align}
where $\{\b{e}_1, \dots, \b{e}_n\}$ are the basis vectors of the coordinate system on the manifold and $\langle \cdot, \cdot \rangle$ denotes inner product.

\begin{corollary}[Metric tensor in a coordinate system with orthonormal bases]
Compare the Eqs. (\ref{equation_basis_orthonormal}) and (\ref{equation_g_e_e}). We conclude that the metric tensor in a coordinate system with orthonormal bases---i.e., the Cartesian coordinate system in Euclidean space (flat space)---is the Kronecker delta:
\begin{align}
g_{ij} =
\delta_{ij} = 
\left\{
    \begin{array}{ll}
        1 & \mbox{if } i = j, \\
        0 & \mbox{if } i \neq j.
    \end{array}
\right.
\end{align}
where $\delta_{ij}$ is the \textit{Kronecker delta} with lower indices, defined in Eq. (\ref{equation_Kronecker_delta}).
\end{corollary}

\begin{remark}[Metric at every point and coordinate formulas for its components]
A Riemannian metric is assigned at every point of the manifold.
More precisely, for every point $\b{p} \in \mathcal{M}$, the metric
$g_{\b{p}}$ is an inner product on the tangent space $T_p\mathcal{M}$:
\[
g_{\b{p}} : T_{\b{p}}\mathcal{M} \times T_{\b{p}}\mathcal{M} \to \mathbb{R}.
\]
Therefore, the metric is not just one fixed matrix in general; it is
a smoothly varying family of inner products, one inner product for
each tangent space.

If a local coordinate system $(x^1,\dots,x^n)$ is chosen around
$\b{p}$, with coordinate basis:
\[
\b{e}_i = \partial_i := \frac{\partial}{\partial x^i},
\]
then, the metric components are:
\[
g_{ij}(\b{p}) := g_p(\b{e}_i,\b{e}_j).
\]
Equivalently, if the point $\b{p}$ has local coordinates
$\b{x}=(x^1,\dots,x^n)$, we often write:
\[
g_{ij}(\b{x}) := g_p(\b{e}_i,\b{e}_j),
\]
to emphasize that the metric components may depend on the
coordinates of the point. Thus, in coordinates, the metric at point $\b{p}$ is represented by the matrix:
\[
\big[g_{ij}(\b{x})\big]_{i,j=1}^n.
\]

In many important manifolds, the metric has additional symmetry or
structure. In such cases, the metric may be described by a formula:
one inserts the coordinates of the point into the formula and obtains
the numerical values of the components $g_{ij}(\b{x})$ at that point.
For example, on a flat Euclidean space with Cartesian coordinates,
the metric components are constant everywhere:
\[
g_{ij}(\b{x}) = \delta_{ij}.
\]
Hence, the same metric matrix is obtained at every point.

On a curved manifold, however, the components $g_{ij}(\b{x})$ generally
vary from point to point. In some structured manifolds, such as many
matrix manifolds or many metrics used in general relativity, the metric may still admit a compact closed-form
expression because of the algebraic structure of the manifold. In
more complicated manifolds, such a closed-form expression may not be
available globally, and the metric may need to be computed locally,
numerically, or through a chosen coordinate chart. In all cases, the
intrinsic object is the family of inner products $g_{\b{p}}$, while the
matrix $[g_{ij}(\b{x})]_{i,j}$ is its coordinate representation in a chosen
chart.
\end{remark}

\begin{proposition}[Squared length of a vector]\label{proposition_squared_length_of_vector}
Consider a vector $\b{V}$. The squared length of the vector, which is $\langle\b{V}, \b{V}\rangle$, is:
\begin{align}
\boxed{
\langle\b{V}, \b{V}\rangle = V^i V^j g_{ij},
}
\end{align}
where $\langle\cdot,\cdot\rangle$ denotes the inner product and $g_{ij}$ is the metric tensor defined in Eq. (\ref{equation_g_e_e}).
\end{proposition}
\begin{proof}
\begin{align*}
\langle\b{V}, \b{V} \rangle\overset{(\ref{equation_contravariant_components})}{=} \langle V^i \b{e}_i, V^j \b{e}_j \rangle \overset{(a)}{=} V^i V^j \langle\b{e}_i, \b{e}_j\rangle \overset{(\ref{equation_g_e_e})}{=} V^i V^j g_{ij},
\end{align*}
where $(a)$ is because $V^i$ and $V^j$ are scalars and can be moved around in the multiplication and dot product. 
\end{proof}

\begin{remark}[Use of metric for calculating length of vector]
According to Proposition \ref{proposition_squared_length_of_vector}, metric $g$ can be used to calculate the squared length---and consequently the length---of a vector in a Riemannian manifold. 
\end{remark}

\subsection{Conversion of Contravariant and Covariant Components Using Metric Tensor}

\begin{lemma}[Conversion of contravariant components to covariant components by index lowering]
The contravariant components can be converted to covariant components by using the metric tensor as:
\begin{align}\label{equation_relation_contravariant_covariant}
\boxed{
V_i = g_{ij} V^j.
}
\end{align}
In other words, we can do index lowering using the metric tensor. 
\end{lemma}
\begin{proof}
\begin{align*}
V_i &\overset{(\ref{equation_covariant_components})}{=} \langle\b{V}, \b{e}_i\rangle \overset{(\ref{equation_contravariant_components})}{=} \langle V^j \b{e}_j, \b{e}_i \rangle \overset{(a)}{=} V^j \langle \b{e}_j, \b{e}_i \rangle  \\
&\overset{(b)}{=} V^j \langle\b{e}_i, \b{e}_j \rangle \overset{(\ref{equation_g_e_e})}{=} V^j g_{ij} \overset{(c)}{=} g_{ij} V^j,
\end{align*}
where $(a)$ is because $V^j$ is a scalar component and comes out of inner product, $(b)$ is because order does not matter in inner product, and $(c)$ is because both $V^j$ and $g_{ij}$ are scalars so their order does not matter in multiplication. 
\end{proof}

The metric tensor can be considered as a matrix where $i$ and $j$ index the dimension of space (manifold) along the row and column of matrix, respectively. 
To understand better, for example, in a three-dimensional space (manifold), Eq. (\ref{equation_relation_contravariant_covariant}) becomes:
\begin{align*}
\begin{bmatrix}
V_1\\
V_2\\
V_3
\end{bmatrix}
=
\begin{bmatrix}
g_{11} & g_{12} & g_{13} \\
g_{21} & g_{22} & g_{23} \\
g_{31} & g_{32} & g_{33}
\end{bmatrix}
\begin{bmatrix}
V^1\\
V^2\\
V^3
\end{bmatrix}.
\end{align*}
This is for going from contravariant to covariant components. 
We can use the inverse of matrix of tensor metric to go from covariant to contravariant components:
\begin{align*}
\begin{bmatrix}
V^1\\
V^2\\
V^3
\end{bmatrix}
=
\begin{bmatrix}
g_{11} & g_{12} & g_{13} \\
g_{21} & g_{22} & g_{23} \\
g_{31} & g_{32} & g_{33}
\end{bmatrix}^{-1}
\begin{bmatrix}
V_1\\
V_2\\
V_3
\end{bmatrix}.
\end{align*}
This gives us the definition of inverse metric tensor. By convention, the inverse of metric tensor $g_{ij}$ is denoted by $g^{ij}$ using superscripts instead of subscripts. 

\begin{definition}[Inverse metric tensor]
Considering that the metric tensor $g_{ij}$ is a matrix, the matrix inverse of metric tensor---also called the \textbf{inverse metric tensor} in short---is denoted by $g^{ij}$ and the multiplication of these two becomes identity matrix. In coordinate-based writing, we have:
\begin{align}\label{equation_metric_inverse}
\boxed{
g_{ik}\, g^{kj} = g^{kj}\, g_{ik} = \delta_i^j,
}
\end{align}
where the index $k$ is contracted and we have noticed that $g_{ij}$ and $g^{ij}$ are scalars (components of matrices) so their order in multiplication does not matter. The $\delta_i^j$ is the Kronecker delta with lower and upper indices, defined in Eq. (\ref{equation_Kronecker_delta_lower_upper}).
\end{definition}

\begin{figure*}[!t]
\centering
\includegraphics[width=6.5in]{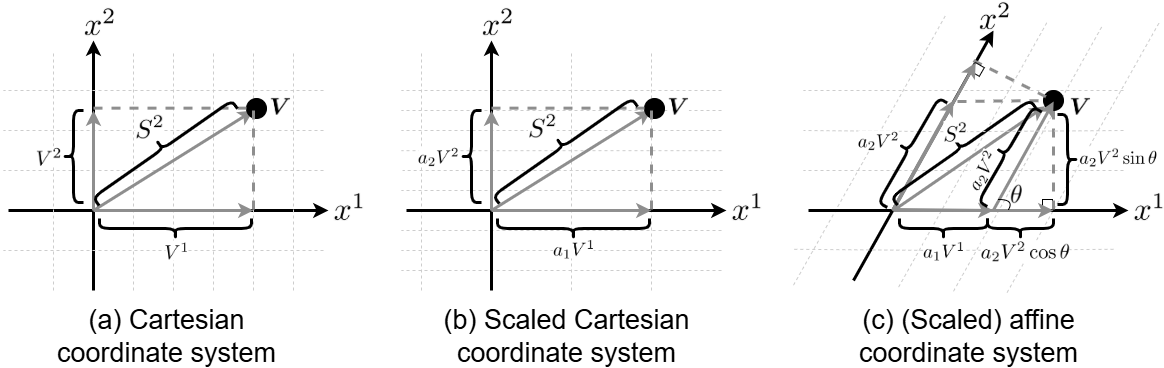}
\caption{Generalized Pythagorean theorem in (a) Cartesian coordinate system, (b) scaled Cartesian coordinate system, and (c) [scaled] affine coordinate system.}
\label{figure_coordinate_system_pythogrean}
\end{figure*}

\begin{corollary}
According to Eqs. (\ref{equation_metric_inverse}) and (\ref{equation_Kronecker_delta_lower_upper}), in an $n$-dimensional manifold, we have:
\begin{align}
g_{ij}\, g^{ij} = g^{ij}\, g_{ij} = \delta_i^i = \underbrace{1 + \dots + 1}_{n \text{ times}} = n,
\end{align}
where the Einstein convention is used, $n$ is the dimensionality of manifold, and we have noticed that $g_{ij}$ and $g^{ij}$ are scalars (components of matrices) so their order in multiplication does not matter. 
\end{corollary}

\begin{lemma}[Conversion of covariant components to contravariant components by index raising]\label{lemma_index_raising_by_metric}
The covariant components can be converted to contravariant components by using the inverse metric tensor as:
\begin{align}\label{equation_relation_covariant_contravariant}
\boxed{
V^i = g^{ij} V_j.
}
\end{align}
In other words, we can do index raising using the inverse metric tensor. 
\end{lemma}
\begin{proof}
\begin{align*}
&V_i \overset{(\ref{equation_relation_contravariant_covariant})}{=} g_{ij} V^j \overset{(a)}{\implies} g^{ki} V_i = g^{ki} g_{ij} V^j \\
&\overset{(\ref{equation_metric_inverse})}{\implies} g^{ki} V_i = \delta_j^k V^j \overset{(\ref{equation_index_substitution_delta})}{\implies} g^{ki} V_i = V^k \\
&\overset{(b)}{\implies} g^{ij} V_j = V^i, 
\end{align*}
where $(a)$ is because of left-multiplying the expression sides by $g^{ki}$ and $(b)$ is because of renaming the dummy variables $k \rightarrow i$ and $i \rightarrow j$.
\end{proof}

\subsection{Generalized Pythagorean Theorem}\label{section_Pythagorean_theorem}


Consider Fig. \ref{figure_coordinate_system_pythogrean}-a which illustrates a vector $\b{V}$ in a two-dimensional Cartesian coordinate system where $V^1$ and $V^2$ denote the contravariant components of the vector $\b{V}$.
According to the Pythagorean theorem, there is:
\begin{align*}
S^2 = (V^1)^2 + (V^2)^2 = V^1 V^1 + V^2 V^2,
\end{align*}
where $S^2$ is the squared length of the vector $\b{V}$. 
By Einstein convention, this equation can be stated as:
\begin{align}
S^2 &= \delta_{ij} V^i V^j \label{equation_s2_cartesian_coordinates} \\
&= \delta_{11} V^1 V^1 + \delta_{12} V^1 V^2 + \delta_{21} V^2 V^1 + \delta_{22} V^2 V^2, \nonumber
\end{align}
where $\delta_{ij}$ is the Kronecker delta with lower indices, defined in Eq. (\ref{equation_Kronecker_delta}).

Now, consider Fig. \ref{figure_coordinate_system_pythogrean}-b which depicts a vector $\b{V}$ in a two-dimensional scaled Cartesian coordinate system with orthogonal coordinates.
The coordinates are scaled as $a_1 V^1$ and $a_2 V^2$ with $a_1$ and $a_2$ are scalars of scaled coordinates. 
According to the Pythagorean theorem, there is:
\begin{align}\label{equation_s2_scaled_cartesian_coordinates}
S^2 = (a_1 V^1)^2 + (a_2 V^2)^2 = a_1^2 V^1 V^1 + a_2^2 V^2 V^2.
\end{align}

This time, consider Fig. \ref{figure_coordinate_system_pythogrean}-c illustrating a vector $\b{V}$ in a two-dimensional affine coordinate system where coordinates are not necessarily perpendicular. Let the angle of axes be denoted by $\theta$.
Again, the coordinates are scaled as $a_1 V^1$ and $a_2 V^2$ with $a_1$ and $a_2$ are scalars of scaled coordinates. 
As shown in Fig. \ref{figure_coordinate_system_pythogrean}-c, a right triangle is formed whose bases are $a_1 V^1 + a_2 V^2 \cos \theta$ and $a_2 V^2 \sin \theta$.
According to the Pythagorean theorem for this right triangle, there is:
\begin{align}
S^2 &= (a_1 V^1 + a_2 V^2 \cos \theta)^2 + (a_2 V^2 \sin \theta)^2 \nonumber \\
&\overset{(a)}{=} a_1^2 (V^1)^2 + a_2^2 (V^2)^2 \cos^2 \theta + 2 a_1 V^1 a_2 V^2 \cos \theta \nonumber \\
&~~~+ a_2^2 (V^2)^2 \sin^2 \theta \nonumber \\
&\overset{(b)}{=} a_1^2 (V^1)^2 + a_2^2 (V^2)^2 + 2 a_1 V^1 a_2 V^2 \cos \theta \nonumber \\
&\overset{(c)}{=} a_1^2 V^1 V^1 + a_1 a_2 (\cos \theta) V^1 V^2 \nonumber \\
&~~~+ a_2 a_1 (\cos \theta) V^2 V^1 + a_2^2 V^2 V^2, \label{equation_s2_affine_coordinates}
\end{align}
where $(a)$ is because of binomial theorem, $(b)$ is because $\sin^2 \theta + \cos^2 \theta = 1$, and $(c)$ is because $V^1$ and $V^2$ are scalars so $V^1 V^2 = V^2 V^1$. 

All Eqs. (\ref{equation_s2_cartesian_coordinates}), (\ref{equation_s2_scaled_cartesian_coordinates}), and (\ref{equation_s2_affine_coordinates}) can be generally stated as:
\begin{align*}
S^2 = g_{11} V^1 V^1 + g_{12} V^1 V^2 + g_{21} V^2 V^1 + g_{22} V^2 V^2,
\end{align*}
which is in two dimensions. In general, for any number of dimensions, we can say using the Einstein convention:
\begin{align}\label{equation_generalized_Pythagorean_theorem}
S^2 = g_{ij} V^i V^j,
\end{align}
where $i$ and $j$ are indices over the dimensions and $g_{ij}$ is the metric tensor. 

The metric tensor $g_{ij}$ can be considered as a matrix or tensor. For example, in Eqs. (\ref{equation_s2_cartesian_coordinates}), (\ref{equation_s2_scaled_cartesian_coordinates}), and (\ref{equation_s2_affine_coordinates}), the metrics are:
\begin{align*}
&\text{Metric in Eq. (\ref{equation_s2_cartesian_coordinates}): \quad } g_{ij} = 
\begin{bmatrix}
1 & 0\\
0 & 1
\end{bmatrix},
\\
&\text{Metric in Eq. (\ref{equation_s2_scaled_cartesian_coordinates}): \quad } g_{ij} = 
\begin{bmatrix}
a_1^2 & 0\\
0 & a_2^2
\end{bmatrix}, 
\\
&\text{Metric in Eq. (\ref{equation_s2_affine_coordinates}): \quad } g_{ij} = 
\begin{bmatrix}
a_1^2 & a_1 a_2 \cos \theta \\
a_2 a_1 \cos \theta & a_2^2
\end{bmatrix}.
\end{align*}

Equation (\ref{equation_generalized_Pythagorean_theorem}) is the generalized Pythagorean theorem stated below.
\begin{proposition}[Generalized Pythagorean theorem]
Assume the vector $\b{V}$ connects two points on a Riemannian manifold.
The squared distance between two points on the manifold is calculated as:
\begin{align}\label{equation_generalized_Pythagorean_theorem_2}
\boxed{
S^2 = g_{ij} V^i V^j.
}
\end{align}
where $V^i$ or $V^j$ denotes the contravariant components of vector $\b{V}$ and $g_{mn}$ denotes the metric tensor. 
\end{proposition}
\begin{proof}
Proof was provided in Section \ref{section_Pythagorean_theorem} with some examples.
\end{proof}

On a curvy manifold, we usually consider infinitesimal distances of close points as differential distances. In this case, the Eq. (\ref{equation_generalized_Pythagorean_theorem_2}) is written as:
\begin{align}\label{equation_generalized_Pythagorean_theorem_differential}
\boxed{
dS^2 = g_{ij} dx^i dx^j,
}
\end{align}
where $dS$, $dx^i$, and $dx^j$ are differentials (infinitesimal measurements) of $S$, $x^i$, and $x^j$, respectively. Here, we have replaced the vector components $V^i$ and $V^j$ with the infinitesimal differentials $dx^i$ and $dx^j$, respectively. 

\begin{remark}[The components of the metric depend on the coordinates of the point on the manifold]
The metric is defined at each point on the manifold, and therefore it may vary from point to point. Hence, Eq. (\ref{equation_generalized_Pythagorean_theorem_differential}) is often written as:
\begin{align}\label{equation_generalized_Pythagorean_theorem_differential2}
\boxed{
dS^2 = g_{ij}(x)\, dx^i dx^j,
}
\end{align}
where $g_{ij}(x)$ emphasizes that the components of the metric depend on the coordinates of the point on the manifold.

Moreover, as we will prove later in Lemma \ref{lemma_tensor_is_metric} and Corollary \ref{corollary_metric_tensor_independent_of_coordinate_system}, the metric tensor itself is independent of the choice of coordinate system, although its components $g_{ij}$ change when we change coordinates. What is intrinsic to the manifold is the geometry encoded by the metric, such as the curvature.
When we move on a curved manifold, the metric generally varies from point to point regardless of the coordinate system used. On a flat manifold, however, it is possible to choose coordinates in which the metric is the same everywhere (for example, Cartesian coordinates in Euclidean space), although in other coordinate systems the components $g_{ij}(x)$ may still vary with position.
\end{remark}

\begin{remark}
The Pythagorean theorem is a special case of the generalized Pythagorean theorem where the Riemannian manifold is the Euclidean space (flat space) with Cartesian coordinate system.
\end{remark}
\begin{proof}
The Pythagorean theorem applies in Cartesian coordinates where the metric tensor is the Kronecker delta defined in Eq. (\ref{equation_Kronecker_delta}). The Pythagorean theorem is Eq. (\ref{equation_s2_cartesian_coordinates}), which is Eq. (\ref{equation_generalized_Pythagorean_theorem_2}) with $g_{ij} = \delta_{ij}$. 
\end{proof}


\subsection{Properties of Metric Tensor}


\begin{lemma}[Symmetry of metric tensor]\label{lemma_metric_is_symmetric}
The metric tensor is symmetric:
\begin{align}\label{equation_metric_is_symmetric}
\boxed{
g_{ij} = g_{ji}.
}
\end{align}
\end{lemma}
\begin{proof}
\begin{align*}
g_{ij} \overset{(\ref{equation_g_e_e})}{=} \langle \b{e}_i, \b{e}_j\rangle \overset{(a)}{=} \langle \b{e}_j, \b{e}_i \rangle \overset{(\ref{equation_g_e_e})}{=} g_{ji}, 
\end{align*}
where $(a)$ is because of commutativity of the inner product.  
\end{proof}

\begin{corollary}[Number of unique elements of metric tensor]
In a metric tensor $g_{ij}$, each of the indices $i$ and $j$ run over $\{1, 2, \dots, n\}$ where $n$ is the dimensionality:
\begin{align*}
&i \in \{1, 2, \dots, n\}, \quad j \in \{1, 2, \dots, n\}. 
\end{align*}
Thus, the matrix of the tensor is an $n \times n$ matrix having $n^2$ elements. However, as the metric tensor is symmetric (see Lemma \ref{lemma_metric_is_symmetric}), the matrix of tensor is a symmetric matrix. Therefore, there are $n(n+1)/2$ unique elements, out of $n^2$ elements, in the metric tensor. In other words, the metric at each point on the manifold is represented by $n(n+1)/2$ numbers. 

For example, in the four-dimensional space-time manifold in general relativity, the metric tensor has $4^2=16$ elements, from which $(4 \times 5)/2 = 10$ elements are unique. Hence, the metric at each point on the space-time manifold is represented by $10$ numbers. 
\end{corollary}

\begin{lemma}[Metric is a tensor]\label{lemma_tensor_is_metric}
The metric is a tensor; that is why the metric is often called the metric tensor. 
\end{lemma}
\begin{proof}
This proof is inspired by \cite{susskind2025general}.
Consider two coordinate systems $X$ and $Y$ on the same manifold. The metric is a function of coordinates, so let the metric be denoted by $g_{ij}(x)$ and $\widetilde{g}_{ij}(y)$ in the coordinate systems $x$ and $y$, respectively. 

On the one hand, the generalized Pythagorean theorem in these two coordinate systems are:
\begin{align}
&dS^2 \overset{(\ref{equation_generalized_Pythagorean_theorem_differential2})}{=} g_{ij}(x)\, dx^i dx^j, \label{equation_generalized_Pythagorean_theorem_differential2_in_proof} \\
&dS^2 \overset{(\ref{equation_generalized_Pythagorean_theorem_differential2})}{=} \widetilde{g}_{ij}(y)\, dy^i dy^j, \label{equation_generalized_Pythagorean_theorem_differential2_in_proof_2}
\end{align}
where $dS^2$ is used in both expressions because the differential of distance, or length of two very close points, is the same regardless of change of coordinates. 

On the other hand, according to the chain rule in derivatives, we have:
\begin{align}
&dx^i = \frac{\partial x^i}{\partial y^j} dy^j, \label{equation_chain_rule_derivative_in_proof_dX_partialX_partialY_dY} \\
&dx^j = \frac{\partial x^j}{\partial y^i} dy^i. \label{equation_chain_rule_derivative_in_proof_dX_partialX_partialY_dY_2}
\end{align}
Substituting Eqs. (\ref{equation_chain_rule_derivative_in_proof_dX_partialX_partialY_dY}) and (\ref{equation_chain_rule_derivative_in_proof_dX_partialX_partialY_dY_2}) in Eq. (\ref{equation_generalized_Pythagorean_theorem_differential2_in_proof}) gives:
\begin{align}
dS^2 &\overset{(\ref{equation_generalized_Pythagorean_theorem_differential2_in_proof})}{=} g_{ij}(x)\, dx^i dx^j \nonumber \\
&\overset{(a)}{=} g_{ij}(x) \left(\frac{\partial x^i}{\partial y^j} dy^j\right) \left( \frac{\partial x^j}{\partial y^i} dy^i \right) \nonumber \\
&\overset{(b)}{=} \left(\frac{\partial x^i}{\partial y^j} \frac{\partial x^j}{\partial y^k} g_{ij}(x) \right) dy^i dy^j, \label{equation_dS2_g_partialX_partialX_dY_dY}
\end{align}
where $(a)$ is because of Eqs. (\ref{equation_chain_rule_derivative_in_proof_dX_partialX_partialY_dY}) and (\ref{equation_chain_rule_derivative_in_proof_dX_partialX_partialY_dY_2}), and $(b)$ is because of rearranging the terms in multiplication. 

Comparing Eqs. (\ref{equation_generalized_Pythagorean_theorem_differential2_in_proof_2}) and (\ref{equation_dS2_g_partialX_partialX_dY_dY}) gives:
\begin{align}\label{equation_g_tilde_partialX_partialX_g}
\widetilde{g}_{ij}(y) = \frac{\partial x^i}{\partial y^j} \frac{\partial x^j}{\partial y^k} g_{ij}(x). 
\end{align}
Equation (\ref{equation_g_tilde_partialX_partialX_g}) is the transformation of coordinates for covariant components of a tensor, stated in Eq. (\ref{equation_transformation_coordinates_covariant_Tik}). Therefore, the metric must be a tensor to satisfy the transformation of coordinates for covariant components of a tensor. Hence, metric is a tensor. 
\end{proof}

\begin{corollary}[Metric tensor object is independent of the choice of coordinate system]\label{corollary_metric_tensor_independent_of_coordinate_system}
The metric tensor, as a geometrical object, is independent of the choice of coordinate system (or the reference frame in the language of physics), according to Remark \ref{remark_tensor_coordinate_independent}. However, the components of the metric, $g_{ij}$, change under a coordinate transformation.
\end{corollary}



\subsection{Musical Isomorphisms and Index Raising/Lowering by Metric Tensor}

We can use the metric tensor to lower an index, i.e., bring an upper index down. Moreover, we can use the metric tensor to raise an index, i.e., bring a lower index up.
These two operations are stated in Eqs. (\ref{equation_index_lower}) and (\ref{equation_index_raise}), respectively. 

Another name for index lowering is flat ($\flat$ as in musical note) and another name for index raising is sharp ($\sharp$ as in musical note). 
This analogy is because, in music, $\flat$ lowers a note by a half-step and $\sharp$ raises a note by a half-step.
Using the notations $\flat$ and $\sharp$ for index lowering and raising is called \textit{musical isomorphism} in differential geometry. 

\begin{definition}[Musical Isomorphisms]
Let $(\mathcal{M}, g)$ be a Riemannian manifold. The metric $g$ induces two natural isomorphisms between the tangent space $T_{\b{p}}\mathcal{M}$ and the cotangent space $T_{\b{p}}^*\mathcal{M}$, collectively known as the \textbf{musical isomorphisms}:
\begin{enumerate}
    \item The \textbf{flat} operator $\flat: T_{\b{p}}\mathcal{M} \to T_{\b{p}}^*\mathcal{M}$ is defined for any vector $\b{V} \in T_{\b{p}}\mathcal{M}$ such that for all $\b{W} \in T_{\b{p}}\mathcal{M}$:
    \begin{equation}\label{equation_V_flat_W}
        \b{V}^\flat(\b{W}) = g(\b{V}, \b{W}).
    \end{equation}
    In local coordinates, if $\b{V} = V^j \partial_j$, then:
    \begin{align}
    \boxed{(\b{V}^\flat)_i = g_{ij} v^j.} 
    \end{align}
    \item The \textbf{sharp} operator $\sharp: T_{\b{p}}^*\mathcal{M} \to T_{\b{p}}\mathcal{M}$ is the inverse of the flat operator, $\sharp = (\flat)^{-1}$. For a covector $\omega \in T_{\b{p}}^*\mathcal{M}$, the vector $\omega^\sharp$ is the unique vector satisfying:
    \begin{equation}
        g(\omega^\sharp, w) = \omega(w), \quad \forall w \in T_{\b{p}}\mathcal{M}.
    \end{equation}
    In local coordinates, if $\omega = \omega_j dx^j$, then:
    \begin{align}
    \boxed{
    (\omega^\sharp)^i = g^{ij} \omega_j.
    }
    \end{align}
\end{enumerate}
\end{definition}

\begin{lemma}[Musical Isomorphisms, and index raising/lowering]
Let $(\mathcal{M}, g)$ be a Riemannian manifold. The metric tensor $g_{ij}$ and its inverse $g^{ij}$ define the musical isomorphisms $\flat: T_{\b{p}}\mathcal{M} \to T_{\b{p}}^*\mathcal{M}$ and $\sharp: T_{\b{p}}^*\mathcal{M} \to T_{\b{p}}\mathcal{M}$ via:
\begin{align}
    &\boxed{V_i = g_{ij} V^j} \quad (\text{Lowering / Flat } \flat) \label{equation_index_lower} \\
    &\boxed{V^i = g^{ij} V_j} \quad (\text{Raising / Sharp } \sharp) \label{equation_index_raise}
\end{align}
\end{lemma}
\begin{proof}
To prove Eq. \eqref{equation_index_lower}, let $\b{V} = V^j \partial_j \in T_{\b{p}}\mathcal{M}$. We define the covector $\b{V}^\flat$ such that for any $\b{W} = W^i \partial_i$, $\b{V}^\flat(\b{W}) = g(\b{V}, \b{W})$ (see Eq. (\ref{equation_V_flat_W})). 
We have:
\begin{align*}
V_i W^i &\overset{(\ref{equation_g_V_W})}{=} g(\b{V}, \b{W}) \overset{(\ref{equation_contravariant_components_2})}{=} g(V^j \partial_j, W^i \partial_i) \\
&\overset{(a)}{=} V^j W^i g(\partial_j, \partial_i) \overset{(\ref{equation_g_components})}{=} V^j W^i g_{ji} \overset{(b)}{=} (V^j g_{ji}) W^i \\
&\overset{(c)}{=} (V^j g_{ij}) W^i,
\end{align*}
where $(a)$ is because $V^j$ and $W^i$ are scalar components so they can come out of inner product in metric, $(b)$ is because of rearranging the terms in multiplication, and $(c)$ is because of Eq. (\ref{equation_metric_is_symmetric}). 
Therefore, we have:
\begin{align*}
V_i W^i = (V^j g_{ij}) W^i \overset{(a)}{\implies} V_i = V^j g_{ij}, 
\end{align*}
where $(a)$ is because this holds for all $W^i$.
This proves the Eq. \eqref{equation_index_lower}. 

To prove Eq. \eqref{equation_index_raise}, we start with Eq. \eqref{equation_index_lower} using a dummy index $k$: $V_k = g_{kj} V^j$. Contract both sides with the inverse metric $g^{ik}$:
\begin{equation*}
    g^{ik} V_k = g^{ik} g_{kj} V^j.
\end{equation*}
According to Eq. (\ref{equation_metric_inverse}), we have $g^{ik} g_{kj} = \delta^i_j$. Substituting this into the above equation gives:
\begin{equation*}
g^{ik} V_k = \delta^i_j V^j \overset{(\ref{equation_index_substitution_delta})}{=} V^i.
\end{equation*}
Thus, $V^i = g^{ij} V_j$ after relabeling the dummy index $k$ to $j$. This proves the Eq. \eqref{equation_index_raise}. 
\end{proof}

\section{Connection and Covariant Derivative}\label{section_connection_covariant_derivative}

\subsection{Intuition of Need for Christoffel Symbol}\label{section_intuion_need_for_Christoffel}

\begin{figure}[!h]
\centering
\includegraphics[width=3.2in]{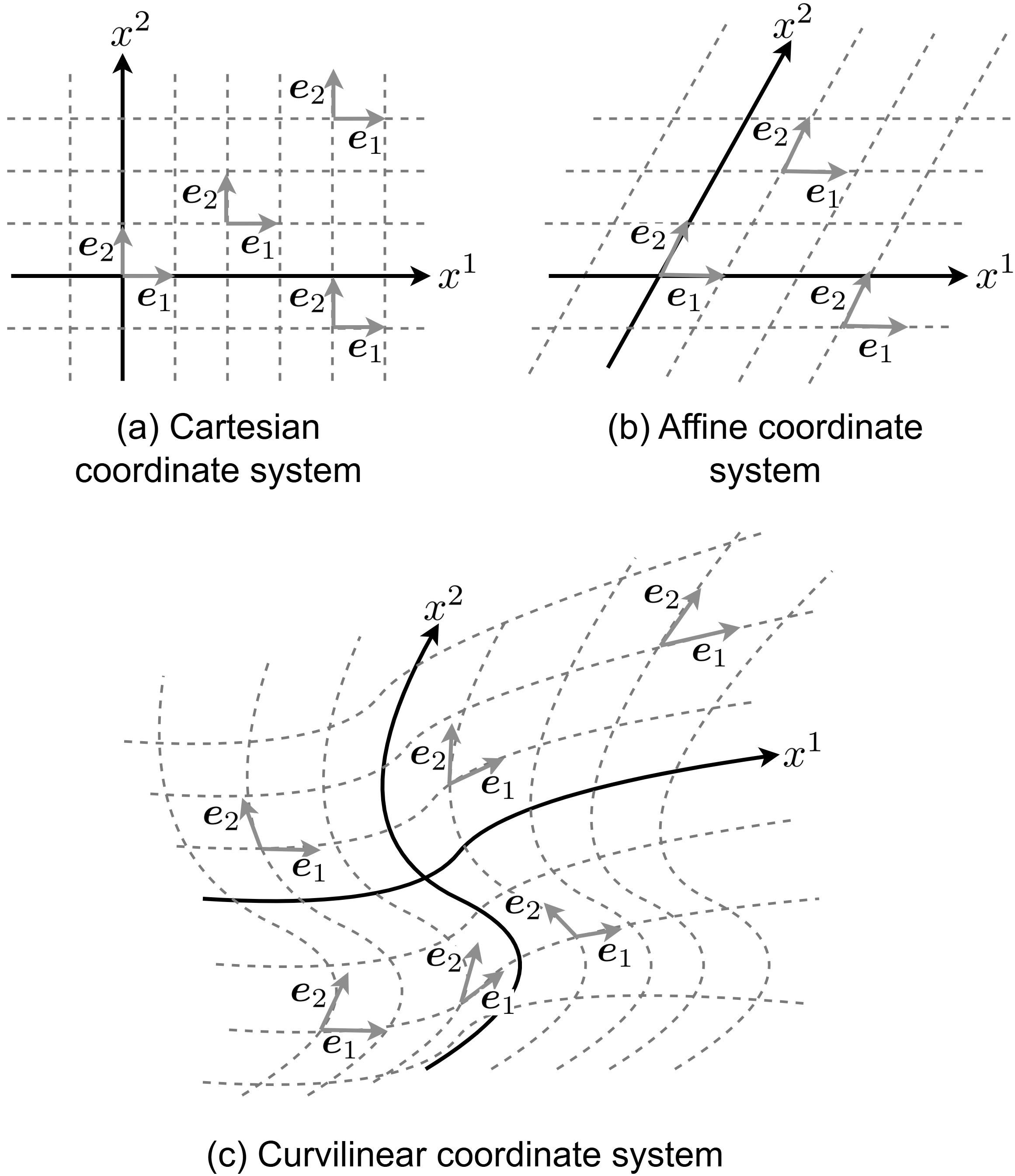}
\caption{Basis vectors $\{\b{e}_j\}_{j=1}^2$ do not change in different coordinates of the (a) Cartesian and (b) affine coordinate systems. However, the basis vectors $\{\b{e}_j\}_{j=1}^2$ change in different coordinates of the (c) curvilinear coordinate systems.}
\label{figure_coordinate_system_basis_vectors}
\end{figure}

Consider Figs. \ref{figure_coordinate_system_basis_vectors}-a and \ref{figure_coordinate_system_basis_vectors}-b with a Cartesian or an affine coordinate system in Euclidean space (flat space). As this figure shows, the basis vectors $\{\b{e}_j\}_{j=1}^n$ are fixed everywhere in both Cartesian and affine coordinate systems. This is because the angles between axes are fixed everywhere in these coordinate systems. 
As the basis vectors are fixed everywhere in Cartesian and affine coordinate systems, the derivative of each basis vector $\b{e}_j$ with respect to each coordinate $x^i$ is zero in these coordinate systems:
\begin{align}\label{equation_partial_e_partial_X}
\frac{\partial \b{e}_j}{\partial x^i} = \b{0}, \quad \forall i,j \in \{1, \dots, n\},
\end{align}
where $\b{0}$ denotes the zero vector.
Note that the derivative of the basis vector $\b{e}_j$ with respect to a scalar element $x^i$ is a vector because changing the element $x^i$ can cause change in each element of the basis vector $\b{e}_j$ \cite{ghojogh2023background}.

As a result, the derivative of metric tensor $g_{ij}$ with respect to each coordinate $x^k$ is also zero in Cartesian and affine coordinate systems:
\begin{align}
\frac{\partial g_{ij}}{\partial x^k} &\overset{(\ref{equation_g_e_e})}{=} \frac{\partial (\langle\b{e}_i, \b{e}_j\rangle)}{\partial x^k} \overset{(a)}{=} \left\langle\frac{\partial \b{e}_i}{\partial x^k}, \b{e}_j\right\rangle + \left\langle\frac{\partial \b{e}_j}{\partial x^k}, \b{e}_i\right\rangle \nonumber \\
&\overset{(\ref{equation_partial_e_partial_X})}{=} \langle\b{0}, \b{e}_j\rangle + \langle\b{0}, \b{e}_i\rangle = 0.
\end{align}

However, as shown in Fig. \ref{figure_coordinate_system_basis_vectors}-c, the basis vectors in a general curvilinear coordinate system change by moving to different points in the coordinate system. Therefore, the derivative of each basis vector $\b{e}_j$ with respect to each coordinate $x^i$ is not zero:
\begin{align}
\frac{\partial \b{e}_j}{\partial x^i} \neq \b{0}, \quad \forall i,j \in \{1, \dots, n\}.
\end{align}
As explained before, the derivative of an $n$-dimensional basis vector $\b{e}_j$ with respect to a scalar number $x^i$ is an $n$-dimensional vector \cite{ghojogh2023background}. Let this derivative vector be denoted by a vector $\b{W}$ which can be represented as:
\begin{align}\label{equation_W_for_Christoffel_symbol}
\b{W} \overset{(\ref{equation_contravariant_components})}{=} W^k \b{e}_k = W^1 \b{e}_1 + W^2 \b{e}_2 + \dots + W^n \b{e}_n,
\end{align}
where $W^k, \forall k \in \{1, \dots, n\}$ are the contravariant components of the vector $\b{W}$ and the Einstein convention is used. 

Therefore, we have:
\begin{align}\label{equation_partial_e_partialX_for_Christoffel_symbol}
\frac{\partial \b{e}_j}{\partial x^i} &= \b{W} \overset{(\ref{equation_W_for_Christoffel_symbol})}{=} W^k \b{e}_k \\
&= W^1 \b{e}_1 + W^2 \b{e}_2 + \dots + W^n \b{e}_n. \nonumber
\end{align}

In the Cartesian coordinate system in Euclidean space, we can denote (rename) $W^k$ in Eq. (\ref{equation_partial_e_partialX_for_Christoffel_symbol}) by $\Gamma^k_{ij}$ to emphasize that it is the $k$-th contravariant component of the derivative of the $j$-th basis vector $\b{e}_j$ with respect to the $i$-th coordinate $x^i$:
\begin{align}\label{equation_partial_e_partialX_Gamma_e}
\frac{\partial \b{e}_j}{\partial x^i} = \Gamma^k_{ij} \b{e}_k = \Gamma^1_{ij} \b{e}_1 + \Gamma^2_{ij} \b{e}_2 + \dots + \Gamma^n_{ij} \b{e}_n.
\end{align}
The $\Gamma^k_{ij}$ is called the \textit{connection coefficient} or \textit{Christoffel symbol}.
Note that, as we will see later, the expression $\partial \b{e}_j / \partial x^i = \Gamma^k_{ij} \b{e}_k$ is only correct for Cartesian coordinate system in euclidean space (flat space). It is not correct for a general curvilinear coordinate system. We will define the Christoffel Symbol, for a general curvilinear coordinate system, in the following. 

\subsection{Connection, Covariant Derivative, and Christoffel Symbol}

\subsubsection{Definitions of Connection and Covariant Derivative}

Recall the pointwise definition of vector field in Definition \ref{definition_poinstwise_vector_field} and Eq. (\ref{equation_poinstwise_vector_field}). It will be used here. 

\begin{definition}[Connection and covariant derivative]\label{definition_connection}
Let $\mathcal{M}$ be a smooth manifold and let $\mathfrak{X}(\mathcal{M})$ denote the space of smooth vector fields on $\mathcal{M}$. 
A \textbf{connection} (also called \textbf{affine connection}) on $\mathcal{M}$ is a map:
\begin{align}
\nabla : \mathfrak{X}(\mathcal{M}) \times \mathfrak{X}(\mathcal{M}) \to \mathfrak{X}(\mathcal{M}),
\end{align}
denoted by:
\begin{align}
\boxed{
(\b{X},\b{Y}) \mapsto \nabla_{\b{X}} \b{Y},
}
\end{align}
satisfying the following properties 
for all $\b{X},\b{Y},\b{Z} \in \mathfrak{X}(\mathcal{M})$ and $f,g \in C^\infty(\mathcal{M})$:
\begin{enumerate}
\item Linearity in the first argument: 
\begin{align}
\nabla_{f\b{X}+g\b{Y}} \b{Z} = f\nabla_{\b{X}} \b{Z} + g\nabla_{\b{Y}} \b{Z}.
\end{align}
\item Linearity in the second argument:
\begin{align}
\nabla_{\b{X}} (\b{Y}+\b{Z}) = \nabla_{\b{X}} \b{Y} + \nabla_{\b{X}} \b{Z}.
\end{align}
\item Leibniz rule (product rule):
\begin{align}\label{equation_Leibniz_rule_connection}
\nabla_{\b{X}} (f\b{Y}) = \b{X}(f)\b{Y} + f\nabla_{\b{X}} \b{Y}.
\end{align}
\end{enumerate}
A connection is a rule that tells how to differentiate vector fields on a manifold. 
A connection provides the operator $\nabla_{\b{X}} \b{Y}$ which means derivative of vector field $\b{Y}$ in the direction of vector field $\b{X}$.

The vector field $\nabla_{\b{X}} \b{Y}$ is called the \textbf{covariant derivative} of $\b{Y}$ in the direction of $\b{X}$.
The covariant derivative $\nabla_{\b{X}} \b{Y}$ means the rate of change of vector field $\b{Y}$ in the direction of vector field $\b{X}$.
\end{definition}

Note that, in some texts of literature, the connection $\nabla$ is denoted by $D$. This notation is used in many physics texts on general relativity; e.g., see \cite{susskind2025general}. 

\begin{remark}[Partial derivative versus covariant derivative]\label{remark_partial_derivative_vs_covariant_derivative}
The regular partial derivative is Eq. (\ref{equation_partial_i}), i.e., $\partial_i := \frac{\partial }{\partial x^i}$. 
This regular derivative is in Cartesian coordinate system in Euclidean space (flat space) because the basis vectors are fixed everywhere and only the vector may change. 

In a general curvilinear coordinate system, the regular partial derivative $\partial_i$ cannot be used because not only may the vector change, but also the basis vector change themselves. Therefore, a new derivative needs to be defined which works on any general curvilinear coordinate system, where it reduces to the regular partial derivatives in Cartesian coordinate system in Euclidean space (flat space). Such a derivative in general curvilinear coordinate system is called the covariant derivative, denoted by $\nabla_i$. 
\end{remark}

\begin{remark}[Covariant derivative with respect to coordinate basis vector]
The covariant derivative $\nabla_{\b{X}} \b{Y}$ captures the rate of change of $\b{Y}$ along vector field $\b{X}$. 
As an example, the covariant derivative $\nabla_{\b{e}_i} \b{Y}$ calculates the rate of change of $\b{Y}$ along the coordinate basis vector $\b{e}_i$ or equivalently along the coordinate $x^i$ (see Definition \ref{definition_coordinate_basis}). According to Eq. (\ref{equation_coordinate_basis_vectors}), there is $\b{e}_i = \partial_i$, so one can write $\nabla_{\partial_i} \b{Y}$ instead of $\nabla_{\b{e}_i} \b{Y}$. In short, people usually write $\nabla_i \b{Y}$ to denote the covariant derivative along the coordinate $x^i$. In summary, we have:
\begin{align}
\boxed{
\nabla_{\b{e}_i} \equiv \nabla_{\partial_i} \equiv \nabla_i
}.
\end{align}
\end{remark}

Covariant derivative was first proposed by \textit{Elwin Bruno Christoffel} (from the former German empire) in 1869 \cite{christoffel1869ueber} and then it was developed by \textit{Gregorio Ricci-Curbastro} (from Italy) in 1887 \cite{ricci1887sulla}.

\begin{lemma}[Covariant derivative of a scalar]
Consider a scalar, or a $(0,0)$ tensor, denoted by $S$.
The covariant derivative of a scalar $S$ is equal to its regular partial derivative:
\begin{align}\label{equation_covariant_derivative_scalar}
\boxed{
\nabla_i S = \partial_i S.
}
\end{align}
\end{lemma}
\begin{proof}
According to the discussion in Remark \ref{remark_partial_derivative_vs_covariant_derivative}, for vectors/tensors, we need Christoffel symbols because the basis changes from point to point.
But for scalars, there is no basis dependence, no indices, and no geometry to ``correct". So, the covariant derivative should reduce to the ordinary derivative.
\end{proof}

\begin{lemma}[Components of the covariant derivative of vector]
Let $\{\b{e}_j\}_{j=1}^n$ be a coordinate basis on a manifold.
The components of the covariant derivative of a vector $\b{V}$, with contravariant components $\{V^1, \dots, V^n\}$ are:
\begin{align}\label{equation_components_of_covariant_derivative}
\boxed{
\nabla_i \b{V} = (\nabla_i V^j) \b{e}_j.
}
\end{align}
\end{lemma}
\begin{proof}
According to Eq. (\ref{equation_contravariant_components}), a vector $\b{V}$ can be stated in terms of linear combination of basis vectors as $\b{V} = V^j \b{e}_j$. Thus:
\begin{align*}
\b{V} \overset{(\ref{equation_contravariant_components})}{=} V^j \b{e}_j \implies \nabla_i \b{V} = \nabla_i (V^j \b{e}_j) \overset{(a)}{=} (\nabla_i V^j) \b{e}_j,
\end{align*}
where $(a)$ is because of linearity of covariant derivative (see definition \ref{definition_connection}).
\end{proof}

\subsubsection{Definition of Christoffel Symbol}

\begin{definition}[Connection coefficients or Christoffel symbols]\label{definition_connection_coefficients_Christoffel_symbols}
Let $\{\b{e}_j\}_{j=1}^n$ be a coordinate basis on a manifold. Then, the covariant derivative of the basis vectors satisfies:
\begin{equation}\label{equation_Christoffel_symbol}
\boxed{
\begin{aligned}
\nabla_i \b{e}_j &= \Gamma^k_{ij} \b{e}_k \\
&= \Gamma^1_{ij} \b{e}_1 + \Gamma^2_{ij} \b{e}_2 + \dots + \Gamma^n_{ij} \b{e}_n,
\end{aligned}
}
\end{equation}
where the coefficients $\Gamma^k_{ij}, \forall i,j,k \in \{1, \dots, n\}$ are \textbf{connection coefficients}, also called \textbf{Christoffel symbols}.
In other words, the Christoffel symbol $\Gamma^k_{ij}$ is the $k$-th contravariant component of the covariant derivative of the $j$-th basis vector $\b{e}_j$ along the $i$-th coordinate $x^i$.
\end{definition}
\begin{proof}
The proof is similar to the discussion in Section \ref{section_intuion_need_for_Christoffel}. 

Let $\{\b{e}_j\}_{j=1}^n$ be the bases of the tangent space $T_{\b{p}}\mathcal{M}$. Since the covariant derivative of a vector field is again a vector field, we have:
\begin{align*}
\nabla_i \b{e}_j \in T_{\b{p}}\mathcal{M}.
\end{align*}
Because $\{\b{e}_j\}_{j=1}^n$ are basis vectors of $T_{\b{p}}\mathcal{M}$, any vector in $T_{\b{p}}\mathcal{M}$ can be expressed as a linear combination of these basis vectors. Therefore, there exist scalar functions $W^k_{ij}$ such that:
\begin{align*}
\nabla_i \b{e}_j = W^k_{ij} \b{e}_k.
\end{align*}
We define $\Gamma^k_{ij} := W^k_{ij}$; therefore:
\begin{align*}
\nabla_i \b{e}_j = \Gamma^k_{ij} \b{e}_k.
\end{align*}
\end{proof}

The \textit{Christoffel symbol} is named after \textit{Elwin Bruno Christoffel} (a mathematician from the former German empire in 1800s). In 1869, Elwin Bruno Christoffel published a fundamental work on covariant derivative in differential geometry \cite{christoffel1869ueber}. 

\begin{remark}
In some texts of literature, Eq. (\ref{equation_Christoffel_symbol}) is stated as:
\begin{align}
\nabla_i \b{e}_j = \Gamma^k_{ji} \b{e}_k \quad \text{ or } \quad \nabla_j \b{e}_i = \Gamma^k_{ij} \b{e}_k.
\end{align}
These equations are also correct if the derivative operator is torsion-free (and this is true in most of the regular cases). This is because, as we will see in Eq. (\ref{remark_Christoffel_sybmol_symmetric}), when the torsion of vector fields is zero on the manifold, the Christoffel symbol is symmetric, so $\Gamma^k_{ij} = \Gamma^k_{ji}$. This will be discussed fully in Section \ref{section_torsion_levi_civita}.
\end{remark}

\subsubsection{Covariant Derivative of Contravariant and Covariant Components}

\begin{proposition}[Covariant derivative of contravariant component]
Let $\b{V} = V^j \b{e}_j$ be a vector field expressed in coordinate bases $\{\b{e}_j\}_{j=1}^n$.
The covariant derivative of contravariant component $V^j$ in the direction $x^i$ is:
\begin{align}\label{equation_covariant_derivative_contravariant}
\boxed{
\nabla_i V^j = \partial_i V^j + \Gamma^j_{ik} V^k,
}
\end{align}
where $\nabla_i V^j$ denotes the $j$-th contravariant component of $\nabla_i \b{V}$:
\begin{align}
\boxed{
\nabla_i V^j \equiv (\nabla_i \b{V})^j.
}
\end{align}
\end{proposition}

\begin{proof}
The vector is expressed as:
\begin{align*}
\b{V} \overset{(\ref{equation_contravariant_components})}{=} V^j \b{e}_j. 
\end{align*}
We take the covariant derivative in the direction $x^i$:
\begin{align*}
\nabla_i \b{V} = \nabla_i (V^j \b{e}_j).
\end{align*}
Using the Leibniz rule (product rule) of the covariant derivative, i.e., Eq. (\ref{equation_Leibniz_rule_connection}), we have:
\begin{align}\label{equation_nabla_V_nabla_V_e_V_nabla_e}
\nabla_i \b{V} = (\nabla_i V^j) \b{e}_j + V^j (\nabla_i \b{e}_j).
\end{align}
The $V^j(X)$, in the first term of Eq. (\ref{equation_nabla_V_nabla_V_e_V_nabla_e}), is just a real-valued function for each fixed $j$.
So, in the first term of Eq. (\ref{equation_nabla_V_nabla_V_e_V_nabla_e}), the component $V^j$ can be treated as a scalar. So, according to Eq. (\ref{equation_covariant_derivative_scalar}), the derivative of the scalar component $V^j$ is the ordinary partial derivative\footnote{Note that Eq. (\ref{equation_covariant_derivative_contravariant_scalar}) does not contradict Eq. (\ref{equation_covariant_derivative_contravariant}) because, here in Eq. (\ref{equation_covariant_derivative_contravariant_scalar}), the $V^j(X)$ is just a real-valued function for each fixed $j$, so it can be treated as a scalar. However, $V^j$ in Eq. (\ref{equation_covariant_derivative_contravariant}) is the component of the covariant derivative of a vector. To highlight this distinction, some literature denotes $V^j$ in Eq. (\ref{equation_covariant_derivative_contravariant}) as $(\nabla_i \b{V})^j$ or $\b{V}^j_{;i}$.}:
\begin{align}\label{equation_covariant_derivative_contravariant_scalar}
\nabla_i V^j \overset{(\ref{equation_covariant_derivative_scalar})}{=} \partial_i V^j.
\end{align}
So, the first term becomes $(\partial_i V^j) \b{e}_j$.

In the second term of Eq. (\ref{equation_nabla_V_nabla_V_e_V_nabla_e}), we have $\nabla_i \b{e}_j = \Gamma^k_{ij} \b{e}_k$, according to Eq. (\ref{equation_Christoffel_symbol}). 
So, the second term becomes $V^j (\Gamma^k_{ij} \b{e}_k)$.

Therefore, Eq. (\ref{equation_nabla_V_nabla_V_e_V_nabla_e}) becomes:
\begin{align}\label{equation_nabla_V_nabla_V_e_V_nabla_e_2}
\nabla_i \b{V} = (\partial_i V^j) \b{e}_j + V^j \Gamma^k_{ij} \b{e}_k.
\end{align}
Now, in the second term, we rename the dummy variable $j$ to $k$ and vice versa ($j \rightarrow k, k \rightarrow j$):
\begin{align*}
V^j \Gamma^k_{ij} \b{e}_k \overset{(a)}{=} V^k \Gamma^j_{ik} \b{e}_j \overset{(b)}{=} \Gamma^j_{ik} V^k \b{e}_j,
\end{align*}
where $(a)$ is because of renaming dummy variables and $(b)$ is because of rearranging. 
Therefore, Eq. (\ref{equation_nabla_V_nabla_V_e_V_nabla_e_2}) becomes:
\begin{align}\label{equation_nabla_V_nabla_V_e_V_nabla_e_3}
\nabla_i \b{V} = (\partial_i V^j) \b{e}_j + \Gamma^j_{ik} V^k \b{e}_j \overset{(a)}{=} (\partial_i V^j + \Gamma^j_{ik} V^k) \b{e}_j,
\end{align}
where $(a)$ is because of factoring out $\b{e}_j$. 

According to Eq. (\ref{equation_components_of_covariant_derivative}), we have $\nabla_i \b{V} = (\nabla_i V^j) \b{e}_j$. Comparing this with Eq. (\ref{equation_nabla_V_nabla_V_e_V_nabla_e_3}) gives:
\begin{align*}
\nabla_i V^j = \partial_i V^j + \Gamma^j_{ik} V^k.
\end{align*}
\end{proof}

\begin{proposition}[Covariant derivative of covariant component]
Let $\b{V}$ be a vector field on the manifold.
The covariant derivative of covariant component $V_j$ in the direction $x^i$ is:
\begin{align}\label{equation_covariant_derivative_covariant}
\boxed{
\nabla_i V_j = \partial_i V_j - \Gamma^k_{ij} V_k,
}
\end{align}
where $\nabla_i V_j$ denotes the $j$-th covariant component of $\nabla_i \b{V}$:
\begin{align}
\boxed{
\nabla_i V_j \equiv (\nabla_i \b{V})_j.
}
\end{align}
\end{proposition}
\begin{proof}
Let $\b{V} = V^j \b{e}_j$ and $\b{W} = W^j \b{e}_j$ be vector fields. Let $V_j$ denote the $j$-th covariant component of vector $\b{V}$. 

On the one hand, using the Leibniz rule (product rule) for the covariant derivative, i.e., Eq. (\ref{equation_Leibniz_rule_connection}), we have:
\begin{align}\label{equation_nablaVW_nablaVW_V_nablaW}
\nabla_i (V_j W^j) &\overset{(\ref{equation_Leibniz_rule_connection})}{=} (\nabla_i V_j) W^j + V_j (\nabla_i W^j).
\end{align}
According to Eq. (\ref{equation_covariant_derivative_contravariant}), the $\nabla_i W^j$ in the second term is:
\begin{align*}
\nabla_i W^j = \partial_i W^j + \Gamma^j_{ik} W^k.
\end{align*}
Thus, Eq. (\ref{equation_nablaVW_nablaVW_V_nablaW}) becomes:
\begin{align}
\nabla_i (V_j W^j) &= (\nabla_i V_j) W^j + V_j (\partial_i W^j + \Gamma^j_{ik} W^k) \nonumber \\
&= (\nabla_i V_j) W^j + V_j \partial_i W^j + V_j \Gamma^j_{ik} W^k. \label{equation_nabla_V_W_1}
\end{align}

On the other hand, consider the scalar $V_j W^j$; it is a scalar because, by Einstein summation convention, we have $\sum_{j=1}^n V_j W^j$.
According to Eq. (\ref{equation_covariant_derivative_scalar}), since it is a scalar, its covariant derivative equals its ordinary derivative:
\begin{align}
\nabla_i (V_j W^j) \overset{(\ref{equation_covariant_derivative_scalar})}{=} \partial_i (V_j W^j) \overset{(a)}{=} (\partial_i V_j) W^j + V_j (\partial_i W^j), \label{equation_nabla_V_W_2}
\end{align}
where $(a)$ is because of the product rule in derivatives. 

Compare Eqs. (\ref{equation_nabla_V_W_1}) and (\ref{equation_nabla_V_W_2}) where they both are expressions for $\nabla_i (V_j W^j)$, so they are equal to each other:
\begin{align*}
&(\nabla_i V_j) W^j + \cancel{V_j \partial_i W^j} + V_j \Gamma^j_{ik} W^k = \\
&\quad\quad\quad\quad\quad\quad\quad\quad\quad\quad(\partial_i V_j) W^j + \cancel{V_j (\partial_i W^j)}.
\end{align*}
We cancel the common term $V_j \partial_i W^j$ on both sides to obtain:
\begin{align*}
(\nabla_i V_j) W^j + V_j \Gamma^j_{ik} W^k = (\partial_i V_j) W^j.
\end{align*}
We rename the dummy indices in the second term as:
\begin{align*}
V_j \Gamma^j_{ik} W^k \overset{(a)}{=} V_k \Gamma^k_{ij} W^j \overset{(a)}{=} \Gamma^k_{ij} V_k W^j,
\end{align*}
where $(a)$ is because of renaming dummy variables ($j \rightarrow k$ and $k \rightarrow j$) and $(b)$ is because of rearranging the terms in multiplication.
Thus:
\begin{align*}
&(\nabla_i V_j) W^j + \Gamma^k_{ij} V_k W^j = (\partial_i V_j) W^j. \\
&\implies (\nabla_i V_j) W^j + \Gamma^k_{ij} V_k W^j - (\partial_i V_j) W^j = 0.
\end{align*}
We factor $W^j$ out:
\begin{align*}
\left( \nabla_i V_j + \Gamma^k_{ij} V_k - \partial_i V_j \right) W^j = 0.
\end{align*}
Since this holds for all $W^j$, we conclude:
\begin{align*}
\nabla_i V_j + \Gamma^k_{ij} V_k - \partial_i V_j = 0.
\end{align*}
Therefore:
\begin{align*}
\nabla_i V_j = \partial_i V_j - \Gamma^k_{ij} V_k.
\end{align*}
\end{proof}

\begin{remark}[Interpretation of the covariant derivative]
The covariant derivatives in Eqs. (\ref{equation_covariant_derivative_contravariant}) and (\ref{equation_covariant_derivative_covariant}) are both interpreted as follows. 
As discussed in Remark \ref{remark_partial_derivative_vs_covariant_derivative}, in a curvilinear coordinate system, not only may the vector change, but also the basis vectors may change. If the basis vectors were fixed, we could use the regular partial derivative of vector with respect to the basis vectors to obtain the rate of changes of vector along the basis vectors. However, we need to add a term---the term $(+ \Gamma^j_{ik} V^k)$ or $(- \Gamma^k_{ij} V_k)$---to make the correction for the change of basis vectors. 
\end{remark}

\begin{lemma}[Christoffel symbols are not tensors]
The Christoffel symbols $\Gamma^k_{ij}$ do not define a tensor.
In fact, Christoffel symbols depend on the choice of coordinates.
\end{lemma}
\begin{proof}
Consider the covariant derivative of a vector field. i.e., Eq. (\ref{equation_covariant_derivative_contravariant}):
\begin{align*}
\nabla_i V^j = \partial_i V^j + \Gamma^j_{ik} V^k.
\end{align*}
Prearranging the terms, we obtain:
\begin{align*}
\Gamma^j_{ik} V^k = \nabla_i V^j - \partial_i V^j.
\end{align*}
The term $\nabla_i V^j$ is a tensor, while $\partial_i V^j$ is not a tensor. Therefore, their difference is not a tensor. Hence, $\Gamma^j_{ik} V^k$ is not a tensor.

Now, if $\Gamma^j_{ik}$ were a tensor, then its contraction with a vector $V^k$ would also be a tensor. Since $\Gamma^j_{ik} V^k$ is not a tensor, it follows that $\Gamma^j_{ik}$ itself cannot be a tensor.
\end{proof}

\begin{proposition}[Covariant derivative of a general $(r,s)$ tensor]\label{proposition_covariant_derivative_general_tensor}
Let $\b{T}$ be a tensor field of type $(r, s)$ expressed in coordinate bases as:
\begin{align*}
\b{T} = T^{i_1 \dots i_r}_{j_1 \dots j_s} \b{e}_{i_1} \otimes \dots \otimes \b{e}_{i_r} \otimes \b{e}^{j_1} \otimes \dots \otimes \b{e}^{j_s},
\end{align*}
where $\{\b{e}_{i_1}, \dots, \b{e}_{i_r}\}$ and $\{\b{e}^{j_1}, \dots, \b{e}^{j_s}\}$ are basis vectors and dual basis vectors, respectively. 

The covariant derivative of the component $T^{i_1 \dots i_r}_{j_1 \dots j_s}$ in the direction $x^k$ is:
\begin{align}\label{equation_covariant_derivative_general_tensor}
\boxed{
\begin{aligned}
\nabla_{\ell} T^{i_1 \dots i_r}_{j_1 \dots j_s} = \partial_{\ell} T^{i_1 \dots i_r}_{j_1 \dots j_s} &+ \sum_{a=1}^r \Gamma^{i_a}_{\ell k} T^{i_1 \dots k \dots i_r}_{j_1 \dots j_s} \\
&- \sum_{b=1}^s \Gamma^k_{\ell j_b} T^{i_1 \dots i_r}_{j_1 \dots k \dots j_s},
\end{aligned}
}
\end{align}
where the notation $T^{\dots k \dots}$ indicates that the index at the $a$-th upper position has been replaced by the summation index $k$, and $T_{\dots k \dots}$ indicates that the index at the $b$-th lower position has been replaced by the summation index $k$.

Note that the rule in Eq. (\ref{equation_covariant_derivative_general_tensor}) is that for every contravariant (upper) index, a term $\Gamma_{\ell k}^i T_{j_1, \dots, j_s}^{i_1 \dots k \dots i_r}$ is added (with positive sign), and for every covariant (lower) index, a $\Gamma_{\ell j}^k T_{j_1, \dots, j_s}^{i_1 \dots k \dots i_r}$ is subtracted (with negative sign).
\end{proposition}

\begin{proof}
Consider the tensor $\b{T}$ written as a tensor product of its components and basis elements. We apply the operator $\nabla_{\ell}$ and utilize the Leibniz rule (product rule), which distributes across the tensor product:
\begin{align*}
\nabla_{\ell} \b{T} &= \nabla_{\ell} \left( T^{i_1 \dots i_r}_{j_1 \dots j_s} \b{e}_{i_1} \otimes \dots \otimes \b{e}_{i_r} \otimes \b{e}^{j_1} \otimes \dots \otimes \b{e}^{j_s} \right) \\
&= (\partial_{\ell} T^{i_1 \dots i_r}_{j_1 \dots j_s}) \b{B} \\
&\quad + T^{i_1 \dots i_r}_{j_1 \dots j_s} \sum_{a=1}^r \b{e}_{i_1} \otimes \dots \otimes (\nabla_{\ell} \b{e}_{i_a}) \otimes \dots \otimes \b{e}^{j_s} \\
&\quad + T^{i_1 \dots i_r}_{j_1 \dots j_s} \sum_{b=1}^s \b{e}_{i_1} \otimes \dots \otimes (\nabla_{\ell} \b{e}^{j_b}) \otimes \dots \otimes \b{e}^{j_s},
\end{align*}
where $\b{B}$ represents the full basis tensor product:
\begin{align*}
\b{B} := \b{e}_{i_1} \otimes \dots \otimes \b{e}_{i_r} \otimes \b{e}^{j_1} \otimes \dots \otimes \b{e}^{j_s}.
\end{align*}
According to Eq. (\ref{equation_Christoffel_symbol}), the connection coefficients are $\nabla_{\ell} \b{e}_{i_a} = \Gamma^k_{\ell i_a} \b{e}_k$ and likewise, we can say $\nabla_{\ell} \b{e}^{j_b} = -\Gamma^{j_b}_{\ell k} \b{e}^k$ for dual basis. 
Substituting these into the expression gives:
\begin{align*}
\nabla_{\ell} \b{T} &= (\partial_{\ell} T^{i_1 \dots i_r}_{j_1 \dots j_s}) \b{B} \\
&\quad + \sum_{a=1}^r \Gamma^k_{\ell i_a} T^{i_1 \dots i_a \dots i_r}_{j_1 \dots j_s} (\dots \otimes \b{e}_k \otimes \dots) \\
&\quad - \sum_{b=1}^s \Gamma^{j_b}_{\ell k} T^{i_1 \dots i_r}_{j_1 \dots j_b \dots j_s} (\dots \otimes \b{e}^k \otimes \dots).
\end{align*}
By relabeling the dummy indices (swapping $k$ with $i_a$ in the first sum and $k$ with $j_b$ in the second sum) to restore the original basis indices $i_1 \dots i_r$ and $j_1 \dots j_s$ for the basis vectors, we factor out the common basis product:
\begin{align*}
\nabla_{\ell} \b{T} = \Big( \partial_{\ell} T^{i_1 \dots i_r}_{j_1 \dots j_s} &+ \sum_{a=1}^r \Gamma^{i_a}_{\ell k} T^{i_1 \dots k \dots i_r}_{j_1 \dots j_s} \\
&- \sum_{b=1}^s \Gamma^k_{\ell j_b} T_{j_1 \dots k \dots j_s}^{i_1 \dots j_r} \Big) \b{B}.
\end{align*}
The quantity inside the parentheses is the component $\nabla_\ell T^{i_1 \dots i_r}_{j_1 \dots j_s}$.
\end{proof}

\begin{corollary}[Covariant derivative of a $(1,1)$ tensor component]
Let $\b{T} = T^i_j \b{e}_i \otimes \b{e}^j$ be a tensor field of type $(1,1)$. 
According to Proposition \ref{proposition_covariant_derivative_general_tensor}, the covariant derivative of the component $T^i_j$ in the direction $x^\ell$ is:
\begin{align}\label{equation_covariant_derivative_tensor_1_1}
\boxed{
\nabla_\ell T^i_j = \partial_\ell T^i_j + \Gamma^i_{\ell k} T^k_j - \Gamma^k_{\ell j} T^i_k,
}
\end{align}
where $\nabla_\ell T^i_j$ denotes the $(1,1)$ component of the tensor $\nabla_\ell \b{T}$:
\begin{align}
\boxed{
\nabla_\ell T^i_j \equiv (\nabla_\ell \b{T})^i_j.
}
\end{align}
\end{corollary}

\begin{corollary}[Covariant derivative of a $(0,2)$ tensor component]
Let $\b{T} = T_{ij} \b{e}_i \otimes \b{e}_j$ be a tensor field of type $(0,2)$. 
According to Proposition \ref{proposition_covariant_derivative_general_tensor}, the covariant derivative of the component $T_{ij}$ in the direction $x^\ell$ is:
\begin{align}\label{equation_covariant_derivative_tensor_0_2}
\boxed{
\nabla_\ell T_{ij} = \partial_\ell T_{ij} - \Gamma^k_{\ell i} T_{kj} - \Gamma^k_{\ell j} T_{ik},
}
\end{align}
where $\nabla_\ell T_{ij}$ denotes the $(0,2)$ component of the tensor $\nabla_\ell \b{T}$:
\begin{align}
\boxed{
\nabla_\ell T_{ij} \equiv (\nabla_\ell \b{T})_{ij}.
}
\end{align}
\end{corollary}

\subsection{Torsion and the Levi-Civita Connection}\label{section_torsion_levi_civita}

\begin{definition}[Coordinate-free definition of torsion tensor]\label{definition_torsion}
Let $\b{X}$ and $\b{Y}$ be vector fields on a manifold. 
The \textbf{torsion tensor} of a connection $\nabla$ is defined as:
\begin{align}\label{equation_torsion}
\boxed{
T(\b{X},\b{Y}) := \nabla_{\b{X}} \b{Y} - \nabla_{\b{Y}} \b{X} - [\b{X},\b{Y}],
}
\end{align}
where $[\b{X},\b{Y}]$ is the Lie bracket of the vector fields (see Definition \ref{definition_Lie_bracket}). 
\end{definition}

\begin{proposition}[Coordinate-based expression of torsion tensor]
Consider coordinate basis vectors $\{\b{e}_i\}_{i=1}^n$ on a manifold.
The torsion tensor of the coordinate basis vectors is:
\begin{align}
\boxed{
T(\b{e}_i,\b{e}_j) = \nabla_i\b{e}_j - \nabla_j\b{e}_i.
}
\end{align}

Therefore, the coordinate expression of the torsion tensor of the coordinate basis vectors is:
\begin{align}\label{equation_torsion_Christoffel_symbol}
\boxed{
T_{ij}^k = \Gamma_{ij}^k - \Gamma_{ji}^k,
}
\end{align}
where $i,j,k \in \{1, \dots, n\}$.
This implies that the torsion measures the failure of the Christoffel symbol (connection coefficient) to be symmetric.
\end{proposition}
\begin{proof}
If we consider $\b{e}_i$ and $\b{e}_j$ as vector fields $X$ and $Y$ in Eq. (\ref{equation_torsion}), we have:
\begin{align*}
T(\b{e}_i,\b{e}_j) = \nabla_i\b{e}_j - \nabla_j\b{e}_i - [\b{e}_i,\b{e}_j],
\end{align*}
where $[\b{e}_i,\b{e}_j]$ is the Lie bracket of basis vectors.
According to Eq. (\ref{equation_Lie_bracket_coordinate_basis}), $[\b{e}_i,\b{e}_j] = 0$. Thus:
\begin{align}\label{equation_T_e_e_partiale_partiale}
T(\b{e}_i,\b{e}_j) = \nabla_i\b{e}_j - \nabla_j\b{e}_i. 
\end{align}

According to Eq. (\ref{equation_contravariant_components}), we have $\b{V} = V^k \b{e}_k$. Therefore, we can state the tensor $T(\b{e}_i,\b{e}_j)$ as:
\begin{align}\label{equation_T_e_e_T_e}
T(\b{e}_i,\b{e}_j) = T_{ij}^k \b{e}_k,
\end{align}
where the lower indices $i,j$ come from the indices of two input vectors $\b{e}_i$ and $\b{e}_j$. 

Moreover, according to Eq. (\ref{equation_Christoffel_symbol}), we have:
\begin{align}\label{equation_nablae_Gammae_nablae_Gammae}
\nabla_i\b{e}_j = \Gamma_{ij}^k \b{e}_k, \quad \nabla_j\b{e}_i = \Gamma_{ji}^k \b{e}_k.
\end{align}
Substituting Eqs. (\ref{equation_T_e_e_T_e}) and (\ref{equation_nablae_Gammae_nablae_Gammae}) in Eq. (\ref{equation_T_e_e_partiale_partiale}) gives $T_{ij}^k = \Gamma_{ij}^k - \Gamma_{ji}^k$. 
\end{proof}

\begin{definition}[Levi-Civita connection --- coordinate-free definition]\label{definition_Levi_Civita_connection_coordinate_free}
Let $(\mathcal{M}, g)$ be a Riemannian manifold. The \textbf{Levi-Civita connection} is the unique linear connection $\nabla$ on $\mathcal{M}$ that satisfies the following two conditions for all vector fields $\b{X}, \b{Y}, \b{Z} \in \mathfrak{X}(\mathcal{M})$:
\begin{enumerate}
    \item \textbf{Metric compatibility:} The connection preserves the Riemannian metric:
    \begin{align}\label{equation_metric_compatibility}
    \boxed{
    \b{X}(g(\b{Y}, \b{Z})) = g(\nabla_{\b{X}} \b{Y}, \b{Z}) + g(\b{Y}, \nabla_{\b{X}} \b{Z}),
    }
    \end{align}
    where $g(\b{Y}, \b{Z}): \mathcal{M} \rightarrow \mathbb{R}$ is the metric as a smooth scalar function on manifold, and $\b{X}(g(\b{Y}, \b{Z}))$ is the directional derivative of scalar function $g(\b{Y}, \b{Z})$ in the direction of $\b{X}$ (see Definition \ref{definition_directional_derivative_function_along_vector_field}). 
    Equation (\ref{equation_metric_compatibility}) is basically the product rule for derivatives. It implies that the connection $\nabla$ preserves the geometry of the manifold as we move along a curve on the manifold.
    \item \textbf{Torsion-free (zero torsion):} The torsion tensor $T$ vanishes identically:
    \begin{align}
    T(\b{X},\b{Y})& \overset{(\ref{equation_torsion})}{=} \nabla_{\b{X}} \b{Y} - \nabla_{\b{Y}} \b{X} - [\b{X},\b{Y}] = \b{0}  \nonumber \\
    &\implies \boxed{\nabla_{\b{X}} \b{Y} - \nabla_{\b{Y}} \b{X} = [\b{X}, \b{Y}],} \label{equation_torsion_free}
    \end{align}
    where $[\b{X}, \b{Y}]$ denotes the Lie bracket of $\b{X}$ and $\b{Y}$.
    According to Eq. (\ref{equation_torsion_Christoffel_symbol}), the Eq. (\ref{equation_torsion_free}) implies that the Christoffel symbols (connection coefficients) are symmetric.
\end{enumerate}
\end{definition}

The \textit{Levi-Civita connection} is named after \textit{Tullio Levi-Civita} (from Italy), although it was originally discovered by \textit{Elwin Bruno Christoffel} (from the former German empire) in 1869 \cite{christoffel1869ueber}. Levi-Civita has important works on parallel transport \cite{levi1916nozione}, which will be introduced in Section \ref{section_parallel_transport}.

\begin{remark}[Symmetry of Christoffel symbols in torsion-free connection]\label{remark_Christoffel_sybmol_symmetric}
According to Eq. (\ref{equation_torsion_Christoffel_symbol}), the Eq. (\ref{equation_torsion_free}) implies that the Christoffel symbols (connection coefficients) are symmetric.
In a coordinate basis $\{\b{e}_i\}_{i=1}^n$, the torsion-free property of Levi-Civita connection implies that the Christoffel symbols are symmetric in their lower indices:
\begin{align}\label{equation_symmtery_Christoffel}
\boxed{
\Gamma_{ij}^k = \Gamma_{ji}^k.
}
\end{align}
Note that in many regular cases, torsion is zero and the Christoffel symbols are symmetric. 

\end{remark}

\begin{definition}[Levi-Civita connection --- coordinate-based definition]
The \textbf{Levi-Civita connection} is a special affine connection determined uniquely by the metric.
It satisfies two conditions:
\begin{itemize}
\item \textbf{Metric compatibility}:
\begin{align}\label{equation_metric_compatibility_nabla_g_zero}
\boxed{
\nabla_k\, g_{ij} = 0.
}
\end{align}
In other words, the covariant derivative of the metric tensor itself is zero. This ensures that the metric tensor acts like a constant during differentiation, which is what allows lengths and angles to be preserved under parallel transport. In other words, it means the metric is preserved under parallel transport (we will introduce parallel transport in Section \ref{section_parallel_transport}).
\item \textbf{Torsion-free (zero torsion)}:
\begin{align}\label{equation_zero_torsion_Christoffel_symmetric}
\boxed{
\Gamma_{ij}^k = \Gamma_{ji}^k.
}
\end{align}
According to Remark \ref{remark_Christoffel_sybmol_symmetric}, the torsion-free condition simplifies to the symmetry of the Christoffel symbols in their lower indices.
\end{itemize}
\end{definition}

\begin{proposition}[Partial derivative of metric tensor]
In a Levi-Civita connection, where metric compatibility or Eq. (\ref{equation_metric_compatibility_nabla_g_zero}) is satisfied:
\begin{align*}
\nabla_k\, g_{ij} = 0,
\end{align*}
then the partial derivative of metric tensor is: 
\begin{align}\label{equation_partial_derivative_metric}
\boxed{
\partial_\ell g_{ij} = \Gamma^k_{\ell i} g_{kj} + \Gamma^k_{\ell j} g_{ik}.
}
\end{align}
\end{proposition}
\begin{proof}
According to Eq. (\ref{equation_covariant_derivative_tensor_0_2}), the covariant derivative of the metric tensor $g_{ij}$ is:
\begin{align*}
\nabla_\ell g_{ij} = \partial_\ell g_{ij} - \Gamma^k_{\ell i} g_{kj} - \Gamma^k_{\ell j} g_{ik}.
\end{align*}
\end{proof}

\begin{corollary}[Number of unique elements of Christoffel symbols]
In the Christoffel symbols $\Gamma^k_{ij}$, each of the indices $i,j,k$ runs
over $\{1,2,\dots,n\}$ where $n$ is the dimensionality:
\begin{align*}
i,j,k \in \{1,2,\dots,n\}.
\end{align*}
Hence, there are $n^3$ elements in total. However, the Christoffel symbols are
symmetric in the lower indices (i.e., $\Gamma^k_{ij}=\Gamma^k_{ji}$).
Therefore, for every fixed $k$, the indices $(i,j)$ form a symmetric
$n\times n$ matrix with $n(n+1)/2$ unique elements. Consequently, the total
number of unique Christoffel symbols is:
\begin{align}
n \times \frac{n(n+1)}{2} = \frac{n^2(n+1)}{2}.
\end{align}
For example, in the four-dimensional space-time manifold in general relativity, there are $4^3=64$ Christoffel symbols in total, of which $4\times(4\times5)/2=40$ are unique.
\end{corollary}

\subsection{Second-Order Covariant Derivative}

\begin{proposition}[Second-order covariant derivative  of contravariant component]
Let $\mathcal{M}$ be a Riemannian manifold with a Levi-Civita connection $\nabla$. For a smooth vector field $V$, the components of the second-order covariant derivative $\nabla^2 V$ in a local coordinate system $(x^1, \dots, x^n)$ are given by:
\begin{align}\label{equation_second_derivative_contravariant}
\boxed{
\nabla_\ell (\nabla_i V^j) = \partial_\ell (\nabla_i V^j) + \Gamma^j_{\ell k} (\nabla_i V^k) - \Gamma^k_{\ell i} (\nabla_k V^j),
}
\end{align}
where $\nabla_i V^j = \partial_i V^j + \Gamma^j_{ik} V^k$, according to Eq. (\ref{equation_covariant_derivative_contravariant}).
\end{proposition}
\begin{proof}
Let $\b{T} = T^j_i \b{e}_j \otimes \b{e}^i$ be the $(1,1)$ tensor defined by the first covariant derivative of $\b{V}$, such that:
\begin{align}\label{equation_T_j_i_nabla_V}
T^j_i := \nabla_i V^j.
\end{align}
According to Eq. (\ref{equation_covariant_derivative_tensor_1_1}), the covariant derivative of a $(1,1)$ tensor with components $T^j_i$ along the $\ell$-th coordinate direction is:
\begin{equation*}
\nabla_\ell T^j_i = \partial_\ell T^j_i + \Gamma^j_{\ell k} T^k_i - \Gamma^k_{\ell i} T^j_k,
\end{equation*}
According to Eq. (\ref{equation_T_j_i_nabla_V}), we have $T^j_i = \nabla_i V^j$, $T^k_i = \nabla_i V^k$, and $T^j_k = \nabla_k V^j$.
Substituting these back into the expression gives:
\begin{equation*}
\nabla_\ell (\nabla_i V^j) = \partial_\ell (\nabla_i V^j) + \Gamma^j_{\ell k} (\nabla_i V^k) - \Gamma^k_{\ell i} (\nabla_k V^j).
\end{equation*}
\end{proof}

\begin{proposition}[Second-order covariant derivative of covariant component]
Let $\mathcal{M}$ be a Riemannian manifold with a Levi-Civita connection $\nabla$. For a smooth vector field $\b{V}$, the components of the second-order covariant derivative $\nabla^2 \b{V}$ for the covariant components $V_j$ in a local coordinate system $(x^1, \dots, x^n)$ are given by:
\begin{equation}\label{equation_second_derivative_covariant}
\boxed{
\nabla_\ell (\nabla_i V_j) = \partial_\ell (\nabla_i V_j) - \Gamma^k_{\ell i} (\nabla_k V_j) - \Gamma^k_{\ell j} (\nabla_i V_k),
}
\end{equation}
where $\nabla_i V_j = \partial_i V_j - \Gamma^k_{ij} V_k$, according to Eq. (\ref{equation_covariant_derivative_covariant}).
\end{proposition}
\begin{proof}
Let $\b{S}$ be a $(0,2)$ tensor field defined by the first covariant derivative of the covariant components of $\b{V}$. We denote its components by:
\begin{align}\label{equation_S_ij_nabla_V}
S_{ij} := \nabla_i V_j.
\end{align} 

According to Eq. (\ref{equation_covariant_derivative_tensor_0_2}), the covariant derivative of a $(0,2)$ tensor with components $S_{ij}$ along the $\ell$-th coordinate direction is:
\begin{equation*}
\nabla_\ell S_{ij} = \partial_\ell S_{ij} - \Gamma^k_{\ell i} S_{kj} - \Gamma^k_{\ell j} S_{ik}.
\end{equation*}
According to Eq. (\ref{equation_S_ij_nabla_V}), we have $S_{ij} = \nabla_i V_j$, $S_{kj} = \nabla_k V_j$, and $S_{ik} = \nabla_i V_k$.
Substituting these back into the expression gives:
\begin{equation*}
\nabla_\ell (\nabla_i V_j) = \partial_\ell (\nabla_i V_j) - \Gamma^k_{\ell i} (\nabla_k V_j) - \Gamma^k_{\ell j} (\nabla_i V_k).
\end{equation*}
\end{proof}

\subsection{Two Kinds of Christoffel Symbols}

In Riemannian geometry, there are two kinds of Christoffel symbols, which are related by the metric tensor $g_{ij}$ and its inverse $g^{kl}$.
These two kinds are called Christoffel symbols of the first kind and second kind. 
The Christoffel symbols of the second kind are the connection coefficients of the Levi-Civita connection; we already defined the second kind in Definition \ref{definition_connection_coefficients_Christoffel_symbols}.
The Christoffel symbols of the first kind are new and we define them in the following. 

\begin{definition}[Christoffel symbols of the \underline{first} kind]\label{definition_Christoffel_symbol_first_kind}
The \textbf{Christoffel symbols of the first kind}, denoted by $\Gamma_{ijk}$, are defined as:
\begin{align}\label{equation_christoffel_first_kind}
\boxed{
\Gamma_{ijk} := g_{i\ell} \Gamma_{jk}^\ell,
}
\end{align}
where $\Gamma_{jk}^\ell$ are the Christoffel symbols of the second kind, i.e., the connection coefficients of the Levi-Civita connection (see Definition \ref{definition_Christoffel_symbol_second_kind}).
In other words, the Christoffel symbols of the first kind are obtained by lowering the index of the second kind.
\end{definition}

\begin{definition}[Christoffel symbols of the \underline{second} kind]\label{definition_Christoffel_symbol_second_kind}
The \textbf{Christoffel symbols of the second kind}, denoted by $\Gamma^k_{ij}$, are the connection coefficients of the Levi-Civita connection. They are obtained by raising the index of the first kind (see Eq. (\ref{equation_index_raise})):
\begin{align}\label{equation_christoffel_second_kind}
\boxed{
\Gamma^k_{ij} := g^{k\ell} \Gamma_{ij\ell},
}
\end{align}
where $\Gamma_{ij\ell}$ are the Christoffel symbols of the second kind (see Definition \ref{definition_Christoffel_symbol_first_kind}).
In other words, the Christoffel symbols of the second kind are obtained by raising the index of the first kind.

As we discussed in Definition \ref{definition_connection_coefficients_Christoffel_symbols} and Eq. (\ref{equation_Christoffel_symbol}), the Christoffel symbols of the second kind are defined as $\nabla_i \b{e}_j = \Gamma^k_{ij} \b{e}_k$, where $\{\b{e}_j\}_{j=1}^n$ are coordinate basis vectors on the manifold.
\end{definition}

\subsection{Christoffel Symbol in Terms of Metric Tensor}

\begin{proposition}[Christoffel symbols of second kind in terms of metric tensor]
When Levi-Civita connection---which satisfies the metric compatibility and Christoffel symbol symmetry---is used, the Christoffel symbols of second kind can be stated in terms of metric tensor:
\begin{align}\label{equation_Christoffel_secondType_metric}
\boxed{
\Gamma^k_{ij} = \frac{1}{2} g^{k\ell} \left( \partial_i g_{j\ell} + \partial_j g_{i\ell} - \partial_\ell g_{ij} \right).
}
\end{align}
\end{proposition}
\begin{proof}
The following proof is inspired by the proof provided in \cite{susskind2025general}.

When we have metric compatibility as in Eq. (\ref{equation_metric_compatibility_nabla_g_zero}), then according to Eq. (\ref{equation_partial_derivative_metric}), we have:
\begin{align}\label{equation_partialg_Gammag_Gammag_1}
\partial_\ell g_{ij} - \Gamma^k_{\ell i} g_{kj} - \Gamma^k_{\ell j} g_{ik} = 0. 
\end{align}
We rename (swap) the dummy variables $\ell \to i$ and $i \to \ell$ in Eq. (\ref{equation_partialg_Gammag_Gammag_1}):
\begin{align*}
\partial_i g_{\ell j} - \Gamma^k_{i \ell} g_{kj} - \Gamma^k_{i j} g_{\ell k} = 0. 
\end{align*}
According to Eq. (\ref{equation_zero_torsion_Christoffel_symmetric}), we have $\Gamma^k_{i \ell} = \Gamma^k_{\ell i}$, so we can replace $\Gamma^k_{i \ell}$ with $\Gamma^k_{\ell i}$ in the second term:
\begin{align}\label{equation_partialg_Gammag_Gammag_2}
\partial_i g_{\ell j} - \Gamma^k_{\ell i} g_{kj} - \Gamma^k_{i j} g_{\ell k} = 0. 
\end{align}
Similarly, we rename (swap) the dummy variables $j \to i$ and $i \to j$ in Eq. (\ref{equation_partialg_Gammag_Gammag_2}):
\begin{align*}
\partial_j g_{\ell i} - \Gamma^k_{j \ell} g_{ki} - \Gamma^k_{j i} g_{\ell k} = 0. 
\end{align*}
According to Eq. (\ref{equation_zero_torsion_Christoffel_symmetric}), we have $\Gamma^k_{ji} = \Gamma^k_{ij}$, so we can replace $\Gamma^k_{ji}$ with $\Gamma^k_{ij}$ in the second term:
\begin{align}\label{equation_partialg_Gammag_Gammag_3}
\partial_j g_{\ell i} - \Gamma^k_{j \ell} g_{ki} - \Gamma^k_{ij} g_{\ell k} = 0. 
\end{align}
Now, consider the three Eqs. (\ref{equation_partialg_Gammag_Gammag_1}), (\ref{equation_partialg_Gammag_Gammag_2}), and (\ref{equation_partialg_Gammag_Gammag_3}) altogether:
\begin{align*}
&\partial_\ell g_{ij} - \Gamma^k_{\ell i} g_{kj} - \Gamma^k_{\ell j} g_{ik} = 0. \\
&\partial_i g_{\ell j} - \Gamma^k_{\ell i} g_{kj} - \Gamma^k_{i j} g_{\ell k} = 0. \\
&\partial_j g_{\ell i} - \Gamma^k_{j \ell} g_{ki} - \Gamma^k_{ij} g_{\ell k} = 0. 
\end{align*}
We add the second and third equations and subtract the first equation:
\begin{align}
&\partial_i g_{\ell j} - \Gamma^k_{\ell i} g_{kj} - \Gamma^k_{i j} g_{\ell k} + \partial_j g_{\ell i} - \Gamma^k_{j \ell} g_{ki} - \Gamma^k_{ij} g_{\ell k} \nonumber \\
&- \partial_\ell g_{ij} + \Gamma^k_{\ell i} g_{kj} + \Gamma^k_{\ell j} g_{ik} = 0 \overset{(a)}{\implies } \nonumber \\
&\partial_i g_{\ell j} - \cancel{\Gamma^k_{\ell i} g_{kj}} - \Gamma^k_{i j} g_{\ell k} + \partial_j g_{\ell i} - \cancel{\Gamma^k_{j \ell} g_{ki}} - \Gamma^k_{ij} g_{\ell k} \nonumber \\
&- \partial_\ell g_{ij} + \cancel{\Gamma^k_{\ell i} g_{kj}} + \cancel{\Gamma^k_{j \ell} g_{ik}} = 0 \implies \nonumber \\
&\partial_i g_{\ell j} + \partial_j g_{\ell i} - \partial_\ell g_{ij} = 2\Gamma^k_{ij} g_{\ell k}, \implies \nonumber \\
&\Gamma^k_{ij} g_{\ell k} = \frac{1}{2} (\partial_i g_{\ell j} + \partial_j g_{\ell i} - \partial_\ell g_{ij}),
\label{equation_partialg_partialg_partialg_2Gammag}
\end{align}
where $(a)$ is because $\Gamma^k_{\ell j} = \Gamma^k_{j \ell}$ according to Eq. (\ref{equation_zero_torsion_Christoffel_symmetric}). 

Multiplying the sides of Eq. (\ref{equation_partialg_partialg_partialg_2Gammag}) by inverse metric $g^{m \ell}$ gives:
\begin{align*}
&\Gamma^k_{ij} g_{\ell k} g^{m \ell} = \frac{1}{2} (\partial_i g_{\ell j} + \partial_j g_{\ell i} - \partial_\ell g_{ij}) g^{m \ell}, \\
&\overset{(\ref{equation_metric_inverse})}{\implies} \Gamma^k_{ij} \delta_k^m = \frac{1}{2} (\partial_i g_{\ell j} + \partial_j g_{\ell i} - \partial_\ell g_{ij}) g^{m \ell} \\
&\overset{(\ref{equation_index_substitution_delta})}{\implies} \Gamma^m_{ij} = \frac{1}{2} (\partial_i g_{\ell j} + \partial_j g_{\ell i} - \partial_\ell g_{ij}) g^{m \ell}.
\end{align*}
Relabeling the dummy index $m \to k$ gives:
\begin{align*}
\Gamma^k_{ij} = \frac{1}{2} g^{k\ell} \left( \partial_i g_{j\ell} + \partial_j g_{i\ell} - \partial_\ell g_{ij} \right).
\end{align*}
\end{proof}

\begin{proposition}[Christoffel symbols of first kind in terms of metric tensor]
When Levi-Civita connection---which satisfies the metric compatibility and Christoffel symbol symmetry---is used, the Christoffel symbols of first kind can be stated in terms of metric tensor:
\begin{align}\label{equation_Christoffel_firstType_metric}
\boxed{
\Gamma_{ijk} = \frac{1}{2} \left( \partial_i g_{jk} + \partial_j g_{ik} - \partial_k g_{ij} \right).
}
\end{align}
\end{proposition}
\begin{proof}
Consider Eq. (\ref{equation_christoffel_second_kind}), where we rename (swap) the dummy indices $k \to \ell$ and $\ell \to k$:
\begin{align*}
\Gamma^\ell_{ij} = g^{\ell k} \Gamma_{ijk}.
\end{align*}
Also, consider the Eq. (\ref{equation_Christoffel_secondType_metric}), where we rename (swap) the dummy indices $k \to \ell$ and $\ell \to k$:
\begin{align*}
\Gamma^\ell_{ij} = \frac{1}{2} g^{\ell k} \left( \partial_i g_{jk} + \partial_j g_{ik} - \partial_k g_{ij} \right).
\end{align*}
Comparing these two equations gives:
\begin{align*}
g^{\ell k} \Gamma_{ijk} = \frac{1}{2} g^{\ell k} \left( \partial_i g_{jk} + \partial_j g_{ik} - \partial_k g_{ij} \right).
\end{align*}
As this holds for all $g^{\ell k}$, we have:
\begin{align*}
\Gamma_{ijk} = \frac{1}{2} \left( \partial_i g_{jk} + \partial_j g_{ik} - \partial_k g_{ij} \right).
\end{align*}
\end{proof}

\section{Ricci Calculus Notation: Comma and Semicolon Derivatives}\label{section_ricci_calculus_notation}

In the study of Riemannian geometry, it is common to use a shorthand notation for derivatives. This notation, often referred to as \textit{Ricci calculus} (also called the \textit{absolute differential calculus}), distinguishes between the standard partial derivative and the coordinate-invariant covariant derivative.
Ricci calculus was developed by \textit{Gregorio Ricci-Curbastro} and later popularized with his student \textit{Tullio Levi-Civita}.

\subsection{Comma and Semicolon Notations}

\begin{definition}[Comma and Semicolon Notations]\label{definition_comma_semicolon_notations}
Let $V^j$ be the contravariant components of vector $\b{V}$ and $V_j$ be the covariant components of vector (or equivalently, let Let $V^j$ be the components of a vector field and $V_j$ be the components of a one-form). We define:
\begin{enumerate}
    \item \textbf{Comma notation (partial derivative):} A comma followed by an index denotes the partial derivative with respect to that coordinate:
    \begin{equation}
    \boxed{
    \begin{aligned}
    & V^j_{,i} \equiv \partial_i V^j = \frac{\partial V^j}{\partial x^i}  \\
    & V_{j,i} \equiv \partial_i V_j = \frac{\partial V_j}{\partial x^i}
    \end{aligned}
    }
    \end{equation}
    \item \textbf{Semicolon notation (covariant derivative):} A semicolon followed by an index denotes the covariant derivative with respect to the connection $\nabla$:
    \begin{equation}
    \boxed{
    \begin{aligned}
    &V^j_{;i} \equiv \nabla_i V^j = (\nabla_i \b{V})^j \\
    &V_{j;i} \equiv \nabla_i V_j = (\nabla_i \b{V})_j
    \end{aligned}
    }
    \end{equation}
    where the three notations are:
    \begin{itemize}
        \item $V^j_{;i}$ (\textbf{Ricci calculus}): Extremely compact for long derivations, especially when taking second derivatives (like the Hessian).
        \item $\nabla_i V^j$ (\textbf{Component-wise notation}): The standard "workhorse" notation you have used for most of your propositions.
        \item $(\nabla_i \b{V})^j$ (\textbf{Operator notation}): The most rigorous form, as it explicitly shows that $\nabla_i$ is an operator acting on the vector field $\b{V}$, and we are then looking at the $j$-th component of the result.
    \end{itemize}
\end{enumerate}
The rule for this notation is:
\begin{itemize}
    \item Indices before the delimiter (comma or semicolon) identify the components of the tensor field.
    \item Indices after the delimiter indicate the coordinates with respect to which you are differentiating.
    \item The order after the delimiter matters. $V_{;ik}$ means we first take the covariant derivative with respect to $x^i$, and then we take the covariant derivative of that resulting tensor with respect to $x^k$.
\end{itemize}
Some examples for the notation are:
\begin{align*}
&V^{jl}_{pq,ik} \equiv \partial_k (\partial_i V^{jl}_{pq}) \\
&V^{jl}_{pq;ik} \equiv \nabla_k (\nabla_i V^{jl}_{pq})
\end{align*}
Note that some texts in the literature denote $V^{jl}_{pq,ik}$ as $V^{jl}_{pq,i,k}$ and denote $V^{jl}_{pq;ik}$ as $V^{jl}_{pq;i;k}$.
\end{definition}

\subsection{Examples of Usage of Comma and Semicolon Notations}

Using the Ricci calculus notation, some of the formulas for the covariant derivatives that we derived previously can be written more compactly:
\begin{itemize}
    \item Covariant derivative of contravariant components:
    \begin{align}
    \boxed{
    V^j_{;i} = V^j_{,i} + \Gamma^j_{ik} V^k
    }
    \end{align}
    
    \item Covariant derivative of covariant components:
    \begin{align}
    \boxed{
    V_{j;i} = V_{j,i} - \Gamma^k_{ij} V_k
    }
    \end{align}
    
    \item Metric compatibility:
    The condition for a Levi-Civita connection ($\nabla_k g_{ij} = 0$) is expressed elegantly as:
    \begin{align}
    \boxed{
    g_{ij;k} = 0
    }
    \end{align}

    \item Equations (\ref{equation_second_derivative_contravariant}) and (\ref{equation_second_derivative_covariant}) for second-order covariant derivatives can be stated by the Ricci calculus notation:
    \begin{equation}
    \boxed{
    \begin{aligned}
    &V^j_{;i\ell} = \partial_\ell V^j_{;i} + \Gamma^j_{\ell k} V^k_{;i} - \Gamma^k_{\ell i} V^j_{;k} \\
    &V_{j;i\ell} = \partial_\ell V_{j;i} - \Gamma^k_{\ell i} V_{j;k} - \Gamma^k_{\ell j} V_{k;i}
    \end{aligned}
    }
    \end{equation}
\end{itemize}

\section{Riemannian Curvature}\label{section_riemannian_curvature_quantities}

Riemannian curvature refers to calculating the curvature of manifold at different points of manifold. There are different types of Riemannian curvature.
In this section, we introduce different types of curvature in an intuitive order:
\begin{enumerate}
\item Riemann Curvature Tensor (the ``parent" object): It contains all information about the curvature (e.g., $20$ independent components in a four-dimensional manifold).
\item Sectional Curvature (the ``geometric essence"): It contains the same amount of information as the Riemann tensor (you can reconstruct Riemann curvature from all sectional curvatures), but it is expressed as a real-valued function on $2$-planes.
\item Ricci Curvature (the ``trace/average"): A compression of the Riemann tensor (e.g., $10$ independent components in a four-dimensional manifold). It loses some information but is easier to work with in optimization (e.g., Ricci flow).
\item Scalar Curvature (the ``total trace"): It is a single number at each point of manifold.
\end{enumerate}

\subsection{Riemann Curvature Tensor}


\subsubsection{Intuition of Riemann Curvature}\label{section_idea_Riemannian_curvature}

\begin{figure}[!h]
\centering
\includegraphics[width=2.5in]{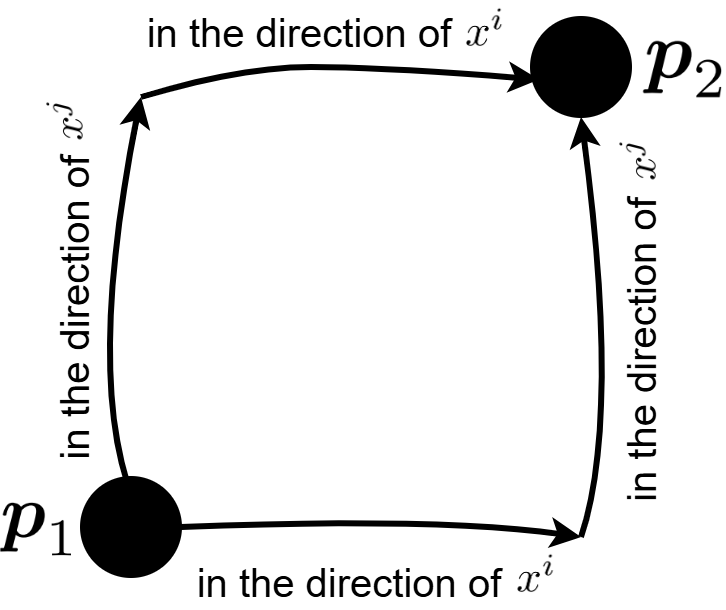}
\caption{Going from point $\b{p}_1$ to point $\b{p}_2$ in two different infinitesimal paths: (1) first in the direction of $dx^i$ and then in the direction of $dx^j$, or (2) first in the direction of $dx^j$ and then in the direction of $dx^i$.}
\label{figure_Riemann_curvature}
\end{figure}

Consider Fig. \ref{figure_Riemann_curvature} where we can go from point $\b{p}_1$ to $\b{p}_2$ through two infinitesimal paths. On the one hand, we can first go from $\b{p}_1$ along coordinate $x^i$ and then go along coordinate $x^j$ to the point $\b{p}_2$ with infinitesimal steps. On the other hand, we can first go from $\b{p}_1$ along coordinate $x^j$ and then go along coordinate $x^i$ to the point $\b{p}_2$ with infinitesimal steps.

Going from $\b{p}_1$ along coordinate $x^i$ and then going along coordinate $x^j$ to the point $\b{p}_2$ with infinitesimal steps can be modeled as:
\begin{itemize}
\item $\nabla_i V^k$ denotes considering the $k$-th component of vector $\b{V}$ going along coordinate $x^i$.
\item $\nabla_j (\nabla_i V^k)$ denotes considering the previous $\nabla_i V^k$ for going along coordinate $x^j$.
\end{itemize}

Going from $\b{p}_1$ along coordinate $x^j$ and then going along coordinate $x^i$ to the point $\b{p}_2$ with infinitesimal steps can be modeled as:
\begin{itemize}
\item $\nabla_j V^k$ denotes considering the $k$-th component of vector $\b{V}$ going along coordinate $x^j$.
\item $\nabla_i (\nabla_j V^k)$ denotes considering the previous $\nabla_j V^k$ for going along coordinate $x^i$.
\end{itemize}

\begin{figure}[!h]
\centering
\includegraphics[width=3.2in]{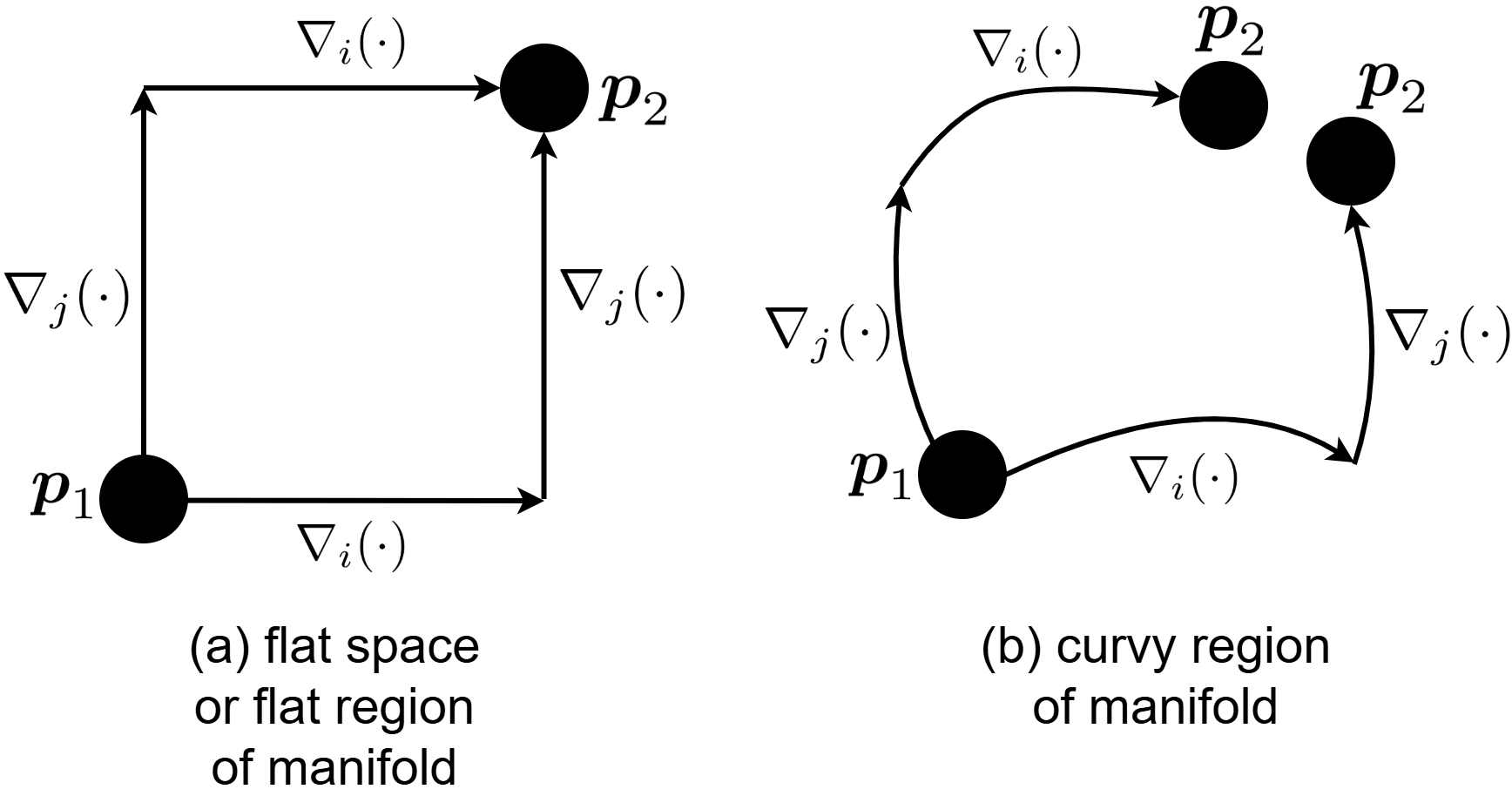}
\caption{Going from point $\b{p}_1$ to point $\b{p}_2$ in two different infinitesimal paths: (a) In a flat space or a flat region of manifold, we end up with the same point $\b{p}_2$ but (b) in a curvy region of manifold, we will end up in different points $\b{p}_2$ which is because of effect of curvature.}
\label{figure_Riemann_curvature2}
\end{figure}

On a flat manifold, without curvature, going along these two paths end up to the same point $\b{p}_2$. However, if the manifold is not flat---meaning that it has curvature---going along these two paths do not end up to the same point. This behavior is illustrated in Fig. \ref{figure_Riemann_curvature2}.
This is the idea of curvature proposed by Riemann. That is why this measure of curvature is also called \textit{Riemannian curvature}.

According to the discussion above, the curvature can be modeled as the difference between the covariant derivatives $\nabla_j (\nabla_i V^k)$ and $\nabla_i (\nabla_j V^k)$:
\begin{align}\label{equation_curvature_idea}
\boxed{
\text{Curvature} \propto \nabla_j \nabla_i V^k - \nabla_i \nabla_j V^k.
}
\end{align}

\subsubsection{Definition of Riemann Curvature}

\begin{definition}[Riemann curvature tensor --- coordinate-free definition \cite{Riemann1854uber}]
Let $(\mathcal{M},\nabla)$ be a smooth manifold equipped with an affine connection $\nabla$.
The \textbf{Riemann curvature tensor} is the $(1,3)$ tensor $R$ defined by:
\begin{equation}\label{equation_Riemann_curvature_tensor_coordinate_free}
\boxed{
\begin{aligned}
R(\b{X},\b{Y})\b{Z} &:= [\nabla_{\b{X}}, \nabla_{\b{Y}}]\b{Z} - \nabla_{[\b{X},\b{Y}]}\b{Z} \\
&= \nabla_{\b{X}} \nabla_{\b{Y}} \b{Z} \;-\; \nabla_{\b{Y}} \nabla_{\b{X}} \b{Z} \;-\; \nabla_{[\b{X},\b{Y}]} \b{Z},
\end{aligned}
}
\end{equation}
for all vector fields $\b{X},\b{Y},\b{Z} \in \mathfrak{X}(\mathcal{M})$.
\end{definition}

Equation (\ref{equation_Riemann_curvature_tensor_coordinate_free}) determines how much the result of ``moving $\b{Z}$" depends on the order of moving along $\b{X}$ and $\b{Y}$. 
As discussed in Section \ref{section_idea_Riemannian_curvature}, $\nabla_{\b{X}} \nabla_{\b{Y}} \b{Z} \;-\; \nabla_{\b{Y}} \nabla_{\b{X}} \b{Z}$ measures this difference. 
However, what is the last term  $-\nabla_{[\b{X},\b{Y}]} \b{Z}$ for? 
We know that, in general, vector fields do not commute, i.e., $[\b{X}, \b{Y}] \neq \b{0}$. So, we need to subtract the ``fake effect"; this removes the part coming from the non-commutativity of directions, and not the curvature. 

The \textit{Riemann curvature tensor} was proposed by \textit{Bernhard Riemann} \cite{Riemann1854uber}. He presented it in a lecture in 1854 and it was published posthumously in 1868. 
Riemann introduced curvature via a quadratic form expansion of the metric.

\subsubsection{Derivation of Riemann Curvature in Coordinates}

\begin{proposition}[Riemann curvature tensor in coordinates]
Consider an $n$-dimensional manifold. 
The Riemann curvature tensor, denoted by $R^\ell_{ijk}$, is a $(1,3)$ tensor which calculates the curvature of manifold as:
\begin{equation}\label{equation_Riemann_curvature_coordinates}
\boxed{
R^\ell_{ijk} := \partial_j \Gamma^\ell_{ki} - \partial_k \Gamma^\ell_{ji} + \Gamma^\ell_{jp} \Gamma^p_{ki} - \Gamma^\ell_{kp} \Gamma^p_{ji},
}
\end{equation}
where $i,j,k,\ell \in \{1, \dots, n\}$.
This expression allows for the direct numerical implementation of Riemann curvature given the Christoffel symbols and their partial derivatives.
\end{proposition}
\begin{proof}
We work in a coordinate basis $\{\b{e}_i\}_{i=1}^n$, where 
$[\b{e}_i,\b{e}_j] = 0, \forall i,j$ according to Eq. (\ref{equation_Lie_bracket_coordinate_basis}). Hence, the last term is zero in coordinate bases:
\begin{align*}
\nabla_{[\b{X},\b{Y}]} \b{Z} = \nabla_{[\b{e}_i,\b{e}_j]} \b{Z} = \nabla_{0} \b{Z} = \b{0}.
\end{align*}
Therefore, the curvature operator reduces to:
\begin{align*}
R(\b{e}_j,\b{e}_i)V^k
=
\nabla_j \nabla_i V^k - \nabla_i \nabla_j V^k.
\end{align*}

According to Eq. (\ref{equation_curvature_idea}), the curvature is $\nabla_j \nabla_i V^k - \nabla_i \nabla_j V^k$. 
Recall that $\nabla_i V^k$ is a $(1,1)$ tensor. According to Eq. (\ref{equation_covariant_derivative_tensor_1_1}), applying the covariant derivative $\nabla_j$ to this tensor yields:
\begin{align}
\nabla_j (\nabla_i V^k) \overset{(\ref{equation_covariant_derivative_tensor_1_1})}{=} \partial_j (\nabla_i V^k) - \Gamma^p_{ji} (\nabla_p V^k) + \Gamma^k_{jp} (\nabla_i V^p). \label{equation_second_cov_deriv_vec}
\end{align}
According to Eq. (\ref{equation_covariant_derivative_contravariant}), we have $\nabla_i V^k = \partial_i V^k + \Gamma^k_{i\ell} V^\ell$. 
Substituting this into Eq. \eqref{equation_second_cov_deriv_vec}, we obtain:
\begin{align*}
\nabla_j \nabla_i V^k &= \partial_j (\partial_i V^k + \Gamma^k_{i\ell} V^\ell) - \Gamma^p_{ji} (\partial_p V^k + \Gamma^k_{p\ell} V^\ell) \\
&~~~~~~ + \Gamma^k_{jp} (\partial_i V^p + \Gamma^p_{i\ell} V^\ell) \\
&= \partial_j \partial_i V^k + (\partial_j \Gamma^k_{i\ell}) V^\ell + \Gamma^k_{i\ell} \partial_j V^\ell - \Gamma^p_{ji} \nabla_p V^k \\
&~~~~~~ + \Gamma^k_{jp} \partial_i V^p + \Gamma^k_{jp} \Gamma^p_{i\ell} V^\ell.
\end{align*}
To find $\nabla_j \nabla_i V^k - \nabla_i \nabla_j V^k$, we subtract the same expression with indices $i$ and $j$ swapped. Assuming a torsion-free connection ($\Gamma^p_{ji} = \Gamma^p_{ij}$), the terms $\Gamma^p_{ji} \nabla_p V^k$ and $\Gamma^p_{ij} \nabla_p V^k$ cancel out, and the second-order partial derivatives $\partial_j \partial_i V^k$ and $\partial_i \partial_j V^k$ cancel out. 
In other words, we have:
\begin{align*}
&\nabla_j \nabla_i V^k - \nabla_i \nabla_j V^k \\
&= \cancel{\partial_j \partial_i V^k} + (\partial_j \Gamma^k_{i\ell}) V^\ell + \Gamma^k_{i\ell} \partial_j V^\ell - \cancel{\Gamma^p_{ji} \nabla_p V^k} \\
&~~~~~ + \Gamma^k_{jp} \partial_i V^p + \Gamma^k_{jp} \Gamma^p_{i\ell} V^\ell \\
&~~~~~ - \cancel{\partial_i \partial_j V^k} - (\partial_i \Gamma^k_{j\ell}) V^\ell - \Gamma^k_{j\ell} \partial_i V^\ell + \cancel{\Gamma^p_{ij} \nabla_p V^k} \\
&~~~~~ - \Gamma^k_{ip} \partial_j V^p - \Gamma^k_{ip} \Gamma^p_{j\ell} V^\ell \\
&= (\partial_j \Gamma^k_{i\ell}) V^\ell + \Gamma^k_{i\ell} \partial_j V^\ell + \Gamma^k_{jp} \partial_i V^p + \Gamma^k_{jp} \Gamma^p_{i\ell} V^\ell \\
&~~~~~ - (\partial_i \Gamma^k_{j\ell}) V^\ell - \Gamma^k_{j\ell} \partial_i V^\ell - \Gamma^k_{ip} \partial_j V^p - \Gamma^k_{ip} \Gamma^p_{j\ell} V^\ell \\
&\overset{(a)}{=} (\partial_j \Gamma^k_{i\ell}) V^\ell + \cancel{\Gamma^k_{i\ell} \partial_j V^\ell} + \cancel{\Gamma^k_{j\ell} \partial_i V^\ell} + \Gamma^k_{jp} \Gamma^p_{i\ell} V^\ell \\
&~~~~~ - (\partial_i \Gamma^k_{j\ell}) V^\ell - \cancel{\Gamma^k_{j\ell} \partial_i V^\ell} - \cancel{\Gamma^k_{i\ell} \partial_j V^\ell} - \Gamma^k_{ip} \Gamma^p_{j\ell} V^\ell \\
&= (\partial_j \Gamma^k_{i\ell} - \partial_i \Gamma^k_{j\ell} + \Gamma^k_{jp} \Gamma^p_{i\ell} - \Gamma^k_{ip} \Gamma^p_{j\ell}) V^\ell,
\end{align*}
where $(a)$ is because of relabeling dummy indices $p \to \ell$ in the third and seventh terms. 

Factoring out the $V^k$ on the left hand side and the $V^\ell$ on the right hand side give:
\begin{align*}
(\nabla_j \nabla_i - &\nabla_i \nabla_j) V^k = \\
&\left( \partial_i \Gamma^\ell_{jk} - \partial_j \Gamma^\ell_{ik} + \Gamma^p_{jk} \Gamma^\ell_{ip} - \Gamma^p_{ik} \Gamma^\ell_{jp} \right) V^\ell.
\end{align*}

We define the \textit{Riemann curvature tensor} $R^k_{\ell ji}$ such that:
\begin{align*}
\nabla_j \nabla_i V^k - \nabla_i \nabla_j V^k = R^k_{\ell ji} V^\ell.
\end{align*}
Comparing the components, the explicit form of the tensor is:
\begin{equation*}
R^k_{\ell ji} = \partial_j \Gamma^k_{i\ell} - \partial_i \Gamma^k_{j\ell} + \Gamma^k_{jp} \Gamma^p_{i\ell} - \Gamma^k_{ip} \Gamma^p_{j\ell}.
\end{equation*}
Relabeling the dummy indices $k \to \ell$, $\ell \to i$, $j \to j$, and $i \to k$ gives:
\begin{equation*}
R^\ell_{ijk} = \partial_j \Gamma^\ell_{ki} - \partial_k \Gamma^\ell_{ji} + \Gamma^\ell_{jp} \Gamma^p_{ki} - \Gamma^\ell_{kp} \Gamma^p_{ji}.
\end{equation*}
\end{proof}

\begin{definition}[The $(0, 4)$ fully covariant Riemann curvature tensor]
By using the metric tensor $g$, we can lower the upper index of the $(1,3)$ Riemann curvature tensor to make it a $(0,4)$ tensor. 
Its coordinate-free definition is:
\begin{align}
R(\b{X}, \b{Y}, \b{Z}, \b{W}) = g(R(\b{X}, \b{Y})\b{Z}, \b{W}).
\end{align}
In components, it is obtained as (see Eq. (\ref{equation_index_lower})):
\begin{align}
\boxed{
R_{mijk} = g_{m\ell} R^\ell_{ijk}.
}
\end{align}
This tensor is called the \textbf{$(0, 4)$ fully covariant Riemann curvature tensor}.
\end{definition}

\subsubsection{Properties of Riemann Curvature Tensor}

\begin{lemma}[Anti-symmetry of $(1,3)$ Riemann curvature]\label{lemma_riemann_curvature_symmetry_1_3}
The Riemann curvature tensor of type $(1,3)$, denoted by $R^\ell{}_{ijk}$, is anti-symmetric in its last two covariant indices. That is, the components satisfy the identity:
\begin{align}
\boxed{
R^\ell_{ijk} = -R^\ell_{ikj}.
}
\end{align}
\end{lemma}
\begin{proof}
According to Eq. (\ref{equation_Riemann_curvature_coordinates}):
\begin{equation*}
R^\ell_{ijk} = \partial_j \Gamma^\ell_{ki} - \partial_k \Gamma^\ell_{ji} + \Gamma^\ell_{js} \Gamma^s_{ki} - \Gamma^\ell_{ks} \Gamma^s_{ji}
\end{equation*}
To verify the anti-symmetry property, we compute $R^\ell{}_{ikj}$ by interchanging the indices $j$ and $k$. By direct substitution, we have:
\begin{align*}
R^\ell_{ikj} &= \partial_k \Gamma^\ell_{ji} - \partial_j \Gamma^\ell_{ki} + \Gamma^\ell_{ks} \Gamma^s_{ji} - \Gamma^\ell_{js} \Gamma^s_{ki} \\
&= -\partial_j \Gamma^\ell_{ki} + \partial_k \Gamma^\ell_{ji} - \Gamma^\ell_{js} \Gamma^s_{ki} + \Gamma^\ell_{ks} \Gamma^s_{ji} \\
&= -\left( \partial_j \Gamma^\ell_{ki} - \partial_k \Gamma^\ell_{ji} + \Gamma^\ell_{js} \Gamma^s_{ki} - \Gamma^\ell_{ks} \Gamma^s_{ji} \right).
\end{align*}
Recognizing the term inside the parentheses as the original definition of $R^\ell_{ijk}$, we obtain:
\[
R^\ell_{ikj} = -R^\ell_{ijk},
\]
which proves that the tensor is identically anti-symmetric in its last two indices.
\end{proof}

\begin{lemma}[Symmetry and anti-symmetry in $(0,4)$ Riemann curvature tensor]\label{lemma_riemann_curvature_symmetry_0_4}
Let $R_{mijk}$ denote the components of the $(0,4)$ Riemann curvature tensor. Then the following algebraic identities hold:
\begin{align}
&\text{Last-pair anti-symmetry: } \boxed{R_{mijk} = -R_{mikj}}, \\
&\text{First-pair anti-symmetry: } \boxed{R_{mijk} = -R_{imjk}}, \\
&\text{Interchange symmetry: } \boxed{R_{mijk} + R_{mjki} + R_{mkij} = 0}, \\
&\text{First Bianchi Identity: } \boxed{R_{mijk} = R_{jkmi}}. 
\end{align}
\end{lemma}

\subsection{Sectional Curvature}\label{section_sectional_curvature}

\subsubsection{Intuition of Sectional Curvature}

Imagine you are standing at a point $\b{p}$ on an $n$-dimensional manifold $\mathcal{M}$. To understand the curviness at $\b{p}$, you pick two linearly independent tangent vectors, $\b{X}$ and $\b{Y}$. These two vectors span a 2D plane (called ``section") in the tangent space $T_{\b{p}}\mathcal{M}$.
The \textit{sectional curvature} $K(\b{X}, \b{Y})$ is defined as the Gaussian curvature of the small two-dimensional surface formed by all the geodesics starting at $\b{p}$ whose tangent vectors lie in that plane.

\subsubsection{Definition of Sectional Curvature}

\begin{definition}[Sectional curvature]
Let $(\mathcal{M},g)$ be a Riemannian manifold, $\b{p} \in \mathcal{M}$, and $\b{X},\b{Y} \in T_{\b{p}} \mathcal{M}$ two linearly independent tangent vectors. 
Consider the two-dimensional plane spanned by $\b{X}$ and $\b{Y}$, i.e., $\mathrm{span}\{\b{X},\b{Y}\}$.
The \textbf{sectional curvature} of the plane spanned by $\b{X}$ and $\b{Y}$ is defined as:
\begin{equation}
\boxed{
K(\b{X},\b{Y}) := \frac{\langle R(\b{X},\b{Y})\b{Y}, \b{X} \rangle}{\langle \b{X},\b{X} \rangle \langle \b{Y},\b{Y} \rangle - \langle \b{X},\b{Y} \rangle^2}.
}
\end{equation}
\end{definition}

\subsubsection{Derivation of Sectional Curvature in Coordinates}

\begin{proposition}[Sectional curvature in coordinates]
Let $\b{X} = X^i \partial_i$ and $\b{Y} = Y^i \partial_i$ be tangent vectors in local coordinates with metric $g_{ij}$.
Then the sectional curvature of the plane spanned by $X$ and $Y$ is
\begin{equation}
\boxed{
K(\b{X},\b{Y}) = \frac{X^\ell Y^i Y^k X^j g_{\ell m} R^m_{ijk}}{(X^i X^j g_{ij})(Y^k Y^\ell g_{k\ell}) - (X^i Y^j g_{ij})^2},
}
\end{equation}
where $R^m_{ijk} := \partial_j \Gamma^m_{ki} - \partial_k \Gamma^m_{ji} + \Gamma^m_{jp} \Gamma^p_{ki} - \Gamma^m_{kp} \Gamma^p_{ji}$ is the Riemann curvature tensor.
\end{proposition}
\begin{proof}
The sectional curvature $K(\b{X}, \b{Y})$ of the plane spanned by $\b{X}, \b{Y} \in T_{\b{p}} \mathcal{M}$ is defined by the formula:
\begin{equation}\label{equation_sectional_curvature_in_proof}
K(\b{X}, \b{Y}) = \frac{\langle R(\b{X}, \b{Y})\b{Y}, \b{X} \rangle}{\langle \b{X}, \b{X} \rangle \langle \b{Y}, \b{Y} \rangle - \langle \b{X}, \b{Y} \rangle^2}.
\end{equation}
We evaluate this expression strictly in local coordinates where $\b{X} = X^i \partial_i$ and $\b{Y} = Y^j \partial_j$.

-- Derivation of the \underline{denominator}:
Using the linearity of the metric inner product, we have:
\begin{align*}
\langle \b{U}, \b{V} \rangle \overset{(\ref{equation_contravariant_components_2})}{=} \langle U^i \partial_i, V^j \partial_j \rangle \overset{(a)}{=} U^i V^j \langle \partial_i, \partial_j \rangle \overset{(\ref{equation_g_components})}{=} g_{ij} U^i V^j,
\end{align*}
where $(a)$ is because $U^i$ and $V^j$ are scalar components coming out of the inner product. Therefore, we can expand the terms in the denominator of Eq. (\ref{equation_sectional_curvature_in_proof}):
\begin{align*}
&\langle \b{X}, \b{X} \rangle = g_{ij} X^i X^j, \\
&\langle \b{Y}, \b{Y} \rangle = g_{k\ell} Y^k Y^\ell, \\
&\langle \b{X}, \b{Y} \rangle = g_{ij} X^i Y^j.
\end{align*}
Substituting these into the denominator yields:
\begin{align*}
&\langle \b{X},\b{X} \rangle \langle \b{Y},\b{Y} \rangle - \langle \b{X},\b{Y} \rangle^2 \\
&~~~~~~~ = (g_{ij} X^i X^j)(g_{k\ell} Y^k Y^\ell) - (g_{ij} X^i Y^j)^2.
\end{align*}

-- Derivation of the \underline{numerator}:
The numerator involves the Riemann curvature operator $R(\b{X}, \b{Y})\b{Y}$. By the multilinearity of the Riemann tensor, we write:
\begin{align*}
R(\b{X}, \b{Y})\b{Y} &\overset{(\ref{equation_contravariant_components_2})}{=} R(X^j \partial_j, Y^k \partial_k)(Y^i \partial_i) \\
&\overset{(a)}{=} X^j Y^k Y^i R(\partial_j, \partial_k) \partial_i,
\end{align*}
where $(a)$ is because multilinearity of the Riemann tensor.

Using the component definition $R(\partial_j, \partial_k) \partial_i = R^m_{ijk} \partial_m$, we obtain:
\begin{align*}
R(\b{X}, \b{Y})\b{Y} = (X^j Y^k Y^i R^m_{ijk}) \partial_m.
\end{align*}
Now, taking the inner product with $X = X^\ell \partial_\ell$ gives:
\begin{align*}
\langle R(\b{X}, \b{Y})\b{Y}, \b{X} \rangle &= \langle (X^j Y^k Y^i R^m_{ijk}) \partial_m, X^\ell \partial_\ell \rangle \\
&\overset{(a)}{=} X^\ell X^j Y^k Y^i R^m_{ijk} \langle \partial_m, \partial_\ell \rangle \\
&\overset{(\ref{equation_g_components})}{=} X^\ell X^j Y^k Y^i R^m_{ijk} g_{m\ell},
\end{align*}
where $(a)$ is because of bringing the scalar components out of the inner product. 

Rearranging (relabeling) the dummy indices to match the proposition, we have:
\begin{align*}
\langle R(\b{X}, \b{Y})\b{Y}, \b{X} \rangle = X^\ell Y^i Y^k X^j g_{\ell m} R^m_{ijk}.
\end{align*}
\end{proof}

\subsection{Ricci Curvature}

\subsubsection{Intuition of Ricci Curvature}

While the Riemann curvature tensor provides a complete description of the manifold's curvature at a point, the \textit{Ricci curvature}\footnote{Ricci is pronounced as REE-chee
or \textipa{/ritSi/} in simplified international phonetic alphabet.} offers a ``summarized" view of how the manifold curves.

Ricci curvature is named after \textit{Gregorio Ricci-Curbastro} (from Italy), who did important work on developing covariant derivative, in 1887 \cite{ricci1887sulla}.

Intuitively, Ricci curvature describes how the volume of an infinitesimal cone of geodesics deviates from that in Euclidean space. 
Intrepretation of the sign of Ricci curvature is as follows:
\begin{itemize}
\item If the Ricci curvature is positive in a certain direction, geodesics tend to converge, and the volume of a small region of the manifold grows more slowly than it would in flat space. 
\item However, if the Ricci curvature is negative in a certain direction, geodesics tend to diverge, and the volume of a small region of the manifold grows faster than it would in flat space.
\end{itemize}

\subsubsection{Definition of Ricci Curvature}

\begin{definition}[Ricci curvature]

Let $R$ be the Riemann curvature tensor of type $(1, 3)$. The \textbf{Ricci curvature tensor} $\mathrm{Ric}$ (often denoted by $R_{ij}$ in coordinate notation), is a $(0, 2)$ tensor, and is defined at each point $\b{p} \in \mathcal{M}$ as follows:
\begin{equation}\label{equation_Ricci_curvature}
\boxed{
R_{ij} := R^k_{ikj}, 
}
\end{equation}
where $R^k_{ikj}$ are the components of the Riemann curvature tensor, and note that $R^k_{ikj}$ is summed over $k$ according to Einstein summation convention. 
The Ricci curvature tensor is also called the \textbf{Ricci tensor} in short. 

In coordinate-free notation, for vector fields $X$ and $Y$, the Ricci tensor is defined as the trace of the Riemann curvature tensor over its first and third indices:
\begin{equation}
\boxed{
\mathrm{Ric}(\b{X}, \b{Y}) := \text{tr}(\b{Z} \mapsto R(\b{Z}, \b{X})\b{Y}),
}
\end{equation}
where the trace is taken over the map that sends a vector $\b{Z}$ to the Riemann curvature operator acting on $\b{Y}$ and $\b{Z}$.
\end{definition}

\subsubsection{Derivation of Ricci Curvature in Coordinates}


\begin{proposition}[Ricci curvature tensor in coordinates]
The components of the Ricci curvature tensor are given by:
\begin{equation}\label{equation_Ricci_curvature_coordinates}
\boxed{
R_{ij} = \partial_k \Gamma^k_{ji} - \partial_j \Gamma^k_{ki} + \Gamma^k_{kp} \Gamma^p_{ji} - \Gamma^k_{jp} \Gamma^p_{ki},
}
\end{equation}
where $\Gamma^k_{ij}$ are the Christoffel symbols of the second kind.
This expression allows for the direct numerical implementation of Ricci curvature given the Christoffel symbols and their partial derivatives.
\end{proposition}
\begin{proof}
Recall the coordinate expression for the Riemann curvature tensor, i.e., Eq. (\ref{equation_Riemann_curvature_coordinates}):
\begin{equation}\label{equation_Riemann_curvature_coordinates_in_proof_ricci}
    R^\ell_{ijk} = \partial_j \Gamma^\ell_{ki} - \partial_k \Gamma^\ell_{ji} + \Gamma^\ell_{jp} \Gamma^p_{ki} - \Gamma^\ell_{kp} \Gamma^p_{ji}.
\end{equation}
To obtain the Ricci tensor $R_{ik}$, we contract the upper index $\ell$ with the second lower index (the index in the $j$-position):
\begin{equation*}
    R_{ik} = R^m_{imk}.
\end{equation*}
Substituting $m$ for $\ell$ and $j$ in Eq. (\ref{equation_Riemann_curvature_coordinates_in_proof_ricci}) gives:
\begin{equation*}
    R_{ik} = R^m_{imk} = \partial_m \Gamma^m_{ki} - \partial_k \Gamma^m_{mi} + \Gamma^m_{mp} \Gamma^p_{ki} - \Gamma^m_{kp} \Gamma^p_{mi}.
\end{equation*}
Relabeling the dummy indices $k \to j$ and $m \to k$, we have:
\begin{equation*}
    R_{ij} = \partial_k \Gamma^k_{ji} - \partial_j \Gamma^k_{ki} + \Gamma^k_{kp} \Gamma^p_{ji} - \Gamma^k_{jp} \Gamma^p_{ki}
\end{equation*}
\end{proof}

\subsection{Scalar Curvature (Ricci Scalar)}

\subsubsection{Intuition of Scalar Curvature}
While the Riemann curvature tensor and the Ricci curvature provide multi-dimensional information about the manifold's geometry, the \textit{scalar curvature} (also known as the \textit{Ricci scalar}) provides a single real number at each point $\b{p} \in \mathcal{M}$. Intuitively, it represents the simplest possible measure of intrinsic curvature. Geometrically, for a small $n$-dimensional ball in the manifold, the scalar curvature describes how the volume of that ball deviates from the volume of a standard ball in Euclidean space of the same radius. 

\subsubsection{Definition of Scalar Curvature}

\begin{definition}[Scalar curvature]
Let $(\mathcal{M}, g)$ be a Riemannian manifold, $R_{ij}$ be the Ricci curvature tensor, and $R^\ell_{ijk}$ be the Riemann curvature tensor. The \textbf{scalar curvature} $R$ (sometimes denoted by $S$) is the trace of the Ricci tensor with respect to the metric:
\begin{equation}\label{equation_scalar_curvature_def}
\boxed{
R := \text{tr}_g(\mathrm{Ric}) = g^{ij}R_{ij} \overset{(\ref{equation_Ricci_curvature})}{=} g^{ij} R^k_{ikj}, 
}
\end{equation}
where $g^{ij}$ are the components of the inverse metric tensor.
The scalar curvature is a $(0,0)$ tensor, i.e., a scalar number.
\end{definition}

\subsubsection{Derivation of Scalar Curvature in Coordinates}

\begin{definition}[Scalar curvature in coordinates]
In a local coordinate system $\{x^1, \dots, x^n\}$, the \textbf{scalar curvature} is calculated by contracting the indices of the Ricci tensor. Its coordinate-based expression is:
\begin{equation}\label{equation_scalar_curvature_christoffel}
\boxed{
R = g^{ij} \left( \partial_k \Gamma^k_{ji} - \partial_j \Gamma^k_{ki} + \Gamma^k_{kp} \Gamma^p_{ji} - \Gamma^k_{jp} \Gamma^p_{ki} \right).
}
\end{equation}
\end{definition}
\begin{proof}
The derivation follows directly from the definition of a trace in Riemannian geometry. Since the Ricci tensor is a $(0, 2)$ tensor, we must use the inverse metric $g^{ij}$ to raise an index before we can perform the final contraction (summation) to obtain a scalar:
\begin{enumerate}
    \item Start with the Riemann curvature tensor $R^l_{ijk}$.
    \item Contract the first and third indices to obtain the Ricci tensor: $R_{ik} = R^j_{ijk}$.
    \item Contract the remaining two indices using the inverse metric: $R = g^{ik}R_{ik}$.
\end{enumerate}
In terms of the Christoffel symbols, this can be expanded as:
\begin{align*}
R &\overset{(\ref{equation_scalar_curvature_def})}{=} g^{ij}R_{ij} \\
&\overset{(\ref{equation_Ricci_curvature_coordinates})}{=} g^{ij} \left( \partial_k \Gamma^k_{ji} - \partial_j \Gamma^k_{ki} + \Gamma^k_{kp} \Gamma^p_{ji} - \Gamma^k_{jp} \Gamma^p_{ki} \right). \label{eq:scalar_curvature_christoffel}
\end{align*}
\end{proof}

The scalar curvature is essential for various applications in Riemannian optimization, such as regularizing objective functions based on the ``total'' curvature of the search space.

\subsection{Gaussian Curvature in Two Dimensions}

\subsubsection{Intuition of Gaussian Curvature}

For a two-dimensional Riemannian manifold (a surface), the \textit{Gaussian curvature} $K$ is the most fundamental measure of intrinsic curvature. Historically, \textit{Carl Friedrich Gauss}'s \textit{Theorema Egregium} \cite{gauss1828disquisitiones} proved in 1828 that this curvature can be determined entirely by measurements (like angles and distances) made \textit{within} the surface, without any reference to how the surface is embedded in a higher-dimensional space. This makes it a primary example of \textit{intrinsic curvature} (see Section \ref{section_intrinsic_extrinsic_curvature}).

\textit{Bernhard Riemann} was the student of \textit{Carl Friedrich Gauss}. Gauss proposed the curvature of a two-dimensional manifold in 1827. This curvature is named Gaussian curvature, named after Gauss. Then, later, Riemann extended curvature to any $n$-dimensional manifold in a lecture in 1854, which was published posthumously in 1868. 

\subsubsection{Definition of Gaussian Curvature}

\begin{definition}[Gaussian curvature \cite{gauss1828disquisitiones}]
For a two-diemnsional Riemannian manifold $(\mathcal{M}, g)$, the \textbf{Gaussian curvature} $K$ at a point $\b{p}$ is defined as the sectional curvature of the tangent space $T_{\b{p}}\mathcal{M}$ (which is itself a 2-plane). It is also equal to half of the scalar curvature $R$ in two dimensions:
\begin{equation}\label{equation_gaussian_curvature_def}
\boxed{
K := \frac{R}{2} \overset{(\ref{equation_scalar_curvature_def})}{=} \frac{g^{ij}R_{ij}}{2}. 
}
\end{equation}
\end{definition}

\subsubsection{Derivation of Gaussian Curvature in Coordinates}

\begin{proposition}[Gaussian curvature in coordinates]
In a local coordinate system $\{x^1, x^2\}$, let $g$ be the determinant of the metric tensor matrix $[g_{ij}]$. The Gaussian curvature can be expressed explicitly as:
\begin{equation}\label{eq:gaussian_curvature_coord}
\boxed{
K = \frac{R_{1212}}{g} = \frac{R_{1212}}{g_{11}g_{22} - (g_{12})^2}, 
}
\end{equation}
where $R_{1212}$ is the only independent component of the $(0, 4)$ Riemann curvature tensor for a two-dimensional manifold.
\end{proposition}
\begin{proof}
We know that the Riemann curvature tensor has rank $4$.
In two dimensions, the Riemann curvature tensor has $2^4 = 16$ components, but due to the symmetries and anti-symmetries established in Lemma \ref{lemma_riemann_curvature_symmetry_1_3} and Lemma \ref{lemma_riemann_curvature_symmetry_0_4}, only one component is algebraically independent: $R_{1212}$.
\begin{enumerate}
    \item From the first-pair and last-pair anti-symmetry: $R_{1212} = -R_{2112} = -R_{1221} = R_{2121}$.
    \item Any component with repeated indices in the first or last pair (e.g., $R_{1112}$ or $R_{1222}$) vanishes by anti-symmetry.
    \item The scalar curvature $R$ is given by $g^{ij}R_{ij}$. In two dimensions, using the relation $R_{ij} = g^{kl}R_{kilj}$, we have:
    \begin{align*}
    R &= g^{11}R_{11} + g^{12}R_{12} + g^{21}R_{21} + g^{22}R_{22} \\
    &= 2(g^{11}g^{22} - (g^{12})^2) R_{1212}.
    \end{align*}
\end{enumerate}
Since the determinant of the inverse metric is $\det(g^{ij}) = 1/g$, and $g^{11}g^{22} - (g^{12})^2 = \det(g^{ij})$, we have $R = 2 R_{1212}/g$. Thus, $K = R/2 = R_{1212}/g$.
\end{proof}

This explicit coordinate form is particularly useful for optimization on 2D surfaces, where the metric tensor is easily defined and the curvature serves as a local measure of the ``difficulty" of the optimization landscape.

\subsection{Other Advanced Curvature Measurements}

Note that there are also more advanced curvature measurements, which we briefly cover. Some of them are Cotton tensor, Weyl Curvature tensor, and Schouten tensor. 




\subsubsection{Schouten Tensor}

\begin{definition}[Schouten tensor \cite{besse1987einstein}]
Let $(\mathcal{M},g)$ be an $n$-dimensional Riemannian manifold with $n \geq 3$. 
The \textbf{Schouten tensor} is a $(0,2)$ tensor defined as:
\begin{align}
\boxed{
P_{ij} := \frac{1}{n-2} \left( R_{ij} - \frac{R}{2(n-1)} g_{ij} \right),
}
\end{align}
where $R_{ij}$ is the Ricci tensor, $R$ is the scalar curvature, and $g_{ij}$ is the metric tensor. 
\end{definition}

\begin{remark}[Intuition of Schouten tensor]
The Schouten tensor is a modified version of the Ricci tensor that is adapted to conformal geometry. 
It appears naturally in the decomposition of the Riemann curvature tensor.
\end{remark}

\subsubsection{Weyl Curvature Tensor}

\begin{definition}[Weyl curvature tensor \cite{weyl1918reine}]
Let $(\mathcal{M},g)$ be an $n$-dimensional Riemannian manifold with $n \geq 3$. 
The \textbf{Weyl curvature tensor} is a $(0,4)$ tensor defined as:
\begin{align}
\boxed{
W_{ijkl} 
:= R_{ijkl} 
- \left( g_{ik} P_{jl} - g_{il} P_{jk} - g_{jk} P_{il} + g_{jl} P_{ik} \right).
}
\end{align}
\end{definition}


\begin{remark}[Information definition of conformal flatness]
Conformal flatness means that the metric is locally a scalar multiple of a flat metric. 
Equivalently, the geometry agrees with Euclidean geometry up to a pointwise scaling 
that preserves angles but not necessarily lengths.
\end{remark}

\begin{remark}[Intuition of Weyl tensor]
The Weyl tensor is the trace-free part of the Riemann curvature tensor. 
It can be proved that, in dimension $n \geq 4$, $W = 0$ if and only if the manifold is locally conformally flat.
\end{remark}

\subsubsection{Cotton Tensor}

\begin{definition}[Cotton tensor \cite{cotton1899varietes}]
Let $(\mathcal{M},g)$ be an $n$-dimensional Riemannian manifold with $n \geq 3$. 
The \textbf{Cotton tensor} is a $(0,3)$ tensor defined as:
\begin{align}
\boxed{
C_{ijk} := \nabla_k P_{ij} - \nabla_j P_{ik}.
}
\end{align}
\end{definition}

\begin{definition}[Codazzi tensor]
Let $(\mathcal{M},g)$ be a Riemannian manifold with Levi-Civita connection $\nabla$. 
A $(0,2)$ tensor $T_{ij}$ is called a Codazzi tensor if:
\begin{align*}
\nabla_k T_{ij} = \nabla_j T_{ik}.
\end{align*}
\end{definition}

\begin{remark}[Intuition of Cotton tensor]
The Cotton tensor measures the failure of the Schouten tensor to be Codazzi. 
It can be proved that, in dimension $n = 3$, the Weyl tensor vanishes identically, and the Cotton tensor completely characterizes conformal flatness.
\end{remark}


\section{Ricci Flow}\label{section_ricci_flow}

\textit{Ricci flow} is an important concept in differential geometry, which makes use of Ricci curvature. It has recently gained attention from computational scientists in different fields. Here, we briefly introduce it so computational scientists can use it in their algorithms. 

\subsection{Intuition of Ricci Flow}

\begin{figure*}[!t]
\centering
\includegraphics[width=6.5in]{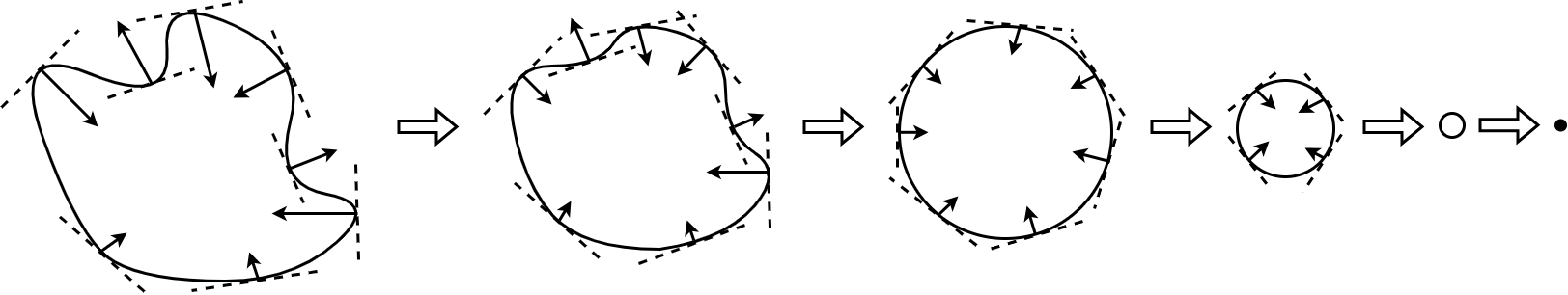}
\caption{Ricci flow on an example manifold where the manifold gradually becomes smoother as a sphere and then disappears. The flow's direction at each point is inverse of the sign of curvature; that is why Ricci flow has negative sign in Eq. (\ref{eq:ricci_flow_pde}). The magnitude of flow is proportional to the amount of curvature, as also obvious in Eq. (\ref{eq:ricci_flow_pde}).}
\label{figure_Ricci_flow}
\end{figure*}

The modern development of Ricci flow was initiated by \textit{Richard S. Hamilton} in 1982 \cite{hamilton1982three} and was profoundly advanced by \textit{Grigori (Grisha) Perelman} in 2002 and 2003 \cite{perelman2002entropy,perelman2003ricci,perelman2003finite}.
The Ricci flow is a geometric process that deforms the Riemannian metric $g_{ij}$ over time in a way that is proportional to the Ricci curvature tensor $R_{ij}$. Intuitively, it can be viewed as a heat equation for the metric, where regions of positive curvature contract and regions of negative curvature expand to ``smooth out'' the manifold's geometry. Such a behavior is depicted in Fig. \ref{figure_Ricci_flow}.

\subsection{The Ricci Flow Equation}

\begin{definition}[The Ricci flow equation \cite{hamilton1982three}]
Let $(\mathcal{M}, g(t))$ be a family of Riemannian manifolds indexed by time $t$. The \textbf{Ricci flow} is defined by the following partial differential equation:
\begin{equation}\label{eq:ricci_flow_pde}
\boxed{
\frac{\partial}{\partial t} g_{ij}(x, t) = -2 R_{ij}(x, t),
}
\end{equation}
where $g_{ij}$ is the metric tensor and $R_{ij}$ is the Ricci curvature tensor at time $t$.
\end{definition}

\begin{proposition}[Ricci flow in local coordinates]
In a local coordinate system $(x^1, \dots, x^n)$, the evolution of the metric under Ricci flow can be expressed by substituting the coordinate derivation of the Ricci tensor (see Eq. (\ref{equation_Ricci_curvature_coordinates})):
\begin{equation}
\boxed{
    \frac{\partial}{\partial t} g_{ij} = -2 \left( \partial_k \Gamma^k_{ji} - \partial_j \Gamma^k_{ki} + \Gamma^k_{kp} \Gamma^p_{ji} - \Gamma^k_{jp} \Gamma^p_{ki} \right),
}
\end{equation}
where the Christoffel symbols $\Gamma^k_{ij}$ are themselves time-dependent since they are functions of $g_{ij}(t)$.
\end{proposition}

\subsection{Intuition for Riemannian Optimization}\label{section_ricci_flow_intuition_for_riemannian_optimization}

In the context of Riemannian optimization (see Section \ref{section_manifold_valued_optimization}), the Ricci flow serves as a theoretical foundation for understanding how the underlying geometry of a matrix manifold changes under certain deformations. For computational scientists, this is relevant in:
\begin{itemize}
    \item \textbf{Geometric Smoothing:} Regularizing noisy manifold data by ``flowing'' the metric toward a more uniform curvature state.
    \item \textbf{Metric Learning:} Adjusting the metric $g_{ij}$ to better suit the convergence properties of Riemannian Gradient Descent (see Section \ref{section_riemannian_gradient_descent}).
\end{itemize}

\subsection{Short-time Existence and Convergence}
For any smooth, compact Riemannian manifold, a unique solution to the Ricci flow equation exists for a short time interval $[0, \epsilon)$. Depending on the initial curvature, the flow may:
\begin{enumerate}
    \item \textbf{Shrink to a point:} Occurs in manifolds with strictly positive curvature (e.g., a sphere $S^n$)\footnote{Recall Eq. (\ref{equation_n_sphere}) in Definition \ref{definition_n_sphere} for definition of $S^n$.}.
    \item \textbf{Expand infinitely:} Occurs in manifolds with negative curvature.
    \item \textbf{Converge to a steady state:} Occurs in ``Ricci-flat'' manifolds where $R_{ij} = 0$.
\end{enumerate}

The Ricci flow, introduced by \textit{Richard S. Hamilton} \cite{hamilton1982three}, plays a central role in the proof of the \textit{Poincaré conjecture}. The conjecture---proposed by \textit{Henri Poincaré} in 1904 \cite{poincare1904cinquieme}---states that every closed, connected, simply connected three-dimensional topological manifold is homeomorphic to the three-dimensional sphere $S^3$ (see Eq. (\ref{equation_n_sphere}) in Definition \ref{definition_n_sphere} for definition of $3$-sphere). The intuition of this conjecture is easy: take any locally three-dimensional topology, having no holes, and mold it in a way to become like a three-dimensional sphere. However, its proof was hard puzzling researchers for almost a century. 

In 1982, \textit{Richard S. Hamilton} introduced the Ricci flow and proved that \cite{hamilton1982three}:
\begin{itemize}
\item For certain initial metrics (e.g., positive Ricci curvature), the solution exists only for a finite time.
\item As time approaches this finite value, the curvature becomes unbounded. This blow-up of curvature is called \textit{singularity} in the Ricci flow.
\end{itemize}

Building on Hamilton's work, \textit{Grigori (Grisha) Perelman} proved the conjecture in a series of groundbreaking papers in 2002 and 2003 \cite{perelman2002entropy,perelman2003ricci,perelman2003finite}. His work developed new analytic tools for the Ricci flow, including entropy formulas (called \textit{Perelman's entropy}) and a detailed analysis of \textit{singularities}, allowing the flow with surgery to be carried out rigorously.
He handled singularities via \textit{surgery}, a procedure in which regions of high curvature that develop under the flow are excised and replaced with standard geometric pieces, after which the Ricci flow is continued on the modified manifold.

\section{Curves (Paths), Parallel Transport, and Geodesics}\label{section_curves_parallel_transport_geodesics}

\subsection{Absolute Differential}\label{section_absolue_differential}

\begin{definition}[Absolute differential]
Let $\b{V} = V^i \b{e}_i$ be a vector field expressed in coordinate bases $\{\b{e}_i\}_{i=1}^n$.
Let $\{x^i\}_{i=1}^n$ be the coordinate system on the manifold. 
The \textbf{absolute differential} (also called \textbf{covariant differential} or \textbf{covariant change}) of $V^i$ is the change in the $i$-th component of vector going from one point to a neighboring point on the trajectory. 
It is denoted by $D V^i$ and is defined as:
\begin{align}\label{equation_absolute_differential}
\boxed{
D V^i := (\nabla_j V^i)\, dx^j,
}
\end{align}
where $\nabla_j V^i$ is the covariant derivative of $V^i$ with respect to coordinate $x^j$ and $d x^j$ is the differential (infinitesimal change) in the direction of coordinate $x^j$. 
\end{definition}

Note that Eq. (\ref{equation_absolute_differential}) is simple to understand. According to chain rule in derivatives, $\nabla_j V^i$ behaves like covariant derivative with respect to small change in $x^j$, i.e., $dx^j$ so multiplying it by $dx^j$ gives $D V^i$.

\begin{proposition}[Absolute differential equation]
Let $\b{V} = V^j \b{e}_i$ be a vector field expressed in coordinate bases $\{\b{e}_j\}_{j=1}^n$.
The absolute differential (also called covariant differential or covariant change) of $V^j$ is equal to:
\begin{align}
\boxed{
D V^i = d V^i + \Gamma^i_{jk} V^k dx^j.
}
\end{align}
\end{proposition}
\begin{proof}
According to Eq. (\ref{equation_covariant_derivative_contravariant}), we have:
\begin{align*}
\nabla_j V^i = \partial_j V^i + \Gamma^i_{jk} V^k \overset{(\ref{equation_partial_i})}{=} \frac{\partial V^i}{\partial x^j} + \Gamma^i_{jk} V^k.
\end{align*}
Multiplying the sides of equation by $dx^j$ gives:
\begin{align*}
&(\nabla_j V^i) dx^j = \frac{\partial V^i}{\partial x^j} dx^j + \Gamma^i_{jk} V^k dx^j \\
&\overset{(a)}{\implies} D V^i = d V^i + \Gamma^i_{jk} V^k dx^j,
\end{align*}
where $(a)$ is because $(\nabla_j V^i) dx^j = D V^i$ according to Eq. (\ref{equation_absolute_differential}) and $(b)$ is because of $ \frac{\partial V^i}{\partial x^j} dx^j = d V^i$ according to chain rule in derivatives. 
\end{proof}

\subsection{Curves (Paths), Parameterization, and Velocity Vector}\label{section_curves}

\begin{figure}[!h]
\centering
\includegraphics[width=3.2in]{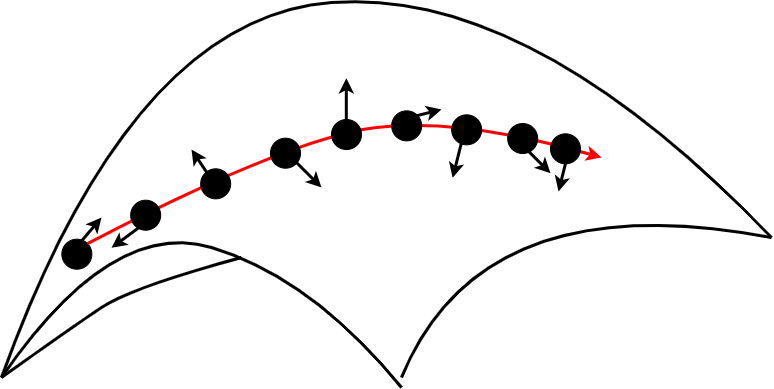}
\caption{A vector field changing along a path on the manifold. The path is illustrated in red and the vectors of vector field are colored black.}
\label{figure_vector_field_along_path}
\end{figure}

Vector fields can change along different paths (or curves) on the manifold. For example, Fig. \ref{figure_vector_field_along_path} illustrates a vector field changing along a path on the manifold. 
To analyze how vector fields change along a specific path, we must first formally define a curve and its associated parameter.

\begin{definition}[Smooth curve or path]
A \textbf{smooth curve} (or \textbf{path}) in an $n$-dimensional manifold $M$ is denoted by $\b{\gamma}$ and is a differentiable mapping from an interval of the real line into the manifold:
\begin{equation}
\boxed{
\begin{aligned}
&\b{\gamma}: I \to \mathcal{M}, \\
&\b{\gamma}: t \mapsto \b{\gamma}(t), \quad \forall t \in I,
\end{aligned}
}
\end{equation}
where $I \subseteq \mathbb{R}$ is an interval of the real line and $t \in I$ is the \textbf{parameter} (or \textbf{time}) of the curve.
The path or curve is denoted by $\b{\gamma}(t)$.
\end{definition}

\begin{definition}[Parameter of the curve]\label{definition_parameter_of_curve}
The \textbf{parameter $t$ of the curve $\b{\gamma}(t)$} can be thought of as a ``coordinate" along the one-dimensional worldline of the path. 
The parameter $t$ increases as we move along the curve. 
In physical applications, $t$ often represents time or the arc length $s$.
\end{definition}

\begin{remark}[Curve in local coordinates]
In a local coordinate system $\{x^i\}_{i=1}^n$, the curve $\b{\gamma}(t)$ is expressed as a set of functions $\{x^i(t)\}_{i=1}^n$:
\begin{equation}
\boxed{x^i = x^i(t)}, \quad \forall t \in I, \quad i \in \{1, \dots, n\},
\end{equation}
where $x^i(t)$ is the $i$-th coordinate of the point on curve (path) at time $t$. Note that the $i$-th coordinate of the point on curve is a function of parameter (time), so it is written as $x^i(t)$. 

A curve $\b{\gamma}(t)$ is technically a point on the manifold $\mathcal{M}$. The list of functions $[x^1(t), \dots, x^n(t)]$ is the image of that point in $\mathbb{R}^n$ under a coordinate chart $\varphi$.
\begin{equation}
\boxed{
\b{\gamma}(t) = \varphi^{-1}\left(x^1(t), x^2(t), \dots, x^n(t)\right).
}
\end{equation}
\end{remark}


\begin{definition}[Velocity vector or tangent vector]\label{definition_velocity_vector}
At any point $p = \b{\gamma}(t)$ on the curve (path), the \textbf{velocity vector} (or \textbf{tangent vector})---denoted by $\dot{\b{\gamma}}(t)$ or $\b{U}$---is a vector tangent to both the curve and the manifold:
\begin{align}\label{equation_velocity_vector2}
\boxed{
\dot{\b{\gamma}}(t) := \frac{d}{dt} \b{\gamma}(t) \in T_{\b{p}}\mathcal{M}.
}
\end{align}
The components of the velocity vector (or tangent vector) $\b{U}$ to the curve $\b{\gamma}(t)$ at any point are:
\begin{equation}\label{equation_velocity_vector}
\boxed{
U^i := \frac{dx^i}{dt} = \dot{x}^i,
}
\end{equation}
where the vector itself is given by:
\begin{align}
\boxed{
\b{U} \overset{(\ref{equation_contravariant_components})}{=} U^i \b{e}_i \overset{(\ref{equation_velocity_vector})}{=} \dot{x}^i \b{e}_i,
}
\end{align}
where $\{\b{e}_i\}_{i=1}^n$ is the coordinate basis and $\dot{x}^i$ denotes the $i$-th coordinate of the velocity (tangent) vector. 
The notation of dot on top of a variable means the derivative of that variable with respect to $t$.

Note that some texts in the literature denote the velocity vector by $\dot{\b{\gamma}}(t)$, and its $i$-th coordinate by $\dot{\gamma}^i(t)$ acting as $\dot{x}^i(t)$:
\begin{equation}\label{equation_gamma_dot_i}
\boxed{
\begin{aligned}
& \dot{\b{\gamma}} := \dot{\gamma}^i \b{e}_i = \b{U} = U^i \b{e}_i = \dot{x}^i \b{e}_i, \\
& \dot{\gamma}^i := \dot{x}^i = \frac{dx^i}{dt}.
\end{aligned}
}
\end{equation}
\end{definition}

\subsection{Absolute Derivative (Intrinsic Derivative) Along a Curve}

Using the definition of the absolute differential $D V^i$ from Section \ref{section_absolue_differential}, we can define the rate of change of a vector field $V^i$ as we move along $\b{\gamma}(t)$. 

\begin{definition}[Absolute derivative or intrinsic derivative along a curve]
By dividing the absolute differential by the parameter increment $dt$, we obtain the \textbf{absolute derivative} (also called the \textbf{intrinsic derivative}):
\begin{equation}\label{equation_absolute_derivative2}
\boxed{
\frac{DV^i}{dt} := \frac{dV^i}{dt} + \Gamma^i_{jk} V^k \frac{dx^j}{dt}.
}
\end{equation}
Substituting the velocity components $\dot{x}^j$ from Eq. (\ref{equation_velocity_vector}), this is often written as:
\begin{equation}\label{equation_absolute_derivative}
\boxed{ 
\frac{DV^i}{dt} = \dot{V}^i + \Gamma^i_{jk} V^k \dot{x}^j, 
}
\end{equation}
where:
\begin{align}
\boxed{\dot{V}^i := \frac{dV^i}{dt}}, \quad \boxed{\dot{x}^i := \frac{dx^i}{dt}}.
\end{align}
The absolute derivative (intrinsic derivative) is the rate of change of a vector field $V^i$ as we move along the path $\b{\gamma}(t)$. 
\end{definition}

\begin{remark}
The absolute derivative (intrinsic derivative)---which is the covariant derivative along a curve---has two terms according to Eq. (\ref{equation_absolute_derivative}). The first term captures the local change and the second term is the change due to space curving:
\begin{align}
\boxed{\frac{DV^i}{dt} = \underbrace{\dot{V}^i}_{\text{local change}} + \underbrace{\Gamma^i_{jk} V^k \dot{x}^j}_{\text{change due to space curving}}}.
\end{align}
\end{remark}

\subsection{Parallel Transport}\label{section_parallel_transport}


\subsubsection{Definition and Equation of Parallel Transport}

\begin{definition}[Parallel transport]
\textbf{Parallel transport} refers to moving (transporting) a vector on the manifold in a way that it remains parallel to itself along the path of movement. 
The condition for parallel transport of a vector $V^i$ along the curve $\b{\gamma}(t)$ is that its absolute derivative is zero along the path of movement:
\begin{equation}\label{equation_parallel_transport_1}
\boxed{
\frac{DV^i}{dt} = 0.
}
\end{equation}
This implies the vector is kept parallel to itself in the curved geometry as it moves along the path.
According to Eq. (\ref{equation_absolute_derivative}), the Eq. (\ref{equation_parallel_transport_1}) for parallel transport can be stated as:
\begin{align}\label{equation_parallel_transport_2}
\boxed{
\dot{V}^i + \Gamma^i_{jk} V^k \dot{x}^j = 0.
}
\end{align}
According to Eq. (\ref{equation_gamma_dot_i}), some texts in the literature denote this equation of parallel transport as:
\begin{align}
\boxed{
\dot{V}^i + \Gamma^i_{jk} V^k \dot{\gamma}^j = 0.
}
\end{align}
\end{definition}

In 1916, \textit{Tullio Levi-Civita} (from Italy)---who was a student of \textit{Gregorio Ricci-Curbastro}---employed Christoffel symbols to formulate the concept of parallel transport and to investigate how parallel transport is connected to curvature \cite{levi1916nozione}.

\subsubsection{Holonomy: The Closed Loop Paradox}

Parallel transport is path-dependent. 
When the manifold is flat with zero curvature, parallel transporting a vector along a closed path---having the same starting and ending points---will give the same vector at the end of path. 
However, when the manifold is curved, the vector at the end of a closed path will not be the same as the vector at the start of path, although we had parallel transporting the vector. 
In other words, parallel transporting is affected by curvature. Therefore, it depends on the path because the path may pass the parts of the manifold having curvature. 

\begin{figure}[!h]
\centering
\includegraphics[width=1.6in]{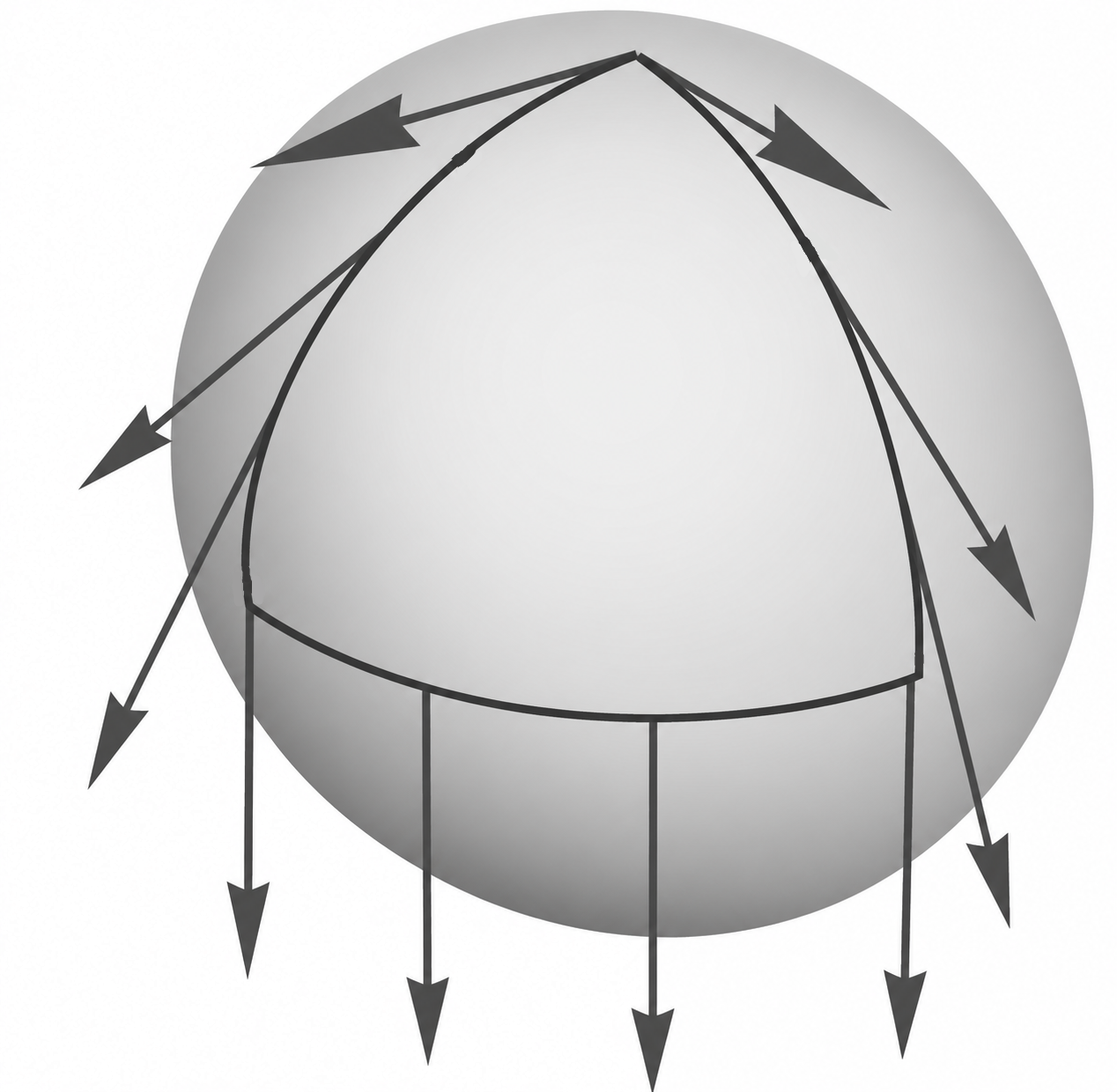}
\caption{Holonomy (a closed loop) where we parallel transport a vector along a closed loop on a curvy manifold but the vector does not return as it started. This behavior is because of curvature of the manifold.}
\label{figure_holonomy}
\end{figure}

A classic way to visualize path-dependence is to consider a closed loop (a \textit{holonomy}).
As illustrated in Fig. \ref{figure_holonomy}, imagine starting at the North Pole of a sphere with a vector pointing toward the equator.
Parallel transport it down to the equator.
Move it along the equator by 90 degrees.
Transport it back to the North Pole.
Upon returning to the start, the vector will no longer point in its original direction. The angle of deviation is directly proportional to the amount of curvature enclosed by the path.

This behavior is because, in a curved manifold, the geometry warps the vector as it moves. Because the Christoffel symbols $\Gamma^i_{jk}$ (which represent the connection) vary across the manifold, the integration of the parallel transport equation---i.e., Eq. (\ref{equation_parallel_transport_1})---yields different results for different paths.

\begin{remark}[Geometrical intuition of torsion]
Torsion, defined in Definition \ref{definition_torsion}, determines whether parallel transport around an infinitesimal parallelogram closes.
If torsion is nonzero, the parallelogram does not close.
\end{remark}

\subsection{Geodesics}\label{section_geodesics}


In a flat Euclidean space, which does not have any curvature, the shortest path between two points is the straight line between them. 
However, in a curvy manifold, the shortest path (or shortest curve) between two points on the manifold is not necessarily the straight line because of the curvature. The shortest curve on the manifold is called \textit{geodesic}. 

\begin{figure}[!h]
\centering
\includegraphics[width=3.2in]{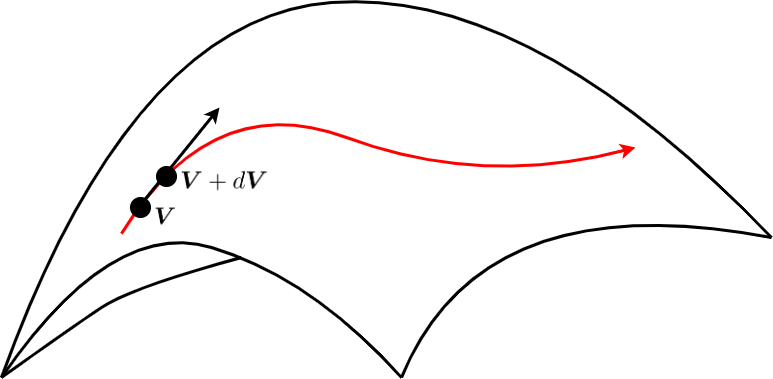}
\caption{A tangent vector along a curve (path). The path is illustrated in red and the vectors of vector field are colored black.}
\label{figure_geodesic_tangent_vector}
\end{figure}

Consider Fig. \ref{figure_geodesic_tangent_vector} which depicts a tangent vector along a curve (path). 
Mathematically, on a geodesic path, the tangent vector (velocity vector) stays parallel to itself along the path. In other words, the covariant derivative (or covariant change) of the tangent vector is zero along the entire geodesic. This means that the tangent vector along the geodesic curve is constant and does not change. 

\begin{definition}[Geodesic]
Let $(\mathcal{M}, \nabla)$ be a smooth manifold equipped with a connection $\nabla$. A smooth curve $\gamma: I \to \mathcal{M}$ is called a \textbf{geodesic} if the covariant derivative of the velocity vector (the covariant derivative of the tangent vector) is zero along the geodesic path.
In other words, for a geodesic curve, the velocity vector field $\b{V} = \dot{\b{\gamma}}$ is parallel transported along the curve itself. In terms of the absolute derivative (intrinsic derivative), this condition is expressed as:
\begin{align}\label{equation_geodesic_absolute}
\frac{DV^i}{dt} = 0
, \quad \forall i \in \{1, \dots, n\}.
\end{align}
As the vector $\b{V} = V^i \b{e}_i$ here is the velocity vector, this equation can be stated as:
\begin{align}\label{equation_geodesic_absolute2}
\boxed{
\frac{D\dot{x}^i}{dt} = 0
}, \quad \forall i \in \{1, \dots, n\},
\end{align}
or equivalently the absolute differential (covariant change) is zero along a geodesic curve:
\begin{align}\label{equation_geodesic_absolute3}
\boxed{
D\dot{x}^i = 0
}, \quad \forall i \in \{1, \dots, n\},
\end{align}

In coordinate-free notation, the defining condition of a
geodesic can be written as:
\begin{align}\label{equation_geodesic_coordinate_free}
\boxed{\nabla_{\dot{\b{\gamma}}(t)} \dot{\b{\gamma}}(t) = \b{0},}
\qquad \forall t \in I,
\end{align}
where $\nabla$ denotes the covariant derivative (see Definition \ref{definition_connection}) and \(\dot{\b{\gamma}}(t)\) is the velocity vector of the
curve \(\b{\gamma}(t)\), defined in Eq. (\ref{equation_velocity_vector2}) in Definition \ref{definition_velocity_vector}. 
This means that a smooth curve on a Riemannian manifold is a geodesic if and only if \textbf{its covariant acceleration vanishes}.
In other words, the velocity vector is parallel transported along the curve itself, or equivalently, the curve has zero intrinsic acceleration.
\end{definition}

\begin{figure}[!h]
\centering
\includegraphics[width=3.2in]{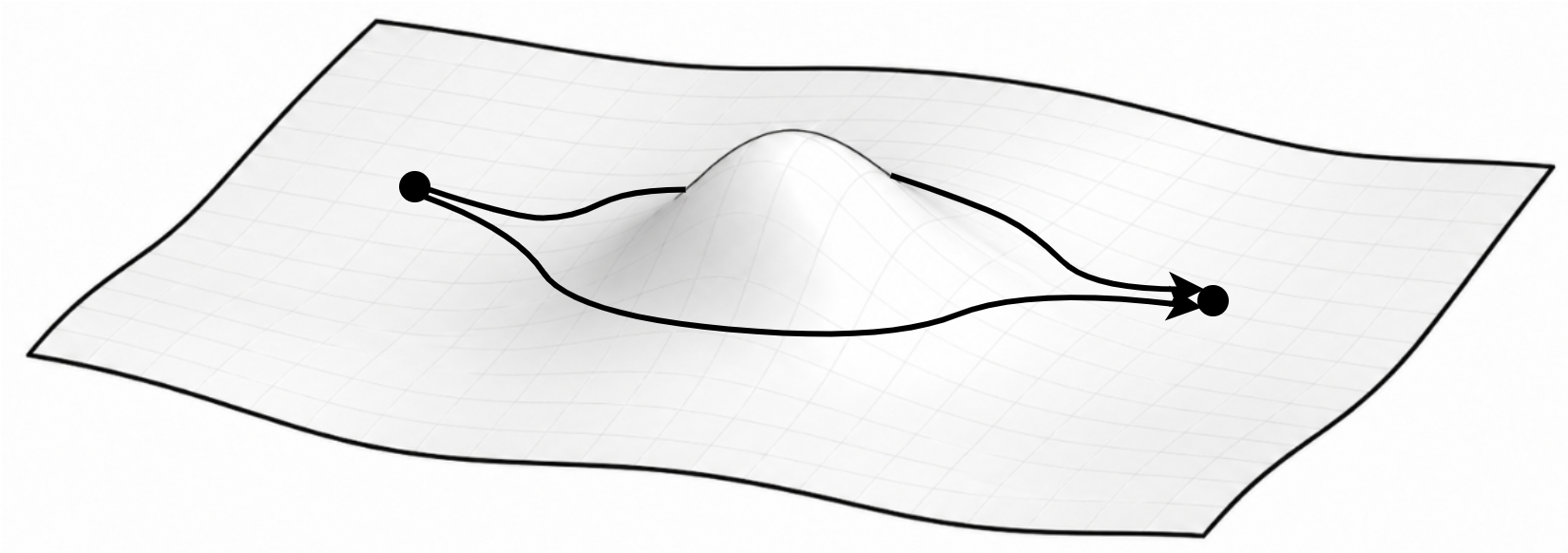}
\caption{Two geodesics, with shortest length, between two points on the manifold.}
\label{figure_geodesic_non_unique}
\end{figure}

\begin{remark}[Non-uniqueness of geodesics]
A geodesic between two points on a Riemannian manifold is not necessarily unique. While a geodesic is locally distance-minimizing, there can exist multiple distinct geodesics of equal (or different) lengths connecting the same two points. Figure \ref{figure_geodesic_non_unique} illustrates a case where two points are connected by multiple geodesics of the same length.
\end{remark}

\begin{proposition}[The geodesic equation]
The coordinate functions $x^k(t)$ of a geodesic $\b{\gamma}(t)$ satisfy the following system of second-order non-linear ordinary differential equations:
\begin{equation}\label{equation_geodesic_ode}
\boxed{
\ddot{x}^k + \Gamma^k_{ij} \dot{x}^i \dot{x}^j = 0}, \quad \forall k \in \{1, \dots, n\},
\end{equation}
where $\Gamma^k_{ij}$ are the Christoffel symbols of the connection, and:
\begin{equation}\label{equation_x_dot_x_ddot}
\boxed{
\begin{aligned}
& \dot{x}^i = \frac{dx^i}{dt}, \\
& \ddot{x}^k = \frac{d}{dt} \left(\frac{dx^k}{dt}\right) = \frac{d^2 x^k}{dt^2}.
\end{aligned}
}
\end{equation}
According to Eq. (\ref{equation_gamma_dot_i}), some texts in the literature denote Eq. (\ref{equation_geodesic_ode}) as:
\begin{equation}\label{equation_geodesic_ode2}
\boxed{
\ddot{\gamma}^k + \Gamma^k_{ij} \dot{\gamma}^i \dot{\gamma}^j = 0
}, \quad \forall k \in \{1, \dots, n\}.
\end{equation}
\end{proposition}
\begin{proof}
Recall Eq. (\ref{equation_absolute_derivative2}) for the definition of the absolute derivative for a vector field along a curve:
\begin{equation*}
\frac{DV^i}{dt} = \frac{dV^i}{dt} + \Gamma^i_{jk} V^k \frac{dx^j}{dt}.
\end{equation*}
According to Eq. (\ref{equation_geodesic_absolute}), the geodesic satisfies:
\begin{equation*}
\frac{DV^i}{dt} = 0,
\end{equation*}
where the vector $\b{V} = V^i \b{e}_i$ is the velocity vector $\dot{x}^i \b{e}_i$ (in other words, $V^i = \dot{x}^i$). 
By combining these two equations, we have:
\begin{equation}\label{equation_proof_DVdt_dVdt_GammaVdxdt}
\frac{dV^i}{dt} + \Gamma^i_{jk} V^k \frac{dx^j}{dt} = 0.
\end{equation}
Substituting $V^i = \dot{x}^i$ into Eq. (\ref{equation_proof_DVdt_dVdt_GammaVdxdt}) gives:
\begin{equation*}
\frac{D\dot{x}^i}{dt} = \frac{d\dot{x}^i}{dt} + \Gamma^i_{jk} \dot{x}^k \frac{dx^j}{dt} \overset{(\ref{equation_x_dot_x_ddot})}{\implies} \ddot{x}^i + \Gamma^i_{jk} \dot{x}^k \dot{x}^j = 0.
\end{equation*}
Relabeling the dummy indices $i \to k$ and $k \to i$ gives:
\begin{align*}
\ddot{x}^k + \Gamma^k_{ji} \dot{x}^i \dot{x}^j = 0.
\end{align*}
Assuming a torsion-free connection, the Christoffel symbols are symmetric (see Eq. (\ref{equation_symmtery_Christoffel})), so $\Gamma^k_{ji} = \Gamma^k_{ij}$. Therefore, $\ddot{x}^k + \Gamma^k_{ij} \dot{x}^i \dot{x}^j = 0$. 
\end{proof}

\begin{remark}[Physics interpretation of the geodesic equation \cite{susskind2025general}]
Equation (\ref{equation_geodesic_ode}) can be interpreted as follows. According to Definition \ref{definition_parameter_of_curve}, assume the parameter of curve is time $t$. 
According to Eq. (\ref{equation_geodesic_ode}):
\begin{equation}\label{equation_xddot_minus_Gamma_xdot_xdot}
\ddot{x}^k = - \Gamma^k_{ij} \dot{x}^i \dot{x}^j.
\end{equation}
The $\ddot{x}^k$ can be interpreted as acceleration, which is the second derivative with respect to time, and $\dot{x}^i$ and $\dot{x}^j$ are components of velocity vectors or tangent vectors (see Definition \ref{definition_velocity_vector}). 
The Christoffel symbol $\Gamma^k_{ij}$ depends on the metric, according to Eq. (\ref{equation_Christoffel_secondType_metric}).
In general relativity, we know that metric is the gravitational field \cite{susskind2025general}. 

Therefore, the left-hand side of Eq. (\ref{equation_xddot_minus_Gamma_xdot_xdot}) is like acceleration, while the right-hand side depends on metric (or gravitational field) and components of the velocity vector. It is similar to Newton's equation $a = F/m$ or equivalently $F = ma$, where $a$ is the acceleration, $F$ is force, and $m$ is mass. Thus:
\begin{equation}
\boxed{
\ddot{x}^k = - \Gamma^k_{ij} \dot{x}^i \dot{x}^j \iff a = \frac{F}{m}.
}
\end{equation}
This shows that the geodesic equation (\ref{equation_geodesic_ode}) can be interpreted as motion of particle in a gravitational field. 
\end{remark}

\section{Riemannian Optimization}\label{section_manifold_valued_optimization}









While the preceding sections established the rigorous foundations of Riemannian geometry, this section transitions those abstract definitions into the computational framework required for \textit{Riemannian optimization}.
Riemannian optimization is also called \textit{optimization on smooth manifolds} or \textit{manifold-valued optimization} because the optimization variables are on a smooth manifold. 

\subsection{Euclidean Optimization versus Riemannian Optimization}

\begin{figure*}[!t]
\centering
\includegraphics[width=5in]{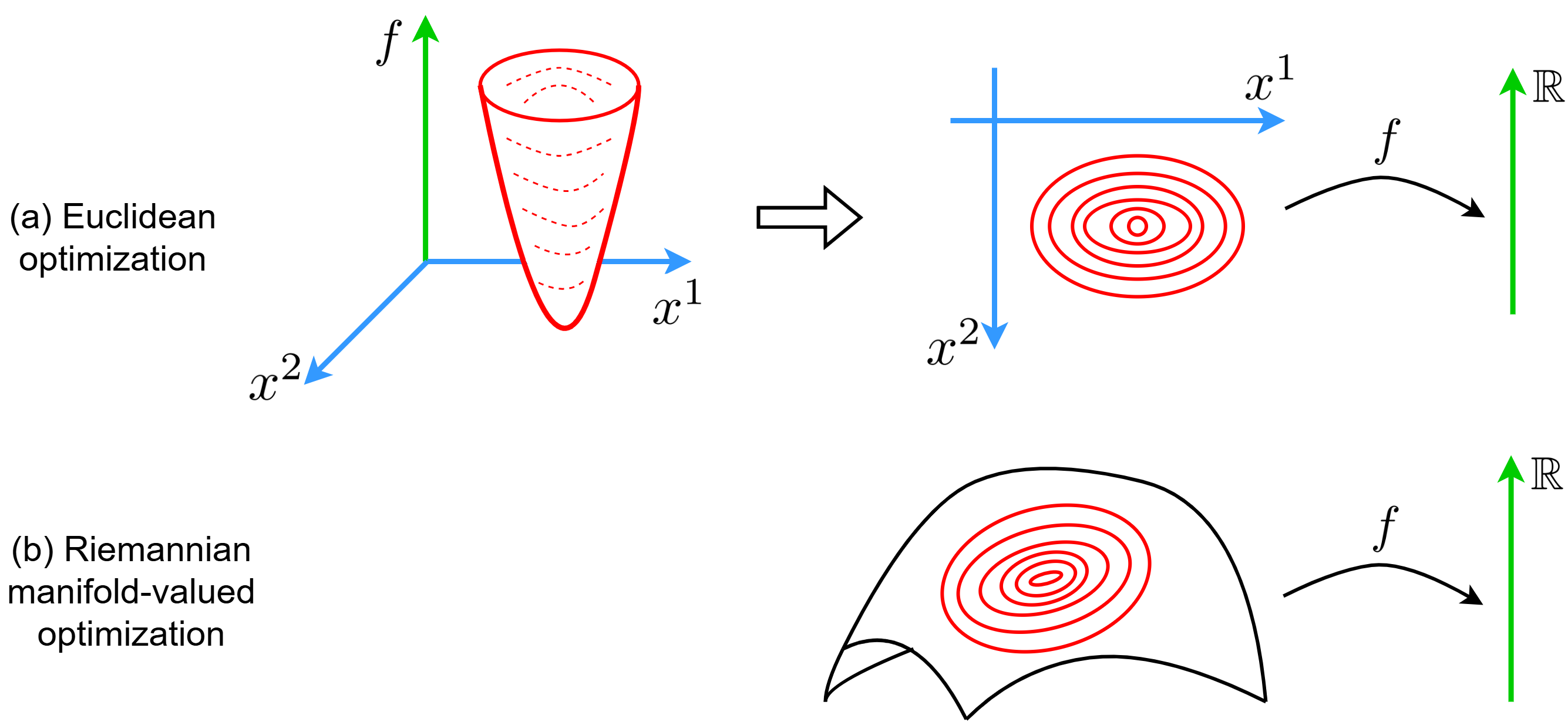}
\caption{(a) Euclidean optimization versus (b) Riemannian optimization. The cost function (contours of cost function) is colored in red, the coordinates are colored in blue, and the real axis for the cost function is colored in green. The red contours show the value levels of cost function in the space or manifold.}
\label{figure_Euclidean_vs_Riemannian_optimization}
\end{figure*}

In \textit{Euclidean optimization}---which is the regular optimization approach---we typically seek a point $x \in \mathbb{R}^n$ that minimizes a cost function $f(x)$. 
In Euclidean optimization, the cost function is a function from the Euclidean space to a scalar: 
\begin{equation}
\begin{aligned}
&f: \mathbb{R}^n \rightarrow \mathbb{R}, \\
&f: \b{x} \mapsto f(\b{x}).
\end{aligned}
\end{equation}
The optimization problem is:
\begin{equation}\label{equation_Euclidean_optimization}
\boxed{
\begin{aligned}
& \underset{\b{x} \in \mathbb{R}^n}{\text{minimize}}
& & f(\b{x}),
\end{aligned}
}
\end{equation}
or equivalently:
\begin{equation}\label{equation_Euclidean_optimization2}
\boxed{
\begin{aligned}
& \underset{\b{x}}{\text{minimize}}
& & f(\b{x}) \\
& \text{subject to}
& & \b{x} \in \mathbb{R}^d.
\end{aligned}
}
\end{equation}
An example Euclidean optimization problem is illustrated in Fig. \ref{figure_Euclidean_vs_Riemannian_optimization}-a where the two dimensional variable $\b{x} \in \mathbb{R}^2$ is mapped to a real-valued cost function. 

If the optimization problem is constrained:
\begin{equation}\label{equation_Euclidean_optimization3}
\boxed{
\begin{aligned}
& \underset{\b{x}}{\text{minimize}}
& & f(\b{x}) \\
& \text{subject to}
& & \b{x} \in \mathcal{S},
\end{aligned}
}
\end{equation}
where $\mathcal{S}$ is the feasibility set (constraint set). 

The Euclidean optimization methods can be revised to have optimization on (possibly curvy) smooth manifolds. 
In \textit{Riemannian optimization} or \textit{manifold-valued optimization} \cite{absil2008optimization,boumal2023introduction}, we generalize this to finding a point\footnote{In Riemannian optimization, we denote the variable $\b{x}$ as the point $\b{p}$ in the manifold.} 
$\b{p}$ on a smooth manifold $\mathcal{M}$ that minimizes a smooth real-valued function $f: \mathcal{M} \rightarrow \mathbb{R}$.
Riemannian optimization optimizes a cost function while the variable lies on a smooth manifold $\mathcal{M}$: 
\begin{equation}
\begin{aligned}
&f: \mathcal{M} \rightarrow \mathbb{R}, \\
&f: \b{p} \mapsto f(\b{p}).
\end{aligned}
\end{equation}
An example Riemannian optimization problem is illustrated in Fig. \ref{figure_Euclidean_vs_Riemannian_optimization}-b where the two dimensional variable $\b{p} \in \mathcal{M}$ belongs to the manifold $\mathcal{M}$ and it is mapped to a real-valued cost function. 

The optimization variable $\b{p}$ in the Riemannian optimization is usually matrix rather than vector; hence, Riemannian optimization is also called \textit{optimization on matrix manifolds}.

In Riemannian optimization, the optimization problem is:
\begin{equation}\label{equation_Riemannian_optimization_1}
\boxed{
\begin{aligned}
& \underset{\b{p} \in \mathcal{M}}{\text{minimize}}
& & f(\b{p}),
\end{aligned}
}
\end{equation}
or equivalently:
\begin{equation}\label{equation_Riemannian_optimization}
\boxed{
\begin{aligned}
& \underset{p}{\text{minimize}}
& & f(\b{p}) \\
& \text{subject to}
& & \b{p} \in \mathcal{M}.
\end{aligned}
}
\end{equation}

\begin{remark}[Converting a constrained Euclidean optimization problem to Riemannian optimization]
If the optimization problem is constrained:
\begin{equation}
\begin{aligned}
& \underset{\b{x}}{\text{minimize}}
& & f(\b{x}) \\
& \text{subject to}
& & \b{x} \in \mathcal{S},
\end{aligned}
\end{equation}
where $\mathcal{S}$ is the feasibility set (constraint set), it is possible to use a trick to convert it to a Riemannian optimization. 
If the constraint $\mathcal{S}$ can be modeled as a smooth manifold, we can make it a manifold-value optimization where the variable $\b{p}$ belongs to the manifold of constraint:
\begin{equation}
\begin{aligned}
& \underset{p}{\text{minimize}}
& & f(\b{p}) \\
& \text{subject to}
& & \b{p} \in \mathcal{M},
\end{aligned}
\end{equation}
where:
\begin{align}
\boxed{
\b{x} \in \mathcal{S} \iff \b{p} \in \mathcal{M}.
}
\end{align}

In other words, we may define the constraint as the matrix manifold of that constraint. For example, if the variable is a matrix and the feasibility set $\mathcal{S}$ is orthogonal matrices, we can use Riemannian optimization on the Stiefel manifold\footnote{We will introduce Stiefel manifold later in Section \ref{section_important_Riemannian_matrix_manifolds}.}, which is the manifold of orthogonal matrices, i.e., $\b{X}^\top \b{X} = \b{I}$ where the matrix $\b{X}$ is the point $\b{p} \in \mathcal{M}$ here and $\b{I}$ denotes the identity matrix. Then, we can use optimization on Stiefel manifold to solve the problem. 
\end{remark}

\subsection{Optimization Path}\label{section_optimization_path}


The primary challenge in transition from Euclidean optimization to Riemannian optimization is that the manifold $\mathcal{M}$ is not a flat Euclidean space. Consequently, standard additive updates of the form $\b{x}_{\nu+1} = \b{x}_\nu + \Delta \b{x}$ (with $\nu$ denoting the iteration index) are generally invalid, as the resulting point may not lie on $\mathcal{M}$. To resolve this, we utilize the tangent space $T_{\b{p}}\mathcal{M}$ and geometric maps to ensure all iterates remain on the manifold.

An optimization algorithm on a manifold generates a sequence of points on manifold, $\{{p}_0, {p}_1, {p}_2, \dots\} \subset \mathcal{M}$, where each point has a local coordinate $\{x^i\}_{i=1}^n$ of chart $\varphi$ around it. This discrete sequence can be viewed as a sampling of a continuous \textit{optimization path} (or \textit{optimization curve}) $\b{\gamma}: [0, T] \rightarrow \mathcal{M}$ that traverses the manifold toward a local minimum.

\begin{figure}[!h]
\centering
\includegraphics[width=3.2in]{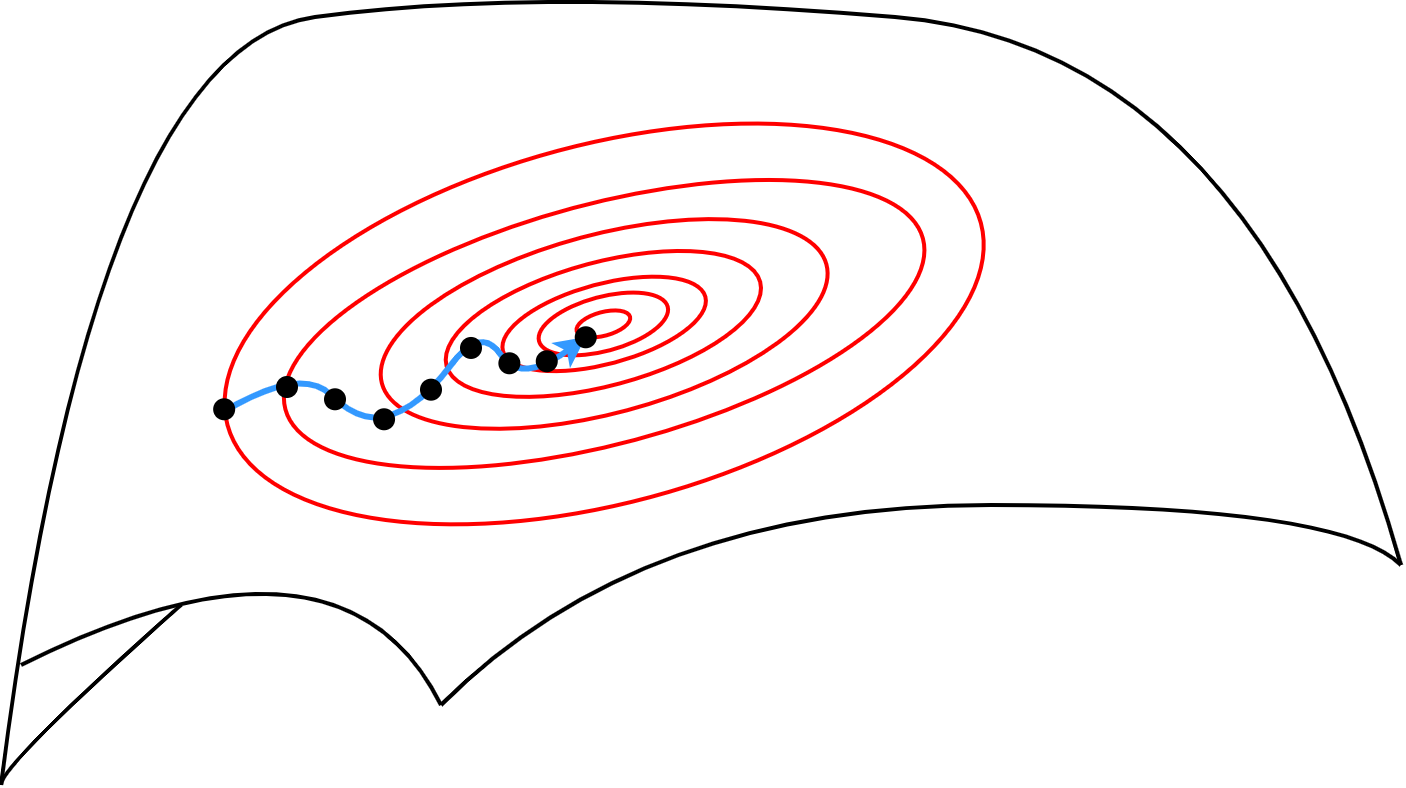}
\caption{Optimization path on the manifold toward the local minimum of the cost function. The contours of cost function are colored in red.}
\label{figure_optimization_path}
\end{figure}

\begin{definition}[Optimization path]
Let $f \in C^\infty(\mathcal{M})$ be a cost function. An \textbf{optimization path} is a smooth curve $\b{\gamma}(t)$:
\begin{align}
\boxed{
\b{\gamma}(t): [0, T] \rightarrow \mathcal{M},
}
\end{align}
that traverses the manifold toward a local minimum.
In other words, the optimization path is a smooth curve $\b{\gamma}(t)$ such that the value of the function decreases along the path, so its derivative with respect to the curve parameter $t$ is non-positive:
\begin{equation}
\boxed{
\frac{d}{dt} f\left(\b{\gamma}(t)\right) \leq 0.
}
\end{equation}
\end{definition}
An example optimization path on, a smooth manifold, toward a local minimum of the cost function, is illustrated in Fig. \ref{figure_optimization_path}.

\begin{definition}[Search direction]
At any point $\b{p} = \b{\gamma}(t)$ on the optimization path, the velocity vector $\dot{\b{\gamma}}(t)$ is a tangent vector in $T_{\b{p}}\mathcal{M}$ (see Definition \ref{definition_velocity_vector}). In the context of optimization, the tangent vector (or velocity vector) represents the \textbf{search direction}.
\end{definition}

Using the coordinate basis $\{\partial_i\}_{i=1}^n$, the optimization path can be expressed locally via the coordinates $x^i(t)$. The velocity vector components are given by $\dot{x}^i = \frac{dx^i}{dt}$.

Ideally, the \textit{straightest} optimization path between two points on the manifold is a \textit{geodesic}. As derived in Section \ref{section_geodesics} and Eq. (\ref{equation_geodesic_ode}), a geodesic path satisfies the equation $\ddot{x}^k + \Gamma^k_{ij} \dot{x}^i \dot{x}^j = 0$. In practical algorithms, we often approximate these paths using \textit{retractions} (to be discussed in Section \ref{section_retraction}) to reduce the computational overhead of solving the second-order geodesic differential equation.

By defining the optimization process as a movement along a path $\b{\gamma}(t)$, we ensure that the intrinsic geometry of the space---specifically the metric $g_{ij}$ and the connection $\Gamma^k_{ij}$---is respected during every iteration of the optimization solver.

\subsection{Riemannian Gradient and Hessian}



In Riemannian optimization, the gradient and Hessian of a scalar function must be defined such that they respect the metric $g$. This section provides the coordinate-level derivations necessary for numerical implementation.

\subsubsection{The Riemannian Gradient}\label{section_Riemannian_gradient}

\begin{definition}[Riemannian gradient --- coordinate-free definition]
Let $f \in C^\infty(\mathcal{M})$ be a smooth function on a Riemannian manifold $(\mathcal{M}, g)$. 
For a smooth function $f \in C^\infty(\mathcal{M})$, the \textbf{Riemannian gradient}, $\operatorname{grad} f$, is the vector field dual to the differential 1-form $df$ with respect to the metric $g$. It is defined by the identity:
\begin{equation}\label{equation_riemannian_gradient}
\boxed{
\begin{aligned}
g(\operatorname{grad} f, \b{X}) &\overset{(\ref{equation_g_inner_product})}{=} \langle \operatorname{grad} f, \b{X} \rangle \\
&= df(\b{X}) = \b{X}(f), \quad \forall \b{X} \in \mathfrak{X}(\mathcal{M}),
\end{aligned}
}
\end{equation}
where $\mathfrak{X}(\mathcal{M})$ denotes the tangent bundle of the manifold and $\langle \cdot, \cdot \rangle$ denotes the inner product. 

Equivalently, at each point $\b{X} \in \mathcal{M}$, the Riemannian gradient $\operatorname{grad} f(\b{X}) \in T_{\b{X}}\mathcal{M}$ is the unique tangent vector satisfying:
\begin{equation}\label{equation_definition_Riemannian_gradient_directional_derivative}
\boxed{
\begin{aligned}
Df(\b{X})[\b{\xi}] = g_{X}(\operatorname{grad} f(\b{X}), \b{\xi}), \quad \forall \b{\xi} \in T_{\b{X}}\mathcal{M},
\end{aligned}
}
\end{equation}
where $Df(\b{X})[\b{\xi}]$ denotes the directional derivative of the function $f$ along the tangent vector $\b{\xi}$ at point $\b{X} \in \mathcal{M}$ (see Eq. (\ref{equation_directional_derivative_function_along_vector_field}) in Definition \ref{definition_directional_derivative_function_along_vector_field}).
\end{definition}

\begin{remark}[Relation of Riemannian gradient and directional derivative]
As established in Definition \ref{definition_directional_derivative_function_along_vector_field}, a vector field $\b{X}$ acts as a differential operator on a smooth function $f: \mathcal{M} \to \mathbb{R}$. The expression $\b{X}(f)$ represents the directional derivative of $f$ in the direction of the vector field $\b{X}$.
It is the weighted sum of the partial derivatives of $f$ along each coordinate direction $x^i$, where the weights are the components $X^i$ of the vector field. The coordinate-based expression of directional derivative is Eq. (\ref{equation_X_f_X_partial_f_partial_x}). 

The directional derivative $\b{X}(f) = X^i \frac{\partial f}{\partial x^i}$ (see Eq. (\ref{equation_X_f_X_partial_f_partial_x})) is equivalent to the inner product of the Riemannian gradient and the vector field:
\begin{align}
\boxed{
\langle \operatorname{grad} f, \b{X} \rangle = df(\b{X}) = \b{X}(f).
}
\end{align}
\end{remark}

\begin{proposition}[Riemannian gradient --- coordinate-based definition]
Consider local coordinate system $\{x^i\}_{i=1}^n$ on a Riemannian manifold, with basis vector $\{\partial_i\}_{i=1}^n$ (see Eq. (\ref{equation_bases_for_vectors_2})). 
The components of the\textbf{ Riemannian gradient}, denoted by $(\operatorname{grad} f)^j$, are given in local coordinates $\{x^i\}_{i=1}^n$ by:
\begin{equation}\label{equation_Riemannian_gradient_components}
\boxed{
(\operatorname{grad} f)^j = g^{ij} \frac{\partial f}{\partial x^i} \overset{(\ref{equation_partial_i})}{=} g^{ij} \partial_i f.
}
\end{equation}
In other words, the Riemannian gradient is a $(1,0)$-tensor (vector field):
\begin{equation}\label{equation_Riemannian_gradient_}
\boxed{
\operatorname{grad} f = \left( g^{ij} \frac{\partial f}{\partial x^i} \right) \partial_j \overset{(\ref{equation_partial_i})}{=} \left( g^{ij} \partial_i f \right) \partial_j.
}
\end{equation}
By Ricci calculus notation (see Definition \ref{definition_comma_semicolon_notations}), we can denote $f_{,i} = \frac{\partial f}{\partial x^i}$ and $f^{;j} = (\operatorname{grad} f)^j$, so Eqs. (\ref{equation_Riemannian_gradient_components}) and (\ref{equation_Riemannian_gradient_}) can be stated as:
\begin{equation}\label{equation_Riemannian_gradient_Ricci}
\boxed{
\begin{aligned}
& f^{;j} = g^{ij} f_{;i} = g^{ij} f_{,i} \\
& \operatorname{grad} f = f^{;j} \partial_j = (g^{ij} f_{,i}) \partial_j. 
\end{aligned}
}
\end{equation}
\end{proposition}
\begin{proof}
On the one hand, according to Eq. (\ref{equation_contravariant_components_2}), the gradient vector can be stated with its components as:
\begin{align}\label{equation_gradf_graf_partiali}
\operatorname{grad} f \overset{(\ref{equation_contravariant_components_2})}{=} (\operatorname{grad} f)^i \partial_i.
\end{align}
According to Eq. (\ref{equation_g_e_e}), we have:
\begin{align}
\langle \operatorname{grad} f, \partial_j \rangle &\overset{(\ref{equation_gradf_graf_partiali})}{=} \langle (\operatorname{grad} f)^i \partial_i, \partial_j \rangle = (\operatorname{grad} f)^i \langle  \partial_i, \partial_j \rangle \nonumber \\
&\overset{(\ref{equation_g_e_e})}{=} (\operatorname{grad} f)^i g_{ij}. \label{equation_inner_gradf_partialj_1}
\end{align}

On the other hand, by the definition of the Riemannian gradient, i.e., Eq. (\ref{equation_riemannian_gradient}), we have:
\begin{align}\label{equation_inner_gradf_partialj_2}
\langle \operatorname{grad} f, \partial_j \rangle \overset{(\ref{equation_riemannian_gradient})}{=} df(\partial_j) \overset{(a)}{=} \partial_j f,
\end{align}
where $(a)$ is because of definition of the differential $df(\b{X}) := \b{X}(f)$ in which $\b{X} = \partial_j$ is used. 

Equating the two expressions for $\langle \operatorname{grad} f, \partial_j \rangle$, i.e., Eqs. (\ref{equation_inner_gradf_partialj_1}) and (\ref{equation_inner_gradf_partialj_2}), we obtain:
\begin{align}\label{equation_grad_f_i_g_ij}
(\operatorname{grad} f)^i g_{ij} = \partial_j f.
\end{align}
To solve for the components $(\operatorname{grad} f)^k$, we multiply both sides of Eq. (\ref{equation_grad_f_i_g_ij}) by the inverse metric $g^{kj}$:
\begin{align*}
(\operatorname{grad} f)^i g_{ij} g^{kj} = g^{kj} \partial_j f.
\end{align*}
Using the property $g_{ij} g^{kj} = \delta_i^k$ (see Eq. (\ref{equation_metric_inverse})), where $\delta_i^k$ is the Kronecker delta, the equation becomes:
\begin{align*}
(\operatorname{grad} f)^i \delta_i^k = g^{kj} \partial_j f.
\end{align*}
By the property of the Kronecker delta, only the term where $i=k$ survives the summation (see Eq. (\ref{equation_index_substitution_delta})):
\begin{align}\label{equation_grad_f_k_components}
(\operatorname{grad} f)^k = g^{kj} \partial_j f = g^{kj} \frac{\partial f}{\partial x^j}.
\end{align}
Finally, substituting Eq. (\ref{equation_grad_f_k_components}) back into Eq. (\ref{equation_gradf_graf_partiali}),where relabeling the dummy variable $k \to i$, we arrive at the explicit coordinate expression for the Riemannian gradient:
\begin{align}
\operatorname{grad} f = \left( g^{ij} \frac{\partial f}{\partial x^i} \right) \partial_j.
\end{align}
\end{proof}

\subsubsection{The Riemannian Hessian}\label{section_Riemannian_Hessian}

\begin{definition}[Riemannian Hessian --- Coordinate-free definition]
The \textbf{Riemannian Hessian} of $f$ is the symmetric $(0,2)$-tensor field defined as the covariant derivative of the gradient (or equivalently, the covariant derivative of the differential $df$). For vector fields $\b{X}, \b{Y}$, it is given by:
\begin{equation}
\boxed{
\operatorname{Hess} f(\b{X}, \b{Y}) = g(\nabla_{\b{X}} \operatorname{grad} f, \b{Y}) \overset{(\ref{equation_g_inner_product})}{=} \langle \nabla_{\b{X}} \operatorname{grad} f, \b{Y} \rangle.
}
\end{equation}
\end{definition}


\begin{proposition}[Riemannian Hessian --- Coordinate-based definition] \label{proposition_hessian_coordinates}
Let $f \in C^\infty(\mathcal{M})$ be a smooth function on a Riemannian manifold $(\mathcal{M}, g)$. 
Consider local coordinate system $\{x^i\}_{i=1}^n$ on a Riemannian manifold, with basis vector $\{\partial_i\}_{i=1}^n$ (see Eq. (\ref{equation_bases_for_vectors_2})). 
The components of the Riemannian Hessian, denoted by $f_{;ij}$, are given in local coordinates $\{x^i\}_{i=1}^n$ by:
\begin{equation}\label{equation_Riemannian_Hessian_components}
\boxed{
\begin{aligned}
f_{;ij} &:= \operatorname{Hess} f(\partial_i, \partial_j) = \nabla_j (\partial_i f) \\
&= \frac{\partial^2 f}{\partial x^i \partial x^j} - \Gamma^k_{ij} \frac{\partial f}{\partial x^k} \\
&\overset{(\ref{equation_partial_i})}{=} \partial_i \partial_j f - \Gamma^k_{ij} \partial_k f,
\end{aligned}
}
\end{equation}
where $\Gamma^k_{ij}$ are the Christoffel symbols of the Levi-Civita connection and $f_{;ij}$ uses the Ricci calculus notation (see Definition \ref{definition_comma_semicolon_notations}). 
The Riemannian Hessian is a $(0,2)$-tensor:
\begin{equation}
\boxed{
\begin{aligned}
\operatorname{Hess} f &= f_{;ij}\, \b{e}^i \otimes \b{e}^j = f_{;ij}\, dx^i \otimes dx^j \\
&\overset{(\ref{equation_Riemannian_Hessian_components})}{=} \left( \frac{\partial^2 f}{\partial x^i \partial x^j} - \Gamma^k_{ij} \frac{\partial f}{\partial x^k} \right) dx^i \otimes dx^j
\end{aligned}
}
\end{equation}

By Ricci calculus notation (see Definition \ref{definition_comma_semicolon_notations}), we can denote $f_{,i;j} = \nabla_j (\partial_i f)$ and $f_{,ij} = \frac{\partial^2 f}{\partial x^i \partial x^j}$, so Eq. (\ref{equation_Riemannian_Hessian_components}) can be stated as:
\begin{equation}\label{equation_Riemannian_Hessian_Ricci}
\boxed{
\begin{aligned}
f_{;ij} &= f_{,i;j} \\
&= f_{,ij} - \Gamma^k_{ji} f_{,k}
\end{aligned}
}
\end{equation}
\end{proposition}
\begin{proof}
The Riemannian Hessian is a symmetric $(0,2)$-tensor field defined as the covariant derivative of the differential $df$. By definition, for any two coordinate basis vector fields $\partial_i$ and $\partial_j$, the components are:
\begin{align}\label{equation_fij_Hess_nabla_partial_f}
f_{;ij} := \operatorname{Hess} f(\partial_i, \partial_j) = \nabla_j (\partial_i f).
\end{align}
Recall that the covariant derivative of a covector field $\omega$ with components $\omega_i$ is given by (see Eq. (\ref{equation_covariant_derivative_covariant})):
\begin{align}\label{equation_nabla_w_partial_w_Gamma_w}
\nabla_j \omega_i \overset{(\ref{equation_covariant_derivative_covariant})}{=} \partial_j \omega_i - \Gamma^k_{ji} \omega_k.
\end{align}
In the case of the Hessian, the covector field is the differential $\omega = df$, whose components are the partial derivatives $\omega_i = \partial_i f = \frac{\partial f}{\partial x^i}$. Substituting these components into the formula for the covariant derivative, we obtain:
\begin{align*}
f_{;ij} \overset{(\ref{equation_fij_Hess_nabla_partial_f})}{=} \nabla_j (\partial_i f) &\overset{(\ref{equation_nabla_w_partial_w_Gamma_w})}{=} \partial_j (\partial_i f) - \Gamma^k_{ji} (\partial_k f) \\
&\overset{(\ref{equation_partial_i})}{=} \frac{\partial }{\partial x^j} \left(\frac{\partial }{\partial x^i} f\right) - \Gamma^k_{ji} \left( \frac{\partial }{\partial x^k} f \right) \\
&= \frac{\partial^2 f}{\partial x^j \partial x^i} - \Gamma^k_{ji} \frac{\partial f}{\partial x^k}.
\end{align*}
Since the Levi-Civita connection is torsion-free, the Christoffel symbols are symmetric in their lower indices, i.e., $\Gamma^k_{ji} = \Gamma^k_{ij}$ (see Eq. (\ref{equation_symmtery_Christoffel})). Furthermore, for smooth functions, partial derivatives commute according to Clairaut's theorem ($\frac{\partial^2 f}{\partial x^j \partial x^i} = \frac{\partial^2 f}{\partial x^i \partial x^j}$). Thus, we arrive at:
\begin{align*}
f_{;ij} = \frac{\partial^2 f}{\partial x^i \partial x^j} - \Gamma^k_{ij} \frac{\partial f}{\partial x^k}.
\end{align*}
\end{proof}

Equation (\ref{equation_Riemannian_Hessian_components}) confirms that the Riemannian Hessian equals the Euclidean Hessian matrix of second partial derivatives minus a correction term involving the Christoffel symbols, which accounts for the geometry of the manifold.
Moreover, $f_{;ij} = f_{,i;j}$ in Eq. (\ref{equation_Riemannian_Hessian_Ricci}) identifies that the Hessian is the covariant derivative of a covector, and $f_{;ij} = f_{,ij} - \Gamma^k_{ji} f_{,k}$ in Eq. (\ref{equation_Riemannian_Hessian_Ricci}) is the explicit coordinate formula required for numerical implementation.

Note that, unlike the Euclidean case, the Riemannian Hessian accounts for the manifold's curvature through the Christoffel term $\Gamma^k_{ij} \partial_k f$. Recall that, according to Eq. (\ref{equation_parallel_transport_2}), parallel transport also contains Christoffel symbols. Therefore, the second-order information (Riemannian Hessian) is consistent with the parallel transport defined by the metric. 

\subsection{Definition of Exponential and Logarithm Maps}\label{section_exponential_logarithm_map}

In Euclidean optimization, the fundamental operations are vector addition and subtraction. To perform optimization on a Riemannian manifold $\mathcal{M}$, we require generalizations of these operations that respect the manifold's curvature and constraints. These generalizations are provided by the exponential and logarithm maps.

\subsubsection{The Exponential Map: Generalizing Addition}\label{section_exponential_map_generalizing_addition}

In Euclidean space, we update a position by adding a velocity vector: $\b{q} = \b{p} + \b{\Delta}$. On a manifold, simply adding a tangent vector to a point would generally result in a point outside the manifold. The exponential map solves this by ``wrapping'' the vector onto the manifold along a geodesic.

\begin{definition}[Exponential Map]\label{definition_exponential_map}
The \textbf{exponential map} at point $\b{p} \in \mathcal{M}$ is a mapping from the tangent space $T_{\b{p}}\mathcal{M}$ to the manifold $\mathcal{M}$:
\begin{equation}
\boxed{
\begin{aligned}
& \mathrm{Exp}_{\b{p}} : T_{\b{p}}\mathcal{M} \to \mathcal{M}, \\
&\mathrm{Exp}_{\b{p}}(\b{\xi}) = \b{q}, \quad \text{where } \b{\xi} \in T_{\b{p}}\mathcal{M} \text{ and } \b{q} \in \mathcal{M}.
\end{aligned}
}
\end{equation}
Specifically, $\mathrm{Exp}_{\b{p}}(\b{\xi})$ is the point reached by following the geodesic $\b{\gamma}(t)$ starting at $\b{p}$ with initial velocity $\b{\xi}$ for one unit of time ($t=1$).
In other words, the exponential map follows a geodesic from a point $\b{p}$ to another point $\b{q}$, where $t=0$ corresponds to the starting point $\b{p}$ and $t=1$ corresponds to the end point $\b{q}$.
\end{definition}

\subsubsection{The Logarithm Map: Generalizing Subtraction}\label{section_logarithm_map}

In Euclidean space, the difference between two points is a vector: $\b{\Delta} = \b{q} - \b{p}$. This operation maps two points to a vector. On a manifold, the ``difference'' between two points $\b{p}, \b{q} \in \mathcal{M}$ is represented by a tangent vector $\b{\xi}$ in the tangent space $T_{\b{p}}\mathcal{M}$.

\begin{definition}[Logarithm map]\label{definition_logarithm_map}
The \textbf{logarithm map} at point $\b{p}$ is a mapping from the manifold to the tangent space at $\b{p}$:
\begin{equation}
\boxed{
\begin{aligned}
& \mathrm{Log}_{\b{p}} : \mathcal{M} \to T_{\b{p}}\mathcal{M}, \\
& \mathrm{Log}_{\b{p}}(q) = \b{\xi}, \quad \text{where } \b{\xi} \in T_{\b{p}}\mathcal{M}.
\end{aligned}
}
\end{equation}
The vector $\b{\xi}$ is the initial velocity vector of the unique geodesic $\b{\gamma}(t)$ such that $\b{\gamma}(0) = \b{p}$ and $\b{\gamma}(1) = \b{q}$. Intuitively, the vector $\b{\xi}$ represents the initial velocity at $\b{p}$ of the geodesic that goes from $\b{p}$ to $\b{q}$.
\end{definition}

\begin{remark}[Relation of exponential and logarithm maps]
The exponential and logarithm maps are inverse operators, so the relation of exponential and logarithm maps is:
\begin{equation}
\boxed{
\begin{aligned}
& \mathrm{Exp}_{\b{p}}(\mathrm{Log}_{\b{p}}(\b{q})) = \b{q}, \\
& \mathrm{Log}_{\b{p}}(\mathrm{Exp}_{\b{p}}(\b{\xi})) = \b{\xi}.
\end{aligned}
}
\end{equation}
In other words, we have:
\begin{align}
\boxed{
q = \mathrm{Exp}_{\b{p}}(\b{\xi}) \iff \mathrm{Log}_{\b{p}}(\b{q}) = \b{\xi}. 
}
\end{align}
Both $\b{q} = \mathrm{Exp}_{\b{p}}(\b{\xi})$ and $\mathrm{Log}_{\b{p}}(q) = \b{\xi}$ represent a geodesic from a point $\b{p}$ to another point $\b{q}$, where $t=0$ corresponds to the starting point $\b{p}$ and $t=1$ corresponds to the end point $\b{q}$.
\end{remark}

\subsection{Numerical Calculation of Exponential Map}\label{section_exponential_map_numerical}

For numerical implementation, it is necessary to express the exponential maps in terms of local coordinates $\{x^i\}_{i=1}^n$. These derivations rely on the geodesic equations established in Section \ref{section_geodesics}.

\subsubsection{Numerical Calculation of Exponential Map by Ordinary Differential Equations}

Recall that the exponential map is a geodesic from a point $\b{p} \in \mathcal{M}$ to another point $\b{q} \in \mathcal{M}$, where $t=0$ corresponds to the starting point $\b{p}$ and $t=1$ corresponds to the end point $\b{q}$.

Given a point $\b{p} \in \mathcal{M}$ and a tangent vector $\b{\xi} \in T_{\b{p}}\mathcal{M}$, the exponential map $\mathrm{Exp}_{\b{p}}(\b{\xi})$ is found by solving an Initial Value Problem (IVP) in Ordinary Differential Equations (ODE). 
Let $\b{\gamma}(t)$ be a curve such that its coordinate representation is $\{x^k(t)\}_{k=1}^n$.
Let $\{p^k\}_{k=1}^n$ denote the local coordinates of chart around point $\b{p}$ and $\{q^k\}_{k=1}^n$ denote the local coordinates of chart around point $\b{q}$.

Consider that the exponential map is a geodesic. Recall from Eq. (\ref{equation_geodesic_ode}) that a geodesic satisfies:
\begin{equation}
\boxed{
\frac{d^2 x^k(t)}{dt^2} + \Gamma^k_{ij} \frac{dx^i(t)}{dt} \frac{dx^j(t)}{dt} = 0, \quad \forall k \in \{1, \dots, n\}.
}
\end{equation}
These are $n$ equations which are $n$ ODE problems. 
For IVP, we need initial conditions:
\begin{align}
&\b{\gamma}(0) = \b{p} \implies \boxed{x^k(0) = p^k}, \label{equation_initial_conditions_ode_exponential_1} \\
&\dot{\b{\gamma}}(0) = \b{\xi} \implies \boxed{\frac{dx^k(t)}{dt}\bigg|_{t=0} = \xi^k}, \label{equation_initial_conditions_ode_exponential_2}
\end{align}
to solve the ODE problems, where $\{x^k(t)\}_{k=1}^n$ are the coordinates of the points on the geodesic, $\{p^k\}_{k=1}^n$ are the local coordinates of chart around the point $\b{p}$, and $\{\xi^k\}_{k=1}^n$ are the components of the vector $\b{\xi}  \in T_{\b{p}}\mathcal{M}$ in the coordinate basis $\{\partial_i\}_{i=1}^n$.

After solving the $n$ ODE problems, we find the $n$ solutions $\{x^k(t)\}_{k=1}^n$ which are the coordinates of the points on the geodesic of exponential map, where each $t$ corresponds to a point on the geodesic, with $t=0$ for the starting point $\b{p}$ and $t=1$ for the end point $\b{q}$.
Therefore, the coordinates of the resulting point $\b{q} = \mathrm{Exp}_{\b{p}}(\b{\xi})$ are given by the solution of the ODE at $t=1$:
\begin{equation}\label{equation_solution_ode_exponential_1}
\boxed{
q^k = [\mathrm{Exp}_{\b{p}}(\b{\xi})]^k = x^k(1), \quad \forall k \in \{1, \dots, n\}.
}
\end{equation}
The coordinate $[q^1, \dots, q^n]^\top$, as the local coordinate of chart around point $\b{q} \in \mathcal{M}$ is the solution of the ODE problem for exponential map. The exponential map yields the point $\b{q}$. 

\subsubsection{Local Power Series Solution of Exponential Map}

In general Riemannian manifolds, the geodesic ODEs are non-linear and coupled, which typically precludes a closed-form solution. 
However, a local approximation of exponential map can be derived via a Taylor expansion of the coordinates $x^k(t)$ around $t=0$, as explained in the following. 

\begin{proposition}[Numerical calculation of exponential map]
The $k$-th component of the exponential map can be approximately calculated as:
\begin{align}\label{eq_exp_map_approx}
\boxed{
q^k = [\mathrm{Exp}_{\b{p}}(\b{\xi})]^k \approx p^k + \xi^k - \frac{1}{2} \Gamma^k_{ij}(\b{p}) \xi^i \xi^j,
} 
\end{align}
where $p^k$ is the $k$-th component of local coordinate of point $\b{p}$, and $q^k$ is the $k$-th component of local coordinate of point $\b{q}$, and $\xi^k$ is the $k$-th component of vector $\b{\xi}$, and $\Gamma^k_{ij}(\b{p})$ is the Christoffel symbol as a function of the local coordinates of point $\b{p}$.
The coordinate $[q^1, \dots, q^n]^\top$ for the point $\b{q} \in \mathcal{M}$ is the output of the exponential map.
\end{proposition}
\begin{proof}
Consider a local approximation can be derived via a Taylor expansion of the coordinates $x^k(t)$ around $t=0$:
\begin{equation}\label{equation_x_x_dx_dt_t_half_d2x_dt2_O_t3}
x^k(t) = x^k(0) + \frac{dx^k}{dt}\bigg|_{t=0} t + \frac{1}{2}\frac{d^2 x^k}{dt^2}\bigg|_{t=0} t^2 + \mathcal{O}(t^3).
\end{equation}

According to Eq. (\ref{equation_geodesic_ode}) for the equation of a geodesic, we have:
\begin{align*}
&\frac{d^2 x^k(t)}{dt^2} + \Gamma^k_{ij} \frac{dx^i(t)}{dt} \frac{dx^j(t)}{dt} = 0 \nonumber \\
&\implies
\frac{d^2 x^k(t)}{dt^2} = - \Gamma^k_{ij} \frac{dx^i(t)}{dt} \frac{dx^j(t)}{dt}.
\end{align*}
By using this to replace the second-order derivative in Eq. (\ref{equation_x_x_dx_dt_t_half_d2x_dt2_O_t3}), the Eq. (\ref{equation_x_x_dx_dt_t_half_d2x_dt2_O_t3}) becomes:
\begin{align*}
x^k(t) =\, &x^k(0) + \frac{dx^k}{dt}\bigg|_{t=0} t \\
&- \frac{1}{2} \Gamma^k_{ij} \frac{dx^i(t)}{dt} \frac{dx^j(t)}{dt} \bigg|_{t=0} t^2 + \mathcal{O}(t^3).
\end{align*}

By substituting the Eqs. (\ref{equation_initial_conditions_ode_exponential_1}) and (\ref{equation_initial_conditions_ode_exponential_2}), which are the initial conditions $x^k(0) = x^k$ and $\dot{x}^k(0) = \xi^k$, we obtain:
\begin{equation}\label{equation_xk_xk_xi_t_half_Gamma_xi_xi_t2}
x^k(t) \approx x^k + \xi^k t - \frac{1}{2} \Gamma^k_{ij}(\b{p}) \xi^i \xi^j t^2.
\end{equation}

According to Eq. (\ref{equation_solution_ode_exponential_1}), evaluating this expansion at $t=1$ provides the second-order coordinate approximation for the exponential map:
\begin{align*}
&y^k = [\mathrm{Exp}_{\b{p}}(\b{\xi})]^k = x^k(1) \nonumber \\
&\overset{(\ref{equation_xk_xk_xi_t_half_Gamma_xi_xi_t2})}{\implies} y^k = [\mathrm{Exp}_{\b{p}}(\b{\xi})]^k \approx x^k + \xi^k - \frac{1}{2} \Gamma^k_{ij}(\b{p}) \xi^i \xi^j.
\end{align*}
\end{proof}

This Eq. (\ref{eq_exp_map_approx}) can be used for numerical implementation of exponential map. This equations shows that the first-order approximation ($x^k + \xi^k$) ignores the manifold's curvature, while the second-order term $\frac{-1}{2} \Gamma^k_{ij}(\b{p}) \xi^i \xi^j$ introduces the geometry via the Christoffel symbols. This series serves as the foundational justification for the \textit{retraction map} discussed in Section \ref{section_retraction}.

\subsection{Numerical Calculation of Logarithm Map}

For numerical implementation, it is necessary to express the logarithm maps in terms of local coordinates $\{x^i\}_{i=1}^n$. These derivations rely on the geodesic equations established in Section \ref{section_geodesics}.

\subsubsection{Numerical Calculation of Logarithm Map by Ordinary Differential Equations}

Recall that the logarithm map is a geodesic from a point $\b{p} \in \mathcal{M}$ to another point $\b{q} \in \mathcal{M}$, where $t=0$ corresponds to the starting point $\b{p}$ and $t=1$ corresponds to the end point $\b{q}$.

Given two points $\b{p}, \b{q} \in \mathcal{M}$ and a tangent vector $\b{\xi} \in T_{\b{p}}\mathcal{M}$, the logarithm map $\mathrm{Log}_{\b{p}}(q)$ is found by solving a Boundary Value Problem (BVP) in Ordinary Differential Equations (ODE). 
Let $\b{\gamma}(t)$ be a curve such that its coordinate representation is $\{x^k(t)\}_{k=1}^n$.
Let $\{p^k\}_{k=1}^n$ denote the local coordinates of chart around point $\b{p}$ and $\{q^k\}_{k=1}^n$ denote the local coordinates of chart around point $\b{q}$.

Consider that the logarithm map is a geodesic. Recall from Eq. (\ref{equation_geodesic_ode}) that a geodesic satisfies:
\begin{equation}
\boxed{
\frac{d^2 x^k(t)}{dt^2} + \Gamma^k_{ij} \frac{dx^i(t)}{dt} \frac{dx^j(t)}{dt} = 0, \quad \forall k \in \{1, \dots, n\}.
}
\end{equation}
These are $n$ equations which are $n$ ODE problems. 
For BVP, we need boundary conditions:
\begin{align}
&\b{\gamma}(0) = \b{p} \implies \boxed{x^k(0) = p^k}, \\
&\b{\gamma}(1) = \b{q} \implies \boxed{x^k(1) = q^k},
\end{align}
to solve the ODE problems, where $\{x^k(t)\}_{k=1}^n$ are the coordinates of the points on the geodesic, $\{p^k\}_{k=1}^n$ are the coordinates of the point $\b{p}$, and $\{q^k\}_{k=1}^n$ are the coordinates of the point $\b{q}$.

After solving the $n$ ODE problems, we find the $n$ solutions $\{x^k(t)\}_{k=1}^n$ which are the coordinates of the points on the geodesic of logarithm map, where each $t$ corresponds to a point on the geodesic, with $t=0$ for the starting point $\b{p}$ and $t=1$ for the end point $\b{q}$.

Once the geodesic $\b{\gamma}(t)$ connecting $\b{x}$ and $\b{y}$ is determined, the components of the tangent vector $\b{\xi} \in T_{\b{p}}\mathcal{M}$ are the initial velocities:
\begin{equation}
\boxed{
\begin{aligned}
\xi^k &= [\mathrm{Log}_{\b{p}}(\b{q})]^k \\
&= \frac{dx^k(t)}{dt}\bigg|_{t=0} = \dot{x}^k(0), \quad \forall k \in \{1, \dots, n\}.
\end{aligned}
}
\end{equation}
The logarithm map yields the vector:
\begin{align}
\boxed{
\b{\xi} = \xi^k \partial_k,
}
\end{align}
where Einstein summation convention is used. 
The vector $\b{\xi} = \xi^k \partial_k$ is the solution of the ODE problem for logarithm map. 
This vector $\b{\xi}$, in the ``flat" tangent space $T_{\b{p}}\mathcal{M}$, points from $\b{p}$ toward $q$ along the shortest path (geodesic).

\subsubsection{Local Power Series Solution of Logarithm Map}

A local approximation of logarithm map can be derived, as explained in the following. 

\begin{proposition}[Numerical calculation of logarithm map]
The $k$-th component of the logarithm map can be approximately calculated as:
\begin{equation}\label{eq_log_map_approx}
\boxed{
\begin{aligned}
\xi^k &= [\mathrm{Log}_{\b{p}}(\b{q})]^k \\
&\approx
(q^k - p^k) + \frac{1}{2} \Gamma^k_{ij}(\b{p}) (q^i - p^i) (q^j - p^j),
\end{aligned}
}
\end{equation}
where $p^k$ is the $k$-th component of local coordinate of point $\b{p}$, and $q^k$ is the $k$-th component of local coordinate of point $\b{q}$, and $\xi^k$ denotes the $k$-th component of the tangent vector $\b{\xi}$, and $\Gamma^k_{ij}(\b{p})$ is the Christoffel symbol as a function of the local coordinates of point $\b{p}$.
The vector $\b{\xi} = \xi^k \partial_k$ is then obtained as the output of the logarithm map. 
\end{proposition}
\begin{proof}
The logarithm map $\mathrm{Log}_{\b{p}}(\b{q}) = \b{\xi}$ is the inverse of the exponential map. To find its coordinate-level power series, we seek the components $\xi^k$ in terms of the coordinate difference:
\begin{align}\label{equation_Deltax_y_x}
\Delta x^k = q^k - p^k.
\end{align}

We derive by \textit{method of undetermined coefficients}.
We begin by recalling the second-order Taylor expansion of the exponential map $\mathrm{Exp}_{\b{p}}(\b{\xi}) = q$ from Eq. \eqref{eq_exp_map_approx}:
\begin{align}
&q^k \overset{(\ref{eq_exp_map_approx})}{=} p^k + \xi^k - \frac{1}{2} \Gamma^k_{ij}(\b{p}) \xi^i \xi^j + \mathcal{O}(\Vert \Delta x \Vert^3) \nonumber \\
&\implies q^k - p^k = \xi^k - \frac{1}{2} \Gamma^k_{ij}(\b{p}) \xi^i \xi^j + \mathcal{O}(\Vert \Delta x \Vert^3) \nonumber \\
&\overset{(\ref{equation_Deltax_y_x})}{\implies} \Delta x^k = \xi^k - \frac{1}{2} \Gamma^k_{ij}(\b{p}) \xi^i \xi^j + \mathcal{O}(\Vert \Delta x \Vert^3). \label{equation_Deltax_forward_expansion}
\end{align}
To invert this mapping, we assume a solution for $\xi^k$ of the form:
\begin{equation}\label{eq_log_ansatz}
\xi^k = \Delta x^k + A^k_{ij} \Delta x^i \Delta x^j + \mathcal{O}(\Vert \Delta x \Vert^3).
\end{equation}
Substituting this Eq. (\ref{eq_log_ansatz}) back into the forward expansion (\ref{equation_Deltax_forward_expansion}), we obtain:
\begin{align*}
\Delta x^k &= \left( \Delta x^k + A^k_{ij} \Delta x^i \Delta x^j \right) \\
&- \frac{1}{2} \Gamma^k_{ij}(\b{p}) \left( \Delta x^i + \cdots \right) \left( \Delta x^j + \cdots \right) + \mathcal{O}(\Vert \Delta x \Vert^3).
\end{align*}
Keeping only terms up to second order, the equation simplifies to:
\begin{align}
&\Delta x^k = \Delta x^k + A^k_{ij} \Delta x^i \Delta x^j - \frac{1}{2} \Gamma^k_{ij}(\b{p}) \Delta x^i \Delta x^j \nonumber \\
&\implies \Delta x^k = \Delta x^k + \Big( A^k_{ij} - \frac{1}{2} \Gamma^k_{ij}(\b{p}) \Big) \Delta x^i \Delta x^j \nonumber \\
&\implies A^k_{ij} - \frac{1}{2} \Gamma^k_{ij}(\b{p}) = 0 \nonumber \\
&\implies A^k_{ij} = \frac{1}{2} \Gamma^k_{ij}(\b{p}). \label{equation_A_half_Gamma}
\end{align}
Therefore, Eq. (\ref{eq_log_ansatz}) becomes:
\begin{align*}
\xi^k &\overset{(\ref{eq_log_ansatz})}{\approx} \Delta x^k + A^k_{ij} \Delta x^i \Delta x^j \\
&\overset{(\ref{equation_A_half_Gamma})}{=} \Delta x^k + \frac{1}{2} \Gamma^k_{ij}(\b{p}) \Delta x^i \Delta x^j \\
&\overset{(\ref{equation_Deltax_y_x})}{=} (q^k - p^k) + \frac{1}{2} \Gamma^k_{ij}(\b{p}) (q^i - p^i) (q^j - p^j).
\end{align*}
This $\xi^k = [\mathrm{Log}_{\b{p}}(\b{q})]^k$ is the $k$-th component of the logarithm map.
\end{proof}

Equation (\ref{eq_log_map_approx}) highlights that while the first-order term $\Delta x^k = (q^k - p^k)$ coincides with the standard Euclidean subtraction used in flat-space optimization, the second-order term $\frac{1}{2} \Gamma^k_{ij}(\b{p}) (q^i - p^i) (q^j - p^j)$ explicitly incorporates the local geometry of the manifold through the Christoffel symbols at the base point $\b{p}$. This series is essential for establishing error bounds of numerical approximations and for understanding the relationship between the Riemannian distance $d(\b{p}, \b{q})$ and the coordinate differences.

\subsection{Retraction Map}\label{section_retraction}

In Riemannian optimization, the exponential map $\mathrm{Exp}_{\b{p}}(\b{\xi})$---introduced in Sections \ref{section_exponential_logarithm_map} and \ref{section_exponential_map_numerical}---is the canonical way to map a tangent vector $\b{\xi} \in T_{\b{p}}\mathcal{M}$ back to the manifold while following a geodesic. However, computing the exponential map often requires solving second-order Ordinary Differential Equations (ODEs), as discussed in Section \ref{section_exponential_map_numerical}, which can be computationally prohibitive for high-dimensional matrix manifolds. To address this, we use a more general and computationally efficient operator called a \textit{retraction map}, or just  \textit{retraction} in short.

\subsubsection{Definition and Properties of Retraction}

\begin{figure}[!h]
\centering
\includegraphics[width=3.2in]{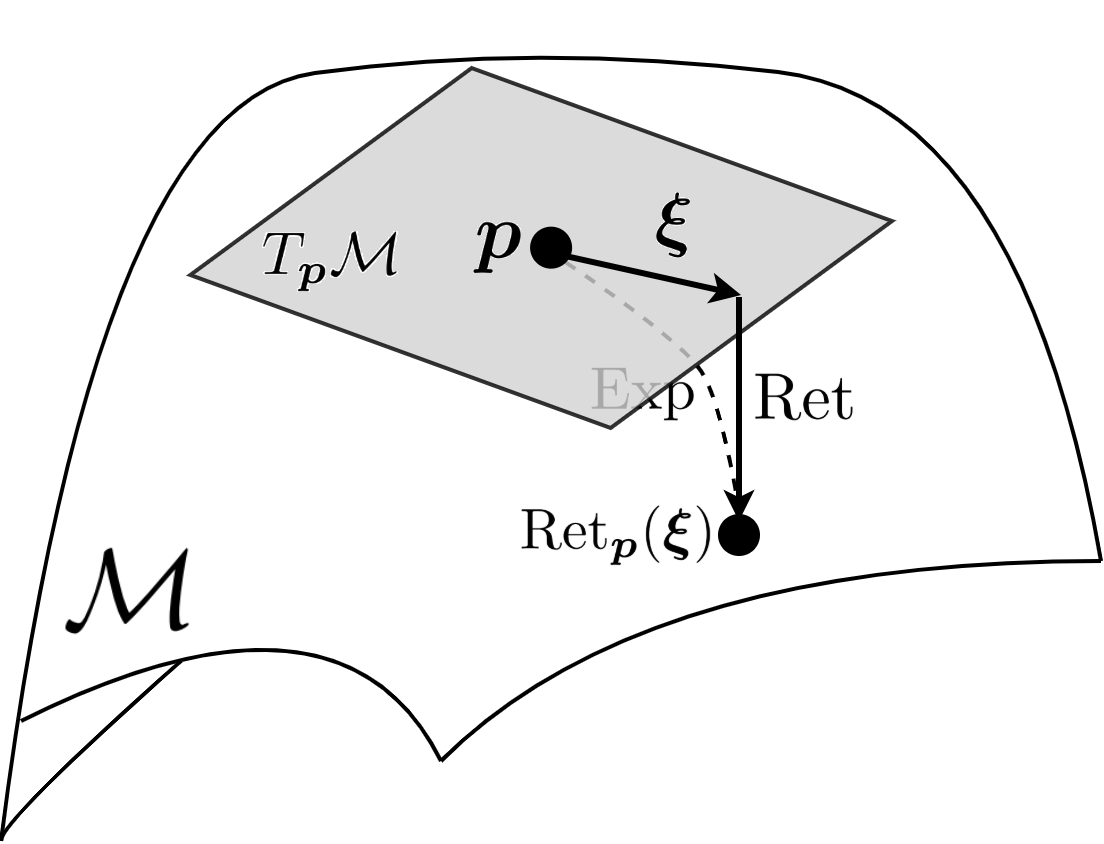}
\caption{Retraction map versus exponential map: the exponential map is depicted by a dashed curve on the manifold, while the retraction map moves from point along the tangent vector and then projects it back onto the manifold. The retraction map approximates the exponential map.}
\label{figure_retraction}
\end{figure}

A retraction is a smooth mapping that approximates the exponential map to at least first order, ensuring that the update step remains a valid descent direction on the manifold.
The retraction map and the exponential map are compared visually in Fig. \ref{figure_retraction}.
As this figure shows, exponential map is a geodesic on the manifold but the retraction map moves by a vector in the tangent space and then it is projected from the tangent space to the manifold. As you see in the figure, the results of exponential map and retraction map are roughly the same, so retraction map can be used instead of exponential map, with fewer calculations. 

\begin{definition}[Retraction Map]\label{definition_retraction}
A \textbf{retraction map} (also called \textbf{retraction} in short) at a point $\b{p} \in \mathcal{M}$ is a smooth mapping from the tangent space at the point to the manifold:
\begin{align}
\boxed{
\mathrm{Ret}_{\b{p}}(\b{\xi}): T_{\b{p}}\mathcal{M} \rightarrow \mathcal{M},
}
\end{align}
which satisfies the following two conditions:
\begin{enumerate}
    \item \textbf{Identity:} The retraction of zero vector in the tangent space at a point is the point itself:
    \begin{align}
    \boxed{
    \mathrm{Ret}_{\b{p}}(\b{0}_{\b{p}}) = \b{p},
    }
    \end{align}
    where $\b{0}_{\b{p}}$ is the zero vector in the tangent space $T_{\b{p}}\mathcal{M}$ at point $\b{p}$. This behavior is obvious according to Fig. \ref{figure_retraction}; assume the vector in the tangent space of a point is zero vector with zero length. Then, projection of vector onto the manifold is the point itself.  
    \item \textbf{Local Rigidity:} The differential of $\mathrm{Ret}_{\b{p}}$ at $\b{0}_{\b{p}}$ is the identity map on $T_{\b{p}}\mathcal{M}$. That is, for any $\b{\xi} \in T_{\b{p}}\mathcal{M}$:
    \begin{equation}
    \boxed{
    \left. \frac{d}{dt} \mathrm{Ret}_{\b{p}}(t\b{\xi}) \right|_{t=0} = \b{\xi}.
    }
    \end{equation}
\end{enumerate}
\end{definition}

\subsubsection{Numerical Calculation of Retraction}\label{section_retraction_numerical}

For numerical implementation, retractions are typically defined through algebraic projections or factorizations rather than ODEs. This allows for high-performance optimizers without a need to calculate ODEs.


\begin{proposition}[Numerical calculation of retraction]\label{proposition_retraction_numerical}
The $k$-th component of the retraction map of vector $\b{\xi}$ at point $\b{p}$ can be approximately calculated as:
\begin{equation}\label{equation_retraction_numerical}
\boxed{
[\mathrm{Ret}_{\b{p}}(\b{\xi})]^k \approx p^k + \xi^k,
}
\end{equation}
where $p^k$ is the $k$-th component of local coordinate of point $\b{p}$ and $\xi^k$ denotes the $k$-th component of the tangent vector $\b{\xi}$.
\end{proposition}
\begin{proof}
Just as the exponential map can be approximated by a power series, the retraction can be viewed as a first-order approximation of the geodesic path. Recall the second-order Taylor expansion of the exponential map from Eq. (\ref{eq_exp_map_approx}):
\begin{equation*}
[\mathrm{Exp}_{\b{p}}(\b{\xi})]^k \approx p^k + \xi^k - \frac{1}{2}\Gamma^k_{ij}(\b{p})\xi^i\xi^j.
\end{equation*}
A retraction $\mathrm{Ret}_{\b{p}}(\b{\xi})$ is numerically implemented such that it matches the first two terms of this expansion:
\begin{equation*}
[\mathrm{Ret}_{\b{p}}(\b{\xi})]^k = p^k + \xi^k + \mathcal{O}(\Vert\b{\xi}\Vert^2) \approx p^k + \xi^k. 
\end{equation*}
\end{proof}

As discussed in the proof of Proposition \ref{proposition_retraction_numerical}, retraction can be viewed as a first-order approximation of the geodesic path. It ignores the second-order term $\frac{-1}{2}\Gamma^k_{ij}(\b{p})\xi^i\xi^j$ of Eq. (\ref{eq_exp_map_approx}), involving the Christoffel symbols $\Gamma^k_{ij}$. By ignoring this term, we achieve a ``flat" update that is subsequently projected back onto the manifold.

According to Eq. (\ref{equation_retraction_numerical}), we have $[\mathrm{Ret}_{\b{p}}(\eta \b{\xi})]^k \approx p^k + \eta \xi^k$. 
This first-order approximation property ensures that for a small step size $\eta$, the retraction $\mathrm{Ret}_{\b{p}}(\eta \b{\xi})$ captures the local geometry of the manifold effectively enough for optimization convergence without the overhead of exact geodesic computation. In other words, as shown in Fig. \ref{figure_retraction}, when the $\eta \b{\xi}$ is small enough, the curvature of the manifold can be neglected so the $x^k + \eta \xi^k$ approximates exponential map along the tangent vector $\eta \xi^k$. 







\subsection{Vector Transport}\label{section_vector_transport}

In Riemannian optimization, many algorithms (e.g., Conjugate Gradient or Quasi-Newton methods) require comparing or combining tangent vectors located at different points on the manifold. While parallel transport $\mathcal{P}_{\gamma}$ along a geodesic is the canonical way to achieve this (see Section \ref{section_parallel_transport}), it is often computationally prohibitive as it requires solving a system of Ordinary Differential Equations (ODEs) involving Christoffel symbols. 

\textit{Vector transport} is a generalized, computationally efficient alternative that relaxes the requirements of parallel transport while remaining compatible with the chosen retraction.
So, vector transport, with fewer computation\footnote{Vector transport is simpler than parallel transport. There are even transformations which simplify vector transport to identity map, reducing the cost of computation significantly \cite{godaz2021vector}.}, can be used instead of parallel transport. 

\subsubsection{Definition of Vector Transport}\label{section_definition_vector_transport}

\begin{definition}[Vector transport]\label{definition_vector_transport}
A \textbf{vector transport} $\mathcal{T}$ on a manifold $\mathcal{M}$ is a smooth mapping:
\begin{equation}
\boxed{
\begin{aligned}
&\mathcal{T}: T_{\b{p}}\mathcal{M} \to T_{\mathrm{Ret}_{\b{p}}(\b{\eta})}\mathcal{M}, \\
&\mathcal{T}: \b{\xi} \mapsto \mathcal{T}_{\b{\eta}}(\b{\xi}), 
\end{aligned}
}
\end{equation}
that associates a vector $\b{\eta} \in T_{\b{p}}\mathcal{M}$ and a vector $\b{\xi} \in T_{\b{p}}\mathcal{M}$ to a vector $\mathcal{T}_{\b{\eta}}(\b{\xi}) \in T_{\mathrm{Ret}_{\b{p}}(\b{\eta})}\mathcal{M}$. 

It must satisfy the following properties (axioms):
\begin{enumerate}
\item \textbf{Well-defined mapping: } The map $\b{\xi} \mapsto \mathcal{T}_{\b{\eta}}(\b{\xi})$ is a well-defined mapping satisfying:
\begin{align}
\boxed{
\mathcal{T}_{\b{\eta}}(\b{\xi}) \in T_{\mathrm{Ret}_{\b{p}}(\b{\eta})}\mathcal{M}.
}
\end{align}
\item \textbf{Linearity in $\b{\xi}$:} The map $\b{\xi} \mapsto \mathcal{T}_{\b{\eta}}(\b{\xi})$ is a linear transformation from $T_{\b{p}}\mathcal{M}$ to $T_{\mathrm{Ret}_{\b{p}}(\b{\eta})}\mathcal{M}$:
\begin{align}
\boxed{
\mathcal{T}_{\b{\eta}}(a\b{\xi}_1 + b\b{\xi}_2) = a \mathcal{T}_{\b{\eta}}(\b{\xi}_1) + b \mathcal{T}_{\b{\eta}}(\b{\xi}_2),
}
\end{align}
where $a$ and $b$ are scalars and $\b{\eta}, \b{\xi}_1, \b{\xi}_2 \in T_{\b{x}}\mathcal{M}$.
\item \textbf{Consistency at zero:} For the zero vector $\b{0}_{\b{x}} \in T_{\b{x}}\mathcal{M}$, the transport is the identity: 
\begin{align}
\boxed{
\mathcal{T}_{\b{0}_{\b{p}}}(\b{\xi}) = \b{\xi}, \quad \forall \b{\xi} \in T_{\b{p}}\mathcal{M}. 
}
\end{align}
\end{enumerate}
\end{definition}

\begin{figure}[!h]
\centering
\includegraphics[width=3.2in]{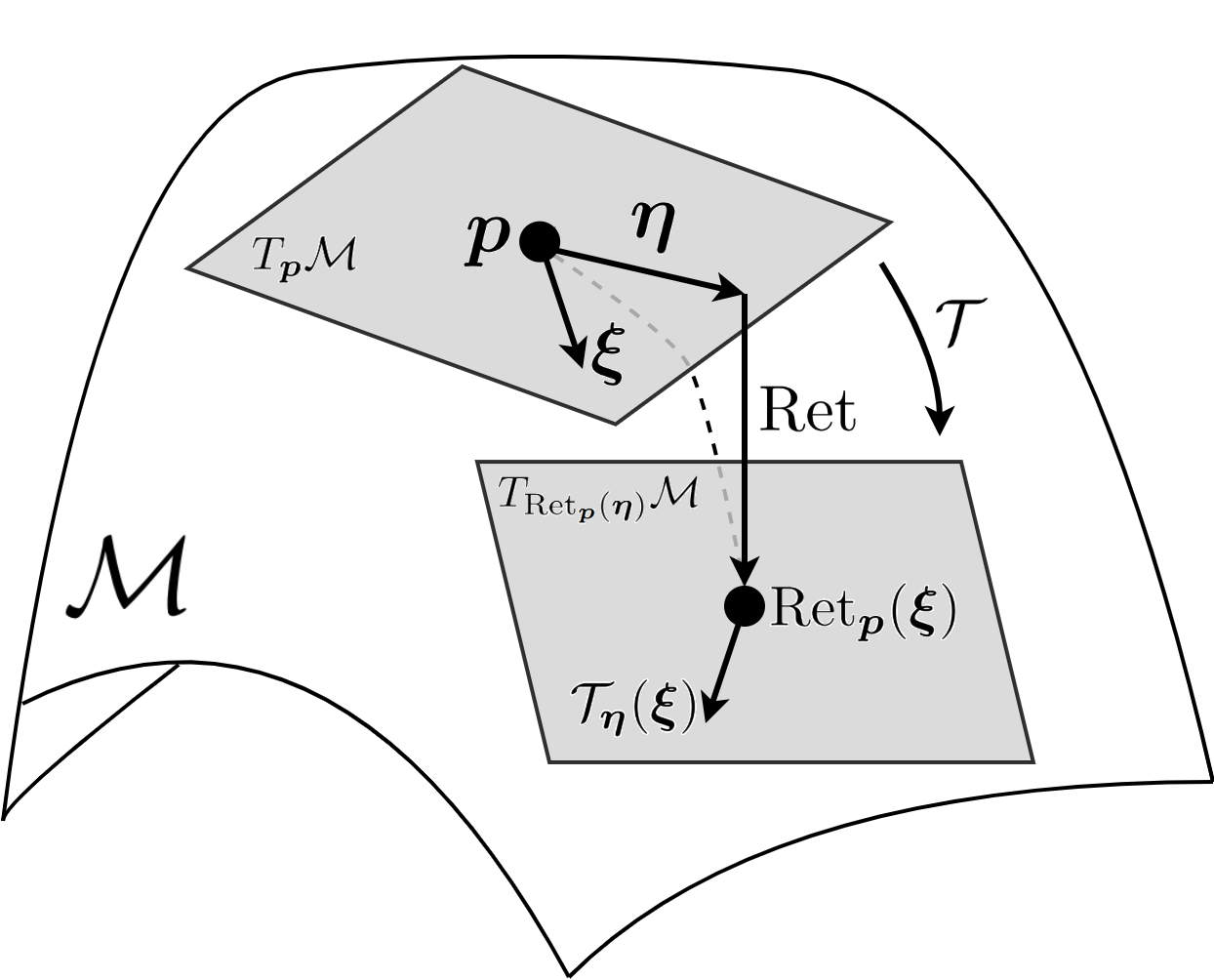}
\caption{Vector transport first performs the retraction $\mathrm{Ret}_{\b{p}}(\b{\eta})$ on the tangent vector $\b{\eta} \in T_{\b{p}}\mathcal{M}$. Thus, it obtains the point $\mathrm{Ret}_{\b{p}}(\b{\eta})$ on the manifold. This new point has a tangent space itself, namely $T_{\mathrm{Ret}_{\b{p}}(\b{\eta})}\mathcal{M}$.
Then, the vector transport $\mathcal{T}_{\b{\eta}}(\b{\xi})$ transforms the vector $\b{\xi}$ from tangent space $T_{\b{p}}\mathcal{M}$ to the tangent space $T_{\mathrm{Ret}_{\b{p}}(\b{\eta})}\mathcal{M}$ without change of its relative direction in the tangent space.}
\label{figure_vector_transport}
\end{figure}

Equivalently, vector transport $\mathcal{T}_{\b{\eta}}(\b{\xi})$ does the following steps:
\begin{enumerate}
\item Given a point $\b{x} \in \mathcal{M}$ and a vector $\b{\eta}$ in its tangent space, i.e., $\b{\eta} \in T_{\b{p}}\mathcal{M}$,
\item It calculates retraction $\mathrm{Ret}_{\b{p}}(\b{\eta})$ to obtain the point $\b{q} = \mathrm{Ret}_{\b{p}}(\b{\eta})$ on the manifold, i.e., $\b{q} \in \mathcal{M}$.
\item It considers a vector $\b{\xi}$ in the tangent space at point $\b{p}$, i.e., $\b{\xi} \in T_{\b{p}}\mathcal{M}$. 
\item It transforms the vector $\b{\xi}$ from the tangent space $T_{\b{p}}\mathcal{M}$ to the tangent space $T_{\b{q}}\mathcal{M}$.
\end{enumerate}
This procedure of vector transport is illustrated in Fig. \ref{figure_vector_transport}.
This shows that vector transport is a mapping from a tangent space to another tangent space on manifold.
You can see it as moving a tangent vector in a tangent space to the corresponding tangent vector in another tangent space. 

\subsubsection{Numerical Calculation of Vector Transport by Differentiated Retraction}\label{section_differentiated_retraction}

The most natural construction of a vector transport is the \textit{differentiated retraction}, which is directly induced by the retraction map.
In other words, we can choose a retraction $\mathrm{Ret}$ and then differentiate it to get a vector transport $\mathrm{D}\mathrm{Ret}$.

Vector transport is the general concept or abstract map, defined by the axioms in Definition \ref{definition_vector_transport}.
Differentiated retraction is one specific construction of such a vector transport. 
Differentiated retraction is not a different concept from vector transport; rather, it is a particular construction of vector transport obtained by differentiating a chosen retraction map.

\begin{definition}[Differentiated retraction]\label{definition_differentiated_retraction}
Let $\mathrm{Ret}_{\b{p}} : T_{\b{p}}\mathcal{M} \to \mathcal{M}$ be a retraction. The \textbf{differentiated retraction} defines a vector transport:
\begin{equation}\label{equation_differential_retraction}
\boxed{
\mathcal{T}_{\b{\eta}}(\b{\xi}) := \mathrm{D}\mathrm{Ret}_{\b{p}}(\b{\eta})[\b{\xi}],
}
\end{equation}
where $\b{\xi} \in T_{\b{p}}\mathcal{M}$, and we use the natural identification\footnote{Since $T_{\b{p}}\mathcal{M}$ is already a vector space, its tangent space at any point is canonically the same vector space.} $\mathcal{T}_{\b{\eta}} (T_{\b{p}}\mathcal{M}) \cong T_{\b{p}}\mathcal{M}$. Hence,
\[
\mathcal{T}_{\b{\eta}} : T_{\b{p}}\mathcal{M} \to T_{\mathrm{Ret}_{\b{p}}(\b{\eta})}\mathcal{M}.
\]

In a local coordinate chart $\varphi$ around $\mathrm{Ret}_{\b{p}}(\b{\eta})$, the $i$-th component of the transported vector is given by:
\begin{equation}\label{equation_differential_retraction_coordinates}
\boxed{
\left(\mathcal{T}_{\b{\eta}}(\b{\xi})\right)^i
=
\sum_{j=1}^d
\frac{\partial (\mathrm{Ret}_{\b{p}}(\b{\eta}))^i}{\partial \eta^j}
\, \xi^j.
}
\end{equation}
\end{definition}

Equation (\ref{equation_differential_retraction_coordinates}) is particularly useful for numerical implementation because the derivative of the retraction is often available via automatic differentiation or simple matrix algebra.

\begin{proposition}[Differentiated retraction defines a vector transport]\label{proposition_differentiated_retraction_vector_transport}
Let $\mathrm{Ret}_{\b{p}} : T_{\b{p}}\mathcal{M} \to \mathcal{M}$ be a retraction. 
The differentiated retraction $\mathcal{T}_{\b{\eta}}(\b{\xi})$ is a vector transport, i.e.,
\[
\mathcal{T}_{\b{\eta}} : T_{\b{p}}\mathcal{M} \to T_{\mathrm{Ret}_{\b{p}}(\b{\eta})}\mathcal{M}
\]
is well-defined, linear in $\b{\xi}$, and satisfies $\mathcal{T}_{\b{0}}(\b{\xi}) = \b{\xi}$.
\end{proposition}

\begin{proof}
Since $\mathrm{Ret}_{\b{p}} : T_{\b{p}}\mathcal{M} \to \mathcal{M}$ is smooth, its differential at $\b{\eta}$ is:
\[
\mathrm{D}\mathrm{Ret}_{\b{p}}(\b{\eta}) :
T_{\b{\eta}}(T_{\b{p}}\mathcal{M}) \to
T_{\mathrm{Ret}_{\b{p}}(\b{\eta})}\mathcal{M}.
\]

\textbf{Checking being well-defined}:
Since $T_{\b{p}}\mathcal{M}$ is already a vector space, its tangent space at any point is canonically the same vector space. Thus, we have canonical identification $T_{\b{\eta}}(T_{\b{p}}\mathcal{M}) \cong T_{\b{p}}\mathcal{M}$. 
Using the canonical identification, we obtain:
\[
\mathcal{T}_{\b{\eta}}(\b{\xi}) \in T_{\mathrm{Ret}_{\b{p}}(\b{\eta})}\mathcal{M},
\]
so the map is well-defined.

\textbf{Checking linearity in } $\b{\xi}$:
Linearity in $\b{\xi}$ follows from linearity of the differential:
\begin{align*}
\mathcal{T}_{\b{\eta}}(a\b{\xi}_1 + b\b{\xi}_2)
&=
\mathrm{D}\mathrm{Ret}_{\b{p}}(\b{\eta})[a\b{\xi}_1 + b\b{\xi}_2]
\\
&=
a\,\mathcal{T}_{\b{\eta}}(\b{\xi}_1)
+
b\,\mathcal{T}_{\b{\eta}}(\b{\xi}_2).
\end{align*}

\textbf{Checking consistency at zero}:
Finally, by the defining property of a retraction, we have:
\[
\mathrm{Ret}_{\b{p}}(\b{0}) = \b{x}
\quad \text{and} \quad
\mathrm{D}\mathrm{Ret}_{\b{p}}(\b{0}) = \mathrm{Id}_{T_{\b{p}}\mathcal{M}},
\]
hence:
\[
\mathcal{T}_{\b{0}}(\b{\xi})
=
\mathrm{D}\mathrm{Ret}_{\b{p}}(\b{0})[\b{\xi}]
=
\b{\xi}.
\]
Therefore, $\mathcal{T}$ satisfies the axioms of a vector transport. Therefore, it is a vector transport. 
\end{proof}



\subsubsection{Comparison of Vector Transport and Parallel Transport}

Parallel transport can be viewed as a specific type of vector transport where the retraction is the exponential map ($\mathrm{Ret} = \mathrm{Exp}$) and the transport is defined by the Levi-Civita connection. In the coordinate-level detail required for implementation:
\begin{itemize}
\item \textbf{Parallel Transport:} According to Eq. (\ref{equation_parallel_transport_2}), parallel transport of the vector $\b{\xi} = \xi^j \partial_j$ is calculated as:
\begin{align}
\boxed{
\frac{d\xi^i}{dt} + \Gamma^i_{jk} \xi^j \frac{dx^k}{dt} = 0.
}
\end{align}
\item \textbf{Vector Transport:} According to Eq. (\ref{equation_differential_retraction}), it requires only the evaluation of $\mathcal{T}_{\b{\eta}}(\b{\xi}) = \mathrm{D}\mathrm{Ret}_{\b{p}}(\b{\eta})[\b{\xi}]$, using Eq. (\ref{equation_differential_retraction_coordinates}), bypassing the explicit calculation of Christoffel symbols $\Gamma^i_{jk}$.
\end{itemize}


\subsection{First-order Riemannian Optimization: Riemannian Gradient Descent}
\label{section_riemannian_gradient_descent}

\subsubsection{The Euclidean Gradient Descent Update}

In Euclidean space, the first-order update for minimizing a cost function $f: \mathbb{R}^n \to \mathbb{R}$ is the standard (Euclidean) Gradient Descent (GD) algorithm:
\begin{equation}\label{section_Euclidean_gradient_descent}
\b{x}_{\nu+1} = \b{x}_\nu - \lambda \nabla f(\b{x}_\nu),
\end{equation}
where $\nu$ denotes the iteration index, $\lambda > 0$ is the learning rate, and $\nabla f(\b{x}_\nu)$ is the Euclidean gradient. As discussed in Sections \ref{section_optimization_path} and \ref{section_exponential_map_generalizing_addition}, this additive update is generally invalid on a Riemannian manifold $\mathcal{M}$ because the sum $\b{x}_\nu - \lambda \nabla f(\b{x}_\nu)$ or $\b{x}_\nu + (- \lambda \nabla f(\b{x}_\nu))$ typically leaves the manifold. To generalize this update, we must utilize the geometric foundations of the tangent space and the exponential (or retraction) map.

\subsubsection{The Riemannian Gradient Descent (RGD) Update}\label{section_RGD_update}

The Riemannian Gradient Descent (RGD) replaces the Euclidean subtraction with a movement along the manifold in the direction of the steepest descent using Riemannian gradient (see Section \ref{section_Riemannian_gradient}). 
Let $\b{p}_\nu \in \mathcal{M}$ be the point of optimization path on the manifold at the iteration $\nu$, and let $\operatorname{grad} f(\b{p}_\nu) \in T_{\b{p}_\nu}\mathcal{M}$ be the Riemannian gradient of the cost function $f$ at point $\b{p}_\nu$. The updated point $\b{p}_{\nu+1}$ is obtained by applying the exponential map to the negative gradient scaled by the learning rate \cite{bonnabel2013stochastic}:
\begin{equation}\label{equation_RGD}
\boxed{
\b{p}_{\nu+1} = \mathrm{Exp}_{\b{p}_\nu}\big(\!-\lambda\, \operatorname{grad} f(\b{p}_\nu)\big),
}
\end{equation}
where $\operatorname{grad} f(\b{p}_\nu)$ is the Riemannian gradient at point $\b{p}$ (see Section \ref{section_Riemannian_gradient}), and addition of $- \lambda \nabla f(\b{x}_\nu)$ in Eq. (\ref{section_Euclidean_gradient_descent}) has been replaced with the exponential map in Eq. (\ref{equation_RGD}); see Section \ref{section_exponential_map_generalizing_addition}.

\subsubsection{Numerical Implementation of RGD using Retractions}\label{section_numerical_RGD_retraction}

For computational efficiency, the exponential map is often replaced by a retraction map $\mathrm{Ret}_{\b{p}_\nu}$ (see Section \ref{section_retraction}), which provides a first-order approximation of the geodesic path while ensuring the iterate remains on the manifold \cite{bonnabel2013stochastic}:
\begin{equation}\label{equation_RGD_Ret}
\boxed{
\b{p}_{\nu+1} = \mathrm{Ret}_{\b{p}_\nu}\big(\!-\lambda\, \operatorname{grad} f(\b{p}_\nu)\big).
}
\end{equation}
Following the coordinate-level detail required for numerical implementation established in Section \ref{section_retraction_numerical}, the $k$-th component of the updated iterate $\b{p}_{\nu+1}$ can be approximated as:
\begin{equation}\label{equation_RGD_approximate}
\boxed{
p^k_{\nu+1} \approx p_\nu^k - \lambda\, [\operatorname{grad} f(\b{p}_\nu)]^k, 
}
\end{equation}
where $p^k_{\nu+1}$ and $p_\nu^k$ denote the $k$-th component of $\b{p}_{\nu+1}$ and $\b{p}_\nu$, respectively, in local coordinates. 

This first-order approximation in Eq. (\ref{equation_RGD_approximate}) effectively ignores the manifold's curvature during the step itself; however, the retraction in Eq. (\ref{equation_RGD_Ret}) is more accurate by ensuring that the final result is projected back onto $\mathcal{M}$.

\subsubsection{Convergence and the Role of Curvature}

The convergence properties of RGD are intrinsically linked to the manifold's curvature, particularly the sectional curvature (see Section \ref{section_sectional_curvature}). As noted in Section \ref{section_ricci_flow_intuition_for_riemannian_optimization}, adjusting the Riemannian metric $g_{ij}$ can be used to optimize the convergence rate by ``smoothing" the underlying geometry toward a more uniform curvature state, e.g., using Ricci flow.

\subsubsection{Riemannian Stochastic Gradient Descent (RSGD)}

In large-scale machine learning applications where the cost function is a sum of many terms, we use Riemannian Stochastic Gradient Descent (RSGD). The update rule follows the same logic but uses a stochastic estimate of the gradient:
\begin{equation}
\boxed{
\b{p}_{\nu+1} = \mathrm{Ret}_{\b{p}}\big(\!-\lambda\, \widehat{\mathrm{grad}} f(\b{p}_\nu)\big).
}
\end{equation}

To calculate the stochastic Riemannian gradient $\widehat{\mathrm{grad}} f(\b{p}_\nu)$, we first compute the Euclidean stochastic gradient and then project it onto the tangent space of the manifold at point $p_\nu$ using the metric tensor.
Based on the derivations in Section \ref{section_Riemannian_gradient} and according to Eq. (\ref{equation_Riemannian_gradient_components}), the components of the Riemannian gradient are obtained by raising the index of the partial derivatives using the inverse metric tensor $g^{ij}$.
For a stochastic estimate, we have:
\begin{equation}
\boxed{
\big(\widehat{\mathrm{grad}} f(\b{p}_\nu)\big)^j = g^{ij}(\b{p}_\nu) \frac{\partial \hat{f}}{\partial x^i}(\b{p}_\nu),
}
\end{equation}
where $\frac{\partial \hat{f}}{\partial x^i}$ is the stochastic estimate of the partial derivative, typically calculated using a mini-batch of data as explained in the following. 

In the context of stochastic optimization, let the objective function $f(\b{p})$ be defined as the expectation of a loss function $\ell$ over a data distribution $\mathcal{D}$, or a finite sum over a dataset $\{ \b{\zeta}_n \}_{n=1}^N$:
\begin{equation}
f(\b{p}) = \mathbb{E}_{\b{\zeta} \sim \mathcal{D}} [\ell(\b{p}; \b{\zeta})] \approx \frac{1}{N} \sum_{n=1}^N \ell(\b{p}; \b{\zeta}_n).
\end{equation}
A stochastic approximation $\hat{f}$ is typically constructed using a mini-batch $\mathcal{B}$ of size $m$, where $\mathcal{B} \subset \{1, \dots, N\}$:
\begin{equation}
\hat{f}(\b{p}) = \frac{1}{m} \sum_{n \in \mathcal{B}} \ell(\b{p}; \b{\zeta}_n).
\end{equation}

The stochastic partial derivatives with respect to the local coordinates $\{x^i\}_{i=1}^n$ at the current iterate $\b{p}_\nu$ are computed as:
\begin{equation}
\boxed{
\frac{\partial \hat{f}}{\partial x^i} (\b{p}_\nu) = \frac{1}{m} \sum_{n \in \mathcal{B}} \frac{\partial \ell(x; \b{\zeta}_n)}{\partial x^i} \bigg|_{x = \varphi(p_\nu)},
}
\end{equation}
where $\varphi(p_\nu)$ represents the coordinate map of the point $p_\nu$ in a local chart. 


Note that, in computational frameworks utilizing automatic differentiation, the terms $\frac{\partial \hat{f}}{\partial x^i}$ are computed via backpropagation. 

\subsection{Second-order Riemannian Optimization: Riemannian Newton’s Method}\label{section_riemannian_second_order_optimization}

While first-order methods like Riemannian Gradient Descent utilize only the steepest descent direction, second-order methods incorporate the curvature of the manifold $\mathcal{M}$ and the cost function $f$ through the Riemannian Hessian (see Section \ref{section_Riemannian_Hessian}). This approach yields quadratic convergence in the neighborhood of a local optimum.

\subsubsection{The Riemannian Newton Equation}
In Euclidean space, the Newton step is defined by:
\begin{align}
\nabla^2 f(\b{x}) \Delta \b{x} = -\nabla f(\b{x}),
\end{align}
where $\nabla f(\b{x})$ and $\nabla^2 f(\b{x})$ are the Euclidean gradient and Hessian, respectively. 
On a Riemannian manifold, we replace the Euclidean gradient and Hessian with their Riemannian counterparts.

\begin{definition}[Riemannian Newton step]
At a point $\b{p} \in \mathcal{M}$, the \textbf{Riemannian Newton step} $\b{\eta} \in T_{\b{p}}\mathcal{M}$ is defined as the solution to the linear system:
\begin{equation}\label{equation_riemannian_newton_step}
\boxed{
\operatorname{Hess} f(\b{p})[\b{\eta}] = -\operatorname{grad} f(\b{p}),
}
\end{equation}
where $\operatorname{Hess} f$ is the Riemannian Hessian operator.
\end{definition}


\subsubsection{Coordinate-based Derivation of Riemannian Newton Equation}

To bridge the gap to numerical implementation, we express the Newton equation in local coordinates $\{x^i\}_{i=1}^n$. 

\begin{proposition}[Coordinate-based equation of Riemannian Newton equation]
Let $\b{\eta} = \eta^i \partial_i$ be the Riemannian Newton step where $\{\partial_i\}_{i=1}^n$ are the basis vectors. 
The components of the Riemannian Newton step equation, in local coordinates $\{x^i\}_{i=1}^n$, are calculated as:
\begin{equation}\label{equation_Riemannian_Newwton_equation_coordinates}
\boxed{
\left( \frac{\partial^2 f}{\partial x^i \partial x^j} - \Gamma^k_{ij} \frac{\partial f}{\partial x^k} \right) \eta^i = -\frac{\partial f}{\partial x^j},
}
\end{equation}
or according to Eq. (\ref{equation_partial_i}), it can be stated as:
\begin{equation}
\boxed{
\left( \partial_i \partial_j f - \Gamma^k_{ij} \partial_k f \right) \eta^i = -\partial_j f.
}
\end{equation}
\end{proposition}
\begin{proof}
Let $\b{\eta} = \eta^i \partial_i$ be the Newton step. 
Recall the component-wise definition of the Riemannian gradient from Eq. \eqref{equation_Riemannian_gradient_components}:
\begin{equation}\label{equation_inproof_gradf_j_gij_partiali_f}
(\operatorname{grad} f)^j = g^{ij} \partial_i f.
\end{equation}
Also, recall the component-wise definition of the Riemannian Hessian from Eq. \eqref{equation_Riemannian_Hessian_components}:
\begin{align}\label{equation_inproof_fij_partial2f_partialx_partialx_minus_Gammaij_partialm_f}
f_{;ij} = \frac{\partial^2 f}{\partial x^i \partial x^j} - \Gamma^m_{ij} \partial_m f.
\end{align}

The Riemannian Newton equation, i.e., Eq. (\ref{equation_riemannian_newton_step}), is a vector equality. In components, the left-hand side of Eq. (\ref{equation_riemannian_newton_step}) involves the contraction of the $(0,2)$ Hessian tensor $f_{;ij}$ with the vector $\eta^i$, which must then be raised by the inverse metric $g^{jk}$ (see Lemma \ref{lemma_index_raising_by_metric}) to match the contravariant nature of the gradient:
\begin{equation}\label{equation_gjk_Hessian_eta_i_minus_gradf_k}
\boxed{
g^{jk} f_{;ij} \eta^i = -(\operatorname{grad} f)^k.
}
\end{equation}
Substituting Eq. (\ref{equation_inproof_gradf_j_gij_partiali_f}) into Eq. (\ref{equation_gjk_Hessian_eta_i_minus_gradf_k}) gives:
\begin{align*}
g^{jk} f_{;ij} \eta^i = -g^{jk} \partial_j f.
\end{align*}
Multiplying the sides by $g_{kl}$ gives:
\begin{align*}
&g_{kl} g^{jk} f_{;ij} \eta^i = - g_{kl} g^{jk} \partial_j f \\
&\overset{(\ref{equation_metric_inverse})}{\implies} \delta_l^j f_{;ij} \eta^i = - \delta_l^j \partial_j f \\
&\overset{(\ref{equation_index_substitution_delta})}{\implies} f_{;il}\, \eta^i = - \partial_l f \\
&\overset{(a)}{\implies} f_{;ij}\, \eta^i = - \partial_j f \\
&\overset{(b)}{\implies} 
\left( \frac{\partial^2 f}{\partial x^i \partial x^j} - \Gamma^m_{ij} \frac{\partial f}{\partial x^m} \right) \eta^i = -\frac{\partial f}{\partial x^j} \\
&\overset{(c)}{\implies} 
\left( \frac{\partial^2 f}{\partial x^i \partial x^j} - \Gamma^k_{ij} \frac{\partial f}{\partial x^k} \right) \eta^i = -\frac{\partial f}{\partial x^j},
\end{align*}
where $(a)$ is because of renaming the dummy index $l \to j$, and $(b)$ is because of Eq. (\ref{equation_inproof_fij_partial2f_partialx_partialx_minus_Gammaij_partialm_f}) for the left-hand side and Eq. (\ref{equation_partial_i}) for the right-hand side, and $(c)$ is because of renaming the dummy index $m \to k$. 
\end{proof}

Solving Eq. (\ref{equation_Riemannian_Newwton_equation_coordinates}), as a linear system of equations, for $\eta^i$ provides the search direction. The update is then performed via a retraction $\mathrm{Ret}_{\b{p}}: T_{\b{p}_\nu}\mathcal{M} \to \mathcal{M}$:
\begin{equation}
\boxed{
\b{p}_{\nu+1} = \mathrm{Ret}_{\b{p}_\nu}(\b{\eta}_\nu),
}
\end{equation}
where $\b{\eta}_\nu$ denotes the vector of Riemannian Newton step $\b{\eta} = \eta^i \partial_i$ at iteration $\nu$. 

\subsubsection{Riemannian Quasi-Newton Methods: RBFGS}

In Euclidean second-order optimization, the \textit{quasi-Newton methods} approximate the inverse Hessian matrix to facilitate computation of the Hessian matrix. 
One of the effective quasi-Newton methods is the \textit{BFGS (Broyden-Fletcher-Goldfarb-Shanno)} algorithm, which was independently and almost simultaneously proposed in 1970 by four researchers Charles George Broyden \cite{broyden1970convergence,broyden1970convergence2}, Roger Fletcher \cite{fletcher1970new}, Donald Goldfarb \cite{goldfarb1970family}, and David Shanno \cite{shanno1970conditioning}. 

Likewise, in Riemannian second-order optimization, calculating the full Hessian and Christoffel symbols $\Gamma^k_{ij}$ is often computationally prohibitive. \textit{Riemannian BFGS (RBFGS)} approximates the Hessian using gradient information from successive iterates. 

A fundamental challenge in RBFGS is that the gradients $\operatorname{grad} f(\b{p}_\nu)$ and $\operatorname{grad} f(\b{p}_{\nu+1})$ reside in different tangent spaces. To perform a valid update, we must use a vector transport $\mathcal{T}_{\b{\eta}_\nu}: T_{\b{p}_\nu}\mathcal{M} \to T_{\b{p}_{\nu+1}}\mathcal{M}$.

Let:
\begin{align}
& \boxed{ \b{s}_\nu := \mathcal{T}_{\b{\eta}_\nu}(\b{\eta}_\nu) }, \\
& \boxed{ \b{y}_\nu := \operatorname{grad} f(\b{p}_{\nu+1}) - \mathcal{T}_{\b{\eta}_\nu}(\operatorname{grad} f(\b{p}_\nu)) }.
\end{align}
Let $\mathcal{B}_{\nu}$ denote the approximate Hessian operator at iteration $\nu$.
The approximate Hessian operator at iteration $\mathcal{B}_{\nu+1}$ is updated as:
\begin{equation}\label{equation_RBFGS}
\boxed{
\mathcal{B}_{\nu+1} = \widetilde{\mathcal{B}}_\nu - \frac{\widetilde{\mathcal{B}}_\nu \b{s}_\nu \b{s}_\nu^\top \widetilde{\mathcal{B}}_\nu}{\b{s}_\nu^\top \widetilde{\mathcal{B}}_\nu \b{s}_\nu} + \frac{\b{y}_\nu \b{y}_\nu^\top}{\b{y}_\nu^\top \b{s}_\nu},
}
\end{equation}
where $\widetilde{\mathcal{B}}_\nu = \mathcal{T}_{\b{\eta}_\nu} \mathcal{B}_\nu \mathcal{T}_{\b{\eta}_\nu}^{-1}$ is the transported approximation from the previous tangent space. In practice, a ``cautious" update is applied to ensure $\mathcal{B}$ remains positive definite, updating only if $\b{y}_\nu^\top \b{s}_\nu > \epsilon \Vert \b{s}_\nu \Vert^2$.
Equation (\ref{equation_RBFGS}) is similar to its counterpart in Euclidean BFGS; for its proof, refer to the proof of update equation in Euclidean BFGS.

As vector transport is computationally expensive in RBFGS, \textit{cautious RBFGS} was proposed \cite{huang2016riemannian} which ignores the curvature condition in the Wolfe conditions \cite{wolfe1969convergence} and only checks the Armijo condition \cite{armijo1966minimization} (see our other paper \cite{ghojogh2021kkt} for explanation of Wolfe and Armijo conditions). Since the curvature condition guarantees that the approximation of Hessian remains positive definite, it compensates by checking a cautious condition \cite{li2001global} before updating the approximation of Hessian.
This cautious RBFGS has been used in the Manopt optimization toolbox \cite{boumal2014manopt} (see Section \ref{section_important_toolboxes_riemannian_optimization} for the Manopt optimization toolbox).

\subsubsection{Riemannian Quasi-Newton Methods: RLBFGS}

It is noteworthy that the quasi-Newton methods, including BFGS, approximate the inverse Hessian matrix by a dense $(n \times n)$ matrix, where $n$ is the dimensionality. For large $n$, storing this matrix is very memory-consuming.
Therefore, the Euclidean \textit{Limited-memory BFGS (LBFGS)}---which uses much less memory than BFGS---was proposed, by Nocedal et al. in the 1980s \cite{nocedal1980updating,liu1989limited}. 
For better understanding of Euclidean LBFGS, refer to Nocedal's book {\citep[Chapter 6]{nocedal2006numerical}}.

The Riemannian counterparts of Euclidean BFGS and Euclidean LBFGS are \textit{Riemannian BFGS}---proposed by Qi in 2010 \cite{qi2010riemannian}---and \textit{Riemannian LBFGS}---proposed by Hosseini and Sra in 2020 \cite{hosseini2020alternative}---respectively. 
The convergence of RBFGS has been analyzed and proved by Ring \cite{ring2012optimization} and Huang \cite{huang2015broyden}.
The properties of RBFGS have been analyzed by Seibert \cite{seibert2013properties}. 
Some other direct extensions of Euclidean BFGS to Riemannian spaces are provided in \citep[Chapter 7]{ji2007optimization}.
A vector-transport-free version of RLBFGS has also been proposed in \cite{godaz2021vector}.

\section{Important Riemannian Matrix Manifolds}
\label{section_important_Riemannian_matrix_manifolds}

Building upon the general optimization framework established in Section \ref{section_manifold_valued_optimization}, we now apply these concepts to specific matrix manifolds. These manifolds are frequently encountered in computational science, where constraints such as orthogonality or positive definiteness are essential. We treat these as embedded submanifolds of Euclidean space $\mathbb{R}^{n \times n}$ or $\mathbb{R}^{n \times d}$. 
In the following, we first introduce the required preliminary background, and then we introduce the important Riemannian matrix manifolds. 

Here, we focus on the Stiefel, Grassmann, and symmetric
positive definite matrix manifolds. Besides these matrix
manifolds studied in this section, there exist several other
important matrix manifolds in optimization and machine
learning. Examples include fixed-rank matrix manifolds
\cite{vandereycken2013low,mishra2014fixed},
the orthogonal group and the special orthogonal group
$\mathrm{SO}(n)$ \cite{edelman1998geometry,boumal2023introduction},
oblique manifolds \cite{absil2006joint,absil2008optimization},
and more general matrix Lie groups
\cite{hall2013lie,boumal2023introduction}. In this monograph,
we focus on the Stiefel, Grassmann, and SPD manifolds because
they are among the most fundamental and frequently used
examples, and because they illustrate the main geometric
constructions in a concrete way.

\subsection{Preliminary Background}

\subsubsection{Group, Lie Group, General Linear Group, Matrix Lie Group, and Orthogonal Group}

\begin{definition}[Group \cite{hall2013lie}]
A \textbf{group} is a set \(\mathcal{G}\) together with a
binary operation:
\[
\cdot : \mathcal{G} \times \mathcal{G} \to \mathcal{G},
\qquad
(\b{X},\b{Y}) \mapsto \b{X}\cdot \b{Y},
\]
such that:
\begin{enumerate}
    \item \textbf{closure:}
    \(\b{X}\cdot\b{Y}\in\mathcal{G}\) for all
    \(\b{X},\b{Y}\in\mathcal{G}\),
    \item \textbf{associativity:}
    \[
    (\b{X}\cdot\b{Y})\cdot\b{Z}
    =
    \b{X}\cdot(\b{Y}\cdot\b{Z}),
    \qquad
    \forall \b{X},\b{Y},\b{Z}\in\mathcal{G},
    \]
    \item \textbf{identity element:} there exists
    \(\b{e}\in\mathcal{G}\) such that:
    \[
    \b{e}\cdot\b{X}=\b{X}\cdot\b{e}=\b{X},
    \qquad
    \forall \b{X}\in\mathcal{G},
    \]
    \item \textbf{inverse element:} for every
    \(\b{X}\in\mathcal{G}\), there exists
    \(\b{X}^{-1}\in\mathcal{G}\) such that
    \[
    \b{X}\cdot\b{X}^{-1}
    =
    \b{X}^{-1}\cdot\b{X}
    =
    \b{e}.
    \]
\end{enumerate}
\end{definition}

\begin{definition}[Lie group \cite{hall2013lie, boumal2023introduction}]
A \textbf{Lie group} is a smooth manifold \(\mathcal{G}\)
equipped with a group structure such that the group
multiplication map:
\begin{align}
m : \mathcal{G} \times \mathcal{G} \to \mathcal{G},
\qquad
m(\b{X},\b{Y}) := \b{X}\b{Y},
\end{align}
and the inversion map:
\begin{align}
\operatorname{inv} : \mathcal{G} \to \mathcal{G},
\qquad
\operatorname{inv}(\b{X}) := \b{X}^{-1},
\end{align}
are smooth.
\end{definition}

The theory of Lie groups originates in the work of \textit{Sophus Lie}---a Norwegian mathematician---on continuous transformation groups; see the classical treatise of Lie and Engel \cite{lieengel1888, lie1890, lieengel1893}.

\begin{definition}[General linear group \cite{hall2013lie}]
The \textbf{general linear group}, denoted by
\(\mathrm{GL}(n)\), is the set of all invertible
\(n\times n\) real matrices:
\begin{align}
\boxed{
\mathrm{GL}(n)
:=
\{
\b{X}\in\mathbb{R}^{n\times n}
\mid
\det(\b{X}) \neq 0
\}.
}
\end{align}
Its group operation is matrix multiplication.
\end{definition}

\begin{definition}[Matrix Lie group \cite{hall2013lie, lee2013smooth}]
A \textbf{matrix Lie group} is a Lie group whose elements are
matrices and whose group operation is matrix multiplication.
Equivalently, it is a subgroup of \(\mathrm{GL}(n)\) which is
also a smooth embedded manifold.
\end{definition}

\begin{remark}[Why Lie groups are relevant for matrix manifolds]
Many important matrix manifolds are closely related to matrix
Lie groups \cite{hall2013lie, boumal2023introduction}. For example, the \textbf{orthogonal group}
\(\mathrm{O}(n)\):
\begin{equation}\label{equation_On_orthogonal_group}
\boxed{
\mathrm{O}(n) := \{ \b{Q} \in \mathbb{R}^{n \times n} \mid \b{Q}^\top \b{Q} = \b{I}_n \},} 
\end{equation}
and the \textbf{special orthogonal group}
\(\mathrm{SO}(n)\):
\begin{equation}
\boxed{
\mathrm{SO}(n) := \{ \b{Q}\in\mathbb{R}^{n\times n}
\mid
\b{Q}^{\top}\b{Q}=\b{I}_n,\ \det(\b{Q})=1 \},} 
\end{equation}
are matrix Lie groups. These groups act
naturally on matrix manifolds and often help explain their
geometry, symmetries, and quotient constructions \cite{hall2013lie}.

However, not every matrix manifold is a Lie group. For
example, the Stiefel manifold \(\mathrm{St}(n,d)\) is a matrix
manifold but, in general, it is not a group under matrix
multiplication because its elements are not square when
\(n \neq d\). Likewise, the Grassmann manifold
\(\mathrm{Gr}(n,d)\) is not a Lie group; rather, it is more
naturally viewed as a quotient manifold \cite{edelman1998geometry, boumal2023introduction}. The SPD manifold
\(\mathbb{S}_{++}^{n}\) is a smooth matrix manifold, but it is
not a Lie group under the usual matrix multiplication because
the product of two symmetric positive definite matrices is not
necessarily symmetric.
\end{remark}


\subsubsection{Embedded Submanifold and Inclusion Map}

In the following, we define embedded submanifold.
The important Riemannian matrix manifolds, i.e., Stiefel, Grassmannian, and SPD matrix manifolds, can be considered as embedded submanifolds of Euclidean space $\mathbb{R}^{n \times d}$. 

\begin{definition}[Submersion, immersion, and embedding {\citep[Chapter 4]{lee2013smooth}}]
Let $\mathcal{M}$ and $\mathcal{N}$ be two smooth manifolds. 
\begin{itemize}
\item A smooth map $F: \mathcal{M} \rightarrow \mathcal{N}$ is a \textbf{smooth submersion} if its differential is surjective, i.e., $\text{rank}(F) = \dim(\mathcal{N})$, where $\dim(\cdot)$ denotes the local dimensionality of manifold. In submersion, we have $\dim(\mathcal{M}) \geq \dim(\mathcal{N})$.
\item A smooth map $F: \mathcal{M} \rightarrow \mathcal{N}$ is a \textbf{smooth immersion} if its differential is injective, i.e., $\text{rank}(F) = \dim(\mathcal{M})$. In immersion, we have $\dim(\mathcal{M}) \leq \dim(\mathcal{N})$.
\item A smooth map $F: \mathcal{M} \rightarrow \mathcal{N}$ is a \textbf{topological embedding} if it is a homeomorphism to its image $F(\mathcal{M}) \subseteq \mathcal{N}$ in the subspace topology. 
\item A smooth map $F: \mathcal{M} \rightarrow \mathcal{N}$ is a \textbf{smooth embedding} if it is both a smooth immersion and a topological embedding. 
\item Let $\mathcal{M} \subseteq \mathcal{A}$ be a subset of the smooth manifold $\mathcal{A}$. The \textbf{inclusion map} is the map $\iota: \mathcal{M} \hookrightarrow \mathcal{A}$ defined by:
\begin{align}
\boxed{
\iota(\b{p}) = p, \quad \forall \b{p} \in \mathcal{M}.
}
\end{align}
\end{itemize}
\end{definition}



\begin{definition}[Embedded submanifold {\citep[Chapter 5]{lee2013smooth}}]
Let $\mathcal{A}$ be a smooth manifold. An \textbf{embedded submanifold} of $\mathcal{A}$ is a subset $\mathcal{M} \subseteq \mathcal{A}$ which is itself a manifold endowed with a smooth structure where the inclusion map $\iota: \mathcal{M} \hookrightarrow \mathcal{A}$ is a smooth embedding.
The quantity $\dim(\mathcal{A}) - \dim(\mathcal{M})$ is called the codimension of $\mathcal{M}$ in $\mathcal{A}$.
\end{definition}

\subsubsection{Ambient Space}

\begin{definition}[Frobenius inner product]
For the ambient space $\mathbb{R}^{n \times d}$, the standard inner product is the \textbf{Frobenius inner product}. For any two matrices $\b{A}, \b{B} \in \mathbb{R}^{n \times d}$, it is defined as:
\begin{equation}\label{equation_Frobenius_inner_product}
\langle \b{A}, \b{B} \rangle_F := \text{tr}(\b{A}^\top \b{B}) = \sum_{i=1}^n \sum_{j=1}^d A_{ij} B_{ij}, 
\end{equation}
where $\text{tr}(.)$ denotes the trace of matrix, and $A_{ij}$ and $B_{ij}$ denote the $(i,j)$-th element of matrix $\b{A}$ and $\b{B}$, respectively. 
\end{definition}

The Frobenius inner product generalizes the dot product to matrices by summing the products of all corresponding entries. In the context of embedded manifolds, it provides the ``Euclidean" baseline for measuring lengths and angles before accounting for the manifold's curvature.

\begin{definition}[Ambient space]
The \textbf{ambient space} for a Riemannian manifold is the surrounding Euclidean vector space in which the manifold is embedded. Intuitively, the ambient space is the large, ``easy" space. For matrix manifolds, the ambient space is the Euclidean space $\mathbb{R}^{n \times d}$, i.e., the set of all $n \times d$ real matrices:
\begin{equation}
\mathcal{E} = \mathbb{R}^{n \times d}.
\end{equation}
It is a vector space where we already know how to compute dot products (the Frobenius inner product).
\end{definition}

In Riemannian optimization, we usually deal with matrix manifolds which consist of matrices.
In such manifolds, we have:
\begin{itemize}
\item The \textit{ambient space} $\mathcal{A}$: This is the large, easy space. For matrices, this is $\mathbb{R}^{n \times d}$. It is a vector space where we already know how to compute dot products (the Frobenius inner product).
\item The \textit{manifold} $\mathcal{M}$: This is the constrained surface sitting inside the ambient space (e.g., the set of matrices which are orthogonal satisfying $\b{X}^\top \b{X} = \b{I}$).
\item The \textit{inclusion map} $\iota$: This is just a formal way of saying ``take a point on the manifold and treat it as a point in the ambient space": $\iota(\b{X}) = \b{X}$.
\end{itemize}

\subsubsection{Pushforward and Pullback}

\begin{definition}[Pushforward and pullback]
Let:
\begin{align*}
\phi: \mathcal{M} \to \mathcal{N},
\end{align*}
be a smooth map from manifold $\mathcal{M}$ to manifold $\mathcal{N}$.
\begin{itemize}
\item The \textbf{pushforward} or \textbf{differential}, denoted by $\phi_*$ or $d\phi_{\b{p}}$, is a map:
\begin{align}
\boxed{
\phi_*: T_{\b{p}}\mathcal{M} \to T_{\phi(\b{p})}\mathcal{N}, 
}
\end{align}
which maps tangent vectors from the source manifold to the target manifold. 
In other words, \underline{pushforward moves vectors forward}.

For an embedding, pushforward identifies a tangent vector on the manifold with its representation in the ambient space.

\item The \textbf{pullback}, denoted by $\phi^*$, is a map:
\begin{align}
\boxed{
\phi^*: T_{\phi(\b{p})}\mathcal{N} \to T_{\b{p}}\mathcal{M}, 
}
\end{align}
which moves tensors (like metrics) from the target manifold $\mathcal{N}$ back to the source manifold $\mathcal{M}$. 
\underline{Pullback moves functions and tensors backward}.

For a metric tensor $g$ on $\mathcal{N}$, the pullback metric tensor $\phi^*g$ on $\mathcal{M}$ is:
\begin{equation}\label{equation_pullback_metric_1}
\boxed{
(\phi^* g)_{\b{p}}(\b{\eta}, \b{\zeta}) := g_{\phi(\b{p})}(\phi_* \b{\eta}, \phi_* \b{\zeta}),
}
\end{equation}
where $\b{p} \in \mathcal{M}$ is a point on $\mathcal{M}$, and $\b{\eta}, \b{\zeta} \in T_{\b{p}}\mathcal{M}$ are tangent vectors on $\mathcal{M}$, and $\phi_* \b{\eta}, \phi_* \b{\zeta} \in T_{\phi(\b{p})}\mathcal{N}$ are tangent vectors on $\mathcal{N}$.
\end{itemize}
\end{definition}

\subsubsection{Pullback Metric}\label{section_pullback_metric}

Imagine you have a flat sheet of rubber (the ambient space) and you have a ruler that only works on flat surfaces. Now, you wrap that rubber sheet over a ball (the manifold).

If you want to measure the distance between two points on the ball, you don't have a ``curved ruler". Instead, you look at how the rubber sheet was stretched over the ball and use the flat ruler you already had for the rubber. You are \textit{pulling back} the measurement system from the flat sheet onto the surface of the ball.

A Riemannian metric is a way to calculate the inner product of two tangent vectors.
On a manifold, the tangent vectors $\b{\eta}$ and $\b{\zeta}$ live in the tangent space $T_{\b{X}}\mathcal{M}$. However, because the matrix manifold $\mathcal{M}$ is embedded in the ambient space (i.e., $\mathbb{R}^{n \times d}$), these tangent vectors are also just matrices in the ambient space.

We do not want to invent a new complicated formula for the inner product on the manifold if we do not have to. The pullback metric says since these tangent vectors are already matrices, let us just use the matrix inner product we already have.

\begin{definition}[Pullback metric]
Let $\iota: \mathcal{M} \hookrightarrow \mathcal{A}$ be the inclusion of the manifold $\mathcal{M}$ into the ambient space $\mathcal{A}$. Let $\langle \cdot, \cdot \rangle_{\mathcal{A}}$ be the standard inner product in the ambient space (the Frobenius inner product).
According to Eq. (\ref{equation_pullback_metric_1}), the \textbf{pullback metric} $g$ on the manifold $\mathcal{M}$ is:
\begin{align}\label{equation_pullback_metric_2}
\boxed{
g_{\b{p}}(\b{\eta}, \b{\zeta}) = \langle d\iota_p(\b{\eta}), d\iota_p(\b{\zeta}) \rangle_{\mathcal{A}},
}
\end{align}
where $\eta, \zeta$ are tangent vectors on manifold $\mathcal{M}$, and $d\iota_p$ is the ``differential" or ``pushforward" which takes a vector tangent to the manifold and maps it to the corresponding vector in the ambient space. For embedded manifolds, $d\iota_p(\b{\eta})$ is simply $\b{\eta}$ viewed as an $n \times d$ matrix.

The pullback metric pulls the inner product from the (e.g., flat) ambient space down onto the tangent space of the (e.g., curved) manifold.
If the ambient space is the Euclidean space $\mathbb{R}^{n \times d}$, the $\langle \cdot, \cdot \rangle_{\mathcal{A}}$ is the Frobenius inner product $\text{tr}(\b{\eta}^\top \b{\zeta})$.
\end{definition}

\subsubsection{Smooth Local Extension on the Ambient Space}\label{section_smooth_local_extension_ambient_space}

\begin{definition}[Smooth local extension on the ambient space]
Let $\mathcal{M}$ be a smooth manifold embedded in an ambient Euclidean space
$\mathbb{R}^N$, let $f : \mathcal{M} \to \mathbb{R}$ be a smooth function, and
let $\b{p} \in \mathcal{M}$. A \textbf{smooth local extension} of $f$ around $\b{p}$ is a smooth
function:
\begin{align}
\bar{f} : U \subset \mathbb{R}^N \to \mathbb{R},
\end{align}
defined on an open neighborhood $U$ of $\b{p}$, such that:
\begin{align}\label{equation_smooth_local_extension_general}
\boxed{
\bar{f}(\b{q}) = f(\b{q}), \qquad \forall \b{q} \in U \cap \mathcal{M}.
}
\end{align}
In other words, $\bar{f}$ agrees with $f$ on the points of the manifold near
$\b{p}$, but it is defined on an open subset of the ambient Euclidean space so that
its Euclidean derivatives can be computed.
\end{definition}

We will use the smooth local extension on the ambient space for deriving the Riemannian gradient and Riemannian Hessian in matrix manifolds. We need the smooth local extension

\subsubsection{Euclidean Gradient}

\begin{definition}[Euclidean gradient]\label{definition_euclidean_gradient}
Let $\mathcal{M}$ be a Riemannian submanifold of $\mathbb{R}^{n \times d}$.
Let $f : \mathcal{M} \to \mathbb{R}$ be a smooth function and
let $\bar{f}$ be smooth local extension of $f$.
The \textbf{Euclidean gradient} of the smooth cost function $\bar{f}$ is denoted by $\nabla \bar{f}(\b{X})$.
The $(i,j)$-th component of $\nabla \bar{f}(\b{X})$ is the partial derivative of the function $f$ with respect to the $(i,j)$-th element of the matrix $\b{X}$:
\begin{align}
\boxed{
(\nabla \bar{f}(\b{X}))_{ij} := \frac{\partial \bar{f}}{\partial \b{X}_{ij}},
}
\end{align}
in which $\b{X}_{ij}$ denotes the $(i,j)$-th element of the matrix $\b{X}$.
\end{definition}

\subsubsection{Ambient Directional Derivative}\label{section_ambient_directional_derivative_matrix_manifolds}

Recall Section \ref{section_directional_derivative} for directional derivative along a vector field. Here, we introduce the directional derivative in the ambient space for matrix manifolds. 

\begin{definition}[Ambient directional derivative]\label{definition_ambient_directional_derivative}
Let $\mathcal{M}$ be a matrix manifold embedded in an ambient
Euclidean space $\mathbb{R}^{n\times d}$. Let
$f:\mathcal{M}\to\mathbb{R}$ be a smooth function, and let:
\[
\bar f:U\subset\mathbb{R}^{n\times d}\to\mathbb{R},
\]
be a smooth local extension of $f$ on an open neighborhood $U$ of
$\b{X}\in\mathcal{M}$, as in Section~\ref{section_smooth_local_extension_ambient_space}.
For $\b{X}\in\mathcal{M}$ and $\b{\Delta}\in T_{\b{X}}\mathcal{M}$, the
\textbf{ambient directional derivative} of $f$ at $\b{X}$ along $\b{\Delta}$
is defined by:
\begin{equation}\label{equation_ambient_directional_derivative_matrix_manifold}
\boxed{
\begin{aligned}
Df(\b{X})[\b{\Delta}]
&:=
D\bar f(\b{X})[\b{\Delta}]
\\
&=
\lim_{h\to 0}
\frac{\bar f(\b{X}+h\b{\Delta})-\bar f(\b{X})}{h}.
\end{aligned}
}
\end{equation}
If $\bar f$ is differentiable, then:
\begin{equation}\label{equation_ambient_directional_derivative_matrix_manifold_inner_product}
\boxed{
D\bar f(\b{X})[\b{\Delta}]
=
\langle \nabla \bar f(\b{X}), \b{\Delta}\rangle_F
=
\text{tr}\big((\nabla \bar f(\b{X}))^\top \b{\Delta}\big).
}
\end{equation}
\end{definition}

\begin{remark}[Compatibility with the general directional derivative]
Section~\ref{section_directional_derivative} defined the directional
derivative of a smooth function on a general manifold. The ambient
directional derivative, introduced here, is the corresponding
Euclidean-space realization for embedded matrix manifolds. Because the curve $\b{X}+h\b{\Delta}$ generally leaves the manifold for $h\neq 0$, one uses a smooth local extension $\bar f$ of $f$ to the ambient space and computes the derivative there.
\end{remark}

\begin{remark}[Expression of ambient directional derivative in matrix manifolds]
Let $f:\mathcal{M}\to\mathbb{R}$ be a smooth function on a matrix
manifold $\mathcal{M}$ embedded in $\mathbb{R}^{n\times d}$, and let:
\[
\bar f:U\subset\mathbb{R}^{n\times d}\to\mathbb{R},
\]
be a smooth local extension of $f$ around $\b{X}\in\mathcal{M}$.
If $\b{\Delta}$ is a tangent vector field on $\mathcal{M}$, then the
directional derivative of $f$ along $\b{\Delta}$ at $\b{X}$ is\footnote{This equation is the counterpart of Eq. (\ref{equation_directional_derivative_function_along_vector_field}) for matrix manifolds, but with the notations used here.}:
\begin{equation}\label{equation_directional_derivative_matrix_2}
\boxed{
\b{\Delta}(\bar{f})(\b{X})
=
D \bar{f}(\b{X})[\b{\Delta}(\b{X})],
}
\end{equation}
Here, $\b{X}$ denotes a point on the manifold, while $\b{\Delta}(\b{X})\in T_{\b{X}}
\mathcal{M}$ is the tangent vector assigned by the vector field
$\b{\Delta}$ at the point $\b{X}$.
\end{remark}

\begin{remark}[Intrinsic and ambient viewpoints for directional derivative on matrix manifolds]\label{remark_intrinsic_ambient_vewpoints_directional_derivative_matrix_manifold}
Let $\mathcal{M}$ be a matrix manifold, let
$f:\mathcal{M}\to\mathbb{R}$ be a smooth function, let
$\b{X} \in \mathcal{M}$, and let
$\b{\Delta} \in T_X\mathcal{M}$.

The directional derivative of $f$ at $\b{X}$ along $\b{\Delta}$ can be defined \textbf{intrinsically} exactly as on any smooth manifold. Namely, if:
\[
\Gamma:(-\epsilon,\epsilon)\to\mathcal{M},
\]
is a smooth curve such that
\[
\Gamma(0)=\b{X},
\qquad
\dot{\Gamma}(0)=\b{\Delta},
\]
then:
\begin{equation}
\boxed{
Df(\b{X})[\b{\Delta}]
=
\frac{d}{dt}f(\Gamma(t))\Big|_{t=0}.
}
\end{equation}
In this intrinsic viewpoint, the curve $\Gamma(t)$ lies on the manifold
$\mathcal{M}$, so no ambient extension is needed.

Because $\mathcal{M}$ is embedded in an \textbf{ambient} Euclidean space, one may
also compute the same directional derivative \textbf{extrinsically}. Let:
\[
\bar f:U\subset\mathbb{R}^{n\times d}\to\mathbb{R},
\]
be a smooth local extension of $f$ on an open neighborhood $U$ of $\b{X}$.
Then, the same directional derivative can be written as:
\begin{equation}
\boxed{
Df(\b{X})[\b{\Delta}]
=
D\bar f(\b{X})[\b{\Delta}]
=
\frac{d}{dt}\bar f(\b{X}+t\b{\Delta})\Big|_{t=0}.
}
\end{equation}
Here, the curve:
\[
t \mapsto \b{X}+t\b{\Delta},
\]
is generally a curve in the ambient Euclidean space rather than on the
manifold itself. Therefore, the ambient extension $\bar f$ is needed in
this extrinsic viewpoint.

Hence, the intrinsic and ambient viewpoints define the same directional
derivative:
\begin{align}
\boxed{
\frac{d}{dt}f(\Gamma(t))\Big|_{t=0}
=
D\bar f(\b{X})[\b{\Delta}].
}
\end{align}
The intrinsic viewpoint is the general manifold definition, while the ambient viewpoint is often more convenient for explicit matrix
calculations.
\end{remark}

According to Remark \ref{remark_intrinsic_ambient_vewpoints_directional_derivative_matrix_manifold}, for a matrix manifold, the directional derivative can be defined
intrinsically exactly as on any smooth manifold, using curves that lie
on the manifold. However, because many matrix manifolds are embedded in
an ambient Euclidean space, we often compute the same derivative
extrinsically through a smooth local extension $\bar f$ to the ambient space. The intrinsic and extrinsic viewpoints represent the same
directional derivative, but the extrinsic one is often more convenient
for explicit matrix calculations.




\begin{definition}[Ambient directional derivative of a vector field on a matrix manifold]
Let \(\mathcal{M}\) be a matrix manifold embedded in an
ambient Euclidean space \(\mathbb{R}^{n\times d}\). Let
\(\b{V}\) be a smooth tangent vector field on \(\mathcal{M}\),
so that:
\[
\b{V}(\b{X}) \in T_{\b{X}}\mathcal{M},
\qquad
\forall \b{X}\in\mathcal{M}.
\]
Let \(\b{U}\) be another smooth tangent vector field on
\(\mathcal{M}\). Since \(\mathcal{M}\) is embedded in
\(\mathbb{R}^{n\times d}\), we may view \(\b{V}\) locally as a
matrix-valued map. Equivalently, we may use a smooth local
extension of \(\b{V}\) to an open neighborhood of \(\b{X}\) in
the ambient Euclidean space.

The ambient directional derivative of \(\b{V}\) at \(\b{X}\)
along \(\b{U}(\b{X})\) is defined by:
\begin{align}
\boxed{
D\b{V}(\b{X})[\b{U}(\b{X})]
:=
\frac{d}{dt}
\b{V}\big(\b{X}+t\b{U}(\b{X})\big)
\Big|_{t=0},
}
\end{align}
whenever this expression is computed using a smooth local
ambient extension of \(\b{V}\). Equivalently:
\begin{align}
\boxed{
D\b{V}(\b{X})[\b{U}(\b{X})]
=
\lim_{t\to 0}
\frac{
\b{V}\big(\b{X}+t\b{U}(\b{X})\big)
-
\b{V}(\b{X})
}{t}.
}
\end{align}
In general, we have:
\[
D\b{V}(\b{X})[\b{U}(\b{X})]
\in \mathbb{R}^{n\times d},
\]
and it does not necessarily belong to the tangent space
\(T_{\b{X}}\mathcal{M}\).
\end{definition}

\begin{lemma}[Projection formula for the Levi-Civita connection under the induced Euclidean metric]\label{lemma_levi_civita_ambient_directional_derivative_relation}
For a Riemannian submanifold \(\mathcal{M}\subset \mathbb{R}^{n\times d}\)
equipped with induced Euclidean metric (the metric induced by the ambient Frobenius inner product),
the Levi-Civita connection is the orthogonal projection of the ambient
Euclidean directional derivative onto the tangent space:
\begin{align}\label{equation_levi_civita_connection_ambient_projection}
\boxed{
(\nabla_{\b{U}} \b{V})(\b{X})
=
\Pi_{\b{X}}
\big(
D\b{V}(\b{X})[\b{U}(\b{X})]
\big),
}
\end{align}
where $\Pi_{\b{X}}$ is the orthogonal projection onto $T_{\b{X}}\mathcal{M}$.
\end{lemma}
\begin{proof}
Let \(\mathcal{M}\) be an embedded submanifold of the
ambient Euclidean space \(\mathbb{R}^{n\times d}\), and suppose
that \(\mathcal{M}\) is equipped with the Riemannian metric
induced from the ambient Frobenius inner product:
\[
g_{\b{X}}(\b{\Delta}_1,\b{\Delta}_2)
=
\langle \b{\Delta}_1,\b{\Delta}_2\rangle_F
=
\operatorname{tr}(\b{\Delta}_1^{\top}\b{\Delta}_2),
\]
for $\b{\Delta}_1,\b{\Delta}_2 \in T_{\b{X}}\mathcal{M}$.
Let \(\b{U}\) and \(\b{V}\) be smooth tangent vector fields on
\(\mathcal{M}\). Since \(\mathcal{M}\) is embedded in
\(\mathbb{R}^{n\times d}\), we can view \(\b{V}\) locally as a
matrix-valued map. Therefore, its ambient Euclidean
directional derivative at \(\b{X}\) along \(\b{U}(\b{X})\) is:
\[
D\b{V}(\b{X})[\b{U}(\b{X})] \in \mathbb{R}^{n\times d}.
\]
However, this ambient derivative does not necessarily lie
in the tangent space \(T_{\b{X}}\mathcal{M}\). Since
\(\mathbb{R}^{n\times d}\) is equipped with the Frobenius
inner product, every ambient matrix can be orthogonally
decomposed into a tangent part and a normal part:
\begin{align*}
D\b{V}(\b{X})[\b{U}(\b{X})]
=\,
&\Pi_{\b{X}}\big(D\b{V}(\b{X})[\b{U}(\b{X})]\big)
\\
&+
\big(
I-\Pi_{\b{X}}
\big)
\big(D\b{V}(\b{X})[\b{U}(\b{X})]\big),
\end{align*}
where:
\begin{align*}
&\Pi_{\b{X}}\big(D\b{V}(\b{X})[\b{U}(\b{X})]\big)
\in T_{\b{X}}\mathcal{M}, \\
& \big(
I-\Pi_{\b{X}}
\big)
\big(D\b{V}(\b{X})[\b{U}(\b{X})]\big)
\in N_{\b{X}}\mathcal{M}.
\end{align*}
Here, \(N_{\b{X}}\mathcal{M}\) denotes the normal space of
\(\mathcal{M}\) at \(\b{X}\).

We define:
\[
(\widetilde{\nabla}_{\b{U}}\b{V})(\b{X})
:=
\Pi_{\b{X}}
\big(
D\b{V}(\b{X})[\b{U}(\b{X})]
\big).
\]
We show that \(\widetilde{\nabla}\) is the Levi-Civita
connection on \(\mathcal{M}\). By the uniqueness of the
Levi-Civita connection, it suffices to show that
\(\widetilde{\nabla}\) is torsion-free and compatible with
the induced metric (see Definition \ref{definition_Levi_Civita_connection_coordinate_free}).

First, we prove torsion-freeness. Since the ambient
Euclidean connection is flat and torsion-free, for smooth
ambient extensions of \(\b{U}\) and \(\b{V}\), we have:
\[
D\b{V}(\b{X})[\b{U}(\b{X})]
-
D\b{U}(\b{X})[\b{V}(\b{X})]
=
[\b{U},\b{V}](\b{X}).
\]
Because \(\b{U}\) and \(\b{V}\) are tangent vector fields on
\(\mathcal{M}\), their Lie bracket \([\b{U},\b{V}]\) is also
tangent to \(\mathcal{M}\). Hence:
\[
\Pi_{\b{X}}\big([\b{U},\b{V}](\b{X})\big)
=
[\b{U},\b{V}](\b{X}).
\]
Therefore:
\[
\begin{aligned}
&(\widetilde{\nabla}_{\b{U}}\b{V})(\b{X})
-
(\widetilde{\nabla}_{\b{V}}\b{U})(\b{X}) \\
&=
\Pi_{\b{X}}
\big(
D\b{V}(\b{X})[\b{U}(\b{X})]
\big)
-
\Pi_{\b{X}}
\big(
D\b{U}(\b{X})[\b{V}(\b{X})]
\big)
\\
&=
\Pi_{\b{X}}
\big(
D\b{V}(\b{X})[\b{U}(\b{X})]
-
D\b{U}(\b{X})[\b{V}(\b{X})]
\big)
\\
&=
\Pi_{\b{X}}\big([\b{U},\b{V}](\b{X})\big)
\\
&=
[\b{U},\b{V}](\b{X}).
\end{aligned}
\]
Thus, \(\widetilde{\nabla}\) is torsion-free.

Second, we prove metric compatibility. Let \(\b{W}\) be another
smooth tangent vector field on \(\mathcal{M}\). Since the
metric is induced from the ambient Frobenius inner product,
we have:
\[
g_{\b{X}}(\b{V}(\b{X}),\b{W}(\b{X}))
=
\langle \b{V}(\b{X}),\b{W}(\b{X})\rangle_F.
\]
Taking the directional derivative along \(\b{U}\), we obtain:
\[
\b{U}\big(g(\b{V},\b{W})\big)(\b{X})
=
D\Big(
\langle \b{V},\b{W}\rangle_F
\Big)(\b{X})[\b{U}(\b{X})].
\]
Using the product rule for the Frobenius inner product gives:
\begin{align*}
\b{U}\big(g(\b{V},\b{W})\big)(\b{X})
=\,
&\big\langle
D\b{V}(\b{X})[\b{U}(\b{X})],
\b{W}(\b{X})
\big\rangle_F
\\
&+
\big\langle
\b{V}(\b{X}),
D\b{W}(\b{X})[\b{U}(\b{X})]
\big\rangle_F.
\end{align*}
Now, we decompose the ambient derivatives into tangent and
normal components. Since \(\b{W}(\b{X})\in T_{\b{X}}\mathcal{M}\)
is orthogonal to every normal vector, we have:
\begin{align*}
\big\langle
D\b{V}(\b{X})&[\b{U}(\b{X})],
\b{W}(\b{X})
\big\rangle_F
\\
&=
\big\langle
\Pi_{\b{X}}
\big(
D\b{V}(\b{X})[\b{U}(\b{X})]
\big),
\b{W}(\b{X})
\big\rangle_F.
\end{align*}
Similarly, since \(\b{V}(\b{X})\in T_{\b{X}}\mathcal{M}\), we have:
\begin{align*}
\big\langle
\b{V}(\b{X}),
D\b{W}&(\b{X})[\b{U}(\b{X})]
\big\rangle_F
\\
&=
\big\langle
\b{V}(\b{X}),
\Pi_{\b{X}}
\big(
D\b{W}(\b{X})[\b{U}(\b{X})]
\big)
\big\rangle_F.
\end{align*}
Therefore:
\[
\begin{aligned}
&\b{U}\big(g(\b{V},\b{W})\big)(\b{X})
\\
&=
\big\langle
(\widetilde{\nabla}_{\b{U}}\b{V})(\b{X}),
\b{W}(\b{X})
\big\rangle_F
+
\big\langle
\b{V}(\b{X}),
(\widetilde{\nabla}_{\b{U}}\b{W})(\b{X})
\big\rangle_F
\\
&=
g_{\b{X}}\big(
(\widetilde{\nabla}_{\b{U}}\b{V})(\b{X}),\b{W}(\b{X})
\big)
\\
&\quad\quad\quad\quad\quad\quad\quad\quad+
g_{\b{X}}\big(
\b{V}(\b{X}),(\widetilde{\nabla}_{\b{U}}\b{W})(\b{X})
\big).
\end{aligned}
\]
Hence, \(\widetilde{\nabla}\) is compatible with the induced
Riemannian metric.

We have shown that \(\widetilde{\nabla}\) is torsion-free and
metric-compatible. By the uniqueness of the Levi-Civita
connection, \(\widetilde{\nabla}\) is the Levi-Civita connection
on the embedded Riemannian submanifold \(\mathcal{M}\).
Therefore:
\[
(\nabla_{\b{U}}\b{V})(\b{X})
=
\Pi_{\b{X}}
\big(
D\b{V}(\b{X})[\b{U}(\b{X})]
\big),
\]
which proves Eq. (\ref{equation_levi_civita_connection_ambient_projection}).
\end{proof}

\begin{remark}[Scope of the projection formula for the Levi-Civita connection]
Lemma \ref{lemma_levi_civita_ambient_directional_derivative_relation} applies to embedded submanifolds equipped with the
Riemannian metric induced by the ambient Euclidean inner product.
In that case, the Levi-Civita connection is obtained by taking the
ambient Euclidean directional derivative and projecting it onto the
tangent space.

For other Riemannian metrics on the same embedded manifold, the
Levi-Civita connection may contain additional correction terms. For
example, on the Stiefel manifold endowed with the canonical metric,
the connection contains an extra term depending on \(\b{X}\), as shown
later in Lemma \ref{lemma_levi_civita_stiefel_canonical}. Thus, Lemma \ref{lemma_levi_civita_ambient_directional_derivative_relation} and Lemma \ref{lemma_levi_civita_stiefel_canonical} are not
contradictory; they correspond to different choices of Riemannian
metric.
\end{remark}

\subsubsection{Characterizing Identity of the Gradient}

While the Riemannian gradient is defined in Eq. (\ref{equation_Riemannian_gradient_}) via the inverse metric components, $\operatorname{grad} f = (g^{ij} \partial_i \bar{f}) \partial_j$, in the context of embedded matrix manifolds, it is often more practical to use the identity $g(\operatorname{grad} f, \b{\Delta}) = \operatorname{tr}(\nabla \bar{f}(\b{X})^\top \b{\Delta})$ for a tangent vector $\b{\Delta}$. We prove this equivalence using the coordinate definitions established in this work. This equation is referred to as ``\textit{characterizing identity of the gradient}" or ``\textit{defining property of the Riemannian gradient}" or ``\textit{gradient compatibility condition}" in the literature. 

\begin{proposition}[Characterizing identity of the gradient]\label{proposition_characterizing_identity_of_gradient}
Let $(\mathcal{M}, g)$ be a Riemannian submanifold of $\mathbb{R}^{n \times d}$. For any smooth function $f: \mathcal{M} \to \mathbb{R}$ and tangent vector $\b{\Delta} \in T_{\b{X}} \mathcal{M}$, the coordinate-based gradient $\operatorname{grad} f = (g^{ij} \partial_i f) \partial_j$ satisfies:
\begin{equation}\label{equation_characterizing_identity_of_gradient}
\boxed{
g(\operatorname{grad} f, \b{\Delta}) = \operatorname{tr}(\nabla \bar{f}(\b{X})^\top \b{\Delta}).
}
\end{equation}
This equation is equivalent to:
\begin{equation}\label{equation_directional_derivative_inner_product}
\boxed{
\begin{aligned}
Df(\b{X})[\b{\Delta}] &= D\bar{f}(\b{X})[\b{\Delta}] \\
&= \langle \nabla \bar{f}(\b{X}), \b{\Delta} \rangle_F \overset{(\ref{equation_Frobenius_inner_product})}{=} \operatorname{tr}(\nabla \bar{f}(\b{X})^\top \b{\Delta}),
\end{aligned}
}
\end{equation}
where $Df(\b{X})[\b{\Delta}]$ or $D\bar{f}(\b{X})[\b{\Delta}]$ denotes classical directional derivative in the ambient Euclidean space $\mathbb{R}^{n \times d}$.
\end{proposition}
\begin{proof}
Let $\{ \partial_1, \dots, \partial_m \}$ be a local basis for the tangent space $T_{\b{X}} \mathcal{M}$. Any tangent vector $\b{\Delta}$ can be expressed as $\b{\Delta} = \Delta^k \partial_k$.
We have:
\begin{align*}
g(\operatorname{grad} f, \b{\Delta}) &\overset{(\ref{equation_Riemannian_gradient_})}{=} g\left( (g^{ij} \partial_i \bar{f}) \partial_j, \Delta^k \partial_k \right) \\
&\overset{(a)}{=} g^{ij} (\partial_i \bar{f}) \Delta^k g(\partial_j, \partial_k) \\
&\overset{(\ref{equation_g_components})}{=} g^{ij} (\partial_i \bar{f}) \b{\xi}^k g_{jk} = (\partial_i \bar{f}) \Delta^k (g^{ij} g_{jk}) \\
&\overset{(\ref{equation_metric_inverse})}{=} (\partial_i \bar{f}) \Delta^k \delta^i_k \overset{(\ref{equation_index_substitution_delta})}{=} (\partial_i \bar{f}) \Delta^i,
\end{align*}
where $(a)$ is because of linearity of the Riemannian metric.

The term $(\partial_i \bar{f}) \Delta^i$ is the directional derivative $D\bar{f}(\b{X})[\b{\Delta}]$. In the ambient Euclidean space $\mathbb{R}^{n \times d}$, this derivative is calculated as the Frobenius inner product with the Euclidean gradient $\nabla \bar{f}(\b{X})$:
\begin{equation*}
D\bar{f}(\b{X})[\b{\Delta}] = \langle \nabla \bar{f}(\b{X}), \b{\Delta} \rangle_F \overset{(\ref{equation_Frobenius_inner_product})}{=} \operatorname{tr}(\nabla \bar{f}(\b{X})^\top \b{\Delta}).
\end{equation*}
Thus, $g(\operatorname{grad} f, \b{\Delta}) = \operatorname{tr}(\nabla \bar{f}(\b{X})^\top \b{\Delta})$ holds for any Riemannian metric $g$.
\end{proof}

\subsection{Stiefel Manifold}
\label{section_stiefel_manifold}

The Stiefel\footnote{Stiefel is pronounced as SHTEE-fel
or \textipa{/Stifel/} in simplified international phonetic alphabet.} manifold is the set of matrices with orthonormal columns. It is a central object in dimensionality reduction \cite{ghojogh2023elements} (e.g., principal component analysis \cite{ghojogh2023principal}) and orthogonal neural networks.
The Stiefel manifold was originally proposed and defined by the Swiss mathematician \textit{Eduard Stiefel} in his 1935 doctoral thesis \cite{stiefel1935richtungsfelder}.
Many of the geometric, computational, and optimization-related
characteristics of the Stiefel manifold are analyzed in
\cite{edelman1998geometry,absil2008optimization}.

\subsubsection{Definition of Stiefel Manifold}\label{section_definition_Stiefel_manifold}

\begin{definition}[Stiefel manifold]\label{definition_stiefel_manifold}
The \textbf{Stiefel manifold}, denoted by $St(n, d)$, is defined as the set of all $n \times d$ real matrices with orthonormal columns, where $n \geq d$:
\begin{equation}\label{equation_Sitefel_manifold_definition}
\boxed{
St(n, d) := \{ \b{X} \in \mathbb{R}^{n \times d} : \b{X}^\top \b{X} = \b{I}_d \},
}
\end{equation}
where $\b{I}_d$ is the $d \times d$ identity matrix.
In other words, the Stiefel manifold is the set of orthogonal matrices where the columns of matrix are orthonormal. 
\end{definition}

\begin{lemma}[Points of the Stiefel manifold are orthonormal matrix frames]
Let \(\b{X} \in \mathbb{R}^{n\times d}\). If
\(\b{X} \in \mathrm{St}(n,d)\), then \(\b{X}\) is a point of
the Stiefel manifold. Conversely, every point
\(\b{p} \in \mathrm{St}(n,d)\) can be written as $\b{p} = \b{X}$
for some \(\b{X} \in \mathrm{St}(n,d)\) satisfying $ \b{X}^\top \b{X} = \b{I}_d$. Thus:
\begin{align}
\boxed{
\b{p} = \b{X} \quad \text{such that} \quad \b{X}^\top \b{X} = \b{I}_d.
}
\end{align}

In other words, points of the Stiefel manifold are exactly
the \(n\times d\) matrices with orthonormal columns.
\end{lemma}
\begin{proof}
By definition, the Stiefel manifold is:
\[
\mathrm{St}(n,d)
:=
\left\{
\b{X} \in \mathbb{R}^{n\times d}
\;\middle|\;
\b{X}^{\top}\b{X} = \b{I}_{d}
\right\}.
\]
Hence, every element of \(\mathrm{St}(n,d)\) is an
\(n\times d\) matrix \(\b{X}\) satisfying
\(\b{X}^{\top}\b{X} = \b{I}_{d}\).

Therefore, if \(\b{X} \in \mathrm{St}(n,d)\), then
\(\b{X}\) is an element of the set \(\mathrm{St}(n,d)\), and
thus it is a point of the Stiefel manifold:
\[
\b{X} \in \mathrm{St}(n,d).
\]

Conversely, let \(\b{p} \in \mathrm{St}(n,d)\). Since
\(\mathrm{St}(n,d)\) is a set of matrices, every element of
it is, by definition, some matrix
\(\b{X} \in \mathbb{R}^{n\times d}\) satisfying
\(\b{X}^{\top}\b{X} = \b{I}_{d}\). Therefore, there exists
\(\b{X} \in \mathrm{St}(n,d)\) such that:
\[
\b{p} = \b{X}.
\]

Hence, points of the Stiefel manifold are exactly matrices
with orthonormal columns.
\end{proof}

\subsubsection{Dimension of Stiefel Manifold}

\begin{remark}[Dimension and embedding of Stiefel manifold]\label{remark_dimension_Stiefel_manifold}
The Stiefel manifold is an embedded submanifold of $\mathbb{R}^{n \times d}$. Its dimensionality is determined by the $d \times d$ symmetric constraint $\b{X}^\top \b{X} = \b{I}_d$. Since a symmetric $d \times d$ matrix has $d(d+1)/2$ independent constraints, the dimension is:
\begin{equation}\label{equation_dimension_Stiefel_manifold}
\boxed{
\dim(St(n, d)) = nd - \frac{d(d+1)}{2}.
}
\end{equation}
\end{remark}

\subsubsection{Tangent and Normal Spaces of Stiefel Manifold}

The tangent space characterizes the first-order geometry of the manifold and defines the set of all feasible directions for optimization.

\begin{proposition}[Tangent space of Stiefel manifold]
The \textbf{tangent space} at a point $\b{X} \in St(n, d)$ is given by:
\begin{equation}\label{equation_stiefel_tangent_space}
\boxed{
T_{\b{X}}St(n, d) = \{ \b{\Delta} \in \mathbb{R}^{n \times d} \mid \b{X}^\top \b{\Delta} + \b{\Delta}^\top \b{X} = \b{0} \},
}
\end{equation}
where $\b{0}$ denotes the zero matrix. 
The $\b{\Delta} \in T_{\b{X}}St(n, d)$ is a tangent vector, but it is a matrix (because the manifold is a matrix manifold). So, we can call it \textbf{tangent vector} or \textbf{tangent matrix}. 
\end{proposition}
\begin{proof}
Consider a smooth curve $\b{X}(t): I \to St(n, d)$ such that:
\begin{align*}
\b{X}(0) = \b{X}, \quad \dot{\b{X}}(0) = \b{\Delta}.
\end{align*}
Since the curve lies on the manifold, it must satisfy the orthogonality constraint for all $t$:
\begin{equation*}
\b{X}(t)^\top \b{X}(t) = \b{I}_d.
\end{equation*}
Differentiating both sides with respect to $t$ yields:
\begin{equation*}
\dot{\b{X}}(t)^\top \b{X}(t) + \b{X}(t)^\top \dot{\b{X}}(t) = \b{0}.
\end{equation*}
Evaluating at $t=0$ with $\b{X}(0) = \b{X}$ and $\dot{\b{X}}(0) = \b{\Delta}$, we obtain:
\begin{equation*}
\b{\Delta}^\top \b{X} + \b{X}^\top \b{\Delta} = \b{0}.
\end{equation*}
\end{proof}

\begin{definition}[Symmetric and skew-symmetric matrices]\label{definition_symmetric_skew_symmetric_matrices}
A square matrix $A$ is said to be \textbf{symmetric} if it is equal to its transpose:
\begin{equation}\label{equation_symmetric_matrix}
\boxed{
\b{A}^\top = \b{A}.
}
\end{equation}
In terms of its entries, this means $\b{A}_{ij} = \b{A}_{ji}$ for all $i, j$, where $\b{A}_{ij}$ denotes the $(i,j)$-th element of matrix $\b{A}$.

A square matrix $\b{A}$ is said to be \textbf{skew-symmetric} if its transpose is equal to its negative:
\begin{equation}\label{equation_skew_symmetric_matrix}
\boxed{
\b{A}^\top = -\b{A}.
}
\end{equation}
In terms of its entries, this means $\b{A}_{ij} = -\b{A}_{ji}$ for all $i, j$. Notably, the diagonal elements must satisfy $\b{A}_{ii} = 0$.
\end{definition}

\begin{lemma}[$\b{X}^\top \b{\Delta}$ is skew-symmetric in Stiefel manifold]\label{lemma_skew_symmetric}
Suppose $\b{\Delta} \in T_{\b{X}}St(n, d)$ is a tangent vector (matrix) in the tangent space of Stiefel manifold. 
The condition $\b{X}^\top \b{\Delta} + \b{\Delta}^\top \b{X} = \b{0}$ in Eq. (\ref{equation_stiefel_tangent_space}) implies that the matrix $\b{X}^\top \b{\Delta}$ is skew-symmetric, satisfying:
\begin{align}\label{equation_skew_symmetric}
\boxed{
(\b{X}^\top \b{\Delta})^\top = - \b{X}^\top \b{\Delta}.
}
\end{align}
\end{lemma}
\begin{proof}
\begin{align*}
&\b{X}^\top \b{\Delta} + \b{\Delta}^\top \b{X} = \b{0} \implies \b{X}^\top \b{\Delta} = - \b{\Delta}^\top \b{X} \\
&\overset{(a)}{\implies} (\b{X}^\top \b{\Delta})^\top = - (\b{\Delta}^\top \b{X})^\top = - \b{X}^\top \b{\Delta},
\end{align*}
where $(a)$ is because of taking transpose from the sides of equation.
\end{proof}

\begin{definition}[Normal space of Stiefel manifold]
The \textbf{normal space} $T_{\b{X}}^\perp \text{St}(n, d)$ consists of all matrices $\b{N} \in \mathbb{R}^{n \times d}$ that are orthogonal to every $\b{\Delta} \in T_{\b{X}} \text{St}(n, d)$. From the theory of constrained optimization (Lagrange multipliers), the normal space to a constraint $h(\b{X}) = \b{0}$ is spanned by the gradients of the constraints. For $\b{X}^\top \b{X} - \b{I} = \b{0}$ in the Stiefel manifold (see Eq. (\ref{equation_Sitefel_manifold_definition})), the normal space is:
\begin{align}\label{equation_normal_space_Stiefel_manifold}
\boxed{
T_{\b{X}}^\perp \text{St}(n, d) = \{ \b{X S} \mid \b{S} \in \mathbb{R}^{d \times d}, \b{S} = \b{S}^\top \}.
}
\end{align}
\end{definition}
\begin{proof}
$\langle \b{XS}, \b{\Delta} \rangle_F = \text{tr}(\b{S}^\top \b{X}^\top \b{\Delta}) = \text{tr}(\b{S} \b{X}^\top \b{\Delta})$. Since $\b{S}$ is symmetric and $\b{X}^\top \b{\Delta}$ is skew-symmetric, the trace of their product is zero. As their inner product is zero, $\b{XS}$ and $\b{\Delta}$ are orthogonal. Thus, the space consisting $\b{XS}$ is normal to the tangent space consisting the tangent vector (matrix) $\b{\Delta}$. 
\end{proof}

\begin{lemma}[Symmetric and skew-symmetric decomposition]\label{lemma_sym_skew_expressions}
Let $\b{A} \in \mathbb{R}^{n \times n}$ be a square matrix. The matrix can be uniquely decomposed into its symmetric part $\operatorname{sym}(\b{A})$ and its skew-symmetric part $\operatorname{skew}(\b{A})$, where:
\begin{equation}\label{equation_sym_skew_expressions}
\boxed{
\begin{aligned}
&\operatorname{sym}(\b{A}) := \frac{1}{2}(\b{A} + \b{A}^\top), \\
&\operatorname{skew}(\b{A}) := \frac{1}{2}(\b{A} - \b{A}^\top).
\end{aligned}
}
\end{equation}
Furthermore, it holds that:
\begin{align}
\boxed{
\b{A} = \operatorname{sym}(\b{A}) + \operatorname{skew}(\b{A}).
}
\end{align}
\end{lemma}
\begin{proof}
Consider the decomposition $\b{A} = \b{S} + \b{K}$ where 
$\b{S} = \frac{1}{2}(\b{A} + \b{A}^\top)$ and 
$\b{K} = \frac{1}{2}(\b{A} - \b{A}^\top)$. 
Direct calculation shows:
\begin{itemize}
    \item $\b{S}^\top = \frac{1}{2}(\b{A}^\top + \b{A}) = \b{S}$, so $\b{S}$ is symmetric according to Eq. (\ref{equation_symmetric_matrix}).
    \item $\b{K}^\top = \frac{1}{2}(\b{A}^\top - \b{A}) = -\b{K}$, so $\b{K}$ is skew-symmetric according to Eq. (\ref{equation_skew_symmetric_matrix}).
    \item $\b{S} + \b{K} = \frac{1}{2}(2\b{A}) = \b{A}$.
\end{itemize}
To show uniqueness, let $\b{A} = \b{S}' + \b{K}'$. Then $\b{A}^\top = \b{S}' - \b{K}'$. 
Solving for $\b{S}'$ and $\b{K}'$ yields the original definitions.
\end{proof}

\begin{remark}[Symmetrization in tangent constraint of the Stiefel and manifold]
On the Stiefel manifold, according to Eq. (\ref{equation_stiefel_tangent_space}), the tangent constraint is:
\[
\text{Tangent on Stiefel: } \b{X}^\top \b{\Delta} + \b{\Delta}^\top \b{X} = \b{0},
\]
so the projection involves the symmetrization operator, defined in Eq. (\ref{equation_sym_skew_expressions}):
\begin{equation}
\boxed{
\begin{aligned}
&\text{Tangent on Stiefel: } \operatorname{sym}(\b{X}^\top \b{\Delta}) = \b{0}, \text{ or, } \\
&\text{Tangent on Stiefel: } \operatorname{sym}(\b{\Delta}^\top \b{X}) = \b{0}.
\end{aligned}
}
\end{equation}
\end{remark}

\begin{lemma}[Projection of matrix onto tangent space of Stiefel manifold] \label{lemma_stiefel_projection}
For any matrix $\b{Z} \in \mathbb{R}^{n \times d}$, the orthogonal projection of $\b{Z}$ onto the tangent space $T_{\b{X}} \text{St}(n, d)$ at a point $\b{X} \in \text{St}(n, d)$ is denoted by $\Pi_{\b{X}}^{\mathrm{St}}(\b{Z})$ and is obtained as:
\begin{equation} \label{equation_stiefel_proj_formula}
\boxed{
\Pi_{\b{X}}^{\mathrm{St}}(\b{Z}) = \b{Z} - \b{X} \operatorname{sym}(\b{X}^\top \b{Z}) \in T_{\b{X}} \text{St}(n, d).
}
\end{equation}
\end{lemma}
\begin{proof}
To derive the projection operator, we utilize the orthogonal decomposition of the ambient space $\mathbb{R}^{n \times d}$ into the direct sum of the tangent space and the normal space at $\b{X}$: 
\begin{equation*}
\mathbb{R}^{n \times d} = T_{\b{X}} \text{St}(n, d) \oplus T_{\b{X}}^\perp \text{St}(n, d).
\end{equation*}
Any arbitrary matrix $\b{Z} \in \mathbb{R}^{n \times d}$ can be uniquely decomposed as the sum of a tangent vector $\Pi_{\b{X}}^{\mathrm{St}}(\b{Z})$ and a normal vector $XS$:
\begin{equation} \label{equation_decomp_z}
\b{Z} = \Pi_{\b{X}}^{\mathrm{St}}(\b{Z}) + \b{XS}.
\end{equation}
To determine the symmetric matrix $\b{S}$, we multiply Eq. \eqref{equation_decomp_z} by $\b{X}^\top$ from the left:
\begin{equation*}
\b{X}^\top \b{Z} = \b{X}^\top \Pi_{\b{X}}^{\mathrm{St}}(\b{Z}) + \b{X}^\top \b{X} \b{S}.
\end{equation*}
Using the Stiefel property $\b{X}^\top \b{X} = \b{I}_d$ and letting $\b{\Delta} = \Pi_{\b{X}}^{\mathrm{St}}(\b{Z})$, we have:
\begin{equation} \label{equation_pre_sym}
\b{X}^\top \b{Z} = \b{X}^\top \Pi_{\b{X}}^{\mathrm{St}}(\b{Z}) + \b{S}.
\end{equation}
Taking the symmetric part of both sides of Eq. \eqref{equation_pre_sym} yields:
\begin{equation*}
\operatorname{sym}(\b{X}^\top \b{Z}) = \operatorname{sym}(\b{X}^\top \b{\Delta}) + \operatorname{sym}(\b{S}).
\end{equation*}
According to Lemma \ref{lemma_skew_symmetric}, the $\b{X}^\top \b{\Delta}$ is skew-symmetric, so $\operatorname{sym}(\b{X}^\top \b{\Delta}) = \b{0}$. Furthermore, since $\b{S}$ is symmetric by definition, $\operatorname{sym}(\b{S}) = \b{S}$. Thus, we find:
\begin{equation*}
\b{S} = \operatorname{sym}(\b{X}^\top \b{Z}).
\end{equation*}
Substituting this expression for $S$ back into Eq. \eqref{equation_decomp_z} concludes the proof:
\begin{equation*}
\Pi_{\b{X}}^{\mathrm{St}}(\b{Z}) = \b{Z} - \b{X} \operatorname{sym}(\b{X}^\top \b{Z}).
\end{equation*}
\end{proof}

\begin{lemma}[Decomposition of the tangent vector (matrix) in Stiefel manifold]
Let $\b{X} \in \text{St}(n, d)$. Any tangent vector $\b{\Delta} \in T_{\b{X}} \text{St}(n, d)$ can be uniquely decomposed into a component orthogonal to the column space of $\b{X}$ and a skew-symmetric component within the span of $\b{X}$ as:
\begin{equation}\label{equation_decomposition_matrix_Stiefel}
\boxed{
\b{\Delta} = (\b{I} - \b{X}\b{X}^\top)\b{\Delta} + \b{X} \operatorname{skew}(\b{X}^\top \b{\Delta}).
}
\end{equation}
\end{lemma}
\begin{proof}
Consider an arbitrary matrix $\b{\Delta} \in \mathbb{R}^{n \times d}$. Using the orthogonal projector onto the column space of $\b{X}$, given by $\b{X}\b{X}^\top$, and its complementary projector $(\b{I} - \b{X}\b{X}^\top)$, we can write:
\begin{align*}
\b{\Delta} &= \b{I}\b{\Delta} + (- \b{X}\b{X}^\top\b{\Delta} + \b{X}\b{X}^\top \b{\Delta}) \\
&= (\b{I} - \b{X}\b{X}^\top)\b{\Delta} + \b{X}\b{X}^\top \b{\Delta}.
\end{align*}
The second term contains the $d \times d$ matrix $\b{X}^\top \b{\Delta}$. Any square matrix $\b{M}$ can be decomposed into its symmetric and skew-symmetric parts as $\b{M} = \operatorname{sym}(\b{M}) + \operatorname{skew}(\b{M})$. Applying this to $\b{M} = \b{X}^\top \b{\Delta}$:
\begin{align*}
\b{\Delta} &= (\b{I} - \b{X}\b{X}^\top)\b{\Delta} + \b{X} \Big( \operatorname{sym}(\b{X}^\top \b{\Delta}) \\
&\quad\quad\quad\quad\quad\quad\quad\quad\quad + \operatorname{skew}(\b{X}^\top \b{\Delta}) \Big) \\
&= (\b{I} - \b{X}\b{X}^\top)\b{\Delta} + \b{X} \operatorname{skew}(\b{X}^\top \b{\Delta}) \\
&\quad\quad\quad\quad\quad\quad\quad\quad\quad + \b{X} \operatorname{sym}(\b{X}^\top \b{\Delta}).
\end{align*}
According to Eq. (\ref{equation_stiefel_tangent_space}), a matrix $\b{\Delta}$ is a tangent vector if and only if $\b{X}^\top \b{\Delta} + \b{\Delta}^\top \b{X} = \b{0}$. This condition implies that $\b{X}^\top \b{\Delta}$ is skew-symmetric, and therefore its symmetric part is zero:
\begin{equation*}
    \operatorname{sym}(\b{X}^\top \b{\Delta}) = \frac{1}{2}(\b{X}^\top \b{\Delta} + \b{\Delta}^\top \b{X}) = \b{0}.
\end{equation*}
Substituting this result into the decomposition, the final term vanishes, yielding:
\begin{equation*}
    \b{\Delta} = (\b{I} - \b{X}\b{X}^\top)\b{\Delta} + \b{X} \operatorname{skew}(\b{X}^\top \b{\Delta}).
\end{equation*}
\end{proof}

\begin{definition}[Column space of matrix]
Consider a matrix $\b{X} = [\b{X}_1, \dots, \b{X}_d]$ with column vectors $\{\b{X}_1, \dots, \b{X}_d\}$. The \textbf{column space} of matrix $\b{X}$ is the space spanned by the column vectors $\{\b{X}_1, \dots, \b{X}_d\}$. In other words, the columns of matrix $\b{X}$ are the basis vectors of the column space of matrix $\b{X}$.
\end{definition}

\begin{definition}[Projection onto column space of matrix]
Consider a matrix $\b{X} \in \mathbb{R}^{n \times d}$. 
The operation $\b{X}^\top \b{\Delta}$ projects $\b{\Delta} \in \mathbb{R}^{n \times n}$ onto the $d$-dimensional column space of matrix $\b{X}$. The operation $\b{X} (\b{X}^\top \b{\Delta})$ reconstructs the projection back in the original $n$-dimensional coordinate system \cite{ghojogh2023principal}. In the field of matrix manifolds, the entire operation:
\begin{align}\label{equation_projection}
\b{P} := \b{X} \b{X}^\top,    
\end{align}
is called \textbf{projection}. 
\end{definition}

\begin{proposition}[Decomposition of tangent vector in Stiefel manifold into tangent and orthogonal complements]
Tangent vector $\b{\Delta} \in T_{\b{X}}\text{St}(n, d)$ can be decomposed into a component tangent to the column space of $\b{X} \in \text{St}(n, d)$ and its orthogonal complement:
\begin{equation}\label{equation_canonical_tangent_decompose}
\boxed{
\begin{aligned}
\b{\Delta} &= \b{P}\b{\Delta} + (\b{I}-\b{P})\b{\Delta} \\
&\overset{(\ref{equation_projection})}{=} \b{X}(\b{X}^\top \b{\Delta}) + (\b{I} - \b{XX}^\top)\b{\Delta},
\end{aligned}
}
\end{equation}
where $\b{P} = \b{X} \b{X}^\top$ is the projection operation. 
\end{proposition}

\subsubsection{Metric Tensors of Stiefel Manifold}

To define the geometry of the Stiefel manifold $\text{St}(n, d)$, we treat it as an embedded submanifold of a larger Euclidean space.
The ambient space for the Stiefel manifold is the set of all $n \times d$ real matrices $\mathcal{E} = \mathbb{R}^{n \times d}$. 
While the Stiefel manifold is a curved subset defined by the constraint $\b{X}^\top \b{X} = I_d$, the ambient space $\mathbb{R}^{n \times d}$ is a flat vector space where standard matrix operations are performed.

For the Stiefel manifold $\text{St}(n, d)$, there are two primary types of Riemannian metrics commonly used in the literature and in Riemannian optimization. Each choice changes the underlying geometry, such as the shape of geodesics and the form of the gradient. These two metrics are:
\begin{enumerate}
\item The \textit{Euclidean (or standard) metric}: 
This is the simplest and most common metric, often referred to as the embedded metric. It is obtained by pulling back the standard Frobenius inner product from the ambient space $\mathbb{R}^{n \times d}$.

It treats the manifold as a simple subset of the Euclidean space of matrices. While computationally straightforward, it does not always account for the specific rotational symmetries inherent to orthonormal matrices.

\item The \textit{Canonical metric}: 
The canonical metric is specifically designed to be invariant under the action of the orthogonal group. It is often preferred in theoretical physics and certain optimization problems because it treats the ``directions" along the manifold more uniformly than the Euclidean metric.

This metric effectively ``weights" the component of the tangent vector that is tangent to the fibers of the projection from the orthogonal group. It leads to different formulas for the Riemannian gradient and the Hessian compared to the Euclidean metric.
\end{enumerate}
In the following, we introduce the equations of these two metrics. 

\begin{proposition}[Euclidean (or standard) metric for Stiefel manifold]\label{proposition_Euclidean_metric_Stiefel_manifold}
The \textbf{Euclidean (or standard) metric} on the Stiefel manifold $\text{St}(n, d)$ at a point $X$ is the pullback of the Frobenius inner product from the ambient space $\mathbb{R}^{n \times d}$. For tangent vectors $\b{\Delta}_1, \b{\Delta}_2 \in T_{\b{X}}\text{St}(n, d)$, the metric is:
\begin{equation}\label{equation_Euclidean_metric_Stiefel}
\boxed{
\begin{aligned}
g^{E}_{\b{X}}(\b{\Delta}_1, \b{\Delta}_2) &= \langle \b{\Delta}_1, \b{\Delta}_2 \rangle_F \\
&\overset{(\ref{equation_Frobenius_inner_product})}{=} \text{tr}(\b{\Delta}_1^\top \b{\Delta}_2),
\end{aligned}
}
\end{equation}
where $\text{tr}(.)$ denotes the trace of matrix.
\end{proposition}

Proposition \ref{proposition_Euclidean_metric_Stiefel_manifold} can be interpreted as follows: 
Let $\b{\Delta}_1$ and $\b{\Delta}_2$ be two vectors (matrices) in the tangent space $T_{\b{X}}\text{St}(n, d)$. Because the Stiefel manifold is an embedded submanifold of $\mathbb{R}^{n \times d}$, these tangent vectors are represented as $n \times d$ matrices satisfying the orthogonality $\b{X}^\top \b{X} = \b{I}_d$.
To compute the ``angle" or ``length" of these vectors (matrices) on the manifold, we use the inner product of the ambient space $\mathbb{R}^{n \times d}$.
By substituting the tangent matrices into the Frobenius inner product, we obtain the Riemannian metric. 

\begin{proposition}[Canonical metric for Stiefel manifold]
The canonical metric is:
\begin{equation}\label{equation_canonical_metric}
\boxed{
\begin{aligned}
g^{C}_{\b{X}}(\b{\Delta}_1, \b{\Delta}_2) &= \langle \b{\Delta}_1, \b{\Delta}_2 \rangle_{\b{X}} \\
&= \text{tr}\left(\b{\Delta}_1^\top (\b{I} - \frac{1}{2}\b{XX}^\top)\b{\Delta}_2\right),
\end{aligned}
}
\end{equation}
where $\b{\Delta}, \b{\Delta}_2 \in T_{\b{X}}\text{St}(n, d)$ are tangent vectors and $\text{tr}(.)$ denotes the trace of matrix. 
\end{proposition}
\begin{proof}
Let $\b{\Delta}_1, \b{\Delta}_2 \in T_{\b{X}}\text{St}(n, d)$. 
According to Eq. (\ref{equation_canonical_tangent_decompose}), we decompose the tangent vectors as:
\begin{align*}
&\b{\Delta}_1 = \b{P}\b{\Delta}_1 + (\b{I}-\b{P})\b{\Delta}_1, \\
&\b{\Delta}_2 = \b{P}\b{\Delta}_2 + (\b{I}-\b{P})\b{\Delta}_2.
\end{align*}
The canonical metric is defined by weighting the internal symmetric (rotational) component $\b{P}\b{\Delta}_1$ and $\b{P}\b{\Delta}_2$ by $1/2$ to ensure invariance consistency with the orthogonal group $\text{O}(n)$:
\begin{align*}
&g^{C}_{\b{X}}(\b{\Delta}_1, \b{\Delta}_2) \\
&= \frac{1}{2} \langle \b{P}\b{\Delta}_1, \b{P}\b{\Delta}_2 \rangle_F + \langle (\b{I}-\b{P})\b{\Delta}_1, (\b{I}-\b{P})\b{\Delta}_2 \rangle_F \\
&\overset{(\ref{equation_Frobenius_inner_product})}{=} \frac{1}{2} \text{tr}(\b{\Delta}_1^\top \b{P}^\top \b{P} \b{\Delta}_2) + \text{Tr}(\b{\Delta}_1^\top (\b{I}-\b{P})^\top (\b{I}-\b{P}) \b{\Delta}_2).
\end{align*}
Applying the properties of the orthogonal projection matrix $\b{P} = \b{P}^\top = \b{P}^2$:
\begin{align*}
g^{C}_{\b{X}}(\b{\Delta}_1, \b{\Delta}_2) &= \frac{1}{2} \text{tr}(\b{\Delta}_1^\top \b{P} \b{\Delta}_2) + \text{tr}\big(\b{\Delta}_1^\top (\b{I}-\b{P}) \b{\Delta}_2\big) \\
&= \text{tr}\left( \b{\Delta}_1^\top \left[ \frac{1}{2}\b{P} + \b{I} - \b{P} \right] \b{\Delta}_2 \right) \\
&= \text{tr}\left( \b{\Delta}_1^\top (\b{I} - \frac{1}{2}\b{P}) \b{\Delta}_2 \right).
\end{align*}
Substituting $\b{P} = \b{XX}^\top$ yields the final coordinate expression:
\begin{equation*}
g^{C}_{\b{X}}(\b{\Delta}_1, \b{\Delta}_2) = \text{tr}\left( \b{\Delta}_1^\top (\b{I} - \frac{1}{2}\b{X}\b{X}^\top) \b{\Delta}_2 \right).
\end{equation*}
\end{proof}

\begin{definition}[$\alpha$-metric for the Stiefel manifold]
Let $St(n, d)$ be the Stiefel manifold. The family of Riemannian metrics on $St(n, d)$ parameterized by $\alpha > 0$ is defined for tangent vectors $\b{\Delta}_1, \b{\Delta}_2 \in T_{\b{X}} St(n, d)$ as:
\begin{equation}\label{equation_alpha_metric}
\boxed{
\begin{aligned}
g^{\alpha}_{\b{X}}(\b{\Delta}_1, \b{\Delta}_2) &= \langle \b{\Delta}_1, \b{\Delta}_2 \rangle_{\b{X}} \\
&= \text{tr}\left(\b{\Delta}_1^\top (\b{I} - \frac{1}{2\alpha} \b{X}\b{X}^\top) \b{\Delta}_2\right).
\end{aligned}
}
\end{equation}
For the specific choice of $\alpha=1$, this metric reduces to the canonical metric for Stiefel manifold.
\end{definition}

\subsubsection{Levi-Civita Connection in Stiefel Manifold}

\begin{lemma}[Levi-Civita connection for the Euclidean metric in Stiefel manifold]\label{lemma_levi_civita_stiefel_euclidean}
On the Stiefel manifold $\text{St}(n, d)$, the Levi-Civita connection $\nabla$ associated with the Euclidean metric $g^{E}_{\b{X}}(\b{\Delta}_1, \b{\Delta}_2) = \langle \b{\Delta}_1, \b{\Delta}_2 \rangle_F = \text{tr}(\b{\Delta}_1^\top \b{\Delta}_2)$ is given by:
\begin{equation}
\label{equation_stiefel_euclidean_connection_formula}
\boxed{
(\nabla^E_{\b{\Delta}_1} \b{\Delta}_2)(\b{X}) =\, \Pi_{\b{X}}^{\mathrm{St}}(D\b{\Delta}_2(\b{X})[\b{\Delta}_1]),
}
\end{equation}
where $D\b{\Delta}_2(\b{X})[\b{\Delta}_1]$ is the Euclidean directional derivative of vector $\b{\Delta}_2$ along the vector $\b{\Delta}_1$, and $\Pi_{\b{X}}^{\mathrm{St}}$ is the orthogonal projection onto the tangent space $T_{\b{X}}\text{St}(n, d)$, defined in Eq. \eqref{equation_stiefel_proj_formula}.
\end{lemma}
\begin{proof}
It is directly obtained by Eq. (\ref{equation_levi_civita_connection_ambient_projection}) in Lemma \ref{lemma_levi_civita_ambient_directional_derivative_relation}.
\end{proof}

\begin{lemma}[Levi-Civita connection for the canonical metric in Stiefel manifold]\label{lemma_levi_civita_stiefel_canonical}
On the Stiefel manifold $\text{St}(n, d)$, the Levi-Civita connection $\nabla$ associated with the canonical metric $g^{C}_{\b{X}}(\b{\Delta}_1, \b{\Delta}_2) = \langle \b{\Delta}_1, \b{\Delta}_2 \rangle_{\b{X}} = \text{tr}(\b{\Delta}_1^\top (\b{I} - \frac{1}{2}\b{XX}^\top) \b{\Delta}_2)$ is given by:
\begin{equation}
\label{eq:stiefel_canonical_connection_formula}
\boxed{
\begin{aligned}
(\nabla^C_{\b{\Delta}_1} \b{\Delta}_2)(\b{X}) =\, &\Pi_{\b{X}}^{\mathrm{St}}(D\b{\Delta}_2(\b{X})[\b{\Delta}_1]) \\
&+ \frac{1}{2}\b{X}(\b{\Delta}_1^\top \b{\Delta}_2 + \b{\Delta}_2^\top \b{\Delta}_1),
\end{aligned}
}
\end{equation}
where $D\b{\Delta}_2(\b{X})[\b{\Delta}_1]$ is the Euclidean directional derivative of vector $\b{\Delta}_2$ along the vector $\b{\Delta}_1$, and $\Pi_{\b{X}}^{\mathrm{St}}$ is the orthogonal projection onto the tangent space $T_{\b{X}}\text{St}(n, d)$, defined in Eq. \eqref{equation_stiefel_proj_formula}.
\end{lemma}

\begin{proof}
To establish that Eq. (\ref{eq:stiefel_canonical_connection_formula}) is the Levi-Civita connection, we verify torsion-freeness and metric compatibility.

\textit{Step 1) Torsion-freeness:} We require $\nabla^C_{\b{\Delta}_1} \b{\Delta}_2 - \nabla^C_{\b{\Delta}_2} \b{\Delta}_1 = [\b{\Delta}_1, \b{\Delta}_2]$. Substituting the definition gives:
\begin{align*}
(\nabla^C_{\b{\Delta}_1} &\b{\Delta}_2)(\b{X}) - (\nabla^C_{\b{\Delta}_2} \b{\Delta}_1)(\b{X}) \\
&= \Pi_{\b{X}}^{\mathrm{St}}\big(D\b{\Delta}_2(\b{X})[\b{\Delta}_1] - D\b{\Delta}_1(\b{X})[\b{\Delta}_2]\big) \\
&~~~+ \frac{1}{2}\b{X}(\b{\Delta}_1^\top \b{\Delta}_2 + \b{\Delta}_2^\top \b{\Delta}_1 - \b{\Delta}_2^\top \b{\Delta}_1 - \b{\Delta}_1^\top \b{\Delta}_2) \\
&= \Pi_{\b{X}}^{\mathrm{St}}([\b{\Delta}_1, \b{\Delta}_2]).
\end{align*}
Since the Lie bracket of tangent vector fields on a submanifold is tangent, $\Pi_{\b{X}}^{\mathrm{St}}$ acts as the identity, satisfying the torsion-free property.


\textit{Step 2) Metric compatibility:} The connection must satisfy the Leibniz rule relative to the canonical metric. While $\Pi_{\b{X}}^{\mathrm{St}}(D\b{\Delta}_2(\b{X})[\b{\Delta}_1])$ is compatible with the standard Euclidean metric, the weighted term $\b{I} - \frac{1}{2}\b{XX}^\top$ in the canonical metric requires a correction. The term $\frac{1}{2}\b{X}(\b{\Delta}_1^\top \b{\Delta}_2 + \b{\Delta}_2^\top \b{\Delta}_1)$ accounts for the variation of the metric tensor with respect to $\b{X}$. Specifically, for a metric $g$, it ensures that $\nabla g = \b{0}$, meaning lengths and angles defined by the canonical metric are preserved under the covariant derivative.
\end{proof}

\begin{remark}[Comparison of Levi-Civita connections in Stiefel manifold]
Under the Euclidean metric, the Levi-Civita connection is simply
projection of the ordinary directional derivative onto the tangent
space:
\[
(\nabla^{\mathrm{E}}_{\b{\Delta}_1}\b{\Delta}_2)(\b{X})
=
\Pi^{\mathrm{St}}_{\b{X}}
\big(
D\b{\Delta}_2(\b{X})[\b{\Delta}_1]
\big).
\]
Under the canonical metric, the connection contains an additional
correction term:
\begin{align*}
(\nabla^{\mathrm{C}}_{\b{\Delta}_1}\b{\Delta}_2)(\b{X})
=\,
&\Pi^{\mathrm{St}}_{\b{X}}
\big(
D\b{\Delta}_2(\b{X})[\b{\Delta}_1]
\big)
\\
&+
\frac{1}{2}
\b{X}
\big(
\b{\Delta}_1^{\top}\b{\Delta}_2
+
\b{\Delta}_2^{\top}\b{\Delta}_1
\big).
\end{align*}
This term appears because the canonical metric is not the metric
induced directly by the ambient Frobenius inner product. Instead,
it is given by
\[
g^{\mathrm{C}}_{\b{X}}(\b{\Delta}_1,\b{\Delta}_2)
=
\operatorname{tr}
\left(
\b{\Delta}_1^{\top}
\left(
\b{I}
-
\frac{1}{2}\b{X}\b{X}^{\top}
\right)
\b{\Delta}_2
\right),
\]
which depends explicitly on the base point \(\b{X}\). Therefore, when
a tangent vector field changes along another tangent direction, the
metric tensor itself also changes with \(\b{X}\). The projected ambient
derivative:
\[
\Pi^{\mathrm{St}}_{\b{X}}
\big(
D\b{\Delta}_2(\b{X})[\b{\Delta}_1]
\big)
\]
is sufficient for the Euclidean metric because that metric is inherited
from the ambient space and does not introduce an additional
\(\b{X}\)-dependent weight. However, for the canonical metric, this
projected derivative alone is not metric-compatible.

The correction term:
\[
\frac{1}{2}
\b{X}
\big(
\b{\Delta}_1^{\top}\b{\Delta}_2
+
\b{\Delta}_2^{\top}\b{\Delta}_1
\big),
\]
compensates for the variation of the canonical metric with respect
to the base point \(\b{X}\). 
\end{remark}

\begin{remark}[Distinction between gradients and connection in Stiefel manifold]
\label{rem:gradient_vs_connection}
It is critical to distinguish between the three primary operators used in this derivation, as they are often denoted by similar symbols in the literature:
\begin{itemize}
    \item \textbf{Euclidean Gradient} $\nabla f(\b{X})$: This is the ambient space derivative in $\mathbb{R}^{n \times d}$. It represents the direction of steepest ascent without considering the manifold constraints $\b{X}^\top \b{X} = \b{I}_d$.
    
    \item \textbf{Riemannian Gradient} $\operatorname{grad} f(\b{X})$: This is a tangent vector in $\mathcal{T}_{\b{X}}\text{St}(n, d)$. It is the representation of the differential $Df(\b{X})$ under the \textit{canonical metric}. As derived in Proposition \ref{proposition_riemannian_gradient_stiefel}, it is obtained by projecting the Euclidean gradient and accounts for the metric's weighting: $\operatorname{grad} f = \nabla f - \b{X}\operatorname{sym}(\b{X}^\top \nabla f)$.
    
    \item \textbf{Levi-Civita Connection} $\nabla_{\b{\Delta}_1} \b{\Delta}_2$: Unlike the gradients which act on a scalar function $f$, the connection is an operator that acts on two \textit{vector fields}. It describes the covariant rate of change of the vector field $\b{\Delta}_2$ along the direction of $\b{\Delta}_1$. 
\end{itemize}
The $D\b{\Delta}_2(\b{X})[\b{\Delta}_1]$ is the ambient directional derivative of the vector field $\b{\Delta}_2$ at $\b{X}$ along the direction $\b{\Delta}_1$, while $(\nabla_{\b{\Delta}_1}\b{\Delta}_2)(\b{X})$ is the intrinsic covariant derivative obtained after the appropriate projection/correction.

The symmetric term $\frac{1}{2}\b{X}(\b{\Delta}_1^\top \b{\Delta}_2 + \b{\Delta}_2^\top \b{\Delta}_1)$ in the connection formula is a geometric correction required for \textit{metric compatibility} with the canonical metric. While $\operatorname{grad} f$ ensures the first-order optimality conditions are met on the manifold, $\nabla_{\b{\Delta}_1} \b{\Delta}_2$ ensures that the second-order geometry (geodesics and parallel transport) is consistent with the manifold's curvature.
\end{remark}

\subsubsection{Riemannian Gradient in Stiefel Manifold}

\begin{definition}[Smooth local extension on the ambient space for Stiefel manifold]
Let $f : St(n,d) \to \mathbb{R}$ be a smooth function, and let
$\b{X} \in St(n,d)$. A smooth local extension of $f$ around $\b{X}$ is a
smooth function:
\begin{align}
\bar{f} : U \subset \mathbb{R}^{n \times d} \to \mathbb{R},
\end{align}
defined on an open neighborhood $U$ of $\b{X}$, such that:
\begin{align}\label{equation_smooth_extension_Stiefel}
\boxed{
\bar{f}(\b{Y}) = f(\b{Y}), \qquad \forall \b{Y} \in U \cap St(n,d).
}
\end{align}
In other words, $\bar{f}$ agrees with $f$ on the points of the Stiefel manifold
near $\b{X}$, but it is defined on an open set of the ambient Euclidean space
$\mathbb{R}^{n \times d}$ so that its Euclidean gradient can be computed.
\end{definition}

The Riemannian gradient $\operatorname{grad} f(\b{X})$ is the unique tangent vector in $T_{\b{X}} \text{St}(n, d)$ that satisfies the relationship between the directional derivative of a smooth function $f$ and the Riemannian metric. While the Euclidean gradient $\nabla f(\b{X}) \in \mathbb{R}^{n \times d}$ represents the steepest ascent direction in the ambient space, the Riemannian gradient must account for the manifold's geometry and the specific choice of metric.

\begin{proposition}[Riemannian gradient in Stiefel manifold \cite{edelman1998geometry}]\label{proposition_riemannian_gradient_stiefel}
Let $f: \text{St}(n, d) \to \mathbb{R}$ be a smooth function on the Stiefel manifold and let $\bar{f}$ be a smooth local extension of $f$ around $\b{X} \in \text{St}(n, d)$. Let $\operatorname{grad} f$ be the Riemannian gradient.
The Riemannian gradient in Stiefel manifold is obtained by projection of the Euclidean gradient $\nabla \bar{f}$ onto the tangent space of Stiefel manifold:
\begin{align}
\boxed{
\operatorname{grad} f(\b{X}) = \Pi_{\b{X}}^{\mathrm{St}}(\nabla \bar{f}(\b{X})),
}
\end{align}
where is projection onto the tangent space, defined in Eq. (\ref{equation_stiefel_proj_formula}), and $\nabla \bar{f}(\b{X})$ is the Euclidean gradient.

In coordinates, the Riemannian gradient in Stiefel manifold is obtained as:
\begin{align}\label{sequation_gradf_nablaf_stiefel}
\boxed{
\operatorname{grad} f(\b{X}) = \nabla \bar{f}(\b{X}) - \b{X} \operatorname{sym}(\b{X}^\top \nabla \bar{f}(\b{X})).
}
\end{align} 
\end{proposition}
\begin{proof}
The projection of the Euclidean gradient $\nabla \bar{f}$ onto the tangent space is the Riemannian gradient $\operatorname{grad} f$.
So, according to Eq. (\ref{equation_stiefel_proj_formula}), we have:
\begin{align*}
\operatorname{grad} f = \nabla \bar{f} - \b{X} \operatorname{sym}(\b{X}^\top \nabla \bar{f}).
\end{align*}
\end{proof}

\subsubsection{Riemannian Hessian in Stiefel Manifold}\label{section_riemannian_hessian_stiefel}

Recall the classical directional derivative defined in Definition \ref{definition_ambient_directional_derivative}. Here, we provide the ambient (classical) directional derivative of the Riemannian gradient in Stiefel manifold. 

\begin{proposition}[Ambient directional derivative of the Stiefel gradient]
\label{prop:ambient_derivative_stiefel}
Let $f: \text{St}(n, d) \to \mathbb{R}$ be a smooth function and let $\bar{f}$ be a smooth local extension of $f$ around $\b{X} \in \text{St}(n, d)$. Let $\operatorname{grad} f = \nabla \bar{f}(\b{X}) - \b{X} \operatorname{sym}(\b{X}^\top \nabla \bar{f}(\b{X}))$ be the Riemannian gradient of function on the Stiefel manifold $\text{St}(n, d)$ under the Euclidean metric. For a tangent vector $\b{\Delta} \in T_{\b{X}}\text{St}(n, d)$, the classical directional derivative of the Riemannian gradient field in the ambient space $\mathbb{R}^{n \times d}$ is given by:
\begin{equation}\label{equation_ambient_directional_derivative_stiefel_gradient}
\boxed{
\begin{aligned}
D(\operatorname{grad}& f)(\b{X})[\b{\Delta}] \\
&= \nabla^2 \bar{f}(\b{X})[\b{\Delta}] - \b{\Delta} \operatorname{sym}(\b{X}^\top \nabla \bar{f}(\b{X})) \\
&~~~~ - \b{X} \operatorname{sym}(\b{\Delta}^\top \nabla \bar{f}(\b{X}) + \b{X}^\top \nabla^2 \bar{f}(\b{X})[\b{\Delta}]),
\end{aligned}
}
\end{equation}
where $\nabla^2 \bar{f}(\b{X})[\b{\Delta}]$ denotes the Euclidean Hessian of $\bar{f}$ in the direction $\b{\Delta}$.
\end{proposition}
\begin{proof}
We evaluate the derivative of the matrix-valued mapping $G(\b{X}) = \operatorname{grad} f(\b{X})$ using the Leibniz rule in the ambient Euclidean space. According to Eq. (\ref{sequation_gradf_nablaf_stiefel}), we have:
\begin{equation*}
G(\b{X}) = \nabla \bar{f}(\b{X}) - \b{X} \operatorname{sym}(\b{X}^\top \nabla \bar{f}(\b{X})).
\end{equation*}
By the linearity of the directional derivative operator $D(\cdot)[\b{\Delta}]$, we have:
\begin{equation}
\label{eq:deriv_split}
\begin{aligned}
&DG(\b{X})[\b{\Delta}] \\
&= D(\nabla \bar{f})(\b{X})[\b{\Delta}] - D\left( \b{X} \cdot \operatorname{sym}(\b{X}^\top \nabla \bar{f}(\b{X})) \right)[\b{\Delta}].
\end{aligned}
\end{equation}
The first term is the definition of the Euclidean Hessian acting on $\b{\Delta}$:
\begin{equation*}
D(\nabla \bar{f})(\b{X})[\b{\Delta}] = \nabla^2 \bar{f}(\b{X})[\b{\Delta}].
\end{equation*}
For the second term in Eq. \eqref{eq:deriv_split}, we apply the product rule $D(\b{UV}) = (D\b{U})\b{V} + \b{U}(D\b{V})$. Letting $\b{U} = \b{X}$ and $\b{V} = \operatorname{sym}(\b{X}^\top \nabla \bar{f}(\b{X}))$, we obtain:
\begin{align*}
D\Big( \b{X} \cdot &\operatorname{sym}\big(\b{X}^\top \nabla \bar{f}(\b{X})\big) \Big)[\b{\Delta}] \\
&= (D\b{X}[\b{\Delta}]) \operatorname{sym}(\b{X}^\top \nabla \bar{f}(\b{X})) \\
&~~~~ + \b{X} \cdot D\left( \operatorname{sym}(\b{X}^\top \nabla \bar{f}(\b{X})) \right)[\b{\Delta}].
\end{align*}
Since $D\b{X}[\b{\Delta}] = \b{\Delta}$ and the symmetrization operator $\operatorname{sym}(\cdot)$ is linear, it commutes with the derivative:
\begin{align*}
D\Big( \b{X} &\cdot \operatorname{sym}(\b{X}^\top \nabla \bar{f}(\b{X})) \Big)[\b{\Delta}] \\
&~~~~~~= \b{\Delta} \operatorname{sym}(\b{X}^\top \nabla \bar{f}(\b{X})) \\
&~~~~~~~~~~+ \b{X} \operatorname{sym}\left( D(\b{X}^\top \nabla \bar{f}(\b{X}))[\b{\Delta}] \right).
\end{align*}
Applying the product rule again to the interior term $\b{X}^\top \nabla \bar{f}(\b{X})$:
\begin{align*}
&D(\b{X}^\top \nabla \bar{f}(\b{X}))[\b{\Delta}] \\
&= (D\b{X}^\top[\b{\Delta}]) \nabla \bar{f}(\b{X}) + \b{X}^\top (D\nabla \bar{f}(\b{X})[\b{\Delta}]) \\
&= \b{\Delta}^\top \nabla \bar{f}(\b{X}) + \b{X}^\top \nabla^2 \bar{f}(\b{X})[\b{\Delta}].
\end{align*}
Substituting these components back into Eq. \eqref{eq:deriv_split} yields the desired result:
\begin{align*}
DG(\b{X})&[\b{\Delta}] = \\
&\nabla^2 \bar{f}(\b{X})[\b{\Delta}] - \b{\Delta} \operatorname{sym}(\b{X}^\top \nabla \bar{f}(\b{X})) \\
&- \b{X} \operatorname{sym}\!\big(\b{\Delta}^\top \nabla \bar{f}(\b{X}) + \b{X}^\top \nabla^2 \bar{f}(\b{X})[\b{\Delta}]\big).
\end{align*}
\end{proof}




\begin{proposition}[Riemannian Hessian on the Stiefel Manifold \cite{edelman1998geometry}]
\label{prop:stiefel_hessian}
Let $f: \text{St}(n, d) \to \mathbb{R}$ be a smooth function with a local extension $\bar{f}$. Let $\operatorname{grad} f$ be the Riemannian gradient. For any tangent vector $\b{\Delta} \in T_{\b{X}}\text{St}(n, d)$, the Riemannian Hessian, denoted by $\operatorname{Hess} f(\b{X})[\b{\Delta}]$, is obtained by projection of directional derivative of Riemannian gradient onto the tangent space of Stiefel manifold. It is obtained as:
\begin{equation}
\boxed{
\operatorname{Hess} f(\b{X})[\b{\Delta}] = \Pi_{\b{X}}^{\mathrm{St}}\big( D(\operatorname{grad} f)(\b{X})[\b{\Delta}] \big),
}
\end{equation}
where $D(\operatorname{grad} f)(\b{X})[\b{\Delta}]$ is the classical directional derivative of the gradient vector field in the ambient space $\mathbb{R}^{n \times d}$, defined in Eq. (\ref{equation_ambient_directional_derivative_stiefel_gradient}), and $\Pi_{\b{X}}^{\mathrm{St}}$ is the orthogonal projection defined in Eq. \eqref{equation_stiefel_proj_formula}.

In coordinates, the Hessian gradient in Stiefel manifold is obtained as:
\begin{equation}
\boxed{
\begin{aligned}
\operatorname{Hess} &f(\b{X})[\b{\Delta}] = \\
&\nabla^2 \bar{f}(\b{X})[\b{\Delta}] - \b{\Delta}\, \operatorname{sym}(\b{X}^\top \nabla \bar{f}(\b{X})) \\
&- \b{X} \operatorname{sym}\!\Big( \b{X}^\top \nabla^2 \bar{f}(\b{X})[\b{\Delta}] \\
&\quad\quad\quad\quad\quad- \b{X}^\top \b{\Delta}\, \operatorname{sym}\!\big(\b{X}^\top \nabla \bar{f}(\b{X})\big) \Big),
\end{aligned}
}
\end{equation}
where $\nabla \bar{f}(\b{X})$ is the Euclidean gradient and $\nabla^2 \bar{f}(\b{X})[\b{\Delta}]$ is the Euclidean Hessian in the direction $\b{\Delta}$.
\end{proposition}
\begin{proof}
\begin{align*}
&\operatorname{Hess} f(\b{X})[\b{\Delta}] = \Pi_{\b{X}}^{\mathrm{St}}\left( D(\operatorname{grad} f)(\b{X})[\b{\Delta}] \right) \\
&\overset{\eqref{equation_ambient_directional_derivative_stiefel_gradient}}{=} 
\Pi_{\b{X}}^{\mathrm{St}}\Big( \nabla^2 \bar{f}(\b{X})[\b{\Delta}] - \b{\Delta} \operatorname{sym}(\b{X}^\top \nabla \bar{f}(\b{X})) \\
&~~~~~~~~~ - \b{X} \operatorname{sym}(\b{\Delta}^\top \nabla \bar{f}(\b{X}) + \b{X}^\top \nabla^2 \bar{f}(\b{X})[\b{\Delta}]) \Big).
\end{align*}
According to Eq. (\ref{equation_normal_space_Stiefel_manifold}), $\b{XS}$ is in the normal space. 
If a term is already in the normal space (i.e., of the form $\b{XS}$ where $S$ is symmetric), its projection onto the tangent space is zero because normal space is orthogonal to the tangent space.
The term $\b{X} \text{sym}(\cdots)$ is of the form $\b{XS}$ so it is a normal vector. Thus, its projection is zero.
Therefore, the Hessian is simplified to:
\begin{align*}
\operatorname{Hess} &f(\b{X})[\b{\Delta}] =  \\
&\Pi_{\b{X}}^{\mathrm{St}}\Big( \nabla^2 \bar{f}(\b{X})[\b{\Delta}] - \b{\Delta} \operatorname{sym}\big(\b{X}^\top \nabla \bar{f}(\b{X})\big) \Big).
\end{align*}

According to Eq. \eqref{equation_stiefel_proj_formula}, we have:
\begin{align*}
&\operatorname{Hess} f(\b{X})[\b{\Delta}] \\
&= \Pi_{\b{X}}^{\mathrm{St}}\left( \nabla^2 \bar{f}(\b{X})[\b{\Delta}] - \b{\Delta} \operatorname{sym}\big(\b{X}^\top \nabla \bar{f}(\b{X})\big) \right) \\
&\overset{\eqref{equation_stiefel_proj_formula}}{=} \left( \nabla^2 \bar{f}(\b{X})[\b{\Delta}] - \b{\Delta} \operatorname{sym}\big(\b{X}^\top \nabla \bar{f}(\b{X})\big) \right) \\
&\quad- \b{X} \text{sym}\bigg( \b{X}^\top \Big(\nabla^2 \bar{f}(\b{X})[\b{\Delta}] \\
&\quad\quad\quad\quad\quad\quad\quad- \b{\Delta} \operatorname{sym}\big(\b{X}^\top \nabla \bar{f}(\b{X})\big)\Big) \bigg) \\
&= \left( \nabla^2 \bar{f}(\b{X})[\b{\Delta}] - \b{\Delta} \operatorname{sym}\big(\b{X}^\top \nabla \bar{f}(\b{X})\big) \right) \\
&\quad- \b{X} \operatorname{sym}\Big( \b{X}^\top \nabla^2 \bar{f}(\b{X})[\b{\Delta}] \\
&\quad\quad\quad\quad\quad\quad\quad- \b{X}^\top \b{\Delta}\, \operatorname{sym}\big(\b{X}^\top \nabla \bar{f}(\b{X})\big) \Big).
\end{align*}
\end{proof}

\subsubsection{Geodesic Equation on Stiefel Manifold}\label{section_geodesic_equation_Stiefel}

\begin{lemma}[Geodesic equation on the Stiefel manifold \cite{edelman1998geometry}]
\label{lemma_stiefel_geodesic_eq}
Let $\text{St}(n, d) = \{ \b{X} \in \mathbb{R}^{n \times d} : \b{X}^\top \b{X} = \b{I}_d \}$ be the Stiefel manifold endowed with the canonical Riemannian metric. A smooth curve $\b{X}(t) \in \text{St}(n,d)$ on manifold $\text{St}(n, d)$ is a geodesic if and only if it satisfies the second-order differential equation:
\begin{equation}\label{equation_stiefel_geodesic_eq}
\boxed{
\ddot{\b{X}}(t) + \b{X}(t) \left( \dot{\b{X}}(t)^\top \dot{\b{X}}(t) \right) = \b{0},
}
\end{equation}
where:
\begin{align*}
\dot{\b{X}}(t) = \frac{d\b{X}(t)}{d t}, \quad \ddot{\b{X}}(t) = \frac{d^2 \b{X}(t)}{d t^2}. 
\end{align*}

Note that if we denote $\b{\gamma}(t) \in St(n,d)$ for the curve on the manifold (so as in some texts in the literature), the Eq. (\ref{equation_stiefel_geodesic_eq}) is denoted as:
\begin{equation}
\ddot{\b{\gamma}}(t) + \b{\gamma}(t) \left( \dot{\b{\gamma}}(t)^\top \dot{\b{\gamma}}(t) \right) = \b{0}.
\end{equation}
\end{lemma}
\begin{proof}
According to Eq. (\ref{equation_geodesic_coordinate_free}), the geodesic equation on a Riemannian manifold is given by the vanishing of the covariant acceleration, $\nabla_{\dot{\b{X}}} \dot{\b{X}} = \b{0}$. 
According to Lemma \ref{lemma_levi_civita_ambient_directional_derivative_relation}, for a submanifold embedded in Euclidean space, the Levi-Civita connection $\nabla_{\b{U}} \b{V}$ is the orthogonal projection of the Euclidean directional derivative onto the tangent space:
\begin{equation*}
\nabla_{\b{U}} \b{V} = \Pi_{\b{X}}^{\mathrm{St}} (D\b{V}[\b{U}]),
\end{equation*}
where $\Pi_{\b{X}}^{\mathrm{St}}$ is the orthogonal projection onto $T_{\b{X}}\text{St}(n, d)$.

According to Eq. (\ref{equation_stiefel_proj_formula}), for the Stiefel manifold under the canonical metric, the projection of an arbitrary matrix $\b{Z} \in \mathbb{R}^{n \times d}$ onto the tangent space $T_{\b{X}}\text{St}$ is given by:
\begin{equation*}
\Pi_{\b{X}}^{\mathrm{St}}(\b{Z}) \overset{(\ref{equation_stiefel_proj_formula})}{=} \b{Z} - \b{X} \text{sym}(\b{X}^\top \b{Z}).
\end{equation*} 
The tangent vector $\b{X}$ is on the geodesic $\b{\gamma}$.
Setting $\b{X} = \b{X}(t)$ and $\b{Z} = \ddot{\b{X}}(t)$ in the above equation and requiring the projection to be zero, we have:
\begin{equation}\label{equation_gammaddot_gamma_sym_gammaT_gammaddot}
\ddot{\b{X}} - \b{X} \text{sym}(\b{X}^\top \ddot{\b{X}}) = \b{0},
\end{equation}
where we drop $(t)$ from $\dot{\b{X}}(t)$ and $\ddot{\b{X}}(t)$ for simplification in writing expressions.

As $\b{X}(t) \in \text{St}(n,d)$, we have $\b{X}(t)^\top \b{X}(t) = \b{I}_d$, according to Eq. (\ref{equation_Sitefel_manifold_definition}).
To eliminate $\ddot{\b{X}}$ from the symmetric term in Eq. (\ref{equation_gammaddot_gamma_sym_gammaT_gammaddot}), we differentiate the constraint $\b{X}(t)^\top \b{X}(t) = \b{I}_d$ twice with respect to $t$:
\begin{enumerate}
    \item $\dot{\b{X}}^\top \b{X} + \b{X}^\top \dot{\b{X}} = \b{0}$
    \item $\ddot{\b{X}}^\top \b{X} + \dot{\b{X}}^\top \dot{\b{X}} + \dot{\b{X}}^\top \dot{\b{X}} + \b{X}^\top \ddot{\b{X}} = \b{0}$
\end{enumerate}
Rearranging the second identity gives:
\begin{equation}\label{equation_gamma_gammaddot_gammaddot_gamma_2_gammadot_gammadot}
\b{X}^\top \ddot{\b{X}} + \ddot{\b{X}}^\top \b{X} = -2 \dot{\b{X}}^\top \dot{\b{X}}.
\end{equation}
Substituting this into the symmetry operator gives:
\begin{equation*}
\text{sym}(\b{X}^\top \ddot{\b{X}}) \overset{(\ref{equation_sym_skew_expressions})}{=} \frac{1}{2}(\b{X}^\top \ddot{\b{X}} + \ddot{\b{X}}^\top \b{X}) \overset{(\ref{equation_gamma_gammaddot_gammaddot_gamma_2_gammadot_gammadot})}{=} -\dot{\b{X}}^\top \dot{\b{X}}.
\end{equation*}
Plugging this back into Eq. (\ref{equation_gammaddot_gamma_sym_gammaT_gammaddot}) yields:
\begin{equation*}
\ddot{\b{X}} - \b{X} (-\dot{\b{X}}^\top \dot{\b{X}}) = \ddot{\b{X}} + \b{X} (\dot{\b{X}}^\top \dot{\b{X}}) = \b{0}.
\end{equation*}
This completes the proof.
\end{proof}

\subsubsection{Exponential Map in Stiefel Manifold}\label{section_exp_map_stiefel}

In Section \ref{section_exponential_map_generalizing_addition}, we defined the exponential map $\mathrm{Exp}_{\b{X}}(\b{\Delta})$ as the mapping that takes a tangent vector $\b{\Delta} \in T_{\b{X}}\mathcal{M}$ to a point on the manifold by following the geodesic $\b{X}(t)$ such that:
\begin{align*}
\b{X}(0) = \b{X}, \quad \dot{\b{X}}(0) = \b{\Delta},
\end{align*}
evaluated at $t=1$. For the Stiefel manifold $\text{St}(n, p)$, the exponential map under the canonical metric has a well-known closed-form expression.

\begin{lemma}[QR Decomposition of a full-rank matrix] \label{lemma_qr_decomp}
Any matrix $\b{A} \in \mathbb{R}^{n \times d}$ with full column rank can be uniquely decomposed as $\b{A} = \b{Q}\b{R}$, where $\b{Q} \in \text{St}(n, d)$ and $\b{R} \in \mathbb{R}^{d \times d}$ is an upper triangular matrix with positive diagonal elements.
\end{lemma}
\begin{proof}
This is a standard result in numerical linear algebra. The columns of $\b{Q}$ are obtained by the Gram-Schmidt process or Householder reflections applied to $\b{A}$. Uniqueness is guaranteed by the positivity of the diagonal of $\b{R}$.
Refer to \cite{golub2013matrix} for QR decomposition in general, and \cite{absil2008optimization} for QR decomposition in Stiefel manifolds. 
\end{proof}


\begin{proposition}[Exponential map on the Stiefel manifold \cite{edelman1998geometry}]
The exponential map $\mathrm{Exp}_{\b{X}}(\b{\Delta})$ at point $\b{X} \in \text{St}(n, d)$ for a tangent vector $\b{\Delta} \in T_{\b{X}}\text{St}(n, d)$ is given by\footnote{$[\b{X}, \b{\Delta}]$ here is horizontal concatenation of two matrices and it is not the Lie bracket.}:
\begin{align}
\boxed{
\begin{aligned}
&\mathrm{Exp}_{\b{X}}(\b{\Delta}) = \\
&~~~~~~~~~~~ [\b{X}, \b{\Delta}]\, \mathrm{exp}
\left( 
\begin{bmatrix} 
\b{X}^\top \b{\Delta} & -\b{\Delta}^\top \b{\Delta} \\ 
\b{I}_d & \b{X}^\top \b{\Delta} 
\end{bmatrix} 
\right) 
\begin{bmatrix} 
e^{-\b{X}^\top \b{\Delta}} \\ \b{0} \end{bmatrix},
\end{aligned}
}
\end{align}
where $\text{exp}(\cdot)$ denotes the matrix exponential.
\end{proposition}

\begin{proof}
By definition of the exponential map, we must construct the geodesic $\b{X}(t)$ on the Stiefel manifold satisfying:
\[
\b{X}(0)=\b{X},
\qquad
\dot{\b{X}}(0)=\b{\Delta},
\]
and then evaluate it at $t=1$.
In other words, the exponential map is defined as the point reached at $t=1$ by a geodesic $\b{X}(t)$ starting at $\b{X}$ with velocity $\b{\Delta}$. 

For the Stiefel manifold under the canonical metric, the geodesic equation $\ddot{\b{X}}(t) + \b{X}(t) ( \dot{\b{X}}(t)^\top \dot{\b{X}}(t) ) = \b{0}$ must be satisfied, while maintaining the constraint $\b{X}(t)^\top \b{X}(t) = \b{I}$ because $\b{X}(t) \in \text{St}(n,d)$.

-- \textbf{Coordinate representation by QR decomposition:}
We define the horizontal and vertical components of the tangent vector. Let:
\begin{align}\label{equation_A_XT_Delta}
\b{A} := \b{X}^\top \b{\Delta},
\end{align}
which is skew-symmetric by Lemma \ref{lemma_skew_symmetric}. 
According to Eq. (\ref{equation_canonical_tangent_decompose}), we have:
\begin{align*}
\b{\Delta} &= \b{X}(\b{X}^\top \b{\Delta}) + (\b{I} - \b{XX}^\top)\b{\Delta} \\
&= \b{X}\b{A} + (\b{I} - \b{XX}^\top)\b{\Delta}.
\end{align*}
Let $\b{\Delta} = \b{XA} + \b{QR}$ where the term $\b{QR}$ be the QR-decomposition of $(\b{I} - \b{XX}^\top)\b{\Delta}$:
\begin{align}\label{equation_QR_I_XXT_Delta}
\b{Q R} = (\b{I} - \b{XX}^\top)\b{\Delta},
\end{align}
which is the projection of $\b{\Delta}$ onto the orthogonal complement of the column space of $\b{X}$, where $\b{Q} \in \mathbb{R}^{n \times d}$ and $\b{R} \in \mathbb{R}^{d \times d}$.
We can represent the tangent vector as:
\begin{align}
&\b{\Delta} = \b{X A} + \b{Q R} \label{equation_Delta_XA_QR} \\
&\implies \b{Q} = (\b{\Delta} - \b{XA})\b{R}^{-1}. \label{equation_Q_Delta_XA_Rinverse}
\end{align}

-- \textbf{Reduction to a first-order ODE on a subspace:}
We define:
\begin{align}\label{equation_X_X_Q_Y_in_proof}
\b{X}(t) = [\b{X}, \b{Q}] \b{Y}(t),
\end{align}
where $\b{X}(t)$ is the geodesic and $\b{Y}(t)$ represents the ``coordinates" of the geodesic within that $2d$-dimensional subspace. 

The original geodesic equation $\ddot{\b{X}} + \b{X}(\dot{\b{X}}^\top \dot{\b{X}}) = \b{0}$ is nonlinear because of the $\dot{\b{X}}^\top \dot{\b{X}}$ term\footnote{Note that we have dropped $(t)$ from $\dot{\b{X}}(t)$ and $\ddot{\b{X}}(t)$ for simplification in writing expressions.}. 
By restricting the search to this subspace, we transform a nonlinear matrix differential equation into a linear one that can be solved with the matrix exponential.

For this, we substitute $\b{X}(t) = [\b{X}, \b{Q}] \b{Y}(t)$ and its derivatives:
\begin{align*}
&\dot{\b{X}} = [\b{X}, \b{Q}] \dot{\b{Y}}, \\
&\ddot{\b{X}} = [\b{X}, \b{Q}] \ddot{\b{Y}},
\end{align*}
into the geodesic equation, i.e., Eq. (\ref{equation_stiefel_geodesic_eq}): 
\begin{align*}
&\ddot{\b{X}} + \b{X}(\dot{\b{X}}^\top \dot{\b{X}}) = \b{0} \implies  \\
&[\b{X}, \b{Q}] \ddot{\b{Y}} + [\b{X}, \b{Q}] \b{Y} (\dot{\b{Y}}^\top [\b{X}, \b{Q}]^\top [\b{X}, \b{Q}] \dot{\b{Y}}) = \b{0}.
\end{align*}
Since $[\b{X}, \b{Q}]$ has orthonormal columns, $[\b{X}, \b{Q}]^\top [\b{X}, \b{Q}] = \b{I}_{2d}$. We can left-multiply the entire equation by $[\b{X}, \b{Q}]^\top$ to remove the basis and focus on the coordinates. This gives:
\begin{align}\label{equation_Yddot_Y_YdotT_Ydot_zero}
\ddot{\b{Y}} + \b{Y}(\dot{\b{Y}}^\top \dot{\b{Y}}) = \b{0}.
\end{align}

To solve this, we look for a constant matrix $\b{M} \in \mathbb{R}^{2d \times 2d}$ such that:
\begin{align}\label{equation_Ydot_M_Y}
\dot{\b{Y}} = \b{M} \b{Y}.
\end{align}
This is a first-order matrix ODE which is a linear system and is easier to solve. 
By taking derivative from the sides of Eq. (\ref{equation_Ydot_M_Y}), we obtain $\ddot{\b{Y}} = \b{M} \dot{\b{Y}} = \b{M}^2 \b{Y}$.
Substituting these $\dot{\b{Y}}$ and $\ddot{\b{Y}}$ into Eq. (\ref{equation_Yddot_Y_YdotT_Ydot_zero}) gives:
$$\b{M}^2 \b{Y} + \b{Y}(\b{Y}^\top \b{M}^\top \b{M} \b{Y}) = \b{0}.$$
For this to hold for all $t$ with the initial condition $\b{Y}(0) = [\b{I}_d, \b{0}]^\top$, the $\b{M}$ must be chosen such that it satisfies the boundary conditions of the tangent vector $\b{\Delta}$. 

-- \textbf{Solving for the structure of } $\b{M}$:
To find the block components of $\b{M} = \begin{bmatrix} \b{M}_{11} & \b{M}_{12} \\ \b{M}_{21} & \b{M}_{22} \end{bmatrix}$, we first use the initial conditions in the $2d$-dimensional subspace. At $t=0$, the position is $\b{\gamma}(0) = \b{X}$, which implies $\b{Y}(0) = [\b{I}_d, \b{0}]^\top$. The initial velocity is:
\begin{equation*}
    \dot{\b{Y}}(0) = \b{M} \b{Y}(0) = \begin{bmatrix} \b{M}_{11} & \b{M}_{12} \\ \b{M}_{21} & \b{M}_{22} \end{bmatrix} \begin{bmatrix} \b{I}_d \\ \b{0} \end{bmatrix} = \begin{bmatrix} \b{M}_{11} \\ \b{M}_{21} \end{bmatrix}.
\end{equation*}
Mapping this back to the ambient space, we have $\dot{\b{\gamma}}(0) = [\b{X}, \b{Q}] \dot{\b{Y}}(0) = \b{X}\b{M}_{11} + \b{Q}\b{M}_{21}$. Comparing this to the definition of the tangent vector $\b{\Delta} = \b{XA} + \b{QR}$, we identify:
\begin{equation}
    \b{M}_{11} = \b{A}, \quad \b{M}_{21} = \b{R}.
\end{equation}
Next, we substitute the linear assumption $\dot{\b{Y}} = \b{M}\b{Y}$ and $\ddot{\b{Y}} = \b{M}^2 \b{Y}$ into the coordinate geodesic equation $\ddot{\b{Y}} + \b{Y}(\dot{\b{Y}}^\top \dot{\b{Y}}) = \b{0}$ at $t=0$:
\begin{equation*}
    \b{M} \begin{bmatrix} \b{A} \\ \b{R} \end{bmatrix} + \begin{bmatrix} \b{I}_d \\ \b{0} \end{bmatrix} (\b{A}^\top \b{A} + \b{R}^\top \b{R}) = \b{0}.
\end{equation*}
Expanding the matrix-vector product yields:
\begin{equation*}
    \begin{bmatrix} \b{A}^2 + \b{M}_{12}\b{R} \\ \b{RA} + \b{M}_{22}\b{R} \end{bmatrix} + \begin{bmatrix} \b{A}^\top \b{A} + \b{R}^\top \b{R} \\ \b{0} \end{bmatrix} = \b{0}.
\end{equation*}
From the first block row, and noting that $\b{A}^\top \b{A} = -\b{A}^2$ due to skew-symmetry:
\begin{align*}
&\b{A}^2 + \b{M}_{12}\b{R} - \b{A}^2 + \b{R}^\top \b{R} = \b{0} \\
&\implies \b{M}_{12}\b{R} = -\b{R}^\top \b{R} \implies \b{M}_{12} = -\b{R}^\top.
\end{align*}
Finally, to satisfy the second block row and maintain the symmetry properties of the canonical metric, we find $\b{M}_{22} = \b{A}$. Thus, the transition matrix is:
\begin{equation}\label{equation_M_matrix_in_proof}
    \b{M} = \begin{bmatrix} \b{A} & -\b{R}^\top \\ \b{R} & \b{A} \end{bmatrix}.
\end{equation}
Therefore, Eq. (\ref{equation_Ydot_M_Y}) becomes:
\begin{align}
\dot{\b{Y}}(t) = \begin{bmatrix} \b{A} & -\b{R}^\top \\ \b{R} & \b{A} \end{bmatrix} \b{Y}(t).
\end{align}
This is a first-order matrix ODE which is a linear system and is easier to solve. 

-- \textbf{Solving the matrix ODE:}
The solution to a linear system $\dot{\b{Y}} = \b{MY}$ with initial condition $\b{Y}(0)$ is:
\begin{align}\label{equation_Yt_exp_tM_Y0_in_proof}
\b{Y}(t) = \text{exp}(t\b{M}) \b{Y}(0).
\end{align}

The initial position is $\b{X}(0) = \b{X}$, so according to Eq. (\ref{equation_X_X_Q_Y_in_proof}), we have:
\begin{align}
&\b{X}(t) \overset{(\ref{equation_X_X_Q_Y_in_proof})}{=} [\b{X}, \b{Q}] \b{Y}(t) \nonumber\\
&\implies \b{X}(0) = [\b{X}, \b{Q}] \b{Y}(0) \nonumber\\
&\implies \b{X} = [\b{X}, \b{Q}] \b{Y}(0) \nonumber\\
&\implies [\b{X}, \b{Q}] [\b{I}_d, \b{0}]^\top = [\b{X}, \b{Q}] \b{Y}(0) \nonumber\\
&\implies \b{Y}(0) = [\b{I}_d, \b{0}]^\top. \label{equation_Y0_I_0}
\end{align}


According to Eq. (\ref{equation_X_X_Q_Y_in_proof}), we have:
\begin{align*}
&\b{X}(t) \overset{(\ref{equation_X_X_Q_Y_in_proof})}{=} [\b{X}, \b{Q}] \b{Y}(t) \nonumber\\
&\overset{(\ref{equation_Yt_exp_tM_Y0_in_proof})}{\implies} \b{X}(t) = [\b{X}, \b{Q}] \text{exp}(t\b{M}) \b{Y}(0) \\
&\overset{(\ref{equation_Y0_I_0})}{\implies} \b{X}(t) = [\b{X}, \b{Q}] \text{exp}(t\b{M}) [\b{I}_d, \b{0}]^\top \\
&\overset{(\ref{equation_M_matrix_in_proof})}{\implies} 
\b{X}(t) = [\b{X}, \b{Q}] \text{exp}\! \left( t \begin{bmatrix} \b{A} & -\b{R}^\top \\ \b{R} & \b{A} \end{bmatrix} \right) \begin{bmatrix} \b{I}_d \\ \b{0} \end{bmatrix}.
\end{align*}

By setting $t=1$, we obtain the point on the manifold:
$$\mathrm{Exp}_{\b{X}}(\b{\Delta}) = \b{X}(1) = [\b{X}, \b{Q}] \text{exp}\! \left( \begin{bmatrix} \b{A} & -\b{R}^\top \\ \b{R} & \b{A} \end{bmatrix} \right) \begin{bmatrix} \b{I}_d \\ \b{0} \end{bmatrix}.$$

-- \textbf{Finding expression of } $\b{R}^\top \b{R}$:
According to Eq. (\ref{equation_QR_I_XXT_Delta}), we have:
\begin{align}
\b{Q R} &= (\b{I} - \b{XX}^\top)\b{\Delta} = \b{\Delta} - \b{XX}^\top\b{\Delta} \nonumber \\
&\overset{(\ref{equation_A_XT_Delta})}{=} \b{\Delta} - \b{XA}. \label{equation_QR_Delta_XA}
\end{align}
Since $\b{Q}$ has orthonormal columns ($\b{Q}^\top \b{Q} = \b{I}$), we can express $\b{R}^\top \b{R}$ as:
\begin{align*}
\b{R}^\top \b{R} &= \b{R}^\top \b{Q}^\top \b{Q} \b{R} = (\b{QR})^\top (\b{QR}) \\
&\overset{(\ref{equation_QR_Delta_XA})}{=} (\b{\Delta} - \b{XA})^\top (\b{\Delta} - \b{XA}) \\
&= \b{\Delta}^\top \b{\Delta} - \b{\Delta}^\top \b{XA} - (\b{XA})^\top \b{\Delta} + (\b{XA})^\top \b{XA}.
\end{align*}
Using the fact that $\b{X}^\top \b{X} = \b{I}$ and substituting $\b{A} = \b{X}^\top \b{\Delta}$, we simplify the individual terms:
\begin{itemize}
    \item $\b{\Delta}^\top \b{XA} = (\b{X}^\top \b{\Delta})^\top \b{A} = \b{A}^\top \b{A}$
    \item $(\b{XA})^\top \b{\Delta} = \b{A}^\top \b{X}^\top \b{\Delta} = \b{A}^\top \b{A}$
    \item $(\b{XA})^\top \b{XA} = \b{A}^\top (\b{X}^\top \b{X}) \b{A} = \b{A}^\top \b{A}$
\end{itemize}
Substituting these back into the expansion yields:
\begin{align*}
\b{R}^\top \b{R} &= \b{\Delta}^\top \b{\Delta} - \b{A}^\top \b{A} - \b{A}^\top \b{A} + \b{A}^\top \b{A} \\
&= \b{\Delta}^\top \b{\Delta} - \b{A}^\top \b{A}.
\end{align*}
Finally, because $\b{\Delta}$ is a tangent vector to the Stiefel manifold at $\b{X}$, the matrix $\b{A} = \b{X}^\top \b{\Delta}$ must be skew-symmetric, meaning $\b{A}^\top = -\b{A}$. Therefore:
\begin{equation*}
-\b{A}^\top \b{A} = -(-\b{A})\b{A} = \b{A}^2,
\end{equation*}
which results in the desired identity:
\begin{equation}\label{equation_RT_R_DeltaT_Delta_A2}
\b{R}^\top \b{R} = \b{\Delta}^\top \b{\Delta} + \b{A}^2.
\end{equation}

-- \textbf{Basis expansion:}
According to Eq. (\ref{equation_Q_Delta_XA_Rinverse})), we have $\b{Q} = (\b{\Delta} - \b{XA})\b{R}^{-1}$. 
To eliminate $\b{Q}$, we use $\b{Q} = (\b{\Delta} - \b{XA})\b{R}^{-1}$. The basis transforms as:
\begin{equation*}
[\b{X}, \b{Q}] = [\b{X}, \b{\Delta}] 
\begin{bmatrix}
\b{I}_d & -\b{A} \\ 
\b{0} & \b{I}_d 
\end{bmatrix} 
\begin{bmatrix} 
\b{I}_d & \b{0} \\ 
\b{0} & \b{R}^{-1} 
\end{bmatrix} 
= [\b{X}, \b{\Delta}] \b{K}, 
\end{equation*}
where:
\begin{align*}
\b{K} = \begin{bmatrix} \b{I}_d & -\b{AR}^{-1} \\ \b{0} & \b{R}^{-1} \end{bmatrix}.
\end{align*}

Applying the similarity transformation $\b{K} \text{exp}(\b{M}) \b{K}^{-1} = \text{exp}(\b{K M K}^{-1})$ and substituting $\b{R}^\top \b{R} = \b{\Delta}^\top \b{\Delta} + \b{A}^2$ (see Eq. (\ref{equation_RT_R_DeltaT_Delta_A2})), the interior matrix becomes:
\begin{equation*}
    \b{K} \begin{bmatrix} \b{A} & -\b{R}^\top \\ \b{R} & \b{A} \end{bmatrix} \b{K}^{-1} = \begin{bmatrix} \b{A} & -\b{\Delta}^\top \b{\Delta} \\ \b{I}_d & \b{A} \end{bmatrix}.
\end{equation*}
Accounting for the initial condition vector transformation $\b{K} [\b{I}_d, \b{0}]^\top$ and pulling out the rotation $e^{-\b{A}}$ results in:
\begin{equation*}
    \mathrm{Exp}_{\b{X}}(\b{\Delta}) = [\b{X}, \b{\Delta}] \text{exp}\! \left( \begin{bmatrix} \b{A} & -\b{\Delta}^\top \b{\Delta} \\ \b{I}_d & \b{A} \end{bmatrix} \right) \begin{bmatrix} e^{-\b{A}} \\ \b{0} \end{bmatrix}.
\end{equation*}

\end{proof}

\subsubsection{Retraction Map in Stiefel Manifold}

As discussed in Section 11.7, a retraction map is a computationally efficient alternative to the exponential map for mapping a tangent vector back to the manifold. For the Stiefel manifold $\text{St}(n, d)$, the QR decomposition provides a standard retraction.

\begin{proposition}[QR-based retraction on the Stiefel manifold \cite{edelman1998geometry}]\label{proposition_qr_based_retraction_stiefel}
Let $\b{X} \in \text{St}(n, d)$ and $\b{\Delta} \in T_{\b{X}}\text{St}(n, d)$. The mapping $\operatorname{Ret}^{QR}_{\b{X}}: T_{\b{X}}\text{St}(n, d) \to \text{St}(n, d)$ defined by:
\begin{equation}\label{eq:stiefel_retraction_qr}
\boxed{
\operatorname{Ret}^{QR}_{\b{X}}(\b{\Delta}) = \operatorname{qf}(\b{X} + \b{\Delta}),
} 
\end{equation}
where $\operatorname{qf}(\cdot)$ denotes the $\b{Q}$ factor of the QR decomposition (see Lemma \ref{lemma_qr_decomp}), is a valid retraction map.
\end{proposition}
\begin{proof}
To be a valid retraction, $\mathrm{Ret}^{QR}_{\b{X}}$ must satisfy the properties in Definition \ref{definition_retraction}:
\begin{enumerate}
    \item \textbf{Consistency:} $\mathrm{Ret}^{QR}_{\b{X}}(\b{0}) = \b{X}$.
    \item \textbf{First-order Agreement:} $\left. \frac{d}{dt} \mathrm{Ret}^{QR}_{\b{X}}(t\b{\Delta}) \right|_{t=0} = \b{\Delta}$.
\end{enumerate}
We prove these in the following. 

\textit{1. Consistency:} For $\b{\Delta} = \b{0}$, $\mathrm{Ret}^{QR}_{\b{X}}(\b{0}) = \text{qf}(\b{X})$. Since $\b{X} \in \text{St}(n, d)$, its columns are already orthonormal ($\b{X}^\top \b{X} = \b{I}_d$ per Eq. (\ref{equation_Sitefel_manifold_definition})). Thus, its QR decomposition is $\b{X} = \b{X}\b{I}_d$, where $\b{I}_d$ is upper triangular. Hence, $\text{qf}(\b{X}) = \b{X}$.

\textit{2. First-order Agreement:} Let $\b{\b{X}}(t) = \mathrm{Ret}^{QR}_{\b{X}}(t\b{\Delta}) = \text{qf}(\b{X} + t\b{\Delta})$. By the definition of QR decomposition:
\begin{equation}
    \b{X} + t\b{\Delta} = \b{\b{X}}(t) \b{R}(t), \label{equation_retraction_deriv_base}
\end{equation}
where $\b{R}(t)$ is upper triangular. At $t=0$, Eq. (\ref{equation_retraction_deriv_base}) becomes:
\begin{align}\label{equation_R_zero_I}
\b{X} = \b{\b{X}}(0)\b{R}(0) \overset{(a)}{\implies} \b{X} = \b{X}\b{R}(0) \implies \b{R}(0) = \b{I}_d,
\end{align}
where $(a)$ is because $\b{X}(0) = \b{X}$.

Differentiating Eq. (\ref{equation_retraction_deriv_base}) with respect to $t$ at $t=0$:
\begin{align}
&\b{\Delta} = \dot{\b{\b{X}}}(0)\b{R}(0) + \b{\b{X}}(0)\dot{\b{R}}(0) \overset{(\ref{equation_R_zero_I})}{=} \dot{\b{\b{X}}}(0) + \b{X}\dot{\b{R}}(0) \nonumber\\
&\implies \dot{\b{\b{X}}}(0) = \b{\Delta} - \b{X}\dot{\b{R}}(0). \label{equation_Xdot_Delta_X_Rdot0}
\end{align}
For $\dot{\b{\b{X}}}(0) \in T_{\b{X}}\text{St}(n, d)$, the condition in Eq. (\ref{equation_stiefel_tangent_space}) must hold: $\b{X}^\top \dot{\b{X}}(0) + \dot{\b{X}}(0)^\top \b{X} = \b{0}$. Substituting $\dot{\b{X}}(0)$ gives:
\begin{align*}
&\b{X}^\top (\b{\Delta} - \b{X}\dot{\b{R}}(0)) + (\b{\Delta} - \b{X}\dot{\b{R}}(0))^\top \b{X} = \b{0} \\
&\implies (\b{X}^\top \b{\Delta} + \b{\Delta}^\top \b{X}) - (\dot{\b{R}}(0) + \dot{\b{R}}(0)^\top) = \b{0}.
\end{align*}
Since $\b{\Delta} \in T_{\b{X}}\text{St}(n, d)$, the first term is zero by Eq. (\ref{equation_stiefel_tangent_space}). This implies $\dot{\b{R}}(0) + \dot{\b{R}}(0)^\top = \b{0}$. Because $\b{R}(t)$ is upper triangular, its derivative $\dot{\b{R}}(0)$ is also upper triangular. The only upper triangular matrix that is also skew-symmetric is the zero matrix. Thus, $\dot{\b{R}}(0) = \b{0}$, which yields:
\begin{align*}
\dot{\b{\b{X}}}(0) \overset{(\ref{equation_Xdot_Delta_X_Rdot0})}{=} \b{\Delta} - \b{X}\dot{\b{R}}(0) \implies \dot{\b{\b{X}}}(0) = \b{\Delta}
\end{align*}
As $\b{X}(t) = \mathrm{Ret}_{\b{X}}(t\b{\Delta})$, we conclude $\left. \frac{d}{dt} \mathrm{Ret}_{\b{X}}(t\b{\Delta}) \right|_{t=0} = \dot{\b{X}}(0) = \b{\Delta}$. 
\end{proof}

There is also a second version of retraction, namely polar-style retraction, on the Stiefel manifold. It is introduced in the following. 

\begin{proposition}[Polar-style retraction on the Stiefel manifold \cite{edelman1998geometry}]
Let \(\b{X} \in St(n,d)\) and let \(\b{\Delta} \in T_{\b{X}}St(n,d)\). Then, the mapping:
\begin{equation}\label{equation_retraction_Stiefel_manifold_polar}
\boxed{
\operatorname{Ret}^{polar}_{\b{X}}(\b{\Delta})
=
(\b{X}+\b{\Delta})
\Big((\b{X}+\b{\Delta})^\top(\b{X}+\b{\Delta})\Big)^{-1/2},
}
\end{equation}
is a valid retraction map on the Stiefel manifold.

Equivalently, since \(\b{X}\in St(n,d)\) and \(\b{\Delta}\in T_{\b{X}}St(n,d)\), we have:
\[
(\b{X}+\b{\Delta})^\top(\b{X}+\b{\Delta})
=
\b{I}_d+\b{\Delta}^\top\b{\Delta}.
\]
Hence, the polar-style retraction in Eq. (\ref{equation_retraction_Stiefel_manifold_polar}) can also be written as:
\begin{align}\label{equation_retraction_Stiefel_manifold_polar_2}
\boxed{
\operatorname{Ret}^{polar}_{\b{X}}(\b{\Delta})
=
(\b{X}+\b{\Delta})
(\b{I}_d+\b{\Delta}^\top\b{\Delta})^{-1/2}.
}
\end{align}
\end{proposition}

\begin{proof}
According to Definition \ref{definition_retraction}, to prove that Eq. \eqref{equation_retraction_Stiefel_manifold_polar} is a valid retraction, we need to verify the following two properties:

\begin{enumerate}
\item Identity:
\[
\operatorname{Ret}^{polar}_{\b{X}}(\b{0})=\b{X}.
\]

\item Local rigidity:
\[
\left.\frac{d}{dt}\operatorname{Ret}^{polar}_{\b{X}}(t\b{\Delta})\right|_{t=0}
=
\b{\Delta}.
\]
\end{enumerate}

\textbf{Step 1: Show that the map lands on the Stiefel manifold.}

We define:
\begin{align}\label{equation_Y_X_plus_Delta_X_DeltaT_X_Delta_half_in_proof}
\b{Y}
:=
(\b{X}+\b{\Delta})
\Big((\b{X}+\b{\Delta})^\top(\b{X}+\b{\Delta})\Big)^{-1/2}.
\end{align}
We show that \(\b{Y}\in St(n,d)\), i.e., according to Eq. (\ref{equation_Sitefel_manifold_definition}), we should show that:
\[
\b{Y}^\top \b{Y}=\b{I}_d.
\]

Let:
\[
\b{A}:=(\b{X}+\b{\Delta})^\top(\b{X}+\b{\Delta}).
\]
Then Eq. (\ref{equation_Y_X_plus_Delta_X_DeltaT_X_Delta_half_in_proof}) becomes:
\[
\b{Y}=(\b{X}+\b{\Delta})\b{A}^{-1/2}.
\]
Since \(\b{A}\) is symmetric positive definite, \(\b{A}^{-1/2}\) exists and is symmetric. Therefore:
\begin{align*}
\b{Y}^\top \b{Y}
&=
\big((\b{X}+\b{\Delta})\b{A}^{-1/2}\big)^\top
\big((\b{X}+\b{\Delta})\b{A}^{-1/2}\big)
\nonumber\\
&=
(\b{A}^{-1/2})^\top(\b{X}+\b{\Delta})^\top(\b{X}+\b{\Delta})\b{A}^{-1/2}
\nonumber\\
&=
\b{A}^{-1/2}\,\b{A}\,\b{A}^{-1/2}
\nonumber\\
&=
\b{I}_d.
\end{align*}
Hence \(\b{Y}\in St(n,d)\).

Comparing Eqs. (\ref{equation_retraction_Stiefel_manifold_polar}) and (\ref{equation_Y_X_plus_Delta_X_DeltaT_X_Delta_half_in_proof}) shows that $\b{Y} = \operatorname{Ret}^{polar}_{\b{X}}(\b{\Delta})$. Therefore:
\[
\operatorname{Ret}^{polar}_{\b{X}}(\b{\Delta}) \in St(n,d).
\]

\textbf{Step 2: Identity property.}

Substituting \(\b{\Delta}=\b{0}\) into Eq. \eqref{equation_retraction_Stiefel_manifold_polar} gives:
\begin{align}
\operatorname{Ret}^{polar}_{\b{X}}(\b{0})
&=
\b{X}\,(\b{X}^\top \b{X})^{-1/2}.
\end{align}
Since \(\b{X}\in St(n,d)\), by Eq. (\ref{equation_Sitefel_manifold_definition}), we have:
\[
\b{X}^\top \b{X}=\b{I}_d.
\]
Hence:
\[
(\b{X}^\top \b{X})^{-1/2}
=
\b{I}_d^{-1/2}
=
\b{I}_d.
\]
Therefore:
\[
\operatorname{Ret}^{polar}_{\b{X}}(\b{0})=\b{X}.
\]
So the identity property holds.

\textbf{Step 3: Local rigidity property.}

Consider the curve:
\begin{align}\label{equation_Yt_X_plus_Delta_X_DeltaT_X_Delta_half_in_proof}
\b{Y}(t)
:=
(\b{X}+t\b{\Delta})
\Big((\b{X}+t\b{\Delta})^\top(\b{X}+t\b{\Delta})\Big)^{-1/2}.
\end{align}
Comparing Eqs. (\ref{equation_retraction_Stiefel_manifold_polar}) and (\ref{equation_Yt_X_plus_Delta_X_DeltaT_X_Delta_half_in_proof}) shows that:
\[
\operatorname{Ret}^{polar}_{\b{X}}(t\b{\Delta})=\b{Y}(t).
\]
We expand \(\b{Y}(t)\) around \(t=0\).

First,
\begin{align}
(\b{X}+&t\b{\Delta})^\top(\b{X}+t\b{\Delta}) \nonumber\\
&=
\b{X}^\top \b{X}
+
t\b{X}^\top \b{\Delta}
+
t\b{\Delta}^\top \b{X}
+
t^2 \b{\Delta}^\top \b{\Delta}
\nonumber\\
&=
\b{I}_d
+
t(\b{X}^\top \b{\Delta}+\b{\Delta}^\top \b{X})
+
t^2 \b{\Delta}^\top \b{\Delta}.
\label{equation_gram_stiefel_polar_retraction}
\end{align}
Because \(\b{\Delta}\in T_{\b{X}}St(n,d)\), according to Eq. (\ref{equation_stiefel_tangent_space}), we have:
\[
\b{X}^\top \b{\Delta}+\b{\Delta}^\top \b{X}=\b{0}.
\]
Therefore, Eq. \eqref{equation_gram_stiefel_polar_retraction} becomes:
\[
(\b{X}+t\b{\Delta})^\top(\b{X}+t\b{\Delta})
=
\b{I}_d+t^2\b{\Delta}^\top \b{\Delta}.
\]
Hence:
\begin{equation}\label{equation_inverse_sqrt_stiefel_polar_retraction}
\Big((\b{X}+t\b{\Delta})^\top(\b{X}+t\b{\Delta})\Big)^{-1/2}
=
(\b{I}_d+t^2\b{\Delta}^\top \b{\Delta})^{-1/2}.
\end{equation}
Now, using the Taylor expansion around \(t=0\), we have:
\[
(\b{I}_d+t^2\b{\Delta}^\top \b{\Delta})^{-1/2}
=
\b{I}_d+\mathcal{O}(t^2),
\]
where $\mathcal{O}(.)$ denotes the big-O complexity notation. 
Therefore, Eq. (\ref{equation_inverse_sqrt_stiefel_polar_retraction}) becomes:
\begin{align*}
\Big((\b{X}+t\b{\Delta})^\top(\b{X}+t\b{\Delta})\Big)^{-1/2} = \b{I}_d+\mathcal{O}(t^2).
\end{align*}

Substituting this into \(\b{Y}(t)\) in Eq. (\ref{equation_Yt_X_plus_Delta_X_DeltaT_X_Delta_half_in_proof}) gives:
\begin{align}
&\b{Y}(t)
=
(\b{X}+t\b{\Delta})(\b{I}_d+\mathcal{O}(t^2))
\nonumber\\
&~~~~~~~~ =
\b{X}+t\b{\Delta}+\mathcal{O}(t^2).
\label{equation_Yt_stiefel_polar_retraction_first_order} \\
&\implies \dot{\b{Y}}(t) = \b{\Delta}.
\end{align}
Therefore:
\[
\b{Y}(0)=\b{X},
\qquad
\dot{\b{Y}}(0)=\b{\Delta}.
\]
Thus:
\[
\left.\frac{d}{dt}\operatorname{Ret}^{polar}_{\b{X}}(t\b{\Delta})\right|_{t=0}
=
\left.\frac{d}{dt}\b{Y}(t)\right|_{t=0}
=
\b{\Delta}.
\]
So the local rigidity property holds.

Since both identity and local rigidity are satisfied, Eq. \eqref{equation_retraction_Stiefel_manifold_polar} is a valid retraction map on the Stiefel manifold.

Since \(\b{X}\in St(n,d)\) and \(\b{\Delta}\in T_{\b{X}}St(n,d)\), we have:
\begin{align*}
\b{X}^\top \b{X} \overset{(\ref{equation_Sitefel_manifold_definition})}{=} \b{I}_d, \quad \b{X}^\top \b{\Delta} + \b{\Delta}^\top \b{X} \overset{(\ref{equation_stiefel_tangent_space})}{=} \b{0}.
\end{align*}
Thus:
\begin{align*}
(\b{X}+\b{\Delta})^\top&(\b{X}+\b{\Delta}) \\
&= \b{X}^\top \b{X} + \b{X}^\top \b{\Delta} + \b{\Delta}^\top \b{X} + \b{\Delta}^\top \b{\Delta} \\
&= \b{I}_d+\b{\Delta}^\top\b{\Delta}.
\end{align*}
Hence, the polar-style retraction in Eq. (\ref{equation_retraction_Stiefel_manifold_polar}) can also be written as:
\begin{align*}
\operatorname{Ret}^{polar}_{\b{X}}(\b{\Delta})
=
(\b{X}+\b{\Delta})
(\b{I}_d+\b{\Delta}^\top\b{\Delta})^{-1/2}.
\end{align*}
\end{proof}

\begin{remark}[Relation of polar-style and QR-based retractions in Stiefel manifold]
The polar-style retraction in Eq. \eqref{equation_retraction_Stiefel_manifold_polar} is an alternative to the QR-based retraction in Eq. (\ref{eq:stiefel_retraction_qr}), on the Stiefel manifold. Both maps return a matrix with orthonormal columns and both agree with the exponential map up to first order. The QR-based retraction uses the \(Q\)-factor of the QR decomposition, while the polar-style retraction uses the symmetric inverse square root of the Gram matrix.

In other words, the QR-based retraction is:
\[
\operatorname{Ret}^{QR}_{\b{X}}(\b{\Delta})
=
\operatorname{qf}(\b{X}+\b{\Delta}),
\]
whereas the polar-style retraction is:
\[
\operatorname{Ret}^{polar}_{\b{X}}(\b{\Delta})
=
(\b{X}+\b{\Delta})
\Big((\b{X}+\b{\Delta})^\top(\b{X}+\b{\Delta})\Big)^{-1/2},
\]
or:
\begin{align*}
\operatorname{Ret}^{polar}_{\b{X}}(\b{\Delta})
=
(\b{X}+\b{\Delta})
(\b{I}_d+\b{\Delta}^\top\b{\Delta})^{-1/2}.
\end{align*}
\end{remark}

\subsubsection{Vector Transport in Stiefel Manifold}

In Riemannian optimization, many algorithms—such as Conjugate Gradient or Quasi-Newton methods—require comparing or combining tangent vectors located at different points on the manifold. While parallel transport along a geodesic is the canonical method for this, it is often computationally prohibitive. Vector transport is a generalized, computationally efficient alternative that relaxes the requirements of parallel transport while remaining compatible with a chosen retraction.

\begin{proposition}[Vector transport in Stiefel manifold] \label{proposition_vector_transport_stiefel}
Let $\b{X} \in \text{St}(n, d)$ and $\b{\Delta}_1 \in T_{\b{X}} \text{St}(n, d)$. Given a direction $\b{\Delta}_2 \in T_{\b{X}} \text{St}(n, d)$ and a retraction $\b{Y} = \mathrm{Ret}_{\b{X}}(\b{\Delta}_2)$, the vector transport $\mathcal{T}_{\b{\Delta}_2}(\b{\Delta}_1) \in T_{\b{Y}} \text{St}(n, d)$ via orthogonal projection is defined as:
\begin{equation}\label{equation_retraction_Stiefel}
\boxed{
\mathcal{T}_{\b{\Delta}_2}(\b{\Delta}_1) = \b{\Delta}_1 - \b{Y} \text{sym}(\b{Y}^\top \b{\Delta}_1).
}
\end{equation}
\end{proposition}
\begin{proof}
The projection-based vector transport is defined by treating the tangent vector $\b{\Delta}_1$ as an element of the ambient Euclidean space $\mathbb{R}^{n \times d}$ and projecting it onto the tangent space at the new point $\b{Y}$.
Using the projection formula in Eq. (\ref{equation_stiefel_proj_formula}), evaluated at $\b{Y}$, we have:
\begin{equation}\label{equation_T_Delta2_Delta1_in_proof}
\mathcal{T}_{\b{\Delta}_2}(\b{\Delta}_1) = \Pi_{\b{Y}}^{\mathrm{St}}(\b{\Delta}_1) \overset{(\ref{equation_stiefel_proj_formula})}{=} \b{\Delta}_1 - \b{Y} \text{sym}(\b{Y}^\top \b{\Delta}_1).
\end{equation}

To verify that $\mathcal{T}_{\b{\Delta}_2}(\b{\Delta}_1) \in T_{\b{Y}} \text{St}(n, d)$, it must satisfy Eq. (\ref{equation_stiefel_tangent_space}), i.e., the symmetry condition. So, we check the symmetry condition:
\begin{align*}
&\b{Y}^\top \mathcal{T}_{\b{\Delta}_2}(\b{\Delta}_1) + (\mathcal{T}_{\b{\Delta}_2}(\b{\Delta}_1))^\top \b{Y}  \\
&\overset{(\ref{equation_T_Delta2_Delta1_in_proof})}{=} 
\b{Y}^\top (\b{\Delta}_1 - \b{Y} \text{sym}(\b{Y}^\top \b{\Delta}_1)) \\
&\quad\quad\quad\quad\quad\quad+ (\b{\Delta}_1 - \b{Y} \text{sym}(\b{Y}^\top \b{\Delta}_1))^\top \b{Y} \\
&=\b{Y}^\top \b{\Delta}_1 - \b{Y}^\top \b{Y} \text{sym}(\b{Y}^\top \b{\Delta}_1) \\
&\quad\quad\quad\quad\quad\quad+ \b{\Delta}_1^\top \b{Y} - \text{sym}(\b{Y}^\top \b{\Delta}_1)^\top \b{Y}^\top \b{Y} \\
&\overset{(a)}{=} \b{Y}^\top \b{\Delta}_1 - \text{sym}(\b{Y}^\top \b{\Delta}_1) \\
&\quad\quad\quad\quad\quad\quad+ \b{\Delta}_1^\top \b{Y} - \text{sym}(\b{Y}^\top \b{\Delta}_1) \\
&= (\b{Y}^\top \b{\Delta}_1 + \b{\Delta}_1^\top \b{Y}) - 2\text{sym}(\b{Y}^\top \b{\Delta}_1) \\
&\overset{(\ref{equation_sym_skew_expressions})}{=} (\b{Y}^\top \b{\Delta}_1 + \b{\Delta}_1^\top \b{Y}) - (\b{Y}^\top \b{\Delta}_1 + \b{\Delta}_1^\top \b{Y}) = \b{0},
\end{align*}
where $(a)$ is because $\b{Y} = \mathrm{Ret}_{\b{X}}(\b{\Delta}_2) \in \text{St}(n, d)$, so $\b{Y}^\top \b{Y} = \b{I}$ according to Eq. (\ref{equation_Sitefel_manifold_definition}).

This shows that satisfy $\mathcal{T}_{\b{\Delta}_2}(\b{\Delta}_1)$ satisfies Eq. (\ref{equation_stiefel_tangent_space}), i.e., the symmetry condition. Thus, it is a valid tangent vector at $\b{Y}$.
\end{proof}

\subsubsection{Vector Transport by Differential QR-Based Retraction in Stiefel Manifold}

As discussed in Section \ref{section_differentiated_retraction}, a natural way to construct a vector transport is by differentiating a chosen retraction.
For the Stiefel manifold, we can use the QR-based retraction
introduced in Eq. (\ref{eq:stiefel_retraction_qr}):
\begin{equation*}
\operatorname{Ret}^{QR}_{\b{X}}(\b{\Delta}) = \operatorname{qf}(\b{X} + \b{\Delta}).
\end{equation*}
We now derive the associated differentiated-retraction vector transport.

\begin{lemma}[Skew-upper triangular decomposition]\label{lemma_skew_upper_triangular_decomposition}
Let $\b{B} \in \mathbb{R}^{d\times d}$. There exist unique matrices
$\b{\Omega}, \b{U} \in \mathbb{R}^{d\times d}$ such that:
\begin{align}
\b{B} = \b{\Omega} + \b{U},
\end{align}
where $\b{\Omega}^\top = -\b{\Omega}$ is skew-symmetric and $\b{U}$ is upper triangular.

If we denote by $\operatorname{low}(\b{B})$ the strictly lower triangular part of $\b{B}$, then:
\begin{align}
\boxed{
\b{\Omega} = \operatorname{low}(\b{B}) - \operatorname{low}(\b{B})^\top,
\qquad
\b{U} = \b{B} - \b{\Omega}.
}
\end{align}
\end{lemma}

\begin{proof}
We first show existence. We define:
\[
\b{\Omega} := \operatorname{low}(\b{B}) - \operatorname{low}(\b{B})^\top.
\]
Then, clearly $\b{\Omega}^\top = -\b{\Omega}$, so $\b{\Omega}$ is skew-symmetric.
Now, we define:
\[
\b{U} := \b{B} - \b{\Omega}.
\]
We show that $\b{U}$ is upper triangular.

By construction, $\operatorname{low}(\b{\Omega}) = \operatorname{low}(\b{B})$.
Hence:
\[
\operatorname{low}(\b{U})
=
\operatorname{low}(\b{B} - \b{\Omega})
=
\operatorname{low}(\b{B}) - \operatorname{low}(\b{\Omega})
=
\b{0}.
\]
Therefore, $\b{U}$ has no entries below the diagonal, so it is upper triangular.

Now, we show uniqueness. Suppose:
\[
\b{B} = \b{\Omega}_1 + \b{U}_1 = \b{\Omega}_2 + \b{U}_2,
\]
where $\b{\Omega}_1,\b{\Omega}_2$ are skew-symmetric and $\b{U}_1,\b{U}_2$ are upper triangular.
Then:
\[
\b{\Omega}_1 - \b{\Omega}_2 = \b{U}_2 - \b{U}_1.
\]
The left-hand side is skew-symmetric, while the right-hand side is upper triangular.
Hence the matrix $\b{\Omega}_1 - \b{\Omega}_2$ is both skew-symmetric and upper triangular.
The only matrix with both properties is the zero matrix. Therefore:
\[
\b{\Omega}_1 - \b{\Omega}_2 = \b{0},
\qquad
\b{U}_2 - \b{U}_1 = \b{0}.
\]
So $\b{\Omega}_1 = \b{\Omega}_2$ and $\b{U}_1 = \b{U}_2$.
This proves uniqueness.
\end{proof}

\begin{proposition}[Vector transport by differential QR-based retraction on the Stiefel manifold]\label{proposition_differential_qr_retraction_stiefel}
Let $\b{X} \in St(n,d)$ and let
$\b{\Delta}_1, \b{\Delta}_2 \in T_{\b{X}}St(n,d)$.
We define the QR-based retraction by Eq. (\ref{eq:stiefel_retraction_qr}):
\[
\operatorname{Ret}^{QR}_{\b{X}}(\b{\Delta})
=
\operatorname{qf}(\b{X}+\b{\Delta}).
\]
Let:
\[
\b{Y}
:=
\operatorname{Ret}^{QR}_{\b{X}}(\b{\Delta}_2)
=
\operatorname{qf}(\b{X}+\b{\Delta}_2),
\]
and let the QR decomposition of $\b{X}+\b{\Delta}_2$ be:
\[
\b{X}+\b{\Delta}_2 = \b{Y}\b{R},
\]
where $\b{R} \in \mathbb{R}^{d\times d}$ is upper triangular and invertible.

The differentiated retraction:
\[
\mathcal{T}^{\operatorname{QR}}_{\b{\Delta}_2}(\b{\Delta}_1) = D\operatorname{Ret}^{QR}_{[\b{X}]}(\b{\Delta}_2)[\b{\Delta}_1]
\in
T_{\b{Y}}St(n,d),
\]
is represented by:
\begin{equation}\label{equation_differentiated_qr_retraction_vector_transport_stiefel}
\boxed{
\mathcal{T}^{\operatorname{QR}}_{\b{\Delta}_2}(\b{\Delta}_1)
=
\b{Y}\b{\Omega}
+
(\b{I}_n - \b{Y}\b{Y}^\top)\b{\Delta}_1\b{R}^{-1}.}
\end{equation}
where:
\begin{align}
& \b{B} := \b{Y}^\top \b{\Delta}_1 \b{R}^{-1}, \\
& \b{\Omega}
:=
\operatorname{low}(\b{B}) - \operatorname{low}(\b{B})^\top.
\end{align}

The Eq.
\eqref{equation_differentiated_qr_retraction_vector_transport_stiefel} defines the differentiated-retraction vector transport associated with the QR-based retraction on the Stiefel manifold.

\end{proposition}

\begin{proof}
According to Definition \ref{definition_differentiated_retraction}, the differentiated retraction gives the vector transport:
\[
\mathcal{T}^{\operatorname{QR}}_{\b{\Delta}_2}(\b{\Delta}_1)
=
D\operatorname{Ret}^{QR}_{\b{X}}(\b{\Delta}_2)[\b{\Delta}_1].
\]
So, it remains to compute the differential of the QR-based retraction explicitly.

\vspace{0.2cm}
\noindent
\textbf{Step 1: Define a curve through the retraction argument.}

Consider the matrix-valued curve:
\begin{align}\label{equation_A_X_Delta2_tDelta1_in_proof}
\b{A}(t) := \b{X}+\b{\Delta}_2+t\b{\Delta}_1.
\end{align}
For $t$ in a neighborhood of zero, the matrix $\b{A}(t)$ has full column rank, so its QR factorization is well-defined:
\[
\b{A}(t) = \b{Y}(t)\b{R}(t),
\]
where:
\[
\b{Y}(t) \in St(n,d),
\qquad
\b{R}(t) \in \mathbb{R}^{d\times d},
\]
is upper triangular and invertible.
By considering Eq. (\ref{equation_A_X_Delta2_tDelta1_in_proof}) and by definition of the QR-based retraction, i.e., Eq. (\ref{eq:stiefel_retraction_qr}), we have:
\[
\b{Y}(t) = \operatorname{qf}(\b{A}(t))
= \operatorname{Ret}^{QR}_{\b{X}}(\b{\Delta}_2+t\b{\Delta}_1).
\]
At $t=0$, we have:
\begin{align}
&\b{Y}(0)=\b{Y},
\quad
\b{R}(0)=\b{R}, \label{equation_Y0_Y_R0_R_in_proof} \\
&\b{A}(0)=\b{X}+\b{\Delta}_2=\b{Y}\b{R}. \nonumber
\end{align}

\vspace{0.2cm}
\noindent
\textbf{Step 2: Differentiate the QR factorization.}

Differentiating $\b{A}(t)=\b{Y}(t)\b{R}(t)$ with respect to $t$ gives:
\[
\b{\Delta}_1
=
\dot{\b{Y}}(t)\b{R}(t)+\b{Y}(t)\dot{\b{R}}(t).
\]
Evaluating at $t=0$:
\begin{align}
&\b{\Delta}_1
=
\dot{\b{Y}}(0)\b{R}
+
\b{Y}\dot{\b{R}}(0) \label{equation_Delta1_YdotR_Y_Rdot_in_proof}\\
&\implies \dot{\b{Y}}(0)
=
\big(\b{\Delta}_1-\b{Y}\dot{\b{R}}(0)\big)\b{R}^{-1}. \label{equation_Ydot0_Delta1_Y_Rdot0_Rinverse_in_proof}
\end{align}

Since we have:
\[
\b{Y}(t)=\operatorname{Ret}_{\b{X}}(\b{\Delta}_2+t\b{\Delta}_1),
\]
the curve \(\b{Y}(t)\) is the composition of the retraction map with the line:
\[
t \mapsto \b{\Delta}_2+t\b{\Delta}_1,
\]
in the vector space \(T_{\b{X}}St(n,d)\). Therefore, by the definition of the differential, we have:
\begin{align}\label{equation_DRetXDelta2Delta1_Ydot0_in_proof}
D\operatorname{Ret}^{QR}_{\b{X}}(\b{\Delta}_2)[\b{\Delta}_1]
=
\left.\frac{d}{dt}\operatorname{Ret}^{QR}_{\b{X}}(\b{\Delta}_2+t\b{\Delta}_1)\right|_{t=0}
=
\dot{\b{Y}}(0).
\end{align}

As we have:
\[
D\operatorname{Ret}^{QR}_{\b{X}}(\b{\Delta}_2)[\b{\Delta}_1]
=
\dot{\b{Y}}(0),
\]
we need to find $\dot{\b{Y}}(0)$. According to Eq. (\ref{equation_Ydot0_Delta1_Y_Rdot0_Rinverse_in_proof}), the problem reduces to determining $\dot{\b{R}}(0)$.

\vspace{0.2cm}
\noindent
\textbf{Step 3: Use the Stiefel constraint to identify the skew-symmetric part.}

Because $\b{Y}(t)\in St(n,d)$ for all $t$, we have according to Eq. (\ref{equation_Sitefel_manifold_definition}):
\[
\b{Y}(t)^\top \b{Y}(t)=\b{I}_d.
\]
Evaluation at $t=0$:
\begin{align}\label{equation_Y0TY_I_YT_Y_I_in_proof}
\b{Y}(0)^\top \b{Y}(0)=\b{I}_d \overset{(\ref{equation_Y0_Y_R0_R_in_proof})}{\implies} \b{Y}^\top \b{Y}=\b{I}_d.
\end{align}

Differentiating with respect to $t$ yields:
\[
\dot{\b{Y}}(t)^\top \b{Y}(t)+\b{Y}(t)^\top \dot{\b{Y}}(t)=\b{0}.
\]
Evaluating at $t=0$:
\[
\dot{\b{Y}}(0)^\top \b{Y}+\b{Y}^\top \dot{\b{Y}}(0)=\b{0}.
\]
Thus:
\[
\b{\Omega} := \b{Y}^\top \dot{\b{Y}}(0),
\]
is skew-symmetric:
\begin{align}\label{equation_Omega_OmegaT_zero_in_proof}
\b{\Omega}^\top = -\b{\Omega} \implies \b{\Omega} + \b{\Omega}^\top = \b{0}.
\end{align}

Now, we left-multiply Eq. (\ref{equation_Delta1_YdotR_Y_Rdot_in_proof}) by $\b{Y}^\top$:
\begin{align*}
&\b{Y}^\top \b{\Delta}_1
=
\b{Y}^\top \dot{\b{Y}}(0)\b{R}
+
\b{Y}^\top \b{Y}\dot{\b{R}}(0) \\
&\overset{(\ref{equation_Y0TY_I_YT_Y_I_in_proof})}{\implies}
\b{Y}^\top \b{\Delta}_1
=
\b{\Omega}\b{R}
+
\dot{\b{R}}(0).
\end{align*}
Multiplying on the right by $\b{R}^{-1}$ gives:
\begin{align}
&\b{B}
:=
\b{Y}^\top \b{\Delta}_1 \b{R}^{-1}
=
\b{\Omega}
+
\dot{\b{R}}(0)\b{R}^{-1} \label{equation_B_YT_Delta1_Rinv_Omega_Rdot0_Rinv_in_proof} \\
&\implies \b{B} - \b{\Omega} = \dot{\b{R}}(0)\b{R}^{-1}. \label{equation_B_Omega_Rdot0_Rinverse_in_proof}
\end{align}
Because $\b{R}(t)$ is upper triangular for all $t$, its derivative $\dot{\b{R}}(0)$ is upper triangular, and therefore
$\dot{\b{R}}(0)\b{R}^{-1}$ is also upper triangular.
So, $\b{B}$ is written as:
\[
\b{B} = \b{\Omega} + \b{U},
\]
where $\b{\Omega}$ is skew-symmetric and $\b{U}:=\dot{\b{R}}(0)\b{R}^{-1}$ is upper triangular.

By the Lemma \ref{lemma_skew_upper_triangular_decomposition}, this decomposition is unique, and therefore:
\[
\b{\Omega}
=
\operatorname{low}(\b{B}) - \operatorname{low}(\b{B})^\top.
\]

\vspace{0.2cm}
\noindent
\textbf{Step 4: Substitute back into the derivative formula.}

From Eq. (\ref{equation_Ydot0_Delta1_Y_Rdot0_Rinverse_in_proof}), we obtain:
\begin{align*}
\dot{\b{Y}}(0)
&=
\b{\Delta}_1\b{R}^{-1}
-
\b{Y}\dot{\b{R}}(0)\b{R}^{-1} \\
&\overset{(\ref{equation_B_Omega_Rdot0_Rinverse_in_proof})}{=} 
\b{\Delta}_1\b{R}^{-1}
-
\b{Y}(\b{B}-\b{\Omega})
\\
&=
\b{\Delta}_1\b{R}^{-1}
-
\b{Y}\b{B}
+
\b{Y}\b{\Omega} \\
&\overset{(\ref{equation_B_YT_Delta1_Rinv_Omega_Rdot0_Rinv_in_proof})}{=} 
\b{\Delta}_1\b{R}^{-1}
-
\b{Y}\b{Y}^\top \b{\Delta}_1 \b{R}^{-1}
+
\b{Y}\b{\Omega} 
\\
&= 
(\b{I}_n-\b{Y}\b{Y}^\top)\b{\Delta}_1\b{R}^{-1}
+
\b{Y}\b{\Omega}.
\end{align*}

Therefore, according to Eq. (\ref{equation_DRetXDelta2Delta1_Ydot0_in_proof}), we have:
\[
D\operatorname{Ret}_{\b{X}}(\b{\Delta}_2)[\b{\Delta}_1]
=
\b{Y}\b{\Omega}
+
(\b{I}_n-\b{Y}\b{Y}^\top)\b{\Delta}_1\b{R}^{-1},
\]
which is Eq. (\ref{equation_differentiated_qr_retraction_vector_transport_stiefel}).

\vspace{0.2cm}
\noindent
\textbf{Step 5: Verify that the transported vector is tangent at $\b{Y}$.}

Let:
\begin{align}\label{equation_Xi_YOmega_I_YYT_Delta1_Rinv_in_proof}
\b{\Xi}
:=
\b{Y}\b{\Omega}
+
(\b{I}_n-\b{Y}\b{Y}^\top)\b{\Delta}_1\b{R}^{-1}.
\end{align}
We check the tangent constraint in Eq. (\ref{equation_stiefel_tangent_space}), i.e., whether $\b{Y}^\top \b{\Xi} + \b{\Xi}^\top \b{Y} = \b{0}$.

First, we consider:
\[
\b{Y}^\top \b{\Xi}
\overset{(\ref{equation_Xi_YOmega_I_YYT_Delta1_Rinv_in_proof})}{=}
\b{Y}^\top \b{Y}\b{\Omega}
+
\b{Y}^\top(\b{I}_n-\b{Y}\b{Y}^\top)\b{\Delta}_1\b{R}^{-1}
=
\b{\Omega},
\]
because $\b{Y}^\top \b{Y}=\b{I}_d$ and:
\[
\b{Y}^\top(\b{I}_n-\b{Y}\b{Y}^\top)=\b{0}.
\]
Similarly, we have:
\[
\b{\Xi}^\top \b{Y}
=
\b{\Omega}^\top.
\]
Therefore:
\[
\b{Y}^\top \b{\Xi} + \b{\Xi}^\top \b{Y}
=
\b{\Omega}+\b{\Omega}^\top
\overset{(\ref{equation_Omega_OmegaT_zero_in_proof})}{=}
\b{0}.
\]
We proved that $\b{\Xi}$ satisfies the tangent constraint in Eq. (\ref{equation_stiefel_tangent_space}); hence, $\b{\Xi}\in T_{\b{Y}}St(n,d)$.
This proves that the differentiated QR-based retraction gives a valid vector transport from $T_{\b{X}}St(n,d)$ to $T_{\b{Y}}St(n,d)$.
\end{proof}

\begin{remark}[Special case of differentiated QR-based retraction at zero step in Stiefel manifold]
If $\b{\Delta}_2=\b{0}$, then $\b{Y}=\b{X}$ and $\b{R}=\b{I}_d$.
Moreover,
\[
\b{B}=\b{X}^\top \b{\Delta}_1.
\]
Because $\b{\Delta}_1\in T_{\b{X}}St(n,d)$, the matrix $\b{X}^\top \b{\Delta}_1$ is skew-symmetric according to Eq. (\ref{equation_skew_symmetric}). Hence:
\[
\operatorname{low}(\b{B})-\operatorname{low}(\b{B})^\top=\b{B}.
\]
So, the transport formula reduces to:
\[
\mathcal{T}^{\operatorname{QR}}_{\b{0}}(\b{\Delta}_1)
=
\b{X}\b{X}^\top \b{\Delta}_1
+
(\b{I}_n-\b{X}\b{X}^\top)\b{\Delta}_1
=
\b{\Delta}_1.
\]
Therefore, the differentiated QR-based vector transport agrees with the identity transport at zero displacement, as expected.
\end{remark}


\begin{remark}[Interpretation of differentiated QR-based retraction in Stiefel manifold]
In Eq. (\ref{equation_differentiated_qr_retraction_vector_transport_stiefel}), the factor \(\b{R}^{-1}\) appears because the QR-based retraction is defined through the decomposition:
\[
\b{X}+\b{\Delta}_2 = \b{Y}\b{R}.
\]
When the input is perturbed in the direction \(\b{\Delta}_1\), both the orthonormal factor \(\b{Y}\) and the upper triangular factor \(\b{R}\) change. Therefore, the derivative of the retraction depends not only on the perturbation \(\b{\Delta}_1\), but also on how this perturbation is adjusted by the change in \(\b{R}\), which leads to the factor \(\b{R}^{-1}\).

The term $(\b{I}_n-\b{Y}\b{Y}^\top)\b{\Delta}_1\b{R}^{-1}$ in Eq. (\ref{equation_differentiated_qr_retraction_vector_transport_stiefel}) is the component of the variation orthogonal to the current columns of \(\b{Y}\). 

The term $\b{Y}\b{\Omega}$ in Eq. (\ref{equation_differentiated_qr_retraction_vector_transport_stiefel}) is the component tangent to the span of \(\b{Y}\), and it remains tangent because \(\b{\Omega}\) is skew-symmetric. 

Hence, the differentiated QR-based retraction, stated in Eq. (\ref{equation_differentiated_qr_retraction_vector_transport_stiefel}), splits naturally into two tangent parts: one describing motion orthogonal to the current frame, and one describing an infinitesimal rotation within the orthonormal frame itself. Their sum gives the transported tangent vector in \(T_{\b{Y}}St(n,d)\).
\end{remark}

\subsubsection{Vector Transport by Differential Polar-Style Retraction in Stiefel Manifold}\label{section_differentiated_polar_retraction_stiefel}

As discussed in Section \ref{section_differentiated_retraction}, a natural way to construct a vector transport is by differentiating a chosen retraction.
For the Stiefel manifold, besides the QR-based retraction,
we can also use the polar-style retraction introduced in Eq. (\ref{equation_retraction_Stiefel_manifold_polar_2}):
\[
\operatorname{Ret}^{polar}_{\b{X}}(\b{\Delta})
=
(\b{X}+\b{\Delta})(\b{I}_d+\b{\Delta}^{\top}\b{\Delta})^{-1/2}.
\]
We now derive the associated differentiated-retraction vector transport.

\begin{lemma}[Derivative of inverse square root matrix factor]\label{lemma_derivative_inverse_square_root_matrix_factor}
Let:
\begin{equation}
\b{S}(t)
:=
\Big(
\b{I}_d + (\b{\Delta}_2+t\b{\Delta}_1)^\top (\b{\Delta}_2+t\b{\Delta}_1)
\Big)^{-1/2}.
\label{equation_S_t_polar_transport_grassmann}
\end{equation}
Then, $\b{S}(t)$ is smooth in $t$ in a neighborhood of $t=0$.
Moreover, if:
\begin{equation}\label{equation_A_t_polar_transport_grassmann}
\b{A}(t)
:=
\b{I}_d + (\b{\Delta}_2+t\b{\Delta}_1)^\top (\b{\Delta}_2+t\b{\Delta}_1),
\end{equation}
then:
\begin{equation}
\dot{\b{A}}(0)
=
\b{\Delta}_1^\top \b{\Delta}_2 + \b{\Delta}_2^\top \b{\Delta}_1.
\label{equation_A_dot_zero_polar_transport_grassmann}
\end{equation}
\end{lemma}

\begin{proof}
We define:
\[
\b{A}(t)
:=
\b{I}_d + (\b{\Delta}_2+t\b{\Delta}_1)^\top (\b{\Delta}_2+t\b{\Delta}_1),
\]
We expand it:
\begin{align*}
\b{A}(t)
=
&\b{I}_d
+
\b{\Delta}_2^\top \b{\Delta}_2
\\
&+
t(\b{\Delta}_1^\top \b{\Delta}_2 + \b{\Delta}_2^\top \b{\Delta}_1)
+
t^2 \b{\Delta}_1^\top \b{\Delta}_1.
\end{align*}
The derivative with respect to $t$ is:
\begin{align*}
&\dot{\b{A}}(t)
=
(\b{\Delta}_1^\top \b{\Delta}_2 + \b{\Delta}_2^\top \b{\Delta}_1)
+
2t \b{\Delta}_1^\top \b{\Delta}_1 \\
&\implies \dot{\b{A}}(0)
=
\b{\Delta}_1^\top \b{\Delta}_2 + \b{\Delta}_2^\top \b{\Delta}_1.
\end{align*}
Because $\b{A}(0)=\b{I}_d+\b{\Delta}_2^\top\b{\Delta}_2$ is symmetric
positive definite, the matrix inverse square root is well-defined and
smooth in a neighborhood of $t=0$. Hence, according to Eqs. (\ref{equation_S_t_polar_transport_grassmann}) and (\ref{equation_A_t_polar_transport_grassmann}), we have $\b{S}(t)=\b{A}(t)^{-1/2}$
is smooth in $t$.
\end{proof}

\begin{proposition}[Differentiated polar-style retraction vector transport on Stiefel manifold]
Let \(\b{X}\in St(n,d)\). Let:
\[
\b{\Delta}_1,\b{\Delta}_2 \in T_{\b{X}}St(n,d).
\]
Consider the polar-style retraction in Eq. (\ref{equation_retraction_Stiefel_manifold_polar_2}):
\[
\operatorname{Ret}^{polar}_{\b{X}}(\b{\Delta})
=
(\b{X}+\b{\Delta})(\b{I}_d+\b{\Delta}^{\top}\b{\Delta})^{-1/2}.
\]
We define:
\[
\b{Y}
:=
(\b{X}+\b{\Delta}_2)(\b{I}_d+\b{\Delta}_2^{\top}\b{\Delta}_2)^{-1/2}.
\]
The differentiated retraction:
\[
\mathcal{T}^{\operatorname{polar}}_{\b{\Delta}_2}(\b{\Delta}_1)
:=
D\operatorname{Ret}^{polar}_{\b{X}}(\b{\Delta}_2)[\b{\Delta}_1]
\in T_{\b{Y}}St(n,d),
\]
is given by:
\begin{equation}\label{equation_differentiated_polar_retraction_stiefel}
\boxed{
\begin{aligned}
\mathcal{T}^{\operatorname{polar}}_{\b{\Delta}_2}(\b{\Delta}_1)
=
\b{\Delta}_1(\b{I}_d+&\b{\Delta}_2^{\top}\b{\Delta}_2)^{-1/2}
\\
&+
(\b{X}+\b{\Delta}_2)\dot{\b{S}}(0),
\end{aligned}
}
\end{equation}
where:
\begin{align}
& \b{S}(t)
:=
\Big(\b{I}_d+(\b{\Delta}_2+t\b{\Delta}_1)^{\top}(\b{\Delta}_2+t\b{\Delta}_1)\Big)^{-1/2}, \label{equation_S_I_Delta2_tDelta1_Delta2_tDelta1_in_differentiated_polar_retraction_stiefel} \\
& \dot{\b{S}}(0)
=
\left.\frac{d}{dt}\b{S}(t)\right|_{t=0}.
\end{align}
Equation (\ref{equation_differentiated_polar_retraction_stiefel}) defines the vector transport by differentiated polar-style retraction on the Stiefel manifold.
\end{proposition}

\begin{proof}
We compute the differential of the polar-style retraction at
\(\b{\Delta}_2\) in the direction \(\b{\Delta}_1\).

\textbf{Step 1: Perturb the tangent input.}

We define the input curve:
\[
\b{\Delta}(t) := \b{\Delta}_2+t\b{\Delta}_1.
\]
Applying the polar-style retraction, i.e., Eq. (\ref{equation_retraction_Stiefel_manifold_polar_2}), to this curve gives:
\[
\operatorname{Ret}^{polar}_{\b{X}}(\b{\Delta}(t))
=
(\b{X}+\b{\Delta}(t))
\big(\b{I}_d+\b{\Delta}(t)^{\top}\b{\Delta}(t)\big)^{-1/2}.
\]
For convenience, we define the curve:
\[
\b{Y}(t)
:=
(\b{X}+\b{\Delta}_2+t\b{\Delta}_1)\b{S}(t),
\]
where \(\b{S}(t)\) is defined in Eq. (\ref{equation_S_I_Delta2_tDelta1_Delta2_tDelta1_in_differentiated_polar_retraction_stiefel}).
According to the definition of polar-style retraction in Eq. (\ref{equation_retraction_Stiefel_manifold_polar_2}),
the curve on the Stiefel manifold is:
\[
\b{Y}(t)=\operatorname{Ret}^{polar}_{\b{X}}(\b{\Delta}_2+t\b{\Delta}_1).
\]
Therefore, by the definition of differential, we have:
\begin{align}\label{equation_DRet_Delta2_Delta1_Ydot0_in_proof}
D\operatorname{Ret}^{polar}_{\b{X}}(\b{\Delta}_2)[\b{\Delta}_1]
=
\left.\frac{d}{dt}\b{Y}(t)\right|_{t=0}
=
\dot{\b{Y}}(0).
\end{align}

\textbf{Step 2: Differentiate the curve.}

Differentiating \(\b{Y}(t)\) with respect to \(t\) gives:
\[
\dot{\b{Y}}(t)
=
\b{\Delta}_1\b{S}(t)
+
(\b{X}+\b{\Delta}_2+t\b{\Delta}_1)\dot{\b{S}}(t).
\]
Evaluating at \(t=0\) yields:
\[
\dot{\b{Y}}(0)
=
\b{\Delta}_1\b{S}(0)
+
(\b{X}+\b{\Delta}_2)\dot{\b{S}}(0).
\]
According to the definition of \(\b{S}(t)\), i.e., Eq. (\ref{equation_S_I_Delta2_tDelta1_Delta2_tDelta1_in_differentiated_polar_retraction_stiefel}), we have:
\[
\b{S}(0)
=
(\b{I}_d+\b{\Delta}_2^{\top}\b{\Delta}_2)^{-1/2}.
\]
Thus:
\[
\dot{\b{Y}}(0)
=
\b{\Delta}_1(\b{I}_d+\b{\Delta}_2^{\top}\b{\Delta}_2)^{-1/2}
+
(\b{X}+\b{\Delta}_2)\dot{\b{S}}(0).
\]
Hence, according to Eq. (\ref{equation_DRet_Delta2_Delta1_Ydot0_in_proof}), we have:
\begin{align*}
D\operatorname{Ret}^{polar}_{\b{X}}(\b{\Delta}_2)[\b{\Delta}_1]
=
&\b{\Delta}_1(\b{I}_d+\b{\Delta}_2^{\top}\b{\Delta}_2)^{-1/2}
\\
&+
(\b{X}+\b{\Delta}_2)\dot{\b{S}}(0).
\end{align*}

\textbf{Step 3: Verify that the derivative is tangent at \(\b{Y}\).}

Because \(\b{Y}(t)\) is obtained by a retraction on the Stiefel manifold,
we have:
\[
\b{Y}(t)\in St(n,d),
\qquad \forall t
\]
in a neighborhood of \(t=0\). Therefore:
\[
\b{Y}(t)^{\top}\b{Y}(t)=\b{I}_d.
\]
Differentiating with respect to \(t\) gives:
\[
\dot{\b{Y}}(t)^{\top}\b{Y}(t)+\b{Y}(t)^{\top}\dot{\b{Y}}(t)=\b{0}.
\]
Evaluating at \(t=0\), and using \(\b{Y}(0)=\b{Y}\), yields:
\[
\dot{\b{Y}}(0)^{\top}\b{Y}+\b{Y}^{\top}\dot{\b{Y}}(0)=\b{0}.
\]
By Eq. (\ref{equation_stiefel_tangent_space}), this is exactly the tangent constraint for
\(T_{\b{Y}}St(n,d)\). Hence:
\[
\dot{\b{Y}}(0)\in T_{\b{Y}}St(n,d).
\]
Therefore, we have:
\[
\mathcal{T}^{\operatorname{polar}}_{\b{\Delta}_2}(\b{\Delta}_1)
=
D\operatorname{Ret}^{polar}_{\b{X}}(\b{\Delta}_2)[\b{\Delta}_1]
\in T_{\b{Y}}St(n,d).
\]
This proves that the above formula is a valid vector transport
on the Stiefel manifold induced by the differentiated polar-style retraction.
\end{proof}

\begin{remark}[Special case of differentiated polar-style retraction at zero step in Stiefel manifold]
If \(\b{\Delta}_2=\b{0}\), then:
\[
\b{Y}
=
(\b{X}+\b{0})(\b{I}_d+\b{0}^{\top}\b{0})^{-1/2}
=
\b{X}.
\]
Also,
\[
\b{S}(0)=\b{I}_d.
\]
Hence:
\[
\mathcal{T}^{\operatorname{polar}}_{\b{0}}(\b{\Delta}_1)
=
\b{\Delta}_1+\b{X}\dot{\b{S}}(0).
\]
On the other hand, because:
\[
\b{S}(t)
=
(\b{I}_d+t^2\b{\Delta}_1^{\top}\b{\Delta}_1)^{-1/2},
\]
we have:
\[
\dot{\b{S}}(0)=\b{0}.
\]
Therefore:
\[
\mathcal{T}^{\operatorname{polar}}_{\b{0}}(\b{\Delta}_1)
=
\b{\Delta}_1.
\]
Thus, the differentiated polar-style retraction reduces to the
identity map at zero step, as expected.
\end{remark}

\begin{remark}[Interpretation of differentiated polar-style retraction in Stiefel manifold]
The formula in Eq. (\ref{equation_differentiated_polar_retraction_stiefel}) for differentiated polar-style retraction has two terms.
The first term:
\[
\b{\Delta}_1(\b{I}_d+\b{\Delta}_2^{\top}\b{\Delta}_2)^{-1/2},
\]
comes from differentiating the explicit factor \(\b{X}+\b{\Delta}\).
The second term:
\[
(\b{X}+\b{\Delta}_2)\dot{\b{S}}(0),
\]
comes from differentiating the normalization factor:
\[
(\b{I}_d+\b{\Delta}^{\top}\b{\Delta})^{-1/2}.
\]
Hence, the differentiated polar-style retraction consists of
the direct variation of the perturbed matrix together with the
variation required to maintain orthonormality. Their sum gives
the transported tangent vector in \(T_{\b{Y}}St(n,d)\).
\end{remark}

\subsection{Grassmann Manifold (Grassmannian)}





In Section \ref{section_definition_Stiefel_manifold}, we defined the Stiefel manifold $\text{St}(n, d)$ as the set of matrices with orthonormal columns. While the Stiefel manifold treats each specific basis as a distinct point, many applications in subspace learning and optimization care only about the \textit{subspace} spanned by those columns. This leads to the definition of the Grassmann manifold.

The Grassmann manifold is named after the German mathematician
\textit{Hermann Grassmann}, whose 1844 work
\textit{Die lineale Ausdehnungslehre} \cite{grassmann1844lineale} introduced the foundational
ideas that later led to what is now called the Grassmannian.
Many of the geometric and algorithmic characteristics of the
Grassmann manifold are analyzed in
\cite{edelman1998geometry,bendokat2024grassmann}.

\subsubsection{Quotient Operator and Equivalence Classes}

To formally define the Grassmann manifold, we must first establish the general definition of a quotient space and the operator $/$.

\begin{definition}[Equivalence relation]
An \textbf{equivalence relation} $\sim$ on a set $\mathcal{M}$ is a relation that satisfies:
\begin{enumerate}
    \item \textbf{Reflexivity:} $\b{X} \sim \b{X}$ for all $\b{X} \in \mathcal{M}$,
    \item \textbf{Symmetry:} $\b{X} \sim \b{Y} \implies \b{Y} \sim \b{X}$,
    \item \textbf{Transitivity:} $\b{X} \sim \b{Y}$ and $\b{Y} \sim \b{Z} \implies \b{X} \sim \b{Z}$,
\end{enumerate}
for all $\b{X}, \b{Y}, \b{Z} \in \mathcal{M}$.
\end{definition}

\begin{definition}[Quotient set, quotient operator, and quotient map]\label{definition_quotient_set}
The \textbf{quotient operator} $/$ acting on a set $\mathcal{M}$ with respect to an equivalence relation $\sim$ produces the \textbf{quotient set} $\mathcal{M} / \sim$. This set is defined as the collection of all disjoint equivalence classes:
\begin{equation}
\boxed{
\mathcal{M} / \sim \,\, := \{ [\b{X}] \mid \b{X} \in \mathcal{M} \},
}
\end{equation}
where $[\b{X}]$ denotes the equivalence class for $\b{X}$ in $\mathcal{M}$:
\begin{align}
\boxed{
[\b{X}] := \{ \b{Y} \in \mathcal{M} \mid \b{Y} \sim \b{X} \}.
}
\end{align}

Associated with this construction is the \textbf{quotient map} (or \textbf{canonical projection}), which maps each element of the original set to its corresponding equivalence class:
\begin{equation}
\boxed{
\begin{aligned}
& \pi: \mathcal{M} \to \mathcal{M} / \sim, \\
& \pi(\b{X}) = [\b{X}].
\end{aligned}
}
\end{equation}
\end{definition}

\begin{definition}[Quotient space]
A quotient set, introduced in Definition \ref{definition_quotient_set}, is a purely algebraic/set-theoretic construct.
A \textbf{quotient space} is a quotient set $\mathcal{M}/\sim$ that is endowed with a topology (and usually a differential structure).
In this topology, a subset of the quotient space is considered ``open" if its preimage under the projection map $\pi: \mathcal{M} \to \mathcal{M}/\sim$ is open in the original manifold $\mathcal{M}$.

\end{definition}

Recall Eq. (\ref{equation_On_orthogonal_group}) for orthogonal group. We repeat it in the following for convenience. 

\begin{definition}[Orthogonal group]
The \textbf{orthogonal group} of degree $d$, denoted by $\mathrm{O}(d)$, is the Lie group consisting of all $d \times d$ real orthogonal matrices. It is defined as the set:
\begin{equation}
\boxed{
\mathrm{O}(d) := \{ \b{Q} \in \mathbb{R}^{d \times d} \mid \b{Q}^\top \b{Q} = \b{I}_d \},} \label{eq:orthogonal_group_def}
\end{equation}
where $\b{I}_d$ is the $d \times d$ identity matrix. 
\end{definition}

\begin{definition}[Orbit]
An \textbf{orbit} $\text{orb}(\b{X})$ is a specific kind of equivalence class $[\b{X}]$ created when a group $G$, such as the orthogonal group $\mathrm{O}(d)$, ``acts" on the elements of a set $\mathcal{M}$. The orbit is defined as:
\begin{align}
\boxed{
\text{orb}(\b{X}) := \{ \b{X}g \mid g \in G \}.
}
\end{align}
In orbit, there must be a group $G$ and a defined action, such as matrix multiplication.
\end{definition}

\begin{remark}[Orbit as a special case of equivalence class]
If we define the equivalence relation as ``$\b{Y} \sim \b{X}$ if there exists $g \in G$ such that $\b{Y} = \b{X}g$", then the equivalence class $[\b{X}]$ is exactly the orbit $\text{orb}(\b{X})$. In other words:
\begin{align}
\boxed{
(\b{Y} \sim \b{X}) \equiv (\exists g \in G: \b{Y} = \b{X}g) \implies [\b{X}] = \text{orb}(\b{X}).
}
\end{align}
\end{remark}

\subsubsection{Definition of Grassmann Manifold}


\begin{definition}[Equivalence relation on the Stiefel manifold]\label{definition_equivalence_relation_on_Stiefel}
Two orthonormal matrices $\b{X}, \b{Y} \in \text{St}(n, d)$ are said to be \textbf{equivalent}, denoted $\b{X} \sim \b{Y}$, if they span the same $d$-dimensional subspace of $\mathbb{R}^n$. This occurs if and only if there exists an orthogonal matrix $\b{Q} \in \mathrm{O}(d)$ such that:
\begin{equation}
\b{Y} = \b{X}\b{Q}.
\end{equation}
Therefore, in Stiefel manifold, we can denote:
\begin{align}\label{equation_equivalence_Stiefel}
\boxed{
\b{X} \sim \b{Y} \iff \b{Y} = \b{X}\b{Q}, \quad \b{Q} \in \mathrm{O}(d).
}
\end{align}
\end{definition}

\begin{definition}[Quotient space $\text{St}(n, d) / \mathrm{O}(d)$]
Let $\text{St}(n, d)$ be the Stiefel manifold and $\mathrm{O}(d)$ be the orthogonal group. The \textbf{quotient space} $\text{St}(n, d) / \mathrm{O}(d)$ is the set of all orbits under the right action of $\mathrm{O}(d)$ on $\text{St}(n, d)$:
\begin{equation}
\boxed{
\text{St}(n, d) / \mathrm{O}(d) := \{ \text{orb}(\b{X}) \mid \b{X} \in \text{St}(n, d) \},
}
\end{equation}
where the orbit (or equivalence class) is given by:
\begin{equation}
\boxed{
\text{orb}(\b{X}) = \{ \b{X}\b{Q} \mid \b{Q} \in \mathrm{O}(d) \}.
}
\end{equation}
\end{definition}

\begin{definition}[Grassmann manifold]
The \textbf{Grassmann manifold} (also called \textbf{Grassmannian}), denoted by $\text{Gr}(n, d)$, is the set of all $d$-dimensional linear subspaces of $\mathbb{R}^n$, where $n \geq d$. Mathematically, it is the quotient space of the Stiefel manifold $\text{St}(n, d)$ under the equivalence action of the orthogonal group $\mathrm{O}(d)$:
\begin{equation}
\boxed{
\text{Gr}(n, d) \cong \text{St}(n, d) / \mathrm{O}(d), 
}\label{equation_grassmann_quotient}
\end{equation}
where two matrices $\b{X}, \b{Y} \in \text{St}(n, d)$ represent the same point in $\text{Gr}(n, d)$ if there exists an orthogonal matrix $\b{Q} \in \mathrm{O}(d)$ such that $\b{Y} = \b{X}\b{Q}$.
According to Eq. (\ref{equation_grassmann_quotient}), the Grassmann manifold is also called the \textbf{quotient manifold of the Stiefel manifold}. 
\end{definition}

\begin{remark}[Linear dimensionality reduction by Grassmann manifold]
The quotient operator $/$ in $\text{St}(n, d) / \mathrm{O}(d)$ performs a ``dimensionality reduction" \cite{ghojogh2023elements} by identifying all matrices that span the same subspace as a single point. 
That is why the Grassmann manifold represents all the $d$-dimensional linear subspaces of $\mathbb{R}^n$, where $n \geq d$.

In other words, projecting onto a linear subspace can be modeled by the Grassmann manifold. The linear dimensionality reduction methods---such as principal component analysis \cite{ghojogh2023principal} and Fisher's linear discriminant analysis \cite{ghojogh2023fisher}---use linear projection onto a subspace \cite{ghojogh2023elements}; thus, Grassman manifold can be useful for linear dimensionality reduction. 
\end{remark}

\begin{lemma}[Points of the Grassmann manifold are equivalence classes of Stiefel points]
Let \(\b{X} \in \mathbb{R}^{n\times d}\). If \(\b{X} \in St(n,d)\), then its equivalence class:
\begin{align}\label{equation_equivalence_X_in_Grassmann}
\boxed{
[\b{X}] := \{\b{X}\b{Q} \mid \b{Q}\in \mathrm{O}(d)\},
}
\end{align}
is a point of the Grassmann manifold \(Gr(n,d)\). Conversely, every point \(\b{p} \in Gr(n,d)\) can be written as:
\begin{align}
\boxed{
p = [\b{X}],
}
\end{align}
for some \(\b{X}\in St(n,d)\).

In other words, a matrix \(\b{X}\in St(n,d)\) determines a point \([\b{X}] \in Gr(n,d)\). Conversely, every point of \(Gr(n,d)\) admits at least one representative \(\b{X}\in St(n,d)\). Therefore, points of the Grassmann manifold are equivalence classes of points of the Stiefel manifold.
\end{lemma}

\begin{proof}
By Definition \ref{definition_equivalence_relation_on_Stiefel}, the Grassmann manifold is the quotient of the Stiefel manifold by the orthogonal group:
\[
Gr(n,d) := St(n,d)/\mathrm{O}(d).
\]
Hence, by the definition of quotient set in Definition \ref{definition_quotient_set}, the elements of \(Gr(n,d)\) are exactly the equivalence classes:
\begin{align*}
&[\b{X}] = \{\b{X}\b{Q} \mid \b{Q}\in \mathrm{O}(d)\}, \\
&\b{X}\in St(n,d),\ \b{Q}\in \mathrm{O}(d).
\end{align*}

Therefore, if \(\b{X}\in St(n,d)\), then \([\b{X}]\) is one of the equivalence classes in the quotient \(St(n,d)/\mathrm{O}(d)\), so:
\[
[\b{X}] \in Gr(n,d).
\]

Conversely, let \(\b{p} \in Gr(n,d)\). Since \(Gr(n,d)\) is the quotient set \(St(n,d)/\mathrm{O}(d)\), every element of \(Gr(n,d)\) is, by definition, an equivalence class of some matrix in \(St(n,d)\). Therefore, there exists \(\b{X}\in St(n,d)\) such that:
\[
p = [\b{X}].
\]

Hence, points of the Grassmann manifold are exactly equivalence classes of points of the Stiefel manifold.
\end{proof}

\begin{remark}[Mathematical interpretation of Grassmannian points]
An element (point) $\b{p}$ of the Grassmann manifold $\text{Gr}(n, d)$ admits several equivalent interpretations, spanning from set-theoretic to algebraic perspectives:
\begin{itemize}
    \item \textbf{Geometric View:} Each point $\b{p} \in \text{Gr}(n, d)$ represents a unique $d$-dimensional linear subspace $\mathcal{V} \subset \mathbb{R}^n$. Geometrically, this is a $d$-dimensional hyperplane passing through the origin.
    
    \item \textbf{Subspace View:} Given an orthonormal basis $\b{X} \in \text{St}(n, d)$, the point $\b{p}$ is the span of its columns:
    \begin{equation}
    \boxed{
    p = \text{span}(\b{X}) = \{ \b{X}\b{a} \mid \b{a} \in \mathbb{R}^d \}.
    }
    \end{equation}
    
    \item \textbf{Quotient View:} Because the manifold is defined as the quotient $\text{St}(n, d) / \mathrm{O}(d)$, a point $\b{p}$ is formally the equivalence class (or orbit) containing all orthonormal bases for that subspace:
    \begin{equation}\label{equation_p_bracket_X}
    \boxed{
    p = [\b{X}] = \{ \b{X}\b{Q} \mid \b{Q} \in \mathrm{O}(d) \}.
    }
    \end{equation}
    
    \item \textbf{Projector View:} A coordinate-independent representation of $\b{p}$ is the unique orthogonal projector $\b{P} = \b{X}\b{X}^\top$. This representation is invariant to the choice of basis since for any $\b{Q} \in \mathrm{O}(d)$, we have $(\b{X}\b{Q})(\b{X}\b{Q})^\top = \b{X}\b{Q}\b{Q}^\top\b{X}^\top = \b{X}\b{X}^\top$.
\end{itemize}
\end{remark}

\begin{remark}[Comparison of points in the Stiefel and Grassmann manifolds]
A point in Stiefel manifold and a point in Grassmann manifold have the following meanings:
\begin{itemize}
\item A single point in the Stiefel manifold $\text{St}(n, d)$ represents a specific set of $d$ orthonormal vectors (a frame).
\item A single point in the Grassmann manifold $\text{Gr}(n, d)$ represents the entire $d$-dimensional hyperplane itself. 
\end{itemize}

\end{remark}

\begin{remark}[Geometric interpretation of orthogonal group and Grassmann manifold]
Elements of the orthogonal group $\mathrm{O}(d)$ represent linear transformations that preserve the Euclidean inner product, lengths of vectors, and angles between vectors in $\mathbb{R}^d$. This group includes both rotations (where $\det(\b{Q}) = 1$) and reflections (where $\det(\b{Q}) = -1$). In the context of the Grassmann manifold $\text{Gr}(n, d)$, $\mathrm{O}(d)$ acts as the structure group that identifies different orthonormal bases for the same $d$-dimensional subspace.
\end{remark}

\begin{remark}[Decomposition of tangent space of the Stiefel manifold]
From a computational perspective, the quotient structure implies that the tangent space of the Stiefel manifold $T_{\b{X}}\text{St}(n, d)$ can be decomposed into two complementary subspaces:
\begin{enumerate}
    \item \textbf{Vertical Space ($\mathcal{V}_{\b{X}}$):} Tangent vectors that correspond to movements within the equivalence class (rotations of the basis $\b{Q}$). These changes do not change the point on the Grassmann manifold $\text{Gr}(n, d)$.
    \item \textbf{Horizontal Space ($\mathcal{H}_{\b{X}}$):} Tangent vectors that are orthogonal to the vertical space. These represent ``actual" movements from one subspace to another.
\end{enumerate}

This decomposition is crucial for defining the Riemannian metric and the gradient on the Grassmann manifold, as optimization algorithms must move ``horizontally" to effectively navigate the space of subspaces.
\end{remark}

\subsubsection{Dimension of Grassmann Manifold}





\begin{lemma}[Dimension of the orthogonal group]
The dimension of the orthogonal group $\mathrm{O}(d)$, consisting of all $d \times d$ matrices $\b{Q}$ such that $\b{Q}^\top \b{Q} = \b{I}_d$, is:
\begin{equation}
\dim(\mathrm{O}(d)) = \frac{d(d-1)}{2}. \label{equation_dim_orthogonal_group}
\end{equation}
\end{lemma}
\begin{proof}
An orthogonal matrix $\b{Q} \in \mathbb{R}^{d \times d}$ is constrained by the relation $\b{Q}^\top \b{Q} = \b{I}_d$. As noted in Remark \ref{remark_dimension_Stiefel_manifold}, the product $\b{Q}^\top \b{Q}$ is a symmetric $d \times d$ matrix, which contains $d(d+1)/2$ independent components. The total number of degrees of freedom in a general $d \times d$ matrix is $d^2$. Subtracting the number of independent constraints imposed by the orthogonality condition, we have:
\begin{align*}
\dim(\mathrm{O}(d)) &= d^2 - \frac{d(d+1)}{2} \\
&= \frac{2d^2 - d^2 - d}{2} = \frac{d^2 - d}{2} = \frac{d(d-1)}{2}.
\end{align*}
This completes the proof.
\end{proof}

\begin{proposition}[Dimension of the Grassmann Manifold]
The Grassmann manifold $\text{Gr}(n, d)$ is a smooth manifold of dimension:
\begin{equation}
\boxed{
\dim(\text{Gr}(n, d)) = d(n - d).} \label{eq:dim_grassmann}
\end{equation}
\end{proposition}
\begin{proof}
By the quotient manifold theorem (see \citep[Theorem 21.10]{lee2013smooth}), the dimension of a quotient space $\mathcal{M}/\mathcal{G}$ is given by $\dim(\mathcal{M}) - \dim(\mathcal{G})$. Applying this to Eq. \eqref{equation_grassmann_quotient}, we use the dimension of the Stiefel manifold derived in Eq. (\ref{equation_dimension_Stiefel_manifold}) and the dimension of the orthogonal group from Eq. \eqref{equation_dim_orthogonal_group}:
\begin{equation*}
\boxed{
\dim(\text{Gr}(n, d)) = \dim(\text{St}(n, d)) - \dim(\mathrm{O}(d)).
}
\end{equation*}
Substituting the expressions:
\begin{align*}
\dim(\text{Gr}(n, d)) &= \left( nd - \frac{d(d+1)}{2} \right) - \frac{d(d-1)}{2} \\
&= nd - \left( \frac{d^2 + d}{2} + \frac{d^2 - d}{2} \right) \\
&= nd - \left( \frac{2d^2}{2} \right) = nd - d^2.
\end{align*}
Factoring out $d$, we obtain:
\begin{equation*}
\dim(\text{Gr}(n, d)) = d(n - d).
\end{equation*}
This result signifies that specifying a $d$-dimensional subspace requires choosing $d$ vectors in $n$ dimensions, while removing the $d^2$ degrees of freedom associated with the internal choice of basis for that subspace.
\end{proof}

\subsubsection{Tangent and Normal Spaces of Grassmann Manifold}\label{section_tangent_normal_spaces_Grassmann}

The Grassmann manifold $\text{Gr}(n, d)$ is the set of all $d$-dimensional subspaces in $\mathbb{R}^n$. A point in $\text{Gr}(n, d)$ is an equivalence class $[\b{X}]$, where $\b{X} \in \text{St}(n, d)$ is an orthonormal basis for the subspace.

\begin{definition}[Direct sum of subspaces]
Let $\mathcal{U}$ and $\mathcal{V}$ be two subspaces of a vector space $\mathcal{W}$. The \textbf{direct sum} of $\mathcal{U}$ and $\mathcal{V}$, denoted by $\mathcal{W} = \mathcal{U} \oplus \mathcal{V}$, exists if:
\begin{enumerate}
    \item $\mathcal{W} = \{ \b{u} + \b{v} \mid \b{u} \in \mathcal{U}, \b{v} \in \mathcal{V} \}$, and
    \item $\mathcal{U} \cap \mathcal{V} = \{ \b{0} \}$.
\end{enumerate}
This implies that every element $\b{w} \in \mathcal{W}$ can be uniquely decomposed into a sum $\b{w} = \b{u} + \b{v}$.
\end{definition}

\begin{lemma}[Horizontal and vertical decomposition of tangent space of Stiefel manifold]\label{lemma_horizontal_vertical_decomposition_tangent_stiefel}
Let $\b{X} \in \text{St}(n, d)$. The tangent space $T_{\b{X}} \text{St}(n, d)$ can be decomposed into a vertical space $\mathcal{V}_{\b{X}}$ and a horizontal space $\mathcal{H}_{\b{X}}$ such that:
\begin{align}
\boxed{
T_{\b{X}} \text{St}(n, d) = \mathcal{V}_{\b{X}} \oplus \mathcal{H}_{\b{X}}.
}
\end{align}
The horizontal space, which is isomorphic to the tangent space of the Grassmannian $T_{[\b{X}]} \text{Gr}(n, d)$, is given by:
\begin{equation}\label{equation_HX_horizontal_space_Stiefel_tangent}
\boxed{
\mathcal{H}_{\b{X}} = \{ \b{\Delta} \in \mathbb{R}^{n \times d} \mid \b{\Delta}^\top \b{X} = \b{0} \}.
}
\end{equation}

The vertical space is given by:
\begin{equation}\label{equation_VX_vertical_space_Stiefel_tangent}
\boxed{
\mathcal{V}_{\b{X}} = \{ \b{X}\b{\Omega} \mid \b{\Omega} = -\b{\Omega}^\top \in \mathbb{R}^{d \times d} \}.
}
\end{equation}
\end{lemma}
\begin{proof}
\hfill\break
\textbf{Vertical Space:} The vertical space consists of tangents to the fiber $[\b{X}] = \{ \b{X}\b{Q} \mid \b{Q} \in \mathrm{O}(d) \}$. Let $\b{Q}(t)$ be a curve in $\mathrm{O}(d)$ with $\b{Q}(0) = \b{I}_d$ and $\dot{\b{Q}}(0) = \b{\Omega}$. 
As $\b{Q} \in \mathrm{O}(d)$, we have:
\begin{align*}
\b{Q}(t)^\top \b{Q}(t) = \b{I}_d.
\end{align*}
Differentiating this expression with respect to $t$ using the product rule:
\begin{equation*}
\dot{\b{Q}}(t)^\top \b{Q}(t) + \b{Q}(t)^\top \dot{\b{Q}}(t) = \b{0}.
\end{equation*}
Evaluating at $t=0$, and letting $\dot{\b{Q}}(0) = \b{\Omega}$:
\begin{equation*}
\b{\Omega}^\top \b{I}_d + \b{I}_d^\top \b{\Omega} = \b{\Omega}^\top + \b{\Omega} = \b{0}.
\end{equation*}
This confirms that $\b{\Omega}$ is a skew-symmetric matrix:
\begin{align*}
\b{\Omega}^\top + \b{\Omega} = \b{0} \implies \b{\Omega} = -\b{\Omega}^\top,
\end{align*}
which is an element of the Lie algebra $\mathfrak{so}(d)$.
Thus, the vertical space is:
\begin{equation*}
\mathcal{V}_{\b{X}} = \{ \b{X}\b{\Omega} \mid \b{\Omega} = -\b{\Omega}^\top \in \mathbb{R}^{d \times d} \}.
\end{equation*}

\textbf{Stiefel Tangent Space Decomposition:} From Eq. (\ref{equation_stiefel_tangent_space}), any $\b{\Delta}_{St} \in T_{\b{X}}\text{St}(n, d)$ satisfies:
\begin{align}\label{equation_XT_DeltaST_DeltaSTT_X_zero}
\b{X}^\top \b{\Delta}_{St} + \b{\Delta}_{St}^\top \b{X} = \b{0}.
\end{align}
We decompose $\b{\Delta}_{St}$ using the basis $\b{X}$ and $\b{X}_\perp$ where $\b{X}_\perp \in \mathbb{R}^{n \times (n-d)}$ is the orthogonal complement:
\begin{equation}\label{equation_Delta_st_XA_XperpB}
\b{\Delta}_{St} = \b{X}\b{A} + \b{X}_\perp \b{B}.
\end{equation}
Substituting Eq. (\ref{equation_Delta_st_XA_XperpB}) into Eq. (\ref{equation_XT_DeltaST_DeltaSTT_X_zero}), gives:
\begin{align*}
&\b{X}^\top (\b{X}\b{A} + \b{X}_\perp \b{B}) + (\b{X}\b{A} + \b{X}_\perp \b{B})^\top \b{X} = \b{0} \\
&\implies (\b{X}^\top \b{X})\b{A} + \b{X}^\top \b{X}_\perp \b{B} \\
&\quad\quad\quad\quad\quad\quad\quad+ \b{A}^\top (\b{X}^\top \b{X}) + \b{B}^\top \b{X}_\perp^\top \b{X} = \b{0} \\
&\overset{(a)}{\implies} \b{I}_d \b{A} + \b{0} + \b{A}^\top \b{I}_d + \b{0} = \b{0} \implies \b{A} = -\b{A}^\top,
\end{align*}
where $(a)$ is because $\b{X}^\top \b{X} = \b{I}_d$ and $\b{X}^\top \b{X}_\perp = \b{X}_\perp^\top \b{X} = \b{0}$.

\textbf{Orthogonality of Horizontal and Vertical Spaces:} 
We already proved the vertical space $\mathcal{V}_{\b{X}} = \{ \b{X}\b{\Omega} \mid \b{\Omega} = -\b{\Omega}^\top \in \mathbb{R}^{d \times d} \}$.
The horizontal space $\mathcal{H}_{\b{X}}$ contains vectors $\b{\Delta}$ such that $\langle \b{\Delta}, \b{V} \rangle = 0$ for all $\b{V} \in \mathcal{V}_{\b{X}}$ under the metric $\langle \b{A}, \b{B} \rangle = \text{tr}(\b{A}^\top \b{B})$. 
For $\b{V} = \b{X}\b{\Omega} \in \mathcal{V}_{\b{X}}$, we have:
\begin{align*}
&\langle \b{X}\b{A} + \b{X}_\perp \b{B}, \b{X}\b{\Omega} \rangle = \text{tr}((\b{X}\b{A} + \b{X}_\perp \b{B})^\top \b{X}\b{\Omega}) = 0 \\
&\implies \text{tr}((\b{A}^\top \b{X}^\top + \b{B}^\top \b{X}_\perp^\top) \b{X}\b{\Omega}) = 0 \\
&\implies \text{tr}(\b{A}^\top \b{X}^\top \b{X}\b{\Omega} + \b{B}^\top \b{X}_\perp^\top \b{X}\b{\Omega}) = 0 \\
&\overset{(a)}{\implies} \text{tr}(\b{A}^\top \b{\Omega}) = 0,
\end{align*}
where $(a)$ is because $\b{X}^\top \b{X} = \b{I}_d$ and $\b{X}_\perp^\top \b{X} = \b{0}$.

Since this must hold for all skew-symmetric $\b{\Omega}$, and $\b{A}$ is already known to be skew-symmetric from the Stiefel condition, i.e., Eq. (\ref{equation_stiefel_tangent_space}), we must have $\b{A} = \b{0}$.

\textbf{Conclusion:} Setting $\b{A} = \b{0}$, the tangent vector in Eq. (\ref{equation_Delta_st_XA_XperpB}) reduces to $\b{\Delta} = \b{X}_\perp \b{B}$. 
Substituting $\b{\Delta} = \b{X}_\perp \b{B}$ in $\b{\Delta}^\top \b{X}$ gives: 
\begin{equation}
\b{\Delta}^\top \b{X} = (\b{X}_\perp \b{B})^\top \b{X} = \b{B}^\top \b{X}_\perp^\top \b{X} \overset{(a)}{=} \b{0},
\end{equation}
where $(a)$ is because $\b{X}_\perp^\top \b{X} = \b{0}$.
Thus, the horizontal space is:
\begin{align*}
\mathcal{H}_{\b{X}} = \{ \b{\Delta} \in \mathbb{R}^{n \times d} \mid \b{\Delta}^\top \b{X} = \b{0} \}.
\end{align*}
\end{proof}

\begin{definition}[Horizontal lift on Grassmann manifold]\label{definition_horizontal_lift_Grassmann}
Let $[\b{X}] \in \text{Gr}(n, d)$ be a point on the Grassmann manifold and $\b{\Theta} \in T_{[\b{X}]} \text{Gr}(n, d)$ be a tangent vector. For any representative $\b{X}$ in the equivalence class $[\b{X}]$, there exists a unique tangent vector $\b{\Delta} \in \mathcal{H}_{\b{X}} \subset T_{\b{X}} \text{St}(n, d)$ such that:
\begin{align}
(\pi_{\b{X}})_* (\b{\Delta}) = \b{\Theta},
\end{align}
where $(\pi_{\b{X}})_*(.)$ is the pushforward of the quotient map $\pi: \text{St}(n, d) \to \text{Gr}(n, d)$ at $\b{X}$. This unique matrix $\b{\Delta}$ is called the \textbf{horizontal lift} of $\b{\Theta}$ at $\b{X}$.
\end{definition}

\begin{remark}[Computation of tangent vector of Grassmann manifold]\label{remark_computation_tangent_Grassmann}
The tangent vector (matrix) of Grassmann manifold, i.e., $\b{\Theta} \in T_{[\b{X}]} \text{Gr}(n, d)$, is an abstract object. 
In Riemannian optimization, we cannot numerically manipulate the abstract tangent vector (matrix) $\b{\Theta}$ directly. Instead, we perform calculations using its horizontal lift $\b{\Delta}$, which is a concrete $n \times d$ matrix satisfying the constraint $\b{\Delta}^\top \b{X} = \b{0}$ derived in Lemma \ref{lemma_horizontal_vertical_decomposition_tangent_stiefel}.
\end{remark}


\begin{proposition}[Tangent space of Grassmann manifold]\label{prop_grassmann_tangent}
The tangent space $T_{[\b{X}]}\text{Gr}(n, d)$ at a point represented by $\b{X} \in \text{St}(n, d)$ is given by:
\begin{equation}\label{equation_tangent_space_Grassmann}
\boxed{
T_{[\b{X}]}\text{Gr}(n, d) \cong H_X = \{ \b{\Delta} \in \mathbb{R}^{n \times d} \mid \b{\Delta}^\top \b{X} = \b{0} \}.
}
\end{equation}
\end{proposition}
\begin{proof}
According to Lemma \ref{lemma_horizontal_vertical_decomposition_tangent_stiefel}, Definition \ref{definition_horizontal_lift_Grassmann}, and Remark \ref{remark_computation_tangent_Grassmann}, we use the horizontal lift of $\b{\Delta} \in \mathcal{H}_{\b{X}} \subset T_{\b{X}} \text{St}(n, d)$ as the tangent vector (matrix) of Grassmann manifold. The $\b{\Delta}$ must satisfy Eq. (\ref{equation_HX_horizontal_space_Stiefel_tangent}). 
\end{proof}

\begin{remark}[Comparison of tangent constraint in the Stiefel and Grassmann manifolds]
The Grassmann tangent constraint is simpler than the Stiefel tangent constraint. On the Stiefel manifold, according to Eq. (\ref{equation_stiefel_tangent_space}), the tangent constraint is:
\[
\text{Tangent on Stiefel: } \b{X}^\top \b{\Delta} + \b{\Delta}^\top \b{X} = \b{0},
\]
so the projection involves the symmetrization operator, defined in Eq. (\ref{equation_sym_skew_expressions}).
On the Grassmann manifold, according to Eq. (\ref{equation_tangent_space_Grassmann}), the tangent constraint reduces to:
\[
\text{Tangent on Grassmann: } \b{X}^\top \b{\Delta} = \b{0} \,\text{ or }\, \b{\Delta}^\top \b{X} = \b{0}.
\]
\end{remark}

\begin{proposition}[Normal space of Grassmann manifold]
The normal space $N_{[\b{X}]}\text{Gr}(n, d)$ in the ambient space $\mathbb{R}^{n \times d}$ is:
\begin{equation}
\begin{aligned}
N_{[\b{X}]}\text{Gr}(n, d) = &\{ \b{X}\b{S} \mid \b{S} = \b{S}^\top \in \mathbb{R}^{d \times d} \} \oplus \\
&\{ \b{X}\b{\Omega} \mid \b{\Omega} = -\b{\Omega}^\top \in \mathbb{R}^{d \times d} \},
\end{aligned}
\end{equation}
which simplifies to any matrix whose columns lie in the subspace $\text{span}(\b{X})$:
\begin{equation}\label{equation_normal_space_Grassmann}
\boxed{
N_{[\b{X}]}\text{Gr}(n, d) = \{ \b{X}\b{K} \mid \b{K} \in \mathbb{R}^{d \times d} \}.
}
\end{equation}
\end{proposition}
\begin{proof}
The normal space $N_{[\b{X}]}\text{Gr}(n, d)$ is defined as the orthogonal complement of the tangent space $T_{[\b{X}]}\text{Gr}(n, d)$ within the ambient space $\mathbb{R}^{n \times d}$ under the Frobenius inner product $\langle \b{A}, \b{B} \rangle = \text{tr}(\b{A}^\top \b{B})$.

\textbf{Ambient Decomposition:} Let $\b{Z} \in \mathbb{R}^{n \times d}$ be an arbitrary matrix. Using the orthogonal completion $\b{X}$ and $\b{X}_\perp$ where $\b{X} \in \text{St}(n, d)$ and $\b{X}^\top \b{X}_\perp = \b{0}$, we decompose $\b{Z}$ as:
\begin{equation}\label{equation_Z_XK_XperpL}
\b{Z} = \b{X}\b{K} + \b{X}_\perp \b{L},
\end{equation}
where $\b{K} \in \mathbb{R}^{d \times d}$ and $\b{L} \in \mathbb{R}^{(n-d) \times d}$.

\textbf{Orthogonality Condition:} From Proposition \ref{prop_grassmann_tangent}, any tangent vector $\b{\Delta} \in T_{[\b{X}]}\text{Gr}(n, d)$ is of the form $\b{\Delta} = \b{X}_\perp \b{B}$ for some $\b{B} \in \mathbb{R}^{(n-d) \times d}$. For $\b{Z}$ to be in the normal space, we require $\langle \b{Z}, \b{\Delta} \rangle = 0$ for all $\b{B}$:
\begin{align*}
&\text{tr}(\b{Z}^\top \b{\Delta}) = 0 \overset{(\ref{equation_Z_XK_XperpL})}{\implies} \text{tr}((\b{X}\b{K} + \b{X}_\perp \b{L})^\top \b{X}_\perp \b{B}) = 0 \\
&\implies \text{tr}(\b{K}^\top \b{X}^\top \b{X}_\perp \b{B} + \b{L}^\top \b{X}_\perp^\top \b{X}_\perp \b{B}) = 0.
\end{align*}
Since $\b{X}^\top \b{X}_\perp = \b{0}$ and $\b{X}_\perp^\top \b{X}_\perp = \b{I}_{n-d}$, this simplifies to:
\begin{equation*}
\text{tr}(\b{L}^\top \b{B}) = 0.
\end{equation*}
This must hold for all $\b{B} \in \mathbb{R}^{(n-d) \times d}$, which implies $\b{L} = \b{0}$. 

\textbf{Simplification:} Consequently, according to Eq. (\ref{equation_Z_XK_XperpL}), $\b{L} = \b{0}$ implies $\b{Z} = \b{X}\b{K}$ for an arbitrary $\b{K} \in \mathbb{R}^{d \times d}$. This proves the second part of the proposition: $N_{[\b{X}]}\text{Gr}(n, d) = \{ \b{X}\b{K} \mid \b{K} \in \mathbb{R}^{d \times d} \}$.

\textbf{Direct Sum Decomposition:} Any square matrix $\b{K}$ can be uniquely decomposed into a symmetric part $\b{S} = \frac{1}{2}(\b{K} + \b{K}^\top)$ and a skew-symmetric part $\b{\Omega} = \frac{1}{2}(\b{K} - \b{K}^\top)$. Since the intersection of symmetric and skew-symmetric subspaces is $\{ \b{0} \}$, we have:
\begin{align*}
&\b{X}\b{K} = \b{X}\b{S} + \b{X}\b{\Omega} \\
&\implies \{ \b{X}\b{S} \mid \b{S} = \b{S}^\top \} \oplus \{ \b{X}\b{\Omega} \mid \b{\Omega} = -\b{\Omega}^\top \}.
\end{align*}
This completes the proof.
\end{proof}

\begin{lemma}[Projection onto the tangent space of Grassmann manifold]\label{lemma_projection_onto_tangent_Grassmann}
Let $\b{X} \in St(n,d)$ be a representative of $[\b{X}] \in Gr(n,d)$. For any $\b{Z} \in \mathbb{R}^{n\times d}$, the orthogonal projection onto $T_{[\b{X}]}Gr(n,d)$ is:
\begin{equation}\label{equation_projection_onto_tangent_space_Grassmann}
\boxed{
\Pi_{[X]}^{\mathrm{Gr}}(\b{Z}) = (\b{I} - \b{X}\b{X}^\top)\b{Z}.
}
\end{equation}
\end{lemma}

\begin{proof}
We define:
\[
\b{\Delta} := (\b{I} - \b{X}\b{X}^\top)\b{Z}.
\]

\textbf{Step 1: Tangency.}
\begin{align*}
\b{\Delta}^\top \b{X}
&= \b{Z}^\top (\b{I} - \b{X}\b{X}^\top)\b{X} \\
&= \b{Z}^\top(\b{X} - \b{X}\b{X}^\top \b{X}) \overset{(a)}{=} \b{Z}^\top(\b{X} - \b{X}) = \b{0},
\end{align*}
where $(a)$ is because $\b{X}^\top \b{X} = \b{I}$. 
As $\b{\Delta}^\top \b{X} = \b{0}$, we have $\b{\Delta} \in T_{[\b{X}]}Gr(n,d)$ according to Eq. (\ref{equation_tangent_space_Grassmann}).

\textbf{Step 2: Normal component.}
The term $\b{Z} - \Pi_{[\b{X}]}^{\mathrm{Gr}}(\b{Z})$ lies in the normal space:
\[
\b{Z} - \Pi_{[\b{X}]}^{\mathrm{Gr}}(\b{Z}) = \b{X}\b{X}^\top \b{Z} \overset{(a)}{=} \b{X}\b{K},
\]
where $(a)$ is because we define $\b{K} := \b{X}^\top \b{Z}$.
This $\b{X}\b{K}$ lies in the normal space.

\textbf{Step 3: Orthogonality.}
For any $\b{\Delta} \in T_{[\b{X}]}Gr(n,d)$:
\begin{align*}
\langle \b{X}\b{X}^\top \b{Z}, \b{\Delta} \rangle_F 
&= \text{tr}((\b{X}\b{X}^\top \b{Z})^\top \b{\Delta}) = \text{tr}(\b{Z}^\top \b{X}\b{X}^\top \b{\Delta}) \\
&= \text{tr}(\b{Z}^\top \b{X} (\b{X}^\top \b{\Delta})) \overset{(a)}{=} 0,
\end{align*}
where $(a)$ is because $\b{X}^\top \b{\Delta} = \b{\Delta}^\top \b{X} = \b{0}$ according to Eq. (\ref{equation_tangent_space_Grassmann}).

Therefore, the decomposition is orthogonal, and the projection formula holds.
\end{proof}

\begin{lemma}[Projected ambient derivative is tangent on the Grassmann manifold]
Let $\b{X} \in St(n, d)$ be a representative of $[\b{X}] \in Gr(n, d)$, and let $\b{Z} \in \mathbb{R}^{n \times d}$ be any matrix. Then, the projected ambient derivative is tangent on the Grassmann manifold:
\begin{equation}\label{equation_projection_connection_grassmann_aux}
\boxed{
\Pi_{[\b{X}]}^{\mathrm{Gr}}(\b{Z}) = (\b{I}_n - \b{X}\b{X}^\top)\b{Z} \in T_{[\b{X}]}Gr(n, d).
}
\end{equation}
\end{lemma}

\begin{proof}
By Lemma \ref{lemma_projection_onto_tangent_Grassmann}, the orthogonal projection onto the tangent space of the Grassmann manifold is:
\[
\Pi_{[\b{X}]}^{\mathrm{Gr}}(\b{Z}) = (\b{I}_n - \b{X}\b{X}^\top)\b{Z}.
\]
We verify the tangent constraint in Eq. (\ref{equation_tangent_space_Grassmann}):
\begin{align*}
\b{X}^\top \Pi_{[\b{X}]}^{\mathrm{Gr}}(\b{Z})
&= \b{X}^\top (\b{I}_n - \b{X}\b{X}^\top)\b{Z}
\\
&= (\b{X}^\top - \b{X}^\top \b{X} \b{X}^\top)\b{Z} \\
&\overset{(a)}{=} (\b{X}^\top - \b{X}^\top)\b{Z} = \b{0},
\end{align*}
where $(a)$ is because $\b{X}^\top \b{X} = \b{I}_d$ since $\b{X} \in \text{St}(n,d)$; see Eq. (\ref{equation_Sitefel_manifold_definition}).
Hence, $\Pi_{[\b{X}]}^{\mathrm{Gr}}(\b{Z})$ satisfies the condition in Eq. (\ref{equation_tangent_space_Grassmann}). Therefore, $\Pi_{[\b{X}]}^{\mathrm{Gr}}(\b{Z}) \in T_{[\b{X}]}Gr(n,d)$.
\end{proof}

\subsubsection{Metric Tensor of Grassmann Manifold}\label{section_metric_tensor_Grassmann}

The Grassmann manifold $\text{Gr}(n, d)$ is defined as the quotient manifold of the Stiefel manifold $\text{St}(n, d)$ under the action of the orthogonal group $\mathrm{O}(d)$. Because it is a quotient manifold, the metric on $\text{Gr}(n, d)$ is inherited from the metric on the total space $\text{St}(n, d)$. To ensure this metric is well-defined on the quotient space, it must be invariant under the action of $\mathrm{O}(d)$.



\begin{proposition}[Riemannian metric on the Grassmann manifold]\label{proposition_metric_Grassmann}
The standard Riemannian metric on the Grassmann manifold $\text{Gr}(n, d)$ at a point $[\b{X}]$ is the metric induced by the Euclidean metric of the Stiefel manifold. For two tangent vectors representing subspaces, let $\b{\Delta}_1, \b{\Delta}_2 \in \mathcal{H}_{\b{X}}$ be their unique horizontal lifts in the tangent space of the Stiefel manifold. The metric on the Grassmann manifold is:
\begin{equation}\label{equation_metric_Grassmann}
\boxed{
\begin{aligned}
g^E_{[\b{X}]}(\b{\Delta}_1, \b{\Delta}_2) &= \text{tr}(\b{\Delta}_1^\top \b{\Delta}_2) \\
&\overset{(\ref{equation_Frobenius_inner_product})}{=} \langle \b{\Delta}_1, \b{\Delta}_2 \rangle_F.
\end{aligned}
}
\end{equation}
This is the quotient metric induced by the Euclidean inner product on the Stiefel total space.
\end{proposition}
\begin{proof}
A metric on a quotient manifold is well-defined if the metric on the total space is invariant under the group action. Let $\b{Q} \in \mathrm{O}(d)$. According to Eq. (\ref{equation_equivalence_Stiefel}), the action is $\phi_{\b{Q}}(\b{X}) = \b{X}\b{Q}$. The pushforward of a tangent vector $\b{\Delta}$ is $\phi_{\b{Q}*} \b{\Delta} = \b{\Delta} \b{Q}$.

We check if the Euclidean metric is invariant under $\mathrm{O}(d)$:
\begin{align*}
&g^E_{\b{X}\b{Q}}(\b{\Delta}_1 \b{Q}, \b{\Delta}_2 \b{Q}) = \text{tr}((\b{\Delta}_1 \b{Q})^\top (\b{\Delta}_2 \b{Q})) \\
&= \text{tr}(\b{Q}^\top \b{\Delta}_1^\top \b{\Delta}_2 \b{Q}) \overset{(a)}{=} \text{tr}(\b{Q} \b{Q}^\top \b{\Delta}_1^\top \b{\Delta}_2) \\
&\overset{(b)}{=} \text{tr}(\b{\Delta}_1^\top \b{\Delta}_2) \overset{(\ref{equation_Euclidean_metric_Stiefel})}{=} g^E_{\b{X}}(\b{\Delta}_1, \b{\Delta}_2),
\end{align*}
where $(a)$ is because of the cyclic property of the trace and $(b)$ is because of $\b{Q}\b{Q}^\top = \b{I}$ for $\b{Q} \in \mathrm{O}(d)$. 

Since the metric is invariant under the orthogonal group action, it descends to a well-defined Riemannian metric on the quotient manifold $\text{Gr}(n, d)$. 
\end{proof}

\subsubsection{Levi-Civita Connection in Grassmann Manifold}

\begin{proposition}[Levi-Civita connection on the Grassmann manifold]\label{proposition_levi_civita_connection_grassmann}
Let $\b{\Delta}_1,\b{\Delta}_2 \in T_{[\b{X}]}Gr(n, d)$ be tangent vector fields. 
Under the Euclidean metric $g^E_{[\b{X}]}(\b{\Delta}_1, \b{\Delta}_2) = \text{tr}(\b{\Delta}_1^\top \b{\Delta}_2)
= \langle \b{\Delta}_1, \b{\Delta}_2 \rangle_F$, the Levi-Civita connection on the Grassmann manifold $Gr(n,d)$ is:
\begin{equation}\label{equation_levi_civita_connection_grassmann}
\boxed{
\begin{aligned}
(\nabla^E_{\b{\Delta}_1}\b{\Delta}_2)(\b{X})
&=
\Pi_{[\b{X}]}^{\mathrm{Gr}}\!\big(D\b{\Delta}_2(\b{X})[\b{\Delta}_1]\big)
\\
&\overset{(\ref{equation_projection_onto_tangent_space_Grassmann})}{=}
(\b{I}_n - \b{X}\b{X}^\top)\,D\b{\Delta}_2(\b{X})[\b{\Delta}_1].
\end{aligned}
}
\end{equation}
\end{proposition}

\begin{proof}
As the metric is the Euclidean metric, this is directly obtained by Eq. (\ref{equation_levi_civita_connection_ambient_projection}) in Lemma \ref{lemma_levi_civita_ambient_directional_derivative_relation}. However, for the sake of completeness, we also prove it in the following, too. 

We verify torsion-freeness and metric compatibility.

\textbf{Step 1: Torsion-freeness.}
We have:
\begin{align*}
&(\nabla^E_{\b{\Delta}_1}\b{\Delta}_2)(\b{X}) - (\nabla^E_{\b{\Delta}_2}\b{\Delta}_1)(\b{X}) \\
&=
\Pi_{[\b{X}]}^{\mathrm{Gr}}\!\big(D\b{\Delta}_2(\b{X})[\b{\Delta}_1]\big)
-
\Pi_{[\b{X}]}^{\mathrm{Gr}}\!\big(D\b{\Delta}_1(\b{X})[\b{\Delta}_2]\big) \\
&=
\Pi_{[\b{X}]}^{\mathrm{Gr}}\!\big(D\b{\Delta}_2(\b{X})[\b{\Delta}_1]
-
D\b{\Delta}_1(\b{X})[\b{\Delta}_2]\big) \\
&\overset{(\ref{equation_lie_bracket_directional_derivative})}{=} \Pi_{[\b{X}]}^{\mathrm{Gr}}([\b{\Delta}_1,\b{\Delta}_2]) \overset{(a)}{=} [\b{\Delta}_1,\b{\Delta}_2],
\end{align*}
where $(a)$ is because the Lie bracket is tangent so its projection onto the tangent space is itself. 
We proved that $\nabla^E_{\b{\Delta}_1}\b{\Delta}_2 - \nabla^E_{\b{\Delta}_2}\b{\Delta}_1 = [\b{\Delta}_1,\b{\Delta}_2]$.
Therefore, according to Eq. (\ref{equation_torsion_free}), the connection $\nabla$ is torsion-free. 

\textbf{Step 2: Metric compatibility.}
The Grassmann metric is:
\[
g^E_{[\b{X}]}(\b{\Delta}_1,\b{\Delta}_2)
=
\operatorname{tr}(\b{\Delta}_1^\top \b{\Delta}_2).
\]
We compute its directional derivative along $\b{\Delta}_0$:
\begin{align}
\b{\Delta}_0\big(&g^E_{[\b{X}]}(\b{\Delta}_1,\b{\Delta}_2)\big) \nonumber\\
&=
\b{\Delta}_0\!\left(\operatorname{tr}(\b{\Delta}_1^\top \b{\Delta}_2)\right) \nonumber\\
&=
\operatorname{tr}\!\left(\b{\Delta}_0(\b{\Delta}_1^\top \b{\Delta}_2)\right) \nonumber\\
&\overset{(\ref{equation_DYpXp_XYk_p_partialkp})}{=}
\operatorname{tr}\!\left((D\b{\Delta}_1[\b{\Delta}_0])^\top \b{\Delta}_2\right)
+
\operatorname{tr}\!\left(\b{\Delta}_1^\top D\b{\Delta}_2[\b{\Delta}_0]\right) \nonumber\\
&\overset{(\ref{equation_Frobenius_inner_product})}{=}
\langle D\b{\Delta}_1[\b{\Delta}_0], \b{\Delta}_2\rangle_F
+
\langle \b{\Delta}_1, D\b{\Delta}_2[\b{\Delta}_0]\rangle_F. \label{equation_Delta0_g_Delta1_Delta2}
\end{align}

Now, we decompose the ambient derivative into tangent and normal components:
\[
D\b{\Delta}_1[\b{\Delta}_0]
=
\Pi_{[\b{X}]}^{\mathrm{Gr}}\big(D\b{\Delta}_1[\b{\Delta}_0]\big)
+
(\b{I}-\Pi_{[\b{X}]}^{\mathrm{Gr}})\big(D\b{\Delta}_1[\b{\Delta}_0]\big).
\]
Hence:
\begin{align*}
\langle D\b{\Delta}_1[\b{\Delta}_0], \b{\Delta}_2\rangle_F
&=
\Big\langle
\Pi_{[\b{X}]}^{\mathrm{Gr}}\big(D\b{\Delta}_1[\b{\Delta}_0]\big),
\b{\Delta}_2
\Big\rangle_F \\
&+
\Big\langle
(\b{I}-\Pi_{[\b{X}]}^{\mathrm{Gr}})\big(D\b{\Delta}_1[\b{\Delta}_0]\big),
\b{\Delta}_2
\Big\rangle_F.
\end{align*}
Since $\b{\Delta}_2 \in T_{[\b{X}]}Gr(n,d)$ is tangent and $(\b{I}-\Pi_{[\b{X}]}^{\mathrm{Gr}})\big(D\b{\Delta}_1[\b{\Delta}_0]\big)$ is in the normal space:
\[
(\b{I}-\Pi_{[\b{X}]}^{\mathrm{Gr}})\big(D\b{\Delta}_1[\b{\Delta}_0]\big)
\in N_{[\b{X}]}Gr(n,d),
\]
the second term is zero by orthogonality of tangent and normal spaces. Therefore:
\begin{align}\label{equation_DDelta1Delta0_Delta2}
\langle D\b{\Delta}_1[\b{\Delta}_0], \b{\Delta}_2\rangle_F
=
\Big\langle
\Pi_{[\b{X}]}^{\mathrm{Gr}}\big(D\b{\Delta}_1[\b{\Delta}_0]\big),
\b{\Delta}_2
\Big\rangle_F.
\end{align}
Similarly, we have:
\begin{align}\label{equation_Delta1DDelta2Delta0}
\langle \b{\Delta}_1, D\b{\Delta}_2[\b{\Delta}_0]\rangle_F
=
\Big\langle
\b{\Delta}_1,
\Pi_{[\b{X}]}^{\mathrm{Gr}}\big(D\b{\Delta}_2[\b{\Delta}_0]\big)
\Big\rangle_F.
\end{align}
Substituting Eqs. (\ref{equation_DDelta1Delta0_Delta2}) and (\ref{equation_Delta1DDelta2Delta0}) into Eq. (\ref{equation_Delta0_g_Delta1_Delta2}) gives:
\begin{align}
\b{\Delta}_0\big(g^E_{[\b{X}]}(\b{\Delta}_1,\b{\Delta}_2)\big)
=
&\Big\langle
\Pi_{[\b{X}]}^{\mathrm{Gr}}\big(D\b{\Delta}_1[\b{\Delta}_0]\big),
\b{\Delta}_2
\Big\rangle_F \nonumber
\\
&+
\Big\langle
\b{\Delta}_1,
\Pi_{[\b{X}]}^{\mathrm{Gr}}\big(D\b{\Delta}_2[\b{\Delta}_0]\big)
\Big\rangle_F. \label{equation_Delta0_g_Delta1_Delta2_2}
\end{align}

By the definition of the proposed connection in Eq. (\ref{equation_levi_civita_connection_grassmann}), we have:
\begin{align*}
&(\nabla^E_{\b{\Delta}_0}\b{\Delta}_1)(\b{X})
=
\Pi_{[\b{X}]}^{\mathrm{Gr}}\big(D\b{\Delta}_1(\b{X})[\b{\Delta}_0]\big), \\
&(\nabla^E_{\b{\Delta}_0}\b{\Delta}_2)(\b{X})
=
\Pi_{[\b{X}]}^{\mathrm{Gr}}\big(D\b{\Delta}_2(\b{X})[\b{\Delta}_0]\big).
\end{align*}
Therefore, Eq. (\ref{equation_Delta0_g_Delta1_Delta2_2}) becomes:
\begin{align*}
\b{\Delta}_0\big(g^E_{[\b{X}]}(\b{\Delta}_1,\b{\Delta}_2)\big)
=\,
&\langle \nabla^E_{\b{\Delta}_0}\b{\Delta}_1, \b{\Delta}_2\rangle_F
\\
&+
\langle \b{\Delta}_1, \nabla^E_{\b{\Delta}_0}\b{\Delta}_2\rangle_F.
\end{align*}
According to Eq. (\ref{equation_metric_Grassmann}) in Proposition \ref{proposition_metric_Grassmann}, the Grassmann metric is the Frobenius inner product on tangent vectors. Therefore, we conclude that:
\begin{align*}
\b{\Delta}_0\big(g^E_{[\b{X}]}(\b{\Delta}_1,\b{\Delta}_2)\big)
=\,
&g^E_{[\b{X}]}\big(\nabla^E_{\b{\Delta}_0}\b{\Delta}_1,\b{\Delta}_2\big)
\\
&+
g^E_{[\b{X}]}\big(\b{\Delta}_1,\nabla^E_{\b{\Delta}_0}\b{\Delta}_2\big).
\end{align*}
Hence, according to Eq. (\ref{equation_metric_compatibility}), the connection $\nabla$ is metric-compatible.

Therefore, since the connection is both torsion-free and metric-compatible, Eq. \eqref{equation_levi_civita_connection_grassmann} is the Levi-Civita connection according to Definition \ref{definition_Levi_Civita_connection_coordinate_free}.
\end{proof}

\begin{remark}[Projection interpretation of the Grassmann
Levi-Civita connection]
In Proposition \ref{proposition_levi_civita_connection_grassmann}, we proved that:
\[
(\nabla^E_{\b{\Delta}_1}\b{\Delta}_2)(\b{X})
=
(\b{I}_n - \b{X}\b{X}^\top)\,D\b{\Delta}_2(\b{X})[\b{\Delta}_1].
\]
According to Eq. (\ref{equation_projection_onto_tangent_space_Grassmann}), the Levi-Civita connection on the Grassmann manifold can be interpreted as projection of the ambient (Euclidean) directional derivative onto the tangent space of Grassmann manifold. 
This simplicity follows from the Euclidean metric structure.
\end{remark}

\subsubsection{Riemannian Gradient in Grassmann Manifold}


\begin{definition}[Smooth local extension (or local representative) on the ambient space for Grassmann manifold]
Let $f : Gr(n,d) \to \mathbb{R}$ be a smooth function, and let $[\b{X}] \in Gr(n,d)$ be represented by $\b{X} \in St(n,d)$. 
A \textbf{smooth local extension} of $f$ around $\b{X}$ is a smooth function:
\begin{align}
\bar{f} : U \subset \mathbb{R}^{n \times d} \to \mathbb{R},
\end{align}
defined on an open neighborhood $U$ of $\b{X}$, such that:
\begin{align}\label{equation_smooth_extension}
\boxed{
\bar{f}(\b{Y}) = f([\b{Y}]), \qquad \forall \b{Y} \in U \cap St(n,d).
}
\end{align}
In other words, $\bar{f}$ agrees with $f$ on the points of the Grassmann manifold near $[\b{X}]$, but it is defined on an open set of the ambient Euclidean space $\mathbb{R}^{n \times d}$ so that its Euclidean gradient can be computed.
\end{definition}

\begin{proposition}[Riemannian gradient on the Grassmann manifold \cite{edelman1998geometry}]
Let $f : Gr(n,d) \to \mathbb{R}$ be a smooth function on the Grassmann manifold $Gr(n,d)$ and let $\bar{f}: U \subset \mathbb{R}^{n\times d} \to \mathbb{R}$ be a smooth local extension to $\mathbb{R}^{n\times d}$ around a representative
$\b{X} \in St(n,d)$ of $[\b{X}] \in Gr(n,d)$.
Then, the Riemannian gradient on the Grassmann manifold is obtained as:
\begin{equation}\label{equation_Riemannian_gradient_Grassmann_manifold}
\boxed{
\operatorname{grad} f([\b{X}]) = (\b{I}_n - \b{X}\b{X}^\top)\nabla \bar{f}(\b{X}),
}
\end{equation}
where $\b{I}_n$ denotes the $n \times n$ identity matrix. 
\end{proposition}

\begin{proof}
According to Eq. (\ref{equation_definition_Riemannian_gradient_directional_derivative}), the Riemannian gradient here is:
\begin{equation}\label{equation_definition_Riemannian_gradient_directional_derivative_Grassmann}
\begin{aligned}
Df([\b{X}])[\b{\Delta}] = g^E_{[\b{X}]}(\operatorname{grad} &f([\b{X}]), \b{\Delta}), \quad \\
&\forall \b{\Delta} \in T_{[\b{X}]}Gr(n,d).
\end{aligned}
\end{equation}

\textbf{Step 1: Ambient directional derivative.}
\begin{align}\label{equation_DfXDelta_DfbarXDelta_inner_nablafbarX_Delta}
Df([\b{X}])[\b{\Delta}] \overset{(\ref{equation_smooth_extension})}{=} D\bar{f}(\b{X})[\b{\Delta}]
\overset{(\ref{equation_directional_derivative_inner_product})}{=} \langle \nabla \bar{f}(\b{X}), \b{\Delta} \rangle.
\end{align}

\textbf{Step 2: Decomposition of Euclidean gradient.}
\[
\nabla \bar{f}(\b{X})
= (\b{I} - \b{X}\b{X}^\top)\nabla \bar{f}(\b{X})
+ \b{X}\b{X}^\top \nabla \bar{f}(\b{X}).
\]

The inner product of the terms of this equation with $\b{\Delta}$ gives:
\begin{align}
\langle \nabla \bar{f}(\b{X}), \b{\Delta} \rangle
&= \langle (\b{I} - \b{X}\b{X}^\top)\nabla \bar{f}(\b{X}), \b{\Delta} \rangle \nonumber \\
&~~~~ + \langle \b{X}\b{X}^\top \nabla \bar{f}(\b{X}), \b{\Delta} \rangle \nonumber \\
&\overset{(a)}{=} \langle (\b{I} - \b{X}\b{X}^\top)\nabla \bar{f}(\b{X}), \b{\Delta} \rangle, \label{equation_inner_nablafbarX_Delta}
\end{align}
where $(a)$ is because the term $\b{X}\b{X}^\top \nabla \bar{f}(\b{X})$ is normal, so orthogonal to $\b{\Delta}$; therefore, $\langle \b{X}\b{X}^\top \nabla \bar{f}(\b{X}), \b{\Delta} \rangle$ = 0.

Substituting Eq. (\ref{equation_DfXDelta_DfbarXDelta_inner_nablafbarX_Delta}) in Eq. (\ref{equation_inner_nablafbarX_Delta}) gives:
\begin{align*}
Df([\b{X}])[\b{\Delta}] &= \langle (\b{I} - \b{X}\b{X}^\top)\nabla \bar{f}(\b{X}), \b{\Delta} \rangle \\
&\overset{(\ref{equation_g_inner_product})}{=} g^E_{[\b{X}]}\big((\b{I} - \b{X}\b{X}^\top)\nabla \bar{f}(\b{X}), \b{\Delta}\big).
\end{align*}
By comparing this equation with Eq. (\ref{equation_definition_Riemannian_gradient_directional_derivative_Grassmann}), we have:
\begin{align*}
\operatorname{grad} f([X]) = (\b{I} - \b{X}\b{X}^\top)\nabla \bar{f}(\b{X}).
\end{align*}
\end{proof}

\begin{remark}[Interpretation of the Riemannian gradient on Grassmann manifold]
According to Eq. (\ref{equation_Riemannian_gradient_Grassmann_manifold}), the Riemannian gradient on the Grassmann manifold is obtained by removing the normal component $\b{X}\b{X}^\top \nabla \bar{f}(\b{X})$ and keeping only the tangent component.
\end{remark}


\begin{lemma}[The Riemannian gradient satisfies the tangent constraint on the Grassmann manifold]\label{lemma_riemannian_gradient_tangent_constraint_Grassmann}
The Riemannian gradient, stated in Eq. (\ref{equation_Riemannian_gradient_Grassmann_manifold}), satisfies the tangent constraint in Eq. (\ref{equation_tangent_space_Grassmann}). In other words:
\begin{equation}
\boxed{
\big(\operatorname{grad} f([\b{X}])\big)^\top \b{X} = \b{0},
}
\end{equation}
and therefore, it is a tangent vector (matrix) on the Grassmann manifold:
\begin{align}
\boxed{
\operatorname{grad} f([\b{X}]) \in T_{[\b{X}]}Gr(n,d).
}
\end{align}
\end{lemma}

\begin{proof}
We have:
\begin{align*}
\big(\operatorname{grad} f([\b{X}])\big)^\top \b{X}
&\overset{(\ref{equation_Riemannian_gradient_Grassmann_manifold})}{=} \big((\b{I}_n - \b{X}\b{X}^\top)\nabla \bar{f}(\b{X})\big)^\top \b{X} \\
&= \nabla \bar{f}(\b{X})^\top (\b{I}_n - \b{X}\b{X}^\top)^\top \b{X} \\
&\overset{(a)}{=} \nabla \bar{f}(\b{X})^\top (\b{I}_n - \b{X}\b{X}^\top)\b{X} \\
&= \nabla \bar{f}(\b{X})^\top (\b{X} - \b{X}\b{X}^\top \b{X}) \\
&\overset{(b)}{=} \nabla \bar{f}(\b{X})^\top (\b{X} - \b{X}) = \b{0},
\end{align*}
where $(a)$ is because $\b{I}_n - \b{X}\b{X}^\top$ is symmetric and $(b)$ is because $\b{X}^\top \b{X} = \b{I}_d$. 

As $\big(\operatorname{grad} f([\b{X}])\big)^\top \b{X} = \b{0}$, the Riemannian gradient satisfies Eq. (\ref{equation_tangent_space_Grassmann}), so it belongs to the tangent space of Grassmann manifold. 
\end{proof}

\subsubsection{Riemannian Hessian in Grassmann Manifold}

Recall from Section \ref{section_Riemannian_Hessian} that the Riemannian Hessian is the covariant derivative of the Riemannian gradient. In the Grassmann manifold, we represent tangent vectors by their horizontal lifts in the ambient Euclidean space
$\mathbb{R}^{n\times d}$ and use the projection formula in Eq.
(\ref{equation_projection_onto_tangent_space_Grassmann}); also see Eq. (\ref{equation_Riemannian_gradient_Grassmann_manifold}). Therefore, the Riemannian Hessian is obtained by
taking the ambient directional derivative of the Grassmann
gradient and projecting the result back onto the tangent
space.

\begin{lemma}[Ambient directional derivative of the Grassmann
gradient]
Let $f : Gr(n, d) \to \mathbb{R}$ be a smooth
function and let $\bar f : U \subset \mathbb{R}^{n\times d} \to
\mathbb{R}$ be a smooth local extension around a representative
$\b{X} \in St(n,d)$ of $[\b{X}] \in Gr(n,d)$. Recall from Eq. (\ref{equation_Riemannian_gradient_Grassmann_manifold}) that:
\begin{equation*}
\operatorname{grad} f([\b{X}]) = (I - \b{XX}^\top)\nabla \bar f(\b{X}).
\end{equation*}
For a tangent vector $\b{\Delta} \in T_{[\b{X}]}Gr(n,d)$, the classical directional derivative of the Riemannian gradient field in the ambient space $\mathbb{R}^{n\times d}$ is:
\begin{equation}\label{equation_ambient_derivative_grassmann_gradient}
\boxed{
\begin{aligned}
D(\operatorname{grad} f)([\b{X}])[\b{\Delta}]
&= -(\b{\Delta} \b{X}^\top + \b{X}\b{\Delta}^\top)\nabla \bar f(\b{X}) \\
&\quad + (\b{I}_n - \b{XX}^\top)\nabla^2 \bar f(\b{X})[\b{\Delta}],
\end{aligned}
}
\end{equation}
where $\nabla^2 \bar f(\b{X})[\b{\Delta}]$ denotes the Euclidean
Hessian of $\bar f$ acting on the direction $\b{\Delta}$.
\end{lemma}
\begin{proof}
We define:
\[
P(\b{X}) := \b{I} - \b{XX}^\top.
\]
Then Eq. (\ref{equation_Riemannian_gradient_Grassmann_manifold}) can be written as:
\[
\operatorname{grad} f([\b{X}]) = P(\b{X})\nabla \bar f(\b{X}).
\]
We differentiate this matrix-valued mapping in the ambient
space along the direction $\b{\Delta}$. By the product rule, we have:
\begin{equation}\label{equation_product_rule_grassmann_hessian_1}
\begin{aligned}
D(\operatorname{grad} f)([\b{X}])[\b{\Delta}]
=\, &D(P(\b{X}))[\b{\Delta}] \, \nabla \bar f(\b{X})
\\
&+ P(\b{X})\,D(\nabla \bar f(\b{X}))[\b{\Delta}].
\end{aligned}
\end{equation}
We calculate the two terms separately.

Firstly, since $P(\b{X}) = I - \b{XX}^\top$, we have:
\begin{align}
D(P(\b{X}))[\b{\Delta}]
&= -D(\b{XX}^\top)[\b{\Delta}] \nonumber \\
&\overset{(a)}{=} - (\b{\Delta} \b{X}^\top + \b{X}\b{\Delta}^\top), \label{equation_derivative_projection_grassmann}
\end{align}
where $(a)$ is because of product rule for directional derivative of $\b{XX}^\top$. 

Secondly, by definition of the Euclidean Hessian, we have:
\begin{equation}\label{equation_euclidean_hessian_grassmann}
D(\nabla \bar f(\b{X}))[\b{\Delta}] = \nabla^2 \bar f(\b{X})[\b{\Delta}].
\end{equation}

Substituting
Eqs. \eqref{equation_derivative_projection_grassmann} and
\eqref{equation_euclidean_hessian_grassmann} into
Eq. \eqref{equation_product_rule_grassmann_hessian_1}, we obtain:
\begin{align*}
D(\operatorname{grad} f)([\b{X}])[\b{\Delta}]
=\, &-(\b{\Delta} \b{X}^\top + \b{X}\b{\Delta}^\top)\nabla \bar f(\b{X})
\\
&+ (\b{I} - \b{XX}^\top)\nabla^2 \bar f(\b{X})[\b{\Delta}].
\end{align*}
\end{proof}

\begin{proposition}[Riemannian Hessian on the Grassmann
manifold \cite{edelman1998geometry}]\label{proposition_Grassmann_Hessian}
Let $f : Gr(n, d) \to \mathbb{R}$ be a smooth
function and let $\bar f$ be a smooth local extension to
$\mathbb{R}^{n\times d}$. Let $\b{X} \in St(n,d)$ be a representative
of $[\b{X}] \in Gr(n,d)$, and let
$\b{\Delta} \in T_{[\b{X}]}Gr(n,d)$ with $\b{X}^\top \b{\Delta} = \b{0}$.
The Riemannian Hessian operator on the Grassmann manifold is:
\begin{equation}\label{equation_hessian_grassmann_operator_projection}
\boxed{
\begin{aligned}
\operatorname{Hess} f([\b{X}])[\b{\Delta}]
&= \Pi_{[\b{X}]}^{\mathrm{Gr}}\!\big( D(\operatorname{grad} f)([\b{X}])[\b{\Delta}] \big) \\ 
&\overset{(\ref{equation_projection_onto_tangent_space_Grassmann})}{=} (\b{I}_n - \b{X}\b{X}^\top)\,D(\operatorname{grad} f)(\b{X})[\b{\Delta}].
\end{aligned}
}
\end{equation}
Using Eq. \eqref{equation_ambient_derivative_grassmann_gradient}
and the projection formula in Eq. (\ref{equation_projection_onto_tangent_space_Grassmann}), this becomes:
\begin{equation}\label{equation_hessian_grassmann_operator}
\boxed{
\begin{aligned}
\operatorname{Hess} f([\b{X}])[\b{\Delta}]
=\, &(\b{I}_n - \b{XX}^\top)\nabla^2 \bar f(\b{X})[\b{\Delta}] \\
&- \b{\Delta} \big(\b{X}^\top \nabla \bar f(\b{X})\big).
\end{aligned}
}
\end{equation}
\end{proposition}

\begin{proof}
By the definition of the Riemannian Hessian in
Section \ref{section_Riemannian_Hessian}, the Hessian operator is the covariant
derivative of the Riemannian gradient:
\[
\operatorname{Hess} f([\b{X}])[\b{\Delta}]
= \nabla_{\b{\Delta}} \operatorname{grad} f.
\]
Since the Grassmann manifold is endowed with the metric
induced from the ambient Euclidean space and tangent
vectors are represented by horizontal lifts, the Levi-Civita
connection is obtained by orthogonally projecting the ambient
directional derivative onto the tangent space. Therefore:
\[
\operatorname{Hess} f([\b{X}])[\b{\Delta}]
= \Pi_{[\b{X}]}^{\mathrm{Gr}}\!\big( D(\operatorname{grad} f)([\b{X}])[\b{\Delta}] \big),
\]
which proves Eq.
\eqref{equation_hessian_grassmann_operator_projection}.
Also, according to Eq. (\ref{equation_projection_onto_tangent_space_Grassmann}), the above equation is equivalent to:
\begin{align*}
\operatorname{Hess} f([\b{X}])[\b{\Delta}]
&= (\b{I}_n - \b{X}\b{X}^\top)\,D(\operatorname{grad} f)(\b{X})[\b{\Delta}].
\end{align*}

Now, we substitute Eq. (\ref{equation_ambient_derivative_grassmann_gradient}) in this equation:
\begin{align}
\operatorname{Hess} &f([\b{X}])[\b{\Delta}] = \Pi_{[\b{X}]}^{\mathrm{Gr}}\!\big( D(\operatorname{grad} f)([\b{X}])[\b{\Delta}] \big) \nonumber\\
&\overset{(\ref{equation_ambient_derivative_grassmann_gradient})}{=} \Pi_{[\b{X}]}^{\mathrm{Gr}}\Big(\!
-(\b{\Delta} \b{X}^\top + \b{X}\b{\Delta}^\top)\nabla \bar f(\b{X})
\nonumber\\
&~~~~~~~~~~~~~~~ + (\b{I}_n - \b{XX}^\top)\nabla^2 \bar f(\b{X})[\b{\Delta}]
\Big) \nonumber\\
&\overset{(\ref{equation_projection_onto_tangent_space_Grassmann})}{=} (\b{I}_n - \b{XX}^\top)\Big(\!
-(\b{\Delta} \b{X}^\top + \b{X}\b{\Delta}^\top)\nabla \bar f(\b{X}) \nonumber\\
&~~~~~~~ + (\b{I}_n - \b{XX}^\top)\nabla^2 \bar f(\b{X})[\b{\Delta}]
\Big). \label{equation_hessian_grassmann_expand_1}
\end{align}

We simplify the terms one by one.

The first term in Eq. (\ref{equation_hessian_grassmann_expand_1}) is simplified as:
\begin{align*}
&-(\b{I}_n - \b{XX}^\top)\b{\Delta} \b{X}^\top \nabla \bar f(\b{X})
\\
&= -\b{\Delta} \b{X}^\top \nabla \bar f(\b{X}) + \b{XX}^\top \b{\Delta} \b{X}^\top \nabla \bar f(\b{X}) \\
&\overset{(a)}{=} -\b{\Delta} \b{X}^\top \nabla \bar f(\b{X}),
\end{align*}
where $(a)$ is because $\b{\Delta} \in T_{[\b{X}]}Gr(n,d)$ implies $\b{X}^\top \b{\Delta} = \b{0}$. 

The second term in Eq. (\ref{equation_hessian_grassmann_expand_1}) is simplified as:
\begin{align*}
-(\b{I}_n - &\b{XX}^\top)\b{X}\b{\Delta}^\top \nabla \bar f(\b{X}) \\
&= -\b{X}\b{\Delta}^\top \nabla \bar f(\b{X}) + \b{XX}^\top\b{X}\b{\Delta}^\top \nabla \bar f(\b{X}) \\
&\overset{(a)}{=} -\b{X}\b{\Delta}^\top \nabla \bar f(\b{X}) + \b{X}\b{\Delta}^\top \nabla \bar f(\b{X}) = \b{0},
\end{align*}
where $(a)$ is because $\b{X} \in \text{St}(n,d)$ implies $\b{X}^\top\b{X} = \b{I}$. 

The third term in Eq. (\ref{equation_hessian_grassmann_expand_1}) is simplified as:
\begin{align*}
&(\b{I}_n - \b{XX}^\top)(\b{I}_n - \b{XX}^\top)\nabla^2 \bar f(\b{X})[\b{\Delta}]
\\
&\overset{(a)}{=} (\b{I}_n - \b{XX}^\top)\nabla^2 \bar f(\b{X})[\b{\Delta}],
\end{align*}
where $(a)$ is because $\b{I}_n - \b{XX}^\top$ is an idempotent projection:
\[
(\b{I}_n - \b{XX}^\top)^2 = \b{I}_n - \b{XX}^\top.
\]

Substituting these simplifications into
Eq. \eqref{equation_hessian_grassmann_expand_1} gives:
\begin{align*}
\operatorname{Hess} f([\b{X}])[\b{\Delta}]
=
&-\b{\Delta} \b{X}^\top \nabla \bar f(\b{X})
\\
&+ (\b{I}_n - \b{XX}^\top)\nabla^2 \bar f(\b{X})[\b{\Delta}].
\end{align*}
\end{proof}

\begin{remark}[Tangent constraint of the Grassmann Hessian]
The Grassmann Hessian in Eq. \eqref{equation_hessian_grassmann_operator} lies in the tangent space. Thus, it satisfies the condition in Eq. (\ref{equation_tangent_space_Grassmann}):
\begin{align}
\b{X}^\top \operatorname{Hess} f([\b{X}])[\b{\Delta}] = \b{0}.
\end{align}
\end{remark}
\begin{proof}
According to Eq. (\ref{equation_hessian_grassmann_operator}), we have:
\begin{align*}
\b{X}^\top \operatorname{Hess} f([\b{X}])[\b{\Delta}]
\overset{(\ref{equation_hessian_grassmann_operator})}{=}\, &\b{X}^\top (\b{I}_n - \b{XX}^\top)\nabla^2 \bar f(\b{X})[\b{\Delta}] \\
&- \b{X}^\top \b{\Delta} \big(\b{X}^\top \nabla \bar f(\b{X})\big).
\end{align*}

For the first term of Hessian, we have:
\begin{align*}
\b{X}^\top &(\b{I} - \b{XX}^\top)\nabla^2 \bar f(\b{X})[\b{\Delta}] \\
&= \b{X}^\top \nabla^2 \bar f(\b{X})[\b{\Delta}] - \b{X}^\top \b{XX}^\top \nabla^2 \bar f(\b{X})[\b{\Delta}]  \\
&\overset{(a)}{=} \b{X}^\top \nabla^2 \bar f(\b{X})[\b{\Delta}] - \b{X}^\top \nabla^2 \bar f(\b{X})[\b{\Delta}] = \b{0},
\end{align*}
where $(a)$ is because $\b{X}^\top \b{X} = \b{I}$ for $\b{X} \in \text{St}(n,d)$ according to Eq. (\ref{equation_Sitefel_manifold_definition}). 

For the second term of Hessian, we have:
\[
\b{X}^\top \b{\Delta} \big(\b{X}^\top \nabla \bar f(\b{X})\big) = \b{0},
\]
since $\b{X}^\top \b{\Delta} = \b{0}$ for every tangent vector
$\b{\Delta} \in T_{[\b{X}]}Gr(n,d)$, according to Eq. (\ref{equation_tangent_space_Grassmann}).
\end{proof}

\begin{remark}[Interpretation of the correction term in Grassmann Hessian]
The term $- \b{\Delta} \big(\b{X}^\top \nabla \bar f(\b{X})\big)$ in Eq. (\ref{equation_hessian_grassmann_operator}) is the Grassmann analogue of the Christoffel-symbol
correction term in the coordinate formula of the Riemannian
Hessian in Eq. (\ref{equation_Riemannian_Hessian_components}). In Euclidean space, the Hessian
is simply the derivative of the gradient. On the Grassmann
manifold, however, the gradient field must remain tangent to
the manifold, and differentiating the projection
$(\b{I} - \b{XX}^\top)$ introduces the extra curvature/constraint
correction term.
\end{remark}

\begin{remark}[Grassmann Hessian as Levi-Civita connection of Grassmann gradient]
Let $f : Gr(n,d) \to \mathbb{R}$ be a smooth function on the Grassmann manifold. 
According to Eq. (\ref{equation_Riemannian_gradient_Grassmann_manifold}), the Grassmann gradient is:
\[
\operatorname{grad} f([\b{X}])
=
(\b{I}_n - \b{X}\b{X}^\top)\nabla \bar f(\b{X}).
\]
Therefore, by the Levi-Civita connection in Eq. (\ref{equation_levi_civita_connection_grassmann}), the covariant derivative of the Grassmann gradient is:
\[
\nabla_{\b{\Delta}}\,\operatorname{grad} f
=
(\b{I}_n - \b{X}\b{X}^\top)\,D(\operatorname{grad} f)(\b{X})[\b{\Delta}],
\]
which recovers Eq. (\ref{equation_hessian_grassmann_operator_projection}).
Hence, the Grassmann Hessian, introduced in Proposition \ref{proposition_Grassmann_Hessian}, can also be stated as:
\begin{align}
\boxed{
\operatorname{Hess} f([\b{X}])[\b{\Delta}]
=
\nabla_{\b{\Delta}}\,\operatorname{grad} f.
}
\end{align}
Substituting the expression of $\operatorname{grad} f([\b{X}])$ into this equation recovers the explicit formula for the Riemannian Hessian on the Grassmann manifold.
\end{remark}

\subsubsection{Geodesic Equation on Grassmann Manifold}

Recall from Section \ref{section_geodesics} that a smooth curve on a Riemannian manifold is a geodesic if and only if its covariant acceleration vanishes to zero. In the Grassmann manifold, tangent vectors are represented by their horizontal lifts in the ambient Euclidean space $\mathbb{R}^{n\times d}$, and Proposition \ref{proposition_levi_civita_connection_grassmann} showed that the Levi-Civita connection is obtained by projecting the ambient directional derivative onto the tangent space. Using this, we now derive the geodesic equation on the Grassmann manifold.

\begin{lemma}[Differentiation of the tangent constraint along a curve on the Grassmann manifold]
Let $\b{X}(t) \in St(n,d)$ be a smooth representative curve of $[\b{X}(t)] \in Gr(n,d)$. Assume its velocity is tangent to the Grassmann manifold for all $t$, i.e.,
\begin{equation}
\dot{\b{X}}(t) \in T_{[X(t)]}Gr(n,d),
\end{equation}
so by Eq.~(\ref{equation_tangent_space_Grassmann}) we have:
\begin{equation}
\b{X}(t)^\top \dot{\b{X}}(t) = \b{0}.
\end{equation}
Then:
\begin{equation}\label{equation_XTXdd_XdTXd_Grassmann}
\boxed{
\b{X}(t)^\top \ddot{\b{X}}(t) = - \dot{\b{X}}(t)^\top \dot{\b{X}}(t),
}
\end{equation}
where:
\begin{align*}
\dot{\b{X}}(t) = \frac{d\b{X}(t)}{d t}, \quad \ddot{\b{X}}(t) = \frac{d^2 \b{X}(t)}{d t^2}. 
\end{align*}
Note that if we denote $\b{\gamma}(t) \in St(n,d)$ and $[\b{\gamma}(t)] \in Gr(n,d)$ for the curve on the manifold (so as in some texts of literature), the Eq. (\ref{equation_XTXdd_XdTXd_Grassmann}) is denoted as:
\begin{equation}
\b{\gamma}(t)^\top \ddot{\b{\gamma}}(t) = - \dot{\b{\gamma}}(t)^\top \dot{\b{\gamma}}(t).
\end{equation}
\end{lemma}

\begin{proof}
Since $\dot{\b{X}}(t)$ is tangent to the Grassmann manifold, i.e., $\dot{\b{X}}(t) \in T_{[\b{X}]}\text{Gr}(n, d)$, we have according to Eq.~(\ref{equation_tangent_space_Grassmann}):
\[
\b{X}(t)^\top \dot{\b{X}}(t) = \b{0}.
\]
We differentiate both sides with respect to $t$. Using the product rule for matrix-valued functions, we obtain:
\begin{align*}
&\frac{d}{dt}\big(\b{X}(t)^\top \dot{\b{X}}(t)\big)
=
\dot{\b{X}}(t)^\top \dot{\b{X}}(t) + \b{X}(t)^\top \ddot{\b{X}}(t)
=
\b{0} \\
&\implies \b{X}(t)^\top \ddot{\b{X}}(t) = - \dot{\b{X}}(t)^\top \dot{\b{X}}(t).
\end{align*}
\end{proof}

\begin{proposition}[Geodesic equation on the Grassmann manifold \cite{edelman1998geometry}]\label{proposition_geodesic_Grassmann}
Let $[\b{X}(t)] \in Gr(n,d)$ be a smooth curve, and let $\b{X}(t) \in St(n,d)$ be a smooth representative of this curve such that $\dot{\b{X}}(t) \in T_{[\b{X}(t)]}Gr(n,d)$. Then, $[\b{X}(t)]$ is a geodesic on the Grassmann manifold if and only if:
\begin{equation}\label{equation_geodesic_equation_grassmann_projected}
\boxed{
\big(\b{I}_n - \b{X}(t)\b{X}(t)^\top\big)\ddot{\b{X}}(t) = \b{0}.
}
\end{equation}
Equivalently, the geodesic equation can be written as:
\begin{equation}\label{equation_geodesic_equation_grassmann}
\boxed{
\ddot{\b{X}}(t) + \b{X}(t)\big(\dot{\b{X}}(t)^\top \dot{\b{X}}(t)\big) = \b{0},
}
\end{equation}
where:
\begin{align*}
\dot{\b{X}}(t) = \frac{d\b{X}(t)}{d t}, \quad \ddot{\b{X}}(t) = \frac{d^2 \b{X}(t)}{d t^2}. 
\end{align*}

Note that if we denote $\b{\gamma}(t) \in St(n,d)$ and $[\b{\gamma}(t)] \in Gr(n,d)$ for the curve on the manifold (so as in some texts of literature), the Eqs. (\ref{equation_geodesic_equation_grassmann_projected}) and (\ref{equation_geodesic_equation_grassmann}) are denoted as:
\begin{align}
&\big(\b{I}_n - \b{\gamma}(t)\b{\gamma}(t)^\top\big)\ddot{\b{\gamma}}(t) = \b{0}, \\
&\ddot{\b{\gamma}}(t) + \b{\gamma}(t)\big(\dot{\b{\gamma}}(t)^\top \dot{\b{\gamma}}(t)\big) = \b{0}.
\end{align}
\end{proposition}

\begin{proof}
According to Eq. (\ref{equation_geodesic_coordinate_free}), a smooth curve on a Riemannian manifold is a geodesic if and only if its covariant acceleration vanishes:
\begin{equation}
\nabla_{\dot{\b{X}}} \dot{\b{X}} = \b{0}.
\end{equation}
According to Proposition \ref{proposition_levi_civita_connection_grassmann}, the Levi-Civita connection on the Grassmann manifold is:
\begin{equation*}
\nabla_{\b{\Delta}_1}\b{\Delta}_2
=
(\b{I}_n - \b{XX}^\top)\,D\b{\Delta}_2(\b{X})[\b{\Delta}_1].
\end{equation*}
We apply this formula to the vector field along the curve by taking:
\[
\b{\Delta}_1 = \dot{\b{X}}(t), \qquad \b{\Delta}_2 = \dot{\b{X}}(t),
\]
so the above equation becomes:
\begin{align*}
\nabla_{\dot{\b{X}}}\dot{\b{X}}
=
(\b{I}_n - \b{XX}^\top)\,D\dot{\b{X}}(\b{X})[\dot{\b{X}}],
\end{align*}
where we drop $(t)$ from $\dot{\b{X}}(t)$ and $\ddot{\b{X}}(t)$ for simplification in writing expressions.

Since the ambient directional derivative of the velocity field along the curve is the ordinary second derivative, we have:
\[
D\dot{\b{X}}(\b{X})[\dot{\b{X}}] = \ddot{\b{X}}.
\]
Therefore:
\begin{equation*}
\nabla_{\dot{\b{X}}} \dot{\b{X}}
=
(\b{I}_n - \b{XX}^\top)\ddot{\b{X}}.
\end{equation*}
Thus, the geodesic condition $\nabla_{\dot{\b{X}}}\dot{\b{X}}=\b{0}$ is equivalent to:
\[
(\b{I}_n - \b{XX}^\top)\ddot{\b{X}} = \b{0},
\]
which proves Eq.~(\ref{equation_geodesic_equation_grassmann_projected}).

Now, we derive the equivalent explicit form in Eq.~(\ref{equation_geodesic_equation_grassmann}). From Eq.~(\ref{equation_geodesic_equation_grassmann_projected}), we have:
\[
(\b{I}_n - \b{XX}^\top)\ddot{\b{X}} = \b{0} \implies \ddot{\b{X}} = \b{XX}^\top \ddot{\b{X}}.
\]
Hence, $\ddot{\b{X}}$ lies in the normal space of the Grassmann manifold. Therefore, according to Eq. (\ref{equation_normal_space_Grassmann}), it is of the form:
\begin{align}\label{equation_Xdd_XK_in_proof}
\ddot{\b{X}} = \b{XK},
\end{align}
for some matrix $\b{K} \in \mathbb{R}^{d\times d}$.

To determine $\b{K}$, left-multiply both sides by $\b{X}^\top$:
\begin{align}
\b{X}^\top \ddot{\b{X}} = \b{X}^\top \b{X K} &\overset{(a)}{\implies} \b{X}^\top \ddot{\b{X}} = \b{K} \nonumber\\
&\overset{(\ref{equation_XTXdd_XdTXd_Grassmann})}{\implies} \b{K} = - \dot{\b{X}}^\top \dot{\b{X}}, \label{equation_K_XTXdd_in_proof}
\end{align}
where $(a)$ is because $\b{X}^\top \b{X} = \b{I}_d$ for $\b{X} \in St(n,d)$, according to Eq. (\ref{equation_Sitefel_manifold_definition}).

Substituting Eq. (\ref{equation_K_XTXdd_in_proof}) in Eq. (\ref{equation_Xdd_XK_in_proof}) gives:
\[
\ddot{\b{X}} = -\b{X}(\dot{\b{X}}^\top \dot{\b{X}}) \implies \ddot{\b{X}} + \b{X}(\dot{\b{X}}^\top \dot{\b{X}}) = \b{0}.
\]
This proves Eq.~(\ref{equation_geodesic_equation_grassmann}).
\end{proof}

\begin{remark}[Interpretation of the geodesic equation on the Grassmann manifold]
According to Eq. (\ref{equation_projection_onto_tangent_space_Grassmann}), the Eq.~(\ref{equation_geodesic_equation_grassmann_projected}) states that the tangential component of the ambient acceleration is zero. In other words, the acceleration of a geodesic has no component in the tangent space of the Grassmann manifold.

Equation~(\ref{equation_geodesic_equation_grassmann}) shows more explicitly that the acceleration is entirely normal to the manifold, and its normal component is exactly the correction required to preserve the subspace constraint. This is the Grassmann analogue of the geodesic equation on the Stiefel manifold derived in Section \ref{section_geodesic_equation_Stiefel}, but here the formula is simpler because the Levi-Civita connection on the Grassmann manifold is the orthogonal projection of the ambient derivative in Eq.~(\ref{equation_levi_civita_connection_grassmann}).
\end{remark}

\subsubsection{Exponential Map in Grassmann Manifold}

In Section \ref{section_exponential_map_generalizing_addition}, we defined the exponential map on a Riemannian manifold as the point reached at time $t=1$ by the geodesic starting from a point with a prescribed initial tangent vector. For the Grassmann manifold, because geodesics admit a closed-form expression, the exponential map also has a closed-form formula.

We first prove a useful lemma regarding the compact Singular Value Decomposition (SVD) of a tangent vector on the Grassmann manifold.

\begin{lemma}[Singular value decomposition of a tangent vector on Grassmann manifold]\label{lemma_SVD_tangent_vector_Grassmann}
Let $\b{\Delta} \in T_{[\b{X}]}Gr(n,d)$ and suppose its compact singular value decomposition is \cite{ghojogh2019eigenvalue}:
\[
\b{\Delta} = \b{U}\b{\Sigma}\b{V}^\top,
\]
where $\b{U} \in \mathbb{R}^{n\times d}$, $\b{V} \in \mathbb{R}^{d\times d}$ are orthogonal matrices (or, more generally, $\b{U}$ and $\b{V}$ have orthonormal columns in the compact-rank case), and $\b{\Sigma}$ is diagonal with nonnegative singular values. Then:
\[
\b{X}^\top \b{U} = \b{0}.
\]
\end{lemma}

\begin{proof}
Because $\b{\Delta} \in T_{[\b{X}]}Gr(n,d)$, according to Eq. \eqref{equation_tangent_space_Grassmann}, we have:
\[
\b{X}^\top \b{\Delta} = \b{0}.
\]
Substituting the compact SVD of $\b{\Delta}$ gives:
\[
\b{X}^\top \b{U}\b{\Sigma}\b{V}^\top = \b{0}.
\]
Right-multiplying both sides by $\b{V}$ yields:
\[
\b{X}^\top \b{U}\b{\Sigma} \b{V}^\top \b{V} = \b{0} \b{V} \overset{{(a)}}{\implies} \b{X}^\top \b{U}\b{\Sigma} = \b{0},
\]
where $(a)$ is because $\b{V}^\top \b{V} = \b{I}$ because $\b{V}$ is an orthogonal matrix. 

If the compact SVD is taken over the nonzero singular values, then $\b{\Sigma}$ is invertible on its support. Therefore:
\[
\b{X}^\top \b{U} = \b{0}.
\]
Hence, the left singular vectors of a tangent vector are orthogonal to the columns of $\b{X}$.
\end{proof}

\begin{proposition}[Exponential map on the Grassmann manifold \cite{edelman1998geometry}]
Let $[\b{X}] \in Gr(n,d)$ and let $\b{\Delta} \in T_{[\b{X}]}Gr(n,d)$. Suppose the compact singular value decomposition of $\b{\Delta}$ is:
\[
\b{\Delta} = \b{U}\b{\Sigma}\b{V}^\top.
\]
Then, the Grassmann exponential map at $[\b{X}]$ in direction $\b{\Delta}$ is given by:
\begin{equation}\label{equation_exponential_map_Grassmann}
\boxed{
\operatorname{Exp}_{[\b{X}]}(\b{\Delta})
=
\left[
\b{X}\b{V}\cos(\b{\Sigma})\b{V}^\top
+
\b{U}\sin(\b{\Sigma})\b{V}^\top
\right],
}
\end{equation}
where $[\cdot]$ denotes the equivalence class (or orbit), defined in Eq. (\ref{equation_equivalence_X_in_Grassmann}), and $\cos(\b{\Sigma})$ and $\sin(\b{\Sigma})$ are defined by applying cosine and sine elementwise to the diagonal entries of $\b{\Sigma}$.

Equivalently, the geodesic starting from $[\b{X}]$ with initial velocity $\b{\Delta}$ is represented by:
\begin{equation}\label{equation_geodesic_curve_Grassmann_closed_form}
\boxed{
\b{X}(t)
=
\b{X}\b{V}\cos(t\b{\Sigma})\b{V}^\top
+
\b{U}\sin(t\b{\Sigma})\b{V}^\top,
}
\end{equation}
and therefore:
\begin{align}
\boxed{
\operatorname{Exp}_{[\b{X}]}(\b{\Delta}) = [\b{X}(1)],
}
\end{align}
where $[\cdot]$ denotes the equivalence class (or orbit), defined in Eq. (\ref{equation_p_bracket_X}).
\end{proposition}

\begin{proof}
By definition of the exponential map, we must construct the geodesic $\b{X}(t)$ on the Grassmann manifold satisfying:
\[
\b{X}(0)=\b{X},
\qquad
\dot{\b{X}}(0)=\b{\Delta},
\]
and then evaluate it at $t=1$.

We claim that the curve:
\begin{align}\label{equation_geodesic_curve_Grassmann_closed_form_claimed}
\b{X}(t)
=
\b{X}\b{V}\cos(t\b{\Sigma})\b{V}^\top
+
\b{U}\sin(t\b{\Sigma})\b{V}^\top,
\end{align}
is the desired geodesic.

\textbf{Step 1: Verify the initial position.}

Setting $t=0$, and using $\cos(\b{0})=\b{I}_d$ and $\sin(\b{0})=\b{0}$, we obtain:
\begin{align*}
\b{X}(0)
&=
\b{X}\b{V}\cos(\b{0})\b{V}^\top
+
\b{U}\sin(\b{0})\b{V}^\top
\\
&=
\b{X}\b{V}\b{I}_d\b{V}^\top
=
\b{X}\b{V}\b{V}^\top
\overset{(a)}{=}
\b{X}.
\end{align*}
where $\b{V}\b{V}^\top = \b{I}$ because $\b{V}$ is an untruncated orthogonal matrix. 

\textbf{Step 2: Compute the velocity and verify the initial velocity.}

Differentiating $\b{X}(t)$ with respect to $t$ gives:
\begin{align}\label{equation_Xdot_sin_cos_Grassmann_in_proof}
\dot{\b{X}}(t)
=
\b{X}\b{V}\big(\!-\b{\Sigma}\sin(t\b{\Sigma})\big)\b{V}^\top
+
\b{U}\b{\Sigma}\cos(t\b{\Sigma})\b{V}^\top.
\end{align}
Evaluating at $t=0$, we get:
\begin{align*}
\dot{\b{X}}(0)
&=
\b{X}\b{V}\big(\!-\b{\Sigma}\sin(\b{0})\big)\b{V}^\top
+
\b{U}\b{\Sigma}\cos(\b{0})\b{V}^\top
\\
&=
\b{U}\b{\Sigma}\b{V}^\top
=
\b{\Delta}.
\end{align*}

So, the curve has the correct initial conditions.

\textbf{Step 3: Verify that the curve remains on the Grassmann manifold.}


Since a point $[\b{X}(t)]$ is on the Grassmann manifold, we verify that $\b{X}(t)$ is on the Stiefel manifold:
\begin{align*}
&[\b{X}(t)] \in \text{Gr}(n,d) \implies \b{X}(t) \in \text{St}(n,d) \\
&\overset{(\ref{equation_Sitefel_manifold_definition})}{\implies} \b{X}(t)^\top \b{X}(t)=\b{I}_d.
\end{align*}

From the Lemma \ref{lemma_SVD_tangent_vector_Grassmann}, because $\b{\Delta}=\b{U}\b{\Sigma}\b{V}^\top$ is tangent, we have:
\[
\b{X}^\top \b{U} = \b{0}.
\]
Also, since $\b{X}\in St(n,d)$ and $\b{U}$ has orthonormal columns, we have:
\[
\b{X}^\top \b{X} = \b{I}_d,
\qquad
\b{U}^\top \b{U} = \b{I}_d.
\]
Now:
\begin{align*}
\b{X}(t)^\top \b{X}(t)
&=
\left(
\b{V}\cos(t\b{\Sigma})\b{V}^\top \b{X}^\top
+
\b{V}\sin(t\b{\Sigma})\b{U}^\top
\right) \\
&~~~~~~~~\left(
\b{X}\b{V}\cos(t\b{\Sigma})\b{V}^\top
+
\b{U}\sin(t\b{\Sigma})\b{V}^\top
\right) \\
&=
\b{V}\cos(t\b{\Sigma})\b{V}^\top \b{X}^\top \b{X}\b{V}\cos(t\b{\Sigma})\b{V}^\top \\
&\quad +
\b{V}\cos(t\b{\Sigma})\b{V}^\top \b{X}^\top \b{U}\sin(t\b{\Sigma})\b{V}^\top \\
&\quad +
\b{V}\sin(t\b{\Sigma})\b{U}^\top \b{X}\b{V}\cos(t\b{\Sigma})\b{V}^\top \\
&\quad +
\b{V}\sin(t\b{\Sigma})\b{U}^\top \b{U}\sin(t\b{\Sigma})\b{V}^\top.
\end{align*}
Using $\b{X}^\top\b{X}=\b{I}_d$, $\b{X}^\top\b{U}=\b{0}$, $\b{U}^\top\b{X}=\b{0}$, and $\b{U}^\top\b{U}=\b{I}_d$, this simplifies to:
\begin{align*}
\b{X}(t)^\top \b{X}(t)
&=
\b{V}\cos^2(t\b{\Sigma})\b{V}^\top
+
\b{V}\sin^2(t\b{\Sigma})\b{V}^\top \\
&=
\b{V}\big(\cos^2(t\b{\Sigma})+\sin^2(t\b{\Sigma})\big)\b{V}^\top
\\
&\overset{(a)}{=}
\b{V}\b{I}_d\b{V}^\top
= \b{V}\b{V}^\top
\overset{(b)}{=}
\b{I}_d,
\end{align*}
where $(a)$ because $\cos^2(t\b{\Sigma})+\sin^2(t\b{\Sigma}) = \b{I}$, and $(b)$ is because $\b{V}\b{V}^\top = \b{I}$ because $\b{V}$ is an untruncated orthogonal matrix. 

We proved $\b{X}(t)^\top \b{X}(t) = \b{I}$, so according to Eq. (\ref{equation_Sitefel_manifold_definition}),  
we have $\b{X}(t)\in St(n,d)$ for all $t$, so $[\b{X}(t)]\in Gr(n,d)$ for all $t$.

\textbf{Step 4: Verify the tangent condition along the curve.}

We next verify that the velocity remains tangent to the Grassmann manifold:
\[
\b{X}(t)^\top \dot{\b{X}}(t)=\b{0}.
\]
Using the expressions for $\b{X}(t)$ and $\dot{\b{X}}(t)$---i.e., Eqs. (\ref{equation_geodesic_curve_Grassmann_closed_form_claimed}) and (\ref{equation_Xdot_sin_cos_Grassmann_in_proof})---and again the relations
$\b{X}^\top\b{U}=\b{0}$, $\b{U}^\top\b{X}=\b{0}$, $\b{X}^\top\b{X}=\b{I}_d$, and $\b{U}^\top\b{U}=\b{I}_d$, we obtain:
\begin{align*}
&\b{X}(t)^\top \dot{\b{X}}(t) \\
&= \left( \b{X}\b{V}\cos(t\b{\Sigma})\b{V}^\top + \b{U}\sin(t\b{\Sigma})\b{V}^\top \right)^\top  \\
&~~~~~ \left( \b{X}\b{V}\big(\!-\b{\Sigma}\sin(t\b{\Sigma})\big)\b{V}^\top + \b{U}\b{\Sigma}\cos(t\b{\Sigma})\b{V}^\top \right) \\
&= \left( \b{V}\cos(t\b{\Sigma})\b{V}^\top\b{X}^\top + \b{V}\sin(t\b{\Sigma})\b{U}^\top \right)  \\
&~~~~~ \left( \b{X}\b{V}\big(\!-\b{\Sigma}\sin(t\b{\Sigma})\big)\b{V}^\top + \b{U}\b{\Sigma}\cos(t\b{\Sigma})\b{V}^\top \right) \\
&= \b{V}\cos(t\b{\Sigma})\underbrace{\b{V}^\top\underbrace{\b{X}^\top \b{X}}_{=\b{I}}\b{V}}_{=\b{I}}\big(\!-\b{\Sigma}\sin(t\b{\Sigma})\big)\b{V}^\top \\
&~~~~ + \b{V}\cos(t\b{\Sigma})\b{V}^\top\underbrace{\b{X}^\top \b{U}}_{=\b{0}}\b{\Sigma}\cos(t\b{\Sigma})\b{V}^\top \\
&~~~~ + \b{V}\sin(t\b{\Sigma})\underbrace{\b{U}^\top \b{X}}_{=\b{0}}\b{V}\big(\!-\b{\Sigma}\sin(t\b{\Sigma})\big)\b{V}^\top \\
&~~~~ + \b{V}\sin(t\b{\Sigma})\underbrace{\b{U}^\top \b{U}}_{=\b{I}}\b{\Sigma}\cos(t\b{\Sigma})\b{V}^\top \\
&=
\b{V}\cos(t\b{\Sigma})\big(-\b{\Sigma}\sin(t\b{\Sigma})\big)\b{V}^\top
\\
&+
\b{V}\sin(t\b{\Sigma})\b{\Sigma}\cos(t\b{\Sigma})\b{V}^\top \\
&= \b{0}.
\end{align*}
We just proved $\b{X}(t)^\top \dot{\b{X}}(t) = \b{0}$.
Therefore, according to Eq. (\ref{equation_tangent_space_Grassmann}), we have $\dot{\b{X}}(t)\in T_{[\b{X}(t)]}Gr(n,d)$ for all $t$.

\textbf{Step 5: Compute the acceleration.}

Consider Eq. (\ref{equation_Xdot_sin_cos_Grassmann_in_proof}) for $\dot{\b{X}}(t)$:
\[
\dot{\b{X}}(t)
=
-\b{X}\b{V}\b{\Sigma}\sin(t\b{\Sigma})\b{V}^\top
+
\b{U}\b{\Sigma}\cos(t\b{\Sigma})\b{V}^\top.
\]
Recall that $\b{X}$, $\b{U}$, $\b{V}$, and $\b{\Sigma}$ are constant matrices with respect to $t$. Only the matrix functions $\sin(t\b{\Sigma})$ and $\cos(t\b{\Sigma})$ depend on $t$.

Because $\b{\Sigma}$ is diagonal, say:
\[
\b{\Sigma}=\operatorname{diag}(\sigma_1,\dots,\sigma_d),
\]
the sine and cosine of $t\b{\Sigma}$ are defined elementwise:
\[
\sin(t\b{\Sigma})
=
\operatorname{diag}(\sin(t\sigma_1),\dots,\sin(t\sigma_d)),
\]
\[
\cos(t\b{\Sigma})
=
\operatorname{diag}(\cos(t\sigma_1),\dots,\cos(t\sigma_d)).
\]
Therefore, their derivatives are also taken elementwise:
\[
\frac{d}{dt}\sin(t\b{\Sigma}) = \b{\Sigma}\cos(t\b{\Sigma}),
\quad
\frac{d}{dt}\cos(t\b{\Sigma}) = -\b{\Sigma}\sin(t\b{\Sigma}).
\]

Now, we differentiate $\dot{\b{X}}(t)$ term by term. For the first term:
\begin{align*}
\frac{d}{dt}\Big(\!-\b{X}\b{V}\b{\Sigma}\sin(t\b{\Sigma})\b{V}^\top&\Big)
=
-\b{X}\b{V}\b{\Sigma}\frac{d}{dt}\big(\sin(t\b{\Sigma})\big)\b{V}^\top
\\
&=
-\b{X}\b{V}\b{\Sigma}\big(\b{\Sigma}\cos(t\b{\Sigma})\big)\b{V}^\top \\
&= -\b{X}\b{V}\b{\Sigma}^2\cos(t\b{\Sigma})\b{V}^\top.
\end{align*}

For the second term:
\begin{align*}
\frac{d}{dt}\Big(\b{U}\b{\Sigma}\cos(t\b{\Sigma})\b{V}^\top\Big)
&=
\b{U}\b{\Sigma}\frac{d}{dt}\big(\cos(t\b{\Sigma})\big)\b{V}^\top
\\
&=
\b{U}\b{\Sigma}\big(-\b{\Sigma}\sin(t\b{\Sigma})\big)\b{V}^\top \\
&= -\b{U}\b{\Sigma}^2\sin(t\b{\Sigma})\b{V}^\top.
\end{align*}

Combining the two derivatives gives:
\begin{align}\label{equation_X_ddot_XVSigma2cosVT_USigma2SinVT_in_proof}
\ddot{\b{X}}(t)
=
-\b{X}\b{V}\b{\Sigma}^2\cos(t\b{\Sigma})\b{V}^\top
-
\b{U}\b{\Sigma}^2\sin(t\b{\Sigma})\b{V}^\top.
\end{align}

We now factor this expression. We start from the definition of $\b{X}(t)$, i.e., Eq. (\ref{equation_geodesic_curve_Grassmann_closed_form_claimed}):
\[
\b{X}(t)
=
\b{X}\b{V}\cos(t\b{\Sigma})\b{V}^\top
+
\b{U}\sin(t\b{\Sigma})\b{V}^\top.
\]
We multiply $\b{X}(t)$ on the right by $\b{V}\b{\Sigma}^2\b{V}^\top$:
\begin{align*}
&\b{X}(t)\,\b{V}\b{\Sigma}^2\b{V}^\top \\
&=
\left(
\b{X}\b{V}\cos(t\b{\Sigma})\b{V}^\top
+
\b{U}\sin(t\b{\Sigma})\b{V}^\top
\right)\b{V}\b{\Sigma}^2\b{V}^\top \\
&=
\b{X}\b{V}\cos(t\b{\Sigma})\underbrace{\b{V}^\top\b{V}}_{=\b{I}}\b{\Sigma}^2\b{V}^\top
+
\b{U}\sin(t\b{\Sigma})\underbrace{\b{V}^\top\b{V}}_{=\b{I}}\b{\Sigma}^2\b{V}^\top \\
&=
\b{X}\b{V}\cos(t\b{\Sigma})\b{\Sigma}^2\b{V}^\top
+
\b{U}\sin(t\b{\Sigma})\b{\Sigma}^2\b{V}^\top.
\end{align*}
Since $\b{\Sigma}^2$, $\sin(t\b{\Sigma})$, and $\cos(t\b{\Sigma})$ are all diagonal matrices, they commute. Therefore:
\[
\cos(t\b{\Sigma})\b{\Sigma}^2 = \b{\Sigma}^2\cos(t\b{\Sigma}),
\quad
\sin(t\b{\Sigma})\b{\Sigma}^2 = \b{\Sigma}^2\sin(t\b{\Sigma}),
\]
and thus:
\begin{align}
&\b{X}(t)\,\b{V}\b{\Sigma}^2\b{V}^\top
\nonumber \\
&~~~~~~~~~~ =
\b{X}\b{V}\b{\Sigma}^2\cos(t\b{\Sigma})\b{V}^\top
+
\b{U}\b{\Sigma}^2\sin(t\b{\Sigma})\b{V}^\top. \label{equation_XVsigma2VT_in_proof}
\end{align}

Comparing Eqs. (\ref{equation_X_ddot_XVSigma2cosVT_USigma2SinVT_in_proof}) and (\ref{equation_XVsigma2VT_in_proof}) gives:
\begin{align}\label{equation_Xddot_XVsigma2VT}
\ddot{\b{X}}(t)
=
-\b{X}(t)\,\b{V}\b{\Sigma}^2\b{V}^\top.
\end{align}

\textbf{Step 6: Compute $\dot{\b{X}}(t)^\top\dot{\b{X}}(t)$.}

Since $\b{\Sigma}$, $\sin(t\b{\Sigma})$, and $\cos(t\b{\Sigma})$ are all diagonal matrices, they commute. Therefore:
\begin{align}\label{equation_cos_sigma_sin_sigma_commute}
\cos(t\b{\Sigma})\b{\Sigma} = \b{\Sigma}\cos(t\b{\Sigma}),
\quad
\sin(t\b{\Sigma})\b{\Sigma} = \b{\Sigma}\sin(t\b{\Sigma}).
\end{align}

From the expression of $\dot{\b{X}}(t)$, i.e., Eq. (\ref{equation_Xdot_sin_cos_Grassmann_in_proof}), we have:
\begin{align}
&\dot{\b{X}}(t)^\top \dot{\b{X}}(t) \nonumber\\
&\overset{(\ref{equation_Xdot_sin_cos_Grassmann_in_proof})}{=}
\Big( -\b{X}\b{V}\b{\Sigma}\sin(t\b{\Sigma})\b{V}^\top + \b{U}\b{\Sigma}\cos(t\b{\Sigma})\b{V}^\top \Big)^\top \nonumber\\
&~~~~~~~~\Big( -\b{X}\b{V}\b{\Sigma}\sin(t\b{\Sigma})\b{V}^\top + \b{U}\b{\Sigma}\cos(t\b{\Sigma})\b{V}^\top \Big) \nonumber\\
&= \Big( -\b{V}\sin(t\b{\Sigma}) \b{\Sigma} \b{V}^\top \b{X}^\top + \b{V}\cos(t\b{\Sigma})\b{\Sigma}\b{U}^\top \Big) \nonumber\\
&~~~~~~~~\Big( -\b{X}\b{V}\b{\Sigma}\sin(t\b{\Sigma})\b{V}^\top + \b{U}\b{\Sigma}\cos(t\b{\Sigma})\b{V}^\top \Big) \nonumber\\
&= \b{V}\sin(t\b{\Sigma}) \b{\Sigma} \underbrace{\b{V}^\top \underbrace{\b{X}^\top 
\b{X}}_{=\b{I}} \b{V}}_{=\b{I}} \b{\Sigma}\sin(t\b{\Sigma})\b{V}^\top \nonumber\\
&~~~~ - \b{V}\sin(t\b{\Sigma}) \b{\Sigma} \b{V}^\top \underbrace{\b{X}^\top
\b{U}}_{=\b{0}} \b{\Sigma}\cos(t\b{\Sigma})\b{V}^\top \nonumber\\
&~~~~ - \b{V}\cos(t\b{\Sigma})\b{\Sigma}\underbrace{\b{U}^\top 
\b{X}}_{=\b{0}} \b{V}\b{\Sigma}\sin(t\b{\Sigma})\b{V}^\top \nonumber\\
&~~~~ + \b{V}\cos(t\b{\Sigma})\b{\Sigma}\underbrace{\b{U}^\top 
\b{U}}_{=\b{I}}\b{\Sigma}\cos(t\b{\Sigma})\b{V}^\top \nonumber\\
&\overset{(\ref{equation_cos_sigma_sin_sigma_commute})}{=} 
\b{V} \b{\Sigma} \b{\Sigma} \sin(t\b{\Sigma}) \sin(t\b{\Sigma})\b{V}^\top \nonumber\\
&~~~~ + \b{V}\b{\Sigma}\b{\Sigma}\cos(t\b{\Sigma})\cos(t\b{\Sigma})\b{V}^\top \nonumber\\
&=
\b{V}\b{\Sigma}^2\sin^2(t\b{\Sigma})\b{V}^\top
+
\b{V}\b{\Sigma}^2\cos^2(t\b{\Sigma})\b{V}^\top \nonumber\\
&=
\b{V}\b{\Sigma}^2\big(\sin^2(t\b{\Sigma})+\cos^2(t\b{\Sigma})\big)\b{V}^\top \overset{(a)}{=}
\b{V}\b{\Sigma}^2\b{V}^\top, \label{equation_Xdot_Xdot_V_sigma2_VT_in_proof}
\end{align}
where $(a)$ is because $\sin^2(t\b{\Sigma})+\cos^2(t\b{\Sigma}) = \b{I}$.

\textbf{Step 7: Verify the geodesic equation.}

According to Eq. (\ref{equation_geodesic_equation_grassmann}), the geodesic equation on the Grassmann manifold is:
\[
\ddot{\b{X}}(t) + \b{X}(t)\big(\dot{\b{X}}(t)^\top \dot{\b{X}}(t)\big)=\b{0}.
\]
Substituting Eqs. (\ref{equation_Xddot_XVsigma2VT}) and (\ref{equation_Xdot_Xdot_V_sigma2_VT_in_proof}) in this equation of geodesic gives:
\begin{align*}
\ddot{\b{X}}(t) + &\b{X}(t)\big(\dot{\b{X}}(t)^\top \dot{\b{X}}(t)\big) \\
&= -\b{X}(t)\,\b{V}\b{\Sigma}^2\b{V}^\top + \b{X}(t) \b{V}\b{\Sigma}^2\b{V}^\top = \b{0}.
\end{align*}
Hence, $\b{X}(t)$ satisfies the geodesic equation on Grassmann mnaifold, so it is indeed a geodesic on the Grassmann manifold.

\textbf{Step 8: Obtain the exponential map.}

By definition, the exponential map is the endpoint of this geodesic at time $t=1$. Therefore:
\begin{align*}
\operatorname{Exp}_{[\b{X}]}(\b{\Delta})
&=
[\b{X}(1)]
\\
&=
\left[
\b{X}\b{V}\cos(\b{\Sigma})\b{V}^\top
+
\b{U}\sin(\b{\Sigma})\b{V}^\top
\right].
\end{align*}
This proves Eq. \eqref{equation_exponential_map_Grassmann}.
\end{proof}


\subsubsection{QR-Based Retraction Map in Grassmann Manifold}\label{section_QR_retraction_Grassmann}

As discussed in Section \ref{section_retraction}, a retraction map is a computationally efficient alternative to the exponential map for mapping a tangent vector back to the manifold. 
One of the versions of retraction map on the Grassmann manifold is based on QR decomposition. We introduce it in the following. 

\begin{proposition}[QR-based retraction on the Grassmann manifold \cite{edelman1998geometry}]
Let \([\b{X}] \in Gr(n,d)\), where \(\b{X} \in St(n,d)\) is a representative, and let
\(\b{\Delta} \in T_{[\b{X}]}Gr(n,d)\) satisfy:
\[
\b{X}^\top \b{\Delta} = \b{0}.
\]
Then, the mapping:
\begin{equation}\label{equation_retraction_Grassmann_manifold_qr}
\boxed{
\operatorname{Ret}_{[\b{X}]}(\b{\Delta})
=
[\operatorname{qf}(\b{X}+\b{\Delta})],
}
\end{equation}
is a valid retraction map on the Grassmann manifold, where \(\operatorname{qf}(\cdot)\) denotes the \(Q\)-factor of the QR decomposition (see Lemma \ref{lemma_qr_decomp}) and $[\cdot]$ denotes the equivalence class (or orbit), defined in Eq. (\ref{equation_equivalence_X_in_Grassmann}).
\end{proposition}

\begin{proof}
To prove that Eq. \eqref{equation_retraction_Grassmann_manifold_qr} is a valid retraction, according to Definition \ref{definition_retraction}, we need to verify the following two properties:

\begin{enumerate}
\item Identity:
\[
\operatorname{Ret}_{[\b{X}]}(\b{0}) = [\b{X}].
\]

\item Local rigidity:
\[
\left.\frac{d}{dt}\operatorname{Ret}_{[\b{X}]}(t\b{\Delta})\right|_{t=0}
=
\b{\Delta}.
\]
\end{enumerate}

\textbf{Step 1: Identity.}

For \(\b{\Delta}=\b{0}\), Eq. \eqref{equation_retraction_Grassmann_manifold_qr} gives:
\[
\operatorname{Ret}_{[\b{X}]}(\b{0})
=
[\operatorname{qf}(\b{X})].
\]
Since \(\b{X} \in St(n,d)\), its columns are already orthonormal, so its QR decomposition is
\[
\b{X} = \b{X}\b{I}_d.
\]
Hence:
\[
\operatorname{qf}(\b{X}) = \b{X}.
\]
Therefore:
\[
\operatorname{Ret}_{[\b{X}]}(\b{0}) = [\b{X}].
\]
So the identity property holds.

\textbf{Step 2: Local rigidity.}

We define:
\[
\b{Y}(t) := \operatorname{qf}(\b{X}+t\b{\Delta}).
\]
Then, by definition of QR decomposition, there exists an upper triangular matrix \(\b{R}(t)\) such that:
\begin{equation}\label{equation_qr_curve_Grassmann_retraction}
\b{X}+t\b{\Delta} = \b{Y}(t)\b{R}(t).
\end{equation}
Since \(\b{Y}(t)\) is the \(Q\)-factor, we have:
\[
\b{Y}(t)^\top \b{Y}(t)=\b{I}_d.
\]

Also, consider the claimed retraction map in Eq. (\ref{equation_retraction_Grassmann_manifold_qr}):
\begin{align*}
&\operatorname{Ret}_{[\b{X}]}(\b{\Delta})
=
[\operatorname{qf}(\b{X}+\b{\Delta})] \\
&\implies 
\operatorname{Ret}_{[\b{X}]}(t\b{\Delta})
=
[\operatorname{qf}(\b{X}+t\b{\Delta})]
=
[\b{Y}(t)].
\end{align*}
Thus:
\begin{align}\label{equation_RetXtDelta_bracket_Yt_in_proof}
\operatorname{Ret}_{[\b{X}]}(t\b{\Delta}) = [\b{Y}(t)].
\end{align}


At \(t=0\), Eq. \eqref{equation_qr_curve_Grassmann_retraction} becomes:
\[
\b{X} = \b{Y}(0)\b{R}(0).
\]
Since \(\b{X}\) already has orthonormal columns, we may choose:
\[
\b{Y}(0)=\b{X},
\quad
\b{R}(0)=\b{I}_d.
\]

Now, we differentiate Eq. \eqref{equation_qr_curve_Grassmann_retraction} with respect to \(t\) at \(t=0\):
\[
\b{\Delta}
=
\dot{\b{Y}}(0)\b{R}(0)+\b{Y}(0)\dot{\b{R}}(0)
=
\dot{\b{Y}}(0)+\b{X}\dot{\b{R}}(0).
\]
Hence:
\begin{equation}\label{equation_Ydot0_qr_Grassmann_retraction}
\dot{\b{Y}}(0)
=
\b{\Delta}-\b{X}\dot{\b{R}}(0).
\end{equation}

We next show that \(\dot{\b{R}}(0)=\b{0}\). Since \(\b{Y}(t)^\top \b{Y}(t)=\b{I}_d\), differentiating at \(t=0\) gives:
\[
\b{X}^\top \dot{\b{Y}}(0)+\dot{\b{Y}}(0)^\top \b{X}=\b{0}.
\]
Substituting Eq. \eqref{equation_Ydot0_qr_Grassmann_retraction} into this equation gives:
\begin{align}
&\b{X}^\top (\b{\Delta}-\b{X}\dot{\b{R}}(0))
+
(\b{\Delta}-\b{X}\dot{\b{R}}(0))^\top \b{X}
=
\b{0}
\nonumber\\
&\implies 
(\b{X}^\top \b{\Delta}+\b{\Delta}^\top \b{X})
\nonumber\\
&~~~~~~~~~~~~~~~~ -
(\b{X}^\top \b{X} \dot{\b{R}}(0)+ \b{X}^\top \b{X}\dot{\b{R}}(0)^\top)
=
\b{0} \nonumber \\
&\overset{(a)}{\implies} 
(\b{X}^\top \b{\Delta}+\b{\Delta}^\top \b{X})
-
(\dot{\b{R}}(0)+\dot{\b{R}}(0)^\top)
=
\b{0},
\label{equation_Rdot_symmetry_Grassmann_retraction}
\end{align}
where $(a)$ is because $\b{X} \in \text{St}(n,d)$ so $\b{X}^\top \b{X} = \b{I}$ according to Eq. (\ref{equation_Sitefel_manifold_definition}). 

But \(\b{\Delta} \in T_{[\b{X}]}Gr(n,d)\), so according to Eq. (\ref{equation_tangent_space_Grassmann}), we have:
\[
\b{X}^\top \b{\Delta}=\b{0}, \quad \b{\Delta}^\top \b{X}=\b{0}.
\]
Therefore, Eq. \eqref{equation_Rdot_symmetry_Grassmann_retraction} reduces to:
\[
\dot{\b{R}}(0)+\dot{\b{R}}(0)^\top=\b{0}.
\]
Since \(\b{R}(t)\) is upper triangular for all \(t\), its derivative \(\dot{\b{R}}(0)\) is also upper triangular. The only upper triangular matrix that is also skew-symmetric is the zero matrix. Hence:
\[
\dot{\b{R}}(0)=\b{0}.
\]
Substituting this into Eq. \eqref{equation_Ydot0_qr_Grassmann_retraction} gives:
\[
\dot{\b{Y}}(0)=\b{\Delta}.
\]

Now, recall Eq. (\ref{equation_RetXtDelta_bracket_Yt_in_proof}):
\[
\operatorname{Ret}_{[\b{X}]}(t\b{\Delta})=[\b{Y}(t)].
\]
As we have \(\dot{\b{Y}}(0)=\b{\Delta}\) and \(\b{X}^\top \b{\Delta}=0\), the derivative \(\dot{\b{Y}}(0)\) satisfies the tangent condition in Eq. (\ref{equation_tangent_space_Grassmann}). Thus, it is already the horizontal representative of the tangent vector of the quotient curve \([\b{Y}(t)]\). Therefore:
\[
\left.\frac{d}{dt}[\b{Y}(t)]\right|_{t=0}
=
\b{\Delta}.
\]
Hence:
\[
\left.\frac{d}{dt}\operatorname{Ret}_{[\b{X}]}(t\b{\Delta})\right|_{t=0}
=
\left.\frac{d}{dt}[\b{Y}(t)]\right|_{t=0}
=
\b{\Delta}.
\]
So the local rigidity property holds.

Since both identity and local rigidity are satisfied, Eq. \eqref{equation_retraction_Grassmann_manifold_qr} is a valid retraction map on the Grassmann manifold.
\end{proof}

\subsubsection{Polar-Style Retraction Map in Grassmann Manifold}\label{section_polar_retraction_Grassmann}

There is also a second version of retraction, namely polar-style retraction, on the Grassmann manifold. It is introduced in the following. 
For the Grassmann manifold \(Gr(n,d)\), a convenient retraction is obtained by orthonormalizing the columns of \(\b{X}+\b{\Delta}\). Because a point on the Grassmann manifold is an equivalence class \([\b{X}]\), the retraction should return a subspace rather than a specific orthonormal basis.

We first prove a useful lemma showing that the formula is independent of the chosen representative of the equivalence class.

\begin{lemma}[Well-definedness of orthonormalization on Grassmann manifold]
Let \([\b{X}] \in Gr(n,d)\), and let \(\b{\Delta} \in T_{[\b{X}]}Gr(n,d)\). We define:
\[
\b{Y}
:=
(\b{X}+\b{\Delta})
(\b{I}_d + \b{\Delta}^\top \b{\Delta})^{-1/2}.
\]
Then, \(\b{Y} \in St(n,d)\). Moreover, if \(\widetilde{\b{X}} = \b{X}\b{Q}\) is another representative of \([\b{X}]\), with \(\b{Q} \in O(d)\), and \(\widetilde{\b{\Delta}} = \b{\Delta}\b{Q}\) is the corresponding tangent representative, then:
\begin{align*}
\big[
(\widetilde{\b{X}}+\widetilde{\b{\Delta}})
(\b{I}_d &+ \widetilde{\b{\Delta}}^\top \widetilde{\b{\Delta}})^{-1/2}
\big]
\\
&=
\big[
(\b{X}+\b{\Delta})
(\b{I}_d + \b{\Delta}^\top \b{\Delta})^{-1/2}
\big],
\end{align*}
where $[\cdot]$ denotes the equivalence class (or orbit), defined in Eq. (\ref{equation_equivalence_X_in_Grassmann}).

Hence, the resulting point on \(Gr(n,d)\) depends only on the equivalence class \([\b{X}]\) and tangent vector \(\b{\Delta}\), not on the chosen basis \(\b{X}\).
This characteristic is called well-definedness of orthonormalization on Grassmann manifold. 
\end{lemma}

\begin{proof}
First, because \(\b{\Delta} \in T_{[\b{X}]}Gr(n,d)\), according to Eq. (\ref{equation_tangent_space_Grassmann}), we have:
\[
\b{X}^\top \b{\Delta} = \b{0}.
\]
Also, since \(\b{X} \in St(n,d)\), we have according to Eq. (\ref{equation_Sitefel_manifold_definition}):
\[
\b{X}^\top \b{X} = \b{I}_d.
\]
Therefore:
\begin{align}
(\b{X}+\b{\Delta})^\top &(\b{X}+\b{\Delta})
\nonumber\\
&=
\b{X}^\top \b{X}
+
\b{X}^\top \b{\Delta}
+
\b{\Delta}^\top \b{X}
+
\b{\Delta}^\top \b{\Delta}
\nonumber\\
&=
\b{I}_d + \b{\Delta}^\top \b{\Delta}.
\label{equation_XplusDelta_gram_Grassmann_retraction}
\end{align}
Now, we define:
\[
\b{Y}
:=
(\b{X}+\b{\Delta})
(\b{I}_d + \b{\Delta}^\top \b{\Delta})^{-1/2}.
\]
Then:
\begin{align}
&\b{Y}^\top \b{Y}
\nonumber\\
&=
(\b{I}_d + \b{\Delta}^\top \b{\Delta})^{-1/2}
(\b{X}+\b{\Delta})^\top (\b{X}+\b{\Delta}) \nonumber\\
&\qquad\qquad\qquad\qquad\qquad\qquad\qquad 
(\b{I}_d + \b{\Delta}^\top \b{\Delta})^{-1/2}
\nonumber\\
&\overset{\eqref{equation_XplusDelta_gram_Grassmann_retraction}}{=}
(\b{I}_d + \b{\Delta}^\top \b{\Delta})^{-1/2}
(\b{I}_d + \b{\Delta}^\top \b{\Delta})
(\b{I}_d + \b{\Delta}^\top \b{\Delta})^{-1/2}
\nonumber\\
&=
\b{I}_d.
\end{align}
Hence, according to Eq. (\ref{equation_Sitefel_manifold_definition}), we have \(\b{Y} \in St(n,d)\), so \([\b{Y}] \in Gr(n,d)\).

Now, we consider another representative \(\widetilde{\b{X}} = \b{X}\b{Q}\), where \(\b{Q} \in O(d)\). The corresponding tangent representative is \(\widetilde{\b{\Delta}} = \b{\Delta}\b{Q}\). Then:
\begin{align}\label{equation_DeltatildeTDeltatilde_QTDeltaTDeltaQ}
\widetilde{\b{\Delta}}^\top \widetilde{\b{\Delta}}
=
(\b{\Delta}\b{Q})^\top (\b{\Delta}\b{Q})
=
\b{Q}^\top \b{\Delta}^\top \b{\Delta}\,\b{Q}.
\end{align}
Hence:
\begin{align}
\b{I}_d + \widetilde{\b{\Delta}}^\top \widetilde{\b{\Delta}}
&= \b{Q}^\top \b{Q} + \widetilde{\b{\Delta}}^\top \widetilde{\b{\Delta}} \nonumber\\
&\overset{(\ref{equation_DeltatildeTDeltatilde_QTDeltaTDeltaQ})}{=}
\b{Q}^\top (\b{I}_d + \b{\Delta}^\top \b{\Delta}) \b{Q}. \label{equation_I_plus_DeltatildeTDeltatilde}
\end{align}
It is known that orthogonal similarity preserves the positive definite square root:
\begin{align}\label{equation_QTAQinverse_in_proof}
(\b{Q}^\top \b{A}\b{Q})^{-1/2}
=
\b{Q}^\top \b{A}^{-1/2}\b{Q},
\end{align}
for every symmetric positive definite matrix \(\b{A}\). Therefore:
\begin{align}
(\widetilde{\b{X}}+\widetilde{\b{\Delta}})&
(\b{I}_d + \widetilde{\b{\Delta}}^\top \widetilde{\b{\Delta}})^{-1/2}
\nonumber\\
&\overset{(\ref{equation_I_plus_DeltatildeTDeltatilde})}{=}
(\b{X}\b{Q}+\b{\Delta}\b{Q})
\big(\b{Q}^\top (\b{I}_d+\b{\Delta}^\top\b{\Delta})\b{Q}\big)^{-1/2}
\nonumber\\
&\overset{(\ref{equation_QTAQinverse_in_proof})}{=}
(\b{X}+\b{\Delta})\b{Q}
\Big(
\b{Q}^\top (\b{I}_d+\b{\Delta}^\top\b{\Delta})^{-1/2}\b{Q}
\Big)
\nonumber\\
&\overset{(a)}{=}
(\b{X}+\b{\Delta})
(\b{I}_d+\b{\Delta}^\top\b{\Delta})^{-1/2}
\b{Q},
\end{align}
where $(a)$ is because $\b{Q} \b{Q}^\top = \b{I}$ since $\b{Q}$ is an untruncated orthogonal matrix. 

Right multiplication by \(\b{Q} \in O(d)\) does not change the equivalence class in \(Gr(n,d)\). Thus:
\begin{align*}
\big[
(\widetilde{\b{X}}+\widetilde{\b{\Delta}})
&(\b{I}_d + \widetilde{\b{\Delta}}^\top \widetilde{\b{\Delta}})^{-1/2}
\big]
\\
&=
\big[
(\b{X}+\b{\Delta})
(\b{I}_d + \b{\Delta}^\top \b{\Delta})^{-1/2}
\big].
\end{align*}
This proves well-definedness.
\end{proof}

\begin{proposition}[Polar-style retraction on the Grassmann manifold \cite{edelman1998geometry}]
Let \([\b{X}] \in Gr(n,d)\), where \(\b{X} \in St(n,d)\) is a representative, and let
\(\b{\Delta} \in T_{[\b{X}]}Gr(n,d)\) satisfy \(\b{X}^\top \b{\Delta} = \b{0}\).
The following expression:
\begin{equation}\label{equation_retraction_Grassmann_manifold}
\boxed{
\operatorname{Ret}_{[\b{X}]}(\b{\Delta})
=
\left[
(\b{X}+\b{\Delta})
(\b{I}_d + \b{\Delta}^\top \b{\Delta})^{-1/2}
\right],
}
\end{equation}
is a valid retraction map on the Grassmann manifold, where $[\cdot]$ denotes the equivalence class (or orbit), defined in Eq. (\ref{equation_equivalence_X_in_Grassmann}).
\end{proposition}

\begin{proof}
According to Definition \ref{definition_retraction}, to prove that Eq. \eqref{equation_retraction_Grassmann_manifold} is a valid retraction, we need to verify two properties:

\begin{enumerate}
\item Identity:
\[
\operatorname{Ret}_{[\b{X}]}(\b{0}) = [\b{X}].
\]

\item Local rigidity:
\[
\left.\frac{d}{dt}\operatorname{Ret}_{[\b{X}]}(t\b{\Delta})\right|_{t=0}
=
\b{\Delta}.
\]
\end{enumerate}

\textbf{Step 1: Identity.}

Substituting \(\b{\Delta}=\b{0}\) into Eq. \eqref{equation_retraction_Grassmann_manifold} gives:
\begin{align*}
\operatorname{Ret}_{[\b{X}]}(\b{0})
&=
\left[
(\b{X}+\b{0})
(\b{I}_d + \b{0}^\top \b{0})^{-1/2}
\right]
\nonumber\\
&=
\left[
\b{X}\,\b{I}_d^{-1/2}
\right]
=
[\b{X}].
\end{align*}
So the identity property holds.

\textbf{Step 2: Local rigidity.}

We define the curve:
\[
\b{Y}(t)
:=
(\b{X}+t\b{\Delta})
(\b{I}_d + t^2 \b{\Delta}^\top \b{\Delta})^{-1/2}.
\]
Then:
\[
\operatorname{Ret}_{[\b{X}]}(t\b{\Delta}) = [\b{Y}(t)].
\]
We first expand the matrix factor $(\b{I}_d + t^2 \b{\Delta}^\top \b{\Delta})^{-1/2}$.
Since this matrix is smooth in \(t\) and equals \(\b{I}_d\) at \(t=0\), its Taylor expansion at \(t=0\) is:
\begin{equation}\label{equation_inverse_sqrt_taylor_Grassmann_retraction}
(\b{I}_d + t^2 \b{\Delta}^\top \b{\Delta})^{-1/2}
=
\b{I}_d + \mathcal{O}(t^2),
\end{equation}
where $\mathcal{O}(.)$ denotes the big-O complexity notation. 
Therefore:
\begin{align}
&\b{Y}(t)
=
(\b{X}+t\b{\Delta})(\b{I}_d + \mathcal{O}(t^2))
\nonumber\\
&~~~~~~~~ =
\b{X}+t\b{\Delta}+\mathcal{O}(t^2).
\label{equation_Yt_first_order_Grassmann_retraction} \\
&\implies \dot{\b{Y}}(t) = \b{\Delta}.
\end{align}


Hence, we have:
\[
\b{Y}(0)=\b{X},
\qquad
\dot{\b{Y}}(0)=\b{\Delta}.
\]
Now, a point on the Grassmann manifold is an equivalence class \([\b{X}]\), not a specific matrix \(\b{X}\). Hence, the derivative of the quotient curve \([\b{Y}(t)]\) must be interpreted through a tangent representative. According to Eq. (\ref{equation_tangent_space_Grassmann}), tangent vectors on \(Gr(n,d)\) are represented by horizontal matrices \(\b{\Delta}\) satisfying:
\[
\b{X}^\top \b{\Delta}=\b{0}.
\]
Because the derivative \(\dot{\b{Y}}(0)=\b{\Delta}\) already satisfies this horizontal condition, it is exactly the tangent vector represented by the quotient curve \([\b{Y}(t)]\) at \(t=0\). Therefore,
\[
\left.\frac{d}{dt}[\b{Y}(t)]\right|_{t=0}=\b{\Delta}.
\]
Hence,
\[
\left.\frac{d}{dt}\operatorname{Ret}_{[\b{X}]}(t\b{\Delta})\right|_{t=0}
=
\left.\frac{d}{dt}[\b{Y}(t)]\right|_{t=0}
=
\b{\Delta},
\]
which proves the local rigidity property.

Since both identity and local rigidity are satisfied, Eq. \eqref{equation_retraction_Grassmann_manifold} is a valid retraction map on the Grassmann manifold.
\end{proof}

\begin{remark}[First-order agreement with the exponential map]
Let \(\b{\Delta} = \b{U}\b{\Sigma}\b{V}^\top\) be the compact SVD of the tangent vector \(\b{\Delta}\), as in Lemma \ref{lemma_SVD_tangent_vector_Grassmann}. 

On the one hand, according to Proposition \ref{proposition_geodesic_Grassmann}, the exponential map on the Grassmann manifold is:
\[
\operatorname{Exp}_{[\b{X}]}(\b{\Delta})
=
\left[
\b{X}\b{V}\cos(\b{\Sigma})\b{V}^\top
+
\b{U}\sin(\b{\Sigma})\b{V}^\top
\right].
\]
Using the Taylor expansions:
\[
\cos(\b{\Sigma}) = \b{I}_d + \mathcal{O}(\Vert\b{\Sigma}\Vert^2),
\quad
\sin(\b{\Sigma}) = \b{\Sigma} + \mathcal{O}(\Vert\b{\Sigma}\Vert^3),
\]
we obtain:
\[
\operatorname{Exp}_{[\b{X}]}(t\b{\Delta})
=
[\b{X} + t\b{\Delta} + \mathcal{O}(t^2)],
\]
where $\mathcal{O}(.)$ denotes the big-O complexity notation. 

On the other hand, Eq. \eqref{equation_Yt_first_order_Grassmann_retraction} gives:
\[
\operatorname{Ret}_{[\b{X}]}(t\b{\Delta})
=
[\b{X} + t\b{\Delta} + \mathcal{O}(t^2)].
\]
Therefore, the retraction agrees with the exponential map up to first order, which is exactly the defining behavior of a retraction.
\end{remark}

\begin{remark}[Relation of polar-style and QR-based retractions in Grassmann manifold]
The polar-style retraction in Eq. \eqref{equation_retraction_Grassmann_manifold} is an alternative to the QR-based retraction in Eq. \eqref{equation_retraction_Grassmann_manifold_qr}, on the Grassmann manifold. Both maps return a point on the Grassmann manifold, i.e., an equivalence class representing a \(d\)-dimensional subspace of \(\mathbb{R}^n\), and both agree with the exponential map up to first order. The QR-based retraction uses the \(Q\)-factor of the QR decomposition of \(\b{X}+\b{\Delta}\), while the polar-style retraction uses the symmetric inverse square root of the Gram matrix of \(\b{X}+\b{\Delta}\). In other words, the QR-based retraction is:
\[
\operatorname{Ret}_{[\b{X}]}(\b{\Delta})
=
[\operatorname{qf}(\b{X}+\b{\Delta})],
\]
whereas the polar-style retraction is:
\[
\operatorname{Ret}_{[\b{X}]}(\b{\Delta})
=
\left[
(\b{X}+\b{\Delta})
(\b{I}_d+\b{\Delta}^\top\b{\Delta})^{-1/2}
\right],
\]
where $[\cdot]$ denotes the equivalence class (or orbit), defined in Eq. (\ref{equation_equivalence_X_in_Grassmann}).
\end{remark}

\begin{remark}[Relation of the retractions on the Stiefel and Grassmann manifolds]
The formula in Eq. \eqref{equation_retraction_Grassmann_manifold_qr} is the quotient counterpart of the QR-based retraction on the Stiefel manifold in Eq. (\ref{eq:stiefel_retraction_qr}). 
Moreover, The formula in Eq. \eqref{equation_retraction_Grassmann_manifold} is the quotient counterpart of the polar-style retraction on the Stiefel manifold in Eq. (\ref{equation_retraction_Stiefel_manifold_polar_2}). 

On the Stiefel manifold, the retraction returns the orthonormal matrix:
\begin{align*}
&\operatorname{Ret}_{\b{X}}(\b{\Delta}) = \operatorname{qf}(\b{X}+\b{\Delta}), \quad \text{or,} \\
&\operatorname{Ret}_{\b{X}}(\b{\Delta})
=
(\b{X}+\b{\Delta})
(\b{I}_d+\b{\Delta}^\top\b{\Delta})^{-1/2}.
\end{align*}
On the Grassmann manifold, however, only the subspace matters, so we take its equivalence class:
\begin{align*}
&\operatorname{Ret}_{[\b{X}]}(\b{\Delta})
=
[\operatorname{qf}(\b{X}+\b{\Delta})], \quad \text{or,} \\
&\operatorname{Ret}_{[\b{X}]}(\b{\Delta})
=
\left[
(\b{X}+\b{\Delta})
(\b{I}_d + \b{\Delta}^\top \b{\Delta})^{-1/2}
\right].
\end{align*}
\end{remark}

\subsubsection{Vector Transport in Grassmann Manifold}

As discussed in Section \ref{section_vector_transport}, vector transport is used for
moving a tangent vector from one tangent space to another,
typically along a search direction. On the Grassmann manifold,
because tangent spaces at different points are different linear
subspaces of $\mathbb{R}^{n \times d}$, a tangent vector at
$[\b{X}]$ cannot in general be used directly as a tangent vector
at another point $[\b{Y}]$. Therefore, it must be transported.

A simple and natural vector transport on the Grassmann manifold
is obtained by orthogonally projecting the ambient matrix onto
the tangent space at the new point. This is the Grassmann
counterpart of the projection-based vector transport introduced
for the Stiefel manifold.

We first state a useful lemma on the tangent projection.

\begin{lemma}[Projection onto the tangent space of the Grassmann manifold]
Let $[\b{Y}] \in Gr(n,d)$, where $\b{Y} \in St(n,d)$ is a representative.
For any ambient matrix $\b{Z} \in \mathbb{R}^{n \times d}$, the orthogonal
projection of $\b{Z}$ onto the tangent space $T_{[\b{Y}]}Gr(n,d)$ is:
\begin{equation}
\boxed{
\Pi_{[\b{Y}]}(\b{Z})
=
(\b{I}_n - \b{Y}\b{Y}^\top)\b{Z}.
}
\label{equation_projection_tangent_space_grassmann_for_vector_transport}
\end{equation}
\end{lemma}

\begin{proof}
According to the characterization of tangent space in Eq. \eqref{equation_tangent_space_Grassmann},
a matrix $\b{\Xi} \in \mathbb{R}^{n\times d}$ belongs to
$T_{[\b{Y}]}Gr(n,d)$ if and only if:
\begin{equation}
\b{Y}^\top \b{\Xi} = \b{0}.
\label{equation_tangent_constraint_grassmann_vector_transport}
\end{equation}
Now, we define:
\begin{equation}\label{equation_Z_tilde_I_YYT_Z}
\widetilde{\b{Z}}
:=
(\b{I}_n - \b{Y}\b{Y}^\top)\b{Z}.
\end{equation}
We verify that $\widetilde{\b{Z}}$ satisfies the tangent constraint:
\begin{align*}
\b{Y}^\top \widetilde{\b{Z}}
&\overset{(\ref{equation_Z_tilde_I_YYT_Z})}{=}
\b{Y}^\top(\b{I}_n - \b{Y}\b{Y}^\top)\b{Z} \nonumber\\
&=
(\b{Y}^\top - \b{Y}^\top \b{Y}\b{Y}^\top)\b{Z} \nonumber\\
&\overset{(a)}{=}
(\b{Y}^\top - \b{I}_d \b{Y}^\top)\b{Z} = \b{0},
\end{align*}
where $(a)$ uses $\b{Y}^\top \b{Y} = \b{I}_d$ because $\b{Y} \in St(n,d)$, according to Eq. (\ref{equation_Sitefel_manifold_definition}).

As $\b{Y}^\top \widetilde{\b{Z}} = \b{0}$, it satisfies Eq. (\ref{equation_tangent_space_Grassmann}), so $\widetilde{\b{Z}} \in T_{[\b{Y}]}Gr(n,d)$.

It remains to show that this is indeed the orthogonal projection.
We decompose $\b{Z}$ as:
\begin{equation}
\b{Z}
=
(\b{I}_n - \b{Y}\b{Y}^\top)\b{Z}
+
\b{Y}\b{Y}^\top \b{Z}.
\label{equation_decomposition_ambient_grassmann}
\end{equation}
The first term is tangent, as shown above. The second term lies in the orthogonal complement of the tangent space because for any
$\b{\Xi} \in T_{[\b{Y}]}Gr(n,d)$:
\begin{align*}
\langle \b{Y}\b{Y}^\top \b{Z}, \b{\Xi} \rangle_F
&\overset{(\ref{equation_Frobenius_inner_product})}{=}
\operatorname{tr}\!\Big((\b{Y}\b{Y}^\top \b{Z})^\top \b{\Xi}\Big) \nonumber\\
&=
\operatorname{tr}\!\Big(\b{Z}^\top \b{Y}\b{Y}^\top \b{\Xi}\Big) \nonumber\\
&=
\operatorname{tr}\!\Big(\b{Z}^\top \b{Y}(\b{Y}^\top \b{\Xi})\Big)
\overset{\eqref{equation_tangent_constraint_grassmann_vector_transport}}{=} 0.
\end{align*}

Therefore, the decomposition in Eq.
\eqref{equation_decomposition_ambient_grassmann}
is an orthogonal decomposition of $\b{Z}$ into tangent and normal parts.
Hence, the tangent projection is exactly:
\[
\Pi_{[\b{Y}]}(\b{Z})
=
(\b{I}_n - \b{Y}\b{Y}^\top)\b{Z}.
\]
\end{proof}

\begin{proposition}[Vector transport in Grassmann manifold]
Let $[\b{X}] \in Gr(n,d)$ and let
$\b{\Delta}_1 \in T_{[\b{X}]}Gr(n,d)$.
Given a direction
$\b{\Delta}_2 \in T_{[\b{X}]}Gr(n,d)$
and a retraction:
\begin{equation}
[\b{Y}] = \operatorname{Ret}_{[\b{X}]}(\b{\Delta}_2),
\end{equation}
the projection-based vector transport of $\b{\Delta}_1$ from
$T_{[\b{X}]}Gr(n,d)$ to $T_{[\b{Y}]}Gr(n,d)$ is defined as:
\begin{equation}
\boxed{
\mathcal{T}_{\b{\Delta}_2}(\b{\Delta}_1)
=
(\b{I}_n - \b{Y}\b{Y}^\top)\b{\Delta}_1.
}
\label{equation_vector_transport_grassmann_projection}
\end{equation}
In particular,
\begin{equation}
\mathcal{T}_{\b{\Delta}_2}(\b{\Delta}_1)
\in
T_{[\b{Y}]}Gr(n,d).
\end{equation}
\end{proposition}

\begin{proof}
The idea is the same as in the Stiefel case: we regard
$\b{\Delta}_1$ as an element of the ambient Euclidean space
$\mathbb{R}^{n\times d}$ and project it onto the tangent space at
the new point $[\b{Y}]$.

By Eq. (\ref{equation_projection_tangent_space_grassmann_for_vector_transport}),
the orthogonal projection onto $T_{[\b{Y}]}Gr(n,d)$ is:
\[
\Pi_{[\b{Y}]}(\b{Z})
=
(\b{I}_n - \b{Y}\b{Y}^\top)\b{Z},
\qquad
\forall \b{Z} \in \mathbb{R}^{n\times d}.
\]
Applying this projection to $\b{\Delta}_1$ gives:
\begin{equation}\label{equation_TDelta2Delta1_I_YYT_Delta1_in_proof}
\mathcal{T}_{\b{\Delta}_2}(\b{\Delta}_1)
:=
\Pi_{[\b{Y}]}(\b{\Delta}_1)
=
(\b{I}_n - \b{Y}\b{Y}^\top)\b{\Delta}_1,
\end{equation}
which is exactly Eq.
\eqref{equation_vector_transport_grassmann_projection}.

It remains to verify explicitly that this matrix is tangent at $[\b{Y}]$.
We compute:
\begin{align*}
\b{Y}^\top \mathcal{T}_{\b{\Delta}_2}(\b{\Delta}_1)
&\overset{(\ref{equation_TDelta2Delta1_I_YYT_Delta1_in_proof})}{=}
\b{Y}^\top (\b{I}_n - \b{Y}\b{Y}^\top)\b{\Delta}_1 \nonumber\\
&=
(\b{Y}^\top - \b{Y}^\top \b{Y}\b{Y}^\top)\b{\Delta}_1 \nonumber\\
&\overset{(a)}{=}
(\b{Y}^\top - \b{I}_d \b{Y}^\top)\b{\Delta}_1 =
\b{0},
\end{align*}
where $(a)$ uses $\b{Y}^\top \b{Y} = \b{I}_d$ because $\b{Y} \in St(n,d)$, according to Eq. (\ref{equation_Sitefel_manifold_definition}).

As $\b{Y}^\top \mathcal{T}_{\b{\Delta}_2}(\b{\Delta}_1) = \b{0}$, it satisfies Eq. (\ref{equation_tangent_space_Grassmann}), so:
\[
\mathcal{T}_{\b{\Delta}_2}(\b{\Delta}_1)
\in
T_{[\b{Y}]}Gr(n,d),
\]
which proves that Eq.
\eqref{equation_vector_transport_grassmann_projection}
is a valid vector transport from $T_{[\b{X}]}Gr(n,d)$ to
$T_{[\b{Y}]}Gr(n,d)$.
\end{proof}

\begin{remark}[Interpretation of projection-based vector transport on Grassmann manifold]
The vector transport in Eq.
\eqref{equation_vector_transport_grassmann_projection}
has a simple interpretation. The matrix $\b{\Delta}_1$ is first viewed
as an ambient matrix in $\mathbb{R}^{n\times d}$. Since the tangent
space changes from $T_{[\b{X}]}Gr(n,d)$ to $T_{[\b{Y}]}Gr(n,d)$,
the matrix $\b{\Delta}_1$ is generally no longer tangent at the new point.
Multiplication by $(\b{I}_n - \b{Y}\b{Y}^\top)$ removes the component
of $\b{\Delta}_1$ in the span of $\b{Y}$ and keeps only the component
orthogonal to $\b{Y}$, which is exactly the tangent component at $[\b{Y}]$.
\end{remark}

\begin{remark}[Vector transport associated with the Grassmann retraction]
If the retraction in Eq. (\ref{equation_retraction_Grassmann_manifold_qr}) is used, namely:
\[
\operatorname{Ret}_{[\b{X}]}(\b{\Delta}_2)
=
[\operatorname{qf}(\b{X}+\b{\Delta}_2)],
\]
then setting:
\[
\b{Y} = \operatorname{qf}(\b{X}+\b{\Delta}_2),
\]
in Eq. \eqref{equation_vector_transport_grassmann_projection}
gives the explicit transported vector:
\begin{equation}
\boxed{
\mathcal{T}_{\b{\Delta}_2}(\b{\Delta}_1)
=
\big(\b{I}_n - \operatorname{qf}(\b{X}+\b{\Delta}_2)\,\operatorname{qf}(\b{X}+\b{\Delta}_2)^\top\big)\b{\Delta}_1.
}
\label{equation_vector_transport_grassmann_qr_retraction}
\end{equation}
This formula is computationally simple and is widely used in
Riemannian optimization algorithms on the Grassmann manifold.
\end{remark}

\begin{remark}[Comparison of vector transports in the Stiefel and Grassmann manifolds]
The Grassmann vector transport is simpler than the Stiefel vector
transport. On the Stiefel manifold, according to Eq. (\ref{equation_stiefel_tangent_space}), the tangent constraint is:
\[
\text{Tangent on Stiefel: } \b{X}^\top \b{\Delta} + \b{\Delta}^\top \b{X} = \b{0},
\]
so the projection involves the symmetrization operator, defined in Eq. (\ref{equation_sym_skew_expressions}).
On the Grassmann manifold, according to Eq. (\ref{equation_tangent_space_Grassmann}), the tangent constraint reduces to:
\[
\text{Tangent on Grassmann: } \b{X}^\top \b{\Delta} = \b{0},
\]
hence the orthogonal projection is simply:
\[
(\b{I}_n - \b{X}\b{X}^\top)\b{Z}.
\]
That is why the vector transport formula on the Grassmann manifold
takes the simpler form in Eq.
\eqref{equation_vector_transport_grassmann_projection}, compared to the vector transport formula on the Stiefel manifold, stated in Eq. (\ref{equation_retraction_Stiefel}).
Here are the two formulas for easier comparison:
\begin{align*}
\text{Stiefel: } &\mathcal{T}_{\b{\Delta}_2}(\b{\Delta}_1) = \b{\Delta}_1 - \b{Y} \text{sym}(\b{Y}^\top \b{\Delta}_1), \\
\text{Grassmann: } &\mathcal{T}_{\b{\Delta}_2}(\b{\Delta}_1)
=
(\b{I}_n - \b{Y}\b{Y}^\top)\b{\Delta}_1.
\end{align*}
\end{remark}

\subsubsection{Vector Transport by Differential QR-Based Retraction in Grassmann Manifold}

As discussed in Section \ref{section_differentiated_retraction}, a natural way to construct a vector transport is by differentiating a chosen retraction.
For the Grassmann manifold, we can use the QR-based retraction
introduced in Eq. (\ref{equation_retraction_Grassmann_manifold_qr}):
\begin{equation*}
\operatorname{Ret}_{[\b{X}]}(\b{\Delta})
=
[\operatorname{qf}(\b{X}+\b{\Delta})].
\end{equation*}
We now derive the associated differentiated-retraction vector transport.

\begin{lemma}[Horizontal representative of a derivative on Grassmann manifold]\label{Lemma_horizontal_representative_derivative_grassmann}
Let $\b{Y}(t) \in St(n,d)$ be a smooth curve of representatives
for a curve $[\b{Y}(t)] \in Gr(n,d)$.
The tangent vector:
\[
\frac{d}{dt}[\b{Y}(t)]\Big|_{t=0}
\in T_{[\b{Y}(0)]}Gr(n,d),
\]
is represented by the horizontal component of $\dot{\b{Y}}(0)$, namely:
\begin{equation}
\left(\b{I}_n-\b{Y}\b{Y}^\top\right)\dot{\b{Y}}(0),
\qquad \text{where } \b{Y}:=\b{Y}(0).
\label{equation_horizontal_representative_derivative_grassmann}
\end{equation}
\end{lemma}

\begin{proof}
According to the tangent-space characterization of the Grassmann manifold, i.e., Eq. (\ref{equation_tangent_space_Grassmann}), a matrix $\b{\Xi}$ is a tangent representative at $[\b{Y}]$ if and only if:
\[
\b{Y}^\top \b{\Xi} = \b{0}.
\]
Now, we decompose $\dot{\b{Y}}(0)$ into its orthogonal and parallel parts with
respect to the column space of $\b{Y}$:
\begin{equation}\label{equation_Ydot0_I_YYT_Ydot0_YYTYdot0_in_proof}
\dot{\b{Y}}(0)
=
\left(\b{I}_n-\b{Y}\b{Y}^\top\right)\dot{\b{Y}}(0)
+
\b{Y}\b{Y}^\top \dot{\b{Y}}(0).
\end{equation}
The first term in Eq. (\ref{equation_Ydot0_I_YYT_Ydot0_YYTYdot0_in_proof}) is horizontal because:
\begin{align*}
&\b{Y}^\top \left(\b{I}_n-\b{Y}\b{Y}^\top\right)\dot{\b{Y}}(0)
=
\left(\b{Y}^\top-\b{Y}^\top\b{Y}\b{Y}^\top\right)\dot{\b{Y}}(0)
\nonumber\\
&\overset{(a)}{=}
\left(\b{Y}^\top-\b{I}_d\b{Y}^\top\right)\dot{\b{Y}}(0)
= \left(\b{Y}^\top-\b{Y}^\top\right)\dot{\b{Y}}(0) =
\b{0},
\end{align*}
where $(a)$ uses $\b{Y}^\top \b{Y} = \b{I}$ because $\b{Y} \in \text{St}(n,d)$ according to Eq. (\ref{equation_Sitefel_manifold_definition}).  

The second term in Eq. (\ref{equation_Ydot0_I_YYT_Ydot0_YYTYdot0_in_proof}) lies in the vertical space because it is of the form
$\b{Y}\b{\Omega}$ with $\b{\Omega}=\b{Y}^\top \dot{\b{Y}}(0)\in\mathbb{R}^{d\times d}$.
On the Grassmann manifold, vertical components correspond only to changes
of representative inside the same equivalence class, hence they do not change
the tangent vector of the quotient curve. Therefore, the derivative of the
Grassmann curve is represented by its horizontal part, which is Eq.
\eqref{equation_horizontal_representative_derivative_grassmann}.
\end{proof}

\begin{proposition}[Differentiated QR-retraction vector transport on Grassmann manifold]
Let $[\b{X}] \in Gr(n,d)$, where $\b{X}\in St(n,d)$ is a representative.
Let:
\[
\b{\Delta}_1,\b{\Delta}_2 \in T_{[\b{X}]}Gr(n,d),
\quad
\b{X}^\top \b{\Delta}_1=\b{0},
\quad
\b{X}^\top \b{\Delta}_2=\b{0}.
\]
Consider the QR-based retraction in Eq. (\ref{equation_retraction_Grassmann_manifold_qr}):
\[
\operatorname{Ret}_{[\b{X}]}(\b{\Delta})
=
[\operatorname{qf}(\b{X}+\b{\Delta})].
\]
Let the QR decomposition of $\b{X}+\b{\Delta}_2$ be:
\begin{align}
&\boxed{\b{X}+\b{\Delta}_2
=
\b{Y}\b{R},}
\end{align}
where:
\begin{align}
&\boxed{\b{Y}
:=
\operatorname{qf}(\b{X}+\b{\Delta}_2),}
\label{equation_qr_factorization_for_vector_transport_grassmann}
\end{align}
and $\b{R}\in\mathbb{R}^{d\times d}$ is the upper triangular factor of the
QR decomposition.

The differentiated retraction:
\[
\mathcal{T}^{\operatorname{QR}}_{\b{\Delta}_2}(\b{\Delta}_1) = D\operatorname{Ret}_{[\b{X}]}(\b{\Delta}_2)[\b{\Delta}_1]
\in
T_{[\b{Y}]}Gr(n,d),
\]
is represented by:
\begin{equation}
\boxed{
\mathcal{T}^{\operatorname{QR}}_{\b{\Delta}_2}(\b{\Delta}_1)
=
\left(\b{I}_n-\b{Y}\b{Y}^\top\right)\b{\Delta}_1\b{R}^{-1}.}
\label{equation_differentiated_qr_retraction_vector_transport_grassmann}
\end{equation}
The Eq.
\eqref{equation_differentiated_qr_retraction_vector_transport_grassmann}
defines the differentiated-retraction vector transport associated with the
QR-based retraction on the Grassmann manifold.
\end{proposition}

\begin{proof}
We derive the differentiated retraction by perturbing the retraction input
$\b{\Delta}_2$ along the tangent direction $\b{\Delta}_1$.

\textbf{Step 1: Define a perturbed retraction curve.}

For $t$ in a neighborhood of zero, we define:
\begin{equation*}
\b{Y}(t)
:=
\operatorname{qf}(\b{X}+\b{\Delta}_2+t\b{\Delta}_1).
\end{equation*}
According to the definition of QR-based retraction in Eq. (\ref{equation_retraction_Grassmann_manifold_qr}), the corresponding curve on the Grassmann manifold is:
\begin{equation*}
[\b{Y}(t)]
=
\operatorname{Ret}_{[\b{X}]}(\b{\Delta}_2+t\b{\Delta}_1).
\end{equation*}
By definition of differentiated retraction, we want:
\begin{align*}
D\operatorname{Ret}_{[\b{X}]}(\b{\Delta}_2)[\b{\Delta}_1]
&=
\frac{d}{dt}\operatorname{Ret}_{[\b{X}]}(\b{\Delta}_2+t\b{\Delta}_1)\Big|_{t=0}
\\
&=
\frac{d}{dt}[\b{Y}(t)]\Big|_{t=0}.
\end{align*}

\textbf{Step 2: Differentiate the QR factorization.}

Because $\b{Y}(t)$ is the Q-factor, there exists an upper triangular R-factor matrix $\b{R}(t)$ such that:
\begin{equation}
\b{X}+\b{\Delta}_2+t\b{\Delta}_1
=
\b{Y}(t)\b{R}(t).
\label{equation_qr_curve_factorization_grassmann}
\end{equation}
At $t=0$, this becomes:
\[
\b{X}+\b{\Delta}_2
=
\b{Y}\b{R},
\]
because $\b{Y}(0) = \b{Y}$ and $\b{R}(0) = \b{R}$. This equation is Eq. \eqref{equation_qr_factorization_for_vector_transport_grassmann}.

Differentiating Eq.
\eqref{equation_qr_curve_factorization_grassmann}
with respect to $t$ at $t=0$ gives:
\begin{equation}
\b{\Delta}_1
=
\dot{\b{Y}}(0)\b{R}
+
\b{Y}\dot{\b{R}}(0).
\end{equation}
Multiplying both sides on the right by $\b{R}^{-1}$ yields:
\begin{align}
&\b{\Delta}_1 \b{R}^{-1}
= 
\dot{\b{Y}}(0)\b{R} \b{R}^{-1}
+
\b{Y}\dot{\b{R}}(0) \b{R}^{-1} \nonumber\\
&\overset{(a)}{\implies} \dot{\b{Y}}(0)
=
\b{\Delta}_1\b{R}^{-1}
-
\b{Y}\dot{\b{R}}(0)\b{R}^{-1},
\label{equation_dotY_general_grassmann}
\end{align}
where $(a)$ is because $\b{R}\b{R}^{-1} = \b{I}$.

\textbf{Step 3: Remove the vertical component.}

The second term in Eq.
\eqref{equation_dotY_general_grassmann}
is in the span of $\b{Y}$, so it is vertical for the quotient representation.
By Lemma \ref{Lemma_horizontal_representative_derivative_grassmann} and Eq.
\eqref{equation_horizontal_representative_derivative_grassmann},
the tangent vector on the Grassmann manifold is represented by the horizontal part of $\dot{\b{Y}}(0)$:
\begin{align}
\frac{d}{dt}[&\b{Y}(t)]\Big|_{t=0}
\sim
\left(\b{I}_n-\b{Y}\b{Y}^\top\right)\dot{\b{Y}}(0)
\nonumber\\
&\overset{(\ref{equation_dotY_general_grassmann})}{=}
\left(\b{I}_n-\b{Y}\b{Y}^\top\right)
\left(
\b{\Delta}_1\b{R}^{-1}
-
\b{Y}\dot{\b{R}}(0)\b{R}^{-1}
\right)
\nonumber\\
&\overset{(a)}{=}
\left(\b{I}_n-\b{Y}\b{Y}^\top\right)\b{\Delta}_1\b{R}^{-1},
\label{equation_horizontal_part_dotY_grassmann}
\end{align}
where $(a)$ is because:
\[
\left(\b{I}_n-\b{Y}\b{Y}^\top\right)\b{Y}=\b{0}.
\]
Therefore, according to Eq. (\ref{equation_horizontal_part_dotY_grassmann}) and Lemma \ref{Lemma_horizontal_representative_derivative_grassmann}, the $D\operatorname{Ret}_{[\b{X}]}(\b{\Delta}_2)[\b{\Delta}_1]$ is represented by:
\[
\left(\b{I}_n-\b{Y}\b{Y}^\top\right)\b{\Delta}_1\b{R}^{-1}.
\]

\textbf{Step 4: Verify tangency.}

Let:
\[
\b{\Xi}
:=
\left(\b{I}_n-\b{Y}\b{Y}^\top\right)\b{\Delta}_1\b{R}^{-1}.
\]
Then:
\begin{align*}
\b{Y}^\top \b{\Xi}
&=
\b{Y}^\top
\left(\b{I}_n-\b{Y}\b{Y}^\top\right)\b{\Delta}_1\b{R}^{-1}
\nonumber\\
&=
\left(\b{Y}^\top-\b{Y}^\top\b{Y}\b{Y}^\top\right)\b{\Delta}_1\b{R}^{-1}
\nonumber\\
&\overset{(a)}{=}
\left(\b{Y}^\top-\b{I}_d\b{Y}^\top\right)\b{\Delta}_1\b{R}^{-1}
=
\b{0}.
\end{align*}
where $(a)$ uses $\b{Y}^\top \b{Y} = \b{I}$ because $\b{Y} \in \text{St}(n,d)$ according to Eq. (\ref{equation_Sitefel_manifold_definition}). 

We proved $\b{Y}^\top \b{\Xi} = \b{0}$, so $\b{\Xi}$ satisfies Eq. (\ref{equation_tangent_space_Grassmann}). 
Hence:
\[
\b{\Xi}\in T_{[\b{Y}]}Gr(n,d).
\]

Therefore, the differentiated-retraction vector transport associated with the
QR-based retraction is exactly Eq.
\eqref{equation_differentiated_qr_retraction_vector_transport_grassmann}.
\end{proof}

\begin{remark}[Special case of differentiated QR-based retraction at zero step in Grassmann manifold]
If $\b{\Delta}_2=\b{0}$, then $\b{Y}=\b{X}$ and $\b{R}=\b{I}_d$.
Hence, Eq.
\eqref{equation_differentiated_qr_retraction_vector_transport_grassmann}
reduces to:
\[
\mathcal{T}^{\operatorname{QR}}_{\b{0}}(\b{\Delta}_1)
=
(\b{I}_n-\b{X}\b{X}^\top)\b{\Delta}_1.
\]
Since $\b{\Delta}_1 \in T_{[\b{X}]}Gr(n,d)$, the $\b{\Delta}_1$ satisfies
$\b{X}^\top \b{\Delta}_1=\b{0}$ according to Eq. (\ref{equation_tangent_space_Grassmann}). Therefore, we have:
\[
(\b{I}_n-\b{X}\b{X}^\top)\b{\Delta}_1=\b{\Delta}_1 - \b{X}\underbrace{\b{X}^\top\b{\Delta}_1}_{=\b{0}} = \b{\Delta}_1.
\]
So, the differentiated retraction agrees with the identity map at zero step,
as expected.
\end{remark}

\begin{remark}[Interpretation of differentiated QR-based retraction in Grassmann manifold]
The factor $\b{R}^{-1}$ appears because the Q-factor of the QR decomposition
changes not only by orthogonal motion but also by the change of the triangular
factor. The projection
$\left(\b{I}_n-\b{Y}\b{Y}^\top\right)$ removes the vertical component along
the representative $\b{Y}$, leaving the tangent representative on the
Grassmann manifold.
\end{remark}

\subsubsection{Vector Transport by Differential Polar-Style Retraction in Grassmann Manifold}

As discussed in Section \ref{section_differentiated_retraction}, a natural way to construct
a vector transport is by differentiating a chosen retraction.
For the Grassmann manifold, besides the QR-based retraction,
we can also use the polar-style retraction introduced in
Eq. (\ref{equation_retraction_Grassmann_manifold}):
\begin{equation*}
\operatorname{Ret}_{[\b{X}]}(\b{\Delta})
=
\Big[
(\b{X}+\b{\Delta})
(\b{I}_d+\b{\Delta}^\top \b{\Delta})^{-1/2}
\Big].
\end{equation*}
We now derive the associated differentiated-retraction
vector transport.

Recall Lemma \ref{lemma_derivative_inverse_square_root_matrix_factor} from Section \ref{section_differentiated_polar_retraction_stiefel}. This lemma will be used here for Eqs. (\ref{equation_S_definition_polar_transport_grassmann}) and (\ref{equation_S_dot_definition_polar_transport_grassmann}), too. 



\begin{proposition}[Differentiated polar-style retraction on Grassmann manifold]
Let $[\b{X}] \in Gr(n,d)$, where $\b{X}\in St(n,d)$ is a representative.
Let:
\[
\b{\Delta}_1,\b{\Delta}_2 \in T_{[\b{X}]}Gr(n,d),
\quad
\b{X}^\top \b{\Delta}_1=\b{0},
\quad
\b{X}^\top \b{\Delta}_2=\b{0}.
\]
Consider the polar-style retraction in Eq. (\ref{equation_retraction_Grassmann_manifold}):
\[
\operatorname{Ret}_{[\b{X}]}(\b{\Delta})
=
\Big[
(\b{X}+\b{\Delta})
(\b{I}_d+\b{\Delta}^\top \b{\Delta})^{-1/2}
\Big].
\]
We define:
\begin{equation}
\b{Y}
:=
(\b{X}+\b{\Delta}_2)
(\b{I}_d+\b{\Delta}_2^\top \b{\Delta}_2)^{-1/2}.
\label{equation_Y_polar_retraction_grassmann_transport}
\end{equation}
The differentiated retraction:
\[
\mathcal{T}^{\operatorname{polar}}_{\b{\Delta}_2}(\b{\Delta}_1) = D\operatorname{Ret}_{[\b{X}]}(\b{\Delta}_2)[\b{\Delta}_1]
\in
T_{[\b{Y}]}Gr(n,d),
\]
is represented by:
\begin{equation}\label{equation_differentiated_polar_retraction_grassmann}
\boxed{
\begin{aligned}
\mathcal{T}^{\operatorname{polar}}_{\b{\Delta}_2}(\b{\Delta}_1)
= (\b{I}_n-\b{Y}\b{Y}^\top)
&\Big(
\b{\Delta}_1(\b{I}_d+\b{\Delta}_2^\top \b{\Delta}_2)^{-1/2}
\\
&+
(\b{X}+\b{\Delta}_2)\,\dot{\b{S}}(0)
\Big),
\end{aligned}
}
\end{equation}
where:
\begin{align}
&\b{S}(t)
:=
\Big(
\b{I}_d+(\b{\Delta}_2+t\b{\Delta}_1)^\top(\b{\Delta}_2+t\b{\Delta}_1)
\Big)^{-1/2}, \label{equation_S_definition_polar_transport_grassmann} \\
&\dot{\b{S}}(0)
=
\frac{d}{dt}\b{S}(t)\Big|_{t=0}. \label{equation_S_dot_definition_polar_transport_grassmann}
\end{align}
The Eq. \eqref{equation_differentiated_polar_retraction_grassmann}
defines the vector transport by differentiated polar-style retraction
on the Grassmann manifold.
\end{proposition}

\begin{proof}
We compute the differential of the polar-style retraction at
$\b{\Delta}_2$ in the direction $\b{\Delta}_1$.

\textbf{Step 1: Perturb the tangent input.}

Define the input curve:
\begin{equation}
\b{\Delta}(t)
:=
\b{\Delta}_2+t\b{\Delta}_1.
\label{equation_delta_t_polar_transport_grassmann}
\end{equation}
Applying the polar-style retraction, i.e., Eq. (\ref{equation_retraction_Grassmann_manifold}), to this curve gives:
\begin{equation}
\operatorname{Ret}_{[\b{X}]}(\b{\Delta}(t))
=
\Big[
(\b{X}+\b{\Delta}(t))
(\b{I}_d+\b{\Delta}(t)^\top \b{\Delta}(t))^{-1/2}
\Big].
\end{equation}
For convenience, we define the representative curve:
\begin{equation}
\b{Y}(t)
:=
(\b{X}+\b{\Delta}_2+t\b{\Delta}_1)\,\b{S}(t),
\label{equation_Y_t_polar_transport_grassmann}
\end{equation}
where $\b{S}(t)$ is defined in Eq.
\eqref{equation_S_definition_polar_transport_grassmann}.

According to the definition of polar-style retraction in Eq. (\ref{equation_retraction_Grassmann_manifold}), the corresponding curve on the Grassmann manifold is:
\begin{equation}
[\b{Y}(t)]
=
\operatorname{Ret}_{[\b{X}]}(\b{\Delta}_2+t\b{\Delta}_1).
\label{equation_curve_quotient_polar_transport_grassmann}
\end{equation}
Therefore:
\[
D\operatorname{Ret}_{[\b{X}]}(\b{\Delta}_2)[\b{\Delta}_1]
=
\frac{d}{dt}[\b{Y}(t)]\Big|_{t=0}.
\]

\textbf{Step 2: Differentiate the representative curve.}

Differentiating Eq. \eqref{equation_Y_t_polar_transport_grassmann}
with respect to $t$ gives:
\begin{equation}
\dot{\b{Y}}(t)
=
\b{\Delta}_1 \b{S}(t)
+
(\b{X}+\b{\Delta}_2+t\b{\Delta}_1)\dot{\b{S}}(t).
\end{equation}
Evaluating at $t=0$ yields:
\begin{equation}
\dot{\b{Y}}(0)
=
\b{\Delta}_1 \b{S}(0)
+
(\b{X}+\b{\Delta}_2)\dot{\b{S}}(0).
\label{equation_Y_dot_zero_polar_transport_grassmann}
\end{equation}
According to Eq.
\eqref{equation_S_definition_polar_transport_grassmann}, we have:
\[
\b{S}(0)
=
(\b{I}_d+\b{\Delta}_2^\top \b{\Delta}_2)^{-1/2}.
\]
Thus, Eq. \eqref{equation_Y_dot_zero_polar_transport_grassmann} becomes:
\begin{equation}
\dot{\b{Y}}(0)
=
\b{\Delta}_1(\b{I}_d+\b{\Delta}_2^\top \b{\Delta}_2)^{-1/2}
+
(\b{X}+\b{\Delta}_2)\dot{\b{S}}(0).
\label{equation_Y_dot_zero_expanded_polar_transport_grassmann}
\end{equation}

\textbf{Step 3: Pass from the representative derivative to the Grassmann tangent vector.}

By Lemma
\ref{Lemma_horizontal_representative_derivative_grassmann} and Eq.
\eqref{equation_horizontal_representative_derivative_grassmann},
the derivative of the quotient curve $[\b{Y}(t)]$ is represented by
the horizontal component of $\dot{\b{Y}}(0)$:
\begin{align}
&\frac{d}{dt}[\b{Y}(t)]\Big|_{t=0}
\sim
(\b{I}_n-\b{Y}\b{Y}^\top)\dot{\b{Y}}(0)
\nonumber\\
&=
(\b{I}_n-\b{Y}\b{Y}^\top)
\Big(
\b{\Delta}_1(\b{I}_d+\b{\Delta}_2^\top \b{\Delta}_2)^{-1/2}
\nonumber\\
&\qquad\qquad\qquad\qquad\qquad+
(\b{X}+\b{\Delta}_2)\dot{\b{S}}(0)
\Big),
\end{align}
where $\b{Y}=\b{Y}(0)$ is exactly the matrix in Eq.
\eqref{equation_Y_polar_retraction_grassmann_transport}.

Therefore, $D\operatorname{Ret}_{[\b{X}]}(\b{\Delta}_2)[\b{\Delta}_1]$ is represented by Eq.
\eqref{equation_differentiated_polar_retraction_grassmann}.

\textbf{Step 4: Verify tangency.}

Let:
\begin{align*}
\b{\Xi}
:=
(\b{I}_n-\b{Y}\b{Y}^\top)
\Big(
\b{\Delta}_1(\b{I}_d+&\b{\Delta}_2^\top \b{\Delta}_2)^{-1/2}
\\
&+
(\b{X}+\b{\Delta}_2)\dot{\b{S}}(0)
\Big).
\end{align*}
Then:
\begin{align*}
&\b{Y}^\top \b{\Xi}
\\
&=
\b{Y}^\top(\b{I}_n-\b{Y}\b{Y}^\top)
\Big(
\b{\Delta}_1(\b{I}_d+\b{\Delta}_2^\top \b{\Delta}_2)^{-1/2}
\\
&\qquad\qquad\qquad\qquad\qquad\qquad+
(\b{X}+\b{\Delta}_2)\dot{\b{S}}(0)
\Big)
\nonumber\\
&=
(\b{Y}^\top-\b{Y}^\top\b{Y}\b{Y}^\top)
\Big(
\b{\Delta}_1(\b{I}_d+\b{\Delta}_2^\top \b{\Delta}_2)^{-1/2}
\\
&\qquad\qquad\qquad\qquad\qquad\qquad+
(\b{X}+\b{\Delta}_2)\dot{\b{S}}(0)
\Big)
\nonumber\\
&\overset{(a)}{=}
(\b{Y}^\top-\b{I}_d\b{Y}^\top)(\cdots)
=
\b{0},
\end{align*}
where $(a)$ uses $\b{Y}^\top \b{Y} = \b{I}$ because $\b{Y} \in \text{St}(n,d)$ according to Eq. (\ref{equation_Sitefel_manifold_definition}).

We proved $\b{Y}^\top \b{\Xi} = \b{0}$, so $\b{\Xi}$ satisfies Eq. (\ref{equation_tangent_space_Grassmann}). 
Hence:
\[
\b{\Xi}\in T_{[\b{Y}]}Gr(n,d).
\]
This proves that Eq.
\eqref{equation_differentiated_polar_retraction_grassmann}
is a valid tangent representative of the differentiated retraction.
\end{proof}

\begin{remark}[Special case of differentiated polar-style retraction at zero step in Grassmann manifold]
If $\b{\Delta}_2=\b{0}$, then $\b{Y}=\b{X}$ and
$\b{S}(0)=\b{I}_d$. Hence:
\[
\dot{\b{Y}}(0)
=
\b{\Delta}_1+\b{X}\dot{\b{S}}(0).
\]
After horizontal projection at $[\b{X}]$, the vertical term
$\b{X}\dot{\b{S}}(0)$ vanishes, so:
\[
\mathcal{T}^{\operatorname{polar}}_{\b{0}}(\b{\Delta}_1)
=
(\b{I}_n-\b{X}\b{X}^\top)\b{\Delta}_1
=
\b{\Delta}_1,
\]
because $\b{X}^\top \b{\Delta}_1=\b{0}$.
Thus, the differentiated retraction reduces to the identity map at
zero step, as expected.
\end{remark}

\begin{remark}[Interpretation of differentiated polar-style retraction at zero step in Grassmann manifold]
The formula in Eq.
\eqref{equation_differentiated_polar_retraction_grassmann}
is the Grassmann counterpart of differentiated polar-style retraction
on the Stiefel manifold. The first term:
\[
\b{\Delta}_1(\b{I}_d+\b{\Delta}_2^\top \b{\Delta}_2)^{-1/2},
\]
comes from differentiating the explicit factor $\b{X}+\b{\Delta}$,
while the second term:
\[
(\b{X}+\b{\Delta}_2)\dot{\b{S}}(0),
\]
comes from differentiating the normalization factor
$(\b{I}_d+\b{\Delta}^\top \b{\Delta})^{-1/2}$.
Finally, the projection $(\b{I}_n-\b{Y}\b{Y}^\top)$ removes the vertical
part and keeps only the tangent representative on the Grassmann manifold.
\end{remark}

\subsection{Symmetric Positive Definite (SPD) Manifold}




In this subsection, we study the manifold of \textit{Symmetric Positive Definite (SPD)} matrices. 
This manifold plays an important role in Riemannian optimization because many quantities in machine learning and signal processing are naturally SPD matrices. 
Examples include covariance matrices, kernel matrices, diffusion tensors, and positive quadratic forms. 

Unlike the Stiefel and Grassmann manifolds, which are defined by orthogonality or quotient constraints, the SPD manifold is characterized by a positivity condition together with symmetry. 
The positivity condition is an open condition inside the vector space of symmetric matrices. 
Therefore, the SPD manifold has a particularly simple manifold structure: it is an open subset of the space of symmetric matrices. 
This fact makes several geometric constructions on the SPD manifold especially natural.

Many of the geometric, computational, and optimization-related
characteristics of the SPD manifold are analyzed in
\cite{bhatia2009positive}.

\subsubsection{Definition of SPD Manifold}

We first define the ambient linear space for the SPD manifold, which is the space of symmetric matrices. 
Recall that we defined symmetric matrices in Definition \ref{definition_symmetric_skew_symmetric_matrices} and Eq. (\ref{equation_symmetric_matrix}). In the followng, we define the space of symmetric matrices. 

\begin{definition}[Space of symmetric matrices]
The set of all \textbf{real symmetric} $n\times n$ matrices is denoted by:
\begin{align}
\boxed{
\mathbb{S}^n
:=
\{
\b{X}\in\mathbb{R}^{n\times n}
\mid
\b{X}^\top=\b{X}
\}.
}
\end{align}
The symmetry condition $\b{X}^\top=\b{X}$ means that the entries satisfy:
\begin{align}
\boxed{
X_{ij}=X_{ji}, \qquad \forall i,j\in\{1,\dots,n\},
}
\end{align}
where $X_{ij}$ denotes the $(i,j)$-th element of matrix $\b{X}$.
\end{definition}

\begin{definition}[Symmetric positive definite matrix]
A matrix $\b{X}\in\mathbb{S}^n$ is called \textbf{symmetric positive definite} if:
\begin{align}
\boxed{
\b{v}^\top \b{X} \b{v} > 0,\qquad \forall \b{v}\in\mathbb{R}^n\setminus\{\b{0}\}.
}
\end{align}
Equivalently, $\b{X}$ is symmetric positive definite if and only if all eigenvalues of $\b{X}$ are positive.
\end{definition}

\begin{definition}[SPD manifold]
The \textbf{symmetric positive definite (SPD) manifold} is defined as:
\begin{align}
\boxed{
\mathcal{M}
=
\mathbb{S}_{++}^n
:=
\{
\b{X}\in\mathbb{S}^n
\mid
\b{X}\succ \b{0}
\}.
}
\end{align}
Equivalently:
\begin{align}
\boxed{
\mathbb{S}_{++}^n
=
\{
\b{X}\in\mathbb{R}^{n\times n}
\mid
\b{X}^\top=\b{X},
\;
\b{v}^\top \b{X}\b{v}>0,\ \forall \b{v}\neq \b{0}
\}.
}
\end{align}
\end{definition}

In other words, the SPD manifold is the set of all symmetric matrices whose eigenvalues are strictly positive. 
Because all eigenvalues are positive, every matrix in $\mathbb{S}_{++}^n$ is invertible.

Examples of SPD matrices include:
\begin{itemize}
\item Covariance matrices,
\item Matrices defining quadratic forms $\b{x}^\top \b{W}\b{x}$,
\item Matrices used in generalized Mahalanobis distance metric \cite{ghojogh2022spectral}:
\[
(\b{x} - \b{y})^\top \b{W} (\b{x} - \b{y}),
\]
where $\b{x}, \b{y} \in \mathbb{R}^n$ are two vectors and the $n \times n$ matrix $\b{W}$ has to be positive (semi-)definite to satisfy triangle inequality and convexity of the distance metric.
\end{itemize}

Note that some of the characteristics of SPD manifold are discussed and derived by Suvrit Sra and Reshad Hosseini in \cite{sra2015conic,sra2016geometric,hosseini2020alternative}.

\begin{remark}[Why the SPD manifold is a manifold]
The set $\mathbb{S}^n$ is a vector space. 
The positive definiteness condition $\b{X}\succ \b{0}$ is strict, so it is preserved under sufficiently small perturbations inside $\mathbb{S}^n$. 
Hence, $\mathbb{S}_{++}^n$ is an open subset of $\mathbb{S}^n$. 
Therefore, $\mathbb{S}_{++}^n$ inherits a smooth manifold structure from the ambient vector space $\mathbb{S}^n$.
\end{remark}

\begin{proposition}[The SPD manifold is a smooth manifold]\label{proposition_spd_manifold_smooth}
The set:
\[
\mathbb{S}_{++}^n
=
\{
\b{X}\in\mathbb{S}^n
\mid
\b{X}\succ \b{0}
\},
\]
is a smooth manifold. More precisely, it is an open submanifold of the vector space $\mathbb{S}^n$.
\end{proposition}

\begin{proof}
We prove the claim in steps.

\textbf{Step 1: The ambient space $\mathbb{S}^n$ is a finite-dimensional vector space.}

By definition,
\[
\mathbb{S}^n
=
\{\b{X}\in\mathbb{R}^{n\times n}\mid \b{X}^\top=\b{X}\}.
\]
If $\b{X},\b{Y}\in\mathbb{S}^n$ and $a,b\in\mathbb{R}$, then:
\[
(a\b{X}+b\b{Y})^\top
=
a\b{X}^\top+b\b{Y}^\top
=
a\b{X}+b\b{Y}.
\]
Hence, $a\b{X}+b\b{Y}\in\mathbb{S}^n$. Therefore, $\mathbb{S}^n$ is a vector subspace of $\mathbb{R}^{n\times n}$, and thus a finite-dimensional Euclidean space. 
So $\mathbb{S}^n$ is itself a smooth manifold.

\hfill\break
\textbf{Step 2: Positive definiteness is an open condition in $\mathbb{S}^n$.}

Take any $\b{X}\in \mathbb{S}_{++}^n$. Since $\b{X}$ is symmetric positive definite, all its eigenvalues are positive. Let $\lambda_{\min}(\b{X})>0$ denote its smallest eigenvalue.

Now let $\b{E}\in\mathbb{S}^n$ be a symmetric perturbation satisfying:
\[
\Vert\b{E}\Vert_2 < \lambda_{\min}(\b{X}),
\]
where $\Vert\cdot\Vert_2$ denotes the $\ell_2$ norm of matrix. 

For any nonzero vector $\b{v}\in\mathbb{R}^n$, we have:
\[
\b{v}^\top (\b{X}+\b{E})\b{v}
=
\b{v}^\top \b{X}\b{v} + \b{v}^\top \b{E}\b{v}.
\]
Because $\b{X}$ is symmetric positive definite, we have:
\begin{align}\label{equation_vTXvgeqlambdaminCv_norm2_in_proof}
\b{v}^\top \b{X}\b{v}
\geq
\lambda_{\min}(\b{X})\Vert\b{v}\Vert_2^2.
\end{align}
Also, by the definition of operator norm, we have:
\begin{align}\label{equation_vTEvEnormvnorm2_in_proof}
|\b{v}^\top \b{E}\b{v}|
\leq
\Vert\b{E}\Vert_2 \Vert\b{v}\Vert_2^2,
\end{align}
where $|\cdot|$ denotes the absolute value. 
Equation (\ref{equation_vTEvEnormvnorm2_in_proof}) can be restated as:
\begin{align}
-\Vert\b{E}\Vert_2 \Vert\b{v}\Vert_2^2 \leq \b{v}^\top \b{E}\b{v} &\leq \Vert\b{E}\Vert_2 \Vert\b{v}\Vert_2^2 \implies \nonumber\\
&\b{v}^\top \b{E}\b{v} \leq \Vert\b{E}\Vert_2 \Vert\b{v}\Vert_2^2, \nonumber\\
&\b{v}^\top \b{E}\b{v} \geq -\Vert\b{E}\Vert_2 \Vert\b{v}\Vert_2^2. \label{equation_vTEVgeqminusEnormvnorm2_in_proof}
\end{align}
Therefore:
\begin{align*}
\b{v}^\top (\b{X}+\b{E})\b{v} &= \b{v}^\top \b{X}\b{v} + \b{v}^\top \b{E}\b{v}
\\
&\overset{(a)}{\geq}
\left(\lambda_{\min}(\b{X})-\Vert\b{E}\Vert_2\right)\Vert\b{v}\Vert_2^2,
\end{align*}
where $(a)$ is because of Eqs. (\ref{equation_vTXvgeqlambdaminCv_norm2_in_proof}) and (\ref{equation_vTEVgeqminusEnormvnorm2_in_proof}).

Since $\Vert\b{E}\Vert_2 < \lambda_{\min}(\b{X})$, the quantity in parentheses is positive. Hence:
\[
\b{v}^\top (\b{X}+\b{E})\b{v} > 0,\qquad \forall \b{v}\neq \b{0}.
\]
So $\b{X}+\b{E}$ is still symmetric positive definite.

This proves that every $\b{X}\in \mathbb{S}_{++}^n$ has a neighborhood in $\mathbb{S}^n$ entirely contained in $\mathbb{S}_{++}^n$. Hence, $\mathbb{S}_{++}^n$ is open in $\mathbb{S}^n$.

\hfill\break
\textbf{Step 3: Conclude the manifold structure.}

Every open subset of a smooth manifold is itself a smooth manifold with the induced smooth structure. Since $\mathbb{S}^n$ is a smooth manifold and $\mathbb{S}_{++}^n$ is an open subset of $\mathbb{S}^n$, it follows that $\mathbb{S}_{++}^n$ is a smooth manifold.

Therefore, the SPD manifold is an open submanifold of $\mathbb{S}^n$.
\end{proof}

\begin{remark}[Embedded viewpoint]
The manifold $\mathbb{S}_{++}^n$ can be viewed in two equivalent ways:
\begin{itemize}
\item as an open subset of the vector space $\mathbb{S}^n$, or
\item as a subset of $\mathbb{R}^{n\times n}$ satisfying the symmetry constraint together with the positivity condition.
\end{itemize}
In this monograph, it is often convenient to regard $\mathbb{S}_{++}^n$ as a matrix manifold embedded in $\mathbb{R}^{n\times n}$, while remembering that its natural linear ambient space is the symmetric subspace $\mathbb{S}^n$.
\end{remark}

\begin{lemma}[Points of the SPD manifold are symmetric positive definite matrices]
Let \(\b{X} \in \mathbb{R}^{n\times n}\). If
\(\b{X} \in \mathbb{S}_{++}^{n}\), then \(\b{X}\) is a point
of the SPD manifold. Conversely, every point
\(\b{p} \in \mathbb{S}_{++}^{n}\) can be written as $\b{p} = \b{X}$
for some \(\b{X} \in \mathbb{S}_{++}^{n}\) satisfying $\b{X}\in\mathbb{S}^n, \b{X}\succ \b{0}$.
Thus:
\begin{align}
\boxed{
\b{p} = \b{X} \quad \text{such that} \quad \b{X}\in\mathbb{S}^n, \b{X}\succ \b{0}.
}
\end{align}

In other words, points of the SPD manifold are exactly the
symmetric positive definite matrices.
\end{lemma}
\begin{proof}
By definition, the SPD manifold is:
\[
\mathbb{S}_{++}^{n}
:=
\left\{
\b{X} \in \mathbb{S}^{n}
\;\middle|\;
\b{X} \succ 0
\right\}.
\]
Equivalently:
\[
\mathbb{S}_{++}^{n}
=
\left\{
\b{X} \in \mathbb{R}^{n\times n}
\;\middle|\;
\b{X}^{\top} = \b{X},
\;
\b{v}^{\top}\b{X}\b{v} > 0,
\ \forall \b{v}\neq \b{0}
\right\}.
\]

Therefore, if \(\b{X} \in \mathbb{S}_{++}^{n}\), then
\(\b{X}\) is an element of the set
\(\mathbb{S}_{++}^{n}\), and thus it is a point of the SPD
manifold:
\[
\b{X} \in \mathbb{S}_{++}^{n}.
\]

Conversely, let \(\b{p} \in \mathbb{S}_{++}^{n}\). Since
\(\mathbb{S}_{++}^{n}\) is a set of matrices, every element
of it is, by definition, some symmetric positive definite
matrix \(\b{X}\). Therefore, there exists
\(\b{X} \in \mathbb{S}_{++}^{n}\) such that:
\[
\b{p} = \b{X}.
\]

Hence, points of the SPD manifold are exactly symmetric
positive definite matrices.
\end{proof}

\subsubsection{Dimension of SPD Manifold}

Since the SPD manifold is an open subset of the vector space of symmetric matrices, its dimension is the same as the dimension of that ambient vector space.

\begin{proposition}[Dimension of SPD manifold]
The SPD manifold
\[
\mathbb{S}_{++}^n=\{\b{X}\in \mathbb{S}_n\mid \b{X}\succ 0\},
\]
has dimension:
\begin{align}
\boxed{
\dim(\mathbb{S}_{++}^n)=\frac{n(n+1)}{2}.
}
\end{align}
\end{proposition}

\begin{proof}
By Proposition \ref{proposition_spd_manifold_smooth}, the SPD manifold \(\mathbb{S}_{++}^n\) is an open submanifold of the vector space \(\mathbb{S}_n\). Therefore, its dimension is equal to the dimension of \(\mathbb{S}_n\):
\begin{align}
\boxed{
\dim(\mathbb{S}_{++}^n)=\dim(\mathbb{S}_n).
}
\end{align}
So, it remains to compute the dimension of the vector space \(\mathbb{S}_n\).

Let \(\b{X}\in \mathbb{S}_n\). Because \(\b{X}\) is symmetric, we have:
\[
\b{X}^\top=\b{X},
\]
which means that the entries satisfy:
\[
X_{ij}=X_{ji}, \qquad \forall i,j\in\{1,\dots,n\}.
\]
Hence, the entries of \(\b{X}\) are not all independent.

We count the number of independent entries of a symmetric \(n\times n\) matrix:

\textbf{1) Diagonal entries:}
The diagonal entries:
\[
X_{11},X_{22},\dots,X_{nn}
\]
can be chosen arbitrarily. Therefore, they contribute \(n\) degrees of freedom.

\textbf{2) Off-diagonal entries:}
For every pair \(i<j\), the entry \(X_{ij}\) determines \(X_{ji}\) because:
\[
X_{ji}=X_{ij}.
\]
So, for each unordered pair \((i,j)\) with \(i<j\), we have one independent parameter.

The number of such pairs is:
\[
\binom{n}{2}=\frac{n(n-1)}{2}.
\]

Therefore, the total number of independent parameters in a symmetric matrix is:
\[
n+\frac{n(n-1)}{2}
=
\frac{2n+n(n-1)}{2}
=
\frac{n(n+1)}{2}.
\]
Thus:
\[
\dim(\mathbb{S}_n)=\frac{n(n+1)}{2}.
\]

Since \(\mathbb{S}_{++}^n\) is an open subset of \(\mathbb{S}_n\), it has the same dimension as \(\mathbb{S}_n\). Hence:
\[
\dim(\mathbb{S}_{++}^n)=\frac{n(n+1)}{2}.
\]
This proves the claim.
\end{proof}

\begin{remark}[Why there is no reduction in dimension from positive definiteness]
Unlike the Stiefel manifold, where the orthogonality condition:
\[
\b{X}^\top \b{X}=\b{I}_d,
\]
imposes equality constraints and reduces the dimension, the condition:
\[
\b{X}\succ 0,
\]
is an open condition inside the space \(\mathbb{S}_n\). It does not impose additional equality constraints. Therefore, the SPD manifold has the same dimension as the space of symmetric matrices:
\[
\dim(\mathbb{S}_{++}^n)=\dim(\mathbb{S}_n)=\frac{n(n+1)}{2}.
\]
\end{remark}

\subsubsection{Tangent and Normal Spaces of SPD Manifold}\label{section_tangent_normal_spd}

Because the SPD manifold \(\mathbb{S}_{++}^n\) is an open subset of the vector space \(\mathbb{S}^n\), its tangent space is particularly simple:
it is the whole ambient linear space \(\mathbb{S}^n\) itself (we will prove it in Proposition \ref{proposition_tangent_space_spd}).
The normal space depends on which ambient space we use.
If we view \(\mathbb{S}_{++}^n\) as an embedded submanifold of
\(\mathbb{R}^{n\times n}\), then the normal space is the orthogonal
complement of \(\mathbb{S}^n\) in \(\mathbb{R}^{n\times n}\) (we will prove it in Proposition \ref{proposition_normal_space_spd}).

\begin{lemma}[Orthogonal complement of symmetric matrices]\label{lemma_orthogonal_complement_of_symmetric_matrices}
Let:
\[
\mathbb{S}^n := \{\b{X}\in\mathbb{R}^{n\times n}\mid \b{X}^\top=\b{X}\}.
\]
Consider the Frobenius inner product on \(\mathbb{R}^{n\times n}\), i.e., Eq. (\ref{equation_Frobenius_inner_product}):
\[
\langle \b{A},\b{B}\rangle_F := \operatorname{tr}(\b{A}^\top\b{B}).
\]
Then, the orthogonal complement of \(\mathbb{S}^n\) in
\(\mathbb{R}^{n\times n}\) is the set of skew-symmetric matrices:
\begin{align}
\boxed{
(\mathbb{S}^n)^\perp
=
\{\b{N}\in\mathbb{R}^{n\times n}\mid \b{N}^\top=-\b{N}\}.
}
\end{align}
\end{lemma}

\begin{proof}
We prove both inclusions (both sides of ``if and only if").

\textbf{(i) Every skew-symmetric matrix is orthogonal to every symmetric matrix.}

Let \(\b{N}\in\mathbb{R}^{n\times n}\) be skew-symmetric, so according to Eq. (\ref{equation_skew_symmetric_matrix}):
\[
\b{N}^\top=-\b{N}.
\]
Let \(\b{S}\in\mathbb{S}^n\) be symmetric, so according to Eq. (\ref{equation_symmetric_matrix}):
\[
\b{S}^\top=\b{S}.
\]
Then:
\[
\langle \b{N},\b{S}\rangle_F
\overset{(\ref{equation_Frobenius_inner_product})}{=}
\operatorname{tr}(\b{N}^\top\b{S})
=
\operatorname{tr}((-\b{N})\b{S})
=
-\operatorname{tr}(\b{N}\b{S}).
\]
On the other hand, using invariance of trace under transpose, we have:
\begin{align*}
\operatorname{tr}(\b{N}\b{S})
&=
\operatorname{tr}((\b{N}\b{S})^\top)
=
\operatorname{tr}(\b{S}^\top\b{N}^\top)
=
\operatorname{tr}(\b{S}(-\b{N}))
\\
&=
-\operatorname{tr}(\b{S}\b{N})
\overset{(a)}{=}
-\operatorname{tr}(\b{N}\b{S}),
\end{align*}
where $(a)$ is because of the cyclic property of trace.

We showed $\operatorname{tr}(\b{N}\b{S}) = -\operatorname{tr}(\b{N}\b{S})$. Therefore:
\[
\operatorname{tr}(\b{N}\b{S})=0.
\]
Hence:
\[
\langle \b{N},\b{S}\rangle_F = 0,
\qquad \forall \b{S}\in\mathbb{S}^n.
\]
So, every skew-symmetric matrix belongs to \((\mathbb{S}^n)^\perp\).

\hfill\break
\textbf{(ii) Every matrix orthogonal to all symmetric matrices is skew-symmetric.}

Now, let \(\b{N}\in(\mathbb{S}^n)^\perp\). We show that
\(\b{N}^\top=-\b{N}\).

According to Eq. (\ref{equation_sym_skew_expressions}) in Lemma \ref{lemma_sym_skew_expressions}, we decompose \(\b{N}\) into its symmetric and skew-symmetric parts:
\[
\b{N}
=
\frac{\b{N}+\b{N}^\top}{2}
+
\frac{\b{N}-\b{N}^\top}{2}.
\]
where we define:
\[
\b{N}_{\operatorname{sym}}
:=
\frac{\b{N}+\b{N}^\top}{2},
\qquad
\b{N}_{\operatorname{skew}}
:=
\frac{\b{N}-\b{N}^\top}{2}.
\]
Then \(\b{N}_{\operatorname{sym}}\in \mathbb{S}^n\), and
\(\b{N}_{\operatorname{skew}}^\top=-\b{N}_{\operatorname{skew}}\).

Since \(\b{N}\in(\mathbb{S}^n)^\perp\), it is orthogonal to every
symmetric matrix, in particular to \(\b{N}_{\operatorname{sym}}\). Hence:
\[
0
=
\langle \b{N}, \b{N}_{\operatorname{sym}}\rangle_F
=
\left\langle
\b{N}_{\operatorname{sym}}+\b{N}_{\operatorname{skew}},
\b{N}_{\operatorname{sym}}
\right\rangle_F.
\]
By part (i), every skew-symmetric matrix is orthogonal to every
symmetric matrix, so:
\[
\langle \b{N}_{\operatorname{skew}}, \b{N}_{\operatorname{sym}}\rangle_F = 0.
\]
Therefore:
\[
0
=
\langle \b{N}_{\operatorname{sym}},\b{N}_{\operatorname{sym}}\rangle_F
=
\Vert\b{N}_{\operatorname{sym}}\Vert_F^2.
\]
Thus:
\[
\b{N}_{\operatorname{sym}}=\b{0}.
\]
So:
\[
\frac{\b{N}+\b{N}^\top}{2}=\b{0}
\quad\Longrightarrow\quad
\b{N}^\top=-\b{N}.
\]
Hence, according to Eq. (\ref{equation_skew_symmetric_matrix}), \(\b{N}\) is skew-symmetric.

Combining (i) and (ii), we conclude:
\[
(\mathbb{S}^n)^\perp
=
\{\b{N}\in\mathbb{R}^{n\times n}\mid \b{N}^\top=-\b{N}\}.
\]
\end{proof}

\begin{proposition}[Tangent space of SPD manifold]\label{proposition_tangent_space_spd}
Let \(\b{X}\in\mathbb{S}_{++}^n\). The tangent space of the SPD
manifold at \(\b{X}\) is:
\begin{align}
\boxed{
T_{\b{X}}\mathbb{S}_{++}^n
=
\mathbb{S}^n
=
\{\b{\Delta}\in\mathbb{R}^{n\times n}\mid \b{\Delta}^\top=\b{\Delta}\}.
}
\end{align}
\end{proposition}

\begin{proof}
By Proposition \ref{proposition_spd_manifold_smooth}, the manifold \(\mathbb{S}_{++}^n\) is an open
submanifold of the vector space \(\mathbb{S}^n\).
A basic fact for open subsets of a vector space is that the tangent
space at every point is naturally identified with the vector space
itself\footnote{According to \citep[p.~59, Proposition 3.13]{lee2013smooth}, If $\mathcal{M}$ is an open submanifold of a vector space $\b{V}$, we can combine our identifications $T_{\b{p}}\mathcal{M} \cong T_p\b{V} \cong \b{V}$. Thus, since---according to Proposition \ref{proposition_spd_manifold_smooth}---the \(\mathbb{S}_{++}^n\) is an open submanifold of the vector space
\(\mathbb{S}^n\), we have the canonical identification
\(T_{\b{X}}\mathbb{S}_{++}^n \cong \mathbb{S}^n\).}. Therefore:
\[
T_{\b{X}}\mathbb{S}_{++}^n=\mathbb{S}^n.
\]

We can also verify this directly using curves.

\textbf{Step 1: Every tangent vector is symmetric.}

Let \(\b{\Delta}\in T_{\b{X}}\mathbb{S}_{++}^n\). By definition of
tangent vector, there exists a smooth curve:
\[
\b{X}(t)\subset \mathbb{S}_{++}^n
\]
such that:
\[
\b{X}(0)=\b{X},
\qquad
\dot{\b{X}}(0)=\b{\Delta}.
\]
Because \(\b{X}(t)\in\mathbb{S}_{++}^n\subset \mathbb{S}^n\) for all \(t\),
every \(\b{X}(t)\) is symmetric:
\[
\b{X}(t)^\top=\b{X}(t).
\]
Differentiating both sides with respect to \(t\) at \(t=0\) gives:
\[
\dot{\b{X}}(0)^\top=\dot{\b{X}}(0).
\]
Therefore:
\[
\b{\Delta}^\top=\b{\Delta},
\]
so \(\b{\Delta}\in\mathbb{S}^n\). Hence:
\[
T_{\b{X}}\mathbb{S}_{++}^n \subseteq \mathbb{S}^n.
\]

\textbf{Step 2: Every symmetric matrix is a tangent vector.}

Now let \(\b{\Delta}\in\mathbb{S}^n\). Consider the curve:
\[
\b{X}(t):=\b{X}+t\b{\Delta}.
\]
Since both \(\b{X}\) and \(\b{\Delta}\) are symmetric, we have:
\[
\b{X}(t)^\top=\b{X}(t),
\]
so, \(\b{X}(t)\) is symmetric for all \(t\).

Because \(\b{X}\in\mathbb{S}_{++}^n\), the set \(\mathbb{S}_{++}^n\) is
open in \(\mathbb{S}^n\). Therefore, for sufficiently small \(t\), we have:
\[
\b{X}(t)=\b{X}+t\b{\Delta}\in\mathbb{S}_{++}^n.
\]
Hence, \(\b{X}(t)\) is a smooth curve on the SPD manifold, with:
\[
\b{X}(0)=\b{X},
\qquad
\dot{\b{X}}(0)=\b{\Delta}.
\]
Thus, \(\b{\Delta}\in T_{\b{X}}\mathbb{S}_{++}^n\). Therefore:
\[
\mathbb{S}^n \subseteq T_{\b{X}}\mathbb{S}_{++}^n.
\]

Combining the two inclusions yields:
\[
T_{\b{X}}\mathbb{S}_{++}^n
=
\mathbb{S}^n.
\]
\end{proof}

\begin{proposition}[Normal space of SPD manifold in the ambient space \(\mathbb{R}^{n\times n}\)]\label{proposition_normal_space_spd}
Let \(\b{X}\in\mathbb{S}_{++}^n\). Viewing the SPD manifold as an
embedded submanifold of \(\mathbb{R}^{n\times n}\) equipped with the
Frobenius inner product, its normal space at \(\b{X}\) is:
\begin{align}\label{equation_normal_space_spd_in_Rnn}
\boxed{
N_{\b{X}}\mathbb{S}_{++}^n
=
\{\b{N}\in\mathbb{R}^{n\times n}\mid \b{N}^\top=-\b{N}\}.
}
\end{align}
\end{proposition}

\begin{proof}
By Proposition \ref{proposition_tangent_space_spd}, we have:
\[
T_{\b{X}}\mathbb{S}_{++}^n=\mathbb{S}^n.
\]
By definition, the normal space in the ambient Euclidean space
\(\mathbb{R}^{n\times n}\) is the orthogonal complement of the tangent
space:
\[
N_{\b{X}}\mathbb{S}_{++}^n
:=
\big(T_{\b{X}}\mathbb{S}_{++}^n\big)^\perp.
\]
Therefore:
\[
N_{\b{X}}\mathbb{S}_{++}^n
=
(\mathbb{S}^n)^\perp.
\]
By Lemma \ref{lemma_orthogonal_complement_of_symmetric_matrices}, the orthogonal complement of
\(\mathbb{S}^n\) in \(\mathbb{R}^{n\times n}\) is exactly the set of
skew-symmetric matrices. Hence:
\[
N_{\b{X}}\mathbb{S}_{++}^n
=
\{\b{N}\in\mathbb{R}^{n\times n}\mid \b{N}^\top=-\b{N}\}.
\]
\end{proof}

\begin{remark}[Normal space of SPD manifold inside \(\mathbb{S}^n\)]
If we regard \(\mathbb{S}_{++}^n\) as an open submanifold of
\(\mathbb{S}^n\), rather than an embedded submanifold of
\(\mathbb{R}^{n\times n}\), then its normal space is trivial:
\begin{align}
\boxed{
N_{\b{X}}\mathbb{S}_{++}^n=\{\b{0}\}.
}
\end{align}
This is because, in the ambient space \(\mathbb{S}^n\), the tangent
space already equals the whole ambient space:
\[
T_{\b{X}}\mathbb{S}_{++}^n=\mathbb{S}^n.
\]
The nontrivial normal space in Eq. (\ref{equation_normal_space_spd_in_Rnn}) is the normal space
with respect to the larger ambient space \(\mathbb{R}^{n\times n}\).
\end{remark}

\begin{remark}[Interpretation of the tangent and normal spaces of SPD manifold]
Unlike the Stiefel and Grassmann manifolds, where tangent vectors
must satisfy nontrivial linearized constraints, the SPD manifold is
an open set inside \(\mathbb{S}^n\). Therefore, every symmetric
matrix is an allowable tangent direction. The only directions normal
to the manifold in the ambient space \(\mathbb{R}^{n\times n}\) are the
skew-symmetric directions, because these move the matrix away from
symmetry rather than along the manifold.
\end{remark}

\subsubsection{Metric Tensors of SPD Manifold}\label{section_metrics_SPD}

\begin{remark}[Several metrics are used on the SPD manifold]
Unlike the Grassmann manifold, where the standard quotient
metric is usually adopted, the SPD manifold
$\mathbb{S}_{++}^n$ is endowed with several different
metrics in the literature \cite{bhatia2009positive}.

First, since $\mathbb{S}_{++}^n$ is an open subset of the
vector space $\mathbb{S}^n$, the \textbf{ambient Euclidean
(Frobenius) inner product}:
\[
g^E_{\b{X}}(\b{\Delta}_1,\b{\Delta}_2)
=
\text{tr}(\b{\Delta}_1^\top \b{\Delta}_2),
\]
defines a valid Riemannian metric on
$\mathbb{S}_{++}^n$ (we will prove it in Proposition \ref{proposition_euclidean_metric_spd}). This Euclidean metric is used in the
literature, although it is usually viewed as the metric
inherited from the ambient vector space rather than the most
intrinsic geometry of SPD matrices \cite{thanwerdas2022geometry,pennec2020manifold}.

Second, the \textbf{affine-invariant metric} \cite{moakher2005differential}:
\[
g^{AI}_{\b{X}}(\b{\Delta}_1,\b{\Delta}_2)
=
\text{tr}(\b{X}^{-1}\b{\Delta}_1\b{X}^{-1}\b{\Delta}_2),
\]
is one of the most classical and standard intrinsic metrics
on $\mathbb{S}_{++}^n$ (we will prove it in Proposition \ref{proposition_affine_invariant_metric_spd}). It was advocated in the SPD tensor literature by Pennec, Fillard, and Ayache
\cite{pennec2006riemannian}. 
The affine-invariant metric is the most well-known metric for SPD manifold.

Third, the \textbf{log-Euclidean metric} was introduced by Arsigny,
Fillard, Pennec, and Ayache \citep{Arsigny2005LogEuclidean,arsigny2005fast,arsigny2006log} as a computationally simpler
alternative to the affine-invariant geometry. In that framework, the
matrix logarithm transforms SPD matrices into a flat vector
space of symmetric matrices.

Finally, the \textbf{Bures--Wasserstein metric} \cite{bhatia2019bures} is another important
geometry on $\mathbb{S}_{++}^n$, especially in matrix
analysis, quantum information, and optimal transport.

Therefore, on the SPD manifold, several metrics are used in
the literature. Among them, the affine-invariant metric is
often regarded as the classical intrinsic choice, while the
Euclidean, log-Euclidean, and Bures--Wasserstein metrics are
also important in theory and applications.

Note that different metrics on SPD induce different
Riemannian gradients, Hessians, Levi-Civita connections,
geodesics, and exponential maps in SPD manifold.
\end{remark}

Recall from Section \ref{section_tangent_normal_spd} that for every
$\b{X} \in \mathbb{S}_{++}^n$, the tangent space is:
\[
T_{\b{X}}\mathbb{S}_{++}^n = \mathbb{S}^n.
\]
Therefore, for every point $\b{X} \in \mathbb{S}_{++}^n$,
the tangent vectors:
\[
\b{\Delta}, \b{\Delta}_1, \b{\Delta}_2
\in T_{\b{X}}\mathbb{S}_{++}^n
\]
are symmetric matrices.
In the following, we introduce Riemannian metrics on the SPD manifold. 
First, we introduce the Euclidean (Frobenius) metric on SPD manifold.

\begin{proposition}[Euclidean (Frobenius) metric tensor on SPD manifold]\label{proposition_euclidean_metric_spd}
A natural Riemannian metric on the SPD manifold
$\mathbb{S}_{++}^n$ is the Euclidean metric inherited from
the ambient space $\mathbb{S}^n$. For
$\b{\Delta}_1, \b{\Delta}_2 \in T_{\b{X}}\mathbb{S}_{++}^n$,
the metric is:
\begin{equation}
\label{equation_metric_tensor_SPD_manifold_Euclidean}
\boxed{
g_{\b{X}}^{E}(\b{\Delta}_1,\b{\Delta}_2)
=
\text{tr}(\b{\Delta}_1^\top \b{\Delta}_2)
\overset{(\ref{equation_Frobenius_inner_product})}{=}
\langle \b{\Delta}_1,\b{\Delta}_2 \rangle_F.
}
\end{equation}
Because $\b{\Delta}_1$ and $\b{\Delta}_2$ are symmetric,
this can also be written as:
\begin{equation}
\label{equation_metric_tensor_SPD_manifold_Euclidean_symmetric}
g_{\b{X}}^{E}(\b{\Delta}_1,\b{\Delta}_2)
=
\text{tr}(\b{\Delta}_1 \b{\Delta}_2).
\end{equation}
\end{proposition}

\begin{proof}
We show that for every fixed
$\b{X} \in \mathbb{S}_{++}^n$, the mapping
$g_{\b{X}}$ is an inner product on
$T_{\b{X}}\mathbb{S}_{++}^n$.

\textbf{Bilinearity:}
Let $a,b \in \mathbb{R}$ and let
$\b{\Delta}_1,\b{\Delta}_2,\b{\Delta}_3
\in T_{\b{X}}\mathbb{S}_{++}^n$. Then:
\begin{align*}
g_{\b{X}}^{E}(a\b{\Delta}_1+&b\b{\Delta}_2,\b{\Delta}_3)
\overset{(\ref{equation_metric_tensor_SPD_manifold_Euclidean})}{=}
\text{tr}\Big((a\b{\Delta}_1+b\b{\Delta}_2)^\top \b{\Delta}_3\Big) \\
&=
\text{tr}\Big((a\b{\Delta}_1^\top+b\b{\Delta}_2^\top)\b{\Delta}_3\Big) \\
&=
a\,\text{tr}(\b{\Delta}_1^\top \b{\Delta}_3)
+
b\,\text{tr}(\b{\Delta}_2^\top \b{\Delta}_3) \\
&\overset{(\ref{equation_metric_tensor_SPD_manifold_Euclidean})}{=}
a\,g_{\b{X}}^{E}(\b{\Delta}_1,\b{\Delta}_3)
+
b\,g_{\b{X}}^{E}(\b{\Delta}_2,\b{\Delta}_3).
\end{align*}
Similarly, it is linear in the second argument.

\textbf{Symmetry:}
\begin{align*}
g_{\b{X}}^{E}(\b{\Delta}_1,\b{\Delta}_2)
&\overset{(\ref{equation_metric_tensor_SPD_manifold_Euclidean})}{=}
\text{tr}(\b{\Delta}_1^\top \b{\Delta}_2) \overset{(a)}{=}
\text{tr}\big((\b{\Delta}_1^\top \b{\Delta}_2)^\top\big) \\
&\overset{(b)}{=}
\text{tr}(\b{\Delta}_2^\top \b{\Delta}_1) \overset{(\ref{equation_metric_tensor_SPD_manifold_Euclidean})}{=}
g_{\b{X}}^{E}(\b{\Delta}_2,\b{\Delta}_1),
\end{align*}
where $(a)$ is because the trace of a matrix is equal to the trace of its transpose and $(b)$ is because transpose reverses the order of matrix multiplication.

\textbf{Positive definiteness:}
Let
$\b{\Delta} \in T_{\b{X}}\mathbb{S}_{++}^n$
with $\b{\Delta} \neq \b{0}$. Then:
\[
g_{\b{X}}^{E}(\b{\Delta},\b{\Delta})
\overset{(\ref{equation_metric_tensor_SPD_manifold_Euclidean})}{=}
\text{tr}(\b{\Delta}^\top \b{\Delta})
=
\Vert\b{\Delta}\Vert_F^2
>
0.
\]
Also:
\[
g_{\b{X}}^{E}(\b{0},\b{0})=0.
\]
Hence, $g_{\b{X}}$ is positive definite.

Therefore,
Eq.~\eqref{equation_metric_tensor_SPD_manifold_Euclidean}
defines an inner product on
$T_{\b{X}}\mathbb{S}_{++}^n$. Since the formula is smooth
in $\b{X}$, it defines a Riemannian metric on
$\mathbb{S}_{++}^n$.
\end{proof}

In the following, we introduce the affine-invariant metric on SPD manifold. 

\begin{proposition}[Affine-invariant metric tensor on SPD manifold \cite{moakher2005differential}]\label{proposition_affine_invariant_metric_spd}
An important Riemannian metric on
$\mathbb{S}_{++}^n$ is the affine-invariant metric. For
$\b{\Delta}_1, \b{\Delta}_2 \in T_{\b{X}}\mathbb{S}_{++}^n$,
the metric is:
\begin{equation}
\label{equation_metric_tensor_SPD_manifold_affine_invariant}
\boxed{
g_{\b{X}}^{AI}(\b{\Delta}_1,\b{\Delta}_2)
=
\text{tr}(\b{X}^{-1}\b{\Delta}_1\b{X}^{-1}\b{\Delta}_2).
}
\end{equation}
\end{proposition}

\begin{proof}
We show that for every fixed
$\b{X} \in \mathbb{S}_{++}^n$, the mapping in
Eq.~\eqref{equation_metric_tensor_SPD_manifold_affine_invariant}
is an inner product on
$T_{\b{X}}\mathbb{S}_{++}^n$.

Because $\b{X} \in \mathbb{S}_{++}^n$, the matrix $\b{X}$
is invertible and $\b{X}^{-1}$ is symmetric positive
definite.

\textbf{Bilinearity:}
Let $a,b \in \mathbb{R}$ and let
$\b{\Delta}_1,\b{\Delta}_2,\b{\Delta}_3
\in T_{\b{X}}\mathbb{S}_{++}^n$. Then:
\begin{align*}
g_{\b{X}}^{AI}(a&\b{\Delta}_1+b\b{\Delta}_2,\b{\Delta}_3)
\\
&\overset{(\ref{equation_metric_tensor_SPD_manifold_affine_invariant})}{=}
\text{tr}\Big(
\b{X}^{-1}(a\b{\Delta}_1+b\b{\Delta}_2)\b{X}^{-1}\b{\Delta}_3
\Big) \\
&=
a\,\text{tr}(\b{X}^{-1}\b{\Delta}_1\b{X}^{-1}\b{\Delta}_3)
+
b\,\text{tr}(\b{X}^{-1}\b{\Delta}_2\b{X}^{-1}\b{\Delta}_3) \\
&\overset{(\ref{equation_metric_tensor_SPD_manifold_affine_invariant})}{=}
a\,g_{\b{X}}^{AI}(\b{\Delta}_1,\b{\Delta}_3)
+
b\,g_{\b{X}}^{AI}(\b{\Delta}_2,\b{\Delta}_3).
\end{align*}
Similarly, it is linear in the second argument:
\begin{align*}
g_{\b{X}}^{AI}(&\b{\Delta}_1,a\b{\Delta}_2+b\b{\Delta}_3)
=
a\,g_{\b{X}}^{AI}(\b{\Delta}_1,\b{\Delta}_2)
+
b\,g_{\b{X}}^{AI}(\b{\Delta}_1,\b{\Delta}_3).
\end{align*}

\textbf{Symmetry:}
Because $\b{X}^{-1}$, $\b{\Delta}_1$, and $\b{\Delta}_2$
are symmetric, we have:
\begin{align*}
g_{\b{X}}^{AI}(\b{\Delta}_1,\b{\Delta}_2)
&\overset{(\ref{equation_metric_tensor_SPD_manifold_affine_invariant})}{=}
\text{tr}(\b{X}^{-1}\b{\Delta}_1\b{X}^{-1}\b{\Delta}_2) \\
&\overset{(a)}{=}
\text{tr}\Big(
(\b{X}^{-1}\b{\Delta}_1\b{X}^{-1}\b{\Delta}_2)^\top
\Big) \\
&\overset{(b)}{=}
\text{tr}\Big(\b{\Delta}_2^\top(\b{X}^{-1})^\top\b{\Delta}_1^\top(\b{X}^{-1})^\top\Big) \\
&\overset{(c)}{=}
\text{tr}(\b{\Delta}_2\b{X}^{-1}\b{\Delta}_1\b{X}^{-1}) \\
&\overset{(d)}{=}
\text{tr}(\b{X}^{-1}\b{\Delta}_2\b{X}^{-1}\b{\Delta}_1) \\
&\overset{(\ref{equation_metric_tensor_SPD_manifold_affine_invariant})}{=}
g_{\b{X}}^{AI}(\b{\Delta}_2,\b{\Delta}_1),
\end{align*}
where $(a)$ is because the trace of a matrix is equal to the trace of its transpose, $(b)$ is because transpose reverses the order of matrix multiplication, $(c)$ is because $\b{X}^{-1}$, $\b{\Delta}_1$, and $\b{\Delta}_2$ are symmetric, and $(d)$ is because of the cyclic property of trace. 

\textbf{Positive definiteness:}
Let
$\b{\Delta} \in T_{\b{X}}\mathbb{S}_{++}^n$
with $\b{\Delta} \neq \b{0}$, and define:
\begin{align}\label{equation_Y_Xminushals_Delta_Xminushald_in_proof}
\b{Y}
:=
\b{X}^{-1/2}\b{\Delta}\b{X}^{-1/2}.
\end{align}
Then:
\begin{align*}
g_{\b{X}}^{AI}(\b{\Delta},\b{\Delta})
&\overset{(\ref{equation_metric_tensor_SPD_manifold_affine_invariant})}{=}
\text{tr}(\b{X}^{-1}\b{\Delta}\b{X}^{-1}\b{\Delta}) \\
&\overset{(a)}{=}
\text{tr}\!\Big(
\b{X}^{-1/2}\b{X}^{-1/2}\b{\Delta}\b{X}^{-1/2}
\b{X}^{-1/2}\b{\Delta}
\Big) \\
&\overset{(b)}{=}
\text{tr}\!\Big(
\b{X}^{-1/2}\b{\Delta}\b{X}^{-1/2}
\b{X}^{-1/2}\b{\Delta}\b{X}^{-1/2}
\Big) \\
&\overset{(\ref{equation_Y_Xminushals_Delta_Xminushald_in_proof})}{=} \text{tr}(\b{Y} \b{Y}) =
\text{tr}(\b{Y}^2),
\end{align*}
where $(a)$ is because $\b{X}^{-1} = \b{X}^{-1/2}\b{X}^{-1/2}$ and $(b)$ is because of the cyclic property of trace. 

Since $\b{X}^{-1/2}$ and $\b{\Delta}$ are symmetric,
$\b{Y}$ is symmetric. Hence:
\[
\text{tr}(\b{Y}^2) = \text{tr}(\b{Y}\b{Y}) 
\overset{(a)}{=}
\text{tr}(\b{Y}^\top \b{Y})
=
\Vert\b{Y}\Vert_F^2
\geq 0,
\]
where $(a)$ is because $\b{Y} = \b{Y}^\top$ since $\b{Y}$ is symmetric. 

We prove by contradiction in the following. If:
\[
g_{\b{X}}^{AI}(\b{\Delta},\b{\Delta})=0,
\]
then $\Vert\b{Y}\Vert_F^2=0$, so $\b{Y}=\b{0}$. Therefore:
\[
\b{\Delta}
=
\b{X}^{1/2}\b{Y}\b{X}^{1/2}
=
\b{0},
\]
which contradicts $\b{\Delta}\neq \b{0}$. Hence:
\[
g_{\b{X}}^{AI}(\b{\Delta},\b{\Delta})>0,
\qquad
\forall \b{\Delta}\neq \b{0}.
\]

Therefore,
Eq.~\eqref{equation_metric_tensor_SPD_manifold_affine_invariant}
defines an inner product on
$T_{\b{X}}\mathbb{S}_{++}^n$. Since matrix inversion is
smooth on $\mathbb{S}_{++}^n$, the metric depends smoothly
on $\b{X}$. Hence, it is a Riemannian metric on
$\mathbb{S}_{++}^n$.
\end{proof}

\begin{remark}[Affine invariance of the affine-invariant metric]
The metric in
Eq.~\eqref{equation_metric_tensor_SPD_manifold_affine_invariant}
is called affine-invariant because it is invariant under
congruence transformations. Let
$\b{A} \in \mathbb{R}^{n\times n}$ be invertible. Consider a congruence transformations by matrix $\b{A}$:
\[
\Phi(\b{X}) := \b{A}\b{X}\b{A}^\top.
\]
Then, for every $\b{X} \in \mathbb{S}_{++}^n$, we have
$\Phi(\b{X}) \in \mathbb{S}_{++}^n$.

The differential of $\Phi$ at $\b{X}$ applied to a tangent
vector $\b{\Delta} \in T_{\b{X}}\mathbb{S}_{++}^n$ is:
\[
D\Phi(\b{X})[\b{\Delta}]
=
\b{A}\b{\Delta}\b{A}^\top.
\]
Therefore, for
$\b{\Delta}_1,\b{\Delta}_2 \in T_{\b{X}}\mathbb{S}_{++}^n$,
we have:
\begin{align*}
&g_{\Phi(\b{X})}^{AI}
\big(
D\Phi(\b{X})[\b{\Delta}_1],
D\Phi(\b{X})[\b{\Delta}_2]
\big) \overset{(\ref{equation_metric_tensor_SPD_manifold_affine_invariant})}{=} \\
&~~
\text{tr}\Big(
(\b{A}\b{X}\b{A}^\top)^{-1}
(\b{A}\b{\Delta}_1\b{A}^\top)
(\b{A}\b{X}\b{A}^\top)^{-1}
(\b{A}\b{\Delta}_2\b{A}^\top)
\Big).
\end{align*}
Using
\[
(\b{A}\b{X}\b{A}^\top)^{-1}
=
\b{A}^{-\top}\b{X}^{-1}\b{A}^{-1},
\]
we get:
\begin{align*}
&g_{\Phi(\b{X})}^{AI}
\big(
D\Phi(\b{X})[\b{\Delta}_1],
D\Phi(\b{X})[\b{\Delta}_2]
\big) \\
& =
\text{tr}\Big(
\b{A}^{-\top}\b{X}^{-1}\underbrace{\b{A}^{-1}
\b{A}}_{=\b{I}}\b{\Delta}_1\b{A}^\top
\\
&\qquad\qquad\qquad\qquad\b{A}^{-\top}\b{X}^{-1}\underbrace{\b{A}^{-1}
\b{A}}_{=\b{I}}\b{\Delta}_2\b{A}^\top
\Big) \\
&=
\text{tr}\Big(
\b{A}^{-\top}\b{X}^{-1}\b{\Delta}_1
\b{X}^{-1}\b{\Delta}_2
\b{A}^\top
\Big) \\
&\overset{(a)}{=}
\text{tr}\Big(
\b{X}^{-1}\b{\Delta}_1
\b{X}^{-1}\b{\Delta}_2
\underbrace{\b{A}^\top \b{A}^{-\top}}_{=\b{I}}
\Big) \\
&=
\text{tr}\!\big(
\b{X}^{-1}\b{\Delta}_1
\b{X}^{-1}\b{\Delta}_2
\big) \\
&\overset{(\ref{equation_metric_tensor_SPD_manifold_affine_invariant})}{=} 
g_{\b{X}}^{AI}(\b{\Delta}_1,\b{\Delta}_2),
\end{align*}
where $(a)$ is because of cyclic property of trace.

Hence, the metric is invariant under congruence
transformations, which is why it is called the
affine-invariant metric.
\end{remark}

In the following, we introduce the log-Euclidean metric on SPD manifold. 

\begin{proposition}[Log-Euclidean metric tensor on SPD manifold]\label{proposition_log_euclidean_metric_spd}
Let:
\[
\mathrm{Log} : \mathbb{S}_{++}^n \to \mathbb{S}^n,
\qquad
\b{X} \mapsto \mathrm{Log}(\b{X}),
\]
denote the matrix logarithm map (see Definition \ref{definition_logarithm_map}). The log-Euclidean metric
on the SPD manifold is the pullback of the Euclidean
(Frobenius) metric on $\mathbb{S}^n$ by the map $\mathrm{Log}(\cdot)$:
\begin{equation}\label{equation_metric_tensor_SPD_manifold_log_Euclidean_2}
\boxed{
\begin{aligned}
g_{\b{X}}^{LE}(\b{\Delta}_1,&\b{\Delta}_2)
:= \\
&g_{\mathrm{Log}(\b{X})}^{E}
\big(
D\mathrm{Log}(\b{X})[\b{\Delta}_1],
D\mathrm{Log}(\b{X})[\b{\Delta}_2]
\big),
\end{aligned}
}
\end{equation}
for $\b{\Delta}_1, \b{\Delta}_2 \in T_{\b{X}}\mathbb{S}_{++}^n$, where $g_{\mathrm{Log}(\b{X})}^{E}$ is the Euclidean metric, defined in Eq. (\ref{equation_metric_tensor_SPD_manifold_Euclidean}).

Hence, for $\b{\Delta}_1, \b{\Delta}_2 \in T_{\b{X}}\mathbb{S}_{++}^n$,
the metric is:
\begin{equation}
\boxed{
g_{\b{X}}^{LE}(\b{\Delta}_1,\b{\Delta}_2)
=
\text{tr}\Big(
D\mathrm{Log}(\b{X})[\b{\Delta}_1]\,
D\mathrm{Log}(\b{X})[\b{\Delta}_2]
\Big).
}
\label{equation_metric_tensor_SPD_manifold_log_Euclidean}
\end{equation}
\end{proposition}

\begin{proof}
Recall from Section \ref{section_tangent_normal_spd} that:
\[
T_{\b{X}}\mathbb{S}_{++}^n = \mathbb{S}^n,
\qquad \forall \b{X} \in \mathbb{S}_{++}^n.
\]
Also, $\mathbb{S}_{++}^n$ is an open manifold of symmetric
positive definite matrices, while $\mathbb{S}^n$ is a vector
space of symmetric matrices.

The matrix exponential map (see Definition \ref{definition_exponential_map}) and matrix logarithm map (see Definition \ref{definition_logarithm_map}) are smooth
inverse maps between these two spaces:
\[
\mathrm{Exp} : \mathbb{S}^n \to \mathbb{S}_{++}^n,
\qquad
\mathrm{Log} : \mathbb{S}_{++}^n \to \mathbb{S}^n,
\]
with:
\[
\mathrm{Exp}(\mathrm{Log}(\b{X})) = \b{X},
\qquad
\mathrm{Log}(\mathrm{Exp}(\b{S})) = \b{S}.
\]
Therefore, $\mathrm{Log}(\cdot)$ is a diffeomorphism from
$\mathbb{S}_{++}^n$ onto $\mathbb{S}^n$.

Now, the Euclidean metric on the vector space $\mathbb{S}^n$
is the Frobenius inner product:
\begin{align*}
&g_{\b{S}}^{E}(\b{A},\b{B})
\overset{(\ref{equation_metric_tensor_SPD_manifold_Euclidean})}{=}
\text{tr}(\b{A}^\top\b{B}) \overset{(a)}{=} \text{tr}(\b{A}\b{B}), \\
&\b{A}, \b{B} \in T_{\b{S}}\mathbb{S}^n = \mathbb{S}^n,
\end{align*}
where $(a)$ is because $\b{A} = \b{A}^\top$ as $\b{A}$ is symmetric. 

By the definition of pullback metric (see Section \ref{section_pullback_metric}),
the pullback of $g^{E}$ by $\mathrm{Log}(\cdot)$ is:
\begin{align*}
g_{\b{X}}^{LE}(\b{\Delta}_1,&\b{\Delta}_2)
:= \\
&g_{\mathrm{Log}(\b{X})}^{E}
\big(
D\mathrm{Log}(\b{X})[\b{\Delta}_1],
D\mathrm{Log}(\b{X})[\b{\Delta}_2]
\big),
\end{align*}
which is Eq. (\ref{equation_metric_tensor_SPD_manifold_log_Euclidean_2}).
Substituting the Euclidean metric on $\mathbb{S}^n$ gives:
\[
g_{\b{X}}^{LE}(\b{\Delta}_1,\b{\Delta}_2)
=
\text{tr}\Big(
D\mathrm{Log}(\b{X})[\b{\Delta}_1]\,
D\mathrm{Log}(\b{X})[\b{\Delta}_2]
\Big),
\]
which proves Eq. \eqref{equation_metric_tensor_SPD_manifold_log_Euclidean}.

It remains to verify that this is a Riemannian metric.

\textbf{Bilinearity:}
Because $D\mathrm{Log}(\b{X})[\cdot]$ is linear in its argument and
the trace is bilinear, for $a,b \in \mathbb{R}$ we have:
\begin{align*}
g_{\b{X}}^{LE}(&a\b{\Delta}_1+b\b{\Delta}_2,\b{\Delta}_3) \\
&=
\text{tr}\Big(
D\mathrm{Log}(\b{X})[a\b{\Delta}_1+b\b{\Delta}_2]\,
D\mathrm{Log}(\b{X})[\b{\Delta}_3]
\Big) \\
&=
\text{tr}\Big(
\big(aD\mathrm{Log}(\b{X})[\b{\Delta}_1]\\
&~~~~~~~~~~~~~~~~~~~~~ +bD\mathrm{Log}(\b{X})[\b{\Delta}_2]\big)\,
D\mathrm{Log}(\b{X})[\b{\Delta}_3]
\Big) \\
&\overset{(\ref{equation_metric_tensor_SPD_manifold_log_Euclidean})}{=}
a\,g_{\b{X}}^{LE}(\b{\Delta}_1,\b{\Delta}_3)
+
b\,g_{\b{X}}^{LE}(\b{\Delta}_2,\b{\Delta}_3).
\end{align*}
Similarly, it is linear in the second argument.

\textbf{Symmetry:}
Because $\mathrm{Log}(\b{X}(t))$ is symmetric whenever $\b{X}(t)$ is
a smooth curve in $\mathbb{S}_{++}^n$, its derivative is also
symmetric. Hence:
\[
D\mathrm{Log}(\b{X})[\b{\Delta}_1],\,
D\mathrm{Log}(\b{X})[\b{\Delta}_2] \in \mathbb{S}^n.
\]
Therefore:
\begin{align*}
g_{\b{X}}^{LE}(\b{\Delta}_1,\b{\Delta}_2)
&=
\text{tr}\Big(
D\mathrm{Log}(\b{X})[\b{\Delta}_1]\,
D\mathrm{Log}(\b{X})[\b{\Delta}_2]
\Big) \\
&=
\text{tr}\Big(
D\mathrm{Log}(\b{X})[\b{\Delta}_2]\,
D\mathrm{Log}(\b{X})[\b{\Delta}_1]
\Big) \\
&=
g_{\b{X}}^{LE}(\b{\Delta}_2,\b{\Delta}_1),
\end{align*}
where we used the cyclic property of trace.

\textbf{Positive definiteness:}
Let $\b{\Delta} \in T_{\b{X}}\mathbb{S}_{++}^n$ with
$\b{\Delta} \neq \b{0}$. Then:
\begin{align*}
g_{\b{X}}^{LE}(\b{\Delta},\b{\Delta})
&=
\text{tr}\Big(
D\mathrm{Log}(\b{X})[\b{\Delta}]^2
\Big)
\\
&=
\big\Vert
D\mathrm{Log}(\b{X})[\b{\Delta}]
\big\Vert_F^2
\geq 0.
\end{align*}
If this were zero, then:
\[
D\mathrm{Log}(\b{X})[\b{\Delta}] = \b{0}.
\]
But since $\mathrm{Log}$ is a diffeomorphism, its differential
$D\mathrm{Log}(\b{X})$ is invertible. Hence $\b{\Delta} = \b{0}$,
which is a contradiction. Therefore:
\[
g_{\b{X}}^{LE}(\b{\Delta},\b{\Delta}) > 0,
\qquad \forall \b{\Delta} \neq \b{0}.
\]

Hence, Eq.
\eqref{equation_metric_tensor_SPD_manifold_log_Euclidean}
defines a Riemannian metric on $\mathbb{S}_{++}^n$.
\end{proof}

\begin{remark}[Interpretation of the log-Euclidean metric]
The log-Euclidean metric treats the SPD manifold by first
mapping every point $\b{X} \in \mathbb{S}_{++}^n$ to the flat
vector space $\mathbb{S}^n$ through the matrix logarithm.
Then, distances and inner products are measured there
using the ordinary Frobenius geometry. Therefore, this
metric turns the nonlinear SPD manifold into a Euclidean
space in logarithmic coordinates.
\end{remark}

In the following, we introduce the Bures--Wasserstein metric on SPD manifold. 

\begin{lemma}[Symmetric solution of the Sylvester equation
for SPD matrices]\label{lemma_Sylvester_equation_solution}
Let $\b{X} \in \mathbb{S}_{++}^n$ and let $\b{\Delta} \in \mathbb{S}^n$.
Then the \textbf{Sylvester equation} \cite{sylvester1851xxxvii}:
\begin{equation}
\boxed{
\b{X}\b{A} + \b{A}\b{X} = \b{\Delta},
}
\label{equation_Sylvester}
\end{equation}
has a unique symmetric solution $\b{A} \in \mathbb{S}^n$.

Let the eigendecomposition of $\b{X}$ be \cite{ghojogh2019eigenvalue}:
\[
\b{X} = \b{U}\b{\Lambda}\b{U}^\top,
\qquad
\b{\Lambda} = \text{diag}(\lambda_1,\dots,\lambda_n),
\qquad
\lambda_i > 0,
\]
where the columns of $\b{U}$ are eigenvectors of $\b{X}$ and $\{\lambda_1,\dots,\lambda_n\}$ are the eigenvalues of $\b{X}$.

Then the unique symmetric solution $\b{A} \in \mathbb{S}^n$ is:
\begin{equation}
\boxed{
\b{A}
=
\b{U}\,\widetilde{\b{A}}\,\b{U}^\top,
\qquad
\widetilde{A}_{ij}
=
\frac{\widetilde{\Delta}_{ij}}{\lambda_i+\lambda_j},
}
\label{equation_closed_form_Sylvester_solution_BW_metric}
\end{equation}
where:
\begin{align}\label{equation_Deltatilde_UT_Delta_U}
\widetilde{\b{\Delta}}
:=
\b{U}^\top \b{\Delta}\b{U}.
\end{align}

The solution $\b{A}$ can be restated as:
\begin{align}\label{equation_A_U_M_UT_Delta_U_UT}
\boxed{
\b{A}
=
\b{U}
\Big(
\b{M}\odot (\b{U}^\top \b{\Delta}\b{U})
\Big)
\b{U}^\top,
}
\end{align}
where $\odot$ denotes the Hadamard (elementwise) product and the $(i,j)$-th element of matrix $\b{M}$ is defined as:
\begin{align}
M_{ij} := \frac{1}{\lambda_i+\lambda_j}.
\end{align}
\end{lemma}

\begin{proof}
Because $\b{X} \in \mathbb{S}_{++}^n$, it admits an
eigendecomposition \cite{ghojogh2019eigenvalue}:
\[
\b{X} = \b{U}\b{\Lambda}\b{U}^\top,
\]
where $\b{U}$ is orthogonal and
\[
\b{\Lambda} = \text{diag}(\lambda_1,\dots,\lambda_n),
\qquad
\lambda_i > 0.
\]

We define:
\[
\widetilde{\b{A}} := \b{U}^\top \b{A}\b{U},
\qquad
\widetilde{\b{\Delta}} := \b{U}^\top \b{\Delta}\b{U}.
\]
Multiplying Eq.
\eqref{equation_Sylvester}
on the left by $\b{U}^\top$ and on the right by $\b{U}$
gives:
\[
\b{\Lambda}\widetilde{\b{A}}
+
\widetilde{\b{A}}\b{\Lambda}
=
\widetilde{\b{\Delta}}.
\]
Hence, we have entrywise:
\[
(\lambda_i+\lambda_j)\widetilde{A}_{ij}
=
\widetilde{\Delta}_{ij},
\qquad
\forall i,j.
\]
Since $\lambda_i+\lambda_j > 0$, we obtain uniquely:
\[
\widetilde{A}_{ij}
=
\frac{\widetilde{\Delta}_{ij}}{\lambda_i+\lambda_j}.
\]
Therefore, $\widetilde{\b{A}}$ is uniquely determined, and
thus so is $\b{A} = \b{U}\widetilde{\b{A}}\b{U}^\top$.

It remains to show symmetry. Since $\b{\Delta}$ is
symmetric, so is $\widetilde{\b{\Delta}}$. Therefore:
\begin{align}\label{equation_Atildeji_Deltatildeji_over_deltaj_deltai_in_proof}
\widetilde{A}_{ji}
=
\frac{\widetilde{\Delta}_{ji}}{\lambda_j+\lambda_i}
=
\frac{\widetilde{\Delta}_{ij}}{\lambda_i+\lambda_j}
=
\widetilde{A}_{ij}.
\end{align}
Hence, $\widetilde{\b{A}}$ is symmetric, and consequently
$\b{A}$ is symmetric as well.

Now, we define the matrix $\b{M}$ by:
\[
M_{ij} := \frac{1}{\lambda_i+\lambda_j}.
\]
Then, the entrywise formula above is exactly:
\[
\widetilde{\b{A}}
\overset{(\ref{equation_Atildeji_Deltatildeji_over_deltaj_deltai_in_proof})}{=}
\b{M}\odot \widetilde{\b{\Delta}}
\overset{(\ref{equation_Deltatilde_UT_Delta_U})}{=}
\b{M}\odot (\b{U}^\top \b{\Delta}\b{U}),
\]
so:
\[
\b{A}
= \b{U}\,\widetilde{\b{A}}\,\b{U}^\top =
\b{U}
\Big(
\b{M}\odot (\b{U}^\top \b{\Delta}\b{U})
\Big)
\b{U}^\top,
\]
which is Eq. (\ref{equation_A_U_M_UT_Delta_U_UT}).
\end{proof}

\begin{proposition}[Bures--Wasserstein metric tensor on SPD
manifold]\label{proposition_Bures_Wasserstein_metric_spd}
Let $\b{X} \in \mathbb{S}_{++}^n$. For every tangent vector
$\b{\Delta} \in T_{\b{X}}\mathbb{S}_{++}^n = \mathbb{S}^n$, let
$\mathcal{L}_{\b{X}}(\b{\Delta}) \in \mathbb{S}^n$ denote the
unique symmetric solution of Eq.
\eqref{equation_Sylvester}:
\begin{align}\label{equation_Sylvester_BW_metric}
\boxed{
\b{X}\,\mathcal{L}_{\b{X}}(\b{\Delta})
+
\mathcal{L}_{\b{X}}(\b{\Delta})\,\b{X}
=
\b{\Delta}.
}
\end{align}

According to Lemma \ref{lemma_Sylvester_equation_solution}, the unique symmetric solution $\mathcal{L}_{\b{X}}(\b{\Delta})$ can be computed
explicitly from the eigendecomposition of $\b{X}$:
\[
\b{X} = \b{U}\b{\Lambda}\b{U}^\top,
\qquad
\b{\Lambda} = \text{diag}(\lambda_1,\dots,\lambda_n),
\qquad
\lambda_i > 0,
\]
where the columns of $\b{U}$ are eigenvectors of $\b{X}$ and $\{\lambda_1,\dots,\lambda_n\}$ are the eigenvalues of $\b{X}$.

The solution $\mathcal{L}_{\b{X}}(\b{\Delta})$ is given by:
\begin{equation}
\boxed{
\mathcal{L}_{\b{X}}(\b{\Delta})
=
\b{U}
\Big(
\b{M} \odot (\b{U}^\top \b{\Delta}\b{U})
\Big)
\b{U}^\top,
}
\label{equation_Lambda_X_Delta_closed_form_in_BW_metric}
\end{equation}
where $\odot$ denotes the Hadamard (elementwise) product and the $(i,j)$-th element of matrix $\b{M}$ is defined as:
\[
M_{ij} := \frac{1}{\lambda_i+\lambda_j}.
\]

Then, the Bures--Wasserstein metric on the SPD
manifold is:
\begin{equation}
\boxed{
g_{\b{X}}^{BW}(\b{\Delta}_1,\b{\Delta}_2)
:=
\text{tr}\Big(
\mathcal{L}_{\b{X}}(\b{\Delta}_1)\,
\b{X}\,
\mathcal{L}_{\b{X}}(\b{\Delta}_2)
\Big).
}
\label{equation_metric_tensor_SPD_manifold_Bures_Wasserstein}
\end{equation}
Equivalently, using the Sylvester equation, it can be written
as:
\begin{equation}
\boxed{
g_{\b{X}}^{BW}(\b{\Delta}_1,\b{\Delta}_2)
=
\frac{1}{2}
\text{tr}\Big(
\mathcal{L}_{\b{X}}(\b{\Delta}_1)\,\b{\Delta}_2
\Big).
}
\label{equation_metric_tensor_SPD_manifold_Bures_Wasserstein_alt}
\end{equation}
\end{proposition}

\begin{proof}
By Lemma \ref{lemma_Sylvester_equation_solution}, for every $\b{\Delta} \in \mathbb{S}^n$, the
operator $\mathcal{L}_{\b{X}}(\b{\Delta})$ is well-defined and
belongs to $\mathbb{S}^n$.

We first prove the equivalence of Eqs.
\eqref{equation_metric_tensor_SPD_manifold_Bures_Wasserstein}
and
\eqref{equation_metric_tensor_SPD_manifold_Bures_Wasserstein_alt}.
From the Sylvester equation, Eq. (\ref{equation_Sylvester_BW_metric}), for $\b{\Delta}_2$, we have:
\begin{align}\label{equation_Delta2_X_LX_Delta2_LXDelta2_X_in_proof}
\b{\Delta}_2
=
\b{X}\,\mathcal{L}_{\b{X}}(\b{\Delta}_2)
+
\mathcal{L}_{\b{X}}(\b{\Delta}_2)\,\b{X}.
\end{align}

We have:
\begin{align}
\text{tr}\Big(
\mathcal{L}_{\b{X}}(\b{\Delta}_1)\,
&\mathcal{L}_{\b{X}}(\b{\Delta}_2)\,\b{X}
\Big) \nonumber\\
&\overset{(a)}{=}
\text{tr}\Big(
\mathcal{L}_{\b{X}}(\b{\Delta}_2) \b{X}\,\mathcal{L}_{\b{X}}(\b{\Delta}_1)
\Big) \nonumber\\
&\overset{(b)}{=}
\text{tr}\Big(
\big(\mathcal{L}_{\b{X}}(\b{\Delta}_2) \b{X}\,\mathcal{L}_{\b{X}}(\b{\Delta}_1) \big)^\top
\Big) \nonumber\\
&\overset{(c)}{=}
\text{tr}\Big(
\mathcal{L}_{\b{X}}(\b{\Delta}_1)^\top\,
\b{X}^\top\,\mathcal{L}_{\b{X}}(\b{\Delta}_2)^\top
\Big) \nonumber\\
&\overset{(d)}{=}
\text{tr}\Big(
\mathcal{L}_{\b{X}}(\b{\Delta}_1)\,
\b{X}\,\mathcal{L}_{\b{X}}(\b{\Delta}_2)
\Big), \label{equation_tr_LxDeltaLxDelta2_X_tr_LXDelta1_X_LXDelta2_in_proof}
\end{align}
where $(a)$ is by cyclic property of trace, $(b)$ is because trace of a matrix is equal to trace of its transpose, $(c)$ is because transpose reverses the order of matrix multiplication, and $(d)$ is because $\mathcal{L}_{\b{X}}(\b{\Delta}_1)$, $\b{X}$, and $\mathcal{L}_{\b{X}}(\b{\Delta}_2)$ are all symmetric matrices. 

Therefore:
\begin{align}
\text{tr}\Big(
\mathcal{L}_{\b{X}}&(\b{\Delta}_1)\,\b{\Delta}_2
\Big) \nonumber\\
&\overset{(\ref{equation_Delta2_X_LX_Delta2_LXDelta2_X_in_proof})}{=}
\text{tr}\Big(
\mathcal{L}_{\b{X}}(\b{\Delta}_1)\,
\b{X}\,\mathcal{L}_{\b{X}}(\b{\Delta}_2)
\Big)
\nonumber\\
&\quad\quad\quad\quad+
\text{tr}\Big(
\mathcal{L}_{\b{X}}(\b{\Delta}_1)\,
\mathcal{L}_{\b{X}}(\b{\Delta}_2)\,\b{X}
\Big) \nonumber\\
&\overset{(\ref{equation_tr_LxDeltaLxDelta2_X_tr_LXDelta1_X_LXDelta2_in_proof})}{=}
\text{tr}\Big(
\mathcal{L}_{\b{X}}(\b{\Delta}_1)\,
\b{X}\,\mathcal{L}_{\b{X}}(\b{\Delta}_2)
\Big)
\nonumber\\
&\quad\quad\quad\quad+
\text{tr}\Big(
\mathcal{L}_{\b{X}}(\b{\Delta}_1)\,
\b{X}\,\mathcal{L}_{\b{X}}(\b{\Delta}_2)
\Big) \nonumber\\
&=
2\,
\text{tr}\Big(
\mathcal{L}_{\b{X}}(\b{\Delta}_1)\,
\b{X}\,\mathcal{L}_{\b{X}}(\b{\Delta}_2)
\Big).  \label{equation_tr_LXDelta1_Delta2_2_tr_LXDelta1XLXDelta2_in_proof}
\end{align}
Hence:
\[
\text{tr}\Big(
\mathcal{L}_{\b{X}}(\b{\Delta}_1)\,
\b{X}\,
\mathcal{L}_{\b{X}}(\b{\Delta}_2)
\Big)
=
\frac{1}{2}
\text{tr}\Big(
\mathcal{L}_{\b{X}}(\b{\Delta}_1)\,\b{\Delta}_2
\Big).
\]

Now, we verify that $g_{\b{X}}^{BW}$ is an inner
product on $T_{\b{X}}\mathbb{S}_{++}^n$.

\textbf{Bilinearity:}
Because $\mathcal{L}_{\b{X}}(\cdot)$ is linear in $\b{\Delta}$
(the Sylvester equation is linear in the unknown and in the
right-hand side), and because trace is bilinear, we have:
\begin{align*}
g_{\b{X}}^{BW}(&a\b{\Delta}_1+b\b{\Delta}_2,\b{\Delta}_3)
\\
&\overset{(\ref{equation_metric_tensor_SPD_manifold_Bures_Wasserstein})}{=}
\text{tr}\Big(
\mathcal{L}_{\b{X}}(a\b{\Delta}_1+b\b{\Delta}_2)\,
\b{X}\,
\mathcal{L}_{\b{X}}(\b{\Delta}_3)
\Big) \\
&=
\text{tr}\Big(
(a\mathcal{L}_{\b{X}}(\b{\Delta}_1)+
b\mathcal{L}_{\b{X}}(\b{\Delta}_2))\,
\b{X}\,
\mathcal{L}_{\b{X}}(\b{\Delta}_3)
\Big) \\
&\overset{(\ref{equation_metric_tensor_SPD_manifold_Bures_Wasserstein})}{=}
a\,g_{\b{X}}^{BW}(\b{\Delta}_1,\b{\Delta}_3)
+
b\,g_{\b{X}}^{BW}(\b{\Delta}_2,\b{\Delta}_3).
\end{align*}
Similarly, it is linear in the second argument.

\textbf{Symmetry:}
Using Eq.
\eqref{equation_metric_tensor_SPD_manifold_Bures_Wasserstein_alt}, we have:
\begin{align*}
g_{\b{X}}^{BW}&(\b{\Delta}_1,\b{\Delta}_2)
\\
&\overset{(\ref{equation_metric_tensor_SPD_manifold_Bures_Wasserstein_alt})}{=}
\frac{1}{2}
\text{tr}\Big(
\mathcal{L}_{\b{X}}(\b{\Delta}_1)\,\b{\Delta}_2
\Big) \\
&\overset{(\ref{equation_tr_LXDelta1_Delta2_2_tr_LXDelta1XLXDelta2_in_proof})}{=}
\text{tr}\Big(
\mathcal{L}_{\b{X}}(\b{\Delta}_1)\,
\b{X}\,
\mathcal{L}_{\b{X}}(\b{\Delta}_2)
\Big) \\
&\overset{(a)}{=}
\text{tr}\Big(
\big(\mathcal{L}_{\b{X}}(\b{\Delta}_1)\,
\b{X}\,
\mathcal{L}_{\b{X}}(\b{\Delta}_2) \big)^\top
\Big) \\
&\overset{(b)}{=}
\text{tr}\Big(
\mathcal{L}_{\b{X}}(\b{\Delta}_2)^\top\,
\b{X}^\top\,
\mathcal{L}_{\b{X}}(\b{\Delta}_1)^\top
\Big) \\
&\overset{(c)}{=}
\text{tr}\Big(
\mathcal{L}_{\b{X}}(\b{\Delta}_2)\,
\b{X}\,
\mathcal{L}_{\b{X}}(\b{\Delta}_1)
\Big) \\
&\overset{(\ref{equation_metric_tensor_SPD_manifold_Bures_Wasserstein})}{=}
g_{\b{X}}^{BW}(\b{\Delta}_2,\b{\Delta}_1),
\end{align*}
where $(a)$ is because trace of a matrix is equal to trace of its transpose, $(b)$ is because transpose reverses the order of matrix multiplication, and $(c)$ is because $\mathcal{L}_{\b{X}}(\b{\Delta}_1)$, $\b{X}$, and $\mathcal{L}_{\b{X}}(\b{\Delta}_2)$ are all symmetric matrices.

\textbf{Positive definiteness:}
Let:
\[
\b{A} := \mathcal{L}_{\b{X}}(\b{\Delta}) \in \mathbb{S}^n.
\]
Then:
\[
\b{\Delta} = \b{X}\b{A} + \b{A}\b{X}.
\]
Using Eq.
\eqref{equation_metric_tensor_SPD_manifold_Bures_Wasserstein}, we have:
\[
g_{\b{X}}^{BW}(\b{\Delta},\b{\Delta})
=
\text{tr}(\b{A}\b{X}\b{A}).
\]
Because $\b{X}$ is SPD, there exists $\b{X}^{1/2}$ such that
$\b{X} = \b{X}^{1/2}\b{X}^{1/2}$. Hence:
\begin{align*}
\text{tr}(\b{A}&\b{X}\b{A})
=
\text{tr}(\b{A}\b{X}^{1/2}\b{X}^{1/2}\b{A})
\\
&=
\text{tr}\Big(
(\b{X}^{1/2}\b{A})^\top
(\b{X}^{1/2}\b{A})
\Big)
=
\Vert\b{X}^{1/2}\b{A}\Vert_F^2
\geq 0,
\end{align*}
where we used $\b{A}^\top = \b{A}$.

If:
\[
g_{\b{X}}^{BW}(\b{\Delta},\b{\Delta}) = 0,
\]
then:
\[
\Vert\b{X}^{1/2}\b{A}\Vert_F^2 = 0
\quad \Longrightarrow \quad
\b{X}^{1/2}\b{A} = \b{0}.
\]
Since $\b{X}^{1/2}$ is invertible, this implies $\b{A}=\b{0}$.
Therefore:
\[
\b{\Delta} = \b{X}\b{A} + \b{A}\b{X} = \b{0}.
\]
Hence:
\[
g_{\b{X}}^{BW}(\b{\Delta},\b{\Delta}) > 0,
\qquad \forall \b{\Delta} \neq \b{0}.
\]

Therefore,
Eq.
\eqref{equation_metric_tensor_SPD_manifold_Bures_Wasserstein}
defines a Riemannian metric on $\mathbb{S}_{++}^n$.
\end{proof}


\begin{remark}[Interpretation of the Bures--Wasserstein metric]
The Bures--Wasserstein metric does not measure a tangent
vector $\b{\Delta} \in T_{\b{X}}\mathbb{S}_{++}^n$ directly.
Instead, it first represents $\b{\Delta}$ through the unique
symmetric matrix $\b{A}=\mathcal{L}_{\b{X}}(\b{\Delta})$
satisfying:
\begin{align}\label{equation_XA_AX_Delta_in_remark}
\b{X}\b{A}+\b{A}\b{X}=\b{\Delta}.
\end{align}

The matrix $\b{A}$ is interpreted through the factorization:
\begin{align*}
\b{X}=\b{Y}\b{Y}^\top.
\end{align*}
According to this factorization, the matrix $\b{Y}$ is the square-root factor of $\b{X}$.

If $\b{Y}(t)$ is a smooth curve of
matrix factors and:
\[
\dot{\b{Y}}(t)=\b{A}(t)\b{Y}(t),
\]
then for $\b{X}(t)=\b{Y}(t)\b{Y}(t)^\top$, we have:
\[
\dot{\b{X}}(t)
=
\b{A}(t)\b{X}(t)+\b{X}(t)\b{A}(t)^\top.
\]
For symmetric $\b{A}(t)$, this becomes:
\begin{align}\label{equation_Xdot_AX_XA_in_remark}
\dot{\b{X}}(t)
=
\b{A}(t)\b{X}(t)+\b{X}(t)\b{A}(t).
\end{align}

Comparing Eqs. (\ref{equation_XA_AX_Delta_in_remark}) and (\ref{equation_Xdot_AX_XA_in_remark}) shows that the matrix $\b{A}$ describes the infinitesimal motion
of the square-root factor $\b{Y}$ that generates the tangent
vector $\b{\Delta}=\dot{\b{X}}$.

Hence, $\b{A}$ can be interpreted as the velocity in the
square-root representation that generates the tangent vector
$\b{\Delta}$. Therefore, the Bures--Wasserstein metric is
naturally related to matrix square roots, unlike the
affine-invariant metric, which is expressed directly in
terms of $\b{X}^{-1}$.

\end{remark}



\subsubsection{Levi-Civita Connection in SPD Manifold}

Recall from Section \ref{section_tangent_normal_spd} that for every
$\b{X} \in \mathbb{S}_{++}^n$, the tangent space is:
\[
T_{\b{X}}\mathbb{S}_{++}^n = \mathbb{S}^n.
\]
Therefore, tangent vectors and tangent vector fields on the
SPD manifold are symmetric matrices.

In this subsection, we derive the Levi-Civita connection on
the SPD manifold for two important metrics, i.e., the Euclidean metric and the affine-invariant metric.

\begin{proposition}[Levi-Civita connection in SPD manifold under the Euclidean metric]\label{proposition_Levi_Civita_connection_spd_euclidean_metric}
Let
$\b{\Delta}_2 : \mathbb{S}_{++}^n \to \mathbb{S}^n$
be a smooth tangent vector field on the SPD manifold.
Under the Euclidean metric, in Proposition \ref{proposition_euclidean_metric_spd}, the Levi-Civita connection is:
\begin{equation}\label{equation_Levi_Civita_connection_SPD_Euclidean}
\boxed{
(\nabla_{\b{\Delta}_1}^{E}\b{\Delta}_2)(\b{X})
=
D\b{\Delta}_2(\b{X})[\b{\Delta}_1],
}
\end{equation}
for $\b{\Delta}_1\in T_{\b{X}}\mathbb{S}_{++}^n$.
\end{proposition}

\begin{proof}
The manifold $\mathbb{S}_{++}^n$ is an open subset of the
vector space $\mathbb{S}^n$. Hence, under the Euclidean
metric inherited from the ambient vector space, the natural
connection is the ordinary directional derivative.

We first show that
$D\b{\Delta}_2(\b{X})[\b{\Delta}_1]$
is tangent to the SPD manifold. Since
$\b{\Delta}_2(\b{X}) \in \mathbb{S}^n$
for every $\b{X} \in \mathbb{S}_{++}^n$, we have:
\[
\b{\Delta}_2(\b{X})^\top = \b{\Delta}_2(\b{X}).
\]
Differentiating both sides in the direction $\b{\Delta}_1$ gives:
\[
D(\b{\Delta}_2^\top)(\b{X})[\b{\Delta}_1]
=
D\b{\Delta}_2(\b{X})[\b{\Delta}_1].
\]
Because transpose is linear, we have:
\[
D(\b{\Delta}_2^\top)(\b{X})[\b{\Delta}_1]
=
\big(D\b{\Delta}_2(\b{X})[\b{\Delta}_1]\big)^\top.
\]
Therefore:
\[
\big(D\b{\Delta}_2(\b{X})[\b{\Delta}_1]\big)^\top
=
D\b{\Delta}_2(\b{X})[\b{\Delta}_1],
\]
so
$D\b{\Delta}_2(\b{X})[\b{\Delta}_1] \in \mathbb{S}^n
= T_{\b{X}}\mathbb{S}_{++}^n$.

Now, because the metric is the constant Euclidean inner
product on the vector space $\mathbb{S}^n$, the ordinary
directional derivative is torsion-free and metric-compatible.
Hence it is the Levi-Civita connection.

Therefore:
\[
(\nabla_{\b{\Delta}_1}^{E}\b{\Delta}_2)(\b{X})
=
D\b{\Delta}_2(\b{X})[\b{\Delta}_1].
\]
This proves
Eq. \eqref{equation_Levi_Civita_connection_SPD_Euclidean}.
\end{proof}

For the affine-invariant metric, the connection is not the
ordinary directional derivative because the metric depends
on the point $\b{X}$. We first need the derivative of the
matrix inverse.

\begin{lemma}[Derivative of the matrix inverse]
Let
$\b{X} \in \mathbb{S}_{++}^n$
and let
$\b{\Delta}_1 \in \mathbb{S}^n$.
Then:
\begin{equation}
D(\b{X}^{-1})[\b{\Delta}_1]
=
-\b{X}^{-1}\b{\Delta}_1\b{X}^{-1}.
\label{equation_derivative_inverse_SPD}
\end{equation}
\end{lemma}

\begin{proof}
We use the identity:
\[
\b{X}\b{X}^{-1} = \b{I}.
\]
Differentiating both sides in the direction $\b{\Delta}$
gives:
\[
D(\b{X}\b{X}^{-1})[\b{\Delta}_1] = \b{0}.
\]
Applying the product rule gives:
\[
\b{\Delta}_1\b{X}^{-1}
+
\b{X}\,D(\b{X}^{-1})[\b{\Delta}_1]
=
\b{0}.
\]
Multiplying from the left by $\b{X}^{-1}$ gives:
\[
\b{X}^{-1}\b{\Delta}_1\b{X}^{-1}
+
D(\b{X}^{-1})[\b{\Delta}_1]
=
\b{0}.
\]
Therefore:
\[
D(\b{X}^{-1})[\b{\Delta}_1]
=
-\b{X}^{-1}\b{\Delta}_1\b{X}^{-1}.
\]
\end{proof}

\begin{proposition}[Levi-Civita connection in SPD manifold under the affine-invariant metric]\label{proposition_Levi_Civita_connection_spd_affine_invariant_metric}
Let
$\b{\Delta}_2 : \mathbb{S}_{++}^n \to \mathbb{S}^n$
be a smooth tangent vector field on the SPD manifold.
Under the affine-invariant metric, in Proposition \ref{proposition_affine_invariant_metric_spd}, the Levi-Civita connection is:
\begin{equation}\label{equation_Levi_Civita_connection_SPD_affine_invariant}
\boxed{
\begin{aligned}
(\nabla_{\b{\Delta}_1}^{AI}&\b{\Delta}_2)(\b{X})
=
D\b{\Delta}_2(\b{X})[\b{\Delta}_1]
\\
&-
\frac{1}{2}
\Big(
\b{\Delta}_1\b{X}^{-1}\b{\Delta}_2(\b{X})
+
\b{\Delta}_2(\b{X})\b{X}^{-1}\b{\Delta}_1
\Big).
\end{aligned}
}
\end{equation}
\end{proposition}

\begin{proof}
Define the candidate connection:
\begin{align*}
(\widetilde{\nabla}_{\b{\Delta}_1}&\b{\Delta}_2)(\b{X})
:=
D\b{\Delta}_2(\b{X})[\b{\Delta}_1]
\\
&-
\frac{1}{2}
\Big(
\b{\Delta}_1\b{X}^{-1}\b{\Delta}_2(\b{X})
+
\b{\Delta}_2(\b{X})\b{X}^{-1}\b{\Delta}_1
\Big).
\end{align*}
We prove that this connection is:
\begin{enumerate}
\item tangent-valued,
\item torsion-free,
\item metric-compatible.
\end{enumerate}
Hence, by uniqueness of the Levi-Civita connection,
$\widetilde{\nabla}=\nabla^{AI}$.

\textbf{Step 1: Tangent-valuedness.}

Since
$\b{\Delta}_2(\b{X}) \in \mathbb{S}^n$,
we have already shown in the proof of the Euclidean case
that
$D\b{\Delta}_2(\b{X})[\b{\Delta}_1] \in \mathbb{S}^n$.
Also, because
$\b{\Delta}_1$, $\b{X}^{-1}$, and $\b{\Delta}_2(\b{X})$
are symmetric, we have:
\[
\big(\b{\Delta}_1\b{X}^{-1}\b{\Delta}_2(\b{X})\big)^\top
=
\b{\Delta}_2(\b{X})\b{X}^{-1}\b{\Delta}_1.
\]
Hence:
\begin{align*}
\Big(
\b{\Delta}_1\b{X}^{-1}&\b{\Delta}_2(\b{X})
+
\b{\Delta}_2(\b{X})\b{X}^{-1}\b{\Delta}_1
\Big)^\top
\\
&=
\b{\Delta}_1\b{X}^{-1}\b{\Delta}_2(\b{X})
+
\b{\Delta}_2(\b{X})\b{X}^{-1}\b{\Delta}_1.
\end{align*}
So the correction term is symmetric. Therefore,
\[
(\widetilde{\nabla}_{\b{\Delta}_1}\b{\Delta}_2)(\b{X})
\in \mathbb{S}^n
=
T_{\b{X}}\mathbb{S}_{++}^n.
\]

\textbf{Step 2: Torsion-free property.}

Let $\b{\Delta}_1$ and $\b{\Delta}_2$ be two smooth tangent
vector fields. Then:
\begin{align*}
&(\widetilde{\nabla}_{\b{\Delta}_2}\b{\Delta}_1)(\b{X})
-
(\widetilde{\nabla}_{\b{\Delta}_1}\b{\Delta}_2)(\b{X})
\\
&=
D\b{\Delta}_1[\b{\Delta}_2]
-
\frac{1}{2}
\Big(
\b{\Delta}_2\,\b{X}^{-1}\b{\Delta}_1
+
\b{\Delta}_1\,\b{X}^{-1}\b{\Delta}_2
\Big) \\
&\quad -
D\b{\Delta}_2[\b{\Delta}_1]
+
\frac{1}{2}
\Big(
\b{\Delta}_1\,\b{X}^{-1}\b{\Delta}_2
+
\b{\Delta}_2\,\b{X}^{-1}\b{Z}
\Big) \\
&=
D\b{\Delta}_1[\b{\Delta}_2]
-
D\b{\Delta}_2[\b{\Delta}_1].
\end{align*}
Thus, the correction terms cancel, and we obtain:
\[
\widetilde{\nabla}_{\b{\Delta}_2}\b{\Delta}_1
-
\widetilde{\nabla}_{\b{\Delta}_1}\b{\Delta}_2
=
[\b{\Delta}_2,\b{\Delta}_1],
\]
which is exactly the torsion-free condition.

\textbf{Step 3: Metric compatibility.}

We must show that for any smooth tangent vector fields
$\b{\Delta}_1, \b{\Delta}_2, \b{\Delta}_3$, we have:
\begin{align*}
\b{\Delta}_1
\big(
g_{\b{X}}^{AI}(\b{\Delta}_2,\b{\Delta}_3)
\big)
=\,
&g_{\b{X}}^{AI}
\big(
\widetilde{\nabla}_{\b{\Delta}_1}\b{\Delta}_2,
\b{\Delta}_3
\big)
\\
&+
g_{\b{X}}^{AI}
\big(
\b{\Delta}_2,
\widetilde{\nabla}_{\b{\Delta}_1}\b{\Delta}_3
\big).
\end{align*}

Because:
\[
g_{\b{X}}^{AI}(\b{\Delta}_2,\b{\Delta}_3)
=
\text{tr}(\b{X}^{-1}\b{\Delta}_2\b{X}^{-1}\b{\Delta}_3),
\]
differentiating in the direction $\b{\Delta}$ gives:
\begin{align}
D\big(
g_{\b{X}}^{AI}(\b{\Delta}_2,\b{\Delta}_3)
\big)[\b{\Delta}]
&=
\text{tr}\Big(
D(\b{X}^{-1})[\b{\Delta}]\,
\b{\Delta}_2\,
\b{X}^{-1}\b{\Delta}_3
\Big)
\nonumber\\
&+
\text{tr}\Big(
\b{X}^{-1}D\b{\Delta}_2[\b{\Delta}]\,
\b{X}^{-1}\b{\Delta}_3
\Big)
\nonumber\\
&+
\text{tr}\Big(
\b{X}^{-1}\b{\Delta}_2\,
D(\b{X}^{-1})[\b{\Delta}_1]\,
\b{\Delta}_3
\Big)
\nonumber\\
&+
\text{tr}\Big(
\b{X}^{-1}\b{\Delta}_2\,
\b{X}^{-1}D\b{\Delta}_3[\b{\Delta}_1]
\Big).
\label{equation_metric_compatibility_SPD_AI_start}
\end{align}

By Eq. (\ref{equation_derivative_inverse_SPD}), we have:
\[
D(\b{X}^{-1})[\b{\Delta}_1]
=
-\b{X}^{-1}\b{\Delta}_1\b{X}^{-1}.
\]
Substituting into
Eq. \eqref{equation_metric_compatibility_SPD_AI_start},
we get:
\begin{align}
D\big(
g_{\b{X}}^{AI}&(\b{\Delta}_2,\b{\Delta}_3)
\big)[\b{\Delta}_1] \nonumber\\
&=
-\text{tr}\Big(
\b{X}^{-1}\b{\Delta}_1\b{X}^{-1}
\b{\Delta}_2\b{X}^{-1}\b{\Delta}_3
\Big)
\nonumber\\
&+
\text{tr}\Big(
\b{X}^{-1}D\b{\Delta}_2[\b{\Delta}_1]
\b{X}^{-1}\b{\Delta}_3
\Big)
\nonumber\\
&-
\text{tr}\Big(
\b{X}^{-1}\b{\Delta}_2\b{X}^{-1}
\b{\Delta}_1\b{X}^{-1}\b{\Delta}_3
\Big)
\nonumber\\
&+
\text{tr}\Big(
\b{X}^{-1}\b{\Delta}_2\b{X}^{-1}
D\b{Z}[\b{\Delta}_1]
\Big).
\label{equation_metric_compatibility_SPD_AI_expand}
\end{align}

Now we calculate the right-hand side:
\begin{align*}
g_{\b{X}}^{AI}
\big(
\widetilde{\nabla}_{\b{\Delta}_1}\b{\Delta}_2,
\b{Z}
\big)
&=
\text{tr}\Big(
\b{X}^{-1}
\widetilde{\nabla}_{\b{\Delta}_1}\b{\Delta}_2
\b{X}^{-1}\b{Z}
\Big) \\
&=
\text{tr}\Big(
\b{X}^{-1}D\b{\Delta}_2[\b{\Delta}_1]
\b{X}^{-1}\b{\Delta}_3
\Big)
\\
&-
\frac{1}{2}
\text{tr}\Big(
\b{X}^{-1}\b{\Delta}_1\b{X}^{-1}
\b{\Delta}_2\b{X}^{-1}\b{\Delta}_3
\Big)
\\
&-
\frac{1}{2}
\text{tr}\Big(
\b{X}^{-1}\b{\Delta}_2\b{X}^{-1}
\b{\Delta}_1\b{X}^{-1}\b{\Delta}_3
\Big).
\end{align*}
Similarly, we have:
\begin{align*}
g_{\b{X}}^{AI}
\big(
\b{\Delta}_2,
\widetilde{\nabla}_{\b{\Delta}_1}\b{\Delta}_3
\big)
&=
\text{tr}\Big(
\b{X}^{-1}\b{\Delta}_2\b{X}^{-1}
D\b{Z}[\b{\Delta}_1]
\Big)
\\
&-
\frac{1}{2}
\text{tr}\Big(
\b{X}^{-1}\b{\Delta}_2\b{X}^{-1}
\b{\Delta}_1\b{X}^{-1}\b{\Delta}_3
\Big)
\\
&-
\frac{1}{2}
\text{tr}\Big(
\b{X}^{-1}\b{\Delta}_2\b{X}^{-1}
\b{\Delta}_3\b{X}^{-1}\b{\Delta}_1
\Big).
\end{align*}
Using cyclic property of trace, we have:
\begin{align*}
\text{tr}\Big(
\b{X}^{-1}\b{\Delta}_2\b{X}^{-1}
&\b{\Delta}_3\b{X}^{-1}\b{\Delta}_1
\Big)
\\
&=
\text{tr}\Big(
\b{X}^{-1}\b{\Delta}_1\b{X}^{-1}
\b{\Delta}_2\b{X}^{-1}\b{\Delta}_3
\Big).
\end{align*}
Therefore, after adding the two expressions, we obtain:
\begin{align*}
g_{\b{X}}^{AI}
\big(
\widetilde{\nabla}_{\b{\Delta}_1}\b{\Delta}_2,\,
&\b{\Delta}_3
\big)
+
g_{\b{X}}^{AI}
\big(
\b{\Delta}_2,
\widetilde{\nabla}_{\b{\Delta}_1}\b{\Delta}_3
\big) \\
&=
\text{tr}\Big(
\b{X}^{-1}D\b{\Delta}_2[\b{\Delta}_1]
\b{X}^{-1}\b{\Delta}_3
\Big)
\\
&\quad+
\text{tr}\Big(
\b{X}^{-1}\b{\Delta}_2\b{X}^{-1}
D\b{\Delta}_3[\b{\Delta}_1]
\Big)
\\
&\quad-
\text{tr}\Big(
\b{X}^{-1}\b{\Delta}_1\b{X}^{-1}
\b{\Delta}_2\b{X}^{-1}\b{\Delta}_3
\Big)
\\
&\quad-
\text{tr}\Big(
\b{X}^{-1}\b{\Delta}_2\b{X}^{-1}
\b{\Delta}_1\b{X}^{-1}\b{\Delta}_3
\Big).
\end{align*}
This is exactly
Eq. \eqref{equation_metric_compatibility_SPD_AI_expand}.
Hence:
\begin{align*}
D\big(
g_{\b{X}}^{AI}(\b{\Delta}_2,\b{\Delta}_3)
\big)[\b{\Delta}]
=\,
&g_{\b{X}}^{AI}
\big(
\widetilde{\nabla}_{\b{\Delta}_1}\b{\Delta}_2,
\b{\Delta}_3
\big)
\\
&\quad+
g_{\b{X}}^{AI}
\big(
\b{\Delta}_2,
\widetilde{\nabla}_{\b{\Delta}_1}\b{\Delta}_3
\big).
\end{align*}
So, the connection is metric-compatible.

Since $\widetilde{\nabla}$ is torsion-free and
metric-compatible, by uniqueness of the Levi-Civita
connection, it is the Levi-Civita connection of the
affine-invariant metric. Therefore:
\begin{align*}
(\nabla_{\b{\Delta}_1}^{AI}&\b{\Delta}_2)(\b{X})
=
D\b{\Delta}_2(\b{X})[\b{\Delta}_1]
\\
&-
\frac{1}{2}
\Big(
\b{\Delta}_1\b{X}^{-1}\b{\Delta}_2(\b{X})
+
\b{\Delta}_2(\b{X})\b{X}^{-1}\b{\Delta}_1
\Big).
\end{align*}
This proves
Eq. \eqref{equation_Levi_Civita_connection_SPD_affine_invariant}.
\end{proof}

\begin{remark}[Comparison of Levi-Civita connections in SPD manifold]
Under the Euclidean metric, the SPD manifold is an open
subset of the flat vector space $\mathbb{S}^n$. Therefore,
the Levi-Civita connection is simply the ordinary
directional derivative:
\[
(\nabla_{\b{\Delta}_1}^{E}\b{\Delta}_2)(\b{X})
=
D\b{\Delta}_2(\b{X})[\b{\Delta}_1].
\]

Under the affine-invariant metric, the connection contains
an additional correction term:
\begin{align*}
(\nabla_{\b{\Delta}_1}^{AI}&\b{\Delta}_2)(\b{X})
=
D\b{\Delta}_2(\b{X})[\b{\Delta}_1]
\\
&-
\frac{1}{2}
\Big(
\b{\Delta}_1\b{X}^{-1}\b{\Delta}_2(\b{X})
+
\b{\Delta}_2(\b{X})\b{X}^{-1}\b{\Delta}_1
\Big).
\end{align*}
The additional term appears because the affine-invariant metric
depends on the base point $\b{X}$, unlike the Euclidean
metric. Therefore, the geometry of the SPD manifold under
the affine-invariant metric is not flat in the Euclidean
sense.
\end{remark}

\subsubsection{Riemannian Gradient in SPD Manifold}

\begin{definition}[Smooth local extension on the ambient space for SPD manifold]
Let $f : \mathbb{S}_{++}^n \to \mathbb{R}$ be a smooth function, and let
$\b{X} \in \mathbb{S}_{++}^n$. A \textbf{smooth local extension} of $f$ around
$\b{X}$ is a smooth function:
\begin{align}
\bar{f} : U \subset \mathbb{R}^{n \times n} \to \mathbb{R},
\end{align}
defined on an open neighborhood $U$ of $\b{X}$, such that:
\begin{align}\label{equation_smooth_extension_SPD}
\boxed{
\bar{f}(\b{Y}) = f(\b{Y}), \qquad \forall \b{Y} \in U \cap \mathbb{S}_{++}^n.
}
\end{align}
In other words, $\bar{f}$ agrees with $f$ on the points of the SPD manifold
near $\b{X}$, but it is defined on an open set of the ambient Euclidean space
$\mathbb{R}^{n \times n}$ so that its Euclidean gradient can be computed.
\end{definition}

Because $\mathbb{S}_{++}^n$ is an open subset of $\mathbb{S}^n$, one may also
view $f$ locally as a smooth function defined on an open subset of the vector
space $\mathbb{S}^n$. However, for deriving matrix formulas using ambient
Euclidean calculus, it is convenient to use a smooth local extension
$\bar{f}$ defined on an open neighborhood in $\mathbb{R}^{n\times n}$.

Recall from Section \ref{section_tangent_normal_spd} that for every point
$\b{X} \in \mathbb{S}_{++}^n$, the tangent space is:
\[
T_{\b{X}}\mathbb{S}_{++}^n = \mathbb{S}^n.
\]
Hence, every tangent vector on the SPD manifold is a
symmetric matrix.

Also, recall from Proposition \ref{proposition_characterizing_identity_of_gradient} that the Riemannian
gradient is characterized by the identity:
\[
g_{\b{X}}(\operatorname{grad} f(\b{X}), \b{\Delta})
=
\text{tr}\big((\nabla \bar f(\b{X}))^\top \b{\Delta}\big),
\,\,
\forall \b{\Delta} \in T_{\b{X}}\mathbb{S}_{++}^n,
\]
where $\bar f$ is a smooth local extension of
$f : \mathbb{S}_{++}^n \to \mathbb{R}$ to an open neighborhood
of $\mathbb{S}_{++}^n$ in $\mathbb{R}^{n\times n}$, and
$\nabla \bar f(\b{X})$ is the Euclidean gradient in the
ambient space.

Because tangent vectors on $\mathbb{S}_{++}^n$ are symmetric,
only the symmetric part of the ambient Euclidean gradient
contributes to directional derivatives. We first state this
fact explicitly.

\begin{lemma}[Only the symmetric part of the ambient gradient contributes on SPD manifold]
Let $f : \mathbb{S}_{++}^n \to \mathbb{R}$ be a smooth
function and let $\bar f$ be a smooth local extension to an
open neighborhood in $\mathbb{R}^{n\times n}$. For every
$\b{X} \in \mathbb{S}_{++}^n$ and every
$\b{\Delta} \in T_{\b{X}}\mathbb{S}_{++}^n = \mathbb{S}^n$, we have:
\begin{align}\label{equation_DfXDelta_tr_sym_nablafX_Delta}
D\bar f(\b{X})[\b{\Delta}]
=
\text{tr}\big((\nabla \bar f(\b{X}))^\top \b{\Delta}\big)
=
\text{tr}\big(\operatorname{sym}(\nabla \bar f(\b{X}))\,\b{\Delta}\big),
\end{align}
where $\operatorname{sym}(\cdot)$ is defined in Eq. (\ref{equation_sym_skew_expressions}).
\end{lemma}

\begin{proof}
Because $\b{\Delta} \in \mathbb{S}^n$, we have
$\b{\Delta}^\top = \b{\Delta}$. Therefore:
\begin{align*}
\text{tr}\big((\nabla \bar f(\b{X}))^\top \b{\Delta}\big)
=\,
&\frac{1}{2}
\text{tr}\big((\nabla \bar f(\b{X}))^\top \b{\Delta}\big)
\\
&+
\frac{1}{2}
\text{tr}\big((\nabla \bar f(\b{X}))^\top \b{\Delta}\big).
\end{align*}
For the second term, using the fact that the trace of a
matrix equals the trace of its transpose, we get:
\begin{align*}
\text{tr}\big((\nabla \bar f(\b{X}))^\top \b{\Delta}\big)
&=
\text{tr}\Big(
\big((\nabla \bar f(\b{X}))^\top \b{\Delta}\big)^\top
\Big) \\
&=
\text{tr}\big(\b{\Delta}^\top \nabla \bar f(\b{X})\big) \\
&\overset{(a)}{=}
\text{tr}\big(\b{\Delta}\,\nabla \bar f(\b{X})\big) \\
&\overset{(b)}{=}
\text{tr}\big(\nabla \bar f(\b{X})\,\b{\Delta}\big),
\end{align*}
where $(a)$ is because of $\b{\Delta} = \b{\Delta}^\top$ and $(b)$ is because of cyclic property of trace.
Hence:
\begin{align*}
\text{tr}\big((\nabla &\bar f(\b{X}))^\top \b{\Delta}\big)
\\
&=
\frac{1}{2}
\text{tr}\big((\nabla \bar f(\b{X}))^\top \b{\Delta}\big)
+
\frac{1}{2}
\text{tr}\big(\nabla \bar f(\b{X})\,\b{\Delta}\big) \\
&=
\text{tr}\left(
\frac{\nabla \bar f(\b{X}) + (\nabla \bar f(\b{X}))^\top}{2}
\,\b{\Delta}
\right) \\
&=
\text{tr}\big(\operatorname{sym}(\nabla \bar f(\b{X}))\,\b{\Delta}\big).
\end{align*}
This proves the result.
\end{proof}

Now, we derive the Riemannian gradients corresponding to
the Euclidean and affine-invariant metrics, introduced in Section \ref{section_metrics_SPD}.

\begin{proposition}[Riemannian gradient in SPD manifold under the Euclidean metric]\label{proposition_Riemannian_gradient_SPD_Euclidean_metric}
Let $f : \mathbb{S}_{++}^n \to \mathbb{R}$ be a smooth
function and let $\bar f$ be a smooth local extension.
Under the Euclidean (Frobenius) metric on
$\mathbb{S}_{++}^n$, the Riemannian gradient is:
\begin{align}\label{equation_Riemannian_gradient_SPD_Euclidean_metric}
\boxed{
\operatorname{grad}^{E} f(\b{X})
=
\operatorname{sym}(\nabla \bar f(\b{X})).
}
\end{align}
\end{proposition}

\begin{proof}
By Proposition \ref{proposition_euclidean_metric_spd}, the Euclidean metric on
$\mathbb{S}_{++}^n$ is:
\[
g_{\b{X}}^{E}(\b{\Delta}_1,\b{\Delta}_2)
=
\text{tr}(\b{\Delta}_1\b{\Delta}_2),
\qquad
\b{\Delta}_1,\b{\Delta}_2 \in T_{\b{X}}\mathbb{S}_{++}^n.
\]
Let:
\[
\b{G}
:=
\operatorname{grad}^{E} f(\b{X}).
\]
By the defining property of the Riemannian gradient, for
every $\b{\Delta} \in T_{\b{X}}\mathbb{S}_{++}^n$, we must have:
\begin{align}
g_{\b{X}}^{E}(\b{G},\b{\Delta})
&=
\text{tr}\big((\nabla \bar f(\b{X}))^\top \b{\Delta}\big) \nonumber \\
&\overset{(\ref{equation_DfXDelta_tr_sym_nablafX_Delta})}{=} \text{tr}\big(\operatorname{sym}(\nabla \bar f(\b{X}))\,\b{\Delta}\big). \label{equation_gXGDelta_tr_sym_nablafX_Delta_in_proof}
\end{align}

We also have:
\begin{align}\label{equation_gXGDelta_tr_G_Delta_in_proof}
g_{\b{X}}^{E}(\b{G},\b{\Delta}) \overset{(\ref{equation_g_inner_product})}{=} \langle \b{G},\b{\Delta} \rangle_F \overset{(\ref{equation_Frobenius_inner_product})}{=} \text{tr}(\b{G}^\top\b{\Delta}) \overset{(a)}{=} \text{tr}(\b{G}\b{\Delta}),
\end{align}
where $(a)$ is because $\b{G}$ is symmetric (because Riemannian gradient in SPD is in the tangent space of SPD manifold, so it should be a SPD matrix). 

According to Eqs. (\ref{equation_gXGDelta_tr_sym_nablafX_Delta_in_proof}) and (\ref{equation_gXGDelta_tr_G_Delta_in_proof}), we have:
\[
\text{tr}(\b{G}\b{\Delta})
=
\text{tr}\big(\operatorname{sym}(\nabla \bar f(\b{X}))\,\b{\Delta}\big),
\qquad
\forall \b{\Delta} \in \mathbb{S}^n.
\]
Therefore, by uniqueness of the Riemannian gradient, we have:
\[
\b{G}
=
\operatorname{sym}(\nabla \bar f(\b{X})).
\]
Hence:
\[
\operatorname{grad}^{E} f(\b{X})
=
\operatorname{sym}(\nabla \bar f(\b{X})).
\]
\end{proof}

\begin{proposition}[Riemannian gradient in SPD manifold under the affine-invariant metric]\label{proposition_Riemannian_gradient_SPD_affine_invariant_metric}
Let $f : \mathbb{S}_{++}^n \to \mathbb{R}$ be a smooth
function and let $\bar f$ be a smooth local extension.
Under the affine-invariant metric on $\mathbb{S}_{++}^n$,
the Riemannian gradient is:
\begin{equation}
\boxed{
\begin{aligned}
\operatorname{grad}^{AI} f(\b{X})
&=
\b{X}\,\operatorname{sym}(\nabla \bar f(\b{X}))\,\b{X} \\
&\overset{(\ref{equation_sym_skew_expressions})}{=} \frac{1}{2} \b{X} \Big(\nabla \bar f(\b{X}) + \big(\nabla \bar f(\b{X})\big)^\top\Big)\,\b{X}.
\end{aligned}
}
\end{equation}
\end{proposition}

\begin{proof}
By Proposition \ref{proposition_affine_invariant_metric_spd}, the affine-invariant metric is:
\[
g_{\b{X}}^{AI}(\b{\Delta}_1,\b{\Delta}_2)
=
\text{tr}(\b{X}^{-1}\b{\Delta}_1\b{X}^{-1}\b{\Delta}_2).
\]
Let:
\[
\b{G}
:=
\operatorname{grad}^{AI} f(\b{X}).
\]
Then, for every $\b{\Delta} \in T_{\b{X}}\mathbb{S}_{++}^n$,
the defining property of the gradient gives:
\begin{align*}
&\text{tr}(\b{X}^{-1}\b{G}\b{X}^{-1}\b{\Delta})
=
\text{tr}\big((\nabla \bar f(\b{X}))^\top \b{\Delta}\big)
\\
&~~~~~~~~~~~\overset{(\ref{equation_DfXDelta_tr_sym_nablafX_Delta})}{=}
\text{tr}\big(\operatorname{sym}(\nabla \bar f(\b{X}))\,\b{\Delta}\big),
\qquad
\forall \b{\Delta} \in \mathbb{S}^n.
\end{align*}

Since this holds for every symmetric $\b{\Delta}$, we obtain:
\[
\b{X}^{-1}\b{G}\b{X}^{-1}
=
\operatorname{sym}(\nabla \bar f(\b{X})).
\]
Multiplying from left and right by $\b{X}$ gives:
\[
\b{G}
=
\b{X}\,\operatorname{sym}(\nabla \bar f(\b{X}))\,\b{X}.
\]
Therefore:
\[
\operatorname{grad}^{AI} f(\b{X})
=
\b{X}\,\operatorname{sym}(\nabla \bar f(\b{X}))\,\b{X}.
\]
\end{proof}

\begin{remark}[Comparison of the Riemannian gradients in SPD manifold]
The formula of the Riemannian gradient depends on the
choice of metric on $\mathbb{S}_{++}^n$.

Under the Euclidean metric, the gradient is simply the
symmetric part of the ambient Euclidean gradient:
\[
\operatorname{grad}^{E} f(\b{X})
=
\operatorname{sym}(\nabla \bar f(\b{X})).
\]

Under the affine-invariant metric, the geometry weights the
gradient by the point $\b{X}$ from both sides:
\[
\operatorname{grad}^{AI} f(\b{X})
=
\b{X}\,\operatorname{sym}(\nabla \bar f(\b{X}))\,\b{X}.
\]



Therefore, although the same scalar function $f$ is being
optimized, the steepest-ascent direction depends on the
geometry chosen on the SPD manifold.
\end{remark}

\subsubsection{Riemannian Hessian in SPD Manifold}

Recall from Section \ref{section_tangent_normal_spd} that for every
$\b{X} \in \mathbb{S}_{++}^n$, the tangent space is:
\[
T_{\b{X}}\mathbb{S}_{++}^n = \mathbb{S}^n.
\]
Therefore, tangent vectors
$\b{\Delta}, \b{\Delta}_1, \b{\Delta}_2 \in T_{\b{X}}\mathbb{S}_{++}^n$
are symmetric matrices.

Also, recall from Section \ref{section_Riemannian_Hessian} that the Riemannian
Hessian is defined by the covariant derivative of the
Riemannian gradient:
\[
\operatorname{Hess} f(\b{X})[\b{\Delta}]
:=
\nabla_{\b{\Delta}} \operatorname{grad} f(\b{X}).
\]
Because the formula of the Riemannian gradient depends
on the metric, the Riemannian Hessian on the SPD
manifold also depends on the chosen metric.
In this subsection, we derive the Riemannian Hessian for
the Euclidean and affine-invariant metrics.

\begin{proposition}[Riemannian Hessian in SPD manifold under the Euclidean metric]
Let
$f : \mathbb{S}_{++}^n \to \mathbb{R}$
be a smooth function, and let
$\bar{f}$ be a smooth local extension of $f$
to an open neighborhood of $\b{X}$ in $\mathbb{S}^n$.
Under the Euclidean metric on $\mathbb{S}_{++}^n$,
the Riemannian Hessian is:
\begin{equation}\label{equation_hessian_SPD_manifold_Euclidean}
\boxed{
\operatorname{Hess}^{E} f(\b{X})[\b{\Delta}]
=
\operatorname{sym}\!\big(
\nabla^2 \bar{f}(\b{X})[\b{\Delta}]
\big),
}
\end{equation}
for $\b{\Delta} \in T_{\b{X}}\mathbb{S}_{++}^n$.
\end{proposition}

\begin{proof}
By Proposition \ref{proposition_Riemannian_gradient_SPD_Euclidean_metric}, the Riemannian gradient under the
Euclidean metric is:
\[
\operatorname{grad}^{E} f(\b{X})
=
\operatorname{sym}\!\big(\nabla \bar{f}(\b{X})\big).
\]
Since $\mathbb{S}_{++}^n$ is an open subset of the vector
space $\mathbb{S}^n$, the Euclidean Levi-Civita connection
is simply the ordinary directional derivative in the ambient
space (see Eq. (\ref{equation_Levi_Civita_connection_SPD_Euclidean})). Therefore:
\begin{align*}
\operatorname{Hess}^{E} f(\b{X})[\b{\Delta}]
&=
\nabla_{\b{\Delta}} \operatorname{grad}^{E} f(\b{X})
\\
&\overset{(\ref{equation_Levi_Civita_connection_SPD_Euclidean})}{=}
D\big(\operatorname{grad}^{E} f\big)(\b{X})[\b{\Delta}].
\end{align*}
Substituting the formula of the gradient gives:
\[
D\big(\operatorname{grad}^{E} f\big)(\b{X})[\b{\Delta}]
\overset{(\ref{equation_Riemannian_gradient_SPD_Euclidean_metric})}{=}
D\Big(
\operatorname{sym}\!\big(\nabla \bar{f}(\b{X})\big)
\Big)[\b{\Delta}].
\]
Because $\operatorname{sym}(\cdot)$ is linear, its derivative
is itself, so:
\[
D\Big(
\operatorname{sym}\!\big(\nabla \bar{f}(\b{X})\big)
\Big)[\b{\Delta}]
=
\operatorname{sym}\!\Big(
D(\nabla \bar{f})(\b{X})[\b{\Delta}]
\Big).
\]
By the definition of the Euclidean Hessian of the smooth
local extension, we have:
\[
D(\nabla \bar{f})(\b{X})[\b{\Delta}]
=
\nabla^2 \bar{f}(\b{X})[\b{\Delta}].
\]
Hence:
\[
\operatorname{Hess}^{E} f(\b{X})[\b{\Delta}]
=
\operatorname{sym}\!\big(
\nabla^2 \bar{f}(\b{X})[\b{\Delta}]
\big).
\]
This proves Eq.
\eqref{equation_hessian_SPD_manifold_Euclidean}.
\end{proof}

\begin{remark}[Simplification for symmetric local extensions]
If the smooth local extension $\bar{f}$ is chosen on the
ambient vector space $\mathbb{S}^n$, then its Euclidean
gradient and Euclidean Hessian are already symmetric.
In that case, Eq.
\eqref{equation_hessian_SPD_manifold_Euclidean}
simplifies to:
\[
\operatorname{Hess}^{E} f(\b{X})[\b{\Delta}]
=
\nabla^2 \bar{f}(\b{X})[\b{\Delta}].
\]
\end{remark}


\begin{proposition}[Riemannian Hessian in SPD manifold under the affine-invariant metric]
Let
$f : \mathbb{S}_{++}^n \to \mathbb{R}$
be a smooth function, and let
$\bar{f}$ be a smooth local extension of $f$
to an open neighborhood of $\b{X}$ in $\mathbb{S}^n$.
Under the affine-invariant metric on $\mathbb{S}_{++}^n$,
the Riemannian Hessian is:
\begin{equation}\label{equation_hessian_SPD_manifold_affine_invariant}
\boxed{
\begin{aligned}
&\operatorname{Hess}^{AI} f(\b{X})[\b{\Delta}]
\\
&\quad=
\b{X}\,
\operatorname{sym}\!\big(
\nabla^2 \bar{f}(\b{X})[\b{\Delta}]
\big)\,
\b{X}
\\
&\quad+
\frac{1}{2}
\Big(
\b{\Delta}\,\operatorname{sym}\!\big(\nabla \bar{f}(\b{X})\big)\,\b{X}
+
\b{X}\,\operatorname{sym}\!\big(\nabla \bar{f}(\b{X})\big)\,\b{\Delta}
\Big),
\end{aligned}
}
\end{equation}
for $\b{\Delta} \in T_{\b{X}}\mathbb{S}_{++}^n$.
\end{proposition}

\begin{proof}
By Proposition \ref{proposition_Riemannian_gradient_SPD_affine_invariant_metric}, the Riemannian gradient under the
affine-invariant metric is:
\[
\operatorname{grad}^{AI} f(\b{X})
=
\b{X}\,\operatorname{sym}\!\big(\nabla \bar{f}(\b{X})\big)\,\b{X}.
\]
Let:
\[
\b{S}(\b{X})
:=
\operatorname{sym}\!\big(\nabla \bar{f}(\b{X})\big).
\]
Then:
\[
\operatorname{grad}^{AI} f(\b{X})
=
\b{X}\b{S}(\b{X})\b{X}.
\]
By definition of the Riemannian Hessian, we have:
\[
\operatorname{Hess}^{AI} f(\b{X})[\b{\Delta}]
=
\nabla_{\b{\Delta}}^{AI}
\operatorname{grad}^{AI} f(\b{X}).
\]
Using the Levi-Civita connection in
Eq. \eqref{equation_Levi_Civita_connection_SPD_affine_invariant},
we get:
\begin{align}
\operatorname{Hess}^{AI} f(\b{X})[\b{\Delta}]
&=
D\big(\operatorname{grad}^{AI} f\big)(\b{X})[\b{\Delta}]
\nonumber\\
&\quad-
\frac{1}{2}
\Big(
\b{\Delta}\b{X}^{-1}\operatorname{grad}^{AI} f(\b{X})
\nonumber\\
&\quad\quad\,+
\operatorname{grad}^{AI} f(\b{X})\b{X}^{-1}\b{\Delta}
\Big).
\label{equation_hessian_SPD_AI_start}
\end{align}

We first calculate the ambient directional derivative of the
gradient. Since
$\operatorname{grad}^{AI} f(\b{X}) = \b{X}\b{S}(\b{X})\b{X}$,
the product rule gives:
\begin{align}
&D\big(\operatorname{grad}^{AI} f\big)(\b{X})[\b{\Delta}] \nonumber\\
&\quad=
\b{\Delta}\b{S}(\b{X})\b{X}
+
\b{X}\,D\b{S}(\b{X})[\b{\Delta}]\,\b{X}
+
\b{X}\b{S}(\b{X})\b{\Delta}.
\label{equation_directional_derivative_gradient_SPD_AI}
\end{align}

Now, because $\b{S}(\b{X}) = \operatorname{sym}\!\big(\nabla \bar{f}(\b{X})\big)$, we have:
\begin{align*}
D\b{S}(\b{X})[\b{\Delta}]
&=
\operatorname{sym}\!\Big(
D(\nabla \bar{f})(\b{X})[\b{\Delta}]
\Big) \\
&=
\operatorname{sym}\!\big(
\nabla^2 \bar{f}(\b{X})[\b{\Delta}]
\big).
\end{align*}
Therefore, Eq.
\eqref{equation_directional_derivative_gradient_SPD_AI}
becomes:
\begin{align}
&D\big(\operatorname{grad}^{AI} f\big)(\b{X})[\b{\Delta}] \nonumber\\
&\quad=
\b{\Delta}\b{S}(\b{X})\b{X}
+
\b{X}\,\operatorname{sym}\!\big(
\nabla^2 \bar{f}(\b{X})[\b{\Delta}]
\big)\,\b{X}
\nonumber\\
&\quad\quad+
\b{X}\b{S}(\b{X})\b{\Delta}.
\label{equation_directional_derivative_gradient_SPD_AI_2}
\end{align}

Next, we simplify the correction term in
Eq. \eqref{equation_hessian_SPD_AI_start}.
Because of $\operatorname{grad}^{AI} f(\b{X})
=
\b{X}\b{S}(\b{X})\b{X}$,
we have:
\begin{align*}
&\b{X}^{-1}\operatorname{grad}^{AI} f(\b{X})
=
\b{S}(\b{X})\b{X}, \\
&\operatorname{grad}^{AI} f(\b{X})\b{X}^{-1}
=
\b{X}\b{S}(\b{X}).
\end{align*}
Hence:
\begin{align}
&\frac{1}{2}
\Big(
\b{\Delta}\b{X}^{-1}\operatorname{grad}^{AI} f(\b{X})
+
\operatorname{grad}^{AI} f(\b{X})\b{X}^{-1}\b{\Delta}
\Big)
\nonumber\\
&\quad\quad=
\frac{1}{2}
\Big(
\b{\Delta}\b{S}(\b{X})\b{X}
+
\b{X}\b{S}(\b{X})\b{\Delta}
\Big).
\label{equation_connection_correction_term_SPD_AI}
\end{align}

Substituting
Eqs. \eqref{equation_directional_derivative_gradient_SPD_AI_2}
and
\eqref{equation_connection_correction_term_SPD_AI}
into
Eq. \eqref{equation_hessian_SPD_AI_start},
we obtain:
\begin{align*}
&\operatorname{Hess}^{AI} f(\b{X})[\b{\Delta}] \\
&=
\Big(
\b{\Delta}\b{S}(\b{X})\b{X}
+
\b{X}\,\operatorname{sym}\!\big(
\nabla^2 \bar{f}(\b{X})[\b{\Delta}]
\big)\,\b{X}
\\
&\quad+
\b{X}\b{S}(\b{X})\b{\Delta}
\Big) 
-
\frac{1}{2}
\Big(
\b{\Delta}\b{S}(\b{X})\b{X}
+
\b{X}\b{S}(\b{X})\b{\Delta}
\Big) \\
&=
\b{X}\,\operatorname{sym}\!\big(
\nabla^2 \bar{f}(\b{X})[\b{\Delta}]
\big)\,\b{X}
\\
&\quad+
\frac{1}{2}
\Big(
\b{\Delta}\b{S}(\b{X})\b{X}
+
\b{X}\b{S}(\b{X})\b{\Delta}
\Big).
\end{align*}
Recalling that
$\b{S}(\b{X}) = \operatorname{sym}(\nabla \bar{f}(\b{X}))$,
we get:
\begin{align*}
&\operatorname{Hess}^{AI} f(\b{X})[\b{\Delta}]
\\
&\quad=
\b{X}\,
\operatorname{sym}\!\big(
\nabla^2 \bar{f}(\b{X})[\b{\Delta}]
\big)\,
\b{X}
\\
&\quad\quad+
\frac{1}{2}
\Big(
\b{\Delta}\,\operatorname{sym}\!\big(\nabla \bar{f}(\b{X})\big)\,\b{X}
+
\b{X}\,\operatorname{sym}\!\big(\nabla \bar{f}(\b{X})\big)\,\b{\Delta}
\Big).
\end{align*}
This proves Eq.
\eqref{equation_hessian_SPD_manifold_affine_invariant}.
\end{proof}

\begin{remark}[Comparison of the Riemannian Hessians in SPD manifold]
Under the Euclidean metric, the SPD manifold is treated as
an open subset of the vector space $\mathbb{S}^n$, so the
Riemannian Hessian is simply the symmetric part of the
ambient Euclidean Hessian:
\[
\operatorname{Hess}^{E} f(\b{X})[\b{\Delta}]
=
\operatorname{sym}\!\big(
\nabla^2 \bar{f}(\b{X})[\b{\Delta}]
\big).
\]

Under the affine-invariant metric, the Hessian contains two
parts:
\begin{enumerate}
\item
the term:
\[
\b{X}\,
\operatorname{sym}\!\big(
\nabla^2 \bar{f}(\b{X})[\b{\Delta}]
\big)\,
\b{X},
\]
which is the metric-weighted second derivative, and

\item
the correction term:
\[
\frac{1}{2}
\Big(
\b{\Delta}\,\operatorname{sym}\!\big(\nabla \bar{f}(\b{X})\big)\,\b{X}
+
\b{X}\,\operatorname{sym}\!\big(\nabla \bar{f}(\b{X})\big)\,\b{\Delta}
\Big),
\]
which comes from the non-Euclidean Levi-Civita connection
of the affine-invariant geometry.
\end{enumerate}

Therefore, unlike the Euclidean case, the affine-invariant
Hessian is not only the second derivative of the local
extension. It also contains an additional geometric term
caused by the curvature of the metric.
\end{remark}

\subsubsection{Geodesic Equation on SPD Manifold}

Recall from Section \ref{section_geodesics} that a smooth curve on a
Riemannian manifold is a geodesic if and only if its
covariant acceleration vanishes. In the SPD manifold,
the tangent space at every point $\b{X} \in \mathbb{S}_{++}^n$
is the vector space $\mathbb{S}^n$. Therefore, a smooth
curve:
\[
\b{X}(t) \in \mathbb{S}_{++}^n,
\]
has velocity and acceleration matrices:
\[
\dot{\b{X}}(t) = \frac{d\b{X}(t)}{dt},
\qquad
\ddot{\b{X}}(t) = \frac{d^2\b{X}(t)}{dt^2},
\]
with $\dot{\b{X}}(t), \ddot{\b{X}}(t) \in \mathbb{S}^n$.

In this subsection, we derive the geodesic equation on
the SPD manifold for two important metrics, i.e., the
Euclidean metric and the affine-invariant metric.

\begin{lemma}[Directional derivative of the velocity field along a curve on the SPD manifold]
Let $\b{X}(t) \in \mathbb{S}_{++}^n$ be a smooth curve.
Define the tangent vector field along the curve by:
\[
\b{\Delta}(t) := \dot{\b{X}}(t) \in T_{\b{X}(t)}\mathbb{S}_{++}^n.
\]
Then the ordinary directional derivative of the velocity
field along the curve is:
\begin{align}\label{equation_DDeltaXXdot_Xddot_in_lemma}
D\b{\Delta}(\b{X}(t))[\dot{\b{X}}(t)] = \ddot{\b{X}}(t).
\end{align}
\end{lemma}

\begin{proof}
The vector field along the curve is defined by:
\[
\b{\Delta}(t) = \dot{\b{X}}(t).
\]
Differentiating both sides with respect to $t$ gives:
\[
\frac{d}{dt}\b{\Delta}(t) = \ddot{\b{X}}(t).
\]
By the chain rule, the derivative of $\b{\Delta}$ along the
curve $\b{X}(t)$ in the direction $\dot{\b{X}}(t)$ is:
\[
\frac{d}{dt}\b{\Delta}(t)
=
D\b{\Delta}(\b{X}(t))[\dot{\b{X}}(t)].
\]
Therefore:
\[
D\b{\Delta}(\b{X}(t))[\dot{\b{X}}(t)] = \ddot{\b{X}}(t).
\]
This proves the claim.
\end{proof}

\begin{proposition}[Geodesic equation on SPD manifold under the Euclidean metric]\label{proposition_geodesic_spd_euclidean_metric}
Let $\b{X}(t) \in \mathbb{S}_{++}^n$ be a smooth curve on
the SPD manifold endowed with the Euclidean metric in
Proposition \ref{proposition_euclidean_metric_spd}. Then, $\b{X}(t)$ is a geodesic if and only if:
\begin{equation}
\boxed{
\ddot{\b{X}}(t) = \b{0}.
}
\label{equation_geodesic_equation_SPD_manifold_Euclidean}
\end{equation}
\end{proposition}

\begin{proof}
According to Eq. (\ref{equation_geodesic_coordinate_free}), a smooth curve $\b{X}(t)$ is a geodesic if
and only if its covariant acceleration vanishes:
\[
\nabla^{E}_{\dot{\b{X}}(t)} \dot{\b{X}}(t) = \b{0}.
\]
According to Proposition \ref{proposition_Levi_Civita_connection_spd_euclidean_metric}, under the Euclidean metric,
the Levi-Civita connection on the SPD manifold is:
\[
\nabla^{E}_{\b{\Delta}_1}\b{\Delta}_2(\b{X})
=
D\b{\Delta}_2(\b{X})[\b{\Delta}_1].
\]
Now let:
\[
\b{\Delta}(t) := \dot{\b{X}}(t).
\]
Applying the connection formula with
$\b{\Delta}_1 = \b{\Delta}_2 = \dot{\b{X}}(t)$ gives:
\[
\nabla^{E}_{\dot{\b{X}}(t)} \dot{\b{X}}(t)
=
D\b{\Delta}(\b{X}(t))[\dot{\b{X}}(t)]
\overset{(\ref{equation_DDeltaXXdot_Xddot_in_lemma})}{=}
\ddot{\b{X}}(t).
\]
Therefore, the geodesic condition
$\nabla^{E}_{\dot{\b{X}}(t)} \dot{\b{X}}(t) = \b{0}$
is equivalent to:
\[
\ddot{\b{X}}(t) = \b{0}.
\]
This proves Eq.
\eqref{equation_geodesic_equation_SPD_manifold_Euclidean}.
\end{proof}

\begin{remark}[Interpretation of the Euclidean geodesic equation on SPD manifold]
Under the Euclidean metric, the SPD manifold is viewed
as an open subset of the vector space $\mathbb{S}^n$.
Therefore, the Levi-Civita connection is just the ordinary
directional derivative, and the geodesic equation becomes:
\[
\ddot{\b{X}}(t)=\b{0}.
\]
Hence, Euclidean geodesics are exactly straight lines in
the ambient vector space of symmetric matrices:
\[
\b{X}(t)=\b{X}(0)+t\,\dot{\b{X}}(0),
\]
as long as the curve remains inside $\mathbb{S}_{++}^n$ (this is because second-order derivative of $\b{X}(0)+t\,\dot{\b{X}}(0)$ becomes zero).
So, in the Euclidean geometry of SPD matrices, geodesics
do not bend intrinsically; they are ordinary linear motions
restricted to the open cone of positive definite matrices.
\end{remark}

\begin{proposition}[Geodesic equation on SPD manifold under the affine-invariant metric]\label{proposition_geodesic_spd_affine_invariant_metric}
Let $\b{X}(t) \in \mathbb{S}_{++}^n$ be a smooth curve on
the SPD manifold endowed with the affine-invariant metric
in Proposition \ref{proposition_affine_invariant_metric_spd}. Then, $\b{X}(t)$ is a geodesic if and
only if:
\begin{equation}
\boxed{
\ddot{\b{X}}(t) - \dot{\b{X}}(t)\b{X}(t)^{-1}\dot{\b{X}}(t)
= \b{0}.
}
\label{equation_geodesic_equation_SPD_manifold_affine_invariant}
\end{equation}
\end{proposition}

\begin{proof}
A smooth curve $\b{X}(t)$ is a geodesic if and only if:
\[
\nabla^{AI}_{\dot{\b{X}}(t)} \dot{\b{X}}(t) = \b{0}.
\]
According to Proposition \ref{proposition_Levi_Civita_connection_spd_affine_invariant_metric}, under the affine-invariant
metric, the Levi-Civita connection on the SPD manifold is:
\begin{align*}
\nabla^{AI}_{\b{\Delta}_1}\b{\Delta}_2&(\b{X})
=
D\b{\Delta}_2(\b{X})[\b{\Delta}_1]
\\
&-
\frac{1}{2}
\Big(
\b{\Delta}_1 \b{X}^{-1}\b{\Delta}_2(\b{X})
+
\b{\Delta}_2(\b{X}) \b{X}^{-1}\b{\Delta}_1
\Big).
\end{align*}
Now let:
\[
\b{\Delta}(t) := \dot{\b{X}}(t).
\]
Setting
\[
\b{\Delta}_1 = \b{\Delta}_2 = \dot{\b{X}}(t)
\]
in the connection formula gives:
\begin{align*}
&\nabla^{AI}_{\dot{\b{X}}(t)} \dot{\b{X}}(t)
=
D\b{\Delta}(\b{X}(t))[\dot{\b{X}}(t)]
\\
&\quad\quad-
\frac{1}{2}
\Big(
\dot{\b{X}}(t)\b{X}(t)^{-1}\dot{\b{X}}(t)
+
\dot{\b{X}}(t)\b{X}(t)^{-1}\dot{\b{X}}(t)
\Big) \\
&\quad\quad=
D\b{\Delta}(\b{X}(t))[\dot{\b{X}}(t)]
-
\dot{\b{X}}(t)\b{X}(t)^{-1}\dot{\b{X}}(t) \\
&\quad\quad\overset{(\ref{equation_DDeltaXXdot_Xddot_in_lemma})}{=}
\ddot{\b{X}}(t)
-
\dot{\b{X}}(t)\b{X}(t)^{-1}\dot{\b{X}}(t).
\end{align*}
Hence, the geodesic condition
$\nabla^{AI}_{\dot{\b{X}}(t)} \dot{\b{X}}(t) = \b{0}$
is equivalent to:
\[
\ddot{\b{X}}(t)
-
\dot{\b{X}}(t)\b{X}(t)^{-1}\dot{\b{X}}(t)
= \b{0}.
\]
This proves Eq.
\eqref{equation_geodesic_equation_SPD_manifold_affine_invariant}.
\end{proof}

\begin{remark}[Interpretation of the affine-invariant geodesic equation on SPD manifold]
Unlike the Euclidean metric, the affine-invariant metric
depends on the base point $\b{X}$. Therefore, its
Levi-Civita connection is not the ordinary directional
derivative. The correction term:
\[
\dot{\b{X}}(t)\b{X}(t)^{-1}\dot{\b{X}}(t)
\]
captures the intrinsic curvature of the SPD manifold under
the affine-invariant geometry. Hence, the geodesic equation
is no longer linear:
\[
\ddot{\b{X}}(t)
-
\dot{\b{X}}(t)\b{X}(t)^{-1}\dot{\b{X}}(t)
= \b{0}.
\]
This equation shows that affine-invariant geodesics bend
according to the inverse of the current point $\b{X}(t)$.
This is one reason why the affine-invariant geometry is
considered an intrinsic geometry of SPD matrices, in
contrast to the ambient Euclidean geometry.
\end{remark}

\subsubsection{Exponential Map in SPD Manifold}\label{section_exponential_map_spd}

In Section \ref{section_exponential_map_generalizing_addition}, we defined the exponential map on a
Riemannian manifold as the point reached at time $t=1$ by
the geodesic starting from a point with a prescribed initial
tangent vector. On the SPD manifold, the form of the
exponential map depends on the chosen metric. In this
subsection, we derive the exponential map on
$\mathbb{S}_{++}^n$ for the Euclidean metric and the
affine-invariant metric \cite{bhatia2009positive}.

We first prove a useful lemma for the matrix exponential
of a symmetric matrix.

\begin{lemma}[Derivative and positivity of the matrix exponential of a symmetric matrix]
\label{lemma_derivative_matrix_exponential_symmetric_SPD}
Let $\b{S}\in\mathbb{S}^n$ be a symmetric matrix, and
define:
\[
\b{E}(t):=\text{exp}(t\b{S}).
\]
Then, for every $t\in\mathbb{R}$, we have:
\begin{align}
&\b{E}(t)\in\mathbb{S}_{++}^n, \label{equation_exp_tS_is_SPD} \\
&\frac{d}{dt}\text{exp}(t\b{S})
=
\b{S}\text{exp}(t\b{S})
=
\text{exp}(t\b{S})\b{S}. \label{equation_derivative_exp_tS}
\end{align}
\end{lemma}

\begin{proof}
Because $\b{S}\in\mathbb{S}^n$, by eigendecomposition \cite{ghojogh2019eigenvalue}, there
exist an orthogonal matrix $\b{U}$ and a diagonal matrix
$\b{\Lambda}=\operatorname{diag}(\lambda_1,\dots,\lambda_n)$
such that:
\[
\b{S}=\b{U}\b{\Lambda}\b{U}^\top.
\]
Therefore:
\begin{align*}
\text{exp}(t\b{S})
&=
\text{exp}(t\b{U}\b{\Lambda}\b{U}^\top)
=
\b{U}\text{exp}(t\b{\Lambda})\b{U}^\top
\\
&=
\b{U}\operatorname{diag}(e^{t\lambda_1},\dots,e^{t\lambda_n})\b{U}^\top.
\end{align*}
Since $e^{t\lambda_i}>0$ for all $i$, the diagonal matrix
$\operatorname{diag}(e^{t\lambda_1},\dots,e^{t\lambda_n})$
is positive definite. Hence,
$\text{exp}(t\b{S})\in\mathbb{S}_{++}^n$, which proves
Eq. \eqref{equation_exp_tS_is_SPD}.

Now, we prove the derivative formula. Using the power
series definition of matrix exponential, we have:
\[
\text{exp}(t\b{S})
=
\sum_{k=0}^{\infty}\frac{(t\b{S})^k}{k!}
=
\sum_{k=0}^{\infty}\frac{t^k\b{S}^k}{k!}.
\]
Differentiating term by term with respect to $t$ gives:
\[
\frac{d}{dt}\text{exp}(t\b{S})
=
\sum_{k=1}^{\infty}\frac{k\,t^{k-1}\b{S}^k}{k!}
=
\sum_{k=1}^{\infty}\frac{t^{k-1}\b{S}^k}{(k-1)!}.
\]
Re-indexing the summation with $m:=k-1$ yields:
\begin{align*}
\frac{d}{dt}\text{exp}(t\b{S})
&=
\sum_{m=0}^{\infty}\frac{t^m\b{S}^{m+1}}{m!}
\\
&=
\b{S}\sum_{m=0}^{\infty}\frac{t^m\b{S}^m}{m!}
=
\b{S}\text{exp}(t\b{S}).
\end{align*}
Because $\text{exp}(t\b{S})$ is a power series in $\b{S}$, it
commutes with $\b{S}$. Therefore:
\[
\b{S}\text{exp}(t\b{S})=\text{exp}(t\b{S})\b{S}.
\]
This proves Eq. \eqref{equation_derivative_exp_tS}.
\end{proof}

\begin{proposition}[Exponential map on SPD manifold under the Euclidean metric]
\label{proposition_exponential_map_SPD_manifold_Euclidean_metric}
Let $\b{X}\in\mathbb{S}_{++}^n$ and
$\b{\Delta}\in T_{\b{X}}\mathbb{S}_{++}^n=\mathbb{S}^n$.
Under the Euclidean metric on $\mathbb{S}_{++}^n$, the
Riemannian exponential map is:
\begin{equation}
\boxed{
\operatorname{Exp}^{E}_{\b{X}}(\b{\Delta})
=
\b{X}+\b{\Delta},
}
\label{equation_exponential_map_SPD_manifold_Euclidean_metric}
\end{equation}
whenever $\b{X}+\b{\Delta}\in\mathbb{S}_{++}^n$.
\end{proposition}

\begin{proof}
By definition of exponential map, $\operatorname{Exp}^{E}_{\b{X}}(\b{\Delta})$
is the point reached at time $t=1$ by the geodesic
$\b{Y}(t)$ satisfying:
\[
\b{Y}(0)=\b{X},
\qquad
\dot{\b{Y}}(0)=\b{\Delta}.
\]
According to Proposition \ref{proposition_geodesic_spd_euclidean_metric}, under the Euclidean metric,
the geodesic equation on $\mathbb{S}_{++}^n$ is:
\[
\ddot{\b{Y}}(t)=0.
\]
We solve this matrix differential equation step by step.

\textbf{Step 1: Integrate the geodesic equation once.}

Since $\ddot{\b{Y}}(t)=0$, integrating with respect to $t$
gives:
\[
\dot{\b{Y}}(t)=\b{C},
\]
where $\b{C}\in\mathbb{S}^n$ is a constant matrix.

\textbf{Step 2: Use the initial velocity condition.}

At $t=0$, we have:
\[
\dot{\b{Y}}(0)=\b{\Delta}.
\]
Hence:
\[
\b{C}=\b{\Delta}.
\]
Therefore:
\[
\dot{\b{Y}}(t)=\b{\Delta}.
\]

\textbf{Step 3: Integrate once more.}

Integrating again with respect to $t$ gives:
\[
\b{Y}(t)=\b{X}_0+t\b{\Delta},
\]
where $\b{X}_0$ is a constant symmetric matrix.

\textbf{Step 4: Use the initial position condition.}

Because $\b{Y}(0)=\b{X}$, we obtain:
\[
\b{X}=\b{X}_0.
\]
Therefore, the geodesic is:
\begin{equation*}
\b{Y}(t)=\b{X}+t\b{\Delta}.
\end{equation*}

\textbf{Step 5: Evaluate at $t=1$.}

By definition of exponential map:
\[
\operatorname{Exp}^{E}_{\b{X}}(\b{\Delta})
=
\b{Y}(1)
=
\b{X}+\b{\Delta}.
\]
This proves Eq.
\eqref{equation_exponential_map_SPD_manifold_Euclidean_metric}.

Finally, because $\mathbb{S}_{++}^n$ is an open subset of
$\mathbb{S}^n$, the line $\b{X}+t\b{\Delta}$ remains in
$\mathbb{S}_{++}^n$ for sufficiently small $t$, but it need
not remain in $\mathbb{S}_{++}^n$ for all $t$. Therefore,
the Euclidean exponential map is given by Eq.
\eqref{equation_exponential_map_SPD_manifold_Euclidean_metric}
on the domain where $\b{X}+\b{\Delta}$ is still SPD.
\end{proof}

\begin{remark}[Interpretation of the Euclidean exponential map on SPD manifold]
Under the Euclidean metric, the SPD manifold is viewed as
an open subset of the vector space $\mathbb{S}^n$. Hence,
the geodesics are just straight lines:
\[
\b{Y}(t)=\b{X}+t\b{\Delta}.
\]
Therefore, the exponential map is simply ordinary matrix
addition, as long as the endpoint remains inside
$\mathbb{S}_{++}^n$.
\end{remark}

\begin{proposition}[Exponential map on SPD manifold under the affine-invariant metric]
\label{proposition_exponential_map_SPD_manifold_affine_invariant_metric}
Let $\b{X}\in\mathbb{S}_{++}^n$ and
$\b{\Delta}\in T_{\b{X}}\mathbb{S}_{++}^n=\mathbb{S}^n$.
Under the affine-invariant metric on $\mathbb{S}_{++}^n$,
the Riemannian exponential map is:
\begin{equation}\label{equation_exponential_map_SPD_manifold_affine_invariant_metric}
\boxed{
\operatorname{Exp}^{AI}_{\b{X}}(\b{\Delta})
=
\b{X}^{1/2}
\exp\!\Big(
\b{X}^{-1/2}\b{\Delta}\b{X}^{-1/2}
\Big)
\b{X}^{1/2},
}
\end{equation}
where $\exp(\cdot)$ denotes the ordinary matrix exponential.
\end{proposition}

\begin{proof}
By definition of exponential map,
$\operatorname{Exp}^{AI}_{\b{X}}(\b{\Delta})$ is the point
reached at time $t=1$ by the geodesic $\b{Y}(t)$ such that:
\[
\b{Y}(0)=\b{X},
\qquad
\dot{\b{Y}}(0)=\b{\Delta}.
\]
According to Proposition \ref{proposition_geodesic_spd_affine_invariant_metric}, under the affine-invariant
metric, the geodesic equation on $\mathbb{S}_{++}^n$ is:
\begin{equation}
\ddot{\b{Y}}(t)
=
\dot{\b{Y}}(t)\,\b{Y}(t)^{-1}\,\dot{\b{Y}}(t).
\label{equation_geodesic_equation_SPD_affine_invariant_for_exponential}
\end{equation}

We claim that the solution of this initial value problem is:
\begin{equation}
\b{Y}(t)
=
\b{X}^{1/2}
\exp\Big(
t\,\b{X}^{-1/2}\b{\Delta}\b{X}^{-1/2}
\Big)
\b{X}^{1/2}.
\label{equation_candidate_geodesic_SPD_affine_invariant}
\end{equation}

We define:
\begin{equation}
\b{S}
:=
\b{X}^{-1/2}\b{\Delta}\b{X}^{-1/2}.
\label{equation_definition_S_for_affine_invariant_exponential_SPD}
\end{equation}
Because $\b{X}\in\mathbb{S}_{++}^n$, its square root
$\b{X}^{1/2}$ is symmetric positive definite, and since
$\b{\Delta}\in\mathbb{S}^n$, we have:
\[
\b{S}^\top
=
(\b{X}^{-1/2}\b{\Delta}\b{X}^{-1/2})^\top
=
\b{X}^{-1/2}\b{\Delta}\b{X}^{-1/2}
=
\b{S}.
\]
Hence, $\b{S}\in\mathbb{S}^n$.

Using Eq.
\eqref{equation_candidate_geodesic_SPD_affine_invariant},
we now verify the initial conditions and the geodesic
equation.

\textbf{Step 1: Verify that $\b{Y}(t)\in\mathbb{S}_{++}^n$.}

By Lemma
\ref{lemma_derivative_matrix_exponential_symmetric_SPD},
because $\b{S}$ is symmetric, we have
$\exp(t\b{S})\in\mathbb{S}_{++}^n$ for all $t$. Therefore:
\[
\b{Y}(t)
=
\b{X}^{1/2}\exp(t\b{S})\b{X}^{1/2}
\in
\mathbb{S}_{++}^n,
\]
because congruence by the invertible matrix $\b{X}^{1/2}$
preserves positive definiteness.

\textbf{Step 2: Verify the initial position.}

At $t=0$, we have:
\[
\b{Y}(0)
=
\b{X}^{1/2}\text{exp}(\b{0})\b{X}^{1/2}
=
\b{X}^{1/2}\b{I}\b{X}^{1/2}
=
\b{X}.
\]

\textbf{Step 3: Compute the first derivative.}

By Lemma
\ref{lemma_derivative_matrix_exponential_symmetric_SPD}, we have:
\[
\frac{d}{dt}\text{exp}(t\b{S})=\b{S}\text{exp}(t\b{S}).
\]
Hence:
\begin{align}
\dot{\b{Y}}(t)
&=
\b{X}^{1/2}
\frac{d}{dt}\text{exp}(t\b{S})
\b{X}^{1/2}
\nonumber\\
&=
\b{X}^{1/2}\b{S}\text{exp}(t\b{S})\b{X}^{1/2}.
\label{equation_Ydot_affine_invariant_exponential_SPD}
\end{align}
At $t=0$:
\begin{align*}
\dot{\b{Y}}(0)
&=
\b{X}^{1/2}\b{S}\text{exp}(\b{0})\b{X}^{1/2} \\
&=
\b{X}^{1/2}\b{S}\b{X}^{1/2} \\
&\overset{(\ref{equation_definition_S_for_affine_invariant_exponential_SPD})}{=}
\b{X}^{1/2}
(\b{X}^{-1/2}\b{\Delta}\b{X}^{-1/2})
\b{X}^{1/2}
=
\b{\Delta}.
\end{align*}
Therefore, the initial velocity condition is satisfied.

\textbf{Step 4: Compute the second derivative.}

Differentiating Eq.
\eqref{equation_Ydot_affine_invariant_exponential_SPD}
again gives:
\begin{align}
\ddot{\b{Y}}(t)
&=
\b{X}^{1/2}\b{S}^2\text{exp}(t\b{S})\b{X}^{1/2}.
\label{equation_Yddot_affine_invariant_exponential_SPD}
\end{align}

\textbf{Step 5: Compute the inverse of $\b{Y}(t)$.}

Since $\text{exp}(t\b{S})$ is invertible and
$\text{exp}(t\b{S})^{-1}=\text{exp}(-t\b{S})$, we have:
\begin{equation}
\b{Y}(t)^{-1}
=
\b{X}^{-1/2}\text{exp}(-t\b{S})\b{X}^{-1/2}.
\label{equation_inverse_Y_affine_invariant_exponential_SPD}
\end{equation}

\textbf{Step 6: Verify the geodesic equation.}

Using Eqs.
\eqref{equation_Ydot_affine_invariant_exponential_SPD}
and
\eqref{equation_inverse_Y_affine_invariant_exponential_SPD},
we get:
\begin{align}
&\dot{\b{Y}}(t)\b{Y}(t)^{-1}\dot{\b{Y}}(t) \nonumber\\
&=
\Big(
\b{X}^{1/2}\b{S}\exp(t\b{S})\b{X}^{1/2}
\Big)
\Big(
\b{X}^{-1/2}\exp(-t\b{S})\b{X}^{-1/2}
\Big)
\nonumber\\
&\quad\quad\quad\quad\quad\quad\quad\quad\quad\quad\quad\quad\Big(
\b{X}^{1/2}\b{S}\exp(t\b{S})\b{X}^{1/2}
\Big) \nonumber\\
&\overset{(a)}{=}
\b{X}^{1/2}
\b{S}\exp(t\b{S})\exp(-t\b{S})\b{S}\exp(t\b{S})
\b{X}^{1/2} \nonumber\\
&\overset{(b)}{=} \b{X}^{1/2}\b{S}^2\exp(t\b{S})\b{X}^{1/2}, \label{equation_YdYinvYd_XhalfS2exptSXhalf_in_proof}
\end{align}
where $(a)$ is because we have $\b{X}^{1/2} \b{X}^{-1/2} = \b{I}$ and $\b{X}^{-1/2} \b{X}^{1/2} = \b{I}$, and $(b)$ is because we have $\exp(t\b{S})\exp(-t\b{S})=\b{I}$. 

Comparing Eqs.
\eqref{equation_Yddot_affine_invariant_exponential_SPD} and (\ref{equation_YdYinvYd_XhalfS2exptSXhalf_in_proof}) gives:
\[
\dot{\b{Y}}(t)\b{Y}(t)^{-1}\dot{\b{Y}}(t)
=
\ddot{\b{Y}}(t).
\]
Thus, $\b{Y}(t)$ satisfies the affine-invariant geodesic
equation
\eqref{equation_geodesic_equation_SPD_affine_invariant_for_exponential}.

Therefore, $\b{Y}(t)$ is the geodesic starting from
$\b{X}$ with initial velocity $\b{\Delta}$. Evaluating at
$t=1$ gives:
\begin{align*}
\operatorname{Exp}^{AI}_{\b{X}}(\b{\Delta})
&=
\b{Y}(1)
\\
&=
\b{X}^{1/2}
\exp\Big(
\b{X}^{-1/2}\b{\Delta}\b{X}^{-1/2}
\Big)
\b{X}^{1/2}.
\end{align*}
This proves Eq.
\eqref{equation_exponential_map_SPD_manifold_affine_invariant_metric}.
\end{proof}

\begin{remark}[Interpretation of the affine-invariant exponential map on SPD manifold]
Let:
\[
\b{S}:=\b{X}^{-1/2}\b{\Delta}\b{X}^{-1/2}.
\]
The affine-invariant exponential map can be understood in
three steps.

First, the tangent vector $\b{\Delta}$ at $\b{X}$ is
normalized by:
\[
\b{\Delta}\mapsto \b{S}
=
\b{X}^{-1/2}\b{\Delta}\b{X}^{-1/2}.
\]

Second, in these normalized coordinates, one considers the
curve:
\[
\widetilde{\b{Y}}(t)=\exp(t\b{S}),
\]
which starts at the identity matrix:
\[
\widetilde{\b{Y}}(0)=\exp(\b{0})=\b{I}.
\]

Third, this normalized curve is mapped back to the original
base point $\b{X}$ by congruence with $\b{X}^{1/2}$:
\[
\b{Y}(t)
=
\b{X}^{1/2}\widetilde{\b{Y}}(t)\b{X}^{1/2}
=
\b{X}^{1/2}\exp(t\b{S})\b{X}^{1/2}.
\]
This final step transfers the geodesic
from the identity back to the manifold point $\b{X}$.
\end{remark}

\subsubsection{Retraction Map in SPD Manifold}\label{section_retraction_spd}

As discussed in Section \ref{section_retraction}, a retraction map is a smooth
mapping from the tangent space back to the manifold that
approximates the exponential map at least to first order.
On the SPD manifold \(\mathbb{S}_{++}^n\), retractions depend
on the chosen metric. In this subsection, we derive
retraction maps for the Euclidean metric and the
affine-invariant metric.

Recall from Section \ref{section_tangent_normal_spd} that for every
\(\b{X} \in \mathbb{S}_{++}^n\), we have:
\[
T_{\b{X}}\mathbb{S}_{++}^n = \mathbb{S}^n.
\]
Hence, every tangent vector
\(\b{\Delta} \in T_{\b{X}}\mathbb{S}_{++}^n\)
is a symmetric matrix.

We first derive a simple retraction for the Euclidean metric.

\begin{proposition}[Retraction map on SPD manifold under the Euclidean metric]
Let \(\b{X} \in \mathbb{S}_{++}^n\) and
\(\b{\Delta} \in T_{\b{X}}\mathbb{S}_{++}^n = \mathbb{S}^n\).
We define:
\begin{equation}
\boxed{
\operatorname{Ret}^{\mathrm{Euc}}_{\b{X}}(\b{\Delta})
:=
\b{X} + \b{\Delta}.
}
\label{equation_retraction_SPD_manifold_Euclidean}
\end{equation}
Then, for sufficiently small \(\b{\Delta}\), this map is a valid
retraction on \(\mathbb{S}_{++}^n\) under the Euclidean metric.
\end{proposition}

\begin{proof}
To prove that
\(\operatorname{Ret}^{\mathrm{Euc}}_{\b{X}}\)
is a retraction, we verify the two properties in
Definition \ref{definition_retraction}.

\textbf{Step 1: Show that the image remains in
\(\mathbb{S}_{++}^n\) for sufficiently small \(\b{\Delta}\).}

Because \(\b{X} \in \mathbb{S}_{++}^n\), the matrix \(\b{X}\) is
symmetric positive definite. Also,
\(\b{\Delta} \in \mathbb{S}^n\), so
\(\b{X}+\b{\Delta}\) is symmetric.

It remains to show positive definiteness for sufficiently
small \(\b{\Delta}\). Since \(\b{X}\) is positive definite, its
smallest eigenvalue satisfies:
\[
\lambda_{\min}(\b{X}) > 0.
\]
If \(\Vert\b{\Delta}\Vert_2 < \lambda_{\min}(\b{X})\), then for every
nonzero \(\b{v} \in \mathbb{R}^n\), we have:
\begin{align*}
\b{v}^\top (\b{X}+\b{\Delta}) \b{v}
&=
\b{v}^\top \b{X}\b{v}
+
\b{v}^\top \b{\Delta}\b{v} \\
&\geq
\lambda_{\min}(\b{X}) \Vert\b{v}\Vert_2^2
-
|\b{v}^\top \b{\Delta}\b{v}| \\
&\geq
\lambda_{\min}(\b{X}) \Vert\b{v}\Vert_2^2
-
\Vert\b{\Delta}\Vert_2 \Vert\b{v}\Vert_2^2 \\
&=
\left(
\lambda_{\min}(\b{X})-\Vert\b{\Delta}\Vert_2
\right)\Vert\b{v}\Vert_2^2
> 0.
\end{align*}
Therefore, \(\b{X}+\b{\Delta} \in \mathbb{S}_{++}^n\) for all
sufficiently small \(\b{\Delta}\).

\textbf{Step 2: Verify the identity property.}

By Eq. \eqref{equation_retraction_SPD_manifold_Euclidean}, we have:
\[
\operatorname{Ret}^{\mathrm{Euc}}_{\b{X}}(\b{0})
=
\b{X}+\b{0}
=
\b{X}.
\]

\textbf{Step 3: Verify the local rigidity property.}

Consider the curve:
\[
\b{Y}(t)
:=
\operatorname{Ret}^{\mathrm{Euc}}_{\b{X}}(t\b{\Delta})
=
\b{X}+t\b{\Delta}.
\]
Differentiating with respect to \(t\) gives:
\[
\frac{d}{dt}\b{Y}(t) = \b{\Delta}.
\]
Hence:
\[
\left.
\frac{d}{dt}
\operatorname{Ret}^{\mathrm{Euc}}_{\b{X}}(t\b{\Delta})
\right|_{t=0}
=
\b{\Delta}.
\]
Therefore, the differential of
\(\operatorname{Ret}^{\mathrm{Euc}}_{\b{X}}\) at \(\b{0}\)
is the identity map on
\(T_{\b{X}}\mathbb{S}_{++}^n\).

Thus,
\(\operatorname{Ret}^{\mathrm{Euc}}_{\b{X}}\)
is a valid retraction under the Euclidean metric.
\end{proof}

\begin{remark}[Interpretation of the Euclidean retraction on SPD manifold]
Under the Euclidean metric, the SPD manifold is an open
subset of the vector space \(\mathbb{S}^n\). Therefore, the
simplest way to move from \(\b{X}\) along a tangent vector
\(\b{\Delta}\) is just ordinary matrix addition:
\[
\b{X} \mapsto \b{X}+\b{\Delta}.
\]
This is exactly the first-order approximation of the
Euclidean exponential map:
\[
\operatorname{Exp}^{\mathrm{Euc}}_{\b{X}}(\b{\Delta})
=
\b{X}+\b{\Delta}.
\]
Hence, for sufficiently small steps, the Euclidean retraction
coincides with the Euclidean exponential map.
\end{remark}

For the affine-invariant metric, a more intrinsic retraction
is obtained by truncating the affine-invariant exponential
map after the quadratic term.

\begin{lemma}[Positivity of \(\b{I}_n+\b{S}+\frac{1}{2}\b{S}^2\)]\label{lemma_I_S_halfS2_psd}
Let \(\b{S} \in \mathbb{S}^n\). Then:
\begin{equation}
\b{I}_n+\b{S}+\frac{1}{2}\b{S}^2 \in \mathbb{S}_{++}^n.
\label{equation_positivity_I_plus_S_plus_half_S_square}
\end{equation}
\end{lemma}

\begin{proof}
Because \(\b{S}\) is symmetric, by eigendecomposition \cite{ghojogh2019eigenvalue} there
exist an orthogonal matrix \(\b{U}\) and a diagonal matrix
\(\b{\Lambda}=\operatorname{diag}(\lambda_1,\dots,\lambda_n)\)
such that:
\[
\b{S} = \b{U}\b{\Lambda}\b{U}^\top.
\]
Therefore:
\begin{align*}
&\b{I}_n+\b{S}+\frac{1}{2}\b{S}^2
=
\b{U}
\left(
\b{I}_n+\b{\Lambda}+\frac{1}{2}\b{\Lambda}^2
\right)
\b{U}^\top \\
&=
\b{U}\,
\operatorname{diag}\!\left(
1+\lambda_1+\frac{1}{2}\lambda_1^2,\,
\dots,\,
1+\lambda_n+\frac{1}{2}\lambda_n^2
\right)
\b{U}^\top.
\end{align*}
For every \(i\), we have:
\begin{align*}
1+\lambda_i+\frac{1}{2}\lambda_i^2
&=
\frac{1}{2}\left(
\lambda_i^2+2\lambda_i+2
\right)
\\
&=
\frac{1}{2}\left(
(\lambda_i+1)^2+1
\right)
> 0.
\end{align*}
Hence all eigenvalues of
\(\b{I}_n+\b{S}+\frac{1}{2}\b{S}^2\)
are strictly positive. Therefore, this matrix is symmetric
positive definite, which proves
Eq. \eqref{equation_positivity_I_plus_S_plus_half_S_square}.
\end{proof}

\begin{proposition}[Retraction map on SPD manifold under the affine-invariant metric]
Let \(\b{X} \in \mathbb{S}_{++}^n\) and
\(\b{\Delta} \in T_{\b{X}}\mathbb{S}_{++}^n = \mathbb{S}^n\).
We define:
\begin{equation}
\boxed{
\operatorname{Ret}^{AI}_{\b{X}}(\b{\Delta})
:=
\b{X}
+
\b{\Delta}
+
\frac{1}{2}\b{\Delta}\b{X}^{-1}\b{\Delta}.
}
\label{equation_retraction_SPD_manifold_affine_invariant}
\end{equation}
Then,
\(\operatorname{Ret}^{AI}_{\b{X}}\)
is a valid retraction on \(\mathbb{S}_{++}^n\) under the
affine-invariant metric.
\end{proposition}

\begin{proof}
We again verify the two properties in Definition \ref{definition_retraction}.

\textbf{Step 1: Show that the image belongs to
\(\mathbb{S}_{++}^n\).}

First, because \(\b{X}\) and \(\b{\Delta}\) are symmetric and
\(\b{X}^{-1}\) is also symmetric, we have:
\[
(\b{\Delta}\b{X}^{-1}\b{\Delta})^\top
=
\b{\Delta}^\top(\b{X}^{-1})^\top\b{\Delta}^\top
=
\b{\Delta}\b{X}^{-1}\b{\Delta}.
\]
Hence,
\(\operatorname{Ret}^{AI}_{\b{X}}(\b{\Delta})\) in Eq. (\ref{equation_retraction_SPD_manifold_affine_invariant})
is symmetric:
\begin{align*}
\operatorname{Ret}^{AI}_{\b{X}}(\b{\Delta}) = \big( \operatorname{Ret}^{AI}_{\b{X}}(\b{\Delta}) \big)^\top.
\end{align*}

Now, we rewrite every term in Eq. (\ref{equation_retraction_SPD_manifold_affine_invariant}) so that \(\b{X}^{1/2}\) appears on
the left and on the right.

First, for the term \(\b{X}\) in Eq. (\ref{equation_retraction_SPD_manifold_affine_invariant}), we have:
\begin{equation*}
\b{X}
=
\b{X}^{1/2}\b{X}^{1/2}
=
\b{X}^{1/2}\b{I}_n\b{X}^{1/2}.
\end{equation*}

Second, for the term \(\b{\Delta}\) in Eq. (\ref{equation_retraction_SPD_manifold_affine_invariant}), we insert the identity
\(\b{I}_n=\b{X}^{1/2}\b{X}^{-1/2}
=\b{X}^{-1/2}\b{X}^{1/2}\) on both sides:
\begin{align*}
\b{\Delta}
&=
\b{X}^{1/2}\b{X}^{-1/2}\b{\Delta}\b{X}^{-1/2}\b{X}^{1/2} \\
&=
\b{X}^{1/2}
\big(
\b{X}^{-1/2}\b{\Delta}\b{X}^{-1/2}
\big)
\b{X}^{1/2}.
\end{align*}

Third, for the quadratic term in Eq. (\ref{equation_retraction_SPD_manifold_affine_invariant}), since
\(\b{X}^{-1}=\b{X}^{-1/2}\b{X}^{-1/2}\), we have:
\begin{align*}
\b{\Delta}\b{X}^{-1}\b{\Delta}
&=
\b{\Delta}\b{X}^{-1/2}\b{X}^{-1/2}\b{\Delta}.
\end{align*}
Now, we insert \(\b{I}_n=\b{X}^{1/2}\b{X}^{-1/2}\) on the left and
\(\b{I}_n=\b{X}^{-1/2}\b{X}^{1/2}\) on the right:
\begin{align*}
&\b{\Delta}\b{X}^{-1}\b{\Delta} \\
&=
\b{X}^{1/2}
\big(
\b{X}^{-1/2}\b{\Delta}\b{X}^{-1/2}
\b{X}^{-1/2}\b{\Delta}\b{X}^{-1/2}
\big)
\b{X}^{1/2} \\
&=
\b{X}^{1/2}
\big(
\b{X}^{-1/2}\b{\Delta}\b{X}^{-1/2}
\big)^2
\b{X}^{1/2}.
\end{align*}
Therefore:
\begin{equation*}
\frac{1}{2}\b{\Delta}\b{X}^{-1}\b{\Delta}
=
\b{X}^{1/2}
\left(
\frac{1}{2}
\big(
\b{X}^{-1/2}\b{\Delta}\b{X}^{-1/2}
\big)^2
\right)
\b{X}^{1/2}.
\end{equation*}

Substituting these three expressions into
\(
\b{X}
+
\b{\Delta}
+
\frac{1}{2}\b{\Delta}\b{X}^{-1}\b{\Delta}
\),
we obtain:
\begin{align*}
&\operatorname{Ret}^{AI}_{\b{X}}(\b{\Delta})
=
\b{X}
+
\b{\Delta}
+
\frac{1}{2}\b{\Delta}\b{X}^{-1}\b{\Delta} \\
&=
\b{X}^{1/2}\b{I}_n\b{X}^{1/2}
+
\b{X}^{1/2}
\big(
\b{X}^{-1/2}\b{\Delta}\b{X}^{-1/2}
\big)
\b{X}^{1/2} \\
&\quad+
\b{X}^{1/2}
\left(
\frac{1}{2}
\big(
\b{X}^{-1/2}\b{\Delta}\b{X}^{-1/2}
\big)^2
\right)
\b{X}^{1/2}.
\end{align*}
Factoring out \(\b{X}^{1/2}\) from the left and right gives:
\begin{align*}
&\operatorname{Ret}^{AI}_{\b{X}}(\b{\Delta}) =
\b{X}^{1/2}
\Big(
\b{I}_n
+
\b{X}^{-1/2}\b{\Delta}\b{X}^{-1/2}
\\
&\quad\quad\quad\quad\quad\quad\quad+
\frac{1}{2}
\big(
\b{X}^{-1/2}\b{\Delta}\b{X}^{-1/2}
\big)^2
\Big)
\b{X}^{1/2}.
\end{align*}

We define:
\[
\b{S}
:=
\b{X}^{-1/2}\b{\Delta}\b{X}^{-1/2}.
\]
Because \(\b{X}^{-1/2}\) and \(\b{\Delta}\) are symmetric,
\(\b{S}\in\mathbb{S}^n\). Therefore:
\[
\operatorname{Ret}^{AI}_{\b{X}}(\b{\Delta})
=
\b{X}^{1/2}
\left(
\b{I}_n+\b{S}+\frac{1}{2}\b{S}^2
\right)
\b{X}^{1/2}.
\]
By Lemma \ref{lemma_I_S_halfS2_psd},
\(\b{I}_n+\b{S}+\frac{1}{2}\b{S}^2 \in \mathbb{S}_{++}^n\).
Since congruence by \(\b{X}^{1/2}\) preserves positive
definiteness, we conclude:
\[
\operatorname{Ret}^{AI}_{\b{X}}(\b{\Delta})
\in
\mathbb{S}_{++}^n.
\]

\textbf{Step 2: Verify the identity property.}

Substituting \(\b{\Delta}=\b{0}\) into
Eq. \eqref{equation_retraction_SPD_manifold_affine_invariant}
gives:
\[
\operatorname{Ret}^{AI}_{\b{X}}(\b{0})
=
\b{X}
+
\b{0}
+
\frac{1}{2}\b{0}\b{X}^{-1}\b{0}
=
\b{X}.
\]

\textbf{Step 3: Verify the local rigidity property.}

Consider:
\[
\b{Y}(t)
:=
\operatorname{Ret}^{AI}_{\b{X}}(t\b{\Delta})
=
\b{X}
+
t\b{\Delta}
+
\frac{t^2}{2}\b{\Delta}\b{X}^{-1}\b{\Delta}.
\]
Differentiating with respect to \(t\) yields:
\[
\frac{d}{dt}\b{Y}(t)
=
\b{\Delta}
+
t\b{\Delta}\b{X}^{-1}\b{\Delta}.
\]
Evaluating at \(t=0\), we obtain:
\[
\left.
\frac{d}{dt}
\operatorname{Ret}^{AI}_{\b{X}}(t\b{\Delta})
\right|_{t=0}
=
\b{\Delta}.
\]
Hence, the differential of
\(\operatorname{Ret}^{AI}_{\b{X}}\) at \(\b{0}\)
is the identity on \(T_{\b{X}}\mathbb{S}_{++}^n\).

Therefore,
\(\operatorname{Ret}^{AI}_{\b{X}}\)
is a valid retraction under the affine-invariant metric.
\end{proof}

\begin{remark}[Relation of the affine-invariant retraction to the affine-invariant exponential map]
From Section \ref{section_exponential_map_spd}, the affine-invariant exponential map
on the SPD manifold is:
\[
\operatorname{Exp}^{AI}_{\b{X}}(\b{\Delta})
=
\b{X}^{1/2}
\exp\!\left(
\b{X}^{-1/2}\b{\Delta}\b{X}^{-1/2}
\right)
\b{X}^{1/2}.
\]
Using the power series (Taylor series) of the matrix exponential:
\[
\exp(\b{S})
=
\b{I}_n+\b{S}+\frac{1}{2}\b{S}^2+\cdots,
\]
with
\(
\b{S}=\b{X}^{-1/2}\b{\Delta}\b{X}^{-1/2},
\)
we get:
\begin{align*}
\operatorname{Exp}^{AI}_{\b{X}}(\b{\Delta})
&=
\b{X}^{1/2}
\left(
\b{I}_n+\b{S}+\frac{1}{2}\b{S}^2+\cdots
\right)
\b{X}^{1/2} \\
&=
\b{X}
+
\b{\Delta}
+
\frac{1}{2}\b{\Delta}\b{X}^{-1}\b{\Delta}
+\cdots.
\end{align*}
Therefore, the retraction in
Eq. \eqref{equation_retraction_SPD_manifold_affine_invariant}
is the second-order truncation of the affine-invariant
exponential map. It is computationally simpler than the
exact exponential map while preserving the correct
first-order behavior required of a retraction.
\end{remark}

\begin{remark}[Why the Euclidean and affine-invariant retractions are different]
The Euclidean metric treats \(\mathbb{S}_{++}^n\) as an open
subset of the vector space \(\mathbb{S}^n\), so the simplest
retraction is just linear addition:
\[
\b{X}+\b{\Delta}.
\]
In contrast, the affine-invariant metric depends on the base
point \(\b{X}\) through \(\b{X}^{-1}\), so a natural retraction
should also depend on \(\b{X}\). This is why the
affine-invariant retraction contains the quadratic correction
term:
\[
\frac{1}{2}\b{\Delta}\b{X}^{-1}\b{\Delta}.
\]
That correction reflects the intrinsic geometry induced by
the affine-invariant metric.
\end{remark}

\subsubsection{Vector Transport in SPD Manifold}

Recall from Section \ref{section_definition_vector_transport} that a vector transport is a
smooth mapping that moves a tangent vector from one
tangent space to another tangent space on the manifold.
Also, recall from Section \ref{section_tangent_normal_spd} that for every
$\b{X} \in \mathbb{S}_{++}^n$, we have:
\[
T_{\b{X}}\mathbb{S}_{++}^n = \mathbb{S}^n.
\]
Therefore, every tangent vector on the SPD manifold is a
symmetric matrix.

In this subsection, we derive natural vector transports on
the SPD manifold for the Euclidean and affine-invariant
metrics. Under the Euclidean metric, because the SPD
manifold is an open subset of the vector space
$\mathbb{S}^n$, the natural vector transport is simply the
identity map. Under the affine-invariant metric, a natural
choice is the parallel transport along the affine-invariant
geodesic, which gives a canonical vector transport.

\begin{proposition}[Vector transport in SPD manifold under the Euclidean metric]
Let $\b{X} \in \mathbb{S}_{++}^n$ and let
$\b{\Delta}, \b{\Xi} \in T_{\b{X}}\mathbb{S}_{++}^n
= \mathbb{S}^n$.
Using the Euclidean retraction in Eq. (\ref{equation_retraction_SPD_manifold_Euclidean}),
we define:
\begin{equation}
\boxed{
\mathcal{T}^{E}_{\b{\Delta}}(\b{\Xi}) := \b{\Xi} }
\in T_{\operatorname{Ret}^{E}_{\b{X}}(\b{\Delta})}
\mathbb{S}_{++}^n.
\label{equation_vector_transport_SPD_manifold_Euclidean}
\end{equation}
Then, $\mathcal{T}^{E}$ is a valid vector transport on
$\mathbb{S}_{++}^n$ under the Euclidean metric.
\end{proposition}

\begin{proof}
According to Definition \ref{definition_vector_transport}, we need to verify the three
axioms of vector transport.

\textbf{Step 1: Well-defined mapping.}

By Section \ref{section_tangent_normal_spd}, for every point
$\b{Y} \in \mathbb{S}_{++}^n$, the tangent space is:
\[
T_{\b{Y}}\mathbb{S}_{++}^n = \mathbb{S}^n.
\]
Let:
\[
\b{Y} := \operatorname{Ret}^{E}_{\b{X}}(\b{\Delta}).
\]
Because $\operatorname{Ret}^{E}_{\b{X}}(\b{\Delta}) \in
\mathbb{S}_{++}^n$, we have:
\[
T_{\b{Y}}\mathbb{S}_{++}^n = \mathbb{S}^n.
\]
Since $\b{\Xi} \in \mathbb{S}^n$, Eq.
(\ref{equation_vector_transport_SPD_manifold_Euclidean})
gives:
\[
\mathcal{T}^{E}_{\b{\Delta}}(\b{\Xi}) = \b{\Xi}
\in \mathbb{S}^n
=
T_{\operatorname{Ret}^{E}_{\b{X}}(\b{\Delta})}
\mathbb{S}_{++}^n.
\]
Hence, the mapping is well-defined.

\textbf{Step 2: Linearity in $\b{\Xi}$.}

Let $a,b \in \mathbb{R}$ and
$\b{\Xi}_1,\b{\Xi}_2 \in T_{\b{X}}\mathbb{S}_{++}^n$.
Then:
\begin{align*}
\mathcal{T}^{E}_{\b{\Delta}}
(a\b{\Xi}_1 + b\b{\Xi}_2)
&=
a\b{\Xi}_1 + b\b{\Xi}_2 \\
&=
a\,\mathcal{T}^{E}_{\b{\Delta}}(\b{\Xi}_1)
+
b\,\mathcal{T}^{E}_{\b{\Delta}}(\b{\Xi}_2).
\end{align*}
Therefore, the map is linear in the transported vector.

\textbf{Step 3: Consistency at zero.}

For $\b{\Delta}=\b{0}$, we have:
\[
\mathcal{T}^{E}_{\b{0}}(\b{\Xi}) = \b{\Xi}.
\]
Therefore, the transport reduces to the identity map at
zero displacement.

Thus, $\mathcal{T}^{E}$ satisfies all axioms of vector
transport. This proves that Eq.
(\ref{equation_vector_transport_SPD_manifold_Euclidean})
defines a valid vector transport on the SPD manifold under
the Euclidean metric.
\end{proof}

\begin{remark}[Interpretation of the Euclidean vector transport on SPD manifold]
Under the Euclidean metric, the SPD manifold is an open
subset of the vector space $\mathbb{S}^n$. Although the
base point changes from $\b{X}$ to
$\operatorname{Ret}^{E}_{\b{X}}(\b{\Delta})$, all tangent
spaces are identified with the same vector space
$\mathbb{S}^n$. Therefore, no correction or projection is
needed, and the vector transport is simply the same
symmetric matrix viewed in the new tangent space.
\end{remark}

\begin{lemma}[A useful identity for the affine-invariant geodesic]\label{lemma_useful_identity_for_affine_invariant_geodesic}
Let $\b{X} \in \mathbb{S}_{++}^n$ and
$\b{\Delta} \in T_{\b{X}}\mathbb{S}_{++}^n$.
We define:
\begin{align}
\b{S} := \b{X}^{-1/2}\b{\Delta}\b{X}^{-1/2}.
\end{align}
Consider the affine-invariant geodesic:
\begin{equation}
\b{X}(t)
:=
\b{X}^{1/2}\exp(t\b{S})\b{X}^{1/2}.
\label{equation_geodesic_for_vector_transport_SPD_AI}
\end{equation}
Then:
\begin{align}
\dot{\b{X}}(t)
&=
\b{X}^{1/2}\b{S}\exp(t\b{S})\b{X}^{1/2},
\label{equation_Xdot_for_vector_transport_SPD_AI}
\\
\b{X}(t)^{-1}
&=
\b{X}^{-1/2}\exp(-t\b{S})\b{X}^{-1/2}.
\label{equation_Xinv_for_vector_transport_SPD_AI}
\end{align}
Moreover, for every matrix $\b{B}(t)$ of compatible size,
we have:
\begin{align}
\dot{\b{X}}(t)\b{X}(t)^{-1}
\Big(
\b{X}^{1/2}\b{B}(t)\b{X}^{1/2}
\Big)
&=
\b{X}^{1/2}\b{S}\b{B}(t)\b{X}^{1/2},
\label{equation_identity_1_vector_transport_SPD_AI}
\\
\Big(
\b{X}^{1/2}\b{B}(t)\b{X}^{1/2}
\Big)
\b{X}(t)^{-1}\dot{\b{X}}(t)
&=
\b{X}^{1/2}\b{B}(t)\b{S}\b{X}^{1/2}.
\label{equation_identity_2_vector_transport_SPD_AI}
\end{align}
\end{lemma}

\begin{proof}
From Section \ref{section_exponential_map_spd}, the affine-invariant exponential map
is:
\[
\operatorname{Exp}^{AI}_{\b{X}}(\b{\Delta})
=
\b{X}^{1/2}
\exp\!\big(
\b{X}^{-1/2}\b{\Delta}\b{X}^{-1/2}
\big)
\b{X}^{1/2}.
\]
Therefore, the curve in Eq.
(\ref{equation_geodesic_for_vector_transport_SPD_AI}) is
the geodesic starting from $\b{X}$ with initial tangent
vector $\b{\Delta}$.

Differentiating Eq.
(\ref{equation_geodesic_for_vector_transport_SPD_AI})
with respect to $t$ gives:
\[
\dot{\b{X}}(t)
=
\b{X}^{1/2}
\frac{d}{dt}\exp(t\b{S})
\b{X}^{1/2}.
\]
Because $\b{S}$ commutes with $\exp(t\b{S})$, we have:
\[
\frac{d}{dt}\exp(t\b{S}) = \b{S}\exp(t\b{S}),
\]
so:
\[
\dot{\b{X}}(t)
=
\b{X}^{1/2}\b{S}\exp(t\b{S})\b{X}^{1/2}.
\]
This proves Eq.
(\ref{equation_Xdot_for_vector_transport_SPD_AI}).

Also, since:
\[
\b{X}(t)
=
\b{X}^{1/2}\exp(t\b{S})\b{X}^{1/2},
\]
its inverse is:
\[
\b{X}(t)^{-1}
=
\b{X}^{-1/2}\exp(-t\b{S})\b{X}^{-1/2},
\]
which proves Eq.
(\ref{equation_Xinv_for_vector_transport_SPD_AI}).

Now, let:
\[
\b{V}(t) := \b{X}^{1/2}\b{B}(t)\b{X}^{1/2}.
\]
Using Eqs.
(\ref{equation_Xdot_for_vector_transport_SPD_AI}) and
(\ref{equation_Xinv_for_vector_transport_SPD_AI}), we get:
\begin{align*}
&\dot{\b{X}}(t)\b{X}(t)^{-1}\b{V}(t) \\
&=
\b{X}^{1/2}\b{S}\exp(t\b{S})\b{X}^{1/2}
\b{X}^{-1/2}\\
&\quad\quad\quad\quad\quad\quad\quad\exp(-t\b{S})\b{X}^{-1/2}
\b{X}^{1/2}\b{B}(t)\b{X}^{1/2}
\\
&=
\b{X}^{1/2}\b{S}\exp(t\b{S})\exp(-t\b{S})
\b{B}(t)\b{X}^{1/2}
\\
&=
\b{X}^{1/2}\b{S}\b{B}(t)\b{X}^{1/2}.
\end{align*}
This proves Eq.
(\ref{equation_identity_1_vector_transport_SPD_AI}).

Similarly:
\begin{align*}
&\b{V}(t)\b{X}(t)^{-1}\dot{\b{X}}(t) \\
&=
\b{X}^{1/2}\b{B}(t)\b{X}^{1/2}
\b{X}^{-1/2}\\
&\quad\quad\quad\quad\quad\exp(-t\b{S})\b{X}^{-1/2}
\b{X}^{1/2}\b{S}\exp(t\b{S})\b{X}^{1/2}
\\
&=
\b{X}^{1/2}\b{B}(t)\exp(-t\b{S})
\b{S}\exp(t\b{S})\b{X}^{1/2}.
\end{align*}
Because $\b{S}$ commutes with every polynomial in
$\b{S}$, and hence with $\exp(t\b{S})$, we have:
\[
\exp(-t\b{S})\b{S}\exp(t\b{S}) = \b{S}.
\]
Therefore:
\[
\b{V}(t)\b{X}(t)^{-1}\dot{\b{X}}(t)
=
\b{X}^{1/2}\b{B}(t)\b{S}\b{X}^{1/2}.
\]
This proves Eq.
(\ref{equation_identity_2_vector_transport_SPD_AI}).
\end{proof}

\begin{proposition}[Vector transport in SPD manifold under the affine-invariant metric]
Let $\b{X} \in \mathbb{S}_{++}^n$ and let
$\b{\Delta}, \b{\Xi} \in T_{\b{X}}\mathbb{S}_{++}^n
= \mathbb{S}^n$.
We define:
\begin{align}
&\boxed{\b{S} := \b{X}^{-1/2}\b{\Delta}\b{X}^{-1/2},} \\
&\boxed{\b{A} := \b{X}^{-1/2}\b{\Xi}\b{X}^{-1/2}.}
\end{align}
Let:
\begin{equation}
\b{Y}
:=
\operatorname{Exp}^{AI}_{\b{X}}(\b{\Delta})
=
\b{X}^{1/2}\exp(\b{S})\b{X}^{1/2}.
\label{equation_endpoint_vector_transport_SPD_AI}
\end{equation}
We define:
\begin{equation}
\boxed{
\mathcal{T}^{AI}_{\b{\Delta}}(\b{\Xi})
:=
\b{X}^{1/2}
\exp\!\Big(\frac{1}{2}\b{S}\Big)
\b{A}
\exp\!\Big(\frac{1}{2}\b{S}\Big)
\b{X}^{1/2}.
}
\label{equation_vector_transport_SPD_manifold_AI}
\end{equation}
Then, $\mathcal{T}^{AI}_{\b{\Delta}}(\b{\Xi})$ is the
parallel transport of $\b{\Xi}$ along the affine-invariant
geodesic from $\b{X}$ to $\b{Y}$. Hence, it defines a
valid vector transport on $\mathbb{S}_{++}^n$ under the
affine-invariant metric.
\end{proposition}

\begin{proof}
We prove the result in two parts. First, we show that Eq.
(\ref{equation_vector_transport_SPD_manifold_AI}) is
indeed the parallel transport along the affine-invariant
geodesic. Then, we verify the vector transport axioms.

\textbf{Step 1: Define the geodesic and the transported
vector field.}

Consider the affine-invariant geodesic:
\[
\b{X}(t)
=
\b{X}^{1/2}\exp(t\b{S})\b{X}^{1/2},
\qquad t \in [0,1].
\]
Then:
\[
\b{X}(0)=\b{X},
\qquad
\b{X}(1)=\b{Y}.
\]
Now, we define the matrix-valued field along this curve by:
\begin{equation}
\b{V}(t)
:=
\b{X}^{1/2}
\exp\!\Big(\frac{t}{2}\b{S}\Big)
\b{A}
\exp\!\Big(\frac{t}{2}\b{S}\Big)
\b{X}^{1/2}.
\label{equation_Vt_vector_transport_SPD_AI}
\end{equation}
At $t=0$, we have:
\begin{align*}
\b{V}(0)
&=
\b{X}^{1/2}\b{A}\b{X}^{1/2}
\\
&=
\b{X}^{1/2}\b{X}^{-1/2}\b{\Xi}\b{X}^{-1/2}\b{X}^{1/2}
=
\b{\Xi}.
\end{align*}
At $t=1$, Eq.
(\ref{equation_Vt_vector_transport_SPD_AI}) gives:
\[
\b{V}(1)
=
\b{X}^{1/2}
\exp\!\Big(\frac{1}{2}\b{S}\Big)
\b{A}
\exp\!\Big(\frac{1}{2}\b{S}\Big)
\b{X}^{1/2}
=
\mathcal{T}^{AI}_{\b{\Delta}}(\b{\Xi}).
\]

\textbf{Step 2: Show that $\b{V}(t)$ is tangent-valued.}

Because $\b{\Xi} \in \mathbb{S}^n$, the matrix
$\b{A}=\b{X}^{-1/2}\b{\Xi}\b{X}^{-1/2}$ is symmetric.
Also, $\b{S}$ is symmetric because
$\b{\Delta} \in \mathbb{S}^n$, hence
$\exp(\frac{t}{2}\b{S})$ is symmetric.
Therefore:
\[
\b{V}(t)^\top
=
\b{X}^{1/2}
\exp\!\Big(\frac{t}{2}\b{S}\Big)
\b{A}
\exp\!\Big(\frac{t}{2}\b{S}\Big)
\b{X}^{1/2}
=
\b{V}(t).
\]
So:
\[
\b{V}(t) \in \mathbb{S}^n
=
T_{\b{X}(t)}\mathbb{S}_{++}^n.
\]

\textbf{Step 3: Differentiate $\b{V}(t)$.}

Differentiating Eq.
(\ref{equation_Vt_vector_transport_SPD_AI}) gives:
\begin{align*}
\dot{\b{V}}(t)
&=
\b{X}^{1/2}
\frac{d}{dt}
\left(
\exp\!\Big(\frac{t}{2}\b{S}\Big)
\b{A}
\exp\!\Big(\frac{t}{2}\b{S}\Big)
\right)
\b{X}^{1/2}
\\
&=
\frac{1}{2}
\b{X}^{1/2}
\left(
\b{S}\exp\!\Big(\frac{t}{2}\b{S}\Big)
\b{A}
\exp\!\Big(\frac{t}{2}\b{S}\Big)
\right)
\b{X}^{1/2}
\\
&\quad
+
\frac{1}{2}
\b{X}^{1/2}
\left(
\exp\!\Big(\frac{t}{2}\b{S}\Big)
\b{A}
\exp\!\Big(\frac{t}{2}\b{S}\Big)\b{S}
\right)
\b{X}^{1/2}.
\end{align*}
Let:
\[
\b{B}(t)
:=
\exp\!\Big(\frac{t}{2}\b{S}\Big)
\b{A}
\exp\!\Big(\frac{t}{2}\b{S}\Big).
\]
Then:
\begin{equation}
\dot{\b{V}}(t)
=
\frac{1}{2}
\b{X}^{1/2}
\big(
\b{S}\b{B}(t) + \b{B}(t)\b{S}
\big)
\b{X}^{1/2}.
\label{equation_Vdot_vector_transport_SPD_AI}
\end{equation}

\textbf{Step 4: Use the affine-invariant Levi-Civita
connection.}

From Proposition \ref{proposition_Levi_Civita_connection_spd_affine_invariant_metric}, the affine-invariant Levi-Civita
connection on the SPD manifold is:
\begin{align*}
\nabla^{AI}_{\b{\Delta}_1}&\b{\Delta}_2(\b{X})
=
D\b{\Delta}_2(\b{X})[\b{\Delta}_1]
\\
&-
\frac{1}{2}
\Big(
\b{\Delta}_1\b{X}^{-1}\b{\Delta}_2(\b{X})
+
\b{\Delta}_2(\b{X})\b{X}^{-1}\b{\Delta}_1
\Big).
\end{align*}
Therefore, along the curve $\b{X}(t)$, the parallel
transport equation is:
\begin{equation}\label{equation_parallel_transport_equation_SPD_AI}
\begin{aligned}
&\nabla^{AI}_{\dot{\b{X}}(t)}\b{V}(t) \\
&=
\dot{\b{V}}(t)
-
\frac{1}{2}
\Big(
\dot{\b{X}}(t)\b{X}(t)^{-1}\b{V}(t)
+
\b{V}(t)\b{X}(t)^{-1}\dot{\b{X}}(t)
\Big)
\\
&=
\b{0}.
\end{aligned}
\end{equation}

Now, by Lemma \ref{lemma_useful_identity_for_affine_invariant_geodesic}, with $\b{V}(t)=\b{X}^{1/2}\b{B}(t)\b{X}^{1/2}$,
we have:
\begin{align*}
\dot{\b{X}}(t)\b{X}(t)^{-1}\b{V}(t)
&=
\b{X}^{1/2}\b{S}\b{B}(t)\b{X}^{1/2},
\\
\b{V}(t)\b{X}(t)^{-1}\dot{\b{X}}(t)
&=
\b{X}^{1/2}\b{B}(t)\b{S}\b{X}^{1/2}.
\end{align*}
Substituting these and Eq.
(\ref{equation_Vdot_vector_transport_SPD_AI}) into Eq.
(\ref{equation_parallel_transport_equation_SPD_AI}) gives:
\begin{align*}
\nabla^{AI}_{\dot{\b{X}}(t)}\b{V}(t)
&=
\frac{1}{2}
\b{X}^{1/2}
\big(
\b{S}\b{B}(t) + \b{B}(t)\b{S}
\big)
\b{X}^{1/2}
\\
&\quad
-
\frac{1}{2}
\b{X}^{1/2}
\big(
\b{S}\b{B}(t) + \b{B}(t)\b{S}
\big)
\b{X}^{1/2}
\\
&=
\b{0}.
\end{align*}
Hence, $\b{V}(t)$ is parallel along the affine-invariant
geodesic.

Therefore, $\b{V}(1)$ is exactly the parallel transport of
$\b{\Xi}$ from $\b{X}$ to $\b{Y}$, and this is the formula
in Eq. (\ref{equation_vector_transport_SPD_manifold_AI}).

\textbf{Step 5: Verify the axioms of vector transport.}

\textit{Well-defined mapping:}
Since $\b{V}(1)$ is symmetric, we have:
\[
\mathcal{T}^{AI}_{\b{\Delta}}(\b{\Xi})
=
\b{V}(1)
\in
\mathbb{S}^n
=
T_{\b{Y}}\mathbb{S}_{++}^n,
\]
where $\b{Y}=\operatorname{Exp}^{AI}_{\b{X}}(\b{\Delta})$.

\textit{Linearity in $\b{\Xi}$:}
The mapping $\b{\Xi} \mapsto \b{A}
=\b{X}^{-1/2}\b{\Xi}\b{X}^{-1/2}$ is linear, and the map:
\[
\b{A}
\mapsto
\b{X}^{1/2}
\exp\!\Big(\frac{1}{2}\b{S}\Big)
\b{A}
\exp\!\Big(\frac{1}{2}\b{S}\Big)
\b{X}^{1/2}
\]
is also linear. Hence:
\[
\mathcal{T}^{AI}_{\b{\Delta}}
(a\b{\Xi}_1 + b\b{\Xi}_2)
=
a\,\mathcal{T}^{AI}_{\b{\Delta}}(\b{\Xi}_1)
+
b\,\mathcal{T}^{AI}_{\b{\Delta}}(\b{\Xi}_2).
\]

\textit{Consistency at zero:}
If $\b{\Delta}=\b{0}$, then $\b{S}=\b{0}$ and
$\exp(\frac{1}{2}\b{S})=\b{I}_n$. Therefore:
\begin{align*}
\mathcal{T}^{AI}_{\b{0}}(\b{\Xi})
&=
\b{X}^{1/2}\b{I}_n
\b{X}^{-1/2}\b{\Xi}\b{X}^{-1/2}
\b{I}_n\b{X}^{1/2}
\\
&=
\b{\Xi}.
\end{align*}

Thus, $\mathcal{T}^{AI}$ satisfies all axioms of vector
transport. This proves the proposition.
\end{proof}

\begin{remark}[Interpretation of the affine-invariant vector transport on SPD manifold]
The affine-invariant transport first whitens the tangent
vector $\b{\Xi}$ at $\b{X}$ by:
\[
\b{\Xi}
\mapsto
\b{X}^{-1/2}\b{\Xi}\b{X}^{-1/2}.
\]
Then, it moves this whitened tangent vector by half of the
geodesic motion on the left and half on the right through
the factors $\exp(\frac{1}{2}\b{S})$. Finally, it maps the
result back to the ambient matrix coordinates by congruence
with $\b{X}^{1/2}$.

Therefore, unlike the Euclidean case where the same
symmetric matrix can simply be reused at the new point,
the affine-invariant geometry requires a base-point-dependent
conjugation that respects the Levi-Civita connection of the
metric.
\end{remark}

\subsubsection{Vector Transport by Differential Retraction in SPD Manifold}

As discussed in Section \ref{section_differentiated_retraction}, a natural way to construct
a vector transport is by differentiating a chosen retraction.
For the SPD manifold, we derive this construction for the
Euclidean and affine-invariant retractions introduced in
Section \ref{section_retraction_spd}.

Recall that for the SPD manifold, we have:
\[
T_{\b{X}}\mathbb{S}_{++}^n = \mathbb{S}^n,
\]
for every $\b{X}\in \mathbb{S}_{++}^n$. Therefore, every tangent
vector is a symmetric matrix.

\begin{proposition}[Differentiated Euclidean retraction on SPD manifold]
Let $\b{X}\in \mathbb{S}_{++}^n$ and let:
\[
\b{\Xi}, \b{\Delta}\in T_{\b{X}}\mathbb{S}_{++}^n = \mathbb{S}^n.
\]
Consider the Euclidean retraction in Eq. \eqref{equation_retraction_SPD_manifold_Euclidean}:
\[
\operatorname{Ret}_{\b{X}}^{\mathrm{Euc}}(\b{\Delta})
=
\b{X}+\b{\Delta}.
\]
We define:
\[
\b{Y}
:=
\operatorname{Ret}_{\b{X}}^{\mathrm{Euc}}(\b{\Delta})
=
\b{X}+\b{\Delta}.
\]
Then, the vector transport associated with the differentiated
Euclidean retraction is:
\begin{align}
\boxed{
\mathcal{T}_{\b{\Delta}}^{\mathrm{Euc}}(\b{\Xi})
:=
\left.\frac{d}{dt}
\operatorname{Ret}_{\b{X}}^{\mathrm{Euc}}(\b{\Delta}+t\b{\Xi})
\right|_{t=0}
=
\b{\Xi}.
}
\end{align}
Hence, the differentiated-retraction vector transport under
the Euclidean metric is the identity map.
\end{proposition}

\begin{proof}
By definition of differentiated retraction, we consider the
curve:
\[
\b{Z}(t)
:=
\operatorname{Ret}_{\b{X}}^{\mathrm{Euc}}(\b{\Delta}+t\b{\Xi}).
\]
Using Eq. \eqref{equation_retraction_SPD_manifold_Euclidean},
we have:
\[
\b{Z}(t)
=
\b{X}+(\b{\Delta}+t\b{\Xi})
=
\b{X}+\b{\Delta}+t\b{\Xi}.
\]
Differentiating with respect to $t$ gives:
\[
\dot{\b{Z}}(t)=\b{\Xi}.
\]
Therefore:
\[
\left.\frac{d}{dt}
\operatorname{Ret}_{\b{X}}^{\mathrm{Euc}}(\b{\Delta}+t\b{\Xi})
\right|_{t=0}
=
\dot{\b{Z}}(0)
=
\b{\Xi}.
\]
Hence:
\[
\mathcal{T}_{\b{\Delta}}^{\mathrm{Euc}}(\b{\Xi})=\b{\Xi}.
\]

It remains to verify that this is indeed a vector transport.

\textbf{Step 1: Tangent-valuedness.}
Because $\b{\Xi}\in \mathbb{S}^n$, the transported vector
$\mathcal{T}_{\b{\Delta}}^{\mathrm{Euc}}(\b{\Xi})=\b{\Xi}$ is symmetric.
Since
\[
T_{\b{Y}}\mathbb{S}_{++}^n=\mathbb{S}^n,
\]
we have:
\[
\mathcal{T}_{\b{\Delta}}^{\mathrm{Euc}}(\b{\Xi})
\in
T_{\b{Y}}\mathbb{S}_{++}^n.
\]

\textbf{Step 2: Linearity in $\b{\Xi}$.}
For any $a,b\in\mathbb{R}$ and
$\b{\Xi}_1,\b{\Xi}_2\in T_{\b{X}}\mathbb{S}_{++}^n$, we have:
\begin{align*}
\mathcal{T}_{\b{\Delta}}^{\mathrm{Euc}}(a\b{\Xi}_1+b\b{\Xi}_2)
&=
a\b{\Xi}_1+b\b{\Xi}_2 \\
&=
a\,\mathcal{T}_{\b{\Delta}}^{\mathrm{Euc}}(\b{\Xi}_1)
+
b\,\mathcal{T}_{\b{\Delta}}^{\mathrm{Euc}}(\b{\Xi}_2).
\end{align*}

\textbf{Step 3: Identity at zero step.}
If $\b{\Delta}=\b{0}$, then:
\[
\mathcal{T}_{\b{0}}^{\mathrm{Euc}}(\b{\Xi})=\b{\Xi}.
\]

Therefore, the differentiated Euclidean retraction defines a
valid vector transport on $\mathbb{S}_{++}^n$ under the
Euclidean metric.
\end{proof}

\begin{lemma}[Symmetry of the differentiated affine-invariant retraction term]\label{lemma_symmetry_differentiated_affine_invariant_retraction}
Let $\b{X}\in \mathbb{S}_{++}^n$ and let:
\[
\b{\Xi},\b{\Delta}\in \mathbb{S}^n.
\]
Then, the following matrix is symmetric:
\begin{align}
\b{\Xi}\b{X}^{-1}\b{\Delta}
+
\b{\Delta}\b{X}^{-1}\b{\Xi} \in \mathbb{S}^n.
\end{align}
\end{lemma}

\begin{proof}
Because $\b{X}\in \mathbb{S}_{++}^n$, its inverse $\b{X}^{-1}$
is also symmetric. Since $\b{\Xi}$ and $\b{\Delta}$ are symmetric,
we have:
\[
(\b{\Xi}\b{X}^{-1}\b{\Delta})^\top
=
\b{\Delta}^\top(\b{X}^{-1})^\top\b{\Xi}^\top
=
\b{\Delta}\b{X}^{-1}\b{\Xi},
\]
and similarly,
\[
(\b{\Delta}\b{X}^{-1}\b{\Xi})^\top
=
\b{\Xi}\b{X}^{-1}\b{\Delta}.
\]
Therefore:
\begin{align*}
\left(
\b{\Xi}\b{X}^{-1}\b{\Delta}
+
\b{\Delta}\b{X}^{-1}\b{\Xi}
\right)^\top
&=
(\b{\Xi}\b{X}^{-1}\b{\Delta})^\top
+
(\b{\Delta}\b{X}^{-1}\b{\Xi})^\top \\
&=
\b{\Delta}\b{X}^{-1}\b{\Xi}
+
\b{\Xi}\b{X}^{-1}\b{\Delta}.
\end{align*}
This is exactly the original matrix. Hence, it is symmetric.
\end{proof}

\begin{proposition}[Differentiated affine-invariant retraction on SPD manifold]
Let $\b{X}\in \mathbb{S}_{++}^n$ and let:
\[
\b{\Xi}, \b{\Delta}\in T_{\b{X}}\mathbb{S}_{++}^n = \mathbb{S}^n.
\]
Consider the affine-invariant retraction in
Eq. \eqref{equation_retraction_SPD_manifold_affine_invariant}:
\[
\operatorname{Ret}_{\b{X}}^{AI}(\b{\Delta})
=
\b{X}
+
\b{\Delta}
+
\frac{1}{2}\b{\Delta}\b{X}^{-1}\b{\Delta}.
\]
We define:
\[
\b{Y}
:=
\operatorname{Ret}_{\b{X}}^{AI}(\b{\Delta}).
\]
Then, the vector transport associated with the differentiated
affine-invariant retraction is:
\begin{equation}
\boxed{
\begin{aligned}
\mathcal{T}_{\b{\Delta}}^{AI}(\b{\Xi})
&:=
\left.\frac{d}{dt}
\operatorname{Ret}_{\b{X}}^{AI}(\b{\Delta}+t\b{\Xi})
\right|_{t=0}
\\
&=
\b{\Xi}
+
\frac{1}{2}
\left(
\b{\Xi}\b{X}^{-1}\b{\Delta}
+
\b{\Delta}\b{X}^{-1}\b{\Xi}
\right).
\end{aligned}
}
\end{equation}
Hence, the differentiated-retraction vector transport under
the affine-invariant metric is given by the first-order
variation of the quadratic retraction term.
\end{proposition}

\begin{proof}
By definition of differentiated retraction, let:
\[
\b{Z}(t)
:=
\operatorname{Ret}_{\b{X}}^{AI}(\b{\Delta}+t\b{\Xi}).
\]
Using Eq. \eqref{equation_retraction_SPD_manifold_affine_invariant},
we obtain:
\[
\b{Z}(t)
=
\b{X}
+
(\b{\Delta}+t\b{\Xi})
+
\frac{1}{2}
(\b{\Delta}+t\b{\Xi})\b{X}^{-1}(\b{\Delta}+t\b{\Xi}).
\]
We now expand the quadratic term:
\begin{align*}
&(\b{\Delta}+t\b{\Xi})\b{X}^{-1}(\b{\Delta}+t\b{\Xi}) \\
&=
\b{\Delta}\b{X}^{-1}\b{\Delta}
+
t\,\b{\Xi}\b{X}^{-1}\b{\Delta}
+
t\,\b{\Delta}\b{X}^{-1}\b{\Xi}
+
t^2\b{\Xi}\b{X}^{-1}\b{\Xi}.
\end{align*}
Substituting this into $\b{Z}(t)$ gives:
\begin{align*}
\b{Z}(t)
&=
\b{X}
+
\b{\Delta}
+
t\b{\Xi}
+
\frac{1}{2}\b{\Delta}\b{X}^{-1}\b{\Delta} \\
&\quad
+
\frac{t}{2}
\Big(
\b{\Xi}\b{X}^{-1}\b{\Delta}
+
\b{\Delta}\b{X}^{-1}\b{\Xi}
\Big)
+
\frac{t^2}{2}\b{\Xi}\b{X}^{-1}\b{\Xi}.
\end{align*}
Differentiating with respect to $t$ yields:
\[
\dot{\b{Z}}(t)
=
\b{\Xi}
+
\frac{1}{2}
\Big(
\b{\Xi}\b{X}^{-1}\b{\Delta}
+
\b{\Delta}\b{X}^{-1}\b{\Xi}
\Big)
+
t\,\b{\Xi}\b{X}^{-1}\b{\Xi}.
\]
Evaluating at $t=0$, we get:
\begin{align*}
\frac{d}{dt}
\operatorname{Ret}_{\b{X}}^{AI}(\b{\Delta}+&t\b{\Xi})
\Big|_{t=0}
\\
&=
\b{\Xi}
+
\frac{1}{2}
\left(
\b{\Xi}\b{X}^{-1}\b{\Delta}
+
\b{\Delta}\b{X}^{-1}\b{\Xi}
\right).
\end{align*}
Therefore:
\[
\mathcal{T}_{\b{\Delta}}^{AI}(\b{\Xi})
=
\b{\Xi}
+
\frac{1}{2}
\left(
\b{\Xi}\b{X}^{-1}\b{\Delta}
+
\b{\Delta}\b{X}^{-1}\b{\Xi}
\right).
\]

It remains to verify that this is a valid vector transport.

\textbf{Step 1: Tangent-valuedness.}
By Lemma \ref{lemma_symmetry_differentiated_affine_invariant_retraction}, the matrix:
\[
\b{\Xi}\b{X}^{-1}\b{\Delta}
+
\b{\Delta}\b{X}^{-1}\b{\Xi}
\]
is symmetric. Since $\b{\Xi}$ is also symmetric, it follows
that $\mathcal{T}_{\b{\Delta}}^{AI}(\b{\Xi})$ is symmetric.
Because
\[
T_{\b{Y}}\mathbb{S}_{++}^n=\mathbb{S}^n,
\]
we conclude:
\[
\mathcal{T}_{\b{\Delta}}^{AI}(\b{\Xi})
\in
T_{\b{Y}}\mathbb{S}_{++}^n.
\]

\textbf{Step 2: Linearity in $\b{\Xi}$.}
Let $a,b\in\mathbb{R}$ and
$\b{\Xi}_1,\b{\Xi}_2\in T_{\b{X}}\mathbb{S}_{++}^n$.
Then:
\begin{align*}
&\mathcal{T}_{\b{\Delta}}^{AI}(a\b{\Xi}_1+b\b{\Xi}_2) \\
&=
a\b{\Xi}_1+b\b{\Xi}_2 \\
&\quad
+\frac{1}{2}
\Big(
(a\b{\Xi}_1+b\b{\Xi}_2)\b{X}^{-1}\b{\Delta}
+
\b{\Delta}\b{X}^{-1}(a\b{\Xi}_1+b\b{\Xi}_2)
\Big) \\
&=
a\left[
\b{\Xi}_1
+
\frac{1}{2}
\left(
\b{\Xi}_1\b{X}^{-1}\b{\Delta}
+
\b{\Delta}\b{X}^{-1}\b{\Xi}_1
\right)
\right] \\
&\quad
+
b\left[
\b{\Xi}_2
+
\frac{1}{2}
\left(
\b{\Xi}_2\b{X}^{-1}\b{\Delta}
+
\b{\Delta}\b{X}^{-1}\b{\Xi}_2
\right)
\right] \\
&=
a\,\mathcal{T}_{\b{\Delta}}^{AI}(\b{\Xi}_1)
+
b\,\mathcal{T}_{\b{\Delta}}^{AI}(\b{\Xi}_2).
\end{align*}

\textbf{Step 3: Identity at zero step.}
If $\b{\Delta}=\b{0}$, then:
\[
\mathcal{T}_{\b{0}}^{AI}(\b{\Xi})
=
\b{\Xi}
+
\frac{1}{2}
\left(
\b{\Xi}\b{X}^{-1}\b{0}
+
\b{0}\b{X}^{-1}\b{\Xi}
\right)
=
\b{\Xi}.
\]

Therefore, the differentiated affine-invariant retraction
defines a valid vector transport on $\mathbb{S}_{++}^n$ under
the affine-invariant metric.
\end{proof}

\begin{remark}[Interpretation of differentiated retraction on SPD manifold]
For the Euclidean metric, the retraction is ordinary matrix
addition, so its differentiated vector transport is trivial:
the vector does not change at all.

For the affine-invariant metric, the retraction contains the
quadratic correction term
$\frac{1}{2}\b{\Delta}\b{X}^{-1}\b{\Delta}$.
Therefore, differentiating the retraction with respect to the
transported vector $\b{\Xi}$ produces the additional
symmetrized correction:
\[
\frac{1}{2}
\left(
\b{\Xi}\b{X}^{-1}\b{\Delta}
+
\b{\Delta}\b{X}^{-1}\b{\Xi}
\right).
\]
Hence, unlike the Euclidean case, the transported vector is
adjusted according to both the step direction $\b{\Delta}$
and the local geometry encoded by $\b{X}^{-1}$.
\end{remark}

\begin{remark}[Special case at zero step]
If $\b{\Delta}=\b{0}$, both differentiated retractions reduce
to the identity:
\[
\mathcal{T}_{\b{0}}^{\mathrm{Euc}}(\b{\Xi})=\b{\Xi},
\qquad
\mathcal{T}_{\b{0}}^{AI}(\b{\Xi})=\b{\Xi}.
\]
Therefore, the differentiated retraction agrees with the
required identity behavior of vector transport at the base
point.
\end{remark}

\section{Important Software Toolboxes and Textbooks for Riemannian Optimization}\label{section_toolboxes_for_riemannian_optimization}

In Sections \ref{section_manifold_valued_optimization} and \ref{section_important_Riemannian_matrix_manifolds}, we derived the main geometric and
algorithmic ingredients of Riemannian manifold-valued
optimization, including the Riemannian gradient, Riemannian
Hessian, exponential map, retraction, and vector transport.
In practice, software toolboxes for Riemannian optimization
are useful because they implement these ingredients in a
numerically stable and reusable way for many manifolds and
optimization algorithms.

This section briefly explains what a practical Riemannian
optimization toolbox must provide, based on the equations
derived in this monograph, and then lists some important
software toolboxes and references.

\subsection{Important Toolboxes for Riemannian Optimization}\label{section_important_toolboxes_riemannian_optimization}

The previous sections show that Riemannian optimization
is not merely Euclidean optimization with a constraint.
Rather, every iteration requires geometric operations that
respect the manifold structure.

\begin{proposition}[Core ingredients of a practical Riemannian optimization toolbox]\label{proposition_ingredients_riemannian_optimization_toolbox}
\label{proposition_core_ingredients_toolbox_riemannian_optimization}
Let $(\mathcal{M},g)$ be a Riemannian manifold and let
$f:\mathcal{M}\rightarrow\mathbb{R}$ be a smooth cost function.
Suppose we want to implement the first-order and second-order
optimization methods derived in Sections \ref{section_riemannian_gradient_descent} and \ref{section_riemannian_second_order_optimization}.
Then, a practical Riemannian optimization toolbox must
provide, either explicitly or implicitly, the following
ingredients:
\begin{enumerate}
\item a representation of a point $\b{p}\in\mathcal{M}$ and of a
tangent vector $\b{\eta}\in T_{\b{p}}\mathcal{M}$,

\item a way to evaluate the metric
$g_{\b{p}}(\b{\xi},\b{\eta})$ for
$\b{\xi},\b{\eta}\in T_{\b{p}}\mathcal{M}$,

\item a way to compute the Riemannian gradient
$\operatorname{grad}f(\b{p})\in T_{\b{p}}\mathcal{M}$,

\item a way to compute the Riemannian Hessian
$\operatorname{Hess}f(\b{p})[\b{\eta}]$, or at least a
quasi-Newton approximation of it,

\item a way to move from one point to another point on the
manifold along a tangent vector, typically through the
exponential map $\operatorname{Exp}_{\b{p}}(\b{\eta})$ or a
retraction $\operatorname{Ret}_{\b{p}}(\b{\eta})$,

\item a way to move tangent vectors between tangent spaces,
namely a vector transport
$\mathcal{T}_{\b{\eta}}(\b{\xi})
:
T_{\b{p}}\mathcal{M}
\rightarrow
T_{\operatorname{Ret}_{\b{p}}(\b{\eta})}\mathcal{M}$,
whenever the algorithm needs to reuse search directions,
gradients, or Hessian approximations across iterations.
\end{enumerate}

Therefore, a software toolbox for Riemannian optimization
must encode the manifold geometry together with the
optimization algorithms.
\end{proposition}

\begin{proof}
We prove this by inspecting the equations derived earlier in
the monograph.

\textbf{Step 1: Necessity of point and tangent-vector
representations.}
Optimization on a manifold proceeds through iterates
$\b{p}_{\nu}\in\mathcal{M}$ and tangent directions
$\b{\eta}_{\nu}\in T_{\b{p}_{\nu}}\mathcal{M}$.
Hence, any implementation must be able to store and manipulate
both points on the manifold and tangent vectors attached to
those points.

\textbf{Step 2: Necessity of the metric.}
By definition of the Riemannian gradient in Section \ref{section_Riemannian_gradient},
the gradient is characterized by:
\[
Df(\b{p})[\b{\xi}]
=
g_{\b{p}}(\operatorname{grad}f(\b{p}),\b{\xi}),
\qquad
\forall \b{\xi}\in T_{\b{p}}\mathcal{M}.
\]
Therefore, even defining the gradient requires access to the
metric tensor. Furthermore, line-search criteria, descent
tests, and curvature-related quantities also depend on the
metric.

\textbf{Step 3: Necessity of the Riemannian gradient.}
The Riemannian gradient descent update derived in
Section \ref{section_RGD_update} is:
\[
\b{p}_{\nu+1}
=
\operatorname{Exp}_{\b{p}_{\nu}}\!
\big(\!
-\!\lambda\,\operatorname{grad}f(\b{p}_{\nu})
\big).
\]
Its practical retraction-based implementation in
Section \ref{section_numerical_RGD_retraction} is:
\[
\b{p}_{\nu+1}
=
\operatorname{Ret}_{\b{p}_{\nu}}\!
\big(\!
-\!\lambda\,\operatorname{grad}f(\b{p}_{\nu})
\big).
\]
Hence, no first-order Riemannian method can be implemented
without computing $\operatorname{grad}f(\b{p}_{\nu})$.

\textbf{Step 4: Necessity of the Hessian or its approximation.}
The Riemannian Newton method in Section \ref{section_riemannian_second_order_optimization} solves the
Newton equation:
\[
\operatorname{Hess}f(\b{p}_{\nu})[\b{\eta}_{\nu}]
=
-
\operatorname{grad}f(\b{p}_{\nu}).
\]
Thus, second-order methods require the action of the
Riemannian Hessian. Even when the exact Hessian is not used,
quasi-Newton methods such as RBFGS and RLBFGS replace it by
an approximation. Therefore, a toolbox must provide either
the Hessian itself or a mechanism for approximating it.

\textbf{Step 5: Necessity of exponential map or retraction.}
The iterate update must remain on the manifold. As discussed
in Sections \ref{section_exponential_logarithm_map} and \ref{section_retraction}, direct Euclidean addition
$\b{p}_{\nu}+\b{\eta}_{\nu}$ is generally invalid because the
result may leave $\mathcal{M}$. Therefore, the update must be
performed using either the exponential map or a retraction:
\[
\b{p}_{\nu+1}
=
\operatorname{Exp}_{\b{p}_{\nu}}(\b{\eta}_{\nu})
\quad\text{or}\quad
\b{p}_{\nu+1}
=
\operatorname{Ret}_{\b{p}_{\nu}}(\b{\eta}_{\nu}).
\]
Hence, a practical toolbox must implement at least one of
these two operations.

\textbf{Step 6: Necessity of vector transport.}
For methods that reuse tangent vectors across iterations,
such as conjugate-gradient and quasi-Newton methods, vectors
living in different tangent spaces must be compared or
combined. However:
\begin{align*}
&\operatorname{grad}f(\b{p}_{\nu}) \in T_{\b{p}_{\nu}}\mathcal{M}, \\
&\operatorname{grad}f(\b{p}_{\nu+1}) \in T_{\b{p}_{\nu+1}}\mathcal{M},
\end{align*}
and these two tangent spaces are generally different. Hence,
to combine old and new search information, we require a
mapping from one tangent space to another, i.e., vector
transport.

Therefore, all six ingredients are required, either directly
or through internally equivalent implementations, for a
practical Riemannian optimization toolbox.
\end{proof}

\begin{corollary}[Retraction is often sufficient in practice]
\label{corollary_retraction_often_sufficient_in_practice}
Let $\operatorname{Ret}_{\b{p}}:T_{\b{p}}\mathcal{M}\to\mathcal{M}$
be a retraction. If a toolbox provides this retraction and
can differentiate it, then it can construct a valid vector
transport by differentiated retraction (see Eq. (\ref{equation_differential_retraction})):
\[
\mathcal{T}_{\b{\eta}}(\b{\xi})
:=
D\operatorname{Ret}_{\b{p}}(\b{\eta})[\b{\xi}],
\]
using the natural identification
$T_{\b{\eta}}(T_{\b{p}}\mathcal{M})\cong T_{\b{p}}\mathcal{M}$.
Therefore, many practical toolboxes need not implement a
separate closed-form vector transport for every manifold.
\end{corollary}

\begin{proof}
This follows directly from Definition \ref{definition_differentiated_retraction} and Proposition \ref{proposition_differentiated_retraction_vector_transport}
in Section \ref{section_differentiated_retraction}. There, we showed that the differentiated
retraction is well-defined, linear in $\b{\xi}$, and satisfies
the identity condition at zero:
\[
\mathcal{T}_{\b{0}}(\b{\xi}) = \b{\xi}.
\]
Hence, it defines a valid vector transport.
\end{proof}

\begin{remark}[Meaning of the Proposition \ref{proposition_ingredients_riemannian_optimization_toolbox} for software design]
Proposition
\ref{proposition_core_ingredients_toolbox_riemannian_optimization}
shows that a Riemannian optimization toolbox has two
coupled responsibilities:
\begin{enumerate}
\item First, it must provide the
\textbf{geometry layer}, including manifold representation,
metric, gradient conversion, retraction, and vector transport.
\item Second, it must provide the \textbf{optimization layer},
including gradient descent, conjugate-gradient, Newton, and
quasi-Newton routines that use those geometric primitives.
Therefore, these toolboxes are not merely collections of
solvers; they are software realizations of the geometric
equations derived throughout this monograph.
\end{enumerate}
\end{remark}

Some important software toolboxes are the following:

\begin{itemize}
\item \textbf{Manopt} \cite{boumal2014manopt}, mainly for
MATLAB\footnote{\url{https://github.com/NicolasBoumal/manopt}}. It is a toolbox for optimization on manifolds and
linear spaces, with support for manifold definitions,
state-of-the-art solvers, and automatic differentiation
features. It is one of the best-known general-purpose
packages in Riemannian optimization.

\item \textbf{PyManopt} \cite{townsend2016pymanopt}, for Python\footnote{\url{https://github.com/pymanopt/pymanopt}}.
It is a Python toolbox for optimization on Riemannian
manifolds with support for automatic differentiation, which
makes implementation of gradients and Hessians easier.

\item \textbf{GeomStats} \cite{miolane2020geomstats}, for Python\footnote{\url{https://github.com/geomstats/geomstats}}.
It is broader than a pure optimization toolbox. It provides
manifolds, Riemannian metrics, exponential and logarithm
maps, geodesics, and learning/statistical tools on manifolds.

\item \textbf{Geoopt} \cite{kochurov2020geoopt}, for PyTorch\footnote{\url{https://github.com/geoopt/geoopt}}.
It is designed for manifold-aware optimization in deep
learning, especially when one wants Riemannian optimization
inside neural network models.

\item \textbf{ROPTLIB} \cite{huang2017roptlib}, mainly with a
C++ core\footnote{\url{https://github.com/whuang08/ROPTLIB}}. It is a Riemannian manifold optimization library
with object-oriented design and efficient implementation of
generic algorithms.

\item \textbf{StochMan} \cite{softwareStochMan}, for Python\footnote{\url{https://github.com/MachineLearningLifeScience/stochman}}.
It focuses on stochastic manifolds and computations on
manifolds learned from noisy finite data.

\item \textbf{MixEst} \cite{hosseini2015mixest}, mainly
for MATLAB\footnote{\url{https://github.com/utvisionlab/mixest}}. It is designed for mixture-model parameter
estimation and includes Riemannian optimization methods
such as Riemannian LBFGS.
\end{itemize}

\begin{remark}[Which toolbox is suitable for which purpose]
If the goal is classical Riemannian optimization in MATLAB,
Manopt is one of the most standard choices. If the goal is
rapid prototyping in Python with automatic differentiation,
PyManopt is a natural choice. If the goal is manifold-based
statistics and geometric machine learning beyond optimization
alone, GeomStats is very useful. If the goal is integrating
Riemannian optimization into deep learning pipelines, Geoopt
is especially convenient because of its PyTorch interface.
If high computational efficiency in a compiled language is
important, ROPTLIB is an attractive option.
If a Riemannian LBFGS algorithm is needed or mixture distribution optimization is required as an application, MixEst can be used. 
\end{remark}

\subsection{Important Textbooks about Riemannian Geometry and Optimization}

There are many good textbooks and references about Riemannian geometry and Riemannian optimization.
Some of the important books and references are the following:

\begin{itemize}
\item \textbf{John M. Lee's books on topology and manifolds:}
\begin{itemize}
\item \emph{Introduction to Topological Manifolds}
\cite{lee2010topological},
\item \emph{Introduction to Smooth Manifolds}
\cite{lee2013smooth}.
\item \emph{Introduction to Riemannian manifolds} \cite{lee2018riemannian}
\item \emph{Riemannian manifolds: an introduction to curvature} \cite{lee2006riemannian2}
\end{itemize}

\item \textbf{Two important books on Riemannian optimization:}
\begin{itemize}
\item \emph{Optimization Algorithms on Matrix Manifolds}
by Pierre-Antoine Absil, Mahony, and Sepulchre \cite{absil2008optimization}.
\item \emph{An Introduction to Optimization on Smooth
Manifolds} by Nicolas Boumal \cite{boumal2023introduction}.
\end{itemize}

\item \textbf{Important references on matrix manifolds:}
\begin{itemize}
    \item Paper \cite{edelman1998geometry} for Stiefel and Grassmann manifolds.
    \item The book \cite{bhatia2009positive} for SPD manifold. 
    \item The book \cite{golub2013matrix} for matrix calculations in general.
\end{itemize}

\end{itemize}

There is also a survey and brief introduction to manifold optimization by Hu et al. \cite{hu2020brief}.

\begin{remark}[Role of these references in this monograph]
The books by Lee provide background on manifolds and smooth
geometry, while the books by Absil et al. and Boumal are
among the most important references specifically for
Riemannian optimization. The software papers cited in Section \ref{section_important_toolboxes_riemannian_optimization} are useful
because they show how the geometric objects derived in this
monograph are translated into practical numerical tools.
\end{remark}

\section{Conclusion}\label{section_conclusion}

In this monograph, we developed a detailed and self-contained
treatment of the foundations of Riemannian geometry for
Riemannian optimization. Starting from the basic notions
of topology, smooth manifolds, tangent and cotangent spaces,
tensor calculus, metric tensors, and connections, we gradually
built the geometric machinery required for optimization on
nonlinear spaces. We then connected these foundations to
optimization by deriving the Riemannian gradient, the
Riemannian Hessian, geodesics, exponential and logarithm
maps, retractions, and vector transport.

A central goal of this monograph was to bridge the gap between
abstract geometric theory and practical implementation.
Although many of the presented results are classical, they are
often stated in the literature in compact coordinate-free form
or with omitted intermediate steps. In contrast, the focus here
was on transparent derivations, explicit formulas, and
implementation-oriented development. This level of detail is
especially important for researchers and practitioners who wish
to design, analyze, or implement Riemannian optimization
algorithms from first principles.

We further specialized the general theory to several important
matrix manifolds arising in optimization and machine learning.
In particular, we studied the Stiefel manifold, the Grassmann
manifold, and the symmetric positive definite manifold
$\mathbb{S}_{++}^n$. For these manifolds, we derived their
tangent and normal spaces, metric tensors, Levi-Civita
connections, Riemannian gradients, Riemannian Hessians,
geodesic equations, exponential maps, retraction maps, and
vector transport operators. We also discussed differentiated
retractions as practical tools for numerical vector transport.
These manifolds appear in a wide range of applications,
including orthogonality-constrained optimization, subspace
learning, covariance estimation, diffusion models, geometric
signal processing, and geometric machine learning.

Another contribution of this monograph is the unification of
coordinate-based and geometric viewpoints. Throughout the
text, we emphasized both the coordinate-free meaning of the
main objects and their explicit coordinate or matrix
representations. This dual viewpoint is valuable because it
clarifies the geometric interpretation of the formulas while also
making them suitable for numerical computation. In this way,
the monograph may serve both as a theoretical reference and
as a practical guide for implementation.

The material presented here also highlights a broader message:
optimization on manifolds is not merely Euclidean optimization
with constraints added afterward. Rather, the geometry of the
search space fundamentally determines the notions of direction,
distance, curvature, acceleration, and feasible motion. Once the
underlying manifold structure is taken seriously, many
optimization procedures admit natural geometric analogues,
such as Riemannian gradient descent, Riemannian Newton's
method, retraction-based algorithms, and transported search
directions.

There are several natural directions for future development.
One direction is to extend the scope of the monograph to other
important manifolds and quotient manifolds arising in modern
applications \citep{absil2008optimization,boumal2023introduction}.
Another direction is to develop more advanced optimization
methods, including trust-region methods, conjugate-gradient
methods, stochastic variance-reduced methods, and large-scale
second-order methods on manifolds
\citep{absil2007trust,absil2008optimization,
kasai2016riemannian,sato2019riemannian,zhang2016riemannian,huang2015broyden}.
Riemannian coordinate descent \cite{han2024riemannian,hamed2024riemannian,huang2021riemannian} can also be of interest.
It is also of interest to study geometric structures beyond the
Riemannian setting, such as Finsler, symplectic, or information
geometries, when they are relevant to optimization and
learning \citep{francca2021optimization,betancourt2018symplectic,
amari2016information,asanjarani2021finsler}.
In addition, deeper connections with applications in machine
learning, computer vision, robotics, control, and scientific
computing can further motivate new geometric models and
numerical methods \citep{boumal2023introduction,fei2025survey,
saveriano2023learning}.

In conclusion, Riemannian geometry provides a rigorous and
powerful language for optimization on nonlinear spaces, while
Riemannian optimization provides a computationally
meaningful context in which the geometric objects acquire
algorithmic significance. By assembling the fundamental
concepts, derivations, and matrix-manifold specializations into
a unified exposition, this monograph aims to make the subject
more accessible, more explicit, and more useful for both study
and practice.

\section*{Dedication}

I dedicate this paper to my lovely wife, Golbahar Amanpour. I am grateful for having her in my life. 

\section*{Acknowledgment}

Regarding the concepts in this monograph, I thank the following people:
\begin{itemize}
\item I thank Ming Miao. At first, when I self-studied about Riemannian geometry, I was confused and frightened! At that time (I think it was 2019), Ming was a student in applied mathematics in the University of Waterloo. In an informal meeting with him in the university, Ming was the first person who explained to me some of the initial theories of Riemannian geometry and helped me understand some concepts better.

\item I thank Reza Godaz and Prof. Reshad Hosseini who introduced me to the field of Riemannian optimization and concepts such as retraction and vector transport. In 2021, we collaborated and proposed a paper for vector-transport-free Riemannian LBFGS for optimization on SPD manifolds \cite{godaz2021vector}. 
Reza Godaz also has a paper on using Riemannian optimization for fast singular value decomposition \cite{godaz2021accurate}.
Prof. Hosseini, at the University of Tehran, has had important articles on Riemannian optimization especially on SPD manifolds, in collaboration with Prof. Suvrit Sra at MIT (e.g., see \cite{sra2013geometric,sra2015conic,sra2016geometric,hosseini2015matrix,hosseini2020alternative,hosseini2020recent,zadeh2016geometric,sra2016positive}). 
In 2017, Prof. Hosseini was also the external thesis examiner for my master's degree at Sharif University of Technology (The resulting publication from my master's thesis was \cite{ghojogh2017fisherposes}, should you be interested). 

\item I thank Prof. Spiro Karigiannis at the Department of Pure Mathematics in the University of Waterloo. During my post-doctoral fellowship in the University of Waterloo in 2021, he kindly let me sit in his ``smooth manifolds" class. 
One of my initial trials in differential geometry was \cite{ghojogh2022manifold} after taking that class.
I learned a lot about topology and smooth manifolds through his course. During that course, I learned important concepts from the books of Prof. John M. Lee \cite{lee2010topological,lee2013smooth}.

\item I thank Prof. N. J. Wildberger. I learned some differential geometry concepts from his YouTube channel ``Insights into Mathematics": \url{https://www.youtube.com/@njwildberger}

\item I thank Prof. Leonard Susskind. I learned many differential geometry concepts---such as tensor algebra, covariant and contravariant components, geodesics, parallel transport, and Riemannian curvature---from his book \cite{susskind2025general} and his YouTube videos\footnote{\url{https://www.youtube.com/playlist?list=PL6i60qoDQhQGaGbbg-4aSwXJvxOqO6o5e}} (his ``The Theoretical Minimum" courses on physics).

\item I thank Prof. Masud Naaseri. I learned a lot about tensor algebra and covariant and contravariant components from his very useful YouTube videos\footnote{He has useful videos in Persian. For example, two of his videos on covariant and contravariant components are \url{https://www.youtube.com/watch?v=4ToZR1tTunM} and \url{https://www.youtube.com/watch?v=_P26pZFProI}.}. He has also published related books in Persian language \cite{naaseri2025einsten}. 
I was honored to once visit Prof. Naaseri and talk to him in 
Amirkabir University of Technology (Tehran Polytechnic) in 2011, when I started my university studies.

\item I thank Prof. Pierre-Antoine Absil, Prof. Nicolas Boumal, and their coauthors for their valuable books on Riemannian optimization \cite{absil2008optimization,boumal2023introduction}. 
I also acknowledge the paper \cite{edelman1998geometry} for Stiefel and Grassmann manifolds and the book \cite{bhatia2009positive} for SPD manifold. 
\end{itemize}

I have another monograph about Euclidean optimization \cite{ghojogh2021kkt,ghojogh2023background}. The reader may consult that paper to learn Euclidean optimization. I also have tutorial videos on optimization, based on my tutorial paper \cite{ghojogh2021kkt}, on YouTube: \url{https://www.youtube.com/playlist?list=PLPrxGIUWsqP3ZBM4Zy5YqfCh1BqM5sJov}.

I also have a brief video on YouTube about Riemannian optimization (entitled ``Introduction to Riemannian Optimization for Optimization on Riemannian Matrix Manifolds" by Benyamin Ghojogh): \url{https://www.youtube.com/watch?v=NZSW-oqDq9A}. The reader may also consult this video.


\bibliography{References}
\bibliographystyle{icml2016}

\end{document}